\documentclass[11pt,oneside]{article}
\usepackage[a4paper,margin=24mm,headheight=15pt]{geometry}
\usepackage{fontspec}
\usepackage{amsmath,amssymb,amsthm,mathtools,mathrsfs,bm}
\usepackage{booktabs,longtable,tabularx,array,multirow,enumitem}
\usepackage{graphicx,pdflscape,adjustbox,calc}
\usepackage{xcolor}
\usepackage{tikz}
\usetikzlibrary{arrows.meta,positioning,calc}
\usepackage{needspace}
\definecolor{inkblue}{HTML}{23415B}
\usepackage{microtype}
\usepackage{fancyhdr}
\IfFileExists{xurl.sty}{\usepackage{xurl}}{\usepackage{url}}
\usepackage{hyperref}
\hypersetup{colorlinks=true,linkcolor=inkblue,citecolor=inkblue,urlcolor=inkblue,bookmarksdepth=2,
  hypertexnames=false,pdftitle={The Moduli of Canonical Threefolds of Minimum Volume},
  pdfsubject={Equality moduli and two rigid canonical threefolds},
  pdfauthor={Pinxian Bie, Peien Du, Zhengjie Yu}}
\newtheorem{theorem}{Theorem}[section]
\newtheorem{proposition}[theorem]{Proposition}
\newtheorem{lemma}[theorem]{Lemma}
\newtheorem{corollary}[theorem]{Corollary}
\theoremstyle{definition}

\theoremstyle{remark}
\newtheorem{remark}[theorem]{Remark}

\newcommand{\CC}{\mathbb C}\newcommand{\QQ}{\mathbb Q}
\newcommand{\ZZ}{\mathbb Z}\newcommand{\PP}{\mathbb P}
\newcommand{\OO}{\mathcal O}
\newcommand{\Ext}{\operatorname{Ext}}
\newcommand{\Proj}{\operatorname{Proj}}
\newcommand{\Hom}{\operatorname{Hom}}\newcommand{\Tor}{\operatorname{Tor}}

\newcommand{\Frac}{\operatorname{Frac}}

\newcommand{\Aut}{\operatorname{Aut}}
\allowdisplaybreaks[2]
\numberwithin{equation}{section}

\newcommand{\basket}[1]{\textnormal{\texttt{#1}}}
\usepackage{tocloft}
\title{\LARGE\bfseries The Moduli of Canonical Threefolds\\of Minimum Volume}
\author{Pinxian Bie \and Peien Du \and Zhengjie Yu}
\date{}

\begin{document}
\maketitle
\begin{abstract}
We determine the moduli space of complex projective canonical threefolds
of volume $1/420$. It consists of an irreducible rational component of
dimension $33$ and two isolated points, represented by the explicitly
constructed rigid threefolds Vega and Altair.
\end{abstract}
\noindent\textbf{2020 Mathematics Subject Classification.}
14J30, 14J10, 14D15, 14B07.\\
\textbf{Keywords.} Canonical volume, moduli, canonical ring, rigidity.
{\small\setlength{\parskip}{0pt}\tableofcontents}
\clearpage
\part{The moduli space and the two rigid models}
\section{Introduction}\label{ms:introduction}\label{sec:introduction}
The sharp lower bound for the canonical volume of a smooth complex
projective threefold of general type is $1/420$
\cite[Theorem~1.1]{LBVolume}. We determine the canonical models at
this bound and their moduli. Besides the familiar weighted hypersurfaces
$X_{46}\subset\PP(4,5,6,7,23)$, there are precisely two models.
They are isolated and rigid. We call them \emph{Vega} and \emph{Altair};
their first nonzero pluricanonical sections have degrees four and five,
respectively.

Throughout, a canonical model means a normal complex projective
threefold with canonical singularities and ample canonical divisor.
We write
\[
 R(Z,K_Z)=\bigoplus_{m\geq0}H^0(Z,\OO_Z(mK_Z)),
 \qquad P_m(Z)=h^0(Z,\OO_Z(mK_Z)),
\]
where the sheaves are reflexive powers. The volume on this model is
$K_Z^3$. Passing to the canonical model identifies the birational
classification of smooth threefolds at the bound with the isomorphism
classification below. We write $B(Z)$ for the Reid basket of a
crepant terminalization, with $(b,r)$ denoting the contribution of
$\frac1r(1,-1,b)$.

\begin{theorem}[Equality models]\label{ms:classification}\label{thm:classification}
Let $Z$ be a canonical model with $K_Z^3=1/420$. Then
$p_g(Z)=q(Z)=0$, and exactly one of the following holds.
\begin{enumerate}[label=\textup{(\roman*)}]
\item $\chi(\OO_Z)=1$ and $Z$ is a normal canonical hypersurface of
degree $46$ in $\PP(4,5,6,7,23)$.
\item $\chi(\OO_Z)=2$ and $Z\cong Z_4=\mathrm{Vega}$, with
\[
 B(Z)=\mathcal B_4
 =7(1,2)+2(1,3)+3(1,4)+(1,5)+2(2,5)+(3,7).
\]
\item $\chi(\OO_Z)=2$ and $Z\cong Z_5=\mathrm{Altair}$, with
\[
 B(Z)=\mathcal B_5
 =4(1,2)+2(1,3)+(1,5)+2(2,5)+(3,7)+2(3,8).
\]
\end{enumerate}
Conversely, every normal canonical member of the indicated hypersurface
family and both exceptional models have canonical volume $1/420$.
\end{theorem}

The two exceptional models are the principal geometric feature of the
equality locus. We describe them before stating the moduli theorem.

\subsection*{The two rigid models}
For $e=4,5$, write $Z_e$ for Vega and Altair, respectively, and let
$a_e$ be their first nonzero canonical section. The surface
$D_e=(a_e=0)$ is normal, with klt singularities, and its Weil
polarization $H_e=K_{Z_e}|_{D_e}$ satisfies
\[
 K_{D_4}\sim5H_4,\quad H_4^2=\frac1{105},
 \qquad
 K_{D_5}\sim6H_5,\quad H_5^2=\frac1{84}.
\]
In each case the minimal resolution of $D_e$ is obtained by blowing up
the triple point of a type-IV fibre on an elliptic K3 surface.
The marked curves and the polarization distinguish the two surfaces.
The threefold construction specifies how their section algebras extend
in the direction of $a_e$.

For Vega the underlying elliptic K3 surface has the explicit equation
\begin{equation}\label{ms:intro-vega-k3}
 y^2+\frac{9t^2-8}{8}xy
    +\frac{27}{256}t^2(t^2-1)^2y
 =x^3+\frac{27}{256}t^2(t^2-1)x^2.
\end{equation}
Its singular fibres are of types
$\mathrm{IV},I_8,2I_5,2I_1$; the section $P=(0,0)$ has height $1/10$.
The sections $O,P,[5]P$ and the marked fibre components determine the
polarization after the blowup and contraction. For Altair a convenient
explicit description begins with the even part of the surface algebra:
\begin{equation}\label{ms:intro-altair-surface}
\begin{split}
 T=\CC[\theta_6,b_8,z_{14},P_{16},q_{18}]/\bigl(
 &P^2-zq-\tfrac14b^4+\tfrac3{44}\theta^4b,\\
 &q^2-bz^2+\tfrac8{11}\theta zP
       +\tfrac7{11}\theta^2b^3+\tfrac4{121}\theta^6\bigr).
\end{split}
\end{equation}
The full surface algebra is $T\oplus N$, where, as a module over
$A_0=\CC[\theta_6,b_8,z_{14}]$, $N$ has generators
$e_{15},J_{21},K_{23},f_{19},L_{25},M_{27}$ and relations
\[
 ze-bJ+\theta K=0,\qquad zf-bL+\theta M=0.
\]
Its multiplication, as well as Vega's surface algebra, is specified
in Section~\ref{ms:constructions} and Appendix~\ref{app:coefficients}. These are surfaces with rational
singularities whose minimal resolutions have K3 minimal models.

The ambient canonical rings have a parallel form. Let $S_e$ be the
explicit graded algebra, free of rank eight over $A_e$, with parameter
and basis degrees as follows:
\begin{center}
\small
\begin{tabular}{@{}lll@{}}
\toprule
 & $A_e$ & Degrees of an $A_e$-basis of $S_e$\\
\midrule
Vega & $\CC[a_4,u_{10},b_{12},c_{14}]$
 & $0,15,18,20,21,23,26,41$\\
Altair & $\CC[a_5,b_8,u_{12},z_{14}]$
 & $0,15,18,19,21,22,25,40$\\
\bottomrule
\end{tabular}
\end{center}
The products are fixed rational polynomials, not free parameters.
Within these algebras put
\begin{align}
 F_4={}&u^3-\frac{108}{25}a^2ub+\frac{432}{125}a^3q
          -\frac{972}{625}a^4c+\frac{9983169}{100000}a^5u,
          \label{ms:intro-F4}\\
 F_5={}&b^2u+abe+\frac14a^2q+\frac{2207}{4400}a^4b.
          \label{ms:intro-F5}
\end{align}
Here $q$ has degree $18$ in both algebras, and the $e$ in $F_5$ is
the degree-$15$ basis element. Then the precise construction is
\begin{equation}\label{ms:intro-rings}
 R_e=\overline{S_e}^{\,\Frac(S_e)(\sqrt{F_e})},
 \qquad Z_e=\Proj R_e,
\end{equation}
where the bar means integral closure in the indicated quadratic field.
The square roots have degrees $15$ and $14$, respectively.
Section~\ref{ms:constructions} gives the multiplication rules and
the colon-ideal description that makes this normalization explicit.
In particular, \eqref{ms:intro-rings} specifies the entire canonical ring.

Both rings have rank sixteen over their parameter rings, but have
depth two. The restriction of the canonical ring to $D_e$ misses one
section of its full surface section algebra: degree five for Vega and
degree six for Altair. This defect accounts for
\[
 h^2(Z_e,\OO_{Z_e})=1,\qquad \chi(\OO_{Z_e})=2.
\]
It is essential to retain it in the threefold construction.
The local geometries also differ. Altair has only isolated terminal
cyclic quotient singularities, whereas Vega has one additional
nonterminal canonical point. On a Gorenstein chart cover of that point,
the exceptional surface of a crepant extraction is
$\PP^1\times\PP^1$ with normal bundle $\OO(-2,-2)$.
The canonical Cartier indices of Vega and Altair are $420$ and $840$.

\subsection*{The moduli space}
Let $\mathcal M_{1/420}$ denote the coarse moduli space of canonical
models with volume $1/420$, with the usual admissible
$\QQ$-Gorenstein families. The following describes this canonical locus;
no compactification by noncanonical or nonnormal varieties is included.

\begin{theorem}[The minimum-volume moduli space]\label{ms:main-moduli}
There is a disjoint union
\[
 \mathcal M_{1/420}
 =\mathcal M_{46}\ \sqcup\ \{[Z_4]\}\ \sqcup\ \{[Z_5]\}.
\]
The space $\mathcal M_{46}$ parametrizes all normal canonical
hypersurfaces of degree $46$ in $\PP(4,5,6,7,23)$. It is irreducible
and rational of dimension $33$, and has finite quotient singularities.
Each member has a smooth ordinary miniversal deformation base of
dimension $33$, and all such deformations arise by varying its equation.
The other two components are reduced isolated points: every ordinary
deformation of $Z_4$ or $Z_5$ over a local Artinian $\CC$-algebra is
trivial.
\end{theorem}

Thus the hypersurface locus includes all its normal canonical members,
not only the quasismooth ones. The rigidity assertion concerns the
unmarked singular threefolds themselves. In the moduli stack the two
isolated components retain their finite automorphism groups; the
statement about two reduced points concerns the coarse space.

The positive-geometric-genus case has a different minimum volume,
$1/30$, and its equality models form a rational hypersurface moduli
space of dimension $103$ \cite[Theorems~1.1 and~1.4]{PGVolume}.
The two rigid components above occur in the geometric-genus-zero case.
Their stable pluricanonical birationality index is $23$
(Section~\ref{sec:stability-completion}); the bound $14$
for positive geometric genus in \cite[Theorem~1.1]{Yu14} does not apply
to them.

We give the two constructions first, then complete the equality
classification and describe the moduli of the hypersurface family.
All proofs and defining coefficients are included in the paper.

\section{Construction of the two exceptional models}
\label{ms:constructions}\label{sec:portraits}

We construct Vega and Altair as projective spectra of explicit normal
graded rings.  In both cases the minimal resolution of the first
canonical divisor has a K3 minimal model, and an involution makes its three-dimensional extension
manageable.  The invariant algebra is free of rank eight over a
polynomial parameter ring; the full canonical ring is a normal
quadratic extension, of rank sixteen over the parameter ring.  The full ring is not
Cohen--Macaulay.  Keeping this distinction is necessary both for the
construction and for the deformation calculation.

All subscripts on variables in this section denote degrees.  The
complete rational multiplication tables are printed in Appendix~\ref{app:coefficients}.  They are defining data, with rules
recovering every product, rather than numerical approximations.  Here
we give the surfaces, the parameter modules, the quadratic extensions,
and the arguments identifying their projective spectra as canonical
models.

\subsection{A normal quadratic extension}
\label{ms:normalization}

We will use the following elementary construction.  Let $S$ be a
normal finitely generated complex domain, and suppose
\[
 H^2=FG,\qquad F,G\ne0,
\]
where $F,G$ have no common height-one prime and $F$ is not a square in
$\operatorname{Frac}S$.  If $v^2=F$ and $I=(F):_S H$, then
\begin{equation}\label{ms:colon-cover}
 R=S\oplus I/v
\end{equation}
is the entire integral closure of $S$ in
$\operatorname{Frac}(S)(v)$.  Products in its second summand are
$ (x/v)(y/v)=xy/F$.

Indeed, at every height-one prime containing $F$, the element $G$ is
a unit, so $\operatorname{ord}F=2\operatorname{ord}H$.  The colon
condition implies $xy/F\in S$ by normality.  Conversely, trace splits
any integral element into an invariant part in $S$ and an
anti-invariant part $\xi$.  Both $v\xi$ and $(H/v)\xi$ are integral
invariants, and therefore belong to $S$; hence $v\xi\in I$.  This
also proves finiteness and normality.  The extension is \'etale at
height-one primes.  In the graded situation, if
$\omega_S\simeq S(1)$, reflexive pullback gives
$\omega_R\simeq R(1)$.

This statement concerns affine cones.  The singularities and canonical
polarization of $\operatorname{Proj}R$ still have to be checked, and
projective stabilizers must be included in that check.

\subsection{Vega: the marked K3 surface}
\label{ms:vega-surface}

Start with the elliptic K3 surface $T_4$ whose minimal Weierstrass
equation over $\mathbb P^1_t$ is
\begin{equation}\label{ms:vega-k3}
 y^2+\frac{9t^2-8}{8}xy+
 \frac{27}{256}t^2(t^2-1)^2y
 =x^3+\frac{27}{256}t^2(t^2-1)x^2.
\end{equation}
Let $O$ be its zero section and $P=(0,0)$.  Its discriminant is
\[
 -\frac{19683}{68719476736}\,
 t^8(t^2-1)^5(81t^2-80).
\]
The fibres are of types $\mathrm{IV}$ at infinity, $I_8$ at zero,
$I_5$ at $\pm1$, and $I_1$ at the two remaining roots.  The section
$P$ meets component $2$ at each reducible multiplicative fibre and
the identity component at $\mathrm{IV}$.  It generates the
Mordell--Weil group, with height $1/10$.

We specify the polarization, since the unpolarized K3 surface alone
does not determine the construction.  Write $F$ for the fibre class,
$B_1,B_2$ for the nonidentity components of the $\mathrm{IV}$ fibre,
$T_1,\ldots,T_7$ for the $I_8$ root chain, and
$Q_1^\pm,\ldots,Q_4^\pm$ for the two $I_5$ root chains.  Blow up the
triple point of the $\mathrm{IV}$ fibre, with exceptional curve $E$;
on the resulting surface $W_4$ put
\[
 C'=F-B_1-B_2-E,\qquad V_i=B_i-E,\qquad A=[5]P.
\]
Here $[5]P$ denotes multiplication in the Mordell--Weil group, not
five times the divisor $P$. Divisor classes from $T_4$ in the following formulas are pulled back
to $W_4$.  Define
\begin{align*}
 J={}&-3O-6F+4P+3T_1+5T_2+4T_3+4T_4+3T_5+2T_6+T_7\\
 &+3Q_1^++5Q_2^++4Q_3^++2Q_4^+
   +3Q_1^-+5Q_2^-+4Q_3^-+2Q_4^-,\qquad M=J-E.
\end{align*}
Order the nineteen curves $R_i$ in the blocks
\[
 (T_1,T_2,A,C',O);\quad V_1;\quad V_2;\quad
 (T_4,T_5,T_6,T_7);\quad(Q_i^+);\quad(Q_i^-),
\]
and give them respectively the coefficients
\[
 (t_i)=\frac{(4,1,5,2,1)}7;\quad\frac13;\quad\frac13;\quad
 \frac{(4,3,2,1)}5;\quad\frac{(3,1,4,2)}5;\quad
 \frac{(3,1,4,2)}5.
\]
The rational divisor $\mathcal H_4=M-\sum t_iR_i$ satisfies
\begin{equation}\label{ms:vega-polarization}
 \begin{gathered}
 5\mathcal H_4=E+\Delta_W,\qquad
 \Delta_W=\frac{T_1+2T_2+3A+4C'+2O}{7}
                  +\frac{V_1+V_2}{3},\\
 \mathcal H_4^2=\frac1{105},\qquad
 \mathcal H_4E=\frac1{21},\qquad
 \mathcal H_4F=\frac17,\qquad \mathcal H_4R_i=0.
 \end{gathered}
\end{equation}
These equalities follow from the fibre intersection matrices and the
section-height formula; the integral calculation is given in
Section~\ref{sec:exceptional-quartic-surface}.

The divisor $\mathcal H_4$ is nef and big.  To see its precise null
locus, use the effective expression in
\eqref{ms:vega-polarization}.  A further null curve would avoid
$E,C',V_1,V_2$, hence would be vertical since
$F=C'+V_1+V_2+3E$.  The remaining fibre components have positive
intersection, so the null locus consists of the nineteen $R_i$.
The pair $(W_4,\Delta_W)$ is klt, $105\mathcal H_4$ is Cartier, and
\[
 105\mathcal H_4-(K_{W_4}+\Delta_W)=100\mathcal H_4
\]
is nef and big. The basepoint-free theorem therefore contracts them to a normal klt surface
$D_4$.  The pushforward $H_4$ of $M$ is an ample integral Weil
divisor, with
\[
 K_{D_4}\sim5H_4,\qquad H_4^2=1/105.
\]
The exceptional graphs are $[2,2,2,3,2]$, two $[3]$, and three
$A_4$.  In particular these are rational singularities, although
the surface is birational to a K3 surface.  The integral root $H_4$
of $K_{D_4}$ is unique: the class-group calculation gives
$\operatorname{Cl}(D_4)\simeq\mathbb Z^2$ with no torsion.

Let $B_4=\bigoplus_{m\ge0}H^0(D_4,mH_4)$, and let $\theta_5$ be its
canonical section.  The base involution $t\mapsto-t$ descends to
$(D_4,H_4)$.  Choose its graded linearization so that $\theta$ is
anti-invariant.  The invariant section algebra is
\begin{equation}\label{ms:vega-surface-basis}
 U_4=\mathbb C[u_{10},b_{12},c_{14}]
 \langle1,h_{15},q_{18},r_{20},s_{21},k_{23},\ell_{26},t_{41}\rangle,
 \qquad u=\theta^2,
\end{equation}
a free module with $\omega_{U_4}=U_4(5)$.  The brackets indicate a
module basis, not freely adjoined algebra generators.  Some of the
products, displaying the surface structure, are
\[
 \begin{aligned}
 h^2&=bq-6ur-27u^3,&hq&=bs-2uk,\\
 q^2&=b^3-3ubc+2u\ell,&
 \ell^2&=u(c^3+b^2q)-u^2cq,&h\ell&=t-2u^2s.
 \end{aligned}
\]
The remaining products are obtained by setting $a=0$ in the complete
Vega table.  Thus this description retains the whole surface algebra,
including its polarization and involution.

\subsection{Vega: the three-dimensional extension}
\label{ms:vega-ring}

Put $A_4=\mathbb C[a_4,u_{10},b_{12},c_{14}]$.  Define $S_4$ to be
free over $A_4$ on the eight elements in
\eqref{ms:vega-surface-basis}, with the multiplication in Appendix~\ref{app:quartic-exact-coefficients}.  That appendix specifies
the twenty-one products among $h,q,r,s,k,\ell$ and recovers all
products with $t$ from $h\ell=t+\Psi$.  For example,
\[
 th=(h^2)\ell-h\Psi,\qquad
 tx=h(\ell x)-\Psi x\ (x=q,r,s,k),\qquad
 t^2=(th)\ell-t\Psi.
\]
These rules specify a finite algebra without a choice of any
undetermined coefficient.  Its reduction modulo $a$ is $U_4$.
The basis associators vanish and the Frobenius pairing has determinant
one.  Since the top basis element has degree $41$ and the parameter
weights sum to $40$, finite duality gives $\omega_{S_4}=S_4(1)$.
Normality lifts from $S_4/(a)=U_4$; in particular $S_4$ is a normal
domain.

The quadratic field is defined by $v_{15}^2=F_4$, where
\begin{equation}\label{ms:vega-cover-F}
 F_4=u^3-\frac{108}{25}a^2ub+\frac{432}{125}a^3q
       -\frac{972}{625}a^4c+\frac{9983169}{100000}a^5u.
\end{equation}
Equivalently, on setting
\[
 A'=\frac65a,\quad B'=b-\frac{369747}{16000}a^3,
 \quad Q'=q-\frac9{20}ac,
\]
one obtains $F_4=u^3-3(A')^2uB'+2(A')^3Q'$.
For an explicit colon presentation, use the degree-$36$ element
\begin{align*}
 H_4^{\rm alg}={}&u\ell+\frac{977}{240}au^2b
 -\frac{27}{10}a^2c^2-\frac{136323}{22400}a^3b^2
 +\frac{34749891}{512000}a^3uc\\
 &-\frac{18013827951}{40960000}a^4u^2
 +\frac{1295195440041}{2508800000}a^6b
 -\frac{14422585817323809}{4014080000000}a^9\\
 &+\frac65acq-\frac{109779}{4480}a^2uq
 -\frac65abr+\frac{367173}{40000}a^4r.
\end{align*}
The product table gives $(H_4^{\rm alg})^2=F_4G_4$ with
$G_4\in S_{4,42}$ and
$G_4\bmod a=c^3+b^2q-ucq$.  This identity determines $G_4$
uniquely in the domain $S_4$ and can be recovered by exact division
using the defining table.  The reductions of $F_4,G_4$ have no common
prime divisor, and the residue square root is the genuine surface
extension $\sqrt u$.  Thus \eqref{ms:colon-cover} gives
\begin{equation}\label{ms:vega-definition}
 R_4=S_4\oplus\frac{(F_4):H_4^{\rm alg}}{v_{15}},
 \qquad Z_4=\operatorname{Proj}R_4.
\end{equation}
We call $Z_4$ \emph{Vega}.

It is useful to see the complete module, rather than only the
quadratic field.  Its anti-invariant part $N_4$ has twelve minimal
generators of degrees
\[
 \mathcal D_4=(15,17,19,20,21,23,25,26,27,28,30,36)
\]
and resolution
\begin{equation}\label{ms:vega-resolution}
 0\longrightarrow A_4(-41)\longrightarrow
 \bigoplus_{d\in\{27,29,31,37,40\}}A_4(-d)
 \longrightarrow\bigoplus_{d\in\mathcal D_4}A_4(-d)
 \longrightarrow N_4\longrightarrow0.
\end{equation}
The first map is $(-c,b,-u,a,0)^{\mathsf t}$; the five relations
and twelve numerator generators are printed with the multiplication
table.  Reduction modulo $a$ gives
\begin{equation}\label{ms:vega-deficient}
 R_4/(a)=B_4\ominus\mathbb C\theta_5.
\end{equation}
Here the notation means the graded subalgebra obtained by deleting
this one-dimensional degree-five piece.  It is not the complete
section algebra of $D_4$, though its Proj is $D_4$.  Exactness of
\eqref{ms:vega-resolution} and the missing class are proved by
the explicit colon generators (Section~\ref{sec:quartic-module46}).  In particular the rank of
$R_4$ over $A_4$ is sixteen, and its polarization cube is
\[
 \frac{16}{4\cdot10\cdot12\cdot14}=\frac1{420}.
\]

\subsection{Altair: a complete surface algebra}
\label{ms:altair-surface}

The quintic surface has a particularly direct algebraic construction.
Put $A_0=\mathbb C[\theta_6,b_8,z_{14}]$ and
\begin{equation}\label{ms:altair-even-surface}
 \begin{split}
 T&=A_0[P_{16},q_{18}]/(f,g),\\
 f&=P^2-zq-\tfrac14b^4+\tfrac3{44}\theta^4b,\\
 g&=q^2-bz^2+\tfrac8{11}\theta zP
                   +\tfrac7{11}\theta^2b^3+\tfrac4{121}\theta^6.
 \end{split}
\end{equation}
This complete intersection is free over $A_0$ on $1,P,q,Pq$, and
$\omega_T=T(6)$.  To obtain the full surface algebra one must add
an odd module $N$, with generators
\[
 e_{15},J_{21},K_{23},f_{19},L_{25},M_{27}
\]
and exactly two $A_0$-linear relations
\begin{equation}\label{ms:altair-koszul}
 ze-bJ+\theta K=0,\qquad zf-bL+\theta M=0.
\end{equation}
Here the variable $f_{19}$ is distinct from the equation $f$ in
\eqref{ms:altair-even-surface}.  The module is the sum of two
Koszul cokernels, is torsion free of rank four, and is locally free
away from the parameter origin.

For completeness, the following compact rules determine its entire
multiplication.  On $b\ne0$ use the column basis $(e,K,f,M)$ and
let multiplication by $P,q$ be the matrices
\[
 \mathsf P=\begin{pmatrix}
 0&b^3/4&7\theta z/11&3b^2\theta^2/44-z^2\\
 b&0&3\theta^2/11&-7\theta z/11\\
 -\theta z/b&-b\theta^2/4+z^2/b&0&-b^3/4\\
 -\theta^2/b&\theta z/b&-b&0
 \end{pmatrix},
\]
\[
 \mathsf Q=\begin{pmatrix}
 \frac{4\theta^2z}{11b}&-\frac{\theta(b^2\theta^2+4z^2)}{11b}
                         &bz&-7b^3\theta/22\\
 \frac{4\theta^3}{11b}&-\frac{4\theta^2z}{11b}&14b\theta/11&-bz\\
 z&-b^2\theta/2&-\frac{4\theta^2z}{11b}
                         &\frac{\theta(b^2\theta^2+4z^2)}{11b}\\
 2\theta&-z&-\frac{4\theta^3}{11b}&\frac{4\theta^2z}{11b}
 \end{pmatrix}.
\]
Set
\[
 \mathsf H=\begin{pmatrix}0&0&0&-b\\0&0&b&0\\0&b&0&0\\-b&0&0&0\end{pmatrix},
 \qquad h(v,w)=v^{\mathsf t}\mathsf Hw.
\]
The product of two odd vectors is
\begin{equation}\label{ms:altair-odd-product}
 \begin{split}
 vw={}&h(\mathsf P\mathsf Qv,w)+\tfrac8{11}\theta z^2h(v,w)\\
 &+h(\mathsf Qv,w)P+h(\mathsf Pv,w)q+h(v,w)Pq.
 \end{split}
\end{equation}
Substitute $J=(ze+\theta K)/b$ and $L=(zf+\theta M)/b$ for the
six global generators.  Every resulting product is polynomial over
$A_0$, and the two matrices preserve the global module
\eqref{ms:altair-koszul}.  For example
\[
 e^2=b^2z-2b\theta P+\theta^2q.
\]
The identities $[\mathsf P,\mathsf Q]=0$,
$\mathsf P^{\mathsf t}\mathsf H=\mathsf H\mathsf P$,
$\mathsf Q^{\mathsf t}\mathsf H=\mathsf H\mathsf Q$, and
$f(\mathsf P,\mathsf Q)=g(\mathsf P,\mathsf Q)=0$, together with
the odd triple products, verify associativity.  Torsion freeness
extends these identities across $b=0$.  Thus the rules define the
whole graded algebra $B_5=T\oplus N$, not only its localization.

The algebra $T$ is normal by the Jacobian criterion in codimension
one.  The odd product identifies $N$ with
$\operatorname{Hom}_T(N,T)$, so $B_5/T$ is \'etale in codimension
one; both summands satisfy $S_2$.  Consequently $B_5$ is a normal
domain.  In fact its punctured cone is smooth: the singular locus
of the even algebra is resolved by the odd cover, whose transverse
cotangent space has dimension two at each of those points.  The
complete Jacobian calculation, including the loci $\theta=0$ and
$b=0$, is given in Section~\ref{sec:quintic-explicit-surface46}.

The surface $D_5=\operatorname{Proj}B_5$, with
$H_5=\mathcal O_{D_5}(1)$, has rational quotient singularities and
\[
 K_{D_5}\sim6H_5,\qquad H_5^2=1/84.
\]
Its minimal resolution is the blowup of an elliptic K3 surface at
the triple point of a type-$\mathrm{IV}$ fibre.  The exceptional
curve is the polar curve; its degree and marked fractions are
\[
 H_5E=1/14,\qquad (4/7,3/4,3/4).
\]
On the resolution, the relevant quotient chains are $[3,2,2]$ and
two $[2,3,2]$.  The exceptional curve of the blowup meets their
three $(-3)$ centres once and misses their $(-2)$ arms.  Blowing
it down recovers the three components of the type-$\mathrm{IV}$
fibre.  This describes the K3 geometry directly from the surface
algebra (Section~\ref{sec:quintic-polar-sources46}).

There is an involution fixing $b,z,q,e,f,J,L$ and negating
$\theta,P,K,M$.  With $u=\theta^2$ and $k=\theta P$, its invariant
algebra is free on
\begin{equation}\label{ms:altair-surface-invariants}
 U_5=\mathbb C[b_8,u_{12},z_{14}]
       \langle1,q_{18},k_{22},qk,e_{15},f_{19},J_{21},L_{25}\rangle.
\end{equation}
The algebra $B_5\ominus\mathbb C\theta_6$ is a subring; the
anti-invariant part is generated by
\[
 -b\theta,P,\theta^3,\theta z,\theta e,K,\theta q,\theta f,
 \theta J,M,
\]
together with $\theta L,Pq$ as parameter-module
generators.  The omitted class will account for the intermediate
cohomology of the threefold.

\subsection{Altair: the three-dimensional extension}
\label{ms:altair-ring}

Set $A_5=\mathbb C[a_5,b_8,u_{12},z_{14}]$.  Define $S_5$ to be
the free $A_5$-module on
\[
 1,q_{18},k_{22},t_{\flat,40},e_{15},f_{19},j_{21},\ell_{25},
\]
with the products in Appendix~\ref{app:quintic-exact-coefficients}.  The symbol
$t_\flat$ differs from $qk$ by the triangular expression printed
there; it is a basis convention, not an additional algebra
generator.  As for Vega, the products with the top element are
recovered from $qk=t_\flat+\Psi$.  The resulting algebra is
associative and $S_5/(a)=U_5$.  Its Frobenius pairing has top degree
$40$ and determinant one, so $\omega_{S_5}=S_5(1)$; normality lifts
from the central surface.

Define
\begin{equation}\label{ms:altair-cover}
 \begin{split}
 F_5={}&b^2u+abe+\tfrac14a^2q+\tfrac{2207}{4400}a^4b,\\
 H_5^{\rm alg}={}&-bk+\tfrac12a\ell-\tfrac{71}{440}a^2bu
                    -\tfrac1{220}a^3e-\tfrac1{48400}a^6,\\
 G_5={}&zq+\tfrac14b^4-\tfrac3{44}u^2b
       +a(\tfrac1{10}bf-\tfrac2{11}ue)\\
      &+\tfrac1{55}a^2k+\tfrac{9123}{4400}a^2bz
                       +\tfrac{203}{24200}a^4u.
 \end{split}
\end{equation}
The table verifies $(H_5^{\rm alg})^2=F_5G_5$.  On $a=0$ the
first and third expressions become $b^2\theta^2$ and $P^2$.
Their zero divisors have no common component on the even surface;
the square root $-b\theta$ gives its nontrivial quadratic field.
Applying \eqref{ms:colon-cover}, put
\begin{equation}\label{ms:altair-definition}
 w_{14}^2=F_5,\qquad
 R_5=S_5\oplus\frac{(F_5):H_5^{\rm alg}}{w_{14}},
 \qquad Z_5=\operatorname{Proj}R_5.
\end{equation}
We call $Z_5$ \emph{Altair}.  Its anti-invariant module $N_5$ has
generator degrees
\[
 \mathcal D_5=(14,16,18,20,21,23,24,25,27,27,31,34)
\]
and exact resolution
\begin{equation}\label{ms:altair-resolution}
 0\longrightarrow A_5(-40)\longrightarrow
 \bigoplus_{d\in\{26,28,32,35,39\}}A_5(-d)
 \longrightarrow\bigoplus_{d\in\mathcal D_5}A_5(-d)
 \longrightarrow N_5\longrightarrow0.
\end{equation}
The first map is $(z,u,-b,a,0)^{\mathsf t}$.  The numerator
generators and five first relations identify the entire colon
module.  In particular
\[
 R_5/(a)=B_5\ominus\mathbb C\theta_6,
 \qquad \frac{\operatorname{rank}_{A_5}R_5}{5\cdot8\cdot12\cdot14}
          =\frac{16}{6720}=\frac1{420}.
\]
The absence of the one extra class is checked by the injectivity of
$S_{5,20}\oplus S_{5,22}\to S_{5,50}$,
$(x,y)\mapsto H_5^{\rm alg}x-F_5y$; a six-by-six coefficient minor
has determinant $1/2$.  Thus no section is lost by replacing the
normalization with the explicit module presentation.

\subsection{Canonical polarization and singularities}
\label{ms:local-geometry}

For both constructions, dualizing the displayed module resolutions
gives
\[
 \operatorname{Ext}^2_{A_e}(N_e,\omega_{A_e})=\mathbb C(1),
 \qquad \operatorname{Ext}^1_{A_e}(N_e,\omega_{A_e})=\mathbb C.
\]
Once the local charts are checked, these yield
\begin{equation}\label{ms:cohomology}
 h^1(Z_e,\mathcal O_{Z_e}(m))=\mathbf1_{m=1},\qquad
 h^2(Z_e,\mathcal O_{Z_e}(m))=\mathbf1_{m=0}
 \quad(e=4,5).
\end{equation}
Thus $p_g=q=0$ and $\chi(\mathcal O)=2$.  Equivalently, the
homogeneous rings have depth two; this is compatible with the
projective varieties having Cohen--Macaulay singularities, since
the defect is supported at the cone vertex.

For Altair the entire punctured cone is smooth.  It is smooth along
$a=0$ by the surface construction.  On $a=1$ its singular locus
is finite, since a positive-dimensional projective component would
meet the ample boundary.  The global degree-six derivation of the
full normalized ring constructed in Section~\ref{sec:altair-cone-smoothness}, satisfies $D(a)=0$ and so
preserves this slice. It has no zero there, whereas
it must vanish on a finite reduced singular locus.  This proves
smoothness on the remaining chart.  Enumerating the residual
grading stabilizers gives the following complete quotient atlas:
\[
 \begin{array}{c|c|c}
 \text{number}&r&\text{tangent weights modulo }r\\\hline
 4&2&(1,1,1)\\
 2&3&(2,2,1)\\
 1&7&(5,6,2)\\
 2&8&(5,7,3)\\
 2&5&(3,2,4)\\
 1&5&(4,1,4).
 \end{array}
\]
Each is terminal: after reordering and choosing a generator the
weights have the form $(1,-1,b)$ with $(b,r)=1$.  They give
\[
 \mathcal B_5=4(1,2)+2(1,3)+(1,5)+2(2,5)+(3,7)+2(3,8).
\]
Its canonical Cartier index is $840$.

Vega has the ordinary terminal points
\[
 3(1,2)+2(1,3)+3(1,4)+(1,5)+2(2,5)+(3,7)
\]
and one further canonical point.  The latter has completed germ
\begin{equation}\label{ms:vega-node}
 \left(\frac{\mathbb C[[x_{00},x_{01},x_{10},x_{11}]]}
 { (x_{00}x_{11}-x_{01}x_{10})}\right)^{\langle g\rangle},
 \qquad
 g(x)=(ix_{00},-ix_{01},-ix_{10},ix_{11}).
\end{equation}
To explain its geometry, the quotient by $g^2=-I$ is the cone over
$\mathbb P^1\times\mathbb P^1$ embedded by $\mathcal O(2,2)$.
Blowing up its vertex is smooth and crepant, with exceptional
normal bundle $\mathcal O(-2,-2)$.  The residual involution has
four fixed points on that surface and characters $(-1,-1,-1)$ on
their two tangent directions and normal line.  Its quotient
therefore has four terminal $\frac12(1,1,1)$ points.  This is a
crepant projective extraction, so the terminalization basket is
\[
 \mathcal B_4=7(1,2)+2(1,3)+3(1,4)+(1,5)+2(2,5)+(3,7).
\]
The canonical Cartier index is $420$.  The full chart computations,
and the identification of the actual completed germ with
\eqref{ms:vega-node}, are in Sections~\ref{sec:quartic-geometry46} and~\ref{sec:vega-completed-germ}.

Finally, $\omega_{R_e}=R_e(1)$ and the chart calculations identify
$\omega_{Z_e}^{[m]}=\mathcal O_{Z_e}(m)$: there is no divisorial
grading correction.  Hence $K_{Z_e}$ is ample.  Normality gives
depth at least two, so the graded section map identifies
$R_{e,m}=H^0(Z_e,\mathcal O_{Z_e}(m))$ in every degree.  The rings
constructed above are therefore the full canonical rings, and the
rank calculations give $K_{Z_e}^3=1/420$.

\section{Recognition and ordinary rigidity}
\label{ms:recognition-rigidity}

The explicit constructions prove existence.  Two additional steps
are needed for the moduli statement: every realization of either
basket must give the corresponding ring, and its ordinary
deformations must be trivial.  The second statement does not follow
from uniqueness of the complex isomorphism class.

\subsection{Recognition from the first canonical divisor}
\label{ms:recognition}

For any realization of $\mathcal B_4$ or $\mathcal B_5$, complete
restriction and the nonzero two-form recover the first prime normal
canonical divisor, its integral polarization, and its marked
surface algebra.  In the quartic case the marked section-height
calculation forces the K3 equation \eqref{ms:vega-k3}; in the
quintic case the polar multiplication and the two Koszul summands
force \eqref{ms:altair-even-surface}--\eqref{ms:altair-odd-product}.
The symmetry extends to the ambient ring with its canonical
linearization.  These steps retain the full coherent restriction,
including the missing class; their proofs are given in Sections~\ref{sec:exceptional-quartic-surface},
\ref{sec:vega-ambient-symmetry},
\ref{sec:quintic-surface-recognition46}, and
\ref{sec:quintic-involution46}.

We record why the subsequent extension is unique; associativity of
the invariant ring alone would not suffice.  For Vega the negative
tangent dimensions of $U_4$ in degrees $-4,-8,-12$ are $3,1,1$,
and there are no later essential choices.  Compatibility with the
normal surface cover imposes one linear condition at first order.
In the resulting chart write the parameters as
$(\kappa,\delta,\eta)$.  The genuine polar supplies
anti-invariants of weights $15,21$, with squares $F,G$ and product
$H$ satisfying
\[
 H^2=FG,\quad F_0=u^3,\quad H_0=u\ell,\quad
 G_0=c^3+b^2q-ucq,\quad F-u^3\in a^2S_{4,22}.
\]
The coefficient $\beta$ of $a^2ub$ is nonzero.  Solving these
weighted coefficient equations, including the endpoint and zero
first-tangent cases, leaves
\[
 \kappa=\frac1{10},\qquad
 \delta=\frac{81843}{2800},\qquad
 \eta=\frac{1546119}{5600},\qquad
 \beta=-\frac{108}{25}.
\]
The other associative branch forces $\beta=0$ and is excluded by
the actual polar.  The resulting table and quadratic field are
exactly those of \eqref{ms:vega-definition} (Section~\ref{sec:quartic-recognition46}).

For Altair the first two negative tangent spaces of $U_5$, in
degrees $-5,-10$, each have dimension one; all later essential
spaces vanish.  The actual square identity forces the nonzero
first class, normalized by $\mu_1(e,e)=\ell$.  If $c$ denotes the
remaining second-order parameter, the order-five associativity
obstruction is
\[
 \frac9{550}\left(c+\frac{259}{880}\right)=0.
\]
It determines $c$, and all later orders lift uniquely modulo
homogeneous coordinate changes.  Finally $H^2=FG$ is a triangular
system for the seven remaining cover coefficients and gives
exactly \eqref{ms:altair-cover}.  This calculation includes the
vanishing loci of the parameters (Section~\ref{sec:quintic-recognition46}).

In either case, an arbitrary canonical ring is normal and finite
over this identified invariant algebra, with the identified
quadratic fraction field.  Uniqueness of integral closure now
identifies its whole ring with $R_4$ or $R_5$.  Thus the two
constructions exhaust all realizations of their baskets, up to
isomorphism of canonical models.

\subsection{The ordinary deformation spaces}
\label{ms:rigidity}

\begin{theorem}\label{ms:ordinary-rigidity}
For $e=4,5$, the canonical model $Z_e=\operatorname{Proj}R_e$
satisfies
\[
 T^1(R_e)_0=0,\qquad T^1_{Z_e}=H^1(Z_e,T_{Z_e})=0.
\]
Every ordinary deformation over a local Artinian complex algebra
is trivial.
\end{theorem}

\begin{proof}
First consider local deformations.  For an isolated normal complex singularity of
depth at least three, the cotangent tangent space is
$H^1(U,T_U)$ on the punctured spectrum.  A finite cover free there
computes this space by invariants.  Thus the isolated
three-dimensional quotient points in the atlases are locally
rigid.  At Vega's additional point, use the node cover in
\eqref{ms:vega-node}.  The node has a one-dimensional smoothing
space, but $g$ sends its equation to its negative, hence acts by
$-1$ on the constant smoothing direction.  Its invariant tangent
space is zero.  All local cotangent $T^1$ sheaves therefore vanish.

For the global calculation the parameter-module decompositions are
\[
 \begin{split}
 N_4&\simeq\Omega^2_{A_4}(\mathbb C)(-1)
       \oplus\Omega^3_{A_4}(\mathbb C)
       \oplus A_4(-20)\oplus A_4(-21),\\
 N_5&\simeq\Omega^2_{A_5}(\mathbb C)(-1)
       \oplus\Omega^3_{A_5}(\mathbb C)
       \oplus A_5(-16)\oplus A_5(-24).
 \end{split}
\]
Here $\Omega^j$ denotes the $j$th syzygy in the parameter Koszul
resolution.  The common module theorem and the displayed
resolutions give $\operatorname{Ext}^1_{A_e}(R_e,R_e)_0=0$.
Thus a graded first-order algebra deformation, after lifting its
four parameters, is represented by a multiplication deformation
of the full presented module.  One must still allow the parameter
frame to move.

For all pairs of the eight invariant and twelve anti-invariant
generators, allow every degree-compatible coefficient of a
commutative cochain $\phi$, with $\phi(1,x)=0$.  Impose the five
module relations and all linearized associativity equations
\[
 \phi(\mu(x,y),z)+\mu(\phi(x,y),z)
 -\phi(x,\mu(y,z))-\mu(x,\phi(y,z))=0.
\]
Coordinate changes include semilinear maps lifting derivations
$A_e\to R_e$.  In particular Altair has the anti-invariant frame
direction $D(z_{14})=w_{14}$; omitting it would compute a relative
deformation space, not ordinary deformations.

The exact coefficient matrices give
\[
 \begin{array}{c|c|r|r|r|r}
 &\text{parity}&\text{coefficients}&\text{constraint rank at least}
       &\text{kernel at most}&\text{coordinate changes}\\\hline
 Z_4&-&4452&4404&48&48\\
 Z_4&+&4632&4592&40&40\\
 Z_5&-&4188&4149&39&39\\
 Z_5&+&4260&4216&44&44
 \end{array}
\]
The rank bounds are certified by nonzero minors modulo $65521$,
where all relevant denominators are units; hence they are valid
over $\mathbb Q$.  The matching coordinate-change spaces are
constructed as exact rational coboundaries, and satisfy all
equations.  The entire tangent kernel is therefore the coboundary
space in both parities.  This proves $T^1(R_e)_0=0$.  The matrices
are determined by the complete printed multiplication and module
relations; the verification procedure and rational gauges are in
Section~\ref{sec:exceptional-global-rigidity46}.

Every ordinary first-order deformation of $Z_e$ is locally trivial
at its singularities, so its reflexive canonical powers are flat.
For $m\ge2$, canonical vanishing gives
$H^1(Z_e,\omega_{Z_e}^{[m]})=0$.  Degree one remains zero because
$p_g=0$, and degree zero is the scalar base because
$H^1(\mathcal O_{Z_e})=0$.  Taking all canonical sections therefore
gives a flat deformation of the full ring.  The ample canonical
class and relative Proj recover the original projective
deformation; hence $T^1_{Z_e}$ injects into $T^1(R_e)_0$.
Local-to-global deformation theory gives
$T^1_{Z_e}=H^1(Z_e,T_{Z_e})=0$, and induction over small Artinian
extensions proves the assertion.
\end{proof}

The theorem concerns the canonical models themselves and requires
no preserved involution, section, polarization frame, or module
splitting.  In particular the two exceptional moduli points are
reduced and isolated.  It makes no assertion about the deformation
space of a crepant terminalization of the nonterminal model Vega.

\section{Completeness of the equality classification}
\label{ms:sec:completeness}

We prove that the models constructed above exhaust the equality locus.
The detailed equality calculations are carried out in Parts~II and~III;
we give here the reductions and the precise points at which geometric
exclusions enter. A basket always means the Reid basket of a crepant
terminalization. Its constituents need not be actual quotient points.

\begin{theorem}\label{ms:thm:three-baskets}
Let $Z$ be a canonical model with $K_Z^3=1/420$. Then $p_g(Z)=q(Z)=0$,
and its Euler characteristic and terminalization basket are exactly one
of the following:
\[
\begin{array}{c|l}
\chi(\OO_Z)&B\\\hline
1&3(1,2)+(1,4)+(1,6)+(2,5)+(3,7),\\
2&7(1,2)+2(1,3)+3(1,4)+(1,5)+2(2,5)+(3,7),\\
2&4(1,2)+2(1,3)+(1,5)+2(2,5)+(3,7)+2(3,8).
\end{array}
\]
The last two rows are $\mathcal B_4$ and $\mathcal B_5$, respectively.
\end{theorem}

\subsection{The numerical range and the common criteria}
Write $P_m=h^0(Z,\omega_Z^{[m]})$ and $V=K_Z^3$.
The positive-geometric-genus theorem
\cite[Theorem~1.1]{PGVolume} gives $V\ge1/30$ if $p_g>0$.
The irregular-volume theorem \cite[Theorem~1.5]{JiangSeveri}
gives $V\ge3/8$ if $q>0$. Thus $p_g=q=0$ here, and the
holomorphic-form bound \cite[Theorem~2.9]{ChenJiangForms} restricts
$\chi(\OO_Z)$ to $1,2,3$. In particular
\[
 h^{2,0}(Z)=\chi(\OO_Z)-1.
\]
For $\chi=1$, the numerical equality statement
\cite[Theorem~1.2(2)]{CCII} gives the first basket in the theorem.

Since $1/420<1/336$, \cite[Corollary~1.5]{CCIII} gives
$P_2,\ldots,P_7\le1$. Riemann--Roch and multiplication of sections force
$P_2=0$. The canonical packing procedure of \cite[Section~3]{CCI},
with Reid's Riemann--Roch formula, gives a finite equality list.
Its completeness does not depend on a guessed index cutoff: the initial
basket has tail indices at most eight and total index sum at most
$26\chi+8\le86$, and prime packings preserve this sum.
Proposition~\ref{prop:equality-census} and Appendix~\ref{app:census}
give the resulting counts and every basket:
\[
\begin{array}{c|ccc}
\chi&1&2&3\\\hline
\text{numerical equality candidates}&1&65&283.
\end{array}
\]
These are necessary candidates, not assertions of existence.

We use the Hodge-null descent, invariant-pencil estimates and
whole-member trace results of \cite[Sections~4--5 and~8--11]{LBVolume}
under their stated hypotheses. In particular, Hodge-null descent requires
$q=0$, an invariant-pencil estimate requires a nonzero two-form
annihilating the pencil, and component trace uses a component of
multiplicity one. Their proofs are not repeated here.
The general polar inequality of \cite[Proposition~3.5]{Yu14} applies
to a terminal minimal threefold with nef and big canonical divisor,
without assuming $p_g>0$. It gives the useful necessary condition
\begin{equation}\label{ms:eq:polar-filter}
 P_{2e+2}-P_2\ge h^{2,0}(Z)
 \qquad\text{whenever }P_e>0.
\end{equation}

One equality specialization should be made explicit. For
$0\ne a\in R_e$, use the full coherent quotient
\[
 \mathcal Q_m(a)=
 \OO_Z(mK_Z)/a\OO_Z((m-e)K_Z),
\]
including all multiplicities. Vanishing and Serre duality give
\begin{equation}\label{ms:eq:whole-restriction}
 h^0(\mathcal Q_m(a))
 =P_m-P_{m-e}+(\chi-1)\mathbf1_{m=e+1},
 \qquad P_j=0\ (j<0).
\end{equation}
Indeed the sole possible connecting term is
$H^1(Z,\omega_Z)$, of dimension $\chi-1$; it occurs when $m=e+1$.
The quotient is pure of dimension two, and its additional trace
direction cannot be supported on a wholly invisible remainder
\cite[Lemmas~8.1, 8.8 and~9.4]{LBVolume}.
This is the point where replacing the full divisor by its reduced
support would lose information. Together with divisor-character
congruences, these restrictions provide the common tests recorded
in Sections~\ref{sec:restriction-toolbox} and~\ref{sec:uniform-exclusions}
below.

\subsection{The two-form branch and the remaining baskets}
For $\chi=3$ there are two independent two-forms. The equality-specific
polar test \eqref{ms:eq:polar-filter} already removes five baskets:
C3-00017, C3-00030, C3-00031 and C3-00048 have $P_4=P_{10}=1$,
while C3-49963 has $P_3=P_8=1$. In each case the displayed target
plurigenus must instead be at least two. These five rows overlap the
exclusions below and are not counted again. For uniform coverage, the
degeneracy calculation of Section~\ref{sec:chi3} excludes all $283$
candidates in disjoint groups of sizes
\[
 259+4+17+3.
\]
The first group fails the simultaneous character and Hodge tests;
the four exceptional coplanar cases force an effective torsion
translate of the canonical class; the last twenty have a cubic
section and are excluded by invariant pencils or the single-prime
character test. This is an equality calculation: the lower-bound
theorem for $\chi=3$ in \cite[Theorem~7.1]{LBVolume} alone would not
exclude this branch.

For the coplanar cases the small-order cover argument can be shortened.
If the torsion class has order $s$, the connected canonical cover has
$V'=s/420$, $p_g'\ge1$ and $2\chi'\ge s+5$, with $s\ge3$.
For $s\le13$, \cite[Theorem~1.1]{PGVolume} gives $V'\ge1/30$, a
contradiction. For $s=14$, its strict positive-Euler-characteristic
assertion gives $V'>1/30$, again a contradiction.
For $s\ge15$, the holomorphic-form estimate used in
Section~\ref{sec:chi3} gives
\[
 V'\ge\frac{h^{2,0}-3}{54}
 \ge\frac{s-1}{108}>\frac{s}{420}.
\]
Thus all orders are excluded. The exact character lists and the other
three groups are proved in Theorem~\ref{thm:no-chi3-equality46}.

When $\chi=2$, the common tests remove $37$ baskets, and the
visible-prime, torsion and complete-restriction arguments remove a
further $17$; see Proposition~\ref{prop:initial-cuts46} and
Section~\ref{sec:uniform-exclusions}. Hence
\[
 65\longrightarrow28\longrightarrow11.
\]
Two of these eleven are $\mathcal B_4$ and $\mathcal B_5$. For clarity, the other nine
and their geometric obstructions are listed below. The identifiers
refer to the complete multisets in Table~\ref{tab:residual11}.

\begin{center}
\small
\begin{tabular}{@{}lp{.48\textwidth}l@{}}
\toprule
Basket & Obstruction & Detailed proof\\
\midrule
C2-01825 & Adjacent pencils and the crepant graph &
Thm.~\ref{thm:c1825-same-pencil-extension}\\
C2-02075 & Complete moving-surface trace &
Thm.~\ref{thm:c2075-excluded46}\\
C2-02349 & The marked quartic K3 configuration &
Thm.~\ref{thm:2349-exclusion46}\\
C2-02350 & The full conductor of the prime octic &
Thm.~\ref{thm:octic-complete-exclusion46}\\
C2-02352 & The elliptic canonical-root cover &
Thm.~\ref{thm:c2352-elliptic-root-extension}\\
C2-02353 & The quintic polar and its transverse jet &
Thm.~\ref{thm:2353-exclusion46}\\
C2-02354 & Two coupled visible primes &
Thm.~\ref{thm:c2354-coupled34}\\
C2-02355 & Marked curves and full cancellation &
Thm.~\ref{thm:2355-exclusion46}\\
C2-02357 & Descent through the full conductor &
Thm.~\ref{thm:2357-exclusion46}\\
\bottomrule
\end{tabular}
\end{center}

The octic case requires the longest separate analysis.
Its first effective canonical member is a reduced visible
prime; its normalization has a K3 minimal model, but the member is
nonnormal. The conductor and polar bounds leave precisely the states
\[
\begin{gathered}
 (7,1),(7,2),(7,3),(7,4),(12,1),(14,1),(14,2),\\
 (18,1),(21,1),(25,1),(7,1)+(18,1).
\end{gathered}
\]
Here each pair records the degree and delta invariant of a conductor
component. Theorem~\ref{thm:octic-eleven-states46} proves completeness
of this list, with the necessary source markings. The local
source-character and conductor arguments in
Theorem~\ref{thm:octic-complete-exclusion46} exclude every state.
For the higher-delta degree-seven cases the full conductor algebra,
including its filtration and multiplication, is retained; the
intersection matrix alone is not the required obstruction.

\begin{proof}[Proof of Theorem~\ref{ms:thm:three-baskets}]
The numerical range and the complete census cover every equality
model. The preceding reductions remove the entire $\chi=3$ branch
and all but $\mathcal B_4,\mathcal B_5$ in the $\chi=2$ branch. The sole $\chi=1$
basket is the one displayed. Existence and recognition of the
corresponding models are established by the canonical-ring
constructions, rather than inferred from these numerical tests.
\end{proof}

\section{The moduli space at the minimum volume}
\label{ms:sec:moduli}

We consider canonical models over $\CC$: thus $Z$ is normal and
projective, has canonical singularities, and $K_Z$ is ample. Families
are Koll\'ar families, meaning that the reflexive powers of the relative
dualizing sheaf are flat and commute with base change. Write
$\mathcal M_{1/420}$ for the coarse moduli space of these models with
$K_Z^3=1/420$, and $\mathcal M_{46}$ for its locus with
$\chi(\OO_Z)=1$. The canonical condition is part of this notation;
no compactification is intended.

\subsection{The hypersurface component}

Put $A=\CC[x,y,z,t]$, with weights $(4,5,6,7)$, and set
\[
 U^{\mathrm{can}}=
 \{f\in A_{46}: X_f=(u^2=f)\subset\PP(4,5,6,7,23)
                  \text{ is normal canonical}\}.
\]
This is a nonempty open subset of the $37$-dimensional vector space
$A_{46}$. A general member is quasismooth; the parameter space also
includes every normal canonical member which is not quasismooth.

\begin{theorem}\label{ms:thm:hypersurface-moduli}
Every canonical model of volume $1/420$ and Euler characteristic one
is an $X_f$ as above, and its entire canonical ring is
\begin{equation}\label{ms:eq:hypersurface-ring}
 R(X_f,K_{X_f})=A[u]/(u^2-f),\qquad \deg u=23.
\end{equation}
Two such models $X_f$ and $X_g$ are isomorphic if and only if
\begin{equation}\label{ms:eq:isomorphism}
 g(\alpha x,\beta y,\gamma z,\delta t)=c^2f(x,y,z,t)
 \quad\text{for some }\alpha,\beta,\gamma,\delta,c\in\CC^*.
\end{equation}
Their coarse moduli space $\mathcal M_{46}$ is irreducible and rational
of dimension $33$, with finite quotient singularities. Its general
member has automorphism group $\mu_2$, generated by $u\mapsto-u$.
Every member has unobstructed ordinary Artinian deformations and a
smooth miniversal base of dimension $33$.
\end{theorem}

\begin{proof}
We first recall what is needed to pass from a basket to a canonical
ring. The equality basket in the Euler-one case is
\[
 3(1,2)+(1,4)+(1,6)+(2,5)+(3,7)
 \qquad\text{\cite[Theorem~1.2(2)]{CCII}}.
\]
Its first four nonzero section spaces, in degrees $4,5,6,7$, are
one-dimensional. Their generators $x,y,z,t$ are algebraically
independent: the proof uses the connected pencil in degree $10$ and
rules out a first relation by restriction to its general member
(Lemmas~\ref{lem:chi1-first-pencil}--\ref{lem:chi1-cutoff}).
This geometric step is necessary; the Hilbert series alone does not
exclude a relation accompanied by an additional generator.

For the resulting polynomial subring $A\subset R(Z,K_Z)$, let $d$
be the fraction-field degree. The multiplicity inequality gives
\[
 \frac1{420}\geq\frac{d}{4\cdot5\cdot6\cdot7}
                   =\frac{d}{840}.
\]
The case $d=1$ would make $Z$ rational, so $d=2$. Equality in this
inequality makes $R$ finite over $A$
(Lemma~\ref{lem:chi1-saturation}). Normality and the factoriality of
$A$ identify it with a quadratic cover $A[u]/(u^2-f)$
(Lemma~\ref{lem:chi1-quadratic-cover}). If $\deg u=b$, adjunction gives
$K_Z\sim_{\QQ}\OO_Z(b-22)$; the canonical grading then forces $b=23$.
Conversely, adjunction and the hypersurface cohomology calculation give
\[
 K_{X_f}=\OO_{X_f}(1),\qquad
 K_{X_f}^3=\frac{46}{4\cdot5\cdot6\cdot7\cdot23}=\frac1{420},
 \qquad \chi(\OO_{X_f})=1.
\]
They also identify the displayed ring with the full canonical ring.

An isomorphism of the canonical models induces a graded isomorphism
of their canonical rings. Since the pieces of degrees $4,5,6,7$ are
one-dimensional, it scales the four lower-weight generators. It sends
$u$ to $cu+h_{23}$, where $h_{23}\in A_{23}$; comparison of the
linear-$u$ terms forces $h_{23}=0$. This proves
\eqref{ms:eq:isomorphism}. The coordinate group is therefore
\[
 G=(\mathbb G_m)^5/\mathbb G_m^{(4,5,6,7,23)},\qquad \dim G=4.
\]
Its kernel on the coefficient space is the deck group $\mu_2$.
All stabilizers in $G$ are finite, since a positive-dimensional
stabilizer would give a nontrivial torus action on a projective variety
of general type. Thus the coarse quotient has dimension $37-4=33$.
It is the quotient in the separated canonical locus of the moduli
stack of stable varieties \cite{AHStable}.

Rationality can be seen directly. The exponent vectors of
\[
 z^3t^4,\qquad yzt^5,\qquad y^2z^6,\qquad xt^6
\]
form a basis of
$L_{46}=\{a\in\ZZ^4:(4,5,6,7)\cdot a\equiv0\pmod{46}\}$:
their determinant has absolute value $46$, the index of $L_{46}$.
On the open set where their coefficients are nonzero, the effective
torus $G/\mu_2$ normalizes all four coefficients to one. If a remaining
exponent is $a=\sum n_i a_i$, the corresponding quotient coordinate is
\[
 c_a\prod_{i=0}^{3}c_i^{-n_i}.
\]
These give $33$ independent parameters and a dense open subset of
$\mathcal M_{46}$ birational to $\mathbb A^{33}$. They also show that the
general stabilizer is exactly $\mu_2$. The deformation assertions are
proved below.
\end{proof}

\subsection{Ordinary deformations and the coarse local rings}

The coefficient count must describe all ordinary deformations, not
only deformations which preserve a chosen weighted presentation. We
give the argument because it determines the local structure of the
moduli space, including at nonquasismooth members.

At a non-Gorenstein point of $X_f$, the local canonical cover is a
three-dimensional Cohen--Macaulay Gorenstein hypersurface. The
non-Gorenstein locus is finite. Its cover is \'etale away from the
closed point and lifts uniquely over a nilpotent thickening of the
punctured spectrum. Depth three extends this lifted cover to a finite
flat algebra across the point: the required obstruction group is
$H^1(U',\OO_{U'})=0$. Taking the character summands shows that all
reflexive canonical powers are flat and commute with Artinian base
change. At Gorenstein points the same conclusion is immediate.
Thus every ordinary Artinian deformation is a Koll\'ar deformation.
The extension and base-change details are given in
Lemma~\ref{lem:chi1-index-cover-deformation}.

For every integer $m$, the hypersurface ring gives
\[
 H^1(X_f,\OO_{X_f}(m))=H^2(X_f,\OO_{X_f}(m))=0.
\]
Consequently the relative canonical section spaces of an Artinian
deformation are free and commute with reduction. Lift the five
generators of degrees $4,5,6,7,23$. Graded Nakayama and the Hilbert
function give a single relation of degree $46$ with invertible
$u^2$ coefficient. Completing the square writes the whole deformation
as $u^2=f_T$. Conversely every lift of the $37$ coefficients gives
such a deformation. Coefficient lifts exist across every small
extension, so the deformation functor is unobstructed.

The four infinitesimal coordinate directions are independent because
the stabilizer is finite. A formal slice therefore has dimension $33$;
equivalently,
\[
 T^1(X_f)=A_{46}/\langle xf_x,yf_y,zf_z,tf_t\rangle,
 \qquad \dim T^1(X_f)=33.
\]
Only finite automorphisms remain on the slice. Hence
\begin{equation}\label{ms:eq:coarse-local}
 \widehat{\OO}_{\mathcal M_{46},[X_f]}
 \simeq\CC[[s_1,\ldots,s_{33}]]^{\Aut(X_f)}.
\end{equation}
This proves the remaining assertions of
Theorem~\ref{ms:thm:hypersurface-moduli}. The deck involution acts
trivially on the parameters, so the coarse space is smooth at its
general point. Special members can have additional finite diagonal
automorphisms and a singular coarse moduli germ; their miniversal
deformation bases still have dimension $33$.

\subsection{The two isolated components}

For Vega and Altair the ordinary deformation calculation gives
\[
 T^1(Z_4)=H^1(Z_4,T_{Z_4})=0,
 \qquad T^1(Z_5)=H^1(Z_5,T_{Z_5})=0.
\]
Here local rigidity and the calculation for the full canonical ring
are both essential. In particular, the latter allows the polynomial
parameter frame to move and includes deformations which do not
preserve the involution
(Theorem~\ref{thm:exceptional-ordinary-rigidity46}). Induction across
small Artinian extensions gives ordinary formal rigidity. This
establishes reduced isolated moduli points, rather than merely a
classification of their complex isomorphism classes.

Together with the equality classification, we obtain the disjoint
union of coarse moduli spaces
\begin{equation}\label{ms:eq:total-moduli}
 \mathcal M_{1/420}
   =\mathcal M_{46}\ \sqcup\ \{[Z_4]\}\ \sqcup\ \{[Z_5]\}.
\end{equation}
The Euler characteristic is constant in a proper flat family, so the
Euler-one component cannot meet either Euler-two point within this
canonical locus. At the stack level the isolated components retain
their finite automorphism groups: they are
$B\Aut(Z_4)$ and $B\Aut(Z_5)$.

These statements concern canonical models. They neither classify
the noncanonical or nonnormal boundary of a stable compactification
nor assert properness of $\mathcal M_{1/420}$. Likewise, rigidity of
Vega does not assert rigidity of its crepant terminalizations.

\clearpage
\part{Common tools and the finite reduction}

\section{Equality range, divisor characters and the finite census}
\label{sec:inputs}
We now pass from the moduli statement to the numerical and geometric
conditions satisfied by a possible equality model. Throughout Parts~II
and~III, $Z$ denotes such a model. We use
\[
 P_m=h^0(Z,\omega_Z^{[m]}),\quad
 e=\min\{m>0:P_m>0\},\quad
 \delta=\min\{m:P_m\ge2\},\quad
 \chi=\chi(\OO_Z).
\]
Here $\omega_Z^{[m]}$ denotes the reflexive power. The usual definitions
and conventions for canonical, terminal and klt singularities are those
of \cite{KM}. All varieties are over $\mathbb C$.
The invariants above, and the canonical volume, are unchanged on passing
to a smooth projective resolution or a terminal minimal model.
In particular, $p_g=q=0$ implies $h^{2,0}=\chi-1$.

\subsection{Models, terminal baskets and divisor characters}
Choose a projective $\mathbb Q$-factorial terminal minimal model
$g:Y\to Z$, and put $K=K_Y=g^*K_Z$, $V=K^3$. The existence of the
canonical model and of suitable minimal models follows from \cite{BCHM};
crepant terminalization in this setting is also described in
\cite{Kawakita}. A \emph{terminal basket} is a multiset of pairs $(b,r)$,
where $0<b\le r/2$ and $\gcd(b,r)=1$, representing the local contribution
of $\frac1r(1,-1,b)$ to Reid's Riemann--Roch formula. A nonquotient
terminal singularity can contribute several virtual constituents. In
particular a basket is not, by definition, a list of actual quotient
points. We separately prove every use of an actual quotient chart.

Set $I=\operatorname{lcm}\{r:(b,r)\in B(Y)\}$ and $v=IV$. For an
integral Weil divisor $T$ on $Y$ put
\[
 n(T)=IK^2T,\qquad M(T)=IKT^2,\qquad
 h(T,U)=IKTU.
\]
These are integers: intersect with general sufficiently divisible
canonical hyperplanes avoiding the finite terminal singular locus, and
use two consecutive sufficiently large multiples. For a character
$i\in\mathbb Z/r$ define
\[
 \Phi_{b,r}(i)=\frac{\overline{bi}(r-\overline{bi})}{2r},
\]
where a bar is the least nonnegative residue. Reid's formula is
\begin{equation}\label{eq:RR46}
 P_m=\frac{m(m-1)(2m-1)}{12}V-(2m-1)\chi
       +\sum_{(b,r)\in B}\sum_{j=1}^{m-1}\Phi_{b,r}(j),
       \qquad m\ge2.
\end{equation}
The simultaneous Weil difference formulas, in the form established in
\cite[Section~2, equations~(2.1)--(2.4)]{LBVolume}, give
\begin{align}
 n(T)&\equiv\sum (I/r)b^2i_T\pmod I,\label{eq:degree-character46}\\
 M(T)&\equiv-\sum (I/r)\overline{bi_T}(r-\overline{bi_T})\pmod {2I},
       \label{eq:square-character46}\\
 h(T,U)&\equiv\sum (I/r)b^2i_Ti_U\pmod I.\label{eq:mixed-character46}
\end{align}
In a necessary-condition computation, allowing independent characters
at virtual constituents enlarges the set. In every mixed calculation
the character of the same divisor is retained; permutations of repeated
constituents must be matched in all possible ways.

\subsection{Two-forms, null divisors, and complete restriction}
Let $W\to Y$ be a smooth projective resolution. For a space
$\mathcal H\subset H^0(W,\Omega_W^2)$ and a divisor $T$, write
\[
 Q_T(\alpha,\beta)=\int_W[T]\wedge\alpha\wedge\overline\beta
\]
with the conventional positive factor. Pullbacks of rational Cartier
classes are understood. For an effective prime this is positive
semidefinite; its kernel is precisely the kernel of restriction to a
smooth model of the prime. Functorial reflexive pullback \cite{Kebekus}
shows that every canonically contracted prime is invisible. We call a
prime \emph{visible} to a nonzero form if its restriction is nonzero.
$Q_K$ is positive definite, by Kodaira's lemma and positivity.

We use the following precise consequences of the lower-bound proof.
\begin{lemma}\label{lem:common-tools46}\label{lem:complete-canonical-restriction46}
\begin{enumerate}
\item If $q=0$ and an integral Weil divisor $G$ on $Y$ satisfies
$K^2G=KG^2=0$, then $g_*G$ is torsion in $\operatorname{Cl}(Z)$,
Cartier at generic curve points, and $Q_G=0$.
\item If $T$ is visible and $S$ resolves its normalization, then for
$A=K|_S$ one has $A^2=n(T)/I$ and
$IAK_S\le n(T)+M(T)$. Consequently
$-n(T)\le M(T)\le n(T)^2/v$.
\item The full quotient sheaf
$\mathcal Q_m=\mathcal O_Z(mK_Z)/s\mathcal O_Z((m-e)K_Z)$,
for $0\ne s\in R_e$, is Cohen--Macaulay of dimension two. If
$p_g=q=0$, then, with $P_j=0$ for $j<0$,
\begin{equation}\label{eq:complete-restriction46}
 h^0(\mathcal Q_m)=P_m-P_{m-e}+(\chi-1)\mathbf1_{m=e+1}.
\end{equation}
The extra directions are dualizing-trace directions. A direction with
nonzero connecting class cannot be supported entirely on an invisible
divisor, even if that divisor is nonreduced.
\end{enumerate}
\end{lemma}
\begin{proof}
The first assertion is \cite[Proposition~4.1, Corollary~4.2 and Lemma~4.3]{LBVolume}: Hodge index on a
general canonical hyperplane makes $G$ numerically trivial there;
vanishing gives $\operatorname{Pic}^0=0$, and class-group Lefschetz
descends torsion. Negative definiteness of the exceptional Du Val
curves gives the codimension-two Cartier assertion. Normalization
adjunction with its effective conductor gives the second assertion.
For the last assertion we use the full quotient and trace formalism in
\cite[Lemmas~8.1, 8.8 and~9.4; Proposition~11.3]{LBVolume}, and retain
the calculation for arbitrary $h=\chi-1$.
Indeed, the cohomology sequence of
\[
0\longrightarrow\OO_Z((m-e)K_Z)\xrightarrow{\ s\ }
\OO_Z(mK_Z)\longrightarrow\mathcal Q_m\longrightarrow0
\]
has no positive-degree obstruction except at $m-e=1$.
For $j\ge2$, vanishing gives $H^1(Z,\OO_Z(jK_Z))=0$;
for $j\le0$ the same assertion follows from Serre duality and
vanishing, with $H^1(Z,\OO_Z)=0$ at $j=0$. At $j=1$,
Serre duality gives $h^1(Z,\omega_Z)=h^2(Z,\OO_Z)=\chi-1$.
Since $m=e+1\ge2$, the target $H^1(Z,\OO_Z(mK_Z))$ vanishes.
This proves \eqref{eq:complete-restriction46}. Locally on the canonical
cover the quotient is by a nonzerodivisor; taking eigensummands
preserves depth.
The invisibility statement uses the full dualizing trace of the proper
scheme. It does not replace that scheme by its reduced support.
\end{proof}

The canonical ring $R=\bigoplus_{m\ge0}R_m$ is a normal graded domain.
Indeed it is the intersection, inside the graded rational function
field, of the valuation conditions defining pluricanonical sections.
If $a\in R_u,b\in R_w$ satisfy $a^{w/d}=c b^{u/d}$, where
$d=\gcd(u,w)$ and $c\ne0$, Bezout's identity produces a homogeneous
fraction of degree $d$ with a power in $R$. Normality puts it in $R_d$.
We will refer to this as the integral-root argument.

\subsection{The complete equality census}\label{sec:equality-census}
\paragraph{The geometric range.}
The published irregular-volume bound is $V\ge3/8$ when $q>0$
\cite[Theorem~1.5]{JiangSeveri}. The positive-geometric-genus bounds
in \cite[Section~2]{ChenJiangForms}, and the holomorphic-form reduction
\cite[Theorem~2.9]{ChenJiangForms}, give
\[
 V=1/420\quad\Longrightarrow\quad p_g=q=0,
       \qquad \chi\in\{1,2,3\}.
\]
The equality statement in the cited preprint, outside $\chi=2,3$,
also gives the unique numerical
$\chi=1$ basket, in agreement with the published
\cite[Theorem~1.2(2)]{CCII}.

Define the \emph{pluricanonical section index}
$\delta=\min\{m:P_m\ge2\}$. Chen--Chen's published
\cite[Corollary~1.5, equivalently Corollary~4.8]{CCIII} states
$V<1/336\Rightarrow\delta\ge8$. Thus all equality models have
$P_2,\ldots,P_7\le1$. This strict numerical comparison is essential:
the finite search below is not a search over arbitrary threefolds.

\paragraph{Seeds and an intrinsic index bound.}
We recall the canonical unpacking and packing of \cite[\S3]{CCI}.
An initial basket consists of entries $(1,r)$ with multiplicities
$n_{1,r}$; its level-five basket is obtained by replacing
$\epsilon_5$ copies of $(1,2)+(1,3)$ by $(2,5)$.
The initial Riemann--Roch identities imply
\[
 0=-3P_2-P_3+P_4+P_5+P_6-P_7
             -n_{1,5}-2\sum_{r\ge6}n_{1,r}.
\]
If $P_2=1$, the binary bounds force
$P_4=P_5=P_6=1$ and $P_3=P_7=0$, contradicting multiplication
$R_2R_5\subset R_7$. Hence $P_2=0$. Put $n_5=n_{1,5}$.
For each binary tuple $(p_3,p_4,p_5,p_6,p_7)$ the seed data are
\begin{align}
 n_5+2\sum_{r\ge6}n_{1,r}&=-p_3+p_4+p_5+p_6-p_7,\label{eq:seedtail46}\\
 n_{1,2}&=5\chi-4p_3+p_4,\nonumber\\
 n_{1,3}&=4\chi+2p_3-3p_4+p_5,\nonumber\\
 n_{1,4}&=\chi+p_3+2p_4-p_5-\sum_{r\ge5}n_{1,r},\nonumber\\
 \epsilon_5&=2\chi-p_3+2p_5-p_6-\sum_{r\ge5}n_{1,r}.\label{eq:seed46}
\end{align}
Only nonnegative multiplicities after the level-five packing are
retained. The right side of \eqref{eq:seedtail46} is at most three,
so at most one tail constituent has $r\ge6$. Moreover
\begin{equation}\label{eq:seedvolume46}
60V(B^{(5)})=\chi-3p_3+p_5+2p_6-n_5
       +\sum_{r\ge6}(60/r-13)n_{1,r}.
\end{equation}
If a tail index were at least nine, its right side would be at most
$3+1+2+60/9-13=-1/3$. Thus the only tail indices are $6,7,8$.
Finally
\[
 \sum r(B^{(0)})=26\chi+2p_3+p_4-p_5+n_5
             +\sum_{r\ge6}(r-4)n_{1,r}\le26\chi+8\le86.
\]
A prime packing $(b,r)+(c,s)\mapsto(b+c,r+s)$, with
$|bs-cr|=1$, preserves this sum and decreases length. This proves
finiteness without an assumed final index cutoff. The lcm $I$ is
not asserted to be at most $86$.

\begin{proposition}[Equality census]\label{prop:equality-census}
Every canonical realization with $V=1/420$ occurs among the
$349$ baskets in Appendix~\ref{app:census}. There is one with
$\chi=1$, $65$ with $\chi=2$, and $283$ with $\chi=3$.
\end{proposition}
\begin{proof}
Generate the seeds by \eqref{eq:seedtail46}--\eqref{eq:seed46}.
From each, perform every prime packing with $r+s>5$, storing each
multiset once. The canonical sequence theorem of \cite{CCI} ensures
that every geometric basket is encountered. Stop a branch only when
$V<1/420$ or a tested $P_m$ is negative. Both failures persist under
packing. At a node use \eqref{eq:RR46}, the Miyaoka--Reid condition
\[
 \sum_{(b,r)}(r-r^{-1})\ge24\chi,
\]
and Hopf multiplication $P_{a+b}\ge P_a+P_b-1$ whenever
$P_aP_b>0$. Test these necessary conditions through degree $60$.
Failures of the last two conditions filter the node, never its
descendants. The complete counts are
\[
\begin{array}{c|r|r|r|r}
\chi&\text{seeds}&\text{visited with }V\ge1/420&\text{retained}&V=1/420\\\hline
1&7&61&36&1\\
2&25&2897&2401&65\\
3&35&61601&53789&283.
\end{array}
\]
For clarity, the volume at every node can be computed without
recursion as $6\chi-\sum b(r-b)/r$, since $P_2=0$.
The source implements both this whole-basket computation and an
independent traversal with the packing decrement $1/(rs(r+s))$.
Both give the displayed multisets. Appendix~\ref{app:census}
prints every equality basket. Finite tests are used only as necessary
conditions: no inference of existence is made from passing them.
\end{proof}

There is one additional use of the same generator later. With
$\chi=3$, $V=1/210$, $(P_2,\ldots,P_7)=(0,0,1,0,0,1)$ and
$P_{13}=0$, it returns exactly
\begin{align}
 B'_1&=8(1,2)+2(1,3)+4(1,4)+2(2,5)+(2,7)+2(5,12),\label{eq:cover-baskets46}\\
 B'_2&=7(1,2)+3(1,3)+5(1,4)+2(2,5)+2(3,7)+(5,12).\nonumber
\end{align}
The displayed low plurigenera already imply $\delta\ge8$, so this
application needs no volume-based section-index bound.



\section{Forms, visible divisors and pencils}
\label{sec:forms-pencils}

We begin with the geometric constraints that apply before the individual
baskets are considered. A holomorphic two-form detects the components on
which restriction is nonzero; its differential bracket relates those
components to low-degree canonical sections. We combine this observation
with the geometry of a complete pluricanonical system and with marked
adjunction. At the end of the section these tools exclude the entire
Euler-characteristic-three branch.

\subsection{The differential bracket and the obstruction to rational time}
\label{ag:sec:common-differential}

The differential bracket used in the equality argument has a useful
intrinsic formulation. We state it in arbitrary dimension, although its
applications to the classification will involve holomorphic two-forms
on threefolds.

\begin{proposition}[The canonical Poisson bracket]
\label{ag:prop:poisson}
Let $X$ be a normal projective $n$-fold and let
$\alpha\in H^0(X,\Omega_X^{[n-1]})$. Its canonical section ring
has an ordinary Poisson bracket of graded degree one. On a smooth
canonical frame, write
$a=f\nu^e$, $b=g\nu^m$ and $\alpha=\iota_v\nu$. Then
\begin{equation}\label{ag:eq:poisson}
 \{a,b\}_{\alpha}
 =\bigl(ef\,v(g)-mg\,v(f)\bigr)\nu^{e+m+1}.
\end{equation}
For homogeneous $a$, the Hamiltonian derivation
$D_{a,\alpha}=\{a,-\}_{\alpha}$ has degree $e+1$, satisfies
$D_{a,\alpha}(a)=0$, and obeys
\[
 [D_{a,\alpha},D_{b,\alpha}]
 =D_{\{a,b\}_{\alpha},\alpha}.
\]
The construction is linear in $\alpha$, and these Poisson brackets
are compatible.
\end{proposition}

\begin{proof}
Under a change of frame $\nu'=h\nu$, the coefficients become
$f'=h^{-e}f$, $g'=h^{-m}g$ and $v'=h^{-1}v$. Thus
\[
 ef'v'(g')-mg'v'(f')
 =h^{-(e+m+1)}\bigl(efv(g)-mgv(f)\bigr),
\]
because the two terms containing $v(h)$ cancel. The bracket is
therefore intrinsic on the smooth locus and extends as a section of
the indicated reflexive canonical power.

For the Jacobi identity, let $t$ be the fibre coordinate on
$K_X^{\vee}$ in this frame. Sections $f\nu^e$ correspond to
functions $ft^e$, and their bracket is induced by
\[
 \Pi=t^2\partial_t\wedge v=E\wedge(tv),
 \qquad E=t\partial_t.
\]
Since $[E,tv]=tv$, the Schouten square of $\Pi$ vanishes.
Skew-symmetry and Leibniz are immediate, and Jacobi gives the
commutator identity. The same construction applies to every linear
combination of forms, proving compatibility.
\end{proof}

Here the grading is the ordinary canonical grading; it introduces no
super-Poisson sign convention. In dimension three this is precisely
the differential bracket of \cite{LBVolume}. If $0\ne a\in R_e$
and $u\in R_m$ have orders $b,c$ along a prime visible to $\alpha$,
and $D_{a,\alpha}(u)=0$, comparison of the leading coefficient gives
\begin{equation}\label{ag:eq:visible-bracket-valuation}
 ec=mb.
\end{equation}
In odd dimension $n=2r+1$, a two-form of generic rank $2r$ supplies
an $(n-1)$-form $\alpha=\omega^r$ for the proposition. An arbitrary
two-form in higher dimension need not have a rank-one kernel.

\begin{lemma}[A visible component and an invisible remainder]
\label{lem:visible-invisible-bracket-division}
Let $Z$ be a canonical threefold, let $\alpha$ be a nonzero
holomorphic two-form on a smooth model, and suppose an actual section
$a\in R_m$ has divisor $T+E$, where $T$ is a reduced visible prime
and $E\ge0$ is supported on invisible primes, with arbitrary
multiplicities. For $0\ne s\in R_n$ put $l=\operatorname{ord}_T(s)$.
If $ml\ne n$, then $D_a(s)\ne0$ and
\[
 \operatorname{ord}_T D_a(s)=l,\qquad
 E\le\operatorname{div}D_a(s).
\]
If in addition $l\ge1$, then $D_a(s)/a\in R_{n+1}$ is nonzero
and has order $l-1$ along $T$.
\end{lemma}
\begin{proof}
Use the local expression \eqref{ag:eq:poisson}. At the generic point
of $T$, write $a=t u\nu^m$ and $s=t^l h\nu^n$, with $u,h$ units.
Visibility makes $v(t)$ a unit, so the leading coefficient of the
bracket is $(ml-n)uhv(t)t^l$. At an invisible prime, $v(t)\in(t)$;
hence $v$ preserves every power of $(t)$, and both terms of the
bracket vanish to at least the multiplicity of that prime in $E$.
For $l\ge1$ division by $a$ therefore creates no divisorial pole.
Normality extends the quotient as a reflexive canonical section.
\end{proof}

The global origin of the bracket imposes a further restriction.
The following no-rational-time statement is a special case of the
Riccati rigidity lemma \cite[Lemma~5.2]{LBVolume}. We recall its
dimension-independent proof.

\begin{lemma}[No rational time]
\label{ag:lem:no-time}
Let $W$ be a smooth projective $n$-fold, let
$0\ne\sigma\in H^0(W,K_W)$ and
$\alpha\in H^0(W,\Omega_W^{n-1})$, and define the rational vector
field $v$ by $\iota_v\sigma=\alpha$. There is no rational function
$t$ with $v(t)=1$. Consequently, for a rational function $u$,
$v^N(u)=0$ implies $v(u)=0$.
\end{lemma}

\begin{proof}
Holomorphic forms on $W$ are closed. If $v(t)=1$, then away from
the poles of $t$ we have
\[
 \sigma=dt\wedge\alpha=d(t\alpha).
\]
After resolving the boundary to simple normal crossings, degeneration
of the logarithmic Hodge--de Rham spectral sequence injects global
logarithmic top forms into the de Rham cohomology of the complement
\cite[Corollaire~3.2.14]{DeligneHodgeII}. It forbids the nonzero
holomorphic form $\sigma$ from being exact there. For the final
assertion, take the smallest $N\ge2$ with $v^N(u)=0$. The ratio
$v^{N-2}(u)/v^{N-1}(u)$ would be a rational time.
\end{proof}

For $0\ne a\in R_e$, pass to a smooth projective resolution of its
normalized canonical-root cover. The tautological top form $\sigma$
and the pulled-back form $\alpha$ are holomorphic, with
$\sigma^e=a$ in the rational canonical algebra. The tautological
form is already regular on the cyclic cover of the canonical line
bundle of a smooth model, and remains regular after pullback to a
resolution. If $u=g\sigma^m$, then
\[
 D_{a,\alpha}(u)=e\,v(g)\sigma^{m+e+1},
 \qquad \iota_v\sigma=\alpha.
\]
Lemma~\ref{ag:lem:no-time} therefore gives
$D_{a,\alpha}^N(u)=0\Rightarrow D_{a,\alpha}(u)=0$.
It also gives the quadratic consequence used below:
\begin{equation}\label{ag:eq:quadratic-derivation-obstruction}
 b\in R_{e+1},\quad D_{a,\alpha}(b)\in\CC b^2
 \quad\Longrightarrow\quad D_{a,\alpha}(b)=0.
\end{equation}
Indeed, after writing $b=g\sigma^{e+1}$, a nonzero equality
$D_{a,\alpha}(b)=c b^2$ would give $v(1/g)=-c/e$, and hence a
rational time. No extension of $v$ across its poles is required.

\begin{lemma}\label{lem:commuting46}
Suppose the seed $a\in R_e$ and sections of degrees $m_1,\ldots,m_t$
commute for one nonzero form. Put $d=\gcd(e,m_1,\ldots,m_t)$.
If two independent commuting monomials have degree $N$, then
\[
 V\ge\frac1{dN(N+1)}.
\]
In addition, if $j$ is the scaled canonical degree of a connected
general moving surface and $\theta$ the Stein degree, then
\[
 I/d\le j+j^2/v,\qquad 0<\theta j\le N(v-1/d).
\]
\end{lemma}
\begin{proof}
We combine the special-model comparison and invariant-pencil inequalities
of \cite[Proposition~5.6, Theorem~5.7 and Corollary~5.8]{LBVolume}
with the following valuation argument, which applies at equality volume.
A visible component of the seed has multiplicity divisible
by $e/d$. The valuation identity forces every chosen degree-$N$
monomial to have multiplicity at least $N/d$ there. On the special
pencil model, adjunction to the irregular general surface gives
$AL\ge1/d$, while $AL\le L^2+(L^2)^2/V$. The fixed component
has positive scaled degree at least one. Intersecting the fixed and
moving decomposition with $IK^2$ gives the second displayed
inequality, and the first follows from the adjunction comparison.
The corresponding volume estimate is the first assertion.
\end{proof}

\label{sec:finite-criteria46}
\subsection{Connected pencils and complete systems}
Resolve a pencil in $|mK_Y|$, let $F$ be its connected general fibre,
and let $\theta$ be the degree of the Stein base over the parameter line.
The Stein base is $\mathbb P^1$, since $q(Y)=0$. On $Y$ write
\[
 mK_Y\sim\theta F_Y+E,\qquad E\ge0,
\]
where $F_Y$ is the integral moving divisor. On the special terminal
pencil model of \cite[Proposition~5.6]{LBVolume}, put
\[
 L=K_Y|_F,\quad d=L^2,\quad A=\sigma^*K_{F_0},\quad
 k=A^2,\quad a=LA,\quad z=K_YF_Y^2.
\]
Here restrictions mean pullbacks, and $F_0$ is the smooth minimal
surface of $F$. In particular $k$ is a positive integer. The effective
exceptional comparison on that model gives
\begin{equation}\label{eq:pencil-area46}
 a\le LK_F\le d+z\le d+d^2/V,\qquad
 a^2\ge kd,\qquad d\le(m/\theta)V.
\end{equation}
These inequalities are preserved on higher resolutions. They are not
asserted by replacing the special comparison divisor on an arbitrary
resolution. If $LK_F\le u<1$, the fractional-null-curve argument of
\cite[Lemma~13.6 and Corollary~13.7]{LBVolume} strengthens Hodge index to
\begin{equation}\label{eq:fractional-area46}
 u^2\ge d(k+1/3).
\end{equation}
Thus, setting $r=m/\theta$, we have
\[
 k\le r(1+r)^2V,\qquad
 r(1+r)V<1\ \Longrightarrow\ k+\tfrac13\le r(1+r)^2V.
\]
In particular $\theta<m$ here, for $\theta\ge m$ would give
$1\le k\le4V<1$. For integer-grid calculations write $j=Id$ and
$M=Iz$. Then $j,M$ are integers, $M\ge0$, $vM\le j^2$, and,
when $(j+M)/I<1$,
\begin{equation}\label{eq:integer-pencil46}
 (j+M)^2\ge Ij(k+1/3).
\end{equation}

\begin{lemma}[The square of a complete moving surface]
\label{lem:complete-moving-square46}
If $P_a=s\ge3$ and the complete $a$-canonical image has dimension
at least two, a smooth general moving surface has minimal square
$k\ge2s-3$. Consequently, for any complete system with $P_a\ge3$,
\begin{equation}\label{eq:complete-system-volume46}
 V\ge\frac{2P_a-3}{a(a+1)^2}.
\end{equation}
\end{lemma}
\begin{proof}
On a smooth resolution let $M$ be the free moving divisor of the
complete system, $S\in|M|$ general, $C=M|_S$, and
$D=\sigma_*C$ on the minimal surface $S_0$. The integral divisor $D$
is nef. Adjunction gives $K_S-C=K_W|_S$, which is big because it
contains the restriction of the nef and big canonical pullback.
Hence $K_{S_0}-D$ is big and $k>K_{S_0}D$.

If the image has dimension two, $C\equiv cG$, where $G$ is a connected
curve fibre and $c\ge s-2$: restriction gives $h^0(S,C)\ge s-1$, whereas a degree-$c$ line bundle on a curve has
at most $c+1$ sections. More precisely this inequality is applied to
the line bundle on the Stein curve defining $|C|$. Put
$G_0=\sigma_*G$. If $K_{S_0}G_0=1$, then $k>c\ge1$ and Hodge index
forces $G_0^2=0$, contrary to the parity of $G_0^2+K_{S_0}G_0$.
Thus $K_{S_0}G_0\ge2$ and $k\ge2c+1\ge2s-3$.

If the image has dimension three, put $b=M^3>0$. Two restrictions
give $h^0(C,M|_C)\ge s-2$. Adjunction gives
$2g(C)-2=K_WM^2+2b$, with $K_WM^2$ a positive even integer.
The degree-$b$ line bundle $M|_C$ is special; Clifford's inequality
gives $b\ge2s-6$. Moreover $D^2\ge b$ and
$(K_{S_0}-D)D$ is a positive even integer. Thus
$K_{S_0}D\ge2s-4$ and $k\ge2s-3$ again.

Take a general subpencil and apply \eqref{eq:pencil-area46} with
$d\le aV$. It gives $k\le a(a+1)^2V$. If the complete image is a
curve instead, its moving degree satisfies $\theta\ge s-1$, and
\eqref{eq:pencil-area46} gives $V\ge\theta^3/[a(a+\theta)^2]$.
For $a\ge1$, $\theta\ge2$, this is at least
$(2\theta-1)/[a(a+1)^2]$; multiplying out, or minimizing
$(a+1)/(a+\theta)$ at $a=1$, proves the comparison.
\end{proof}

The Hilbert function of a nondegenerate integral surface in
$\mathbb P^{s-1}$ is at least
\begin{equation}\label{eq:hilbert-surface46}
 \binom{b+2}{2}+(s-3)\binom{b+1}{2}\qquad(b\ge1).
\end{equation}
To see this, take two successive general hyperplane sections and use
injectivity of multiplication by their equations. The final
zero-dimensional nondegenerate section has Hilbert function at least
$1,s-2,s-2,\ldots$; summing twice gives the displayed expression.
The same bound holds for a variety of higher dimension by restricting
to an integral surface section. Thus, if $P_{ba}$ is smaller than
\eqref{eq:hilbert-surface46}, the complete $a$-system is a pencil,
and necessarily
\begin{equation}\label{eq:hilbert-pencil-volume46}
 V\ge\frac{(P_a-1)^3}{a(a+P_a-1)^2}.
\end{equation}

\subsection{Fibre geography and full-member trace}
Write $(k,p,q)=(K_{F_0}^2,p_g(F_0),q(F_0))$. The standard surface
inequalities used in the lower-bound proof give
\[
 0\le q\le p\le\lfloor k/2\rfloor+2,
 \qquad q>0\Longrightarrow k\ge2p.
\]
Here $\chi(\mathcal O_{F_0})>0$. If $q=0$, restriction of absolute
two-forms is injective, so $p\ge h$. A nonzero form in the restriction kernel is locally
$dt\wedge\eta$; whether or not $p<h$, the kernel-form estimate
\cite[Theorem~5.7 and Corollary~5.8]{LBVolume} gives
\begin{equation}\label{eq:kernel-area46}
 p<h\Longrightarrow 2\le k\le r(1+2r)V.
\end{equation}
More generally, the same upper bound on $k$ holds whenever that
restriction has a nonzero kernel. Thus $k>r(1+2r)V$ guarantees
injectivity even for irregular fibres.

For completeness, the argument applies at $V=1/420$ as follows.
Let $\alpha$ be a nonzero kernel form and work on a common smooth
refinement $W$ of the special pencil model. On a general fibre write
$\alpha=dt\wedge\eta$, with $0\ne\eta\in H^0(F,\Omega_F^1)$.
If $K_W=\mu^*K_Y+\Delta$, then
\[
 L-A\equiv(K_F-A)-\Delta|_F.
\]
The one-form $\eta$ vanishes on the rational curves supporting
$K_F-A$. A vertical component of $\Delta$ misses the general fibre.
For a horizontal component, functorial reflexive pullback
\cite{Kebekus} makes $\alpha$ vanish on its smooth model, since its
image in $Y$ has dimension at most one. As $dt$ is nonzero there,
$\eta$ vanishes on the normalizations of its intersections with $F$.
Consequently the nonzero nef isotropic class
$B=[\sqrt{-1}\eta\wedge\overline\eta]$ satisfies $B(L-A)=0$.
Hodge index makes $B^\perp$ negative semidefinite, so
$(L-A)^2\le0$. Using \eqref{eq:pencil-area46} gives
\[
 k\le2a-d\le d+2d^2/V\le r(1+2r)V.
\]
Also $q(F)>0$, and the irregular-surface inequality with its
addendum \cite{BeauvilleDebarre1982} gives $k\ge2p\ge2$.
This proves the estimate without any strict assumption
$V<1/420$ and without a period-map assertion. In particular the
profile $(k,p,q)=(1,2,0)$ is subject only to the displayed numerical
and trace tests; it is not deleted by a general fibre prohibition.

If $P_{m+1}=0$, ordinary trace into the full member gives $p\le h$.

For $\ell\ge2$ put
\[
 U_{\ell,m}=P_{\ell(m+1)}-P_{\ell(m+1)-m},\qquad
 B_\ell=\binom\ell2 k+1+p-q.
\]
The pluri-trace and reduced-member defect
\cite[Theorems~8.2 and~8.5]{LBVolume} give
\begin{equation}\label{eq:profile-trace46}
 U_{\ell,m}\ge\theta B_\ell.
\end{equation}
Here reducedness means reducedness of the full member on the
canonical model $Z$, as in \cite[Theorem~8.5 and Corollary~8.6]{LBVolume}.
Thus every repeated prime detected by the estimate survives on $Z$
and has positive scaled degree at least one.
If this full canonical-model member is reduced, they give the stronger bound
\begin{equation}\label{eq:profile-defect46}
 U_{\ell,m}\ge\theta B_\ell+
       \max\{0,P_{m+1}-\theta p+\rho\},
 \qquad \rho=\operatorname{rank}
 [H^2(Y,\mathcal O_Y)\to\textstyle\bigoplus H^2(F_i,\mathcal O_{F_i})].
\end{equation}
In particular one may use $\rho=h$ when $q=0$, and $\rho=0$
without further information. For clarity, the term is $h$, not
$\theta h$: there is one connecting map for the whole member.
One proof is to apply the homogeneous power map from $R_{m+1}$ to
the quotient by the $\ell$-canonical traces of the $\theta$ moving
components. Its zero locus is contained in a space of dimension
$\theta p-\rho$, since a rational canonical form whose $\ell$-th
power is regular is regular. Reducedness identifies its generic
values on all other components, and purity gives equality globally.
The height theorem yields \eqref{eq:profile-defect46}.
If the inequality fails, a fixed prime has multiplicity at least
two and positive scaled degree. Accordingly every area bound must
also hold with
\begin{equation}\label{eq:repeated-cost46}
 r'=(m-2/v)/\theta>0
\end{equation}
in place of $r$. No visibility assumption on that repeated prime
is needed for this degree cost.


\subsection{The normalized polar square and marked adjunction}
\begin{lemma}\label{lem:polar-square46}
Let $a\in R_e$, $e\ge2$, and let $\alpha$ be a holomorphic two-form.
On a Gorenstein chart write $a=f\nu^e$ and
$\alpha=\iota_\partial\nu$. Set $\theta=\partial f$,
$\gamma=\partial^2 f$. There is a section $u\in R_{2e+2}$ with
\begin{equation}\label{eq:normalized-polar46}
 u=\theta^2-\frac{2e}{2e+1}f\gamma+
                         \frac{f^2\kappa}{2e+1},
\end{equation}
for a regular local coefficient $\kappa$ in each canonical frame.
It is unique modulo $a^2R_2$. If a visible prime has multiplicity
$k$ in $(a)$, then its multiplicity in $(u)$ is exactly $2k-2$.
An invisible prime of multiplicity $b$ in $(a)$ has multiplicity
at least $2b$ in $(u)$. In particular, if every visible component
is repeated, then $u/a$ is a nonzero section of $R_{e+2}$.
\end{lemma}
\begin{proof}
This is \cite[Lemma~9.1]{LBVolume}; the construction
is recalled to specify the actual global section. First lift
$\theta^2$ through the complete quotient by $a$, using vanishing
in degree $e+2$. If $u_0-\theta^2=f\sigma$, the coefficient
$(2e+1)\sigma+2e\gamma$ transforms modulo $f$ as a section of
the degree-$e+2$ quotient. The derivative terms in a change
$\nu'=g\nu$ cancel. Lift it using vanishing in degree two and
subtract $a$ times that lift divided by $2e+1$; this gives
\eqref{eq:normalized-polar46}. Depth two extends the quotient
sections across the finite non-Gorenstein locus. The same
transformation calculation gives the stated ambiguity.

At a visible prime write $f=t^kb$ with $b$ and $\partial t$ units.
The coefficient of $t^{2k-2}$ is
\[
 \left(k^2-\frac{2e}{2e+1}k(k-1)\right)b^2(\partial t)^2
 =\frac{k(k+2e)}{2e+1}b^2(\partial t)^2\ne0.
\]
At an invisible prime $\partial$ preserves its ideal, so both
derivatives have order at least $b$. These are every possible
divisorial pole of $u/a$; normality proves the final assertion.
\end{proof}

\begin{corollary}[The cubic covariant]\label{cor:polar-cubic-covariant}
Keep the notation of Lemma~\ref{lem:polar-square46}, and choose its
normalized section $u$. For a section $b$ of degree $d$, write
$D_a(b)=\{a,b\}_\alpha$; its local numerator is
$ef\partial b-d\theta b$. There are sections
\[
 v=-\frac{D_a(u)}{2e+2}\in R_{3e+3},\qquad
 p\in R_{2e+4},\qquad w\in R_{4e+6}
\]
such that
\begin{equation}\label{eq:polar-cubic-covariant}
 D_a(v)=-(3e+3)u^2+a^2p,\qquad
 v^2-u^3=a^2w,\qquad D_a(w)=2vp.
\end{equation}
In the canonical frame of \eqref{eq:normalized-polar46},
\begin{equation}\label{eq:polar-cubic-leading}
 p\equiv\frac{3e^2(2e+3)}{(2e+1)^2}\gamma^2
                    \pmod{(f,\theta)},\qquad
 w\equiv-\frac{\theta^2p}{2e+3}\pmod f.
\end{equation}
The choice of $u$ is unique when $R_2=0$, as in every application
below.
\end{corollary}
\begin{proof}
The bracket makes $v$ a global regular section. To verify the
divisibility in \eqref{eq:polar-cubic-covariant}, put
$S=2e+1$, $T=2e+2$ and
\[
 U=S\theta^2-2ef\gamma+f^2\kappa,
 \qquad V=T\theta U-ef\partial U.
\]
Thus $u=U/S$ and $v=V/(TS)$. Applying $\partial f=\theta$
and $\partial\theta=\gamma$, direct differentiation gives
\[
 S\bigl(ef\partial V-(3e+3)\theta V\bigr)
       +T(3e+3)U^2\in(f^2),
\]
and, after division by $f^2$, its residue modulo $(f,\theta)$ is
$3Te^2(2e+3)\gamma^2$. The same expansion gives
$SV^2-T^2U^3\in(f^2)$ and, on dividing by $T^2S^3f^2$,
the stated restriction of $w$. Finally, differentiating
$v^2-u^3=a^2w$ by the weighted derivation $D_a$ gives
$D_a(w)=2vp$, since $D_a(a)=0$.

This calculation is a polynomial identity in $f$, its derivatives
and the regular coefficient $\kappa$; it does not assume that $f$
is reduced or that its prime components are visible. Consequently
both divided rational sections $p$ and $w$ are regular at every generic point of
the full divisor $(a)$, with its actual multiplicities. They are
regular away from $(a)$ as well. Normality extends them uniquely
as sections of their respective canonical powers across codimension two.
The argument works for every permitted choice of $u$; the final
uniqueness assertion is its ambiguity modulo $a^2R_2$.
\end{proof}

\begin{lemma}[Normalized marked adjunction]\label{lem:full-adjunction46}
Let $C$ be a projective curve that is a log canonical centre of
$(U,B)$ at its generic point, where $B\ge0$, $U$ is normal and klt,
and $K_U+B$ is rational Cartier. On $C^\nu$, let $r_x$ be the
nontrivial inertia orders of its canonical index-cover branches.
Then
\begin{equation}\label{eq:marked-adjunction46}
 2g(C^\nu)-2+\sum_x(1-1/r_x)\le(K_U+B)C.
\end{equation}
The same inequality holds for the inertia of a specified finite
quasi-\'etale cyclic root cover, by the same argument. If its root
line is $L$ and has primitive fraction $k_x/r_x$, then
$LC\equiv\sum_xk_x/r_x\pmod{\mathbb Z}$. In particular,
if the canonical line on a branch has primitive fraction $k_x/r_x$,
then $K_UC\equiv\sum_xk_x/r_x\pmod{\mathbb Z}$.
At an ordinary $\frac1R(1,-1,b)$ point, an inertia-$r$ branch
has normalized arc coordinate orders at least
$(j,r-j,\overline{bj})$, where $j\equiv-b^{-1}k\pmod r$.
\end{lemma}
\begin{proof}
The degree inequality is \cite[Lemma~10.2 and Proposition~10.3]{LBVolume}, based on
normalized-centre adjunction \cite[Theorem~1.2]{FH}.
One may retain the inertia terms by taking a global cyclic cover
whose local cover at the marked points is the canonical cover,
using an auxiliary ample divisor disjoint from them. Adjunction
upstairs and Riemann--Hurwitz downstairs have identical auxiliary
branch terms; cancellation leaves \eqref{eq:marked-adjunction46}.
Additional branch terms are nonnegative. The fractional-degree
identity is the degree formula for the resulting orbifold line
bundle. On a normalized lifted arc, equivariance of its three
coordinates gives the three stated residues and hence their
minimum nonnegative orders. A nonzero prescribed-character
numerator must pay the minimum order of its actual monomials.
\end{proof}

We will also use \cite[Lemma~9.9]{LBVolume}: if a visible prime
has $M=-n$ and $n/I\notin\mathbb Z$, its normalization on the
canonical model is a Du Val K3 surface with an integral ample Weil
class $H$ of square $n/I$. Here a \emph{Du Val K3 surface} means
a normal surface whose only singularities are rational double
points and whose minimal resolution is a K3 surface. The
exceptional ADE root lattice has rank at most nineteen, and
\begin{equation}\label{eq:k3-rr46}
 h^0(mH)=2+\frac{m^2H^2-\sum_Qq_Q(mH)}2\quad(m>0).
\end{equation}
The correction norms are $i(r-i)/r$ for $A_{r-1}$; for $D_r$
the vector and spin norms are $1,r/4$; for $E_6,E_7$ they are
$4/3,3/2$; $E_8$ has trivial discriminant group. Thus
$\sum_Qq_Q(H)-H^2\in2\mathbb Z$. These formulas include
nonprimitive local classes and extra Cartier root factors.

\subsection{From virtual constituents to actual quotient charts}
\begin{proposition}[Source transfer]\label{prop:actual-source-transfer46}
Let $g:Y\to Z$ be a crepant terminalization.
\begin{enumerate}
\item A closed-point crepant centre of index $r>1$ requires the
following subbasket; no such centre has index greater than six:
\[
\begin{array}{c|l}
r&\text{mandatory constituents}\\\hline
2&4(1,2)\\3&3(1,3)\\4&(1,2)+2(1,4)\\
5&(1,5)+(2,5)\\6&(1,2)+(1,3)+(1,6).
\end{array}
\]
Additional constituents with zero auxiliary character are allowed.
\item If $g^{-1}(P)$ has dimension at most one, the canonical
index at $P$ occurs as an actual terminal index somewhere in that
fibre. In general every canonical index is either at most six
or an actual terminal index.
\item Suppose $Q$ is the unique actual point of index $r>6$,
is an ordinary $\frac1r(1,-1,b)$ quotient, and its image has
index exactly $r$. Then $g$ is an isomorphism near its image
unless
\[
 4\mid r,\qquad b\equiv r/2\pm1\pmod r,
 \qquad\text{another actual terminal point has index }r/2.
\]
\end{enumerate}
\end{proposition}
\begin{proof}
For the first assertion, Kawakita constructs an effective exceptional
divisor that detects the index indicator function
\cite[Theorem~2.7]{Kawakita}; the nonzero-character constituents
are classified in his Theorem~4.2, Table~2. Normalizing $b$ to
$\min(b,r-b)$ gives the displayed list. The terminal basket is
unchanged between crepant terminal minimal models related by flops
\cite[footnote~1]{CCIII}, so the mandatory constituents occur in
the given basket.

For the second assertion take the canonical index-one cover at
$P$ and normalize its base change to $Y$. Crepancy identifies the
local covers above a point of terminal index $e$ with $r/e$
components, each with stabilizer of order $e$. For a positive
dimensional fibre, rationality and formal functions give
$H^1(C,\mathcal O_C)=0$ for its connected reduced curve $C$.
Indeed the $H^1$ maps from successive thickenings onto the reduced
fibre are surjective, since the kernels have dimension at most one.
The normalization sequence now forces all components to be rational,
the bipartite incidence graph to be a tree, and the local delta
invariant at a point with $s$ branches to be $s-1$. Equality in
this last statement identifies its local ring with the subring of
the branch normalization rings having equal constant terms. Thus
the branches are smooth with independent tangent directions.
A cyclic group acting on the finite tree fixes a vertex; a fixed
rational component has a fixed point. A full stabilizer therefore
occurs, so $e=r$ at some source point. A zero dimensional connected
fibre gives the same conclusion immediately.

For the third assertion, the first part makes the fibre at most
one dimensional. If it were a curve, its canonical cover tree
would have exactly one full-group fixed point: a second would
give another index-$r$ point. No component is invariant, since
a finite cyclic action on $\mathbb P^1$ has either two fixed
points or is trivial. At the unique fixed point the ambient cover
is smooth of dimension three, so at most three branches meet.
A three-cycle has tangent characters
$j,j+r/3,j+2r/3$; the characters $1,-1$ would then force $r=3$
or $6$. Thus only a two-cycle remains. Its eigencharacters differ
by $r/2$, forcing $b\equiv r/2\pm1$ and $4\mid r$.
The stabilizer of one component is $\mu_{r/2}$, and its action
on that component has a second fixed point. The stabilizer there
is exactly $\mu_{r/2}$, giving the additional half-index point.
When these conditions fail the fibre is zero dimensional; proper
birationality and normality imply a local isomorphism.
\end{proof}

A singleton virtual constituent of the indices used below is an
ordinary point: its axial multiplicity in the terminal
classification is one. The proposition proves unchanged ordinary
neighbourhoods at the following indices:
\[
\begin{array}{c|l}
\text{row}&\text{unchanged singleton indices}\\\hline
02347&7\\02349&7,12\\02350&3,5,12,14\\
02353&7,12\\02355&7,12\\02357&5,7,9,18\\02397&7.
\end{array}
\]
For index twelve the possible tangent exception has $b=5$, but
there is no terminal point of index six. In row 02357 first treat
the index-eighteen point, then the index-nine point; the latter
cannot map into an already unchanged neighbourhood. The other
high indices have no proper multiple among the terminal indices.
For the index-five points of rows 02350 and 02357, the missing
$(1,5)$ excludes a crepant point centre; for the index-three point
of row 02350 the missing third $(1,3)$ does the same. Larger
image indices are excluded by index transfer and the already
unchanged higher-index points; an index-six centre is also
excluded by the missing $(1,6)$. The same rational-tree proof
applies. Five cannot permute two or three branches without a fixed
branch. A three-cycle at an index-three point has a zero tangent
character, whereas $(1,-1,b)$ has none. These fibres are therefore
zero dimensional as well.


\subsection{A small K3 polarization obstruction}
\label{sec:k3-square105-independent}

\begin{lemma}[A Du Val K3 square obstruction]
\label{lem:k3-square105-extension}
No Du Val K3 surface has an integral ample Weil divisor $H$ with
$H^2=1/105$.
\end{lemma}

\begin{proof}
On the smooth K3 resolution, the exceptional ADE root lattice is
orthogonal to the positive class $H$, so its rank is at most nineteen.
The Du Val correction-norm congruence \eqref{eq:k3-rr46} is
\[
 \sum q_P(H)-H^2\in2\mathbb Z.
\]
Multiply by105 and reduce in the localization at seven.  Only
$A_6$ and $A_{13}$ can contribute under this rank bound.  For their
nonzero seven-primary classes the contributions are respectively
$-j^2$ and $-j^2/2$ modulo seven.  Both range through
$\{3,5,6\}$, so one factor cannot supply the required residue one.
At least two factors are required.  The rank bound and the still
required five- and three-primary contributions force exactly two
$A_6$ factors, one $A_4$, and one $A_2$.  At most one rank remains,
which can only add an $A_1$ correction.

Reduction at five forces the $A_4$ class to be $\pm2$, with norm
$6/5$.  Reduction at seven forces the two $A_6$ classes to be
$\pm2$ and $\pm3$, with norms $10/7$ and $12/7$.
The nonzero $A_2$ class has norm $2/3$.  Their sum is
\[
 \frac{10}7+\frac{12}7+\frac65+\frac23
 =5+\frac1{105}.
\]
This violates the even-integer congruence.  A single non-Cartier
$A_1$ adds $1/2$ and does not repair it; a Cartier factor changes
nothing.  Two $A_1$ factors would repair parity but have total
exceptional rank twenty.  This is impossible.
\end{proof}


\subsection{Two forms and the degeneracy divisor}

We isolate the two-form argument in a form that also applies to the
cyclic covers used later. For a rational Cartier divisor $E$ on a
smooth model, write $Q_E$ for the Hermitian pairing on holomorphic
two-forms induced by $c_1(E)$, with the usual positive scalar.

\begin{lemma}[Two-form degeneracy bound]
\label{lem:two-form-degeneracy-bound}
Let $Z$ be a canonical threefold with ample $K_Z$, and let
$Y\to Z$ be a crepant terminal $\mathbb Q$-factorial model.
Suppose
\[
 p_g(Z)=q(Z)=P_3(Z)=0,\qquad h^{2,0}(Z)=2.
\]
Put $K=K_Y$, $V=K^3$, $v=IV$, and $c=4/3$, where $IK$ is Cartier.
Assume
\[
 V<1/14,\qquad c(c+1)V<1,\qquad c(c+1)^2V<1/3.
\]
Then the degeneracy divisor $D$ of the two forms satisfies
\[
 Q_D=Q_K,\qquad 0<K^2D\le\frac43V.
\]
Every visible prime in $D$ has normalization of geometric genus
at most one. The components with any fixed restriction-kernel
line have, with their multiplicities in $D$, total canonical
degree at most $2V/3$.

In particular there are visible primes $T,U$ with distinct
restriction kernels. With
$n=IK^2T$, $m=IK^2U$, $M=IKT^2$, $N=IKU^2$, and $h=IKTU$,
they satisfy
\begin{align}
 1\le n,m&\le\lfloor2v/3\rfloor,&
 -n\le M&<n^2/v,&-m\le N&<m^2/v,\nonumber\\
 0\le h&\le\min(n+M,m+N),&
 m^2M+n^2N&\ne2nmh,\label{eq:two-form-test46}\\
 (n^2-vM)(m^2-vN)&\ge(nm-vh)^2.\nonumber
\end{align}
The simultaneous character congruences
\eqref{eq:degree-character46}--\eqref{eq:mixed-character46}
hold for the same local characters.
\end{lemma}

\begin{proof}
Let $\mathcal H=H^0(\Omega^2)$ on a smooth model. Its generic
evaluation rank is two: rank one with two sections would give
$V\ge1/14$ by \cite[Lemma~2.7]{ChenJiangForms}. The two anchors
$\mathcal O_Y(-K)\to T_Y$ span a saturated rank-two sheaf $G$,
with
\[
 \det G=\mathcal O_Y(D-2K),\qquad
 (T_Y/G)^{**}=\mathcal O_Y(K-D),\qquad D\ge0.
\]
The Frobenius obstruction is a section of $3K$, hence vanishes.
The Bott partial connection on the quotient therefore gives
$Q_D=Q_K$, as in \cite[Theorems~6.1, 6.5 and~6.7]{LBVolume}.
To see that this identity includes the singular locus, take a
line-bundle multiple of $K-D$. Contracting its transition
cocycles with either anchor gives the class
$c_1(K-D)\wedge\alpha$. The partial-connection coefficients
extend across the codimension-two foliation singular set by
Hartogs. At the isolated terminal points the canonical sheaf
has depth three, so passage to the regular locus does not alter
the relevant $H^1$. Grauert--Riemenschneider vanishing then
identifies the group with that on the smooth model. Pairing with
the conjugate forms gives the asserted identity.

The pairing $Q_K$ is positive definite because $K$ is nef and big.
Effective primes give positive semidefinite restriction pairings;
primes of canonical degree zero have zero pairing, since their
images have dimension at most one and reflexive forms pull back
functorially. Thus $K^2D>0$. To apply slope semistability, push
the saturated sheaves to $Z$ and compute degrees on a general
complete canonical-intersection curve avoiding the exceptional
images. The ample-canonical-model semistability theorem
\cite[Theorem A]{Guenancia} gives
\[
 K^2(D-2K)/2\le -V/3,
 \qquad K^2(Z_\gamma-K)\le -V/3,
\]
where $Z_\gamma$ is the divisorial zero part of
$0\ne\gamma\in\mathcal H$. Hence
$K^2D\le4V/3$ and $K^2Z_\gamma\le2V/3$.

We next check the small-surface assertion. Suppose that a prime
$T$ with $K^2T=c'V$, $c'\le4/3$, has a smooth normalization model
$S$ with $p_g(S)\ge2$, and let $A$ be the canonical pullback.
Then $u=AK_S\ge0$, and adjunction and Hodge index give
\[
 d=A^2=c'V,
 \qquad u=AK_S\le c'(c'+1)V<1,
 \qquad u^2/d\le c'(c'+1)^2V<1/3.
\]
The general-type case contradicts Hodge index
against the integral minimal canonical class. In the properly
elliptic case every $A$-null curve is vertical. Indeed a horizontal
null curve with minimal-image canonical degree $s\ge1$ gives,
after writing $A=\sigma^*A_0-\sum a_iE_i^*$ and
$b=A_0K_{S_0}/s$, the adjunction estimate
\[
 d\le b^2(s+2)+2b\sum a_i
   =2bu+b^2(2-s)\le3u^2,
\]
contrary to $u^2/d<1/3$. A general elliptic fibre therefore avoids
the non-Gorenstein locus and has positive integral $A$-degree.
Since $p_g(S)\ge2$, the canonical system moves at least one such
fibre, forcing $u\ge1$. The remaining Kodaira dimensions have
$p_g(S)\le1$. This proves the small-surface assertion; it is the
argument of \cite[Lemmas~6.2--6.3 and Proposition~6.4]{LBVolume} under the explicit
strict inequalities stated here.

Each visible restriction pairing now has rank one. Their sum,
with the multiplicities in $D$, is $Q_K$, so at least two have
distinct kernels. If $T$ occurs with multiplicity $b$ and
$\gamma$ lies in its restriction kernel, choose transverse
coordinates $T=(x=0)$ and a nonvanishing form $dy\wedge dz$.
Write
$\gamma=f\,dy\wedge dz+g\,dz\wedge dx+h\,dx\wedge dy$.
The degeneracy multiplicity gives $g,h\in(x^b)$ and restriction
gives $f|_T=0$. Closedness implies $f_x+g_y+h_z=0$, hence
$f\in(x^{b+1})$. Thus $Z_\gamma\ge bT$. Summing over the
components with the same kernel proves their degree bound.

Visible adjunction gives $M\ge-n$ and $N\ge-m$. Equality in
the single-divisor Hodge inequality would make the corresponding
restriction pairing proportional to the rank-two pairing $Q_K$
by Lemma~\ref{lem:common-tools46}, which is impossible.
A form vanishing along one prime restricts nontrivially to the
other, giving $h\le\min(n+M,m+N)$. The three-class Hodge
determinant gives the last inequality in
\eqref{eq:two-form-test46}. Finally the excluded equality in its
middle line would make $mT-nU$ Hodge-null and therefore
$mQ_T=nQ_U$, contradicting the distinct kernels.
\end{proof}


\subsection{The Euler-characteristic-three branch}\label{sec:chi3}
For $V=1/420$ and $\chi=3$, the preliminary reduction gives
$p_g=q=0$ and $h^{2,0}=2$. If $P_3=0$, all hypotheses of
Lemma~\ref{lem:two-form-degeneracy-bound} hold. We apply its
simultaneous inequalities to the complete census.

\paragraph{The exhaustive finite test}
For each basket enumerate character multisets in
$(\mathbb Z/r)^k$ for every repeated type $k(b,r)$, calculate
the two single-divisor residues, and retain the finite ranges in
\eqref{eq:two-form-test46}. For two such states, enumerate every
matching of their character multisets and compute the mixed
residue. This is an explicit finite enumeration: all character
coordinates lie in $[0,r-1]$, $n,m\le\lfloor2v/3\rfloor$, and
$h\le\min(n+M,m+N)$. Repeated matching is equivalently the
subset recurrence
\[
 D_{j+1}(S\cup\{k\})\supseteq
       \{t+a_jc_k\pmod r:t\in D_j(S)\},\quad D_0(\varnothing)=\{0\}.
\]
It retains all permutations and introduces none outside the
necessary enlargement.

Of the $263$ equality baskets with $P_3=0$, $259$ have no
admissible different-kernel pair. The other four all have $I=4620$,
$v=11$, and the only numerical pairs are
\[
 (n,M,m,N,h)=(5,1,6,2,4),\quad(6,2,7,-7,0).
\]
Both have zero Hodge determinant. The per-basket empty lists are
tabulated in Appendix~\ref{app:census}; the four nonempty lists
are displayed completely here.

When every different-kernel pair is coplanar, choose one such
pair $T,U$. On a smooth general canonical hyperplane, the plane
spanned by $K,T$ has signature $(1,1)$, by the strict single-prime
Hodge inequality. Its orthogonal complement is negative definite.
The zero Gram determinant therefore puts $U$ in this plane.
Every further visible prime has kernel distinct from that of
$T$ or $U$, so the same argument puts it in the same plane.
The map of this plane to Hermitian restriction forms is injective:
evaluate a relation on the kernel of the rank-one prime form.
Since $Q_D=Q_K$, the full visible part $E$ has $E-K$ Hodge-null.
Lemma~\ref{lem:common-tools46} descends it to an actual torsion
class $\tau=E_Z-K_Z$, Cartier in codimension two. This is the
argument of \cite[Proposition~6.9]{LBVolume}, using only the rank-one
restriction, distinct kernels and Hodge descent established above
at equality. Invisible components and all multiplicities are
retained until this step.

For a null character vector put $C=\sum\Phi_{b,r}(i)$. Retain
the zero linear and square residues, and zero mixed residue
against both members of at least one admitted pair. The resulting
counts, by correction mass, are
\[
\begin{array}{c|rrrrr|r}
\text{row}&C=0&C=1&C=2&C=3&C=4&\max C\\\hline
\mathrm{C3\!-!00166}&1&2&0&0&0&1\\
\mathrm{C3\!-!47272}&1&1&0&1&1&4\\
\mathrm{C3\!-!47286}&1&1&0&1&0&3\\
\mathrm{C3\!-!49956}&1&1&0&1&1&4.
\end{array}
\]
The matchings used for different necessary pairwise equations
may vary; this only enlarges the allowed set. In particular the
maximum applies separately to every power of the actual torsion
class, without assuming additivity of sequential virtual characters.

If $s$ is the exact order of $\tau$, its connected canonical
quasi-\'{e}tale cover has
\[
 V'=s/420,\quad p_g'\ge1,\quad
 2\chi'=6s-\sum_{j=1}^{s-1}C_j\ge s+5,
\]
using even the weaker uniform cap $C_j\le5$.
Moreover $s\ge3$, since $sE\sim sK$ and $P_1=P_2=0$.
Thus $h^{2,0}\ge\chi'\ge(s+5)/2$. The form-volume bounds
in \cite[Corollary~2.3]{ChenJiangForms}, also recorded in
\cite[Section~12, equation~(12.5)]{LBVolume}, give $V'\ge1/18$ for
$3\le s\le14$, contradicting $s/420\le1/30$.
For $s\ge15$ they give
$V'\ge(h^{2,0}-3)/54\ge(s-1)/108>s/420$.
All orders are covered. The four coplanar baskets are impossible.

\paragraph{The twenty baskets with a cubic section}
Seventeen of them have $P_5=1$, $P_9\le2$. Choose nonzero
$a\in R_3$, $b\in R_5$. The linear map
$\mathcal H\to R_9/\mathbb C a^3$, $\alpha\mapsto D_{a,\alpha}b$,
has a nonzero kernel. For a kernel form, $D^2b=0$, hence $Db=0$.
The sections $a^5,b^3$ are independent by the integral-root
argument and $P_1=0$. Lemma~\ref{lem:commuting46} gives
$V\ge1/(15\cdot16)>1/420$.

The other three rows are C3-00047, C3-00050 and C3-00060.
They have $I=420,v=1$. For an actual member of $3K$, its
restriction form is $3Q_K$, so at least two rank-one visible
components occur. One has $n\le\lfloor3v/2\rfloor=1$.
The single-divisor congruences with $-n\le M<n^2/v$ have
no solution for $n=1$ in each of these three baskets. This is
the single-column computation printed in the appendix.
It excludes all three.

\begin{theorem}\label{thm:no-chi3-equality46}
There is no canonical threefold with $V=1/420$ and $\chi=3$.
\end{theorem}
\begin{proof}
The disjoint counts are $259+4+17+3=283$.
Proposition~\ref{prop:equality-census} proves their completeness.
\end{proof}

We will also need the following consequence for cyclic covers.
The same single-divisor calculation for either basket in
\eqref{eq:cover-baskets46} gives $I=420,v=2$, $P_3=0$ and
no state with $1\le n\le\lfloor2v/3\rfloor=1$.
Thus neither basket can be realized. This specialized observation
will replace a global nonvanishing theorem in the octic proof.



\section{Complete restriction, actual torsion and the Koszul defect}
\label{sec:restriction-toolbox}

The numerical tests in this section concern actual sections and actual
divisor classes. Their usefulness comes from retaining the cohomology
of the complete restriction sheaf, including its connecting direction.
We first record the possible torsion profiles, then prove the
restriction and Koszul statements that will be used in the exclusions.

\subsection{Actual torsion and opposite section profiles}\label{sec:torsion-local46}
An actual torsion class in this paper is an element of
$\operatorname{Cl}(Z)_{\rm tors}$ that is Cartier at every generic
curve point. Its sequential integral Weil pullback $D$ to a
terminalization satisfies $K^2D=KD^2=0$
\cite[proof of Proposition~3.3; Section~12]{LBVolume}. Indeed its rational pullback is numerically
trivial. Integral pullback compatibility identifies the two pullbacks
at the generic points of curve-centred exceptional divisors. Their
difference is supported over points and therefore has zero intersection
with $K$. This assertion is applied separately to each actual class.

For a tuple of local characters $\mathbf i=(i_Q)$ set
\[
 C(\mathbf i)=\sum_Q\Phi_Q(i_Q),\qquad
 d_2(\mathbf i)=\frac{C(\mathbf i)}2+
       \sum_Q\bigl(\Phi_Q(i_Q+1)-\Phi_Q(1)\bigr).
\]
A \emph{null character tuple} satisfies
\begin{equation}\label{eq:null-characters46}
 \sum_Q(I/r_Q)b_Q^2i_Q=0\pmod I,\qquad
 \sum_Q(I/r_Q)\overline{b_Qi_Q}(r_Q-\overline{b_Qi_Q})=0\pmod{2I}.
\end{equation}
In particular $C(\mathbf i)$ is a nonnegative integer. Put $D_2=0$ and
\begin{equation}\label{eq:twist-differences46}
 D_m(\mathbf i)=P_m-P_2+
       \sum_{j=2}^{m-1}\sum_Q
       \bigl(\Phi_Q(i_Q+j)-\Phi_Q(j)\bigr),\quad m\ge3.
\end{equation}
Weak canonical Riemann--Roch, Serre duality and vanishing, in the form
of \cite[Section~12, equations~(12.1)--(12.3)]{LBVolume}, imply the following rules. For an actual
class $\tau$, choose necessary null tuples $\mathbf i,\mathbf i^-$
for $\tau,-\tau$. They have equal $C$; writing
$H_m(\tau)=h^0(Z,mK_Z+\tau)$, one has
\begin{align}
 H_m(\tau)&=H_2(\tau)+D_m(\mathbf i),\label{eq:twist-H46}\\
 H_2(\tau)+H_2(-\tau)&=2P_2+d_2(\mathbf i)+d_2(\mathbf i^-),
       \label{eq:twist-pair46}\\
 \chi(\mathcal O_Z(\tau))+\chi(\mathcal O_Z(-\tau))
       &=2\chi-C(\mathbf i).\label{eq:twist-euler46}
\end{align}
The weak difference formula has no curve term here, precisely because
the classes are Cartier in codimension two. The two opposite Euler
characteristics need not be equal.

These formulas give a finite necessary list with no bound on the
order of $\tau$. Enumerate all null tuples, group equal pairs
$(C,d_2,(D_m)_{m=2}^{120})$, and match any two states having the same
$C$. If the right side $S$ of \eqref{eq:twist-pair46} is a
nonnegative integer, try $H_2(\tau)=0,\ldots,S$ and
$H_2(-\tau)=S-H_2(\tau)$. Retain only nonnegative dimensions
through degree $120$ and the opposite-character product inequalities
\begin{equation}\label{eq:twist-product46}
 H_a(\tau)+H_b(-\tau)-1\le P_{a+b}
       \quad(2\le a,b\le24, H_aH_b>0).
\end{equation}
Each retained pair is called a \emph{section-dimension state}; it is
only a necessary pair of sequences. The counts for the $28$ rows are
in Table~\ref{tab:remaining28}. They are verified in two ways:
meet-in-the-middle enumeration of grouped tuples, and independent
primary congruence counts weighted by the multinomial number of
ordered tuples. Thus grouping repeated constituents does not lose
any tuple.

Further filters always come from actual multiplication. Examples are
$H_a(\tau)+P_b-1\le H_{a+b}(\tau)$, or squaring sections and
applying the same exhaustive list to the actual class $2\tau$.
We never impose $\mathbf i(2\tau)=2\mathbf i(\tau)$ at the virtual
constituents, nor infer that a retained state is realized. The
trivial class is retained whenever a cohomological argument assumes
that a nontrivial class has no degree-zero section.

\paragraph{Canonical multiplication and integral powers.}
\label{subsec:late-torsion46}
The ordinary canonical sections impose necessary conditions
on every actual torsion class, before a moving-system or
fixed-part role is assigned to that class.  Write
$H_m(\tau)=h^0(Z,mK_Z+\tau)$.  The domain multiplication
inequality gives
\begin{equation}
 H_{m+l}(\tau)\ge H_m(\tau)+P_l-1
 \quad\text{if }H_m(\tau)P_l>0.
 \label{eq:canonical-twist-hopf34}
\end{equation}
The same inequality applies to $-\tau$.  In particular
each nonzero fixed canonical section gives an injective
map $H_m(\tau)\hookrightarrow H_{m+l}(\tau)$.
This condition is different from the cancellation
inequality between $\tau$ and $-\tau$: both are
necessary and may be imposed together.  It uses no
additivity for virtual local basket characters.

For the refinements below, write
$q_\tau(m)=\sum_Q\Phi_Q(i_{\tau,Q}+m)$ and $C_\tau=q_\tau(0)$.
With $H_m^\pm=h^0(mK_Z\pm\tau)$, the preceding identities become
\begin{align}
 H_m^+&=H_2^++P_m-P_2+
       \sum_{j=2}^{m-1}(q_\tau(j)-q_0(j)),\label{eq:late-twist46}\\
 H_2^++H_2^-&=2P_2+C_\tau+
       q_\tau(1)+q_{-\tau}(1)-2q_0(1),\label{eq:late-anchor46}\\
 \chi(\mathcal O_Z(\tau))+\chi(\mathcal O_Z(-\tau))
       &=2\chi-C_\tau.\label{eq:late-euler46}
\end{align}
The difference formula is
\[
 E_\tau(m+1)-E_\tau(m)=m^2K^3/2-2\chi+q_\tau(m),
 \qquad E_\tau(m)=\chi(\mathcal O_Z(mK_Z+\tau)).
\]
Serre duality for the maximal Cohen--Macaulay divisorial sheaves gives
$E_\tau(m)=-E_{-\tau}(1-m)$; hence
\begin{equation}\label{eq:late-reflection46}
 q_\tau(m)=q_{-\tau}(-m)\quad\hbox{for every integer }m.
\end{equation}
This is a reflection of correction functions, not an assertion that
unrelated virtual characters can be added.

These identities also give the individual paired formulas used above,
without appealing to an aggregate cover formula. At $m=0$, the weak
difference formula and $E_\tau(0)=-E_{-\tau}(1)$ give
\eqref{eq:late-euler46}; reflection gives $C_\tau=C_{-\tau}$.
Add the two difference formulas at $m=1$ and use
$2P_2=V-6\chi+2q_0(1)$ to obtain
\eqref{eq:late-anchor46}. Finally, subtract the untwisted difference
formula and sum from $2$ to $m-1$; vanishing identifies these Euler
characteristics with section dimensions and proves
\eqref{eq:late-twist46}. Thus every test concerns actual opposite
classes separately, even though independent virtual tuples enlarge
their necessary numerical possibilities.

Here is a fully specified finite necessary test. Enumerate every
character assignment satisfying the two null congruences; repeated
identical constituents can be grouped, with their multinomial weights
retained. Form the full correction vector of period $I$. Pair only
vectors satisfying \eqref{eq:late-reflection46}, and use
\eqref{eq:late-anchor46} to enumerate the nonnegative integer anchors
$H_2^+,H_2^-$. Keep profiles nonnegative through degree $120$ and
satisfying
\[
 H_a^++H_b^--1\le P_{a+b}
 \quad(2\le a,b\le24, H_a^+H_b^->0).
\]
For a fixed candidate profile $T$, form a directed graph on the
remaining paired profiles. An edge $J\to J'$ is allowed only if the
four Hopf inequalities belonging to
\[
 \sigma+\tau=\sigma',\quad-\sigma-\tau=-\sigma',\quad
 \sigma'-\tau=\sigma,\quad-\sigma'+\tau=-\sigma
\]
hold in every recorded degree through $120$. Remove $T$ if it has
no path to a canonical profile, or no path to its inverse profile;
iterate. Every actual class survives, because its actual finite
sequence of multiples supplies both paths at every deletion stage.
A surviving path does not construct a torsion group.

For later use we record the outputs for four of the remaining baskets:
\[
\begin{array}{c|r|r|r|r}
 B&\text{ordered assignments}&\text{ordered null}&
 \text{paired profiles}&\text{power survivors}\\\hline
 02349&10321920&2016&16&8\\
 02353&6804000&189&7&3\\
 02355&1290240&64&16&6\\
 02357&408240&27&13&3.
\end{array}
\]
All surviving profiles have $H_2^\pm=0$. For $02353$ they also
have $H_3^\pm=H_6^\pm=0$, and for $02355$ they have $H_3^\pm=0$.
For $02357$ the complete low list is
\[
 (C,H_2,H_3,H_4)=(0,0,0,0),(3,0,0,1),(3,0,0,1).
\]
The two noncanonical profiles in this list are inverse, with
degree-six dimensions $(1,0)$ and $(0,1)$. Squaring gives $H_1=0$
for every actual torsion class in all four baskets. These lists are obtained by the finite procedure just described; its inputs are the complete baskets in Appendix~\ref{app:census}.

\subsection{Complete restriction to a prime member}
\label{sec:twisted-prime-toolbox}
\begin{lemma}[Complete restriction to a torsion-twisted prime]
\label{lem:twisted-prime-restriction34}
Let $Z$ be a normal projective canonical threefold with
$K_Z$ ample, $p_g(Z)=q(Z)=0$, and $h=h^{2,0}(Z)$.
Let $\Gamma$ be a finite subgroup of the Weil divisor class
group consisting of torsion classes.  Suppose that
\[
 0\ne s\in H^0(Z,\mathcal O_Z(dK_Z+\tau)),\qquad
 \tau\in\Gamma,\quad d\ge1,
 \qquad \operatorname{div}_Z(s)=T
\]
is a single reduced prime.  Put
\[
 \mathcal Q_{n,\sigma}=
 \mathcal O_Z(nK_Z+\sigma)\big/
 s\mathcal O_Z((n-d)K_Z+\sigma-\tau),
 \qquad
 G_n(\sigma)=h^0(Z,\mathcal Q_{n,\sigma}).
\]
For any two nonzero homogeneous sections of these quotient
sheaves, their product is nonzero.  In particular, whenever
both source spaces are nonzero,
\[
 G_{a+b}(\sigma+\rho)
 \ \ge\ G_a(\sigma)+G_b(\rho)-1.
\tag{\ref{lem:twisted-prime-restriction34}.1}
\]
Writing $H_n(\sigma)=h^0(Z,\mathcal O_Z(nK_Z+\sigma))$,
the dimensions needed below are
\begin{align*}
 G_n(\sigma)
 &=H_n(\sigma)-H_{n-d}(\sigma-\tau)
 &&(n-d\ge2),\\
 G_{d+1}(\tau)
 &=H_{d+1}(\tau)-P_1+h.
\end{align*}
\end{lemma}

\begin{proof}
Every sheaf $\mathcal O_Z(nK_Z+\sigma)$ in the statement
is maximal Cohen--Macaulay.  Locally, the normalized finite
cyclic cover trivializing the relevant rational-Cartier class
is quasi-\'{e}tale and canonical; its Cohen--Macaulay direct
image splits into character summands.  This is the same
local argument used for canonical powers in
Lemma~\ref{lem:complete-canonical-restriction46}. The injective multiplication
by $s$ therefore has a Cohen--Macaulay quotient of pure
dimension two.

The ideal $\mathcal O_Z(-T)$ annihilates each quotient.
Thus $\mathcal Q_{n,\sigma}$ is a torsion-free rank-one
sheaf on the integral surface $T$.  At its generic point,
the ambient local ring is a discrete valuation ring and
$s$ has order one.  After choosing local frames the quotient
is a one-dimensional vector space over $\mathbb C(T)$,
and multiplication of two nonzero generic coefficients is
nonzero.  Purity makes restriction of global sections to
this generic point injective.  Consequently there are no
zero products between nonzero homogeneous pieces.
The displayed dimension inequality is the usual bilinear
Hopf inequality, applied to each fixed pair of pieces.
No domain assertion for the sum over the finite group
$\Gamma$ is required.

For $n-d\ge2$, Kawamata--Viehweg vanishing, with the
numerically trivial torsion twist, gives
$H^1(Z,\mathcal O_Z((n-d)K_Z+\sigma-\tau))=0$.
The quotient exact sequence gives the first dimension
formula.  At $(n,\sigma)=(d+1,\tau)$ its denominator is
$\omega_Z$.  Serre duality gives $h^1(\omega_Z)=h$,
while $H^1(Z,\mathcal O_Z((d+1)K_Z+\tau))=0$.
The same exact sequence gives the second formula.
Equivalently, this distinguished quotient is the actual
dualizing sheaf $\omega_T$.
\end{proof}

For an untwisted prime member, the ordinary quotient $R/(s)$ is also
a domain. Indeed, if $xy$ vanishes on the reduced prime divisor of
$s$, then one factor, say $x$, has positive order there. The quotient
$x/s$ has no divisorial pole and is a canonical section by normality.
Thus $x\in(s)$. We will use the complete restriction when its extra
cohomology direction is needed, and the ordinary quotient when it suffices.

\subsection{The two-twist Koszul defect}
\label{sec:cohomological-quotient34}

The correction in the full restriction formula has a second useful
form. Two effective divisors with no common component give a Koszul
sequence, and its first cohomology imposes a lower bound on the
corresponding quotient. We state the version with both twists, so that
the same argument covers ordinary and torsion-twisted moving systems.

\begin{lemma}[Two-twist Koszul defect]
\label{lem:two-twist-koszul34}
\label{lem:twisted-koszul-defect34}
Let $Z$ be a normal projective canonical threefold with ample $K_Z$.
Let $\lambda,\mu,\sigma$ be numerically trivial $\mathbb Q$-Cartier
Weil divisor classes, and put $H_j(\rho)=h^0(Z,jK_Z+\rho)$.
Suppose that nonzero sections
\[
 s\in H^0(Z,aK_Z+\lambda),\qquad
 t\in H^0(Z,bK_Z+\mu),\qquad a,b>0,
\]
have no common prime in their effective divisors. Then
\begin{align}
 h^2(Z,\mathcal O_Z(\lambda+\mu-\sigma))
 \le{}&H_{a+b+1}(\sigma)-H_{b+1}(\sigma-\lambda)
       -H_{a+1}(\sigma-\mu)+H_1(\sigma-\lambda-\mu).
 \label{eq:two-twist-koszul34}
\end{align}
In particular, for $\lambda=\tau$ and $\mu=0$,
\begin{align}
 h^2(Z,\mathcal O_Z(\tau-\sigma))
 \le{}&H_{a+b+1}(\sigma)-H_{b+1}(\sigma-\tau)
       -H_{a+1}(\sigma)+H_1(\sigma-\tau).
 \label{eq:twisted-koszul-defect34}
\end{align}
Every actual torsion Weil class is $\mathbb Q$-Cartier, so the lemma
applies in particular to the actual torsion classes used below.
\end{lemma}

\begin{proof}
Write $\mathcal L_j(\rho)=\mathcal O_Z(jK_Z+\rho)$. These
divisorial sheaves are maximal Cohen--Macaulay: locally, a cyclic
cover that trivializes a Cartier multiple of the relevant Weil
divisor is quasi-\'etale; its normalization has canonical
singularities, and the sheaf in question is a character summand of
its Cohen--Macaulay direct image. In particular Serre duality applies
to these sheaves.

Let $\mathcal I$ be the image of multiplication by $(s,t)$ in
$\mathcal L_{a+b+1}(\sigma)$. We have an exact sequence
\[
 0\longrightarrow\mathcal L_1(\sigma-\lambda-\mu)
 \xrightarrow{(t,-s)}
 \mathcal L_{b+1}(\sigma-\lambda)\oplus
 \mathcal L_{a+1}(\sigma-\mu)
 \xrightarrow{(s,t)}\mathcal I\longrightarrow0.
\]
Indeed, a relation $sA+tB=0$ identifies $A/t=-B/s$ as a rational
section of the kernel sheaf. At every prime one of $s,t$ has order
zero. Thus the quotient has no divisorial pole, and reflexivity
extends it across codimension two.

Since $a,b>0$ and the twists are numerically trivial,
Kawamata--Viehweg vanishing kills $H^1$ of the two middle sheaves.
Taking sections gives
\begin{align*}
 h^0(\mathcal I)={}&H_{b+1}(\sigma-\lambda)
 +H_{a+1}(\sigma-\mu)-H_1(\sigma-\lambda-\mu)\\
 &+h^1(Z,\mathcal L_1(\sigma-\lambda-\mu)).
\end{align*}
Duality identifies the last term with
$h^2(Z,\mathcal O_Z(\lambda+\mu-\sigma))$.
The inclusion $\mathcal I\subset\mathcal L_{a+b+1}(\sigma)$
now proves the inequality. Only divisorial coprimality was used;
base points in codimension two cause no difficulty.
\end{proof}

\begin{corollary}[A moving torsion translate]
\label{prop:twisted-quotient-correction34}
Assume in addition $p_g(Z)=q(Z)=0$. If $s\in R_l$ and
$u\in H^0(Z,pK_Z+\tau)$ have no common prime, where $l,p>0$
and $\tau$ is a numerically trivial $\mathbb Q$-Cartier Weil class,
then
\begin{equation}
 H_{l+p+1}(\tau)-H_{p+1}(\tau)-P_{l+1}+P_1
 \ge\chi(\mathcal O_Z)-1.
 \label{eq:twisted-quotient-correction34}
\end{equation}
In particular this applies to every nonzero $s$ if the actual
complete system $|pK_Z+\tau|$ has no divisorial fixed part.
\end{corollary}
\begin{proof}
Apply Lemma~\ref{lem:two-twist-koszul34} with
$(a,b,\lambda,\mu,\sigma)=(l,p,0,\tau,\tau)$.
Then $h^2(\mathcal O_Z)=\chi(\mathcal O_Z)-1$.
For the last assertion choose a member avoiding the finitely many
primes in the divisor of $s$.
\end{proof}

\begin{proposition}[The irrelevant torsion of a two-section quotient]
\label{prop:quotient-correction34}
Assume $p_g(Z)=q(Z)=0$, and put
$h=\chi(\mathcal O_Z)-1$ and $\mathfrak m=R_{>0}$.
Let $s\in R_a$, $u\in R_n$ be nonzero, with $a,n>0$ and
divisorially coprime. They form a regular sequence in $R$, and
for $B=R/(s,u)$ one has
\[
 H^0_{\mathfrak m}(B)_j=
 \begin{cases}\mathbb C^h,&j=a+n+1,\\0,&j\ne a+n+1.\end{cases}
\]
Consequently
\begin{equation}
 P_{a+n+1}-P_{n+1}-P_{a+1}+P_1\ge h.
 \label{eq:quotient-correction34}
\end{equation}
\end{proposition}
\begin{proof}
Vanishing and duality give
\begin{equation}
 H^1(Z,\mathcal O_Z(mK_Z))=
 \begin{cases}H^2(Z,\mathcal O_Z)^*,&m=1,\\0,&m\ne1.
 \end{cases}
 \label{eq:quotient-H1-profile34}
\end{equation}
For $m\ge2$ this is vanishing; for $m\le-1$ it follows by
duality from degree $1-m\ge2$; the remaining two degrees use
$p_g=q=0$. Since $R$ is the full normal section ring,
$H^0_{\mathfrak m}(R)=H^1_{\mathfrak m}(R)=0$, and the
Proj correspondence identifies $H^2_{\mathfrak m}(R)$ with
$\mathbb C^h$ in degree one and zero in the other degrees.

The section $s$ is regular because $R$ is a domain. Normal
division, as in the preceding Koszul proof, shows that $u$ is
regular on $C=R/(s)$. The two exact multiplication sequences give
\[
 H^0_{\mathfrak m}(C)=0,\qquad
 H^1_{\mathfrak m}(C)=H^2_{\mathfrak m}(R)(-a),\qquad
 H^0_{\mathfrak m}(B)=H^2_{\mathfrak m}(R)(-a-n).
\]
Here $M(-a)_j=M_{j-a}$; each intervening multiplication on local
cohomology is zero because its source and target have different
single supporting degrees. This proves the asserted torsion.
Finally $H_B(t)=(1-t^a)(1-t^n)H_R(t)$, and $B_{a+n+1}$ contains
that $h$-dimensional subspace, giving the inequality.
\end{proof}

For the baskets with $\chi=2$, a zero coefficient in the critical
degree already contradicts this one-dimensional torsion. In each
application below we therefore produce the coprime sections before
using the numerical inequality.


\subsection{Visible saturation and a single reduced component}
\label{sec:visible-common}
Let $\alpha$ span the holomorphic two-forms and let $\lambda$ be its
Hodge functional. It is nonnegative on effective primes, zero on
invisible primes and on torsion classes, and positive on $K$.
\begin{lemma}[Visible saturation of a canonical member]
\label{lem:visible-saturation34}
Let $D\sim mK_Z$ be an effective pluricanonical divisor.
Suppose every visible prime component $T_i$ of its actual
terminal pullback satisfies a Hodge equality and descends
to a rational canonical multiple with a torsion correction:
\[
 q_iT_{i,Z}\sim p_iK_Z+\tau_i,
 \qquad q_i>0.
\]
Put $a_i=p_i/q_i=K_Z^2T_{i,Z}/K_Z^3$, and let $b_i$
be its multiplicity in $D$.  Then
\[
 \sum_i b_i a_i=m,
\]
and the entire invisible complement of $D$ vanishes
as a divisor on $Z$.  If $v=1$, every relevant degree
$a_i$ is an integer.  More generally, whenever all the
descents have $q_i=1$, any decomposition of $D$ into
two nonempty sums of such prime components of total
degrees $j$ and $m-j$ supplies actual opposite torsion
classes with
\[
 H_j(\tau)>0,
 \qquad H_{m-j}(-\tau)>0.
\]
\end{lemma}

\begin{proof}
Linear equivalence modulo torsion gives
$\lambda(T_i)=a_i\lambda(K)$.  Since the complement
is invisible,
\[
 m\lambda(K)=\lambda(D)
 =\sum_i b_i a_i\lambda(K).
\]
Cancel $\lambda(K)>0$.  The visible part already
has the entire canonical square-degree $mK_Z^3$.
Every nonzero effective prime on $Z$ has strictly
positive square-degree against the ample $K_Z$, so
no invisible component remains downstairs.

When $v=1$, the integral null-divisor relation has
$q_i=1$ and $p_i=n_i=IK_Y^2T_i$.  In the stated
more general integral-descent situation, the two sums of
components therefore have classes $jK_Z+\tau$ and
$(m-j)K_Z+\sigma$.  Their sum is the actual divisor
$D\sim mK_Z$, so $\sigma=-\tau$ as global divisor
classes.  Effectivity supplies the two sections.
No local virtual-character addition is used.
\end{proof}

\begin{lemma}
\label{lem:single-visible-square34}
Let $Y$ be a terminal $\mathbb Q$-factorial minimal threefold of general type
with $q=p_g=0$ and $h^{2,0}=1$, and let $\alpha$ span its
holomorphic two-form space.  Let $s\in R_m$ be nonzero,
$m\ge2$.  Suppose the full visible part of its actual
effective section divisor is a single prime $T$ with
multiplicity one; every other component is invisible to
$\alpha$.  If $P_{m+1}>0$, then
\[
 P_{2m+2}\ge3.
\]
\end{lemma}

\begin{proof}
Put $c=m+1$ and
$\mathcal Q_j=\mathcal O_Y(jK_Y)/s\mathcal O_Y((j-m)K_Y)$.
Complete restriction gives a section
$\vartheta\in H^0(Y,\mathcal Q_c)$ with nonzero connecting
class in $H^1(Y,\omega_Y)$.  Choose $0\ne b\in R_c$.
The full-visible-scheme lemma says that the image of
$\vartheta$ at the generic point of $T$ is nonzero, and
is linearly independent from the generic image of $b$
if that image is nonzero.  This statement uses the full
invisible complement, including its scheme multiplicities.

First suppose $b$ has order zero along $T$.  Since $T$
has coefficient one, its generic quotient algebra is a
field, after a divisor-frame trivialization.  The two
independent generic images of $b$ and $\vartheta$ have
at least three independent quadratic products.  These
products are restrictions of sections of $\mathcal Q_{2c}$.
Thus
\[
 P_{2c}-P_{m+2}=h^0(Y,\mathcal Q_{2c})\ge3,
\]
and the desired bound follows.

Suppose instead $l=\operatorname{ord}_T(b)\ge1$.
The generic image of $\vartheta^2$ is nonzero.
Kawamata--Viehweg vanishing gives
$H^1(Y,(m+2)K_Y)=0$, so $\vartheta^2$ lifts to a
section $z\in R_{2c}$ with $\operatorname{ord}_T(z)=0$.
The weighted differential bracket
\[
 D_s(b)=\{s,b\}_\alpha\in R_{2m+2}
\]
has exact order $l$ along $T$.  Indeed, in a smooth
generic chart write $s=t a$ and $b=t^l b_0$, with
$a,b_0$ units.  Visibility means that the normal
coefficient of the local vector field defining
$\alpha$ is a unit.  The leading bracket coefficient
is a nonzero unit times
\[
 ml-(m+1),
\]
which cannot vanish for integers $m\ge2$, $l\ge1$.
The three sections $z,D_s(b),b^2$ consequently have
distinct orders $0,l,2l$ along $T$, so they are
linearly independent in $R_{2m+2}$.
\end{proof}

\begin{lemma}\label{lem:adjacent-gaps46}
Assume $\chi=2$, $p_g=q=0$, and $a\in R_e\setminus\{0\}$ with
$P_{e+1}=P_{e+2}=0$. Its visible part consists of one prime $T$
with coefficient one. Writing $(a)_Y=T+E$, the full remainder
$E$ is invisible and $\lambda(T)=e$, where $\lambda(K)=1$.
\end{lemma}
\begin{proof}
The complete dualizing space of the member has dimension one.
The ordinary traces from distinct visible primes are nonzero and
independent at their separate generic points, including in the
socle layers of a nonreduced member. Positivity supplies a visible
prime, so it is unique. If its coefficient were at least two,
Lemma~\ref{lem:polar-square46} would supply a nonzero section of
$R_{e+2}$, a contradiction. The Hodge-mass identity follows by
additivity and invisibility of the entire complement.
\end{proof}

\begin{lemma}[A visible-prime floor]
\label{lem:visible-prime-floor}
Assume $p_g=q=0$ and $h^{2,0}=1$. Fix an integer $f>0$.
Suppose that every necessary visible state $(n,M)$ with $0<n<f$
is a Hodge equality $vM=n^2$, that $n/v$ is integral, and that
$H_{n/v}(\tau)=0$ for every actual torsion class supplied by
null-divisor descent. Then every visible prime has degree at least
$f$. Consequently, if $0\ne a\in R_e$ and $ev<2f$, the full
visible part of $(a)_Y$ is one reduced prime. If $ev=f$, its image
is the entire reduced prime divisor of $a$ on $Z$.
\end{lemma}
\begin{proof}
A visible prime of degree less than $f$ would descend to the
effective class $(n/v)K_Z+\tau$, contrary to the asserted vanishing.
A pluricanonical member has positive Hodge value and therefore
contains a visible component. Two such components, or a repeated
one, would have degree at least $2f$. If $ev=f$, the visible prime
exhausts the canonical degree, and ampleness of $K_Z$ excludes
every additional component downstairs.
\end{proof}


\subsection{Crepant separation of zero-contact divisors}
\label{sec:crepant-graph-binary34}

Zero canonical intersection need not make two Weil divisors disjoint
on the canonical model. The following construction separates them on
a crepant model and, at the same time, restricts their local characters.
It will be used both for the adjacent-pencil obstruction and for
cancellation at a zero-contact endpoint.

\begin{lemma}[Crepant graph separation]
\label{lem:crepant-graph-binary34}
\label{lem:crepant-graph-general}
Let $f:Y\to Z$ be a projective crepant terminal
$\mathbb Q$-factorial model of a canonical threefold with ample $K_Z$,
and write $K=K_Y=f^*K_Z$. Let $H,J$ be effective integral Weil
divisors with no common prime. Suppose
\[
 KHJ=0,\qquad H-J-aK\sim P,\qquad a\in\mathbb Z,
\]
where $P$ is supported on $f$-exceptional divisors with point
centres; its coefficients may have either sign. There is a
projective crepant terminal $\mathbb Q$-factorial model
$Y^\dagger\to Z$ carrying an integral Weil divisor $D$ with
\[
 (K_{Y^\dagger})^2D=K^2H,
 \qquad D\sim0\ \text{or}\ aK_{Y^\dagger}
 \quad\text{locally at every point}.
\]
Consequently, for the terminal basket $B=\{(b_Q,r_Q)\}$ and
$I=\operatorname{lcm}_Qr_Q$,
\begin{equation}
 IK^2H\equiv
 a\sum_Q\frac I{r_Q}b_Q^2\varepsilon_Q\pmod I,
 \qquad \varepsilon_Q\in\{0,1\}.
 \label{eq:binary-graph-congruence34}
\end{equation}
For $a=1$, the local characters are thus zero or the canonical
character, even when $Z$ is not $\mathbb Q$-factorial.
\end{lemma}

\begin{proof}
Push down to obtain $H_Z-J_Z\sim aK_Z$. Choose a rational
section $\eta$ of $aK_Z$ with divisor $H_Z-J_Z$, and choose
$r>0$ such that $rK_Z$ is Cartier. Let $G$ be the normalization
of the graph of $\eta^r$ in
$\mathbb P_Z(\mathcal O_Z\oplus\mathcal O_Z(raK_Z))$.
Then $\pi:G\to Z$ is projective and birational, with fibres of
dimension at most one, since normalization is finite.

Because $Y$ is $\mathbb Q$-factorial, $H\cdot J$ is an effective
rational one-cycle. Its components are contracted to points by
$f$, as $KHJ=0$ and $K_Z$ is ample. Let $S$ be the finite union
of these image points and the centres of the components of $P$.
Over $Z\setminus S$, the pulled-back section has divisor
$rH-rJ$ with disjoint supports. It is regular away from $J$,
and its reciprocal is regular away from $H$, by normality.
Thus $Y\setminus f^{-1}(S)$ maps to the graph and factors
through its normalization.

Every $\pi$-exceptional prime has a curve as centre, since the
fibres over points have dimension at most one. Its generic point
is dominated by a prime of $Y\setminus f^{-1}(S)$, whose
discrepancy over $Z$ is zero. Hence $K_G=\pi^*K_Z$ in
codimension one. This makes $K_G$ rational Cartier; discrepancies
on a further resolution are those over $Z$, so $G$ is canonical.

The zero and infinity sections of the projective-line bundle are
disjoint Cartier divisors. Their pullbacks have no exceptional
prime component: a prime over a curve contained in one section
would map to a curve in the unnormalized graph, contrary to the
finiteness of normalization. Thus these pullbacks are precisely
$rH_G$ and $rJ_G$, where the subscripts denote strict transforms.
They are disjoint, and $H_G-J_G\sim aK_G$. Therefore $H_G$ is
locally linearly equivalent to $aK_G$ near its support and to zero
off its support.

Take a crepant terminal $\mathbb Q$-factorialization
$g:Y^\dagger\to G$ and put $D=g^*H_G$. The local descriptions
show that $D$ is integral. Projection gives
$(K_{Y^\dagger})^2D=K_Z^2H_Z=K^2H$.
The terminal basket is unchanged under this crepant change of
minimal model, by its invariance under flops
\cite[footnote~1]{CCIII}. The degree-character congruence for
$D$ now gives \eqref{eq:binary-graph-congruence34}.
\end{proof}


\subsection{Cancellation and full conductors}
\label{sec:late-tools46}
We retain the entire divisor on the terminal model, including every
canonical-null component. This is essential when a complete restriction
forces a section to vanish on the visible component but does not yet
show divisibility by the whole canonical member. We normalize the
Hodge functional by $\lambda(K)=1$.

\paragraph{A half-canonical obstruction}
\begin{lemma}\label{lem:late-no-half46}
If a terminal basket contains a unique constituent of even index
$r$, and no other index is divisible by $r$, then its canonical
model has no integral Weil class $L$ with $2L\sim K_Z$.
\end{lemma}
\begin{proof}
For the sequential pullback of $L$, Serre duality and the integral
shift formula give $q_L(m)=q_L(-m-1)$. Indeed the polynomial part
of the first Euler difference is invariant under this reflection:
projection gives $K^2L=K^3/2$, with no claim about additivity of
the virtual characters of $L$.
Let $\zeta$ be a primitive $r$th root of unity. At this Fourier
frequency only the specified constituent contributes. For
$\Phi_{b,r}$ its coefficient is the nonzero real number
\[
 c=\bigl(\zeta^{b^{-1}}+\zeta^{-b^{-1}}-2\bigr)^{-1};
\]
this follows by taking the second finite difference of
$\bar j(r-\bar j)/(2r)$. If its shift is $i$, reflection gives
$\zeta^i c=\zeta^{1-i}c$, or $2i\equiv1\pmod r$, impossible.
The divisorial sheaves used for duality are maximal
Cohen--Macaulay, as character summands of their canonical cyclic
covers. Thus the argument applies to an integral Weil class,
without a codimension-two Cartier assumption on $L$.
\end{proof}

\paragraph{Full negative parts and a zero-contact endpoint}
\begin{lemma}\label{lem:late-full-cancellation46}
Let $a\in R_m$ have full divisor $T+E$ on $Y$, with $T$ one
reduced visible prime and $E$ invisible. Suppose $z\in R_j$ vanishes
on $T$, and put $h=j-m>0$. If every integral effective invisible
divisor $J$ with $0<n=IK^2J\le IK^2E$ satisfies $IKJ^2+hn<0$,
then $z/a\in R_h$.
\end{lemma}
\begin{proof}
Take the full positive and negative parts $H,J$ of
$\operatorname{div}z-\operatorname{div}a$. Then $J\le E$, with
all multiplicities and canonical-null components retained, and
$H-J\sim hK$. Distinct effective divisors have nonnegative
intersection against the nef $K$, so
\[
 0\le IKHJ=IKJ^2+hIK^2J.
\]
Thus $IK^2J=0$. Every component of $J$ has nondivisorial image on
$Z$, because $K_Z$ is ample. The rational section $z/a$ therefore
has no divisorial pole on the normal canonical model and extends.
In particular the square of $J$ includes all cross terms with
canonical-null components; no sum of prime-square bounds is used.
\end{proof}

An invisible divisor of positive degree has strict Hodge inequality
$M<n^2/v$: equality would give positive Hodge value by null descent.
If $\mathcal R_n\subset\mathbb Z/(2I)$ is the set of square residues
from the same assignment as its degree residue, its exact necessary
upper bound is
\begin{equation}\label{eq:late-invisible-max46}
 U(n)=\max_{r\in\mathcal R_n}
 \left[r+2I\left\lfloor
 {\lfloor(n^2-1)/v\rfloor-r\over2I}\right\rfloor\right].
\end{equation}
There is no lower cutoff on $M$.

\begin{corollary}\label{lem:late-zero-contact46}
If effective integral $H,J$ on a crepant terminal model have
no common prime, $H-J\sim hK$, $KHJ=0$ and $n=IK^2J>0$,
then $\gcd(h,I)\mid n$.
\end{corollary}
\begin{proof}
Apply Lemma~\ref{lem:crepant-graph-general} to $(J,H)$ and $a=-h$.
The local characters of the resulting integral divisor are zero
or $-h$. Its degree congruence therefore makes $n$ divisible by
$\gcd(h,I)$.
\end{proof}

\paragraph{Small negative primes and integral curve budgets}
\begin{lemma}\label{lem:late-ruling46}
Let $F$ be an effective terminal prime, $s=IK^2F>0$, $N=IKF^2$.
If $s+N\le-L<0$, and $r_{\max}$ bounds the terminal indices, then
\[
 W=\left\lfloor3Is/L\right\rfloor<I/r_{\max},\qquad 2W<L
\]
are incompatible with the existence of $F$.
\end{lemma}
\begin{proof}
On a normalization resolution put $A=K|_S$. Adjunction gives
$A^2=s/I$ and $IAK_S\le-L$. For $t<L/s$ sufficiently close to
$L/s$, run the surface program for the non-pseudoeffective
$K_S+tA$. Its birational steps are point blowdowns, preserving
nefness of the pushed-forward polarization. The endpoint is a
ruling with general degree $e<2/t$ or a plane with line degree
$e<3/t$. A general line pencil in the latter case gives, after one
blowup of polarization degree zero, a rational ruling in either
case with $0<w=Ie\le W$.
Since $0<e<1/r_{\max}$, a general ruling curve has at least two
fractional markings. After a finite base change split them into
sections and compare with two disjoint constant sections of
$\mathbb P^1\times B'$. Their strict transforms are $A$-null. If
$x>0$ is the common intersection with the nef pushed-forward
polarization, the blowup weights over each section sum to $x$.
The relative canonical intersection is therefore at least
$-2x+x+x=0$. Relative ramification has nonpositive coefficient
on nonexceptional vertical components, and exceptional components
are $A$-null. Descending gives $AK_S\ge(2g(B)-2)e\ge-2e$.
This is the marked-ruling argument of
\cite[Lemma~21.4]{LBVolume}, including the plane endpoint.
Consequently $w\ge L/2$, a contradiction.
\end{proof}

\begin{lemma}\label{lem:late-intrinsic-semigroup46}
Suppose the entire member $D=(a_m=0)$ is reduced on $Z$, and let
$T$ be a visible component, of degree $n$ and square $M$ on $Y$.
If $n+M>0$, every positive-degree component of its intrinsic
polar cycle has image degree $w=IK_Z\Gamma$ satisfying
\[
 2g(\Gamma^\nu)-2+\sum(1-1/r_x)\le2(m+1)w/I,
 \qquad \sum k_x/r_x\equiv w/I\pmod{\mathbb Z}.
\]
The total $n+M$ is a nonnegative integral combination of these
degrees. The full intersection $T\cdot E$, for
$\operatorname{div}_Y(a_m)=T+E$, obeys the same assertion with
coefficient $m+1$ and total $IKTE$.
\end{lemma}
\begin{proof}
On a smooth normalization resolution the intrinsic cycle is
$\operatorname{div}(\alpha|_S)+C_T$, whose total $A$-degree is
$(n+M)/I$. It is effective and integral, and its components are
counted with multiplicities and mapping degrees. The full polar
restricts locally as $v(te)|_T=e|_T v(t)$, so the full remainder
adds zeros and cannot remove the intrinsic cycle.
At a curve where the reduced member is singular, its generic
threshold is at most one and gives adjunction coefficient at most
$m+1$. If the member is regular there, its transverse Cartier
curve in a Du Val surface forces the ambient space to be regular.
For polar order $j\ge1$, the square lift has order $2j$ on $D$;
inversion of adjunction for
$D+(2j)^{-1}\operatorname{div}(u_{2m+2})$ gives coefficient
$(m+1)(1+1/j)\le2(m+1)$. Apply
Lemma~\ref{lem:full-adjunction46} in either case.
Every curve in the full intersection maps into the singular locus
of the reduced member, including curves from an exceptional
component: a regular transverse member would leave no ambient
exceptional fibre. Thus the stronger coefficient applies to it.
Positive intersection multiplicities and finite mapping degrees are
integers. Summing proves the semigroup statements without dropping
null cross terms.
\end{proof}

\paragraph{Cubic integrality and the full conductor}
\begin{lemma}\label{lem:late-cubic-descent46}
Let $A$ be a reduced complete plane-curve ring in characteristic
zero, with normalization $\bar A=\prod k[[t_i]]$, delta invariant
$\delta$, branch multiplicities $m_i$ and conductor exponents $c_i$.
Suppose $\xi\in\bar A$, $\theta,d\in A$, and
\[
 \xi^3=\theta d,\qquad v_i(\theta)\ge c_i+m_i-1.
\]
Then $\xi^2\in A$ if $\delta\le2$, or if $\delta\le4$ and $d$
is a nonunit. Zero components of $\xi$ are allowed.
\end{lemma}
\begin{proof}
The full conductor is $\prod t_i^{c_i}k[[t_i]]$, contained in $A$,
and $\sum c_i=2\delta$. A smooth branch has
$c_i=\sum_{j\ne i}I(C_i,C_j)\le\delta$.
For a nonunit $d$, $v_i(d)\ge m_i$, so
\[
 v_i(\xi^2)\ge2\left\lceil(c_i+2m_i-1)/3\right\rceil.
\]
For $m_i\ge2$, $c_i\le8$ and this is at least $c_i$;
for $m_i=1$, use $c_i\le4$. If $\delta\le2$, the weaker
$2\lceil(c_i+m_i-1)/3\rceil$ already suffices, since the two
bounds are now $c_i\le4$ and $c_i\le2$. Thus the square lies in
the full conductor, which supplies all branch-gluing conditions.
\end{proof}

The inequality on the actual polar requires no assumption that its
ambient vector field is linear over the function field of the
singular curve. On a normalized transverse branch use coordinates
$(t_i,z)$; both transverse ambient coordinates have order at least
$m_i$. The pullback of a regular ambient two-form has order at
least $m_i-1$. Finite duality identifies
$\omega_{\bar A}=\mathfrak c\omega_A$, where $\mathfrak c$ is the
full conductor, so the residue identity adds precisely $c_i$ to
that order. A finite separable extension can split the branches;
membership descends by faithful flatness. This proves the polar
inequality used in Lemma~\ref{lem:late-cubic-descent46}.

\paragraph{The exact trace on a full complete-intersection curve}
\begin{lemma}\label{lem:late-ci-trace46}
Let two actual canonical members of degrees $a,b$ meet in a reduced
integral curve $\Gamma$ on a canonical threefold. Suppose their
local equations form a regular sequence on the local canonical
covers. Let $C$ be its normalization and $d=a+b+1$. On a cover
branch with inertia $r$ and primitive canonical mark $k/r$, the
normalization trace has an integer order $c\ge0$ satisfying
\[
 c\equiv dk+1\pmod r.
\]
The full trace gives the exact degree identity
\begin{equation}\label{eq:late-ci-trace46}
 d\deg(K_Z|_C)=2g(C)-2+
     \sum_{x\in C}\frac{r_x-1+c_x}{r_x}.
\end{equation}
At an unmarked point take $r_x=1$. If $c_x=0$ at a branch on
its full local cover, that local complete intersection is a single
smooth normal branch.
\end{lemma}
\begin{proof}
A one-dimensional reduced complete intersection $A$ is Gorenstein.
For its finite normalization $B$, finite duality gives
$\omega_B=\operatorname{Hom}_A(B,\omega_A)$. After choosing a
residue generator of $\omega_A$, evaluation at one identifies
this module with the conductor $\operatorname{Hom}_A(B,A)$.
On a normalized cover branch with uniformizer $t$, write the
image of $dt$ as $g\omega_A$. Its order $c=v_t(g)$ is therefore
nonnegative. Normalize the inertia generator to act on $t$ with
character one. The numerator character of the adjunction residue
is $-dk$, so equivariance gives $c-dk\equiv1\pmod r$.

The coarse differential is $du=r t^{r-1}dt$, for $u=t^r$.
Thus the pullback of the full dualizing line has, at that coarse
point, the excess $(c+r-1)/r$ above $\omega_C$. This is an equality
of rational divisor classes: finite duality uses the full
intersection, including every local branch and its conductor.
Taking degrees gives \eqref{eq:late-ci-trace46}.

Finally $c=0$ means that the conductor contains an element which
is a unit on that normalized branch. A local ring has a single
residue value on all its branches; multiplication by the branch
idempotents in $B$ would otherwise give an element of $A$ with
both zero and nonzero residue. Thus there is only one branch,
and the conductor is the unit ideal. Then $A=B$ is a discrete
valuation ring, so the branch is smooth. This also rules out
an additional branch of trivial inertia which would have been
missed by merely counting nontrivial marks.
\end{proof}

\begin{lemma}[Restriction to a complete-intersection curve]
\label{lem:complete-ci-canonical-restriction}
Let two actual canonical sections of degrees $a,b>0$ form a
regular sequence on the local canonical covers of a projective
canonical threefold $Z$ with $K_Z$ ample. Let $\Gamma$ be their
full intersection, and let $\mathcal F_m$ be its degree-$m$
canonical twist, obtained by taking invariants on these covers.
If $m-a$, $m-b$ and $m-a-b$ all lie outside $\{0,1\}$, then
restriction is surjective and
\[
 h^0(\Gamma,\mathcal F_m)
       =P_m-P_{m-a}-P_{m-b}+P_{m-a-b},
 \qquad P_j=0\quad(j<0).
\]
\end{lemma}
\begin{proof}
Taking invariants is exact in characteristic zero, so the local
Koszul complexes give
\[
0\longrightarrow\mathcal O_Z((m-a-b)K_Z)
\longrightarrow\mathcal O_Z((m-a)K_Z)\oplus
                 \mathcal O_Z((m-b)K_Z)
\longrightarrow\mathcal O_Z(mK_Z)
\longrightarrow\mathcal F_m\longrightarrow0.
\]
For every integer $j\notin\{0,1\}$, the first and second
cohomology of $\mathcal O_Z(jK_Z)$ vanish: use canonical
vanishing for $j\ge2$ and Serre duality for $j\le-1$.
Splitting the displayed complex into two short exact sequences
therefore proves surjectivity and the dimension formula.
\end{proof}

\subsection{Elliptic pencils and ambient residue bounds}
\begin{lemma}\label{lem:k3-elliptic-pencil-morphism}
Let $W$ be a smooth projective surface birational to a K3 surface.
A rational pencil on $W$ whose normalized general connected
curve has genus one is already a morphism after Stein factorization.
\end{lemma}
\begin{proof}
Resolve its basepoints by $W'\to W$ and let $F$ be a general
connected fibre. The zero divisor of the pulled-back K3 form is
an effective canonical divisor on $W'$, with positive coefficient
on every exceptional curve of $W'\to W$. A basepoint would
produce a horizontal such curve, giving $K_{W'}F>0$; all other
intersections with the moving fibre are nonnegative. This
contradicts $K_{W'}F=2g(F)-2=0$.
\end{proof}

\begin{lemma}[An ambient residue bound]\label{lem:ambient-residue-bound}
Let a finite group act on a normal hypersurface germ
$F=(f=0)\subset(\mathbb C^3,0)$, freely away from the origin,
and suppose its quotient has rational singularities. Write
$\omega=\operatorname{Res}(dx\wedge dy\wedge dz/f)$.
If $g\omega$ is invariant and an ambient exceptional valuation
$v$ has discrepancy $a(v)$ and satisfies $v(g)<v(f)$, then
\[
                         v(f)\le a(v)+v(g).
\]
\end{lemma}
\begin{proof}
The invariant form descends reflexively to the rational quotient
and extends to its resolution. Finite pullback and resolution
therefore extend $g\omega$ on $F$. The adjoint residue sequence
\cite[Proposition~3.1]{EinLazarsfeldTheta} says that its numerator
has a lift $g+fh$ in
\[
 \operatorname{adj}(F)=
 \pi_*\mathcal O\bigl(K_{\widetilde X/\mathbb C^3}
                     -\pi^*F+\widetilde F\bigr),
\]
where an embedded resolution may be chosen to contain the given
valuation. Thus $v(g+fh)\ge v(f)-a(v)$. Since $h$ is regular
and $v(g)<v(f)$, the left side equals $v(g)$, proving the bound.
\end{proof}



\section{The Euler-one family and its moduli}
\label{sec:euler-one}

\subsection{The numerical basket and the first sections}
\label{sec:chi1-classification}

Assume $\chi(\OO_Z)=1$. We first recover the entire canonical ring
from four low-degree sections; we then determine the resulting family
and its ordinary deformations.

\begin{theorem}[The Euler-characteristic-one equality models]
\label{thm:chi1-classification}
Let $Z$ be a canonical model.  The following conditions are equivalent:
\begin{enumerate}
\item $\chi(\mathcal O_Z)=1$ and $K_Z^3=1/420$;
\item for a weighted homogeneous polynomial
$f\in\mathbb C[x,y,z,t]_{46}$, with
$\deg(x,y,z,t)=(4,5,6,7)$, there is an isomorphism
\[
 Z\simeq X_f:=\{u^2=f(x,y,z,t)\}
       \subset\mathbb P(4,5,6,7,23),
\]
where $X_f$ is normal and has canonical singularities.
\end{enumerate}
In this description the displayed weighted coordinate ring is the
\emph{entire} canonical ring:
\begin{equation}
 R(Z,K_Z)\simeq
 \mathbb C[x,y,z,t,u]/(u^2-f),\qquad \deg u=23.
 \label{eq:chi1-canonical-ring}
\end{equation}
No quasismoothness hypothesis occurs in the implication from
\textup{(1)} to \textup{(2)}.
\end{theorem}

The distinction between a numerical basket and a canonical ring is
essential.  A Hilbert series by itself does not exclude a low-degree
relation accompanied by an extra generator.  We first establish
algebraic independence of four actual sections, and then use equality
in a multiplicity estimate to recover the full ring.

\paragraph{The basket and its torsion profiles.}

Let $g:Y\to Z$ be a projective $\mathbb Q$-factorial terminalization,
so that $K_Y=g^*K_Z$.  Chen--Chen's equality statement
\cite[Theorem~1.2(2)]{CCII} gives the basket
\begin{equation}
 \mathcal B=3(1,2)+(1,4)+(1,6)+(2,5)+(3,7).
 \label{eq:chi1-basket}
\end{equation}
Here and below the entries of a terminal basket may be virtual
constituents of a nonquotient point.  We do not treat them as distinct
geometric points.  Their least common multiple is $I=420$.
Reid's Riemann--Roch formula gives
\begin{equation}
 P_m=\frac{m(m-1)(2m-1)}{12\cdot420}-(2m-1)
       +\sum_{(b,r)\in\mathcal B}\sum_{j=1}^{m-1}
       \frac{\overline{bj}(r-\overline{bj})}{2r},\qquad m\geq2.
 \label{eq:chi1-RR}
\end{equation}
In particular,
\begin{equation}
\begin{array}{c|rrrrrrrrrrrrrrrr}
 m&1&2&3&4&5&6&7&8&9&10&11&12&13&14&15&16\\\hline
 P_m&0&0&0&1&1&1&1&1&1&2&2&3&2&3&3&4.
\end{array}
\label{eq:chi1-low-plurigenera}
\end{equation}
The value $P_1=0$ follows from $P_2=0$ by taking the square of a
canonical section.  A resolution of $Z$ is regular: otherwise
Jiang's irregular-volume bound gives $K_Z^3\geq3/8$
\cite[Theorem~1.5]{JiangSeveri}.  Thus $q(Z)=0$.

For an integral Weil divisor $D$ on $Y$, set
\[
 j(D)=420K_Y^2D,\qquad M(D)=420K_YD^2.
\]
These are integers.  The simultaneous Weil--Riemann--Roch
congruences \cite[Section~2, equations~(2.3)--(2.4)]{LBVolume} state that, for one common
assignment of local characters $i_Q\in\mathbb Z/r_Q$,
\begin{align}
 j(D)&\equiv\sum_Q\frac{420}{r_Q}b_Q^2i_Q\pmod{420},
 \label{eq:chi1-j-congruence}\\
 M(D)&\equiv-\sum_Q\frac{420}{r_Q}
       \overline{b_Qi_Q}(r_Q-\overline{b_Qi_Q})\pmod{840}.
 \label{eq:chi1-M-congruence}
\end{align}
Allowing the characters of virtual constituents to vary independently
only enlarges this necessary set.  Evaluating the $6720$ assignments
in \eqref{eq:chi1-basket} gives, for every $j\pmod{420}$,
\begin{equation}
 M\equiv j^2+210k\pmod{840},\qquad k=0,1,2,3.
 \label{eq:chi1-small-lattice}
\end{equation}
For $j=1,\ldots,10$, it follows that
\[
 0\leq M\leq j^2\quad\Longrightarrow\quad M=j^2.
\]
For $j=1,2,3$, the same conclusion holds with the weaker lower
bound $-2j\leq M$.

The Hodge-null descent of Lemma~\ref{lem:common-tools46} gives
\begin{equation}
 M(D)=j(D)^2\quad\Longrightarrow\quad
 D_Z\sim j(D)K_Z+\tau
 \label{eq:chi1-hodge-torsion}
\end{equation}
for such a torsion class $\tau$.  This implication applies to fixed
effective divisors as well as moving divisors.

Here is the required twisted-section calculation, including its
anchor.  Define
\[
 \Phi_{b,r}(i)=\frac{\overline{bi}(r-\overline{bi})}{2r},\qquad
 C(i)=\sum_Q\Phi_Q(i_Q),
\]
and
\[
 \delta_m(i)=\frac{C(i)}2+
      \sum_{j=1}^{m-1}\sum_Q\bigl(\Phi_Q(i_Q+j)-\Phi_Q(j)\bigr).
\]
Apply the weak difference Riemann--Roch formula to a sequential
pullback of $\tau$.  Addition of $mK_Z$ is compatible with that
pullback, since its crepant pullback is the integral divisor $mK_Y$.
No additivity assertion for two torsion classes is needed.
Writing $\chi_\tau=\chi(\mathcal O_Z(\tau))$, the formula at
$\tau$, together with Serre duality, gives
\[
 \chi_\tau+\chi_{-\tau}=2\chi(\mathcal O_Z)-C(\tau),
 \qquad C(-\tau)=C(\tau).
\]
Summing the canonical-shift increments and applying
Kawamata--Viehweg vanishing for $m\geq2$ gives
\[
 H_m(\tau):=h^0(Z,mK_Z+\tau)
   =P_m+\chi_\tau-\chi(\mathcal O_Z)+\frac{C(\tau)}2
       +\delta_m(i(\tau)).
\]
Thus, with $A_\tau=H_2(\tau)$,
\begin{align}
 A_\tau+A_{-\tau}
   &=2P_2+\delta_2(i(\tau))+\delta_2(i(-\tau)),
 \label{eq:chi1-paired-anchor}\\
 H_m(\tau)&=P_m-P_2+A_\tau+
          \delta_m(i(\tau))-\delta_2(i(\tau)),\qquad m\geq2.
 \label{eq:chi1-twisted-increments}
\end{align}
The codimension-two Cartier condition excludes singular-curve
correction terms in this calculation.

In the order of entries in \eqref{eq:chi1-basket}, the simultaneous
null congruences $j=M=0$ give precisely
\[
 (0,0,0,0,0,0,0),\quad(0,0,1,0,3,0,0),\quad
 (0,1,0,0,3,0,0),\quad(1,0,0,0,3,0,0).
\]
For the zero state, $C=\delta_m=0$.  For each other state, $C=1$
and
\begin{equation}
 \delta_m=\begin{cases}
  \frac12,&m\equiv1,2,3\pmod6,\\
  -\frac12,&m\equiv0,4,5\pmod6.
 \end{cases}
 \label{eq:chi1-period-six}
\end{equation}
Indeed only one index-two character and the index-six character
$3$ contribute; their consecutive increments are
$1,0,0,-1,0,0$ in the six residue classes, and $\delta_2=1/2$.
The inverse state can be any state of the same correction mass.
Equations \eqref{eq:chi1-paired-anchor}--\eqref{eq:chi1-period-six}
therefore imply
\begin{equation}
\begin{array}{c|c|c}
 \text{state}&A_\tau+A_{-\tau}&H_m(\tau),\quad m\geq2\\\hline
 0&0&P_m\\
 \text{nonzero}&1&P_m+A_\tau+\delta_m-\frac12.
\end{array}
\label{eq:chi1-twist-table}
\end{equation}
In the second row $A_\tau\in\{0,1\}$.  In particular,
\begin{equation}
 H_2(\tau)>0\quad\Longrightarrow\quad H_m(-\tau)=0
       \quad(2\leq m\leq6).
 \label{eq:chi1-inverse-vanishing}
\end{equation}

\begin{lemma}[The first pencil]
\label{lem:chi1-first-pencil}
The complete system $|10K_Z|$ has no divisorial fixed part, and
its rational map has connected general fibre.
\end{lemma}

\begin{proof}
Resolve the pencil and take its Stein factorization.  Its base is
$\mathbb P^1$, since $q(Z)=0$.  If $F$ is the pushforward of a
general connected fibre, write
\[
 10K_Z\sim\theta F+E,\qquad E\geq0,\quad\theta\geq1.
\]
Put $j=420K_Z^2F$.  On $Y$, the strict transform of $F$ is moving,
so $M\geq0$; Hodge index gives $M\leq j^2$, and
$1\leq\theta j\leq10$.  By \eqref{eq:chi1-small-lattice},
$M=j^2$, whence $F\sim jK_Z+\tau$.

For $2\leq j\leq6$, \eqref{eq:chi1-twist-table} gives
$H_j(\tau)\leq1$, contrary to movability.  If $j=1$, two
independent sections have three independent quadratic products;
this contradicts $H_2(2\tau)\leq1$.  Here the bound is applied
to the separate actual torsion class $2\tau$.
For $j=7,8,9$, one has $\theta=1$, and $H_j(\tau)\geq2$
forces a nonzero state with $A_\tau=1$.  If $j=7$ or $8$, the
effective fixed part lies in $|3K_Z-\tau|$ or $|2K_Z-\tau|$,
contradicting \eqref{eq:chi1-inverse-vanishing}.  If $j=9$, multiply
the fixed section in $|K_Z-\tau|$ by a nonzero section of
$2K_Z+\tau$ to contradict $P_3=0$.
Thus $j=10$, $\theta=1$, and $K_Z^2E=0$.  Ampleness gives $E=0$.
The equality $\theta=1$ also says that the complete pencil equals
its connected Stein pencil.
\end{proof}

\subsection{Four algebraically independent sections}

Choose nonzero $x\in R_4$, $y\in R_5$, $z\in R_6$, and $t\in R_7$.
Their section divisors on $Y$ or $Z$ will be denoted by
$D_4,D_5,D_6,D_7$, with the model indicated when necessary.

\begin{lemma}
\label{lem:chi1-products}
The products $xz,y^2$ form a basis of $R_{10}$.
On $Z$, $D_5$ has no prime component in common with $D_4+D_6$.
Moreover $xt,yz$ form a basis of $R_{11}$, and no two among
$x^3,z^2,yt$ are proportional.
\end{lemma}

\begin{proof}
Suppose $xz=cy^2$.  On $Y$, the equality $D_4+D_6=2D_5$
allows the componentwise decomposition
\[
 D_4=2A+C,\quad D_5=A+B+C,\quad D_6=2B+C,
\]
where $A,B,C$ are effective integral divisors, $A,B$ have no
common component, and $B-A\sim K_Y$.  The last relation is
obtained from $D_5-D_4$; it does not require dividing a linear
equivalence by two.  Put $a=420K_Y^2A\in\{0,1,2\}$.
If $a=0$, then $A_Z=0$ and $B_Z\sim K_Z$ is effective,
contradicting $P_1=0$.  Otherwise, proper intersection of the
two effective $\mathbb Q$-Cartier divisors $A,B$ gives
\[
 0\leq420K_YAB=M(A)+a,\qquad M(A)\leq a^2.
\]
The lattice \eqref{eq:chi1-small-lattice} forces $M(A)=a^2$.
Consequently
\[
 A_Z\sim aK_Z+\tau,\quad B_Z\sim(a+1)K_Z+\tau,
 \quad C_Z\sim(4-2a)K_Z-2\tau.
\]
If $a=1$, $B_Z$ gives $H_2(\tau)>0$, whereas $B_Z+C_Z$
gives $H_4(-\tau)>0$, contradicting
\eqref{eq:chi1-inverse-vanishing}.  If $a=2$, then $C_Z=0$
and $2\tau\sim0$; effectivity of $A_Z$ contradicts
$H_2(-\tau)=0$.  Thus $xz,y^2$ are independent.

Lemma~\ref{lem:chi1-first-pencil} now gives
$\min(D_4+D_6,2D_5)=0$ on $Z$.
If $xt=cyz$, comparison along $D_5$ makes $D_7-D_5$
effective, contradicting $P_2=0$.  A proportionality involving
$yt$ and either $x^3$ or $z^2$ contradicts the same disjointness.
Finally, $z^2=cx^3$ would make $z/x$ integral over the normal
canonical ring, of degree two, contradicting $R_2=0$.
\end{proof}

\begin{lemma}
\label{lem:chi1-first-three}
The sections $x,y,z$ are algebraically independent.
\end{lemma}

\begin{proof}
The disjoint nonzero divisors $D_4,D_5$ first show that $x,y$
are independent: a homogeneous relation between them has a
binomial factor $x^5-\lambda y^4$, which would give
$5D_4=4D_5$.
Set
\[
 f_0=xz/y^2,\qquad g_0=x^5/y^4\in\mathbb C(Z).
\]
By Lemma~\ref{lem:chi1-first-pencil}, $\mathbb C(f_0)$ is
algebraically closed in $\mathbb C(Z)$.
If $x,y,z$ were dependent, their degree-zero fraction field
would have transcendence degree one.  This assertion can be
checked by adjoining the homogeneous fraction $y/x$ of degree
one, which is transcendental over the degree-zero field.
Thus $g_0=G(f_0)$ for $G\in\mathbb C(T)$.

For every $c\in\mathbb C$, the section $xz-cy^2$ has a nonzero
effective divisor $E_c\sim10K_Z$.  The divisors $E_c$ have no
common components, are disjoint in divisorial support from $D_5$,
and for $c\ne0$ share no component with $D_4$.
The equalities
\[
 \operatorname{div}(f_0)=D_4+D_6-2D_5,\qquad
 \operatorname{div}(g_0)=5D_4-4D_5
\]
therefore show that $G$ has no zero or pole in $\mathbb C^*$
and no pole at zero.  Since $g_0$ is nonconstant,
$G(T)=\lambda T^n$ with $n>0$.  Comparison along $D_5$
gives $n=2$, hence $x^3=\lambda z^2$, already excluded.
\end{proof}

We give the differential estimate needed to control the fourth
section.  For a graded image, $H$ denotes the rational polarization
obtained by dividing a sufficiently divisible Cartier Veronese
polarization by its degree.  This convention avoids treating
ordinary powers of a non-Cartier degree-one sheaf as invertible.

\begin{lemma}[A weighted image estimate]
\label{lem:chi1-weighted-image}
Let $Z$ be a canonical model and suppose four sections of degrees
$a_0,a_1,a_2,a_3$ induce a dominant rational map onto a normal
weighted hypersurface $S_d\subset\mathbb P(a_0,a_1,a_2,a_3)$.
Assume that the evaluated partial derivatives of its equation have
no common divisorial zero on $Z$, and put $c=1+\sum a_i-d$.
Then $c>0$ and
\[
 K_Z^3\geq\frac{2d}{c\,a_0a_1a_2a_3}.
\]
\end{lemma}

\begin{proof}
Pass to a small projective $\mathbb Q$-factorialization $X\to Z$
\cite{BCHM}.  Locally on its smooth locus, write the sections
as $f_i\nu^{a_i}$ in a canonical frame $\nu$.
Contracting the residue form of $F_d=0$ with its weighted Euler
field and evaluating gives the local two-form
\[
 \frac1{F_i(f)}\sum_{j\ne i}(-1)^{\epsilon(i,j)}a_jf_j
                   \bigwedge_{k\notin\{i,j\}}df_k.
\]
The residue identities make the expressions agree.  Under a
change of frame, the derivative terms cancel by horizontality for
the Euler field.  Thus this is a nonzero global section of
$\Omega_X^{[2]}[\otimes]\mathcal O_X((\sum a_i-d)K_X)$:
it is regular at each prime divisor by the gradient hypothesis,
and extends reflexively across the remaining set.
The identification $\Omega_X^{[2]}\simeq T_X[\otimes]\omega_X$
gives a map $\mathcal O_X(-cK_X)\to T_X$ tangent to the curve
fibres.  Its saturation is an algebraically integrable rank-one
foliation $\mathcal F$, with
$cK_X-K_{\mathcal F}$ effective.

If $C$ normalizes a general leaf, the general log-leaf construction
gives an effective $\mathbb Q$-divisor $\Delta$ with
$K_C+\Delta\sim_{\mathbb Q}\nu_C^*K_{\mathcal F}$
\cite[Section~3.2, Definition~3.6]{ADLogLeaf}.
Here the ambient variety need only be normal and the foliation
$\mathbb Q$-Gorenstein; no Fano or canonical-foliation hypothesis
is being used.  A smooth general fibre on a resolution has genus
at least two: by adjunction its canonical degree equals the
intersection with the canonical class of the resolution, and
$K_W=\pi^*K_Z+E_\pi$, $E_\pi\geq0$, has positive degree on
such a moving curve.  Hence
\begin{equation}
 cK_XC\geq K_{\mathcal F}C
      =2g(C)-2+\deg\Delta\geq2.
 \label{eq:chi1-leaf-bound}
\end{equation}
In particular $c>0$.

Resolve the map defined by a divisible Veronese and put
$L=\pi^*K_Z$, $M=h^*H$.  Both divisors are nef and $L-M$
is effective, since it is the divided fixed part of that
subseries.  Thus $L^3\geq LM^2$.
If the Stein base has degree $e\geq1$ over $S_d$, then
\[
 LM^2=eH^2(LC)\geq\frac{2eH^2}{c},\qquad
 H^2=\frac d{a_0a_1a_2a_3}.
\]
This proves the estimate.
\end{proof}

\begin{lemma}
\label{lem:chi1-relations-16}
There is no polynomial relation among $x,y,z,t$ of weighted
degree at most $16$.
\end{lemma}

\begin{proof}
Degrees below $10$ have at most one monomial.  Degrees $10,11$
were handled in Lemma~\ref{lem:chi1-products}.
A relation of degree $12$ must have all three coefficients
nonzero and can be normalized to
$F_{12}=x^3+z^2+yt$.  Its affine gradient
$(3x^2,t,2z,y)$ vanishes only at the vertex, so its weighted
surface is normal.  The map onto that surface is dominant by
Lemma~\ref{lem:chi1-first-three}.  Its evaluated gradient has
no common divisorial zero because $D_4,D_5$ have disjoint supports.
Lemma~\ref{lem:chi1-weighted-image} gives
\[
 K_Z^3\geq\frac{24}{840\cdot11}=\frac1{385}>\frac1{420}.
\]
The degree-$13$ monomials $zt,x^2y$ are independent: a relation
again gives $D_7\geq D_5$ and a section of $2K_Z$.

In degree $14$, the monomials are $t^2,xy^2,x^2z$.
Their pairwise proportionalities respectively reduce to the
degree-$10$ contradiction, or make $t/y$ integral of degree two,
or $t/x$ integral of degree three.  Hence a relation would be
\[
 F_{14}=t^2+a xy^2+b x^2z,
       \qquad ab\ne0.
\]
It is irreducible as a primitive polynomial linear in $z$.
Its affine singular locus is the $z$-axis, of codimension two
in the hypersurface, so the hypersurface is normal.  Its gradient
has no common divisorial zero after evaluation: such a zero
would force $x=y=t=0$, contrary to the disjointness of
$D_4,D_5$.  The same lemma gives
$K_Z^3\geq28/(840\cdot9)=1/270$, again impossible.

In degree $15$ the monomials are $x^2t,xyz,y^3$.
At any prime of $D_4$ the first two vanish and the third does not;
its coefficient in a relation is therefore zero.  Cancellation
of $x$ reduces the remainder to degree $11$.
Finally write a degree-$16$ relation as
\[
 a y^2z+x(bz^2+cyt+dx^3)=0.
\]
If $a=0$, cancellation contradicts the established degree-$12$
independence.  If $a\ne0$, at each prime of $D_4$ put
$A=\operatorname{ord}(x)>0$, $B=\operatorname{ord}(z)$,
$C=\operatorname{ord}(t)$; $y$ is a unit there.
If $B<A$, the first term has order $B$, strictly below the
orders $A+2B,A+C,4A$ of the other terms.  Thus $B\geq A$
at every such prime.  The fraction $z/x$ then extends to a
nonzero section of $2K_Z$, a final contradiction.
\end{proof}

\begin{lemma}[The first-relation cutoff]
\label{lem:chi1-cutoff}
If $t$ is algebraic over $\mathbb C(x,y,z)$, its primitive
irreducible weighted relation has degree at most $16$.
\end{lemma}

\begin{proof}
In a local canonical frame, the polynomial Euler two-form
\[
 4f_4\,df_5\wedge df_6-5f_5\,df_4\wedge df_6
                  +6f_6\,df_4\wedge df_5
\]
gives a nonzero section of
$\Omega^{[2]}[\otimes]\mathcal O(15K)$, hence of
$T[\otimes]\mathcal O(16K)$, on a small factorialization.
There are no denominators.  The corresponding derivation
annihilates the degree-zero field of $\mathbb C[x,y,z]$ and,
in characteristic zero, its algebraic extensions.  It therefore
defines the curve foliation of the map to
\[
 S_d=\operatorname{Proj}\mathbb C[x,y,z,t]/(F_d).
\]
The same log-leaf argument gives $K C\geq1/8$.
Normality of $S_d$ is unnecessary here: its hypersurface Hilbert
series gives $H^2=d/840$, and finite birational normalization
preserves this intersection number.  Resolving a divisible
Veronese as above, with Stein degree $e$, gives
\[
 \frac1{420}=L^3\geq LM^2
       =e\frac d{840}(LC)\geq\frac{ed}{6720}.
\]
Thus $ed\leq16$.
\end{proof}

Lemmas~\ref{lem:chi1-first-three}, \ref{lem:chi1-relations-16}
and \ref{lem:chi1-cutoff} prove that $x,y,z,t$ are algebraically
independent.  This conclusion holds on every canonical equality
model with $\chi=1$.

\subsection{Multiplicity saturation recovers the full ring}

For a positively graded domain $B$ of dimension $n+1$, put
$e(B)=\lim_{s\to1^-}(1-s)^{n+1}H_B(s)$.
For a canonical ring in dimension $n$, this is $K^n$.

\begin{lemma}[Homogeneous field degree and finiteness]
\label{lem:chi1-saturation}
Let $B\subset R$ be finitely generated positive graded domains
with degree-zero part $\mathbb C$, equal dimension, and finite
fraction-field degree $d$.  Then $e(R)\geq d e(B)$.
If $B=A=\mathbb C[x_0,\ldots,x_n]$ is a weighted polynomial
ring and $e(R)=d e(A)$, then $R$ is finite over $A$.
\end{lemma}

\begin{proof}
Choose a field basis of $\operatorname{Frac}R$ over
$\operatorname{Frac}B$ consisting of homogeneous elements
$r_1,\ldots,r_d$ of $R$.  It gives an injection
$\bigoplus_iB(-\deg r_i)\hookrightarrow R$, which proves the
first assertion by comparing Hilbert series for $0<s<1$.

For the second, take any homogeneous $v\in R_b$ and put
$B=A[v]$.  Its relation over $A$ is primitive, irreducible and
weighted homogeneous, because its defining height-one prime in
$A[V]$ is principal.  If its degree in $V$ is $e$ and its
leading coefficient has degree $c\geq0$, its weighted degree
is $eb+c$, and
\[
 e(B)=(e+c/b)e(A).
\]
The first assertion applied to $B\subset R$ gives
$d e(A)\geq(d/e)(e+c/b)e(A)$, forcing $c=0$.
Thus every homogeneous generator of $R$ is integral over $A$.
\end{proof}

\begin{lemma}[A quadratic normal graded cover]
\label{lem:chi1-quadratic-cover}
If a normal graded domain $R$ is finite of field degree two
over a weighted polynomial ring $A$, then
\[
 R=A[w]/(w^2-F),\qquad \deg w=b>0,
\]
where $F\in A_{2b}$ is squarefree.
\end{lemma}

\begin{proof}
$R$ is the integral closure of $A$ in its fraction field.
The quadratic involution preserves $R$ and commutes with its
grading.  Choose a homogeneous anti-invariant $v\ne0$;
$v^2\in A$.  Factor $v^2=h^2F$ in the factorial graded ring
$A$, with $F$ squarefree and the factors homogeneous.
Then $w=v/h$ is integral over $A$, hence belongs to $R$.
The ring $A[w]/(w^2-F)$ is Cohen--Macaulay and regular in
codimension one, since the branch equation is generically simple.
It is therefore normal and is the same integral closure.
\end{proof}

\begin{proof}[Proof of Theorem~\ref{thm:chi1-classification}]
By the independence just proved,
$A=\mathbb C[x,y,z,t]\subset R$ is a polynomial subring.
If its fraction-field degree is $d$, then
\[
 \frac1{420}=e(R)\geq\frac d{4\cdot5\cdot6\cdot7}
                          =\frac d{840}.
\]
Thus $d\leq2$.  Adjoining the homogeneous fraction $y/x$ of
degree one identifies each full graded fraction field with its
degree-zero field with one transcendental variable adjoined.
Consequently $d$ is the degree of the rational map
$Z\dashrightarrow\mathbb P(4,5,6,7)$.
The case $d=1$ would make $Z$ rational, contrary to general type.
Hence $d=2$, equality holds in the multiplicity bound, and
Lemmas~\ref{lem:chi1-saturation}--\ref{lem:chi1-quadratic-cover}
give $R=A[w_b]/(w_b^2-F_{2b})$.

The projection to $\mathbb P(4,5,6,7)$ is finite.  That weighted
space is well formed, its stabilizer locus has codimension at
least two, and the $w$-vertex is absent.  Weighted adjunction
therefore gives $K_Z\sim_{\mathbb Q}\mathcal O_Z(b-22)$.
The original canonical grading gives
$\mathcal O_Z(1)=\mathcal O_Z(K_Z)$ reflexively.
Since this class is ample, $b-22=1$, hence $b=23$.

Conversely, suppose $X_f$ is normal canonical as in (2).
The same finite projection and adjunction give
$K_{X_f}=\mathcal O_{X_f}(1)$, and
\[
 K_{X_f}^3=\frac{46}{4\cdot5\cdot6\cdot7\cdot23}=\frac1{420}.
\]
The graded hypersurface ring is Cohen--Macaulay of dimension four.
Local cohomology identifies its graded pieces with
$H^0(X_f,\mathcal O_{X_f}(m))$ and gives
$H^1(X_f,\mathcal O_{X_f}(m))=
 H^2(X_f,\mathcal O_{X_f}(m))=0$ for every integer $m$.
It is therefore the canonical ring, and its degree-one piece is
zero.  Serre duality yields $\chi(\mathcal O_{X_f})=1$.
\end{proof}

\begin{remark}[Isomorphism and birational scope]
The theorem classifies canonical models up to isomorphism, with
the polynomial varying as described in the following section.
It does not identify all crepant terminalizations as varieties.
Any projective threefold with canonical singularities, of general
type, with $\chi=1$ and volume $1/420$ has one of these canonical
models and is birational to it.  Two members $X_f,X_g$ are
birational if and only if they are isomorphic, because their
canonical rings are birational invariants and their canonical
classes are ample.
\end{remark}


\subsection{The hypersurface family and its moduli}
\label{sec:chi1-moduli}

Put $A=\mathbb C[x,y,z,t]$, with weights $(4,5,6,7)$, and
let
\[
 U^{\mathrm{can}}=
 \{f\in A_{46}: X_f=(u^2=f)\subset\mathbb P(4,5,6,7,23)
                 \text{ is normal canonical}\}.
\]
Let $U^{\mathrm{qs}}$ be its quasismooth locus.  The conclusions
below concern \emph{all} of $U^{\mathrm{can}}$; quasismoothness is
used only to exhibit a nonempty open subset and to describe its
ordinary quotient singularities.

\begin{theorem}[Moduli and deformations of the Euler-one models]
\label{thm:chi1-moduli}
The set $U^{\mathrm{can}}$ is a nonempty open subset of the
$37$-dimensional vector space $A_{46}$.  It parametrizes every
canonical model with $\chi=1$ and volume $1/420$.
Two members are isomorphic precisely when
\begin{equation}
 g(\alpha x,\beta y,\gamma z,\delta t)=c^2 f(x,y,z,t)
 \quad\text{for some }\alpha,\beta,\gamma,\delta,c\in\mathbb C^*.
 \label{eq:chi1-isomorphism}
\end{equation}
Their coarse moduli space is irreducible and rational of
dimension $33$.  The general automorphism group is $\mu_2$,
generated by $u\mapsto-u$.

For every $f\in U^{\mathrm{can}}$, the ordinary Artinian
deformation functor of $X_f$ agrees with its
Koll\'ar/$\mathbb Q$-Gorenstein deformation functor.
It is unobstructed, with a smooth $33$-dimensional miniversal
base.  The completed local ring of the coarse moduli space is
\begin{equation}
 \mathbb C[[s_1,\ldots,s_{33}]]^{\operatorname{Aut}(X_f)}.
 \label{eq:chi1-local-coarse}
\end{equation}
Thus special canonical members may have extra finite
automorphisms and finite quotient singularities in coarse
moduli, but have the same miniversal dimension.
\end{theorem}

\paragraph{Existence and monomial counts.}

The base locus of the degree-$46$ branch linear system on
$\mathbb A^4\setminus\{0\}$ is the union of the four axes.
Indeed the six pairs of coordinates support the monomials
\[
 x^9y^2,\quad x^{10}z,\quad xt^6,\quad y^8z,
 \quad y^5t^3,\quad z^3t^4.
\]
Bertini applies away from the axes.  On the $x,y,z,t$ axes
the respective terms $x^{10}z,y^8z,xz^7,xt^6$ give nonzero
transverse derivatives for a general polynomial.  Hence a
general branch cone has an isolated singularity at the origin,
and $u^2-f=0$ is quasismooth.  This is the classical
$X_{46}$ format \cite{Fletcher}.

For completeness its quasismooth members are terminal.
The four coordinate points have quotient types
\[
 \tfrac14(1,3,3),\quad \tfrac15(4,2,3),\quad
 \tfrac16(5,1,5),\quad \tfrac17(5,6,2),
\]
as is seen by eliminating the variable supplied by the nonzero
axis derivative.  All satisfy the terminal age inequalities.
The $u$-vertex is absent.  On the open part of the remaining
singular line $\mathbb P(4,6)$ the branch equation is
\[
 xz(a x^9+b x^6z^2+c x^3z^4+d z^6),\qquad ad\ne0.
\]
Its three nonzero roots in $z^2/x^3$ are simple in a quasismooth
member, and give three $\frac12(1,1,1)$ points.
This also recovers the basket \eqref{eq:chi1-basket} geometrically.
In particular $U^{\mathrm{qs}}\ne\varnothing$.
The canonical condition is open in this flat weighted family
by deformation invariance of canonical singularities
\cite{KawamataCanonical}; normality is open as well.
Thus $U^{\mathrm{can}}$ is a nonempty open subset of $A_{46}$.

The monomial count can be performed without a Hilbert-series
guess.  If $t$ has exponent $d$, then the remaining exponents
satisfy $4a+5b+6c=46-7d$.  The counts are
\begin{equation}
\begin{array}{c|rrrrrrr|r}
 d&0&1&2&3&4&5&6&\text{total}\\\hline
 \#&13&8&7&4&3&1&1&37.
\end{array}
\label{eq:chi1-monomial-count}
\end{equation}
There are seven monomials of degree $23$:
\[
 yz^3,\quad y^2zt,\quad xz^2t,\quad xyt^2,
 \quad x^2y^3,\quad x^3yz,\quad x^4t.
\]
Thus the space of all weighted degree-$46$ equations in five
variables has dimension $37+7+1=45$.  On the open subset with
nonzero $u^2$ coefficient, scaling the equation and completing
the square remove respectively one coefficient and seven
linear-$u$ coefficients, leaving the $37$ branch coefficients.

Every canonical member has a nonzero square coefficient.  If that coefficient vanished, the
weight-$23$ vertex would lie on the hypersurface.  The local
cyclic cover there has equation $g_{23}+f_{46}=0$ in
$\mathbb A^4_{x,y,z,t}$.  Since $23$ is prime and the four
weights are nonzero modulo $23$, the cover is \emph{\'etale}
off the vertex.  If the quotient were normal canonical, the
cover would be normal canonical as well.  The valuation of
weights $(4,5,6,7)$ has ambient log discrepancy $22$ and
equation order at least $23$.  Hence the hypersurface pair is
not log canonical; inversion of adjunction is a contradiction.
Consequently every normal canonical degree-$46$ hypersurface
in these weights can be put in the square-normalized form.

\paragraph{Isomorphisms and quotient coordinates.}

The canonical-ring pieces in degrees $4,5,6,7$ are all
one-dimensional.  A graded ring isomorphism therefore scales
$x,y,z,t$.  In degree $23$ it sends $u$ to $cu+h_{23}$,
where $c\ne0$ and $h_{23}\in A_{23}$.  Comparing the
linear-$u$ terms of the two square-normalized relations forces
$h_{23}=0$.  This proves \eqref{eq:chi1-isomorphism}, including
its sufficiency.

Weighted scalar multiplication acts trivially on the projective
model.  Thus the relevant change-of-coordinates group is
\[
 G=(\mathbb G_m)^5/
      \mathbb G_m^{(4,5,6,7,23)},\qquad \dim G=4.
\]
Its action on branch coefficients has kernel exactly the deck
group $\mu_2$.  To verify this and rationality explicitly, use
the following reference monomials:
\begin{equation}
 z^3t^4,\qquad yzt^5,\qquad y^2z^6,\qquad xt^6.
 \label{eq:chi1-reference-monomials}
\end{equation}
Their four exponent vectors have determinant of absolute value
$46$.  Each belongs to
\[
 L_{46}=\{a\in\mathbb Z^4:(4,5,6,7)\cdot a\equiv0\pmod{46}\},
\]
which has index $46$ in $\mathbb Z^4$; hence they are a basis
of $L_{46}$.  Every other degree-$46$ exponent therefore has
integral coordinates $n_0,n_1,n_2,n_3$ in this basis, and
$\sum n_i=1$.
The associated characters $(a,-2)$ generate an index-two
sublattice of
$\ker_{\mathbb Z}(4,5,6,7,23)$, proving the assertion about
the kernel.  Equivalently, they form a character basis of the
effective torus $G/\mu_2$.

On the open subset where the four reference coefficients $c_i$
are nonzero, the $33$ expressions
\begin{equation}
 c_a\prod_{i=0}^3c_i^{-n_i}
 \label{eq:chi1-rational-parameters}
\end{equation}
for the remaining monomials are quotient coordinates.  They
normalize the reference coefficients to one.  The residual
normalizing choices act trivially on all branch coefficients.
The normalized equations supply a family over an open subset
of $\mathbb A^{33}$, and two points of this subset represent
the same model only when their displayed parameters agree.
Its automorphism group is exactly $\mu_2$.

Every canonical member has finite stabilizer in $G$.
For example, if the character support of its branch polynomial
had rank less than four, its stabilizer would contain a
nontrivial one-parameter torus acting on the projective
threefold.  The closures of its general orbits would give a
covering family of rational curves, contrary to general type.
The converse rank condition gives a finite subgroup of a torus.
This also shows directly that the only possible exceptional
automorphism groups here are finite diagonal groups.
The quotient therefore has dimension $37-4=33$ everywhere.
The dense open quotient \eqref{eq:chi1-rational-parameters}
proves its rationality and irreducibility.

We take coarse moduli in the separated canonical locus of the
stable-variety moduli stack. The results of
\cite[Theorem~5.3.6 and Proposition~6.1.2, author-version numbering]{AHStable}
identifies the canonical covering-stack formulation with Koll\'ar
families and supplies this locus with a separated coarse moduli space.
The weighted family maps to that locus. The isomorphism criterion above
identifies its geometric fibres with the $G$-orbits, and the
generator-lifting argument below gives the local deformation presentation.
This identifies the coarse quotient and its completed local rings
in the canonical-moduli setting.

As a check using the original equations, the graded ambient
automorphism group has dimension
$1+1+1+1+(1+7)=12$: there are four lower-coordinate scalings,
one scaling of $u$, and seven degree-$23$ translations of $u$.
Removing the projectively trivial scalar leaves dimension $11$.
The same count is therefore $\dim\mathbb P^{44}-11=33$.

\subsection{Ordinary deformations and the index-one cover}

For clarity, the Koll\'ar condition means that every reflexive
power of the relative dualizing sheaf is flat and commutes
with base change.  It is equivalent to deformation of the
canonical covering stack \cite[Section~5]{AHStable}.
Parameter counting alone would not establish this condition
for an arbitrary deformation over a nonreduced base.  The
following local argument supplies the needed implication.

\begin{lemma}[Extending a three-dimensional canonical cover]
\label{lem:chi1-index-cover-deformation}
Let $(S,\mathfrak m)$ be a normal Cohen--Macaulay local
$\mathbb C$-algebra of dimension three, Gorenstein on its
punctured spectrum.  Suppose its index-one cover $B$ is
Cohen--Macaulay and Gorenstein.  Every ordinary flat Artinian
deformation of $S$ is induced by an equivariant flat
deformation of this cover.  Its reflexive canonical powers
are flat and commute with Artinian base change.
\end{lemma}

\begin{proof}
Work in the complete local or analytic category.  Let $T$ be
a local Artinian base, let $S_T$ be a flat deformation, and
let $U_T$ be its punctured spectrum.  The index cover
$U'\to U$ is finite \emph{\'etale}.  Finite \emph{\'etale}
covers and their deck actions lift uniquely across nilpotent
thickenings, giving $U'_T\to U_T$.
Set $B_T=H^0(U'_T,\mathcal O_{U'_T})$.

We justify that this is the required finite flat algebra, not
merely an algebra of punctured sections.  The central covers
satisfy
\[
 H^0(U',\mathcal O_{U'})=B,\qquad
 H^1(U',\mathcal O_{U'})=H^2_{\mathfrak m_B}(B)=0,
\]
by depth three.  In a small extension of Artinian bases with
kernel $J$ annihilated by the maximal ideal, the exact sequence
of the lifted structure sheaves consequently gives
\begin{equation}
 0\longrightarrow J\otimes_{\mathbb C}B
 \longrightarrow B_{T'}\longrightarrow B_T\longrightarrow0.
 \label{eq:chi1-cover-small-extension}
\end{equation}
Induction and the local flatness criterion imply flatness over
the Artinian base and identify the reduction with $B$.
Lifting finitely many generators of the finite central
$S$-module $B$ shows that $B_T$ is finite over $S_T$:
the cokernel is equal to the maximal ideal of $T$ times
itself and is zero by nilpotence.  Localization on a finite
affine cover of the punctured spectrum shows that
$\operatorname{Spec}B_T$ restricts to $U'_T$.
The same exact-sequence argument, or the flatness criterion
with an arbitrary quotient of $T$, proves base-change
compatibility.  On the coarse side, depth three similarly
gives $H^0(U_T,\mathcal O_{U_T})=S_T$.
Taking deck invariants, an exact operation in characteristic
zero, therefore gives $(B_T)^{\mu_r}=S_T$.

The flat deformation $B_T$ has Gorenstein fibres and is a
Gorenstein morphism over $T$, so its relative dualizing module
is invertible.  Choose a generator whose reduction is the
canonical generator of the central index cover.  Projecting
it to the corresponding deck-character eigenspace produces
an eigenvector which is still a generator, by Nakayama.
Its character is the original faithful canonical character.
On the punctured spectrum the quotient is \emph{\'etale}, so
the relative canonical powers of $S_T$ are the corresponding
character summands of $B_T$.  Extending across the closed
point gives the same equality everywhere: both sides are
extensions of their punctured restrictions. Indeed $S_T$ has depth
three. Every double dual over $S_T$ has depth at least two:
a two-term regular sequence on $S_T$ remains regular on
$\operatorname{Hom}_{S_T}(M,S_T)$ for every finite module $M$.
Thus each reflexive canonical power is recovered by extension from
the punctured spectrum. The character summands are flat with
depth-three central fibres, and have the same extension property.
Each summand of $B_T$ is flat
and commutes with base change.  These are precisely the
Koll\'ar conditions.
\end{proof}

For our application the hypotheses are especially concrete.
A canonical threefold is Gorenstein in codimension two, since
its transverse canonical surface singularities are Du Val;
its non-Gorenstein locus is therefore finite.  Canonical
singularities and their index-one covers are rational and
Cohen--Macaulay.  On a weighted hypersurface chart the
index-one cover is a Gorenstein hypersurface (or is smooth),
and its canonical character is faithful because
$K_{X_f}=\mathcal O_{X_f}(1)$.
At Gorenstein points every flat deformation remains relatively
Gorenstein.  Applying Lemma~\ref{lem:chi1-index-cover-deformation}
at the finitely many remaining points proves that every
ordinary Artinian deformation of $X_f$ is admissible.

\paragraph{Completeness and unobstructedness.}

Let $\mathcal X\to\operatorname{Spec}T$ be any such deformation
of $X_f$.  The preceding lemma gives base change and flatness
for all $\omega_{\mathcal X/T}^{[m]}$.
The vanishing
\[
 H^1(X_f,\mathcal O_{X_f}(m))=
 H^2(X_f,\mathcal O_{X_f}(m))=0\qquad(m\in\mathbb Z)
\]
proved from the hypersurface ring implies that
\[
 \mathcal R_m=H^0(\mathcal X,\omega_{\mathcal X/T}^{[m]})
\]
is free over $T$ and commutes with reduction.  Lift the five
central generators in degrees $4,5,6,7,23$.
Graded Nakayama in each degree gives a surjection
\[
 T[x,y,z,t,u]\longrightarrow\mathcal R.
\]
The degree-$46$ kernel is free of rank one.  Its generator
reduces to $u^2-f$, so its square coefficient is a unit.
Normalize that coefficient and complete the square to obtain
$u^2-f_T$, with $f_T\in T[x,y,z,t]_{46}$.
The quotient by this monic polynomial is free in each degree
with the same ranks as $\mathcal R$; the induced graded
surjection is therefore an isomorphism.
Relative Proj recovers $\mathcal X$, since a Cartier multiple
of its relative canonical divisor remains ample.

Conversely every lift $f_T$ produces a flat graded
hypersurface.  Adjunction on the weighted quotient stack
gives the relative canonical line bundle $\mathcal O(1)$;
exactness of finite cyclic invariants gives flatness and
base change of its powers on the coarse family.  Thus all
coefficient lifts are admissible.  In particular they lift
through every small extension, proving unobstructedness.

The $37$-dimensional coefficient space has four independent
infinitesimal coordinate directions.  Indeed a positive
dimensional stabilizer would contradict the finite-stabilizer
argument above.  At any fixed canonical polynomial, choose
four nonzero coefficients whose characters are independent.
Their values can be fixed formally: the associated torus
isogeny is \emph{\'etale} in characteristic zero.  This
constructs a smooth formal slice of dimension $33$.
All residual identifications are the finite automorphisms of
the central model.  This proves
\eqref{eq:chi1-local-coarse} and the deformation assertions of
Theorem~\ref{thm:chi1-moduli}.
The deck involution acts trivially on the coefficient space
and hence on every deformation parameter.  On the reference
open set \eqref{eq:chi1-reference-monomials}, it is the full
automorphism group, so the coarse moduli space is smooth there.
Equivalently, the tangent space is the explicit Jacobian quotient
\[
 T^1(X_f)=
 A_{46}/\langle xf_x,yf_y,zf_z,tf_t\rangle.
\]
The $u$-derivative removes the linear-$u$ directions, and the
Euler identity $46f=4xf_x+5yf_y+6zf_z+7tf_t$ removes the
separate equation-scaling direction.  The four remaining
vectors are independent by finiteness of the stabilizer, so
this quotient has dimension $37-4=33$ at every canonical
member, including those outside the quasismooth locus.

Finally, Theorem~\ref{thm:chi1-classification} makes this a
global classification of the canonical equality locus with
$\chi=1$, not only a local description near its quasismooth
members.  In a flat admissible family with a canonical
central fibre, $\chi$ is constant and the canonical Hilbert
polynomial, hence the volume, is constant.  Any canonical
limit of these equality models is again an $X_f$ in
$U^{\mathrm{can}}$.  Locally near that limit the preceding
generator-lifting proof supplies the weighted equation of the
family itself.  This assertion does not classify the
noncanonical or nonnormal boundary in a stable-variety
compactification, and it makes no claim that the canonical
moduli space is proper.



\section{Uniform exclusions and the remaining eleven baskets}
\label{sec:uniform-exclusions}

We now assume $\chi=2$. The preliminary reduction gives
$p_g=q=0$ and $h^{2,0}=1$. All divisors below are taken on the
terminal model unless a subscript $Z$ indicates their image on the
canonical model. We first apply the general numerical tests and then
use the common lemmas to exclude seventeen further baskets. The
remaining possibilities are collected at the end of the section.

\subsection{The initial geometric reduction}

\begin{proposition}\label{prop:initial-cuts46}
The tests above and the commuting-degree test in
Section~\ref{sec:chi3} exclude $37$ of the $65$ equality
baskets with $\chi=2$. The $28$ rows not excluded by these tests
are exactly those in Table~\ref{tab:remaining28}.
\end{proposition}
\begin{proof}
The complete finite calculation is specified as follows. Use
\eqref{eq:RR46} through degree $120$. Test
\eqref{eq:complete-system-volume46} for $2\le a\le60$ and
\eqref{eq:hilbert-surface46}--\eqref{eq:hilbert-pencil-volume46}
for $2\le a\le30$, $2\le b$, $ab\le60$.

Put $\delta=\min\{m:P_m\ge2\}$. For each $\delta\le m\le20$ with $P_m\ge2$, enumerate
$1\le\theta<m$, $1\le k\le\lfloor r(1+r)^2V\rfloor$,
$0\le p\le\lfloor k/2\rfloor+2$, $0\le q\le p$.
Apply \eqref{eq:pencil-area46}--\eqref{eq:kernel-area46} and
\eqref{eq:profile-trace46} for $\ell=2,3,4$; whenever
\eqref{eq:profile-defect46} fails, apply the same area tests with
\eqref{eq:repeated-cost46}. Emptiness of this explicit integer list
is the required surface-profile obstruction. Table~\ref{tab:initial-cuts46}
provides a witness for every removed row. All comparisons are strict
rational comparisons; the companion source contains the complete
lists, including overlapping witnesses.

There are five differential witnesses. For C2-00015, C2-00017,
and C2-00021 the degree-three and degree-five sections commute,
because $P_9=1$ makes their bracket a multiple of the cube of the
degree-three section, which is killed by the same derivation; the
no-rational-time lemma therefore kills the first derivative;
their degree-fifteen powers are independent by normality and $P_1=0$.
The commuting-section estimate
\cite[Theorem~5.7 and Corollary~5.8]{LBVolume} gives $V\ge1/240$.
For C2-02067 use degrees four and six: $P_{11}=0$ kills their
bracket, and the degree-twelve pencil gives $V\ge1/312$.
For C2-02356 use degrees five and seven: $P_{13}=0$ kills their
bracket, and their independent degree-thirty-five powers would
require $0<\theta j\le35(v-1)=0$, since $v=1$.
These also are impossibilities at $V=1/420$.
\end{proof}


\subsection{Prime members and visible-component bounds}
\paragraph{C2-02346: a torsion-twisted prime.}

\label{sec:twisted-prime-restriction34}
\begin{theorem}
\label{thm:c2346-twisted-prime34}
The basket \texttt{C2-02346} has
no canonical realization.
\end{theorem}

\begin{proof}
Assume that a realization $Z$ exists and take its terminal
$\mathbb Q$-factorial minimal model $Y$.  The relevant data
are $I=420$, $v=1$, $K_Z^3=1/420$, $\chi=2$,
$p_g=q=0$, and $h^{2,0}=1$.  In particular,
\[
 P_1=P_2=P_3=P_5=P_{13}=0,\quad
 P_4=P_6=P_{10}=1,\quad P_{19}=3,\quad P_{20}=5.
\]
The finite paired-torsion calculation
has $51$ necessary profiles.  In every profile and
both orientations, $H_{12}(\sigma)\ge1$.  An effective
$K_Z+\sigma$ would multiply a nonzero section of
$12K_Z-\sigma$ into $R_{13}=0$.  Thus no such effective
degree-one torsion translate exists.

For a visible prime $B$ of scaled degree
$1\le n=420K_Y^2B\le4$, the adjunction bound,
Hodge bound, and simultaneous exact Weil--RR congruences
give exactly
\[
 (n,M)=(1,1),(2,4),(3,9),(4,16),
 \qquad M=420K_YB^2.
\]
All are Hodge equalities.  Their null-divisor descents are
$B_Z\sim nK_Z+\sigma$; degree one has just been excluded.

Choose a nonzero quartic canonical section $s_4$.  Every
visible prime in its actual terminal divisor is among
these Hodge states.  The visible-saturation lemma
(Lemma~\ref{lem:visible-saturation34}) shows that its
entire canonical divisor downstairs is a sum of these
visible components, with no invisible remainder.

The complete degree-five restriction to the full quartic
member has dimension
\[
 h^0(\omega_{D_4})=P_5-P_1+h^{2,0}=1.
\]
Ordinary dualizing trace from a smooth normalization
model of each visible prime contributes a nonzero section
of this space.  Contributions from distinct primes are
independent at their separate generic points, also when
their multiplicities in $D_4$ exceed one, by ordinary dualizing trace
\cite[Lemmas~8.1 and~9.6]{LBVolume}.
Hence there is only one distinct visible prime.
Its integral canonical degree is at least two and its
multiplicity times that degree is four.  Therefore
\[
 D_4=T_4\quad\hbox{or}\quad D_4=2T_2.
\]

The first case is impossible by the untwisted prime
complete-restriction algebra: its degree-five space has
dimension one, while its degree-ten space has dimension
$P_{10}-P_6=0$.  The square of a nonzero degree-five
section cannot vanish on the integral prime.

Consequently $D_4=2T$, where $T$ is reduced and prime,
and the actual Hodge descent gives
\[
 T\sim2K_Z+\tau,\qquad 2\tau=0.
\]
The class is nontrivial because $P_2=0$.  Its defining
section has $H_2(\tau)>0$, and $-\tau=\tau$ also gives
$H_2(-\tau)>0$.  Among the complete $51$ necessary
paired profiles, exactly one has these two positive
dimensions.  Its two orientations coincide.  The only
entries needed are
\[
 H_3(\tau)=0,\qquad H_{17}(\tau)=2,\qquad
 H_{22}(\tau)=5.
\]

Apply Lemma~\ref{lem:twisted-prime-restriction34} to the
actual prime section of $2K_Z+\tau$.  Its complete
restriction pieces satisfy
\[
 G_3(\tau)=H_3(\tau)-P_1+1=1,
\]
\[
 G_{19}(0)=P_{19}-H_{17}(-\tau)=3-2=1,
 \qquad
 G_{22}(\tau)=H_{22}(\tau)-P_{20}=5-5=0.
\]
The product of nonzero homogeneous sections of the
first two spaces would be a nonzero section of the
last one, a contradiction.
\end{proof}


\paragraph{C2-02351: a saturated quintic member.}

\label{sec:hodge-visible-partition34}
\begin{theorem}\label{thm:c2351-visible-partition34}
The basket C2-02351 has no canonical realization.
\end{theorem}
\begin{proof}
Here $I=420$, $v=1$, and
\[
 P_1=P_2=P_3=P_4=P_{13}=0,\qquad
 P_5=P_6=P_9=P_{10}=1,\qquad P_{15}=2.
\]
All forty-three paired torsion profiles have $H_{12}(\pm\tau)\ge1$.
An effective $K_Z+\tau$ would therefore give a nonzero section
of $R_{13}$, so $H_1(\tau)=0$ for every actual class.

The complete visible RR list through degree five is
$(n,M)=(n,n^2)$ for $1\le n\le5$. For $0\ne s_5\in R_5$,
visible saturation (Lemma~\ref{lem:visible-saturation34}) removes
the entire invisible complement on $Z$. Its component degrees
are at least two, so the only partitions are $[5]$ and $[2,3]$.
The second would give one actual class with
$H_2(\tau)H_3(-\tau)>0$, absent from all forty-three profiles.
Thus the divisor of $s_5$ is a single reduced prime.

Its complete restriction algebra has
$G_n=P_n-P_{n-5}+\mathbf1_{n=6}$, in particular
\[
 G_6=2,\qquad G_9=1,\qquad G_{15}=1.
\]
The homogeneous product inequality of
Lemma~\ref{lem:twisted-prime-restriction34} gives
$G_{15}\ge G_6+G_9-1=2$, a contradiction.
\end{proof}


\paragraph{C2-00005 and C2-02074: one reduced visible component.}

\label{sec:single-visible-adjunction34}
\begin{theorem}
\label{thm:c205-single-visible34}\label{thm:single-visible-pair34}
Neither C2-00005 nor C2-02074 has a canonical realization.
\end{theorem}
\begin{proof}
The relevant data and the complete small-degree visible states are
\[
\begin{array}{c|rr|r|rr|l|r}
 B&I&v&e&P_{e+1}&P_{2e+2}&(n,M)\text{ for }n<f&f\\\hline
 00005&1260&3&4&1&2&(3,3),(6,12)&7\\
 02074&420&1&6&1&2&(1,1),(2,4),(3,9)&4.
\end{array}
\]
There are respectively one and two necessary paired torsion profiles.
Each gives $H_2(\tau)=H_3(\tau)=0$ for every actual class.
Squaring a hypothetical section of $K_Z+\tau$ and applying the
degree-two vanishing to $2\tau$ also gives $H_1(\tau)=0$.
Every state in the table is a Hodge equality and descends to one
of these forbidden effective translates. Lemma~\ref{lem:visible-prime-floor}
therefore gives the indicated floor $f$.

In each case $P_e=1$ and $ev<2f$, namely $12<14$ and $6<8$.
The full visible part of the first member is consequently one
reduced prime. Since $P_{e+1}=1$, Lemma~\ref{lem:single-visible-square34}
gives $P_{2e+2}\ge3$, contradicting the table.
\end{proof}


\paragraph{C2-01826: canonical multiplication on torsion translates.}

\label{sec:canonical-twist-filter34}
\begin{theorem}\label{thm:c1826-canonical-twist34}
The basket C2-01826 has no canonical realization.
\end{theorem}
\begin{proof}
Here $I=420$, $v=1$,
\[
 P_1=P_2=P_3=P_4=P_5=0,\qquad
 P_6=P_7=P_8=1,\qquad P_{14}=2.
\]
The paired-torsion enumeration gives three profiles. Besides the
canonical profile, each has an orientation $\sigma=\pm\tau$
with $H_2(\sigma)=1$ and $H_8(\sigma)=0$. Multiplication by
$0\ne a_6\in R_6$ excludes both. Thus every actual torsion
class has $H_j(\tau)=0$ for $2\le j\le5$; squaring gives the
same conclusion for $j=1$.

The visible RR states with $1\le n\le5$ are exactly
$(n,M)=(n,n^2)$. Null descent and the preceding vanishings
therefore give the floor $f=6$ in
Lemma~\ref{lem:visible-prime-floor}. The divisor of $a_6$ is one
reduced prime on $Z$. Since $P_7=1$, the single-visible-component
lemma gives $P_{14}\ge3$, a contradiction.
\end{proof}


\paragraph{C2-02078: a prime quintic and its ordinary quotient.}

\label{sec:c2078-prime-section34}
Suppose that C2-02078 has a canonical realization. Here $I=420$,
$v=1$, and
\[
 P_2=P_3=P_4=0,\qquad P_5=P_7=P_8=1,
 \qquad P_{10}=P_{15}=2.
\]
The torsion enumeration leaves a single paired profile, with
$H_j(\tau)=0$ for $j=2,3,4$. Squaring gives $H_1(\tau)=0$
for every actual class. The visible RR states through degree four
are exactly $(1,1),(2,4),(3,9),(4,16)$, so
Lemma~\ref{lem:visible-prime-floor} gives the floor five.
Hence $0\ne s\in R_5$ has a single reduced prime divisor on $Z$.

The ordinary quotient $R/(s)$ is a domain by the normal-division
argument in Section~\ref{sec:twisted-prime-toolbox}. Its relevant
dimensions are
\[
 \dim(R/(s))_7=P_7-P_2=1,\qquad
 \dim(R/(s))_8=P_8-P_3=1,\qquad
 \dim(R/(s))_{15}=P_{15}-P_{10}=0.
\]
The product of the first two nonzero homogeneous pieces gives a
contradiction. This excludes C2-02078.


\paragraph{C2-02348: linked visible partitions.}

\label{sec:c2348-visible-partition34}

\begin{theorem}
\label{thm:c2348-visible-partition34}
The basket \texttt{C2-02348} has no
canonical realization.
\end{theorem}

\begin{proof}
Its data are
\[
 I=840,\qquad v=2,\qquad K^3=1/420,\qquad\chi=2,
\]
\[
 P_1=P_2=P_3=P_4=P_7=0,
 \qquad P_5=P_6=P_{13}=1,\qquad P_{12}=2.
\]
Tiny volume gives $q=0$, and therefore $h^{2,0}=1$.
Choose $0\ne s_5\in R_5$.

For a visible prime component $T$ of its actual divisor
on the terminal model, the scaled degree satisfies
$1\le n=IK^2T\le10$.  Adjunction and Hodge give
$-n\le M=IKT^2\le n^2/2$.  The exact simultaneous
RR congruences permit precisely
\[
 (n,M)=(2,2),(4,8),(6,18),(8,32),(10,50).
\]
Every state is a Hodge equality, and $n$ is even.
Consequently the actual null-divisor descent has
the integral form
\[
 T_Z\sim (n/2)K_Z+\tau.
\]
Lemma~\ref{lem:visible-saturation34} applies.  The
visible canonical degrees exhaust five and the
invisible complement vanishes downstairs.

If $D_5=\operatorname{div}_Z(s_5)$ were not one
reduced prime, its positive integral component
degrees could be grouped into $1+4$ or $2+3$.
This grouping also covers repeated components.
The corresponding two effective sums have actual
opposite torsion classes, since their sum is
linearly equivalent to $5K_Z$.

The torsion enumeration has $78$ necessary paired
profiles.  None has
$H_2(\tau)>0$ and $H_3(-\tau)>0$, so the split
$2+3$ is impossible.  For the split $1+4$, write
the actual effective classes as
$K_Z+\tau$ and $4K_Z-\tau$.  Exactly $22$ of the
necessary profiles have $H_4(-\tau)>0$, and
each has $H_{12}(-\tau)\ge2$.  Multiplication by
the same effective degree-one translate injects
\[
 H^0(Z,12K_Z-\tau)\hookrightarrow R_{13},
\]
contradicting $P_{13}=1$.  This is a linked-pair
argument; it does not assert that every effective
degree-one torsion translate is absent by itself.

Thus $D_5$ is one reduced prime.  Its actual
complete restriction ring is a domain with
\[
 G_n=P_n-P_{n-5}+\mathbf1_{n=6}.
\]
Here $G_6=2$ and $G_{12}=2$.  A two-dimensional
subspace of a domain has square span of dimension
at least three.  The square of the degree-six
piece therefore contradicts the degree-twelve
dimension two.
\end{proof}


\subsection{Adjacent gaps and the full invisible remainder}

\label{sec:two-anchors46}
\begin{theorem}\label{thm:c2071-c2077-excluded46}
Neither C2-02071 nor C2-02077 has a canonical realization.
\end{theorem}
\begin{proof}
Both have $I=420$, $v=1$. For C2-02071 use $e=4$; its two
adjacent gaps are $P_5=P_6=0$. The simultaneous visible RR list
through degree four is exactly
\[
 (n,M)=(1,1),(2,4),(3,9),(4,16).
\]
Each is a Hodge equality and hence has Hodge mass $n$ by null
descent. Lemma~\ref{lem:adjacent-gaps46} forces mass four, so
the visible prime has degree four and exhausts the entire member
on $Z$. Complete restriction to this reduced prime gives
\[
 (G_7,G_{10},G_{17})=(1,2,1).
\]
The domain product inequality requires $G_{17}\ge G_7+G_{10}-1=2$,
which is impossible.

For C2-02077 the basket is
\[
 4(1,2)+(1,3)+(1,5)+2(2,5)+(5,12)+(5,14).
\]
Use $e=5$ and $P_6=P_7=0$. The raw visible RR list through
degree five is
\[
 (1,1),(2,4),(3,9),(4,2),(4,16),(5,11),(5,25).
\]
Hodge mass five removes the equality states of degrees at most
four. Thus in $(a_5)_Y=T+E$ the visible prime has $n(T)\ge4$
and the full invisible remainder has $IK^2E\le1$.
Let $\mathcal Q_j$ be the full quotient by $a_5$. Its
degree-six connecting direction is nonzero generically on $T$,
by Lemma~\ref{lem:common-tools46}, whereas
\[
 h^0(\mathcal Q_6)=1,\qquad h^0(\mathcal Q_{15})=0,
 \qquad P_9>0.
\]
Multiplication by a nonzero $z\in R_9$ therefore shows that
$z$ vanishes along $T$.

For an arbitrary effective invisible integral divisor $J$ with
$IK^2J=1$, Hodge equality would force nonzero Hodge mass. Thus
$M(J)<1$. The exact degree-one RR residues have largest
representative below one equal to $-13$; consequently
\begin{equation}\label{eq:c2077-invisible46}
 IKJ^2\le-13.
\end{equation}
No lower adjunction bound is imposed on an invisible divisor in
this calculation. This distinction is essential.

The bound applies to every integral effective invisible divisor of
degree at most $IK^2E\le1$, with its full square and all null
components retained. Since $-13+4<0$, the full-cancellation lemma
(Lemma~\ref{lem:late-full-cancellation46}) applies with $h=9-5=4$.
It gives $z/a_5\in R_4=0$, a contradiction.
\end{proof}


\subsection{Moving systems and linked torsion classes}

\begin{lemma}[The first degree-ten pencil]
\label{lem:first-degree-ten-common}
For any of C2-02068, C2-02073, C2-02076 and C2-02080, the first
complete pencil $|10K|$ has Stein degree one. If $T$ is its integral
moving divisor, put $j=420K^2T$, $M=420KT^2$, and let
$(k,p,q)$ be the invariants of its smooth minimal surface. Then
$k=1$ or $2$, and restriction of the ambient two-form is injective.
The possible moving degrees are $j=8,9,10$; for $j=8,9$ one has
$M=j^2$ and the actual moving and fixed parts satisfy
\[
 T_Z\sim jK_Z+\tau,\qquad E_Z\sim(10-j)K_Z-\tau.
\]
For C2-02068 and C2-02076, more precisely, the only possibilities are
\[
 (k,p,q,j,M)=(1,1,0,8,64)
 \quad\text{or}\quad (k,p,q,j)=(2,1,1,10).
\]
The first case gives an actual fixed-free subsystem of $|8K_Z+\tau|$
of dimension one; in the second, the complete $|10K_Z|$ has no
divisorial fixed part.
\end{lemma}
\begin{proof}
Each basket has $I=420$, $v=1$, $P_{10}=2$ and $P_m\le1$ for
$m<10$. Stein degree at least two would give
$k\le5\cdot6^2/420=3/7<1$. Thus the degree is one, and
$k\le10\cdot11^2/420=121/42<3$. A kernel in the restriction
of the ambient form would give $k\le10\cdot21/420=1/2$, also
impossible.

Fractional adjunction gives
$3(j+M)^2\ge420j(3k+1)$, with $0\le M\le j^2$.
For $k=1$ even the Hodge maximum fails when $j\le7$; for
$k=2$ it fails when $j\le9$, since
$3(9+81)^2=24300<26460$. In each of the four baskets, the
simultaneous RR congruences leave only $M=64$ at $j=8$ and
$M=81$ at $j=9$. Null descent gives the displayed linked
classes. If $j=10$, ampleness kills the canonical fixed part.

For C2-02068 and C2-02076 one also has $P_{11}=1$ and
$P_{22}-P_{12}=3$. Thus $k+1+p-q\le3$; together with
injective restriction and surface geography, this leaves only
$(1,1,0)$ and $(2,1,1)$. In the regular case the full-member
trace defect is $\max(0,P_{11}-p+1)=1$, while the available
second-plurigenus difference is exactly three. The member must
therefore contain a repeated fixed prime. Its degree cost is at
least two, giving $j\le8$ and hence $(j,M)=(8,64)$.
The irregular case has $j=10$ by the preceding area bound.
\end{proof}


\paragraph{C2-01823 and C2-01827: the first degree-nine pencil.}

\label{sec:degree-nine-equality34}

The two baskets C2-01823 and C2-01827 have
$I=420$, $V=1/420$, $\chi=2$, $P_9=2$, and $P_m\le1$ for
$m<9$. We use the first-system adjunction, exact Weil
RR lattice, and null-divisor descent on the canonical model $Z$.
Their hypotheses hold here: in particular $q(Z)=p_g(Z)=0$.
Write $T$ for the integral moving first-surface divisor on a
terminal model, $j=420K^2T$ and $M=420KT^2$.

The first pencil has multiplicity one. Indeed multiplicity at
least two gives $t\le9/2$ and $t(t+1)^2V<1$, contradicting
$k=K_{S_0}^2\ge1$. Also $t(t+1)V<1$, so the fractional
adjunction estimate gives $k+1/3\le15/7$, hence $k=1$.
Consequently $j\le9$, $M\le j^2$, and
\[
 3(j+M)^2\ge4j\cdot420.
\]
Even the Hodge maximum fails for $j\le7$, since
$3\cdot7(7+1)^2=1344<1680$. Thus $j=8$ or $9$.
The simultaneous RR congruences leave only $M=64$ at $j=8$.
In that case null-divisor descent supplies a single linked actual
torsion class with
\begin{equation}
 T_Z\sim8K_Z+\tau,\qquad E_Z\sim K_Z-\tau,
 \qquad H_8(\tau)\ge2.
 \label{eq:degree-nine-linked34}
\end{equation}
The subsystem defining $T_Z$ has no divisorial fixed part. If
$j=9$, then $K_Z^2E_Z=0$ and ampleness gives $E_Z=0$; this
argument does not require Hodge equality for $T$.

Both baskets violate the same necessary quotient inequality for
a complete fixed-free degree-nine system. For a nonzero section
$s_3$ and a generic $u_9$ coprime to its divisor, normality gives
\[
 \frac{R_{14}}{s_3R_{11}}
 \lhook\joinrel\longrightarrow
 \frac{R_{23}}{s_3R_{20}}.
\]
The source has dimension $3-1=2$, whereas the target has
dimension $7-6=1$ for C2-01823, and $8-7=1$ for C2-01827.
In particular $j=9$ is impossible for either basket.

\paragraph{C2-01823: section transport excludes the smaller moving part}

For the $j=8$ alternative, the finite enumeration of
the $41$ necessary inverse profiles finds exactly three with
$H_8(\tau)\ge2$. Their degree-twelve dimensions are respectively
$3,2,3$. Thus every such profile has $H_{12}(\tau)\ge2$.
Multiplication by the same actual effective divisor
$E_Z\sim K_Z-\tau$ in \eqref{eq:degree-nine-linked34} injects
\[
 H^0(Z,12K_Z+\tau)\hookrightarrow R_{13}.
\]
This contradicts $P_{13}=1$. Hence C2-01823 has no canonical
realization.

\paragraph{C2-01827: quotient failure forces a prime into the fixed part}

Here the $j=8$ alternative cannot be removed by the preceding
degree-twelve comparison alone. Instead the effective fixed
part $E_Z$ has $420K_Z^2E_Z=1$. Every prime divisor $D$ on $Z$
has $420K_Z^2D$ a positive integer: take its strict transform on
the terminal model for integrality, and use ampleness
on $Z$ for positivity. Therefore $E_Z$ is one reduced prime.

Let $D_3=\operatorname{div}(s_3)$. If $E_Z$ were not a component
of $D_3$, choose a member of the moving subsystem in
\eqref{eq:degree-nine-linked34} avoiding all primes of $D_3$.
Its product with the fixed section of $E_Z$ is an actual
degree-nine section coprime to $D_3$. The quotient injection
above would give $2\le1$, a contradiction. Hence $E_Z$ is a
component of $D_3$. Since $E_Z$ is reduced,
\[
 D_3-E_Z\ge0,\qquad
 D_3-E_Z\sim2K_Z+\tau,
 \qquad H_2(\tau)\ge1.
\]
On the other hand multiplication by $E_Z$ sends
$H^0(3K_Z+\tau)$ into $R_4=0$, so $H_3(\tau)=0$.
The finite enumeration of all $87$ necessary inverse profiles
shows
\[
 H_8(\tau)\ge2\quad\hbox{and}\quad H_3(\tau)=0
 \quad\Longrightarrow\quad H_2(\tau)=0.
\]
Exactly three profiles satisfy the antecedent; all three
have $H_2(\tau)=0$. This contradicts the effective residual
$D_3-E_Z$. Hence C2-01827 has no canonical realization.


\paragraph{C2-02073 and C2-02080: fixed-free degree-ten pencils.}

\label{sec:c2073-c2080-fixed-free46}

By Lemma~\ref{lem:first-degree-ten-common}, only $j=8,9,10$
can occur. For $j=8,9$ the moving and fixed classes are linked by
$T_Z\sim jK_Z+\tau$ and $E_Z\sim(10-j)K_Z-\tau$.

The torsion enumeration gives the following implications:
\[
\begin{array}{c|c|c|c}
\text{basket}&\text{admissible profiles}&j=8&j=9\\\hline
\text{C2-02073}&79&H_8(\tau)\ge2\Rightarrow H_2(-\tau)=0
 &H_9(\tau)\ge2\Rightarrow H_3(\tau)\ge1\\
\text{C2-02080}&17&H_8(\tau)\ge2\Rightarrow H_2(-\tau)=0
 &H_9(\tau)\ge2\Rightarrow H_5(\tau)\ge1
\end{array}
\]
At $j=8$, effectivity of the same $E_Z\sim2K_Z-\tau$ contradicts
the inverse degree-two vanishing. At $j=9$, multiplication by
the same effective section of $K_Z-\tau$ contradicts $P_4=0$
for C2-02073 or $P_6=0$ for C2-02080, respectively.
Only $j=10$ remains. Its effective canonical fixed part has
degree zero and therefore vanishes. This establishes divisorial
fixed-freeness of the complete $10$-canonical system in both cases.

Apply Proposition~\ref{prop:quotient-correction34} for $\chi=2$
and $q=p_g=0$.
For a section of degree $l$ and a generic degree-ten section
coprime to its divisor, that lemma requires
\[
 P_{l+11}-P_{11}-P_{l+1}+P_1\ge1.
\]
For C2-02073 take $l=5$: the left side is
$P_{16}-P_{11}-P_6=3-2-1=0$. For C2-02080 take $l=7$:
it is $P_{18}-P_{11}-P_8=3-2-1=0$.
Both contradict the one-dimensional cohomology defect. Hence
neither basket has a canonical realization.


\paragraph{C2-02076: a moving torsion translate.}

\label{sec:c2076-twisted-quotient46}

The quotient argument also applies to the actual moving torsion
translate when the original complete canonical system has a fixed
divisor. Let $u\in H^0(pK_Z+\tau)$ be chosen from a subsystem
with no divisorial fixed part, and let $s\in H^0(lK_Z)$ be nonzero.
Choose $u$ not vanishing at any prime of $\operatorname{div}(s)$.
Then multiplication by $u$ induces an injection
\begin{equation}
 \frac{H^0(bK_Z-\tau)}{sH^0((b-l)K_Z-\tau)}
 \hookrightarrow
 \frac{H^0((b+p)K_Z)}{sH^0((b+p-l)K_Z)}.
 \label{eq:twisted-quotient-injection34}
\end{equation}
Indeed, if $ut$ is divisible by $s$, the absence of a zero of $u$
at every prime of $s$ implies that $t/s$ has no pole in
codimension one. Normality extends it to the indicated reflexive
section. This is a statement about the actual moving subsystem,
not about a hypothetical torsion class chosen independently.

For \texttt{C2-02076}, $I=420$, $V=1/420$, $\chi=2$, and
\[
 P_4=P_5=P_7=P_{11}=1,\quad P_{10}=2,\quad
 P_{12}=P_{16}=P_{17}=3,\quad P_{22}=6.
\]
Lemma~\ref{lem:first-degree-ten-common} leaves the regular profile
$(k,p_g,q,j,M)=(1,1,0,8,64)$ and the irregular profile
$(k,p_g,q,j)=(2,1,1,10)$. In the regular case
\[
 T_Z\sim8K_Z+\tau,\qquad E_Z\sim2K_Z-\tau,
\]
and the actual moving subsystem has rank two and no divisorial
fixed part.

The paired-torsion calculation gives $215$ necessary
inverse profiles. The simultaneous requirements
\[
 H_8(\tau)\ge2,\qquad H_2(-\tau)\ge1
\]
leave exactly one profile, with
\[
 H_4(-\tau)=0,\qquad H_8(-\tau)=1.
\]
Apply \eqref{eq:twisted-quotient-injection34} with $p=b=8$ and
$l=4$. Its source has dimension one, while its target has dimension
$P_{16}-P_{12}=0$. This excludes the regular profile.

In the irregular case, $j=10$ and the canonical fixed part is zero.
The ordinary quotient injection with anchor $10$,
$l=5$ and $b=7$ would require
\[
 P_{17}-P_{12}\ge P_7-P_2,
 \qquad 0\ge1,
\]
a contradiction. Both complete first-surface profiles are
excluded. Therefore C2-02076 has no canonical realization.


\paragraph{C2-02068: the opposite cohomology dimensions.}

\label{sec:moving-twisted-koszul-c2068-34}

\begin{theorem}
\label{thm:c2068-moving-koszul34}
The equality basket \texttt{C2-02068} has no
canonical realization.
\end{theorem}

\begin{proof}
Here $I=420$, $v=1$, $K^3=1/420$, $\chi=2$, and
\[
 P_2=P_7=0,\qquad
 P_4=P_5=P_6=P_9=P_{11}=P_{13}=1,
 \qquad P_{10}=2,
\]
\[
 P_{12}=P_{16}=3,\qquad P_{22}=6.
\]
As before, $p_g=q=0$ and $h^{2,0}=1$.

By Lemma~\ref{lem:first-degree-ten-common}, the first moving
surface has either $(k,p,q,j,M)=(1,1,0,8,64)$ or
$(k,p,q,j)=(2,1,1,10)$. In the second case $|10K_Z|$ has no
divisorial fixed part. Normal division with an anchor of degree ten
and $(a,b)=(4,6)$ would require
\[
 P_{16}-P_{12}\ge P_6-P_2,\qquad 0\ge1.
\]
We are therefore in the regular case, with the actual linked classes
\[
 T_Z\sim8K_Z+\tau,\qquad E_Z\sim2K_Z-\tau,
\]
and a fixed-free rank-two subsystem of $H^0(8K_Z+\tau)$.

\emph{The linked torsion profile.}
The inverse-torsion enumeration has $474$ necessary
paired profiles.  The actual moving and fixed divisors require
\[
 H_8(\tau)\ge2,\qquad H_2(-\tau)\ge1.
\]
These conditions leave five profiles.  Choose a moving
$u\in H^0(8K_Z+\tau)$ avoiding every prime of a nonzero
$s_4\in R_4$.  Twisted quotient injection requires
\[
 H_8(-\tau)-H_4(-\tau)\le P_{16}-P_{12}=0.
\]
Only one of the five states satisfies this.  Its needed dimensions are
\[
\begin{array}{c|rrrrrrrr}
 j&2&4&5&6&7&8&9&13\\\hline
 H_j(\tau)&1&0&0&1&1&2&0&2\\
 H_j(-\tau)&1&1&0&1&1&1&1&1
\end{array}
\]
and its compatible local torsion mass is exactly $C(\tau)=2$.
Moreover $H_1(\tau)=H_1(-\tau)=0$: either hypothetical
degree-one translate would multiply a nonzero opposite
degree-six section to $R_7=0$.

\emph{The actual double class.}
Put $J_j=H_j(2\tau)$ and $J_j^-=H_j(-2\tau)$.
The selected profile has nonzero degree-two and degree-seven
sections with twist $\tau$, and a nonzero degree-two section
with twist $-\tau$.  Therefore
\[
 J_9\ge1,\qquad J_4^-\ge1.
\]
The exhaustive list of $474$ profiles, applied to
the actual class $2\tau$, has the following finite implication:
\begin{equation}
 J_9\ge1,\ J_4^-\ge1\quad\Longrightarrow\quad J_{13}-J_9\le1.
 \label{eq:c2068-small-double-gap34}
\end{equation}
Only degree $2+7$ and the inverse degree-two square are used;
no high-degree cycle or additive virtual-character assumption
is involved.

\emph{The opposite cohomology dimensions.}
Apply Lemma~\ref{lem:twisted-koszul-defect34} to the actual
coprime pair $(u,s_4)$ of degrees $(8,4)$.
At $\sigma=0$ it gives
\[
 h^2(\mathcal O_Z(\tau))
 \le P_{13}-H_5(-\tau)-P_9+H_1(-\tau)
 =1-0-1+0=0.
\]
At $\sigma=2\tau$ it gives
\[
 h^2(\mathcal O_Z(-\tau))
 \le J_{13}-H_5(\tau)-J_9+H_1(\tau)
 =J_{13}-J_9\le1.
\]
The class $\tau$ is nontrivial because $H_2(\tau)=1>P_2=0$.
Thus the degree-zero cohomology of both torsion sheaves vanishes
by ampleness, and their degree-three cohomology vanishes by
Serre duality and the just-proved $H_1(\pm\tau)=0$.
The torsion RR identity now gives the contradiction
\[
 2=2\chi(\mathcal O_Z)-C(\tau)
 =\chi(\mathcal O_Z(\tau))+\chi(\mathcal O_Z(-\tau))
 \le h^2(\mathcal O_Z(\tau))+h^2(\mathcal O_Z(-\tau))\le1.
\]
Both actual first-surface profiles are excluded.
\end{proof}


\subsection{Coprime pairs and the two-twist defect}
\paragraph{C2-02079: producing two actual visible primes.}

\label{sec:twisted-koszul-c2079-34}
\begin{theorem}
\label{thm:c2079-twisted-koszul34}
The basket \texttt{C2-02079} has no canonical
realization.
\end{theorem}

\begin{proof}
Here $I=420$, $v=1$, $K^3=1/420$, $\chi=2$, and
\[
 B=2(1,2)+(1,3)+3(2,5)+(3,7)+(3,10)+(5,12).
\]
Its relevant canonical dimensions are
\[
 P_1=P_2=P_3=P_4=P_6=P_8=0,\qquad
 P_5=P_7=P_9=1,\qquad P_{10}=P_{15}=2.
\]
Thus $p_g=0$; tiny volume gives $q=0$, and consequently
$h^{2,0}=1$.

\emph{Step 1: actual doubling gives a visible-prime floor.}
The finite paired-torsion calculation has
eleven necessary profiles.  In either orientation it gives
\[
 H_3(\rho)=0,\qquad
 H_2(\rho)>0\Longrightarrow H_6(\rho)>0,\qquad
 H_4(\rho)>0\Longrightarrow H_8(\rho)=0.
\]
If $H_2(\tau)>0$, the square of a degree-two section and its
product with a degree-six section imply
$H_4(2\tau)>0$ and $H_8(2\tau)>0$, a contradiction.
Therefore $H_2(\tau)=0$ for every actual torsion class.  Its
square also excludes every effective $K_Z+\tau$.

For a visible prime, the actual-prime adjunction and Hodge bounds
are $-n\le M\le n^2$.  Through $n=4$, exact simultaneous RR
permits only
\[
 (1,1),(2,4),(3,9),(4,2),(4,16).
\]
The first three states are Hodge equalities and descend to
forbidden effective degree-one, two or three torsion translates.
Every visible prime thus has $n\ge4$.

\emph{Step 2: two actual visible primes.}
Choose $0\ne a\in R_5$.  Its full visible part is one prime $T$
of coefficient one, because $5<2\cdot4$.  Write
\[
 \operatorname{div}_Y(a)=T+E,\qquad\lambda(T)=5,
\]
where $E$ is the full invisible terminal remainder and
$\lambda(K)=1$.  If the entire canonical member were a reduced
prime, its complete restriction domain would have
\[
 G_6=1,\quad G_9=1,\quad G_{15}=0,
\]
contradicting multiplication.  Thus $n_T<5$, so $n_T=4$.
The Hodge-equality square $16$ would force $\lambda(T)=4$.
Consequently
\[
 (n_T,M_T)=(4,2).
\]

Take $0\ne e\in R_9$.  The exact plateau gives
$R_{15}=aR_{10}$.  At the visible prime $T$, the weighted bracket
has exact order
\[
 \operatorname{ord}_T(D_a e)=\operatorname{ord}_T(e)=t,
\]
because $5t-9\ne0$ for integral $t$.  The plateau forces
$t\ge1$, while $5t\le\lambda(e)=9$ gives $t=1$.
The remainder has scaled degree five and Hodge mass four.  It
therefore contains exactly one further visible prime, with
coefficient one:
\[
 \operatorname{div}_Y(e)=T+U+J,
 \qquad T\ne U,\quad\lambda(U)=4,\quad4\le n:=IK^2U\le5,
\]
where $J$ is fully invisible.

Put $M=IKU^2$ and $\ell=IKTU$.  Proper intersections of the
effective terminal remainders give
\[
 0\le\ell\le\min\{5n,34,9n-M\},
 \qquad(4n-\ell)^2\le14(n^2-M).
\]
Only on the equality boundary of the latter inequality is the
additional Hodge-mass equation
\[
 14(4-n)=4n-\ell
\]
imposed.  The coupled individual and mixed RR congruences, using
the same ordered local characters, leave exactly
\[
 (n,M,\ell)=(4,16,16).
\]
The simultaneous calculation considers all $1\,260\,000$ ordered
individual assignments.  There are $120$
first-prime assignments and $120$ surviving labelled pairs;
no strict Hodge-inequality pair survives.

Individual null descent now gives
\[
 U_Z\sim4K_Z+\tau.
\]
Let $u\in H^0(Z,4K_Z+\tau)$ have divisor exactly $U_Z$.
Its divisor is coprime to $\operatorname{div}_Z(a)=T_Z+E_Z$:
$U$ and $T$ are distinct visible primes, and $E$ is invisible.

\emph{Step 3: the two opposite cohomology dimensions.}
Among the eleven necessary paired profiles, the condition
$H_4(\tau)>0$ selects a unique orientation.  It satisfies
\[
 H_6(\tau)=0,\qquad H_6(-\tau)=1,
 \qquad H_5(\tau)=0,
\]
and has local torsion mass $C(\tau)=2$.  Every one of the eleven
profiles also satisfies
\begin{equation}
 H_{10}(\rho)-H_5(\rho)\le1.
 \label{eq:c2079-profile-gap34}
\end{equation}

Apply Lemma~\ref{lem:twisted-koszul-defect34} to the actual
coprime pair $(u,a)$ of degrees $(4,5)$.  With $\sigma=0$ it gives
\[
 h^2(\mathcal O_Z(\tau))
 \le P_{10}-H_6(-\tau)-P_5+H_1(-\tau)
 =2-1-1+0=0.
\]
With $\sigma=2\tau$, it gives
\[
 h^2(\mathcal O_Z(-\tau))
 \le H_{10}(2\tau)-H_6(\tau)-H_5(2\tau)+H_1(\tau)
 \le1,
\]
by \eqref{eq:c2079-profile-gap34} applied to the actual class
$2\tau$.

The torsion class is nontrivial, since $H_4(\tau)>0$ but $P_4=0$.
Ample canonical degree therefore gives $H^0(\mathcal O_Z(\pm\tau))=0$;
Serre duality and the absence of every effective $K_Z\pm\tau$
give $H^3(\mathcal O_Z(\pm\tau))=0$.  The torsion RR
identity yields
\[
 \chi(\mathcal O_Z(\tau))+\chi(\mathcal O_Z(-\tau))
 =2\chi(\mathcal O_Z)-C(\tau)=2.
\]
On the other hand, the two cohomological bounds give
\[
 \chi(\mathcal O_Z(\tau))+\chi(\mathcal O_Z(-\tau))
 \le h^2(\mathcal O_Z(\tau))+h^2(\mathcal O_Z(-\tau))\le1,
\]
a contradiction.
\end{proof}


\paragraph{C2-02069: opposite components of a quartic.}

\label{sec:two-twist-koszul-c2069-34}
\begin{theorem}
\label{thm:c2069-two-twist34}
The basket \texttt{C2-02069} has
no canonical realization.
\end{theorem}

\begin{proof}
Suppose $Z$ is a realization and take a terminal
$\mathbb Q$-factorial minimal model $Y$.  Here
$I=420$, $v=1$, $K^3=1/420$, $\chi=2$, and $p_g=0$.
The tiny-volume irregularity bound gives $q=0$, hence
$h^{2,0}=1$.  The required canonical dimensions are
\[
 P_1=P_2=P_3=P_6=0,\qquad
 P_4=P_5=P_{13}=1,\qquad P_{10}=2.
\]
The paired-torsion enumeration has
$28$ necessary profiles.  All satisfy
$H_{12}(\rho)\ge2$ in both orientations.  Thus
$H_1(\rho)=0$ for every actual torsion class: multiplication
by a nonzero section of $K_Z+\rho$ would inject at least
two sections of $12K_Z-\rho$ into $R_{13}$, whose dimension
is one.

For visible primes of scaled degree at most two, the
actual-prime adjunction bound, Hodge bound, and exact
simultaneous Weil--RR congruences allow only
\[
 (n,M)=(1,1),(2,4),\qquad
 n=420K_Y^2T,\quad M=420K_YT^2.
\]
The first is a Hodge equality whose null descent would
give an effective $K_Z+\rho$, already excluded.  Hence
every visible prime has scaled degree at least two.

Choose $0\ne a_4\in R_4$.  Its visible part is nonempty
and has total degree at most four.  There are only three
possibilities: one prime of coefficient one, one prime
of coefficient two, or two distinct primes each of
coefficient one.  The first is impossible by
Lemma~\ref{lem:single-visible-square34}: $P_5=1$ would
force $P_{10}\ge3$, whereas $P_{10}=2$.

The second possibility is also impossible.  The normalized
adjunction-square construction of \cite[Lemma~9.1]{LBVolume}
says that if every visible
component of a degree-$m$ member is repeated, division
of its normalized square by that member gives a nonzero
actual section of $(m+2)K_Z$.  Its exact visible orders
are $k-2$ and its invisible orders are nonnegative; thus
the assertion includes the entire invisible remainder.
Here $m=4$ and $P_6=0$, a contradiction.

Consequently the divisor of $a_4$ on the canonical model
is
\[
 D_4=T+U
\]
for two distinct reduced visible primes, each of scaled
degree two.  Their degrees already exhaust the degree of
$D_4$, so ampleness leaves no additional canonical prime.
Both have Hodge-equality square four, and actual null
descent gives
\[
 T\sim2K_Z+\tau,\qquad U\sim2K_Z-\tau.
\]
The opposite classes follow from the actual divisor
identity $T+U\sim4K_Z$, not from addition of virtual
basket characters.  Their defining sections are coprime.

Among the complete $28$ necessary paired profiles,
exactly one has both $H_2(\tau)>0$ and $H_2(-\tau)>0$.
It has
\[
 H_5(\tau)=H_5(-\tau)=0,\qquad C(\tau)=2,
\]
where $C$ is the local torsion RR mass.
The class is nontrivial because $P_2=0$.

Apply Lemma~\ref{lem:two-twist-koszul34} with
$a=b=2$, $\lambda=\tau$, and $\mu=-\tau$.
Taking $\sigma=\tau$ yields
\[
 h^2(\mathcal O_Z(-\tau))
 \le H_5(\tau)-P_3-H_3(2\tau)+H_1(\tau)
 =-H_3(2\tau)\le0.
\]
Taking $\sigma=-\tau$ similarly yields
$h^2(\mathcal O_Z(\tau))=0$.
Ample canonical degree and the nontriviality of $\tau$
give $H^0(\mathcal O_Z(\pm\tau))=0$; the previously
proved $H_1(\pm\tau)=0$ and Serre duality give
$H^3(\mathcal O_Z(\pm\tau))=0$.  Consequently
\[
 \chi(\mathcal O_Z(\tau))+
 \chi(\mathcal O_Z(-\tau))\le0.
\]
But the actual opposite torsion RR identity gives
\[
 \chi(\mathcal O_Z(\tau))+
 \chi(\mathcal O_Z(-\tau))
 =2\chi(\mathcal O_Z)-C(\tau)=2,
\]
a contradiction.
\end{proof}


The preceding arguments exclude seventeen distinct rows of
Table~\ref{tab:remaining28}. Each uses either a complete member,+a geometrically produced coprime pair, or a linked moving and fixed
decomposition. Together with Proposition~\ref{prop:initial-cuts46},
they reduce the original sixty-five Euler-characteristic-two baskets
to the eleven possibilities displayed below.


\subsection{The eleven remaining baskets}
\label{sec:eleven-baskets}
The preceding applications give the following endpoint of the common
reduction. The table records full terminalization baskets, with
$I=\operatorname{lcm}(r)$, $e=\min\{m>0:P_m>0\}$ and
$\delta=\min\{m:P_m\ge2\}$. Its last column points to the later
argument; at this stage every row is only a necessary possibility.

\begin{proposition}[The eleven-basket reduction]\label{prop:eleven-baskets}
Every canonical model of volume $1/420$ belongs either to the
Euler-one hypersurface family of Section~\ref{sec:euler-one}, or has
$\chi=2$ and one of the eleven terminalization baskets in
Table~\ref{tab:residual11}.
\end{proposition}
\begin{proof}
The complete census gives $1+65+283$ possibilities. The Euler-three
branch has been excluded in Section~\ref{sec:chi3}, and the Euler-one
branch has been identified in Section~\ref{sec:euler-one}. In the
Euler-two branch, the first uniform tests remove $37$ baskets and the
applications above remove $17$ more. Subtracting their distinct row
identifiers from the complete list gives exactly the eleven rows below.
The intermediate lists and all finite witnesses are retained in
Appendix~\ref{app:census}.
\end{proof}

\begingroup\small
\setlength{\tabcolsep}{3pt}
\begin{longtable}{@{}p{19mm}p{88mm}rccp{27mm}@{}}
\caption{The eleven Euler-two baskets after the common reduction.}\label{tab:residual11}\\
\toprule
Basket & Full multiset & $I$ & $e$ & $\delta$ & Subsequent result\\
\midrule\endfirsthead
\toprule
Basket & Full multiset & $I$ & $e$ & $\delta$ & Subsequent result\\
\midrule\endhead
C2-01825 & $2(1,2)+2(1,3)+(1,4)+(2,5)+(3,7)+(4,13)+(6,13)$ & 5460 & 6 & 9 & Excluded, Thm.~\ref{thm:c1825-same-pencil-extension}\\[3pt]
C2-02075 & $5(1,2)+(1,6)+2(2,5)+(3,7)+(7,20)$ & 420 & 5 & 10 & Excluded, Thm.~\ref{thm:c2075-excluded46}\\[3pt]
C2-02347 & $7(1,2)+2(1,3)+3(1,4)+(1,5)+2(2,5)+(3,7)$ & 420 & 4 & 12 & Vega, Thm.~\ref{thm:quartic-recognition46}\\[3pt]
C2-02349 & $5(1,2)+(1,3)+4(1,4)+(2,5)+(3,7)+(5,12)$ & 420 & 4 & 12 & Excluded, Thm.~\ref{thm:2349-exclusion46}\\[3pt]
C2-02350 & $5(1,2)+(1,3)+2(1,4)+(2,5)+(5,12)+(5,14)$ & 420 & 8 & 12 & Excluded, Thm.~\ref{thm:octic-complete-exclusion46}\\[3pt]
C2-02352 & $3(1,2)+8(1,3)+3(1,4)+(2,5)+(3,7)$ & 420 & 3 & 12 & Excluded, Thm.~\ref{thm:c2352-elliptic-root-extension}\\[3pt]
C2-02353 & $3(1,2)+4(1,3)+(1,5)+2(2,5)+(3,7)+(5,12)$ & 420 & 5 & 12 & Excluded, Thm.~\ref{thm:2353-exclusion46}\\[3pt]
C2-02354 & $2(1,2)+(1,3)+(1,4)+(2,5)+2(3,7)+(8,21)$ & 420 & 5 & 12 & Excluded, Thm.~\ref{thm:c2354-coupled34}\\[3pt]
C2-02355 & $2(1,2)+(1,3)+(1,4)+(2,5)+(3,7)+2(3,8)+(5,12)$ & 840 & 5 & 12 & Excluded, Thm.~\ref{thm:2355-exclusion46}\\[3pt]
C2-02357 & $(1,2)+2(1,3)+(1,4)+(2,5)+(3,7)+(4,9)+(7,18)$ & 1260 & 5 & 12 & Excluded, Thm.~\ref{thm:2357-exclusion46}\\[3pt]
C2-02397 & $4(1,2)+2(1,3)+(1,5)+2(2,5)+(3,7)+2(3,8)$ & 840 & 5 & 14 & Altair, Thm.~\ref{thm:quintic-whole-basket46}\\[3pt]
\bottomrule
\end{longtable}\endgroup

Parts~III and~IV now treat precisely this residual list. The eight
non-octic exclusions are grouped by their geometric mechanism; the
ninth, \basket{C2-02350}, has its own conductor analysis. The two
remaining baskets will be realized and recognized as Vega and Altair.


\clearpage
\part{The difficult basket exclusions}

The nine remaining exclusions use the complete restriction spaces and
the polar operations of Part~II. We first determine the actual visible
components and source points, then compare their section dimensions,
local characters and polar degrees. For the nonnormal octic this
requires the full conductor algebra. Each case below records the
additional geometric obstruction after these common constraints.

\section{Three special geometric obstructions}
\label{sec:special-obstructions}

The next three baskets require different geometric arguments. For
C2-01825, two adjacent pencils are forced to coincide, and a crepant
graph separates the last configuration. For C2-02075 we use the
geography of the complete moving surface. For C2-02352 a differential
identity produces an elliptic canonical-root cover of excessive volume.

\subsection{Exclusion of C2-01825: adjacent pencils and a crepant graph}
\paragraph{The small K3 polarizations and the visible floor.}

\label{sec:c1825-k3-floor-extension}

We first raise the minimum degree of a visible prime. This will force
the adjacent ninth and tenth systems to define the same pencil.

\begin{lemma}[Two forbidden Du Val K3 polarizations]
There is no integral ample Weil divisor $H$ on a Du Val K3 surface
with $H^2=1/780$ or $H^2=1/364$.
\end{lemma}

\begin{proof}
The exceptional ADE lattice on the minimal resolution has rank at most
nineteen. Surface Riemann--Roch gives
$\sum_Pq_P(H)-H^2\in2\mathbb Z$. Consequently the reduced denominator
of $H^2$ divides the least common multiple of the local correction
denominators. For an $A_r$ point the latter divides $r+1$;
$D_r$ points contribute only powers of two of exponent at most two,
and the exceptional types contribute only two or three.

In either proposed square, the factor thirteen requires $A_{12}$:
the next such type, $A_{25}$, already exceeds the rank budget.
This leaves at most seven ranks. For $1/780$, the remaining
denominator $4\cdot3\cdot5$ needs at least
$3+2+4=9$ ranks. Under a seven-rank budget the factor five can only
come from $A_4$; the remaining three ranks cannot supply both four
and three. For $1/364$, the factor seven requires $A_6$ after
$A_{12}$, leaving only one rank, which cannot supply the factor
four. These are contradictions. Extra Cartier factors and
cancellation of local corrections cannot create a missing
prime-power denominator.
\end{proof}

\begin{proposition}
For a hypothetical realization of \texttt{C2-01825}, every visible
prime has scaled degree $IK^2T\ge18$, where $I=5460$ and $IK^3=13$.
\end{proposition}

\begin{proof}
The exact torsion census has five paired profiles.
Canonical-times-twisted multiplication removes two: a section in
degree four multiplied by the two-dimensional degree-nine canonical
space would have to inject into a one-dimensional twisted degree
thirteen space. Of the three remaining profiles, the two
noncanonical possibilities each have, in one common orientation,
degree-two and degree-four spaces whose product has dimension at least two in
degree six of an actual doubled torsion class. Every remaining
profile has degree-six dimension at most one. Thus every actual
torsion class has the canonical section-dimension profile in the
verified degrees. In particular $H_2(\tau)=0$ for every actual
class, and $H_1(\tau)=0$ follows by squaring.

The exact simultaneous Weil--RR congruences, visible adjunction
$M\ge-n$, and Hodge inequality $13M\le n^2$ allow precisely
\[
 (n,M)=(7,-7),\ (13,13),\ (15,-15)
 \qquad (1\le n\le17).
\]
The equality state $(13,13)$ descends to $K+\tau$ and is excluded
by $H_1(\tau)=0$. Saturated adjunction for the other two states
produces a Du Val K3 normalization with an integral ample Weil
class of square $7/5460=1/780$ or $15/5460=1/364$.
The preceding lemma excludes both. The next necessary state is
$(18,-2)$, so the claimed lower bound follows.
\end{proof}

For the first connected $9K$ pencil, fractional adjunction and
the exact RR lattice give $j=IK^2F\ge102$. Hence its fixed
divisor has total scaled degree at most $117-102=15$, strictly
below the floor eighteen. Its entire fixed support is therefore
invisible, even if it is nonreduced. The full fixed-divisor
trace-rank argument of \cite[Theorem~15.2]{LBVolume} applies at equality
as well. Here is the required specialization, retaining the full fixed
scheme. Let $s\in R_9$ define the chosen general member, put
$p=p_g(F)$, and let $r_{10}$ be the rational residue rank of
$R_{10}$ on the general fibre. The full-member trace theorem
\cite[Theorem~8.2]{LBVolume} injects $H^0(F,2K_F)$ into
$R_{20}/sR_{11}$; denote its image by $W_2$, so that
$\dim W_2=P_2(F)$. The homogeneous square map
\[
 R_{10}\longrightarrow (R_{20}/sR_{11})/W_2
\]
has zero locus of dimension at least
$P_{10}-(U_{2,9}-P_2(F))$. A zero has regular canonical residue
on $F$, since its square is regular. The difference between its
ambient section and its component trace is supported on the full
fixed divisor; coherent invisibility makes its connecting class
zero. Its residue therefore lies in a space of dimension $p-1$,
the kernel of the moving trace, since the ambient restriction rank
is one. The inverse image has dimension at most
$P_{10}-r_{10}+p-1$. Comparing dimensions gives
\[
 U_{2,9}-P_2(F)\ge r_{10}-p+1.
\]
This argument has no strict-volume hypothesis. For either remaining
profile $(K_S^2,p_g,q)=(1,1,0)$ or $(1,2,0)$, surface
Riemann--Roch gives $P_2(F)=p+2$. Since
$U_{2,9}=P_{20}-P_{11}=4$, the inequality gives $r_{10}\le1$.
The rank is nonzero on a general moving fibre, so the $10K$ map
factors through the same connected pencil. Its moving multiplicity
is one, since $2j>10\cdot13$.

This argument does not replace the actual residue rank by $P_{10}$.
It also does not charge twice the visible-prime floor to a repeated
fixed prime without proving that prime is visible. The resulting
same-pencil configuration is for the next geometric step.


\paragraph{The adjacent-pencil contradiction.}

\label{sec:c1825-same-pencil-extension}

The sharper visible-prime floor in
Section~\ref{sec:c1825-k3-floor-extension} supplies the initial
reduction.  The final step requires a crepant change of model:
zero canonical intersection does not by itself make Weil divisors
disjoint on a non-$\mathbb Q$-factorial canonical model.

\begin{lemma}[The actual torsion group for C2-01825 is trivial]
\label{lem:c1825-torsion-free-extension}
Every codimension-two-Cartier torsion Weil class on a hypothetical
realization of \texttt{C2-01825} is trivial.
\end{lemma}

\begin{proof}
The five necessary profiles, followed
by the actual multiplication filters used in the floor-eighteen
proof, leave only the canonical section-dimension profile.
Its only compatible correction mass is zero.  In particular every
actual torsion class $\sigma$ has $C(\sigma)=0$ and $H_2(\sigma)=0$.

If $\tau$ had exact order $N\ge2$, its connected canonical
quasi-\'{e}tale cyclic cover would have
\[
 V'=\frac{N}{420},\qquad \chi(\mathcal O_{Z'})=2N,
 \qquad P_2(Z')=0,
 \qquad h^{2,0}(Z')\ge2N-1.
\]
The Euler identity follows from the actual character-cover
Riemann--Roch formula with every correction mass zero;
the plurigenus identity is the eigensheaf decomposition.
No additive law on virtual basket characters is used.

For $N\ge3$, the linear consequence of Chen--Jiang,
Corollary~2.3 in arXiv:2211.00926v3
\cite{ChenJiangForms}, gives
\[
 \frac{N}{420}\ge\frac{h^{2,0}(Z')-3}{60}
 \ge\frac{2N-4}{60},\qquad N\le\frac{28}{13}<3,
\]
a contradiction.  The cited linear consequence holds for
$h^{2,0}\ge4$, which includes this case.

For $N=2$, take a smooth projective model of the cover.
It has at least three independent holomorphic two-forms.
Their evaluation in $\Omega^2$ has rank at most two:
rank three would give a nonzero determinant section of
$\det\Omega^2\simeq2K$, contrary to $P_2=0$.
Chen--Jiang, Lemma~2.2, therefore supplies a subsheaf of
$\Omega^2$ whose determinant has at least two independent
sections.  Their Lemma~2.7 gives
\[
 V'\ge\frac1{14}>\frac1{210}=V',
\]
the final contradiction.
\end{proof}

\begin{theorem}
\label{thm:c1825-same-pencil-extension}
The basket \texttt{C2-01825} has no
canonical realization.
\end{theorem}

\begin{proof}
Assume that $Z$ is a realization and let $h:Y\to Z$ be a
projective crepant terminal $\mathbb Q$-factorial model.
Write $K=K_Y$, $I=5460$ and $v=IK^3=13$.
The first-$9K$ and adjacent-$10K$ conclusions of
Section~\ref{sec:c1825-k3-floor-extension} give the same
connected pencil, with moving multiplicity one in both
systems.  Let $F$ denote a general moving prime on $Y$.
The two actual complete-system decompositions are
\[
 9K\sim F+E_9,\qquad 10K\sim F+E_{10}.
\]
All prime components of $E_9$, including exceptional
ones, are invisible.  Normalize the ambient two-form
functional by $\lambda(K)=1$.  Then
\[
 \lambda(F)=9,\qquad \lambda(E_{10})=1.
\]
Put $j=IK^2F$ and $M=IKF^2$.  The exact first-pencil
constraints give $102\le j\le117$.  In particular the
scaled degree of $E_{10}$ is at most $130-102=28$.
The proved visible-prime floor is eighteen, so $E_{10}$
contains exactly one visible prime $T$, with coefficient one.
Write
\[
 E_{10}=T+J,
\]
where $J$ is the full effective invisible remainder and has
no component $T$.  Thus $\lambda(T)=1$.
Set $n=IK^2T$, $N=IKT^2$ and $\ell=IKFT$.

The actual effective divisors on the terminal
$\mathbb Q$-factorial model give
\begin{align*}
 &18\le n\le130-j,\qquad -n\le N,\qquad 13N\le n^2,\\
 &0\le\ell\le\min\{9n,\,10n-N,\,10j-M\}.
\end{align*}
For example the three upper bounds are respectively
$KE_9T\ge0$, $KJT\ge0$ and $KJF\ge0$.
These positivity statements are made on $Y$, where effective
divisors with no common prime component have effective rational
intersection cycles.

The mixed Hodge inequality is
\[
 (jn-13\ell)^2\le(j^2-13M)(n^2-13N).
\]
If equality holds and $j^2-13M>0$, null-divisor descent and
the actual Hodge masses impose
\[
 (13-n)(j^2-13M)=(117-j)(jn-13\ell).
\]
An individual Hodge equality for $F$ instead forces
$j=117$ from $\lambda(F)=9$; an individual equality
for $T$ forces $n=13$ and is unavailable here.
We impose these conditions together with the simultaneous
local degree, square and mixed Weil--Riemann--Roch congruences.
All $851760$ ordered local assignments are enumerated in the
finite congruence calculation. There remain exactly three numerical types:
\[
\begin{array}{c|rrrrr|r}
 &j&M&n&N&\ell&\text{paired assignments}\\\hline
 \mathrm A&102&768&28&28&252&1\\
 \mathrm B&110&920&20&20&180&2\\
 \mathrm C&103&773&20&20&180&2.
\end{array}
\]
All three have zero mixed Hodge determinant.  Strict Hodge
states are throughout the search and none survives;
the Hodge-mass identity is used only on the equality boundary.

In types A and B, $j+n=130$, so $J_Z=0$ by ampleness of
$K_Z$.  A general complete $10K$ member on the canonical
model is consequently $F_Z+T_Z$, a reduced divisor.
The general moving surface has minimal invariants
$(K_S^2,p_g,q)=(1,1,0)$ or $(1,2,0)$.
The reduced-member trace estimate \eqref{eq:profile-defect46} gives
\[
 P_{22}-P_{12}\ge K_S^2+1+P_{11}+h^{2,0}(Z)=5.
\]
The exact row gives $P_{22}-P_{12}=7-3=4$, a contradiction.

It remains to exclude type C.  The integral divisor
\[
 G=F+2T-11K=T-K-J
\]
has $K^2G=KG^2=0$; these equalities follow directly from
the displayed type C numbers.  Null-divisor descent makes
$G_Z$ a codimension-two-Cartier torsion class, which is
trivial by Lemma~\ref{lem:c1825-torsion-free-extension}.
After subtracting a principal divisor, $G$ is therefore
$h$-exceptional.  Its components centred on curves vanish:
restrict to a general member of a sufficiently divisible
canonical system avoiding the finitely many point centres.
The exceptional curve intersection matrix on that surface
is negative definite, whereas the restriction of $G$ has
square zero because $KG^2=0$.  Thus
\[
 T-J-K\sim P
\]
for an exceptional divisor $P$ supported over points,
with coefficients of either sign allowed.  Moreover
\[
 IK^2T=20,\qquad IK^2J=7,\qquad IKTJ=200-180-20=0.
\]

Apply the crepant graph separation lemma,
Lemma~\ref{lem:crepant-graph-binary34}, to this full
terminal pair $T,J$.  It gives another crepant terminal
minimal model carrying an integral divisor $D$ with
\[
 IK^2D=20,
\]
whose local class at every virtual basket constituent is
either zero or the canonical class one.  The canonical
Reid basket is unchanged under this crepant change of
terminal minimal model.

Reduce its linear Weil--Riemann--Roch congruence modulo
thirteen.  All constituents except $(4,13)$ and $(6,13)$
contribute zero.  For those two constituents the weights are
\[
 \frac{5460}{13}\,4^2\equiv-1\pmod{13},\qquad
 \frac{5460}{13}\,6^2\equiv1\pmod{13}.
\]
The binary local classes therefore allow only
$IK^2D\equiv0,1,-1\pmod{13}$.
But $IK^2D=20\equiv7\pmod{13}$.
This contradiction excludes type C and completes the proof.
\end{proof}


\subsection{Exclusion of C2-02075: the complete moving surface}

\label{sec:c2075-complete46}

This exclusion uses information about a complete moving system that
is absent from numerical surface geography alone.

\begin{lemma}\label{lem:two-forbidden-net-profiles46}
Let $S$ be a smooth minimal surface of general type. Suppose there
is a nonzero integral nef divisor $D$ with $h^0(S,D)\ge2$ and
$K_S-D$ big. Then $(K_S^2,p_g,q)$ is neither $(3,1,1)$ nor
$(4,2,2)$.
\end{lemma}
\begin{proof}
First let $(K_S^2,p_g,q)=(3,1,1)$. Since $K_S-D$ is big,
$0<K_SD<3$. Hodge index and adjunction parity exclude $K_SD=1$;
when $K_SD=2$ they give $D^2=0$. The moving divisor therefore
defines a genus-two pencil, as follows. Write $D=M+E$ with $M$
the moving part and $E$ fixed. Nefness and $D^2=0$ give
$DM=DE=0$; since $M^2,ME\ge0$, both vanish. In particular
$|M|$ has no isolated base point and $M=aF$ for a fibre class
$F$. The positive even integer $K_SF$ is at least two, whereas
$K_SD=2$. Thus $a=1$, $K_SF=2$ and $K_SE=0$. The
$K_S$-null curves have negative-definite intersection matrix;
as $E^2=DE-ME=0$, this forces $E=0$.
Its base is $\mathbb P^1$, because a degree-one line bundle with
two sections has rational base. Polizzi's
Propositions~6.1 and~6.18 \cite{Polizzi} apply to precisely these
invariants and this rational genus-two pencil: the Albanese fibre
$F$ has genus three and $FD=4$. Thus
\[
 (K_S-D)F=(2\cdot3-2)-4=0,
\]
contrary to the positivity of a big divisor against a nonzero nef
fibre class. This application does not assume the Albanese genus
in advance.

Now let $(K_S^2,p_g,q)=(4,2,2)$. If the Albanese image is a curve,
its connected base has genus two. A relatively minimal fibration
with general fibre genus $g\ge2$ has
$K_S^2\ge8(g-1)(2-1)\ge8$, a contradiction. Thus $S$ has maximal
Albanese dimension. The equality theorem of
Barja--Pardini--Stoppino \cite[Theorem~1.1]{BPS} identifies its
canonical model as a double cover of its Albanese surface $A$.
For the Albanese morphism $f:S\to A$ one has
\[
 K_S=f^*L,\qquad L\text{ ample and integral},\qquad L^2=2.
\]
The theorem includes the permitted canonical branch degenerations.

Put $P=f_*D$, $Q=f_*(K_S-D)=2L-P$. The norm of rational functions
preserves linear equivalence under pushforward. The map
$|D|\to|P|$ has finite fibres: over a prescribed effective divisor
the nonexceptional prime components are fixed and have bounded
integral multiplicities; the finitely many exceptional primes also
have bounded multiplicities by intersection with a fixed ample
divisor on $S$. Hence $h^0(A,P)\ge2$. The same argument applied
to high multiples shows that pushforward preserves bigness, so
$Q$ is big. On an abelian surface an integral big class is ample;
therefore $Q^2\ge2$. The effective divisor $P$ is nef.

Hodge index and adjunction parity give $K_SD\ge2$ as above.
Consequently
\[
 LQ=4-K_SD\le2,\qquad
 (LQ)^2\ge L^2Q^2\ge4.
\]
Equality holds throughout. Hodge index gives $Q\equiv L$, hence
$P\equiv L$. Abelian Riemann--Roch and vanishing give
$h^0(A,P)=P^2/2=1$, the required contradiction.
\end{proof}

\begin{theorem}\label{thm:c2075-excluded46}
The equality basket C2-02075 has no canonical realization.
\end{theorem}
\begin{proof}
Its basket is
\[
 5(1,2)+(1,6)+2(2,5)+(3,7)+(7,20),\qquad I=420,
\]
and Riemann--Roch gives
\[
 P_{12}=3,\quad P_{13}=2,\quad P_{14}=4,\quad P_{26}=10.
\]
If the complete $12K$-system has curve image, its moving degree is
at least two and \eqref{eq:pencil-area46} gives
$1\le k\le6\cdot7^2/420=7/10$, impossible.
For a general moving surface in a complete nonpencil,
Lemma~\ref{lem:complete-moving-square46} and the area bound give
\[
 3\le k\le12\cdot13^2/420=169/35<5.
\]
Restriction of the ambient two-form is injective: a kernel form
would give $k\le12(1+24)/420=5/7$. Thus the connecting rank in
\eqref{eq:profile-defect46} is one, including for irregular surfaces.
The bicanonical trace is $U_{2,12}=10-4=6$. Surface geography,
$p\ge1$, and $k+1+p-q\le6$ leave exactly
\[
\begin{array}{c|c|c}
k&(p,q)&k+1+p-q+\max\{0,3-p\}\\\hline
3&(1,0),(2,0),(1,1)&7,7,6\\
4&(1,0),(1,1),(2,1),(2,2)&8,7,7,6.
\end{array}
\]
Every entry greater than six forces a repeated fixed prime.
Its cost $2/I$ changes the area cap to
$10\cdot11^2/420=121/42<3$, excluding that entry.
Only $(k,p,q)=(3,1,1)$ and $(4,2,2)$ remain.
The complete moving surface supplies the integral nef divisor
$D=\sigma_*(M|_S)$ of the proof of
Lemma~\ref{lem:complete-moving-square46}: it has at least two
sections and $K_{S_0}-D$ big. Both remaining profiles contradict
Lemma~\ref{lem:two-forbidden-net-profiles46}.
\end{proof}


\subsection{Exclusion of C2-02352: an elliptic canonical-root cover}

\label{sec:c2352-elliptic-root-extension}

We exclude the basket \texttt{C2-02352} using a
cubic covariant of Corollary~\ref{cor:polar-cubic-covariant} and
the irregular volume bound on a branched canonical-root cover.  No
normality assumption on the prime canonical member, and no assertion that
the root cover has canonical singularities, is made.

\begin{lemma}[Canonical-root volume bound]
\label{lem:canonical-root-volume-extension}
Let $Z$ be a normal projective canonical threefold with ample $K_Z$.
Suppose $s\in H^0(Z,eK_Z)$ has a single reduced prime divisor $T$.
The normalized cover $\pi:W\to Z$ obtained by adjoining an $e$-th
canonical root of $s$ is connected of degree $e$, and
\[
 K_W\sim_{\mathbb Q}\pi^*\left(K_Z+\frac{e-1}{e}T\right)
 \sim_{\mathbb Q}e\pi^*K_Z.
\]
Every smooth projective resolution $\widetilde W$ is of general type and
\[
 \operatorname{vol}(\widetilde W)\le e^4K_Z^3.
\]
In particular, if this cover admits a nonconstant rational map to an elliptic
curve, then $e^4K_Z^3\ge3/8$.
\end{lemma}

\begin{proof}
Choose a rational canonical form $\eta$ and write $s=f\eta^e$.
Adjoin $w$ with $w^e=f$.  The valuation of $f$ at $T$ is congruent to one
modulo $e$, so the field extension has degree $e$.  Normalization in that
field is finite and connected.  In codimension one the only branch divisor
is $T$, with full ramification index $e$; changing the rational canonical
frame changes valuations by multiples of $e$.  The displayed canonical
divisor formula follows from the codimension-one ramification formula,
and both sides extend as rational Cartier divisor classes by normality.
Thus $K_W$ is rational Cartier and ample.

The resolution is of general type because it is generically finite over
the general-type variety $Z$: resolve the rational map to a smooth
resolution of $Z$ and pull back pluricanonical forms.  For sufficiently
divisible $m$, restriction over the smooth locus of the normal $W$ gives
an injection
\[
 H^0(\widetilde W,mK_{\widetilde W})
 \hookrightarrow H^0(W,\mathcal O_W(mK_W)).
\]
Indeed a pluricanonical form upstairs is regular at the generic point of
every prime downstairs, and normality extends its reflexive section.
This does not require canonical singularities on $W$.  Taking leading
asymptotic dimensions gives
$\operatorname{vol}(\widetilde W)\le K_W^3=e^4K_Z^3$.
If there is a nonconstant rational map to an elliptic curve, resolution of
that map shows $q(\widetilde W)>0$.  Jiang's Theorem~1.5 in
\emph{On Severi type inequalities}, arXiv:1901.11207v2
\cite{JiangSeveri}, gives
$\operatorname{vol}(\widetilde W)\ge3/8$.
\end{proof}

\begin{proposition}[Prime cubic differential obstruction]
\label{prop:prime-cubic-elliptic-extension}
Let $Z$ be as above, with $p_g=q=0$ and $h^{2,0}=1$.
Suppose $0\ne a\in R_3$ has a reduced prime divisor $T$ visible
to the ambient two-form, and $P_4=P_{10}=0$. Then
$81K_Z^3\ge3/8$.
\end{proposition}

\begin{proof}
Use the normalized square and cubic covariant of
Lemma~\ref{lem:polar-square46} and
Corollary~\ref{cor:polar-cubic-covariant}. They give
$u\in R_8$, $v\in R_{12}$ and $w\in R_{18}$ with
\[
 D_a(v)=-12u^2,\qquad v^2-u^3=a^2w,\qquad D_a(w)=0,
\]
since the curvature section lies in $R_{10}=0$.
Visibility and reducedness give $\operatorname{ord}_T(u)=0$,
so $u\ne0$ and the first identity also gives $v\ne0$.
If $w=0$, then $(v/u)^2=u$; normality extends $v/u$ to
a nonzero section of $R_4$, contrary to $P_4=0$.

Put $k=\operatorname{ord}_T(w)$. The leading coefficient of
$D_a(w)$ at the visible prime is a nonzero unit times $3k-18$.
Thus $k=6$. Since $\operatorname{div}(a)=T$, the degree-zero
ratio $w/a^6$ has no divisorial pole anywhere on $Z$.
Normality and projectivity make it a nonzero constant $c$, and
\[
                         v^2=u^3+ca^8.
\]
On the cubic canonical-root cover $\pi:W\to Z$, write
$\sigma^3=\pi^*a$. The rational functions
\[
 x=\frac{\pi^*u}{\sigma^8},\qquad
 y=\frac{\pi^*v}{\sigma^{12}}
\]
satisfy $y^2=x^3+c$. Its smooth projective model is elliptic,
and the map is nonconstant: $u$ is a unit along $T$, so $x$ has
a pole of order eight above the fully ramified divisor $T$.
Lemma~\ref{lem:canonical-root-volume-extension}, with $e=3$,
now gives the assertion.
\end{proof}

\begin{theorem}
\label{thm:c2352-elliptic-root-extension}
The basket \texttt{C2-02352} has no canonical
realization.
\end{theorem}

\begin{proof}
Its exact data give $I=420$, $K^3=1/420$, $\chi=2$,
and
\[
 P_3=P_8=1,\qquad P_4=P_5=P_{10}=P_{13}=0.
\]
Jiang's irregular volume bound applied to a smooth model first gives
$q=0$, hence $h^{2,0}=1$.

For completeness, the forced prime cubic can be recovered from the
exact arithmetic.  All thirteen admissible
paired torsion profiles have $H_{12}(\tau)\ge1$ in both orientations.
Thus $H_1(\tau)=0$ for every actual codimension-two-Cartier torsion
class, since an effective $K+\tau$ would multiply a nonzero
$12K-\tau$ section into $R_{13}=0$.
The exact Weil--RR, visible adjunction and Hodge constraints for a
visible prime of scaled degree $n\le3$ allow only
\[
 (n,M)=(1,1),(2,4),(3,9).
\]
Every such prime is a Hodge equality and descends to
$nK_Z+\tau$ by the null-divisor lemma.  The first case is
excluded by $H_1(\tau)=0$.  Apply visible saturation to an actual
degree-three member: its visible components exhaust its full
canonical degree, with integer component degrees at least two.
The only partition of three is therefore a single component of
degree three and multiplicity one.  The invisible complement on
$Z$ vanishes by ampleness.  Thus $\operatorname{div}(s_3)=T$ is the
required visible reduced prime.

Proposition~\ref{prop:prime-cubic-elliptic-extension} would imply
\[
 \frac38\le81K_Z^3=\frac{81}{420}=\frac{27}{140}<\frac38,
\]
a contradiction.
\end{proof}

The companion exact calculation verifies the thirteen necessary
torsion states, the three visible states, the plurigenera and
the strict cover-volume comparison directly from
the basket and \eqref{eq:RR46}.


\section{Coupled primes, marked curves and conductor descent}
\label{sec:coupled-conductor-obstructions}

Here the full remainder of a low-degree canonical member carries
information that is lost on an isolated visible component. We first
use simultaneous characters for two visible primes. We then combine
marked-curve restrictions with cancellation of the complete negative
part, and finally descend a cubic square through the full conductor.

\subsection{Exclusion of C2-02354: two coupled visible primes}

\label{sec:c2354-coupled34}

The exclusion uses two visible primes and their simultaneous intersection
characters. The first step shows that every actual torsion translate has
the same relevant section dimensions as the canonical class; this is the
precise input needed for the prime calculation.

\begin{theorem}
\label{thm:c2354-coupled34}
The basket \texttt{C2-02354} is not geometrically
realizable.
\end{theorem}

\begin{proof}
The exact data are $I=420$, $v=1$, $K^3=1/420$, $\chi=2$, and
\[
 B=2(1,2)+(1,3)+(1,4)+(2,5)+2(3,7)+(8,21).
\]
We have $q=p_g=0$ and $h^{2,0}=1$.  In the degrees used below,
\[
 P_4=P_6=0,\quad P_5=P_7=P_8=P_{11}=1,
 \qquad P_{12}=P_{13}=P_{17}=2.
\]

\emph{Step 1: an actual-doubling filter.}
Write $H_j(\tau)=h^0(Z,jK_Z+\tau)$ for an actual torsion class
Cartier in codimension two.  The independently exhaustive paired
torsion calculation leaves precisely the following three profiles
in the displayed degrees.  The first row has $H_j(\pm\tau)=P_j$ in
every checked degree, and the other two are interchanged by inversion.
\[
\begin{array}{c|cc|cc|cc|cc|cc}
 &\multicolumn{2}{c|}{j=4}&\multicolumn{2}{c|}{j=5}
 &\multicolumn{2}{c|}{j=6}&\multicolumn{2}{c|}{j=8}
 &\multicolumn{2}{c}{j=13}\\
 &H&H^-&H&H^-&H&H^-&H&H^-&H&H^-\\\hline
 0&0&0&1&1&0&0&1&1&2&2\\
 1&1&1&0&0&0&1&0&1&0&1\\
 2&1&1&0&0&1&0&1&0&1&0
\end{array}
\]
Suppose an actual $\tau$ has a noncanonical profile.  Squaring
nonzero sections of $4K_Z+\tau$ and $4K_Z-\tau$ gives
$H_8(2\tau)>0$ and $H_8(-2\tau)>0$.  The same exhaustive list
applied to the actual class $2\tau$ forces its profile to be row~0.
In particular $H_5(2\tau)=1$.  Orient $\tau$ so that its row is
row~1.  The product of nonzero sections of $5K_Z+2\tau$ and
$8K_Z-\tau$ contradicts $H_{13}(\tau)=0$.
Thus only row~0 can occur for an actual class.  In particular
\begin{equation}
\label{eq:c2354-torsion-zero34}
 H_2(\tau)=H_3(\tau)=H_4(\tau)=H_6(\tau)=0
 \quad\hbox{for every actual }\tau.
\end{equation}
There is no effective $K_Z+\tau$ either, since its square would
contradict the degree-two vanishing for the actual class $2\tau$.
Only actual additions of torsion classes have been used; no additive
law is imposed on virtual basket characters.

\emph{Step 2: the first visible prime.}
For a visible prime $T$ on the terminal minimal model put
$n=IK^2T$ and $M=IKT^2$.  The exact simultaneous Weil--RR
congruences with $-n\le M\le n^2$ give
\[
\begin{array}{c|ccccc}
 n&1&2&3&4&5\\\hline
 M&1&4&9&-4,16&-5,5,25.
\end{array}
\]
Hodge equality descends to an effective $nK_Z+\tau$.
Equation~\eqref{eq:c2354-torsion-zero34} removes the equality
states through degree four, so every visible prime has $n\ge4$.
Take $0\ne a\in R_5$.  Since $5<2\cdot4$, its full visible part
is one prime $T$ of coefficient one.  Write its full terminal
divisor as $T+E$, with $E$ invisible, and normalize $\lambda(K)=1$.
Then $\lambda(T)=5$.

The whole canonical divisor of $a$ cannot be a reduced prime:
its complete restriction ring would be a domain with
\[
 G_6=1,\qquad G_{11}=1,\qquad G_{17}=0,
\]
contrary to multiplication.  Ampleness downstairs therefore gives
$n_T<5$.  The only remaining state is
\begin{equation}
\label{eq:c2354-first-prime34}
 (n_T,M_T)=(4,-4),\qquad IK^2E=1.
\end{equation}
No self-intersection is assigned to the strict positive-degree
component of $E$; throughout $E$ denotes the entire terminal
remainder, including exceptional components.

\emph{Step 3: a second visible prime in the eleventh member.}
Take $0\ne e\in R_{11}$. The plateau $R_{17}=aR_{12}$ and
Lemma~\ref{lem:visible-invisible-bracket-division}, applied on the
canonical model with $(m,n)=(5,11)$, give
$\operatorname{ord}_T(D_a e)=t:=\operatorname{ord}_T(e)$, since
$5t-11\ne0$. Hence $t\ge1$.  Hodge positivity gives $t\le2$.  If $t=2$,
the remainder has positive Hodge mass one and canonical scaled
degree $11-2\cdot4=3$, below the visible-prime floor.  Hence $t=1$.
The remaining degree is seven and its Hodge mass is six.  There
is precisely one further visible prime, again of coefficient one:
\[
 \operatorname{div}_Y(e)=T+U+J,\qquad
 T\ne U,\quad \lambda(U)=6,\quad 4\le n:=IK^2U\le7,
\]
where $J$ is the full invisible remainder.

\emph{Step 4: simultaneous prime intersections.}
Set $M=IKU^2$ and $\ell=IKTU$.  Effectiveness of $E,J$, whose
supports contain neither relevant visible prime, yields
\[
 0\le\ell\le\min\{5n,48,11n-M\}.
\]
Hodge index gives
\[
 (4n-\ell)^2\le20(n^2-M).
\]
\emph{Only when equality holds} in this inequality does the
Hodge--Lefschetz comparison additionally require
\[
 20(6-n)=4n-\ell.
\]
The individual RR congruences and their polarization use the same
local characters:
\[
 n_D\equiv\sum_Q\frac I{r_Q}b_Q^2 i_{D,Q}\pmod I,
 \quad
 M_D\equiv-\sum_Q\frac I{r_Q}
       \overline{b_Qi_{D,Q}}(r_Q-\overline{b_Qi_{D,Q}})\pmod{2I},
\]
\[
 \ell\equiv\sum_Q\frac I{r_Q}b_Q^2i_{T,Q}i_{U,Q}\pmod I.
\]
The full enumeration of all $246960$ ordered individual local
assignments, followed by all compatible ordered pairs, leaves
exactly
\[
 (n,M,\ell)=(5,5,0)\quad\hbox{or}\quad(6,36,24).
\]
There are twelve labelled paired assignments in each state.
All strict Hodge-inequality pairs were in the enumeration;
none survive the other necessary conditions.

The second state is a Hodge equality, so $U_Z\sim6K_Z+\tau$,
contradicting~\eqref{eq:c2354-torsion-zero34}.  In the first state
the integral divisor $D=U-T-K$ satisfies
\[
 IK^2D=5-4-1=0,\qquad
 IKD^2=5-4+1-2\cdot0-2\cdot5+2\cdot4=0.
\]
Integral null-divisor ascent gives
$U_Z-T_Z\sim K_Z+\tau$ for an actual codimension-two Cartier
torsion class.  Adding the full effective invisible remainder gives
\[
 U_Z+E_Z\sim6K_Z+\tau,
\]
again contradicting~\eqref{eq:c2354-torsion-zero34}.
\end{proof}


\subsection{Exclusion of C2-02355: marked curves and full cancellation}

\label{sec:2355-exclusion46}

\begin{theorem}\label{thm:2355-exclusion46}
There is no canonical threefold with $\chi=2$, $K^3=1/420$ and
terminalization basket
\[
 B=2(1,2)+(1,3)+(1,4)+(2,5)+(3,7)+2(3,8)+(5,12).
\]
This is the basket denoted C2-02355 in the numerical table.
\end{theorem}

Here $I=840$, $v=2$, and
\[
 P_1=P_2=P_3=P_4=P_6=P_9=0,
 \qquad P_5=P_7=P_8=P_{10}=P_{11}=P_{13}=1,
 \quad P_{12}=2.
\]
The profiles in Section~\ref{subsec:late-torsion46} imply
$H^0(jK_Z+\tau)=0$ for every actual torsion class Cartier in
codimension two and $j=1,2,3$.

\paragraph{The saturated visible state}
\begin{lemma}\label{lem:late-k3-16846}
If $N$ is a Du Val K3 surface and $H$ an integral ample Weil
divisor with $H^2=1/168$, then $H^0(N,6H)=0$.
\end{lemma}
\begin{proof}
This is the K3 Riemann--Roch calculation in
\cite[Section~18, equation~(18.1) and the table following it]{LBVolume}.
\end{proof}

Let $a\in R_5$ be nonzero. The full dualizing trace has dimension
$P_6+1=1$, so its divisor has exactly one distinct visible prime
$T$, possibly with multiplicity. The actual single-divisor
congruences and $-n\le M\le n^2/2$ show that every visible state
with $n<7$ is one of
\[
 (2,2),\quad(4,8),\quad(6,18),\quad(5,-5).
\]
The first three are Hodge equalities. Null descent would give
$T_Z\sim jK_Z+\tau$ for $j=1,2,3$, contrary to the just proved
twisted vanishings. The state $(5,-5)$ requires a separate argument.

If this state occurs with multiplicity two in $(a)$, it exhausts
the total degree ten. Thus $(a)_Z=2T_Z$, and the integral class
$T_Z-2K_Z$ is half-canonical. Lemma~\ref{lem:late-no-half46}
applies to the unique index-twelve constituent and excludes this.
If its multiplicity is one, saturated adjunction
\cite[Lemma~9.9]{LBVolume} identifies the normalization $N$
of $T_Z$ as a Du Val K3, with integral ample $H$ and
$H^2=5/840=1/168$. The full polar of the actual quintic gives a
nonzero section of $6H$: at a generic curve point, if $a=te$,
its value on $T$ is $e\delta(t)$. It is regular and nonzero there,
and finite normalization preserves regularity. Curves contracted
by the polarization map to points, so reflexive extension gives
the section on $N$. The full remainder $e$ adds zeros and creates
no poles. This contradicts Lemma~\ref{lem:late-k3-16846}.

Consequently the unique visible prime of this particular quintic
has degree at least seven. Since the entire quintic has degree
ten, its full divisor on $Y$ is
\[
 (a)=T+E,\qquad T\text{ reduced},\qquad E\text{ invisible},
 \qquad IK^2E\le3.
\]
No global lower bound seven for all visible primes is being asserted.

\paragraph{A complete restriction of dimension zero}
The full invisible-square bounds from
\eqref{eq:late-invisible-max46} are
\[
 U(1)=-17,\qquad U(2)=-68,\qquad U(3)=-13.
\]
In particular $U(n)+2n=-15,-64,-7<0$.
Choose $0\ne z\in R_7$. Multiplication by the actual full polar
$\theta_6\in H^0(\mathcal Q_6)$ gives
\[
 R_7\longrightarrow H^0(\mathcal Q_{13}),\qquad
 z\longmapsto\theta_6z.
\]
The target has dimension $P_{13}-P_8=0$. Since $\theta_6$ is
nonzero at the generic point of the reduced visible prime, $z$
vanishes along $T$. Lemma~\ref{lem:late-full-cancellation46}, with
$h=2$, now gives $z/a\in R_2=0$, a contradiction.
This excludes the entire basket, including reducible or nonreduced
invisible remainders and canonical-null components.

More generally the same argument gives the rank obstruction
\[
 P_j-P_{j-m}\le P_{m+1+j}-P_{j+1}
\]
whenever the full negative-part bound is strict for $h=j-m$.
The geometric statement concerns the complete restriction sheaf,
including its exceptional polar direction, rather than an ambient
canonical section in the missing degree $m+1$.


\subsection{Exclusion of C2-02357: descent through the full conductor}

\label{sec:2357-exclusion46}

\begin{theorem}\label{thm:2357-exclusion46}
There is no canonical threefold with $\chi=2$, $K^3=1/420$ and
terminalization basket
\[
 (1,2)+2(1,3)+(1,4)+(2,5)+(3,7)+(4,9)+(7,18).
\]
This is C2-02357 in the numerical table.
\end{theorem}

\paragraph{Actual torsion.}
\begin{lemma}\label{lem:2357-torsion46}
For basket $02357$, every torsion Weil class Cartier in codimension
two is trivial.
\end{lemma}
\begin{proof}
An element of prime order $p$ gives a connected canonical cyclic
cover of volume $p/420$ and $p_g'=P_2'=P_3'=0$. For odd $p$,
\eqref{eq:late-euler46} and the table give
$\chi'\ge(p+3)/2$, hence $h^{2,0}\ge(p+1)/2$.
The two-form bounds
\[
 h^{2,0}\ge4\Longrightarrow V\ge(h^{2,0}-3)/60,
 \qquad P_2=0, h^{2,0}\ge3\Longrightarrow V\ge1/14
\]
follow from \cite[Section~12, equation~(12.5) and the following paragraph]{LBVolume}.
For the first, use $1/18$ for $4\le h^{2,0}\le6$, $1/10$
for $7\le h^{2,0}\le9$, and $(h^{2,0}-3)/54$ thereafter.
The first excludes $p\ge11$, since it
would imply $p\le7$; the second excludes $p=5,7$.
For $p=2$ the profile must be self-inverse, hence canonical, giving
$\chi'=4$ and the same contradiction. For $p=3$, a canonical
inverse pair gives $\chi'=6$ and is excluded in the same way.
The other pair gives $P_4'=2$. The uniform pencil estimate
\cite[Corollary~13.7]{LBVolume} holds for every canonical
threefold of general type and gives
\[
 P_n\ge2\Longrightarrow V\ge\frac{4}{3n(n+1)^2}.
\]
For this degree-three cover it yields
$V'\ge1/75>1/140=3/420$. Thus no prime order occurs.
\end{proof}

\begin{lemma}\label{lem:late-k3-18046}
An integral ample Weil divisor $H$ on a Du Val K3 surface with
$H^2=1/180$ has $h^0(6H)=h^0(7H)=0$.
\end{lemma}
\begin{proof}
The exceptional root rank is at most nineteen and
$\sum q_P(H)-1/180\in2\mathbb Z$. The nine-primary denominator
forces $A_8$ or $A_{17}$. An $A_{17}$ leaves no room for the
five-primary part, and two $A_8$ factors also cost too much rank.
Thus exactly one $A_8$ occurs, with class $1,2$ or $4$ up to sign.
The five-primary factor is $A_4$ or $A_9$. With $A_9$ the remaining
rank is at most two, and all two-primary denominators are at most
two, so the denominator four cannot occur. Thus $A_4$ is required.
At most seven ranks remain. Enumerate every nonzero local class on
$A_1,\ldots,A_7$, $D_4,\ldots,D_7$, $E_6,E_7$ within that rank,
using the correction norms and their actual multiples. The complete
congruence list, allowing additional Cartier factors, is
\[
\begin{array}{l|c|c}
 \text{non-Cartier profile}&\text{rank}&(h_1,\ldots,h_7)\\\hline
 A_1(1)+A_2(1)+A_3(1)+A_4(2)+A_8(1)&18&(0,0,0,0,0,0,0)\\
 A_2(1)+A_4(2)+A_8(1)+D_5(\mathrm{spin})&19&(0,0,0,0,0,0,0)\\
 A_1(1)+A_3(1)+A_4(2)+A_8(2)&16&(0,0,0,1,1,0,0)\\
 A_4(2)+A_8(2)+D_5(\mathrm{spin})&17&(0,0,0,1,1,0,0).
\end{array}
\]
Here $h_m=2+(m^2/180-\sum q_P(mH))/2$. The two relevant norm
sums are $21/5$ and $769/180$. The complete finite type list and
multiset enumeration are in \texttt{verify\_small\_k3\_46.py}; the
only test is the displayed congruence, with no positivity cutoff.
\end{proof}

\begin{lemma}\label{lem:2357-saturated46}
The first quintic has no visible prime with $(n,M)=(7,-7)$.
\end{lemma}
\begin{proof}
Its total degree is fifteen and the full dualizing trace has
one dimension, so such a prime has multiplicity one or two.
Saturated adjunction, \cite[Lemma~9.9]{LBVolume}, gives a
Du Val K3 normalization with $H^2=7/1260=1/180$.
In multiplicity one the full polar descends to a nonzero section
of $6H$: restrict the actual polar to the finite normalization,
then extend reflexively across the remaining finite set.
In multiplicity two the normalized square
(Lemma~\ref{lem:polar-square46}) has visible order $2k-2=2$ and
at least twice the order of the quintic at every invisible prime.
Dividing it by the quintic gives a global section of $7K_Z$ with
visible order zero. It restricts nontrivially to $7H$ on the
normalization; possible polarization-null curves map to points,
and reflexivity supplies the extension. Both conclusions contradict
Lemma~\ref{lem:late-k3-18046}.
\end{proof}

\paragraph{The exact input and the permissible first states}
Put $K=K_Y$ on a terminal $\mathbb Q$-factorial minimal model over the
canonical model $Z$.  The numerical data are
\[
 I=1260,\qquad IK^3=3,\qquad \chi=2,\qquad K^3=1/420,
\]
\[
 B=(1,2)+2(1,3)+(1,4)+(2,5)+(3,7)+(4,9)+(7,18).
\]
In the degrees used here,
\[
\begin{array}{c|rrrrrrrrrrrrrrrrrrrr}
 m&1&2&3&4&5&6&7&8&9&10&11&12&13&14&15&16&17&18&19&20\\\hline
 P_m&0&0&0&0&1&0&1&1&1&1&1&2&2&2&3&3&3&5&3&5.
\end{array}
\]
The codimension-two-Cartier torsion group is trivial by Lemma~\ref{lem:2357-torsion46}.  In particular the Hodge-equality states that
would give effective $K$, $2K$, $3K$, or $4K$ are impossible.

Choose $0\ne a\in R_5$.  The full-member trace and Lemma~\ref{lem:2357-saturated46} give
\[
 \operatorname{div}_Y(a)=T+E,
\]
where $T$ is one reduced visible prime, $E$ is the full invisible
remainder, and, after normalizing the Hodge functional by
$\lambda(K)=1$, one has $\lambda(T)=5$ and $\lambda(E)=0$.
All canonical-null components of $E$ remain in this notation.
Writing $n=IK^2T$ and $M=IKT^2$, the remaining possibilities are
\[
 (n,M)=(8,-2),(10,10),(11,17),(13,33),(14,42),(15,75).
\tag{E357.1}
\]
Indeed the remaining visible RR states through degree fifteen are these
and $(7,-7)$; the last is excluded by Lemma~\ref{lem:2357-saturated46}.  The intervening Hodge
equalities descend to the forbidden lower canonical degrees.  The
state $(15,75)$ means that the full member is prime on $Z$, not that
its full pullback has no exceptional components.

\paragraph{Small invisible primes and reducedness of the full member}
For an invisible prime $F$ put $s=IK^2F>0$ and $N=IKF^2$.
The simultaneous single-divisor residues and strict invisible Hodge
inequality begin as follows:
\[
\begin{array}{c|l}
 s&\text{largest permissible square values}\\\hline
 1&-23,-93,-233,\ldots\\
 2&-22,-92,-232,\ldots\\
 3&-207,-417,\ldots.
\end{array}
\]
Lemma~\ref{lem:late-ruling46} excludes a prime with $s+N\le-L<0$ if
\[
 W=\left\lfloor{3Is\over L}\right\rfloor<70,
 \qquad 2W<L.
\]
For $(s,N)=(1,-93),(2,-232),(3,-207)$ the values $(L,W)$ are
$(92,41),(230,32),(204,55)$, respectively.  All still smaller square
values are excluded because both inequalities improve as $L$ grows.
Consequently degree-one primes have square $-23$, degree-two primes
have square $-22$ or $-92$, and degree-three invisible primes do not
exist.  These are prime-divisor statements, not assertions about a
reducible invisible negative part.

If such a prime $F$ appears in $E$, put $l=IKTF$.  Nefness and
effectivity, with the full remainder retained, give
\[
 0\le l\le IKTE=5n-M,
\]
and the threefold Hodge index theorem gives
\[
 (ns-3l)^2\le(n^2-3M)(s^2-3N).
\tag{E357.2}
\]
The actual simultaneous local characters also satisfy
\[
 l\equiv\sum_{(b,r)\in B}{I\over r}b^2i_Ti_F\pmod I.
\tag{E357.3}
\]
Here each actual divisor uses the same assignment in its own two
single-divisor RR identities and in (E357.3).  This is the simultaneous identity \eqref{eq:mixed-character46}; no
additional law for adding virtual local characters is imposed.

Enumerating all $408240$ ordered assignments for this basket leaves,
for the first three states in (E357.1), only
\[
\begin{array}{c|c|r|c}
 (n,M)&(s,N)&l&(IK^2(T+F),IK(T+F)^2)\\\hline
 (8,-2)&(1,-23)&26&(9,27)\\
 (10,10)&(2,-22)&30&(12,48)\\
 (11,17)&(1,-23)&27&(12,48).
\end{array}
\]
All three sums are Hodge equalities.  Null descent and torsion
triviality would give a nonzero section of $3K_Z$ or $4K_Z$, which
is impossible.  Thus every positive-degree prime of $E$ has degree
at least four in these states.  Since $IK^2E=7,5,4$, respectively,
its pushforward $E_Z$ is one reduced prime.

For $(n,M)=(13,33)$, the same mixed test permits no degree-one prime
in $E$, whose degree is two; hence $E_Z$ is again one reduced prime.
For $(14,42)$ its degree is one, so this is automatic.  For $(15,75)$
there is no positive-degree remainder.  We have proved that the
\emph{entire} canonical quintic $D=(a=0)$ is reduced in every case.
This conclusion does not discard nonreduced canonical-null components
on $Y$.

\paragraph{Two numerical states disappear before any local classification}
The canonical inertia orders lie in
$\{2,3,4,5,6,7,9,18\}$, and $1260K_Z$ is Cartier by Proposition~\ref{prop:actual-source-transfer46}.  Exact primitive-mark enumeration gives no curve of
scaled degree $1\le w\le9$ satisfying
\[
 -2+\sum(1-1/r_x)\le {12w\over1260},\qquad
 \sum k_x/r_x\equiv w/1260\pmod{\mathbb Z}.
\]
The first arithmetic possibility is $w=10$ with marks $(2/7,13/18)$.
Lemma~\ref{lem:late-intrinsic-semigroup46} applied to $(8,-2)$ would
produce such a curve with $w\le n+M=6$.  Therefore $(8,-2)$ is
impossible.

There is a second useful full-cycle test.  Every positive-degree
component of the complete effective intersection cycle $T\cdot E$
maps to a curve where the reduced full member $D$ is singular.  This
includes intersections with exceptional components of $E$: if $D$
were regular at the generic image point, the transverse Cartier
curve in the canonical Du Val surface would force the ambient space
to be regular, leaving no exceptional component there.  Generic
threshold adjunction therefore has coefficient at most six.  Each
image curve of degree $w$ satisfies the weaker genus-zero conditions
\[
 -2+\sum(1-1/r_x)\le {6w\over1260},\qquad
 \sum k_x/r_x\equiv w/1260\pmod{\mathbb Z}.
\tag{E357.4}
\]
The full intersection on the terminal model has integral positive
multiplicities at generic curve points.  Its finite pushforward
multiplies these by positive integral mapping degrees.  Thus its total
scaled degree is a nonnegative integral combination of the permitted
image degrees.  No individual component or null cross term is dropped.

Through degree forty-two the permitted positive degrees in (E357.4)
are
\[
 \{10,14,20,28,30,35,36,40,42\}.
\tag{E357.5}
\]
For $(13,33)$ the total is $IKTE=32$. It is not a nonnegative
integral combination of $10,14,20,28,30$, which excludes this state.
For later use the possible degree partitions of the other totals are
\[
\begin{array}{c|l}
 28&28;\ 14+14\\
 38&10+28;\ 10+14+14\\
 40&40;\ 10+30;\ 20+20;\ 10+10+20;\ 10+10+10+10.
\end{array}
\tag{E357.6}
\]
Repeated summands may denote distinct curves or a multiple
intersection component; this distinction is in the proof below.

\paragraph{The terminal intersections needed for descent}
We now need only $(n,M)=(10,10),(11,17),(14,42)$ and the prime state.
In the nonprime states let $F$ be the strict transform of the reduced
prime $E_Z$, and write $E=F+N_0$, where $N_0$ is the full effective
canonical-null remainder.  The complete prime ruling bounds and the
mixed tests give
\[
\begin{array}{c|r|l|l}
 (n,M)&s&\text{prime square values before mixing}&
                 \text{surviving }(IKF^2,IKTF)\\\hline
 (10,10)&5&-15,-85,-225&(-15,40),(-225,40)\\
 (11,17)&4&-18,-88&(-18,38)\\
 (14,42)&1&-23&(-23,28).
\end{array}
\tag{E357.7}
\]
For degree four the next value $-228$ is excluded by $(L,W)=(224,67)$;
for degree five the next value $-295$ is excluded by $(L,W)=(290,65)$.
The degree-five value $-225$ \emph{must be retained}: its ruling
bound has $W=85$, so the small-prime lemma does not exclude it.
All surviving rows nevertheless give
\[
                         IKTN_0=0.
\tag{E357.8}
\]
For the $(-225,40)$ row, $IKFN_0=210$ and $IKN_0^2=-210$ are
not ruled out.  We do not claim that every exceptional component is
point-centred in that row.

What (E357.8) does imply is exactly what is needed: $Z$ is regular
at the generic point of every curve contained in $T_Z$.  To see this,
suppose a canonical singular curve $\Gamma$ were contained in $T_Z$.
The transverse canonical surface is Du Val, and its crepant
terminalization is its minimal resolution.  The full canonical
section $a$ vanishes along $\Gamma$, so every exceptional curve of
that transverse resolution has positive coefficient in the full
pullback divisor.  A strict transform of a curve branch of $T_Z$
through the transverse singular point meets the exceptional fibre.
Consequently $T\cdot N_0$ has a positive component dominating
$\Gamma$, and $K_Z\Gamma>0$ contradicts (E357.8).
The same argument applies to the prime state: $E$ is canonical-null
and $IKTE=5n-M=0$.  Thus all generic local curve calculations involving
$T_Z$ take place in a smooth ambient threefold.  Exceptional curves
lying entirely in $F_Z$ away from $T_Z$ are harmless.

\paragraph{Four sections do not vanish on the visible component}
Let $b\in R_7$, $c\in R_8$, $d\in R_9$, and $e\in R_{11}$ be nonzero.
Let $Q_j$ be the global sections of the complete degree-$j$ restriction
sheaf on the full reduced member $D$.  The full polar
$\theta\in Q_6$ is nonzero at the generic point of $T_Z$ and zero
at the generic point of $F_Z$.  Restriction and vanishing give
\[
 \dim Q_{13}=\dim Q_{14}=\dim Q_{17}=\dim Q_{19}=1,
 \qquad Q_5=Q_{10}=0.
\tag{E357.9}
\]
In particular the zero values are complete restriction statements,
not merely a difference of two lower bounds for section dimensions.

First $b$ cannot contain $T$.  Otherwise Hodge-mass positivity gives
$\operatorname{ord}_T(b)=1$, since $\lambda(T)=5$ and
$\lambda(\operatorname{div}b)=7$.  Lemma~\ref{lem:visible-invisible-bracket-division}, with
$(m,n,l)=(5,7,1)$, gives a nonzero $\{a,b\}/a\in R_8$ of
visible order zero.  In the prime case division by $a$ already
contradicts $P_2=0$.  In the two-component case $b$ cannot also contain
$F_Z$, for the same reason.  Thus $b^2$ and
$\theta\,\{a,b\}/a$ give two independent sections of $Q_{14}$:
one is zero on $T_Z$ and nonzero on $F_Z$, while the other is nonzero
on $T_Z$ and zero on $F_Z$.  This contradicts (E357.9).

Next $d$ cannot contain $T$. Its order would again be one by
Hodge mass. The same lemma with $(m,n,l)=(5,9,1)$ gives a nonzero
$\{a,d\}/a\in R_{10}=\mathbb C a^2$ of visible order zero,
whereas $a^2$ has order two.

For $c$ use full negative-part cancellation with $h=8-5=3$.
The exact strict-invisible square maxima in degrees one through five
are
\[
 M_{\max}(1),\ldots,M_{\max}(5)=-23,-22,-207,-18,-15.
\]
Thus $M_{\max}(s)+3s<0$ through degree four, and the only degree-five
endpoint is zero.  The strict cases exclude containment of $T$ in
$c$ for $(11,17)$ and $(14,42)$.  In the state $(10,10)$, such
containment would force the full positive and negative parts
$H-J\sim3K$ of $\operatorname{div}c-\operatorname{div}a$ to satisfy
\[
 IK^2J=5,\qquad IKJ^2=-15,\qquad IKHJ=0.
\]
Lemma~\ref{lem:late-zero-contact46} makes $J$ locally zero or $-3K$ on a
crepant terminal graph model.  Its linear RR congruence therefore
requires $\gcd(3,1260)\mid IK^2J$, contrary to $3\nmid5$.
In the prime case normal division gives the contradiction $P_3=0$
directly.  Consequently $c|_{T_Z}\ne0$ in every surviving state.

It remains to show $e|_{T_Z}\ne0$.  This is immediate by division
in the prime case, and by full cancellation with $h=6$ when
$(n,M)=(14,42)$.  Suppose instead that $e$ contains $T$ in one of
the states $(10,10),(11,17)$.  It cannot contain $F_Z$ as well,
since $P_6=0$.  The nonzero section $\theta^2b$ spans $Q_{19}$.
As $ce$ vanishes on $T_Z$, it is zero in $Q_{19}$; hence $c$ vanishes
on $F_Z$.  Since $b^2$ and $\theta c$ are nonzero on $T_Z$ and span
the same one-dimensional $Q_{14}$, $b$ also vanishes on $F_Z$.
Now $be$ vanishes on both components of $D$, so normal division on
$Z$ gives $be=a h$ for an actual $h\in R_{13}$.  The restriction of
$h$ is zero or a multiple of $\theta b$, and in either case vanishes
on $F_Z$.  Taking the valuation along $F_Z$ gives
\[
                         \operatorname{ord}_{F_Z}b\ge2.
\]
For $(10,10)$ this contradicts
$IKT\operatorname{div}b=70<2IKTF=80$.
For $(11,17)$ it leaves at most one unit of intersection with $T$
after removing $2F$.  The remaining divisor has scaled degree at
most $21-8=13$ and positive Hodge mass seven, so it contains a
visible prime $W\ne T$ of degree between seven and thirteen.
The possible strict states for $W$ are
\[
 (7,-7),(8,-2),(10,10),(11,17),(13,33).
\]
The mixed RR residues with $T=(11,17)$ have respective smallest
nonnegative intersections
\[
                           49,6,60,17,71,
\]
all exceeding the available bound $IKTW\le1$.  This final
contradiction proves that $e$ does not contain $T$.

\paragraph{A forced cubic root on the normalization}
All four restrictions are now nonzero on $T_Z$.  By the
one-dimensional spaces in (E357.9), there are nonzero constants
$A,B,C$ such that on the full member
\[
 b^2=A\theta c,\qquad cd=B\theta e,\qquad ce=C\theta^2b.
\tag{E357.10}
\]
In the function field of $T_Z$ these imply
\[
 be=AC\theta^3,\qquad e^3={AC^2\over B}\theta^4d.
\]
After absorbing a nonzero constant into $d$, put
\[
                        \xi={e\over\theta}.
\]
This is a rational weight-five section on $T_Z$ satisfying
\[
                          \xi^3=\theta d.
\tag{E357.11}
\]
It has no pole on the normalization of $T_Z$, since its cube is
regular.  Extend it by zero on the normalization of $F_Z$, when that
component is present.  Equation (E357.11) remains valid there
because the full polar vanishes on the invisible component.
The resulting section of the normalization of the full member is
nonzero on $T_Z$.

We shall prove that $\xi^2$ descends to the complete weight-ten
restriction sheaf on $D$.  This contradicts $Q_{10}=0$.
At any generic curve point involving $T_Z$ the ambient threefold
is smooth by (E357.8).  The restriction problem is therefore a
reduced plane-curve problem over the curve's function field, with
$\theta=v(f)$ the actual full Jacobian polar.  The canonical frame
makes (E357.11) an ordinary equality of local functions.

We use Lemma~\ref{lem:late-cubic-descent46}:
for a reduced plane curve of delta invariant at most four, a regular
normalization element satisfying $\xi^3=v(f)d$ can fail to have its
square in the original curve ring only when $\delta\ge3$ and $d$
is a unit.  The lemma includes zero values on some normalization
branches and the case in which $v(f)$ is identically zero on those
branches.  It is applied to the actual full polar, not an arbitrary
section with the same divisor degree.

Two elementary global budgets make this local lemma effective.
At an intrinsic nonnormal curve of $T_Z$, the conductor contributes
$2\delta_Tw$ to $n+M$, where $w=IK_Z\Gamma$.  More precisely the
intrinsic polar contributes the conductor plus the zero divisor of
the restricted ambient two-form.  A singular plane branch of
multiplicity $m_i$ has additional two-form vanishing at least $m_i-1$.
Thus a cusp costs at least three units per image degree, while a node
costs two.  At a curve where $T_Z$ meets $F_Z$, the contact multiplicity
is counted in the complete budget $IKTE$; its contribution is not
replaced by the reduced support.

Every such curve has $w\ge10$ by (E357.4).  There is also a useful
section obstruction:
\[
 \text{(E357.4) and }1\le w\le15
 \quad\Longrightarrow\quad d|_\Gamma=0.
\tag{E357.12}
\]
Thus no such curve of degree at most fifteen supports a nonzero
restriction of $d$.  Indeed the only
primitive-mark patterns through degree fifteen are
\[
 (w;\text{marks})=(10;2/7,13/18),\quad(14;2/5,11/18).
\]
Their necessary degree-nine section costs are respectively $15/14$
and $11/10$, whereas $9w/1260$ is $1/14$ and $1/10$.  Thus the
actual ninth-canonical section vanishes on either curve.

\paragraph{The prime state $(15,75)$.}
Here the full intrinsic mass is ninety.  If $\xi^2$ failed to descend
at a curve $\Gamma$, then $2\delta w\le90$ and $w\ge10$, so
$\delta\le4$.  The local lemma would give $\delta\ge3$ and $d$ a
unit.  The same budget gives $w\le15$, contradicting
(E357.12).  Hence the square descends at every generic curve point.
This includes the nonnormal prime case; normality is not assumed.

\paragraph{The states $(10,10)$ and $(11,17)$.}
Their intrinsic masses are twenty and twenty-eight.  A cusp on
$T_Z$ would cost at least $3w\ge30$, and a singularity with delta
at least two costs at least $4w\ge40$.  Thus at any relevant curve
$T_Z$ is smooth or nodal.  If $F_Z$ does not pass through that curve,
the full conductor exponents on visible branches are zero or one,
and the square plainly descends.

Suppose $F_Z$ does pass through it.  On a visible branch write
$c_i=c_{T,i}+j_i$ for the full conductor exponent, where $j_i$ is
its contact with $F_Z$.  If $T_Z$ is smooth, $c_{T,i}=0$ and the
contact budget gives $j_i\le4$.  If $T_Z$ is nodal, each of its
branches meets $F_Z$, so both $j_i\ge1$; their sum is at most four,
and $c_{T,i}=1$, so again $c_i\le4$.
The full polar has order at least $c_i$ on that branch.  If $d$ is
not a unit at the curve, (E357.11) gives
\[
 2\operatorname{ord}\xi\ge
 2\left\lceil{c_i+1\over3}\right\rceil\ge c_i
 \qquad(1\le c_i\le4).
\]
The square then belongs to the full conductor on every visible
branch and is zero on all invisible branches, so it descends.
If $d$ is a unit, the only possible failure of this conductor test
is $c_i=3$.  For a smooth visible branch this requires contact three.
For a nodal visible component it requires contact two on that branch
and at least one on the other.  In either case the complete contact
budget gives $3w\le40$, hence $w\le13$, contradicting
(E357.12).  Thus these two states also give descent everywhere.

\paragraph{The state $(14,42)$.}
The intrinsic mass is fifty-six, so $\delta_T\le2$ at every intrinsic
nonnormal curve.  Away from $F_Z$ the local cubic-root lemma gives
square descent directly.  The full contact budget is twenty-eight;
by (E357.6) the intersection is one degree-twenty-eight component
of multiplicity one, or two degree-fourteen components counted with
multiplicity.  In particular each contact multiplicity is at most two.
If $T_Z$ is smooth at such a curve, the conductor exponent on its
branch is at most two, and integrality of (E357.11) already gives
square descent with zero on the other components.
If $T_Z$ is singular there, contact at most two forces multiplicity
two for $T_Z$, a smooth other component, and total contact exactly
two.  The full plane curve has
\[
                     \delta_D=\delta_T+2\le4.
\]
Failure of square descent would force $d$ to be a unit by the local
lemma.  But the curve has degree fourteen, and (E357.12) excludes
that unit.  This proves descent in the final nonprime state.

\paragraph{Extension through closed points and conclusion}
Along a curve contained only in $F_Z$, the proposed section is zero
and needs no extension argument.  The preceding discussion proves
regularity of $\xi^2$ at every codimension-one point of the full member.
At a closed non-Gorenstein point take the local canonical cover.
The canonical cover is Cohen--Macaulay and its canonical member is a
Cartier hypersurface.  The complete restriction sheaf is the relevant
character summand of this Cohen--Macaulay hypersurface quotient;
in particular it is $S_2$.  Regularity in codimension one therefore
extends across the remaining codimension-two points.  Taking invariants
preserves the extension.  We obtain a nonzero global section of the
\emph{complete} degree-ten restriction sheaf, contrary to $Q_{10}=0$.
All states in (E357.1) have now been excluded, proving the theorem.

The numerical verification \texttt{verify\_c2357\_46.py} enumerates
all labelled local assignments, the simultaneous mixed intersections,
all primitive mark multisets in the proven bounds, and all full-cycle
partitions. Its output explicitly retains the square $-225$. The
normalization and extension arguments above are separate geometric
proofs, not conclusions inferred from a successful computation.


\section{The remaining quartic and quintic configurations}
\label{sec:false-small-surfaces}

The last two non-octic baskets force a small canonical surface with an
elliptic pencil. We determine the entire member and its local sources
before using its polar equation. This order is important: the local
calculation must concern the actual surface, with all conductor and
marking conditions already established.

\subsection{Exclusion of C2-02349: the quartic surface configuration}

\label{sec:2349-preparation46}
The basket in this section is
\[
 B=5(1,2)+(1,3)+4(1,4)+(2,5)+(3,7)+(5,12),
 \quad I=420,\quad v=1,\quad\chi=2.
\]
Its numerical label is C2-02349. The source-transfer theorem
(Proposition~\ref{prop:actual-source-transfer46}) gives the canonical
index set $\{1,2,3,4,5,6,7,12\}$ and preserves the unique ordinary
points of indices five, seven and twelve. The index-five assertion
uses the absence of $(1,5)$ from the mandatory crepant index-five
subbasket; the other two use the high-index bound and the rational-tree
fixed-point argument. Thus all smooth-cover calculations below occur
at actual unchanged sources, not at virtual constituents.

\begin{lemma}\label{lem:2349-prime46}
The divisor of the quartic section is a reduced prime $D$ on $Z$.
Its terminal strict transform has $(n,M)=(4,16)$, and every component
of the full crepant remainder is point-centred.
\end{lemma}
\begin{proof}
Riemann--Roch gives $P_4=1$ and $P_5=P_6=0$.
Lemma~\ref{lem:adjacent-gaps46} gives a unique reduced visible
prime $T$ with $\lambda(T)=4$. The complete visible-prime states
through degree four are $(1,1),(2,4),(3,9),(4,16)$, all Hodge
equalities. Null descent therefore gives $\lambda(T)=n$, so
$(n,M)=(4,16)$. Its degree exhausts the canonical-model member.
On the terminal model the remainder $E=4K-T$ is effective, with
$K^2E=KE^2=0$. Restricting to a general canonical hyperplane, the
negative-definite exceptional intersection form forces $E$ to
have no curve-centred component. In particular the canonical model
is regular at the generic point of every curve in $D$: otherwise
the positive transverse pullback of $D$ on a Du Val resolution
would meet its strict transform and give a curve-centred remainder.
\end{proof}

Write $s\in R_4$ for the quartic section and $T=(s)_Z$.
The further dimensions used here are
$P_7=P_8=P_{10}=1$, $P_{13}=0$ and $P_{12}=3$.
Let $a\in R_7\setminus0$ and let $\theta$ be the coherent
weight-five polar section of $s$ along the normalization of $T$.
Both are nonzero there.  The actual bracket gives a twelfth section
whose restriction is $-7\theta a$.  Since
\[
 \dim\operatorname{im}(R_{12}\to S_T)=P_{12}-P_8=2,
\]
choose $q\in R_{12}$ giving an independent restriction.
The restricted pencil is generated by $\theta a$ and $q|_T$.
On a smooth resolution of the normalization of $T$, resolve this
pencil and take a normalized general connected moving curve $C$.
Put $\ell=420K\cdot C$, and let $d$ be the Stein degree.  Then
\[
 0<d\ell\le12\cdot4=48.
\]
The nonzero ambient two-form excludes a ruled normalization, so
$g(C)\ge1$.

The marked-curve adjunction for the full quartic and a
general twelfth member gives, for primitive canonical markings $k/r$,
\[
 2g(C)-2+\sum(1-1/r)\le\frac{17\ell}{420},
 \qquad \frac\ell{420}\equiv\sum\frac{k}{r}\pmod{\mathbb Z}.
\]
The allowed inertia orders divide the actual basket indices, hence
belong to $\{2,3,4,5,6,7,12\}$.  The restrictions of $\theta$ and
$a$ are regular sections of degrees $5K|_C$ and $7K|_C$.  Their
fractional orders therefore give the two independent budgets
\[
 \sum\frac{\overline{5k}}r\le\frac{5\ell}{420},
 \qquad
 \sum\frac{\overline{7k}}r\le\frac{7\ell}{420},
\]
where least nonnegative residues are used; zero residues are allowed.
Thus no extra positivity at a Cartier marking is assumed.

The exact enumeration, including proper inertia order six and
repeated markings, leaves precisely
\[
\begin{array}{c|c}
 \ell&\text{markings}\\\hline
 35&(1/12)\\
 40&(2/3,3/7)\\
 42&(1/2,3/5)\\
 45&(1/4,6/7)\\
 48&(2/5,5/7).
\end{array}
\]
Both section budgets are saturated in every row.  Moreover
$\ell>24$ and $17\ell/420<2$, so $d=1$ and $g(C)=1$.
The necessary marking at index five, seven or twelve is unique in
its row and hence rational over the pencil's function field.
The relatively minimal elliptic surface consequently has a section.
Ordinary trace and $P_5=0$ give $p_g=1$. The elliptic canonical
bundle formula then makes the smooth minimal model a K3 surface:
the base is $\mathbb P^1$, a section removes multiple fibres,
and $p_g=1$ gives $q=0$, $\chi=2$ and $K=0$.

\paragraph{Three unique-leading-monomial exclusions}
At a cyclic terminal point of type $\frac1r(1,-1,b)$ with canonical
mark $k/r$, the lower arc weights are
\[
 (w_x,w_y,w_z)=(-kb^{-1},kb^{-1},-k)\pmod r,
\]
with positive representatives.  A canonical section of degree $m$
has numerator character $-mb$ modulo $r$.
The unique constituents at indices twelve, seven and five ensure
the relevant sources below are smooth cyclic index-one covers.

The first three rows have the following local data:
\[
\begin{array}{c|c|c|c|c}
\ell&\text{mark}&(w_x,w_y,w_z)&\operatorname{ord}(\theta a)
 &\text{sole leading monomial}\\\hline
35&1/12&(7,5,11)&12&xy\\
40&3/7&(6,1,4)&1&y\\
42&3/5&(1,4,2)&1&x.
\end{array}
\]
In each row the displayed monomial is the only nonconstant
degree-twelve numerator monomial within the exact product-order
budget. A possible constant term of $q$ in the first row vanishes,
since every general curve passes through the marked point.
Thus $q-c\theta a$ has nonzero leading coefficient for general
$c$, and cannot vanish on the marked branch. All three rows are
excluded.

\paragraph{An initial polar equation excludes the fourth row}
For $\ell=45$, the index-seven mark $6/7$ has weights $(5,2,1)$.
The section $a$ is a unit on the branch, while $\theta$ has exact
order two.  The degree-twelve numerator space of weight at most two
is $\langle y,z^2\rangle$.

Write $f$ for the local numerator of the actual quartic section.
Its character is two modulo seven, so the weight-three initial part,
if present, is
\[
 f_3=z(Ay+Bz^2).
\]
Write the pulled-back invariant two-form as
$A_0\,dx\wedge dy+B_0\,dy\wedge dz+C_0\,dz\wedge dx$.
The coefficient characters are $0,5,3$, respectively.  Their
minimum arc weights are $0,4,1$.  Thus the associated vector field
can lower arc weight by at most one.  If $f_3=0$, every monomial
of $f$ has weight at least ten, and its polar derivative has weight
at least nine, contrary to the exact order two of $\theta$.
Consequently $(A,B)\ne(0,0)$.

Let $y_0,z_0$ be the leading arc coefficients at their lower weights;
they are allowed to be zero.  Exact polar order two implies that
$(y_0,z_0^2)\ne(0,0)$.  Since $f$ vanishes on the curve,
\[
 z_0(Ay_0+Bz_0^2)=0.
\]
The projective ratio $[y_0:z_0^2]$ therefore belongs to a fixed
finite set: $[1:0]$, or the single ratio satisfying
$Ay_0+Bz_0^2=0$.  The leading coefficients of $q$ and $\theta a$
are fixed linear forms in these two monomials, and the latter is
nonzero.  Their ratio has only finitely many possible values.
This is incompatible with a general member of the pencil.  Thus
$\ell=45$ is excluded as well.

\paragraph{The unique remaining marking configuration}
Only
\[
 \ell=48,\qquad (2/5,5/7)
\]
remains.  The fifth mark has arc weights $(4,1,3)$; the polar
section is a unit there, and $a$ has exact numerator order four.
The three possible low twelfth monomials are $x,y^4,yz$.
The quartic's initial part is $Az+By^3$, with $A\ne0$ forced
by the polar unit.  Hence the covering surface is smooth there
and the quotient surface has a Du Val $A_4$ point.

At the seventh mark the weights are $(3,4,2)$; the polar section
has exact order four and $a$ is a unit.  The low twelfth monomials
are $y,z^2$, and the quartic's weight-six initial part is
\[
 Ax^2+Byz+Cz^3.
\]
Here $A\ne0$ is necessary: otherwise the same finite-ratio
argument as above would exclude the pencil.  At least one of
$B,C$ is nonzero, because the polar section has exact order four.
These conditions allow a variable leading ratio, so the preceding
local argument stops.  We next determine its conductor and then its normal quotient graphs.

\paragraph{The two marked sections and the full conductor}
\label{sec:c2349-marked-conductor-continuation}

Use the prime quartic and pencil just proved. The terminal strict transform
$T$ has state $(n,M)=(4,16)$, and the remainder of the full quartic
member is point-centred.  On a smooth normalization model put
$A=K|_S$, so $A^2=1/105$.  The restricted twelfth pencil has
connected general elliptic curve $F$ with
\[
 420AF=48,\qquad 12A\sim_{\mathbb Q}F+B,
 \qquad B\ge0,\qquad AB=0.
\]
Its smooth minimal model is an elliptic K3 surface.  The two
primitive markings are $2/5$ and $5/7$, each rational over the
pencil field.  Denote them by $P_5,P_7$ on the generic fibre.
The coherent polar $\theta$ has weight five; the nonzero ambient
section $a$ has weight seven.

\begin{proposition}[Marked four-torsion]
\label{prop:c2349-marked-four-torsion-continuation}
With $P_5$ as origin, the section $P_7$ on this elliptic K3 has
exact order two or four.
\end{proposition}

\begin{proof}
The section budgets are equalities.  As rational zero
divisors on the smooth generic curve, including the canonical
cover fractions, they give
\[
 (\theta|_F)=\frac47P_7,\qquad
 (a|_F)=\frac45P_5.
\]
There are no additional zero contributions: their degrees are
respectively $5AF=4/7$ and $7AF=4/5$.  The two powers
$\theta^7$ and $a^5$ have the same weight thirty-five, hence
their ratio is an honest rational function on the generic curve.
Consequently
\[
 \operatorname{div}_F(\theta^7/a^5)=4P_7-4P_5.
\]
The two points lie over different actual terminal sources and
are distinct.  A divisor $P_7-P_5$ on a smooth genus-one curve
is principal only when the points coincide.  Thus its order
divides four and is not one, proving the assertion in the
Mordell--Weil group.  This argument takes place on the generic
smooth curve and does not require normality of the quartic.
\end{proof}

This torsion is on the elliptic generic fibre; no statement about the threefold torsion group is needed here.

\paragraph{The composite-inertia conductor budget}

The normalized polar square is an actual degree-ten section,
since the complete restriction piece has dimension
$G_{10}=P_{10}-P_6=1$.  The full polar cycle of $\theta$ has
scaled canonical degree twenty.  For a positive polar curve
$\Gamma$, put $w=420K\Gamma$.  Generic marked adjunction gives
\[
 \epsilon(\Gamma):=2g(\Gamma^\nu)-2+\sum(1-1/r_x)
 \le\frac{10w}{420},
 \qquad \frac w{420}\equiv\sum\frac{k_x}{r_x}\pmod{\mathbb Z}.
\]
At a singular curve of $T$ the generic-threshold bound improves
the coefficient ten to five.  In both tests the allowed inertia
orders are the actual set $2,3,4,5,6,7,12$, with repeated marks
allowed.  In particular an order-twelve mark is not replaced
by independent order-four and order-three costs.

The exact singular-curve list for $w\le10$ is
\[
\begin{array}{c|c}
 w&\text{primitive markings}\\\hline
 5&(3/7,7/12)\\
 7&(3/5,5/12)\\
 10&(1/2,2/3,6/7)\ \text{or}\ (1/6,6/7).
\end{array}
\]
The all-polar list gives precisely the following partitions of
total degree twenty, with repetitions representing multiplicities:
\[
 (5,5,5,5),\quad(5,5,10),\quad(5,15),\quad(10,10),\quad(20).
\]
These are necessary lists, not existence assertions.

\begin{proposition}[The only possible positive conductor]
\label{prop:c2349-single-conductor-continuation}
If the canonical quartic is not normal, it has exactly one
positive-degree curve of nonnormality $\Gamma$.  It is rational,
has $420K\Gamma=5$, and has markings $(3/7,7/12)$.
Its geometric generic transverse plane-curve singularity has
$\delta$-invariant one.  The conductor pushes forward to
$2\Gamma$, and
\[
 420A\cdot C_S=10,\qquad420A\cdot K_S=10.
\]
Here $C_S$ denotes the strict conductor on the smooth
normalization model; these equalities concern positive degree
and allow additional $A$-null curves.
\end{proposition}

\begin{proof}
The generic transverse ring is a reduced plane-curve ring in a
smooth terminal chart.  Its finite polar colength is at least
two, so every positive singular curve has $w\le10$ and occurs
in the displayed singular list.  A degree-seven curve cannot
occur in a polar cycle of total degree twenty: its coefficient
would be two, leaving degree six, which is not a sum of allowed
positive polar degrees.

For a reduced plane-curve singularity, the conductor divisor
on its normalization has total degree $2\delta$.  One local
proof uses its Gorenstein duality: if $R\subset\overline R$
and $\mathfrak c$ is the conductor, then
$\ell(R/\mathfrak c)=\ell(\overline R/R)=\delta$;
the exact quotient sequence gives
$\ell(\overline R/\mathfrak c)=2\delta$.
This includes residue extension degrees when the branches are
conjugate.  In particular each singular curve contributes at
least $2w$ to the scaled conductor degree.

Normalization adjunction and the nonzero restricted two-form
give
\[
 420A\cdot(K_S+C_S)=4+16=20,
 \qquad A\cdot K_S\ge0.
\]
Equality $420AC_S=20$ is impossible.  It would give $AK_S=0$.
The saturated-adjunction minimalization argument then contracts
only $A$-null curves and produces a Du Val K3 normalization
with integral ample Weil polarization of square $1/105$.
Lemma~\ref{lem:k3-square105-extension} excludes that polarization.  Equivalently, one can start from
the already proved K3 minimal model and descend $A$ through
all its $A$-null exceptional curves.

Thus the conductor degree is strictly less than twenty.
A degree-ten singular curve, two degree-five singular curves,
or a degree-five curve of $\delta\ge2$ would each contribute
at least twenty.  The only remaining possibility is one
degree-five curve with $\delta=1$, whose conductor degree is
ten.  Its marking inequality also forces genus zero.  The
last two equalities follow from adjunction.  Since the
canonical quartic is Cohen--Macaulay, failure of normality
must occur along a curve; the point-centred remainder identifies
these generic curves with the terminal ones just used.
\end{proof}

\paragraph{The polar cycle in that remaining branch}

\begin{proposition}
\label{prop:c2349-contact-three-continuation}
In the nonnormal alternative of
Proposition~\ref{prop:c2349-single-conductor-continuation},
the first nonzero normal derivative of the quartic along
$\Gamma$ is the third one.  Its positive polar cycle is
therefore either
\[
 4\Gamma\qquad\text{or}\qquad3\Gamma+\Gamma',
 \quad420K\Gamma'=5,\quad\Gamma'\ne\Gamma.
\]

\end{proposition}

\begin{proof}
Let $q\in R_{12}$ be the second restricted pencil member.
The first member is a nonzero multiple of $\theta a$.
Since $AB=0$, $q$ cannot vanish identically on the positive
curve $\Gamma$ where $\theta$ vanishes.  Its degree-twelve
budget on $\Gamma$ is $12\cdot5/420=1/7$.
The marks $(3/7,7/12)$ already cost $1/7$ and zero;
hence its exact numerator order at the seventh point is one.

Use the actual cyclic chart $\frac17(1,-1,3)$.
For this curve its lower arc weights are $(6,1,4)$, so the
order-one numerator is a nonzero multiple of $y$.
In particular $\operatorname{ord}_\Gamma y=1$ at that branch.
The previously established moving-fibre weights are $(3,4,2)$,
and the quartic's weight-six initial form is
\[
 f_6=A_0x^2+B_0yz+C_0z^3,
 \quad A_0\ne0,\quad(B_0,C_0)\ne(0,0).
\]
Since $\Gamma$ is in the singular locus, $f_z$ vanishes
along it.  In the weights $(6,1,4)$ its only possible
order-one term is $B_0y$.  Therefore $B_0=0$ and $C_0\ne0$.

Write the invariant ambient form as $\iota_\delta\eta$ in
a local volume frame.  For the moving-fibre weights
$(3,4,2)$, the three coefficient minima of $\delta$ are
$(8,2,0)$, by their characters $(5,3,0)$ modulo seven.
The exact polar order four consequently forces
$\delta z(0)\ne0$.  It follows that
$\delta^3f(0)=6C_0(\delta z(0))^3\ne0$.

The first derivative vanishes along the singular curve.
Where the lower derivatives vanish, the next derivative is
a well-defined regular canonical jet: its degree is $4+j$
for the $j$th derivative.  This follows by changing a canonical
frame and observing that every extra term contains a lower
derivative.  A nonzero second jet would therefore be a regular
degree-six section on $\Gamma$.  But at its seventh mark
the residue cost is $4/7$, already greater than its whole
degree $6\cdot5/420=1/14$.  Thus the second jet is zero,
whereas the third is nonzero, as just checked locally.

The derivation is transverse to $\Gamma$ near the seventh
point.  In a normal plane choose a formal coordinate with
$\delta=\partial_z$.  For a singular reduced plane curve
whose first nonzero derivative along this direction is the
third, the polar colength is at least three.  If the polar
has linear initial term, eliminate that parameter and the
remaining equation has order at least three; otherwise the
intersection multiplicity is at least four.  The total polar
degree twenty bounds the coefficient of $\Gamma$ by four.
Coefficient four consumes the cycle; coefficient three leaves
degree five, giving exactly the two displayed alternatives.
\end{proof}

\paragraph{The normalized square removes the split polar cycle}

\begin{proposition}[The remaining nonnormal branch is cuspidal]
\label{prop:c2349-cuspidal-conductor-continuation}
In the nonnormal alternative above, the full positive polar cycle
is $4\Gamma$.  The geometric generic transverse singularity is
a cusp, not a node.  The positive canonical zero divisor on the
normalization pushes forward to $2\Gamma$.
\end{proposition}

\begin{proof}
Write $s=s_4$, $D=D_s$, $\theta=\delta s$ and
$\gamma=\delta^2s$. Corollary~\ref{cor:polar-cubic-covariant}
gives actual sections
\[
 u\in R_{10},\qquad v=-Du/10\in R_{15},\qquad
 p\in R_{12},\qquad Dv=-15u^2+s^2p.
\]
In a local canonical frame the normalized formula is
\[
 u=\theta^2-\frac89s\gamma+\frac19s^2\kappa
\]
with $\kappa$ regular.  Applying the actual weighted bracket
and expanding before setting $s=\theta=0$ gives
\[
 p\equiv\frac{176}{27}\gamma^2\pmod{(s,\theta)}.
\tag{**}\label{eq:c2349-square-correction}
\]
More generally, direct weighted differentiation gives the coefficient
$3m^2(2m+3)/(2m+1)^2$ for an anchor of degree $m$.

The preceding proposition proved that $\gamma$ vanishes on
$\Gamma$.  Thus $p|_\Gamma=0$.  The actual space $R_{12}$
has basis
\[
 s^3,\quad D_s(a_7),\quad q_{12}.
\]
Indeed the first vanishes on $T$, whereas the other two restrict
to independent members $-7\theta a_7$ and $q_{12}|_T$.
The first two vanish on $\Gamma$ and $q_{12}|_\Gamma\ne0$.
Consequently $p$ is a constant linear combination of the first
two.  In every Gorenstein chart both belong to $(s,\theta)$,
since $D_s(a_7)=4s\delta a_7-7a_7\theta$.
Equation~\eqref{eq:c2349-square-correction} therefore implies
\[
 \gamma^2\in(s,\theta).
\]

Here is the local polar-colength argument, with its hypotheses
visible.  At the generic singular curve take the normal-plane
regular local ring and put $L=\ell(R/(s,\theta))$.
If $\theta$ belongs to the square of the maximal ideal then
$L\ge4$.  Otherwise use $y=\theta$ as a parameter.
The displayed ideal containment gives $\gamma(0)=0$.
As $\delta$ is transverse to $\Gamma$, a complementary
parameter $x$ has $\delta x(0)\ne0$.
The identity $\delta s=y$ forces the quadratic initial form
of $s$ to be $cxy+dy^2$ with $c\ne0$.
Writing $s(x,0)=bx^L+\cdots$ and differentiating along $y=0$
gives
\[
 \gamma(x,0)=-\frac{Lb\,\delta x(0)}{c}x^{L-2}
                 +\text{higher terms}.
\]
Since $\gamma^2\in(s,y)$, one has $2L-4\ge L$,
and again $L\ge4$.  Total degree twenty and $w=5$ now
force $L=4$ and leave no further positive polar curve.

The transverse singularity has $\delta=1$, so over an
algebraic closure it is a node or a cusp.  A node with the
already proved third-order contact has polar colength three:
write its two smooth transverse branches as $f_1f_2$.
Their orders on a $\delta$-integral arc are one and two;
the polar intersection is
\[
 2I(f_1,f_2)+I(f_1,\delta f_1)+I(f_2,\delta f_2)
 =2+0+1=3.
\]
This contradicts $L=4$.  Hence the generic singularity is
cuspidal.  Its conductor has coefficient two on the unique
normalization branch.  The actual polar-zero identity is the
canonical zero divisor plus this conductor, so the remaining
positive canonical-zero coefficient is $4-2=2$.
\end{proof}

\paragraph{Two blowups and a starred fibre in the nonnormal case}
\begin{proposition}
In the nonnormal quartic branch of C2-02349, its minimal
normalization resolution is obtained from the minimal elliptic K3 by
exactly two successive blowups. The strict conductor has fibre
multiplicity four. The fibre on the minimal K3 is of type
$I_n^*$ or III$^*$. If the difference of the two marked sections
has order four, the fibre is $I_n^*$ with $n$ odd.
\end{proposition}

\begin{proof}
Suppose that the canonical quartic $D$ is nonnormal. Let
$W\to D^\nu$ be the minimal resolution of its normalization and put
$A=K_Z|_W$. The reduction gives a K3 minimal model $S$,
a connected elliptic pencil on $W$, and a unique curve of
nonnormality $\Gamma$ with $420K_Z\Gamma=5$. Its generic transverse
singularity is a cusp with conductor multiplicity two. The full
positive polar cycle of the weight-five adjunction section $\theta$
is $4\Gamma$. The positive canonical zero cycle on the normalization
is $2\Gamma$. Write $R$ for the strict transform of the unique
normalization branch over $\Gamma$. Thus
\[
 A R=\frac1{84},\qquad
 (K_W)_{A\text{-positive}}=2R.
\]
Here $K_W$ is represented by the nonzero restricted holomorphic
2-form. Its divisor is the exceptional canonical divisor of the
minimalization $\sigma:W\to S$.

The elliptic pencil is a morphism on this minimal resolution by
Lemma~\ref{lem:k3-elliptic-pencil-morphism}.

Every $(-1)$-curve exceptional for $\sigma$ and still of square
$-1$ on $W$ has positive $A$-degree. Indeed an $A$-null curve is
exceptional for $W\to D^\nu$, since $K_Z|_D$ is ample; a rational
$(-1)$-curve exceptional for that morphism contradicts minimality
of the resolution. Consequently $R$ is the only final $(-1)$-curve
in the blowup construction of $W$ from $S$.
The canonical coefficient of a newly created exceptional curve is
one plus the sum of the coefficients of the exceptional curves
through its center. A final exceptional curve of coefficient two
can therefore arise only from a second blowup on a coefficient-one
exceptional curve. Any other blowup branch would have a different
final $(-1)$-curve; any later blowup on this branch would give a
final coefficient at least three. It follows that $\sigma$ consists
of exactly two successive blowups. If $E$ is the first exceptional
curve's strict transform, then
\[
 E^2=-2,\quad R^2=-1,\quad E R=1,\quad
 K_W=E+2R,\quad AE=0. \tag{1}
\]

We identify the normalization at the actual seventh point, where
the canonical contraction is a local isomorphism. On its smooth
index-one cover choose coordinates $(x,y,z)$ of characters
$(1,6,3)$. The local proof gives the weighted initial
quartic equation $c x^2+d z^3$, $cd\ne0$, and the conductor branch
has $y$-order one. Equivariant completion of the square writes
$f=X^2+g(y,z)$. The lifted singular curve is smooth, parameterized
by $y$, so an equivariant translation makes it $X=Z=0$.
Then $g=Z^2a(y)+Z^3b(y,Z)$. Its generic transverse type is cuspidal,
so $a(y)$ vanishes identically, while $b(0,0)\ne0$ by the initial
form. Rescaling $Z$ by an invariant analytic cube root yields
\[
                 f=X^2+Z^3.
\]
The normalization is $X=t^3$, $Z=-t^2$, with characters
$(t,y)=(5,6)$. After changing the generator of $\mu_7$ it is exactly
$\frac17(1,4)$.

Its minimal resolution consists of a $(-4)$-curve $V$ and a
$(-2)$-curve $L$, joined once. The strict conductor meets $V$.
The general moving branch has normalized orders $(t,y)=(1,4)$,
so meets $L$ at a general smooth point. To check these assertions
without a convention ambiguity, the fan from the $t=0$ ray to
$y=0$ has intermediate primitive vectors
\[
 (2/7,1/7),\quad(1/7,4/7),
\]
with self-intersections respectively $-4,-2$.
The first is adjacent to the strict $t=0$ divisor; the second is
the moving-branch ray. The exact $y$-order four follows also from
the varying leading ratio in the actual twelfth pencil: its first
member is $\theta a_7$, and a second member must have a nonzero
$y$ initial term in addition to the $z^2$ polar term.
The unique seventh marking therefore makes $L$ a section, whereas
$V$ is vertical.

Adjunction gives $K_W L=0$ and $K_W V=2$.
The local fan identifies $RV=1$. Thus (1) gives $EV=0$.
The second blowup creates $R$ at a point of the first exceptional
curve; since $RV=1$, its center is also on the strict transform
of $V$. The total transform of that first exceptional curve is
$E+R$, whose intersection with $V$ is one. Hence the image of $V$
after the second blowdown meets the first exceptional curve with
multiplicity one. Both blowups have multiplicity one on $V$.
Its image $V_0$ on the K3 is therefore a smooth rational
$(-2)$-curve. This argument rules out a nodal or cuspidal image
without assuming smoothness of the image in advance.
The section $L$ misses the two blowup centers because
$(E+2R)L=0$. It meets $V_0$, so $V_0$ is a multiplicity-one
component of its elliptic fibre.

Let $p$ be the first blowup center on $V_0$. The second center is
where its strict transform meets the first exceptional curve.
The Kodaira configurations give precisely the following
multiplicities for the second exceptional curve:
\[
\begin{array}{l|c}
 p\text{ on the fibre}&\operatorname{mult}(R)\\ \hline
 \text{smooth point of }V_0&2\\
 \text{ordinary node with another multiplicity-one component}&3\\
 \text{type III tangency of two simple components}&4\\
 \text{type IV triple point of three simple components}&4\\
 V_0\text{ meeting a multiplicity-two component}&4
 \end{array}
\]
For the last case, the first exceptional curve has multiplicity
three and the second has multiplicity four. In type III the
second center is also on the other tangent component, producing
$1+1+2=4$; in type IV the first exceptional multiplicity is three.
The last line occurs in $I_n^*$, III$^*$ and II$^*$; the affine
$E_6$ diagram of IV$^*$ has its multiplicity-one ends meeting
double components as well. We retain it at this numerical stage.
These exhaust the possibilities for a smooth multiplicity-one
component in a Kodaira fibre.

On the other hand, the twelfth-pencil reduction gives
$12A\sim_{\mathbb Q}F+B$, where $AB=0$. The actual member
$\theta a_7$ has multiplicity at least four along $R$, since its
polar factor has positive zero cycle $4\Gamma$ and the restriction
of $a_7$ is regular. The fixed divisor has no $R$ component because
$AR>0$. The pencil is connected, so its moving divisor is one
fibre. Its multiplicity along $R$ is at least four. The table
therefore forces exact multiplicity four, and $a_7$ does not
vanish generically along the conductor.

The type III alternative is also impossible, by a separate
intersection argument. Denote its other component after the
two blowups by $V'$. Its full fibre is
\[
                  F_\infty=V+V'+2E+4R.
\]
Since $AV=AE=0$, $AR=5/420$, and $AF=48/420$,
we obtain $AV'=28/420>0$. The positive polar cycle is
exactly $4R$, and the fixed divisor of the twelfth pencil is
$A$-null. Thus the effective rational divisor obtained by
restricting and pulling back the ambient seventh section $a_7$
has coefficient one on $V'$ and coefficient zero on $R$.
This is the full pullback of an effective $\mathbb Q$-Cartier
divisor, hence all its other coefficients are nonnegative.
It is numerically $7A$. As $R$ is not a component, every other
effective component has nonnegative intersection with it.
Consequently
\[
           \frac1{12}=7AR
              =\operatorname{div}(a_7|_W)R
              \ge V'R=1,
\]
a contradiction. No Cartier assumption at the intersection
point is used; retaining the effective full pullback is essential.

The marked-ratio argument above says that the difference of
the marked sections $P_7-P_5$ has exact order two or four.
In characteristic zero specialization of torsion into the component
group is injective at an additive fibre
\cite[Section 7.8.2, Lemma 7.8]{SchuettShioda}.
The component groups in \cite[Section 7.2, Lemma 7.3]{SchuettShioda}
are trivial for II$^*$, cyclic of order three for IV and IV$^*$,
cyclic of order two for III and III$^*$, and respectively
$(\mathbb Z/2)^2$ or $\mathbb Z/4$ for $I_n^*$ with $n$ even or
odd. Thus IV, IV$^*$ and II$^*$ are impossible. Together with
the preceding type III contradiction, the remaining
fibres are exactly those in the statement; order four forces an
odd $I_n^*$.
\end{proof}

\paragraph{Propagation from the degree-seven section}
\begin{proposition}
The fibre III$^*$ is impossible. In the $I_n^*$ alternative,
$n\ge9$, and the connected exceptional configuration of the
normalization above the index-twelve point is exactly the graph
with a central $(-3)$-curve, two short $(-2)$ arms, and a long
arm of ten $(-2)$-curves. Its determinant is 48 and its
discriminant group is cyclic of order 48. The restricted ambient
canonical class has exact order twelve in that group.
\end{proposition}

\begin{proof}
Let $H_7$ be the effective rational divisor on $W$ obtained by
restricting and pulling back the ambient seventh section $a_7$.
It is numerically $7A$ and has no $R$ component. On the normalized
conductor curve, its entire degree $7AR=1/12$ is exhausted by
the primitive index-twelve marking: multiplication of the proved
fractions $(3/7,7/12)$ by seven gives residues $(0,1/12)$.
At the seventh point $a_7$ is a unit. This follows also from its
generic elliptic-fibre divisor $4P_5/5$, which has no seventh zero.

In either remaining starred fibre the curve $R$ meets exactly two
vertical components: $V$ over the seventh point and $E$, the first
blowup exceptional curve. No horizontal exceptional curve over
the twelfth point can meet a general fibre: the exact generic
divisors of $\theta$ and $a_7$ are exhausted by their fifth and
seventh markings. Their product is one pencil generator, so a
horizontal exceptional curve would supply an additional zero
of that generator on the generic fibre, including in trivial inertia.
Consequently the twelfth marking on $R$ lies on $E$.
The full effective pullback and $RV=RE=1$ give
\[
        \operatorname{coeff}_{V}(H_7)=0,\qquad
        \operatorname{coeff}_{E}(H_7)=\frac1{12}.
\]
There are no further vertical components through $R$ in these
starred Kodaira diagrams.

For every vertical component write
$q_C=12\operatorname{coeff}_C(H_7)$. On the exceptional connected
configuration over the twelfth point, all coefficients are positive.
They satisfy the exact intersection equations $H_7C=0$.
Every component of positive $A$-degree in this fibre has its full
fibre multiplicity in $H_7$: the positive polar cycle is exactly
$4R$ and the pencil's fixed divisor is $A$-null.
Thus a component of fibre multiplicity $m_C$ is either exceptional,
where the intersection equation applies, or has
\[
                         q_C=12m_C. \tag{1}
\]
No upper bound $q_C\le12m_C$ is imposed on exceptional components;
the full null fixed divisor is retained.
All positive components of $H_7$ lie in this fibre. Indeed their
scaled $A$-degree is $48-4\cdot5=28$, equal to the entire
$420A H_7=7\cdot4=28$ budget.
The only horizontal component of $H_7$ is the fifth marked
exceptional section, and it cannot meet the twelfth exceptional
configuration, since their images are different points.
Hence the stated fibre intersection equations have no omitted
positive or horizontal terms.

For III$^*$, start from $q_E=1$. The adjacent double component
has square $-3$, and all other unchanged fibre components have
square $-2$. Along the affine $E_7$ diagram the forced coefficients
are
\[
\begin{array}{c|rrrrrrr}
\text{fibre multiplicity}&2&3&4&2\ (\text{side leaf})&3&2&1\\
q&2&5&8&4&7&6&5.
\end{array}
\]
At each step (1) is impossible, so the relevant component must
remain exceptional. At the last leaf its equation is $2\cdot5=6$,
a contradiction. This excludes III$^*$.
In particular the side leaf cannot be positive: its value 24
would exceed the available central sum $2\cdot8-5=11$.
It is therefore exceptional and has the asserted value four.

For $I_n^*$, call the chain of double components
$W_0,\ldots,W_n$, with $V$ one of the two simple ends attached
to $W_0$. Denote the other simple end there by $B$, and the two
simple ends at $W_n$ by $C,D$. For $n=0$ use the single
four-valent double component. After the two blowups,
$E^2=-2$, $W_0^2=-3$, and the unchanged components have square
$-2$. The equations begin
\[
             q_E=1,\qquad q_{W_0}=2,\qquad q_B=1.
\]
Here $B$ cannot be positive: its value twelve would exceed
$3q_{W_0}-q_E=5$, even before adding the other nonnegative
neighbor. Thus its $(-2)$ equation gives $q_B=1$.
For $n=0$ the three other simple leaves would all have coefficient
one, whereas the central equation requires their sum to be five,
so that case is impossible.
For $n\ge1$, induction along the double chain gives
\[
                     q_{W_j}=2j+2
\]
until the first positive component. By (1), a positive double
component can first occur only at $j=11$.

If the whole double chain remains exceptional, its last equation
gives $q_C+q_D=2n+4$. An exceptional simple end has coefficient
$n+1$, and a positive one has coefficient twelve. Both exceptional
ends give $2n+2$, which is impossible. One positive end forces
$n=9$; two positive ends force $n=10$. Thus the alternatives are
exactly: $n=9$ with one positive far end; $n=10$ with both far
ends positive; or $n\ge11$ with first positive double component
$W_{11}$. Continuing past $W_{11}$ as an exceptional component
cannot lead to any later positive double component and cannot
satisfy either far-end equation, so does not give another case.

In all three alternatives the exceptional connected graph over
the twelfth point is the same weighted graph: central $W_0$ of
square $-3$, short arms $E,B$ of square $-2$, and ten successive
$(-2)$ vertices on the long arm. In the $n=9$ case the last
vertex is the exceptional far simple end; otherwise it is
$W_{10}$. Components outside this connected graph are positive
boundary components, so no further exceptional curve over that
point is omitted.

Let $M$ be the positive intersection matrix, ordering the two
short arms first, then the central vertex, then the long arm.
Direct elimination gives
\[
 \det M=48,\qquad
 \operatorname{SNF}(M)=\operatorname{diag}(1,\ldots,1,48).
\]
The vectors for the three effective seventh-section boundary
patterns have denominators exactly twelve. They define the
local class of $7H$, where $H$ is the ambient canonical
polarization. Its order is exactly twelve, because its
coefficient on $E$ is $1/12$ and twelve times every coefficient
is integral. Since the discriminant group has order 48 and
$\gcd(7,48)=1$, the class of $H$ also has
order twelve.
\end{proof}

\paragraph{The normalization cover rules out the cusp}
\begin{lemma}[A one-dimensional normalization defect bounds embedding dimension]
Let $R$ be a two-dimensional complex analytic hypersurface local
domain with embedding dimension at most three, and let $B$ be
its local normalization. Suppose the nonnormal locus is one smooth
curve $\Gamma$, and the geometric generic transverse singularity
is a cusp of $\delta$-invariant one. Suppose $B$ is Cohen--Macaulay.
Then the embedding dimension of $B$ is at most four.
\end{lemma}

\begin{proof}
Put $Q=B/R$. Both $R$ and $B$ have depth two, so the depth lemma
gives $\operatorname{depth}Q\ge1$. The support of $Q$ is exactly
$\Gamma$; finite support at the origin cannot be an extra submodule.
Choose a regular parameter $x$ on $\Gamma$. Its lift to $R$ is
a parameter on the one-dimensional module $Q$, hence acts
injectively. Thus $Q$ is a finite free $\mathbb C[[x]]$-module.
Its generic rank is the transverse normalization defect, namely
one, so $Q\simeq\mathbb C[[x]]$.
The ideal of $\Gamma$ acts trivially generically on $Q$, because
the normalization quotient of an ordinary cusp has length one
and is annihilated by the transverse maximal ideal. Its action
is therefore zero everywhere, by freeness over $\mathbb C[[x]]$.
Choose $u\in B$ lifting a generator of $Q$. Then $B=R+Ru$ as
an $R$-module. In particular three generators of the maximal
ideal of $R$, together with $u-u(0)$, topologically generate
$B$ as a complex analytic algebra. Its embedding dimension
is at most four.
\end{proof}

\begin{proposition}[Normality of the C2-02349 quartic]
The canonical quartic of a hypothetical realization of
C2-02349 is normal.
\end{proposition}

\begin{proof}
Assume the cuspidal alternative. By the preceding
coefficient-propagation proposition, the normalization germ $N$
over the actual index-twelve point has minimal resolution graph
with central weight $-3$, two short arms of weight $-2$, and
a long arm of ten $(-2)$ curves.

This graph is rational. For example, its fundamental cycle is
reduced. Starting with the central curve and adding the two short
arms and then the successive long-arm vertices is a Laufer
computation in which every added curve has intersection one with
the cycle already constructed. The rationality criterion therefore
applies. Equivalently this is the usual reduced-fundamental-cycle
rational graph with each weight at least its valency.
Hence the local divisor class group is finite and the germ is
$\mathbb Q$-Gorenstein. Solving the adjunction equations gives
the negative discrepancies
\[
 \Delta=
 \left(\frac{11}{24},\frac{11}{24},\frac{11}{12},
 \frac56,\frac34,\frac23,\frac7{12},\frac12,
 \frac5{12},\frac13,\frac14,\frac16,\frac1{12}\right).
\]
All are strictly less than one. The normal surface germ is
klt, and therefore is a quotient surface singularity.

The Hirzebruch--Jung fraction for the central vertex and long
arm is
\[
                    [3,\underbrace{2,\ldots,2}_{10}]
                    =\frac{23}{11}.
\]
The dihedral-quotient classification
\cite[\S\S5,7]{WemyssMcKay} identifies the germ with
$\mathbb C^2/G$, where, up to simultaneous conjugacy,
\[
 r=\operatorname{diag}(\zeta_{22},\zeta_{22}^{-1}),
 \qquad
 s=\zeta_{48}\begin{pmatrix}0&1\\-1&0\end{pmatrix}.
\]
This is $D_{23,11}$ in that source. Its order is 528.
Its abelianization is cyclic of order 48; the determinant
character has order 24, agreeing with the discrepancy calculation.
For the reverse identification, complex quotient surface
singularities are taut: their analytic type is determined by the
weighted resolution graph. One may use
\cite[Section 4.2 and Remark 4.4]{LiedtkeMcKay}, which records
Laufer's tautness theorem and the quotient graph classification.
Thus the classification applies to the actual normalization
germ, rather than merely exhibiting an example with its graph.

The ambient canonical twelve-cover, after restriction and
normalization, is the cover of $N$ belonging to the local class
of $H=K_Z|_D$. The preceding proposition proves that this class
has exact order twelve. As the local class group is cyclic of
order 48, all its order-twelve characters have the same kernel.
In the displayed group that kernel is
\[
             \langle r,s^{12}\rangle,\qquad s^{12}=iI.
\]
It is cyclic of order 44. More explicitly
\[
             r^{17}s^{12}
              =\operatorname{diag}(\zeta_{44},\zeta_{44}^{21})
\]
is a generator. Thus the normalization of the surface on the
ambient index-one cover is exactly
\[
                         N'=\frac1{44}(1,21). \tag{1}
\]
This is a connected cover, since the character has exact order
twelve; no unproved identification with the canonical cover
of $N$ is made.

We check the hypotheses of the lemma on this same ambient cover.
Proposition~\ref{prop:actual-source-transfer46} makes the canonical contraction a local
isomorphism at the ordinary index-twelve source. Hence the
ambient index-one cover is smooth, and the lifted quartic has
a hypersurface local ring $R$ of embedding dimension at most
three. It has only the lifted conductor curve. Indeed on the
canonical model every curve has positive $K_Z$-degree, so the
unique-positive-conductor statement excludes all
other conductor curves. A finite cover cannot introduce a
curve mapping to the origin, and away from the origin the
ambient canonical cover is étale. Thus no exceptional
zero-degree conductor curve has been silently omitted.

The conductor has exactly one primitive twelfth marking,
of full inertia twelve. Its lift has one normalized branch.
The degree-$7$ section has exact numerator order one on that
branch, by the saturated degree $1/12$ budget. The corresponding
ordinary chart character is one, so its numerator has a nonzero
linear coordinate term. The lifted conductor is therefore
smooth on the hypersurface side, with that section as a regular
parameter. Its geometric generic transverse defect is still
the cusp defect one. The normalization (1) is Cohen--Macaulay.
All hypotheses of the lemma hold, giving
$\operatorname{embdim}N'\le4$.

But the invariant ring of $\frac1{44}(1,21)$ has at least five
independent indecomposable positive monomials:
\[
                   u^{44},\quad u^{23}v,\quad
                   u^2v^2,\quad uv^{23},\quad v^{44}.
\]
Each has invariant weight zero modulo 44, and none is a product
of two nonconstant invariant monomials. Their residue classes
in the maximal ideal modulo its square are independent, since
distinct invariant monomials are linearly independent.
Thus $\operatorname{embdim}N'\ge5$, a contradiction.
The cusp alternative is impossible.

The canonical quartic is Cohen--Macaulay. Its only possible
failure of normality had been the specified conductor curve,
so its absence proves normality.
\end{proof}


\label{sec:2349-exclusion46}
\begin{theorem}\label{thm:2349-exclusion46}
There is no canonical threefold with $\chi=2$, $K^3=1/420$ and
terminalization basket
$5(1,2)+(1,3)+4(1,4)+(2,5)+(3,7)+(5,12)$.
\end{theorem}

We now use the normal quartic and its marked elliptic K3 pencil
to determine the complete exceptional canonical cycle.

\paragraph{The rational surface and the twelfth slice}
\begin{proposition}
The normal canonical quartic of C2-02349 has
rational singularities. At its actual ordinary index-twelve point,
the index-one surface germ has an equivariant equation
\[
            yz+x^d=0,\qquad d\equiv4\pmod {12},\quad d\ge4,
\]
in coordinates of characters $(1,11,5)$.
Consequently the surface is klt there. Its quotient germ is
$\frac1{12d}(1,11d-1)$.
\end{proposition}

\begin{proof}
The divisor sequence and Serre duality give
$\chi(\mathcal O_D)=\chi(\mathcal O_Z)+P_5=2$.
Its smooth resolution is birational to the K3 model,
so its Euler characteristic is also two. As $D$ is normal,
the Leray formula expresses their difference as the negative
sum of the lengths of $R^1\pi_*\mathcal O_W$. Every length is
nonnegative, so they all vanish. Thus $D$ has rational
singularities.

At the ordinary index-twelve point use the smooth canonical
cover with coordinate characters $(1,11,5)$ and quartic
numerator $f$ of character four. The surface on the cover is
normal: away from the origin the cover is étale over the normal
surface, and the hypersurface is $S_2$. Thus isolated possible
failure of normality is excluded by Serre's criterion.
The only quadratic monomial of character four is $yz$;
there is no linear monomial. All coordinate characters are odd,
so every allowed monomial has even ordinary degree.

Suppose the coefficient of $yz$ were zero. Then
$m=\operatorname{mult}_0 f\ge4$. The residue form $\omega$ has
character $1+11+5-4\equiv1$, so $y\omega$ is invariant.
The ordinary blowup valuation has discrepancy two and $v(y)=1$.
Lemma~\ref{lem:ambient-residue-bound} gives $m\le2+1=3$,
a contradiction. This applies also to a nonreduced tangent cone.

The $yz$ coefficient is therefore nonzero. The equivariant
Morse lemma in the two nondegenerate variables gives
$f=yz+h(x)$. Its character forces the first exponent of $h$
to be congruent to four modulo twelve. It cannot vanish
identically, since the index-one surface is normal and
irreducible. An invariant analytic unit can be absorbed into
$x$, giving $yz+x^d$ with the stated $d$. This is Du Val
type $A_{d-1}$; its finite quotient is klt.

For completeness, uniformize the Du Val surface by
$y=u^d,z=v^d,x=uv$. A lift of the twelfth action to the
uniformization has weights $a,d-a$ modulo $12d$, where
$a\equiv11\pmod {12}$ and $a$ may be chosen coprime to $12d$.
Normalizing the first weight to one gives second weight
$d a^{-1}-1\equiv11d-1\pmod {12d}$.
This proves the stated cyclic quotient type. For $d=4$ its
minimal resolution is the chain
$[-2,-2,-2,-2,-2,-2,-2,-2,-3,-3]$.
\end{proof}

\paragraph{The seventh slice}
\begin{proposition}
In the normal branch of C2-02349, the quartic is rational by the preceding proposition.  At its actual seventh
source, the index-one covering surface is either $A_1$ or $E_8$.
\end{proposition}
\begin{proof}
The ambient smooth cover has coordinate characters $(1,6,3)$ modulo
seven, and the quartic numerator $f$ has character two.  The initial calculation above gives
\[
                Ax^2+Byz+Cz^3,\qquad
                A\ne0,\quad (B,C)\ne(0,0),             \tag{1}
\]
for filtration weights $(3,4,2)$.  If $B\ne0$, its ordinary
quadratic part $Ax^2+Byz$ is nondegenerate, so the covering surface
is $A_1$.

Suppose $B=0$, hence $C\ne0$. The residue form has character
one, so multiplication by $y$ makes it invariant.

Use the primitive ambient valuation of weights $(6,1,4)$.  Its
discrepancy is ten.  Every character-two monomial has weight
congruent to five modulo seven, and the only monomials of weight
five are $yz$ and $y^5$.  If the $y^5$ coefficient were zero,
the valuation of $f$ would be at least twelve.

Since $v(y)=1$, Lemma~\ref{lem:ambient-residue-bound} instead
requires $v(f)\le10+1=11$, a contradiction.

Thus the coefficient of $y^5$ is nonzero.  For weights
$(15,6,10)$, the only allowed monomial of weight below thirty is
$yz$, whose coefficient is zero.  The weight-thirty part consists
exactly of
\[
                         Ax^2+Cz^3+Dy^5,
                         \qquad ACD\ne0.
\]
All remaining terms have strictly greater weight.  This is a
semiquasihomogeneous $E_8$ singularity.  Explicitly the Jacobian
algebra of $x^2+z^3+y^5$ has basis $y^jz^i$ with
$0\le j\le3$ and $0\le i\le1$, whose largest weight is
twenty-eight; the higher-weight terms are removed by the usual
analytic finite-determinacy argument.  Hence the covering surface
is $E_8$.
\end{proof}

\paragraph{The explicit seventh quotient in the $E_8$ alternative}
\label{sec:r38-c2349-e8-graph}
In the $E_8$ alternative at the seventh
point, the index-one surface is
\[
 x^2+z^3+y^5=0,
 \qquad \mu_7:(x,y,z)\longmapsto(\zeta x,\zeta^6y,\zeta^3z).
\]
Its quotient has the following minimal resolution.  The central
vertex $E_0$ has square $-2$; its three arms, in order away from
$E_0$, are
\[
 [E_1]=[-2],\qquad [E_2,E_3]=[-2,-2],\qquad
 [E_4,E_5]=[-2,-3].
\]
In particular its determinant is seven and its negative discrepancy
vector is
\[
 \Delta=\frac17(6,3,4,2,5,4),\qquad -\Delta^2=\frac47.
\]

Here is a direct identification of the quotient, rather than an
inference from its group order.  The star has Seifert data
\[
 (b_0; (\alpha_i,\omega_i))=(2;(2,1),(3,2),(5,3)),
 \qquad e=-\frac7{30}.
\]
The weighted-homogeneous reconstruction and the universal-abelian-cover
formula of \cite[\S2.1, Theorem~2.2.1]{NemethiOkuma}, including its displayed discriminant action,
give cover $x^2+z^3+y^5=0$.  If $M=-(E_iE_j)$, then
\[
 M^{-1}=\frac17
 \begin{pmatrix}
30&15&20&10&18&6\\
15&11&10&5&9&3\\
20&10&18&9&12&4\\
10&5&9&8&6&2\\
18&9&12&6&15&5\\
6&3&4&2&5&4
 \end{pmatrix}.
\]
The endpoint classes $E_1^*,E_3^*,E_5^*$ correspond to $x,z,y$.
The generator $[E_1^*]$ acts on $(x,y,z)$ with characters $(3,4,2)$.
Its fifth power has characters $(1,6,3)$, exactly the actual action.
The three marked points on the central projective line have no
modulus, so this construction identifies the actual quotient germ.
All six vertices have square at most $-2$, hence this is its
minimal resolution.  Solving $M\Delta=(0,0,0,0,0,1)^t$
gives the displayed discrepancies.

The coordinate orders along $E_5$ are $(3,4,2)/7$, obtained from
the three endpoint columns of $M^{-1}$.  They equal the proved arc orders of the general marked fibre.  Thus $E_5$ is the unique
horizontal seventh section; the other five vertices are vertical.

\paragraph{The cyclic seventh quotient and the moving marked section}
\label{sec:r38-c2349-a1-graph}
In the $A_1$ alternative, equivariant Morse coordinates put the
covering equation in the form $x^2+yz=0$, with characters
$(1,6,3)$ modulo seven.  Its smooth double cover is
\[
 x=uv,\qquad y=u^2,\qquad z=-v^2.
\]
The action $\mu_{14}(1,11)$ on $(u,v)$ induces characters
$(6,1,4)$ modulo seven on $(x,y,z)$; replacing its residual
$\mu_7$ generator by its sixth power gives exactly $(1,6,3)$.
The actual quotient is therefore $\frac1{14}(1,11)$.
Its fan lattice is
$\mathbb Z^2+\mathbb Z(1,11)/14$, and, ordered from $(0,1)$
to $(1,0)$, its interior primitive rays are
\[
 \frac1{14}(1,11),\quad\frac1{14}(2,8),\quad
 \frac1{14}(3,5),\quad\frac1{14}(4,2),\quad
 \frac1{14}(9,1).
\]
The relations $v_{i-1}+v_{i+1}=b_i v_i$ give
$(b_1,\ldots,b_5)=(2,2,2,3,2)$.  On the fourth ray the
coordinate orders are
\[
 (v(x),v(y),v(z))=(3,4,2)/7,
\]
so it is the distinguished $(-3)$-curve $H$.
On every other ray $v(x)>3/7$.  Hence an arc with the proved
moving-fibre orders cannot meet another exceptional component
or an intersection of exceptional components.  It meets the
smooth open part of $H$ with multiplicity one.
The same argument applies to the $E_8$ graph: the $x$-orders on
$E_0,\ldots,E_5$ are $(15,11,10,5,9,3)/7$, so $E_5$ is again
the unique possible component and the intersection is transverse.

Finally Lemma~\ref{lem:k3-elliptic-pencil-morphism} makes the
restricted elliptic pencil a morphism on its minimal resolution
$W$. Thus different
fibres are disjoint on $W$.  Their transverse intersections
with $H$ consequently vary along $H$, and the unique seventh
mark gives $FH=1$.  This proves that $H$ is a section and that
its coefficient in the canonical zero divisor is zero; it is
not inferred from the valuation vector alone.

\paragraph{A terminal-component test for effective canonical cycles}
\label{sec:r38-terminal-component-test}
The following elementary observation makes the local graph calculations
independent of a chosen order of blow-downs.

\begin{lemma}[Terminal-component inequality]
Let $\pi:W\to D$ be the minimal resolution of a normal surface,
let $A$ be the pullback of an ample divisor on $D$, and suppose the
smooth minimal model $S$ of $W$ is a K3 surface.  Write $K_W$ for the
zero divisor of the pullback of a nonzero two-form on $S$.
Let $G=\bigcup E_i$ be one complete exceptional fibre of $\pi$.
Assume that every $A$-positive component of $K_W$ meets $G$.
Write
\[
 K_W=\sum_i k_iE_i+K',\qquad b_i=K'E_i.
\]
Then $k_i,b_i$ are nonnegative integers.  If $K_W\ne0$, some
$A$-positive component $C$ has $C^2=-1$.  If its coefficient is $a$
and $m_i=CE_i>0$, then
\[
 a=1+(K_W-aC)C\ge1+k_i m_i,
 \qquad b_i\ge a m_i\ge1+k_i.               \tag{TC}
\]
\end{lemma}
\begin{proof}
The map $W\to S$ is a sequence of point blow-ups.  Its canonical
zero divisor has positive integer coefficient on every exceptional
component and no other components.  At least one last exceptional
curve of this sequence has square $-1$.  Such a curve cannot be
$A$-null: otherwise it would be a $(-1)$-curve contracted by the
minimal resolution $\pi$.  This proves the existence assertion.
Every component outside the complete fibre $G$ which meets it is
$A$-positive, since a curve of $A$-degree zero is $\pi$-exceptional
and distinct fibres of $\pi$ are disjoint.  Adjunction on the
smooth rational $(-1)$-curve $C$ gives $K_WC=-1$ and proves (TC).
\end{proof}

In the present case, every positive component of the canonical zero
divisor meets the seventh exceptional fibre. Indeed the full polar has
scaled degree twenty, and its complete marking enumeration gives only
$(5,5,5,5)$, $(5,5,10)$, $(5,15)$, $(10,10)$, $(20)$; every
component in these partitions has a seventh marking. This concerns
all positive components, without declaring arbitrary invisible curves
polarization-null.

The unique seventh horizontal component $H$ has square $-3$ and
$K_WH=1$.  Since the canonical zero divisor is vertical for the
elliptic pencil, its coefficient at $H$ is zero.  Every other
vertex of the seventh graph has square $-2$, hence canonical
intersection zero.  If $M=-(E_iE_j)$ and $e_H$ is the unit vector
of $H$, the complete nonnegative integral system is
\[
 k_H=0,\qquad b=Mk+e_H\ge0.                  \tag{CY}
\]
Here the general marked fibre meets the smooth open part of $H$
transversely: its restricted valuation is the primitive divisorial
valuation $(3,4,2)/7$, and $H$ is its degree-one horizontal component.
This is the geometric identification needed to impose $k_H=0$.

\paragraph{The $A_1$ index-one alternative.}
Number the chain by $0,1,2,3,4$, with weights $(2,2,2,3,2)$ and
$H=E_3$.  The complete solutions of (CY) are
\[
\begin{array}{c|c|c}
k&b&\text{possible terminal attachment}\\\hline
(0,0,0,0,0)&(0,0,0,1,0)&E_3\\
(0,0,0,0,1)&(0,0,0,0,2)&E_4\\
(1,1,1,0,0)&(1,0,1,0,0)&\text{none}\\
(1,2,1,0,0)&(0,2,0,0,0)&\text{none}\\
(2,2,1,0,0)&(2,1,0,0,0)&\text{none}\\
(3,2,1,0,0)&(4,0,0,0,0)&E_0.
\end{array}                                                    \tag{A1}
\]
This is exhaustive without any numerical cutoff: the equation at
$H$ gives $k_2+k_4\le1$, the equation at $E_2$ gives
$k_1\le2k_2$, and that at $E_1$ gives $k_0\le2k_1-k_2\le3$.
These bounds leave exactly the displayed solutions.  The middle
three have no vertex satisfying $b_i\ge1+k_i$, and are excluded
by (TC).

In the first, second and last rows, any last exceptional curve
must use the unique indicated vertex.  Inequality (TC) and the
listed value of $b_i$ force respectively
\[
 (a,\;CE_i)=(1,1),\quad(2,1),\quad(4,1).
\]
It exhausts all the $b$-contributions.  Every positive canonical
component meets this graph, so $C$ is the only positive component.
Moreover equality holds in the adjunction inequality: its
intersection with $K_W-aC$ is exactly $a-1$, already exhausted
by the displayed seventh null chain.  Thus $C$ meets no further
component of the canonical zero divisor.

For clarity this also excludes hidden null canonical clusters.
The displayed seventh null components have their complete
intersections with the outside canonical divisor recorded by $b$;
no other null curve can meet this complete exceptional fibre.
Any remaining component of $K_W$ is therefore in a separate
connected component of its support.  Each such connected
exceptional cluster of $W\to S$ contains a last $(-1)$-curve,
which is positive for $A$ by the lemma.  There is no second
positive component.  No further canonical component exists.
Consequently
\[
 K_W^2=-a,\qquad N:=-K_W^2=a\in\{1,2,4\}.               \tag{N}
\]
Here the equality also follows directly by multiplying the
canonical-intersection equations by their coefficients: each
listed null component is a $(-2)$-curve with $K_WE_i=0$, while
$K_WC=-1$.

\paragraph{The $E_8$ index-one alternative.}
Use the graph and indexing in Section~\ref{sec:r38-c2349-e8-graph},
with $H=E_5$.  The complete solutions of (CY) are
\[
\begin{array}{c|c|c}
k&b&\text{possible terminal attachment}\\\hline
(0,0,0,0,0,0)&(0,0,0,0,0,1)&E_5\\
(2,1,2,1,1,0)&(0,0,1,0,0,0)&\text{none}\\
(2,1,2,2,1,0)&(0,0,0,2,0,0)&\text{none}.
\end{array}                                                    \tag{E8}
\]
Indeed $k_4\le1$, $k_0\le2k_4$, $k_1+k_2\le2k_0-k_4\le3$,
and $k_3\le2k_2-k_0\le6$.  Substitution in these bounded
inequalities gives exactly the table.  The last two rows fail
(TC).  The first forces one positive coefficient-one curve and
no null canonical component, by the same saturation and terminal
cluster argument as above.  It follows that $N=1$.

\paragraph{The twelfth point forces $a=1$.}
In the $A_1$ case, the preceding argument has shown the full
canonical divisor to consist of $aC$ and the listed seventh
$(-2)$-curves.  The twelfth point has a minimal-resolution
$(-3)$-curve $Z$, as the following continued-fraction calculation
shows uniformly.  Since distinct exceptional fibres are disjoint,
\[
 1=K_WZ=a(CZ).
\]
Both factors on the right are nonnegative integers, so $a=1$
and hence $N=1$.  The $E_8$ case already has $N=1$.

\paragraph{The entire twelfth continued-fraction family}
The analytic form is
$\frac1{12d}(1,11d-1)$ with $d=4+12s$, $s\ge0$.
With $[b_1,\ldots,b_t]$ denoting a negative continued fraction,
\[
\frac{12d}{11d-1}=
\begin{cases}
[2^8,3,3],&s=0,\\
[2^{10},5,2,3],&s=1,\\
[2^{10},3,2^{s-2},4,2,3],&s\ge2.
\end{cases}                                                    \tag{HJ}
\]
These are infinite-family identities, not a finite search.
For $M(b)=\left(\begin{smallmatrix}b&-1\\1&0\end{smallmatrix}\right)$,
\[
 M(2)^n=\begin{pmatrix}n+1&-n\\n&1-n\end{pmatrix},
\]
and the first column of
$M(2)^{10}M(3)M(2)^{s-2}M(4)M(2)M(3)$ is exactly
$(48+144s,43+132s)^t$.  The two exceptional values $s=0,1$
are direct evaluations of their displayed chains.

The negative discrepancy vectors, in the chain order of (HJ), are
\[
\Delta_{12}=
\begin{cases}
(1,2,\ldots,9,7)/12,&s=0,\\
(1,2,\ldots,10,\underbrace{11,\ldots,11}_{s},9,7)/12,&s\ge1.
\end{cases}
\]
Multiplication by the positive intersection matrix gives
$(b_i-2)_i$: interior repeated entries cancel, and the four
junction equalities are immediate.  Thus
\[
 -\Delta_{12}^2=
 \begin{cases}4/3,&s=0,\\10/3,&s\ge1.\end{cases}
\]
Every chain contains a $(-3)$-curve.  The seventh discrepancy
contribution is $4/7$ in both alternatives.  Every other
exceptional discrepancy divisor has nonnegative negative
self-intersection, by negative definiteness.  Since
$K_D^2=25H^2=5/21$, the resolution formula gives
\[
 N=-K_W^2
   =\sum_P(-\Delta_P^2)-K_D^2
   \ge\frac47+\frac43-\frac5{21}
   =\frac53.
\]
This contradicts $N=1$ and excludes C2-02349 and proves the theorem.


\subsection{Exclusion of C2-02353: the quintic surface configuration}

\label{sec:2353-preparation46}
We consider the basket
\[
 B=3(1,2)+4(1,3)+(1,5)+2(2,5)+(3,7)+(5,12),
 \qquad I=420,\quad v=1,\quad\chi=2.
\]
Its numerical label is C2-02353. The source-transfer theorem gives
canonical indices in
\[\{1,2,3,4,5,6,7,12\}.\]
It preserves the unique
ordinary sources
\[Q_7=\tfrac17(1,6,3),\qquad Q_{12}=\tfrac1{12}(1,11,5).\] No assertion that the several
index-five constituents remain at distinct canonical points will
be used. Riemann--Roch gives
\[
 \begin{gathered}
 P_1=P_2=P_3=P_4=P_6=P_7=P_8=P_{11}=P_{13}=0,\\
 P_5=P_9=P_{10}=1,\qquad P_{12}=P_{14}=P_{17}=2,
 \qquad P_{18}=3.
 \end{gathered}
\]

\begin{lemma}\label{lem:2353-low46}
Every actual torsion class $\tau$ Cartier in codimension two is
trivial. Every visible prime has scaled degree at least four.
The full quintic member on the canonical model is either a reduced
prime, or $T+E$, with $T$ a reduced visible prime of terminal state
$(n,M)=(4,2)$ and Hodge mass five, and $E$ a reduced invisible
prime of scaled degree one.
\end{lemma}
\begin{proof}
The exact test of Section~\ref{subsec:late-torsion46} leaves three
necessary profiles. For the canonical profile $C=0$ and $H_m=P_m$.
The other two are inverse, have $C=2$, and their low dimensions are
\[
\begin{array}{c|rrrrrrr}
 &H_2&H_3&H_4&H_5&H_6&H_9&H_{10}\\\hline
 \mathrm{A}&0&0&0&1&0&0&1\\
 \mathrm{B}&0&0&1&0&0&1&1.
\end{array}
\]
In particular $H_2=H_3=0$ for every actual class, and squaring
also gives $H_1=0$. The complete visible-prime states through
degree five are
\[
 (1,1),(2,4),(3,9),(4,2),(4,16),(5,-3),(5,11),(5,25).
\]
The first three are Hodge equalities and descend to the excluded
effective translates of degrees one, two and three. Thus the
visible degree floor is four.

Since $P_6=P_7=0$, Lemma~\ref{lem:adjacent-gaps46} gives a
unique reduced visible prime $T$ with $\lambda(T)=5$. The Hodge state
$(4,16)$ has mass four and cannot occur. Either $n(T)=5$, exhausting
the member on $Z$, or $(n(T),M(T))=(4,2)$. In the latter case the
remaining positive degree is one; integrality of positive scaled
degrees makes its divisor one reduced invisible prime.

Suppose now $H_4(\tau)>0$. Its effective divisor has total Hodge
mass four. The visible degree floor makes it a reduced prime $A$
on $Z$, with no additional component. It differs from $T$, whose
mass is five, and from every invisible component of the quintic.
Its section and $a_5$ are therefore divisorially coprime. The actual
two-section Koszul inequality of Lemma~\ref{lem:two-twist-koszul34},
with degrees $(4,\tau),(5,0)$ and
output twists $0,2\tau$, gives
\begin{align*}
 h^2(\mathcal O_Z(\tau))
 &\le P_{10}-H_6(-\tau)-P_5=0,\\
 h^2(\mathcal O_Z(-\tau))
 &\le H_{10}(2\tau)-H_6(\tau)-H_5(2\tau)\le1.
\end{align*}
The last upper bound uses the maximum of
$H_{10}(\rho)-H_5(\rho)$ over all three profiles, and is valid
for the actual class $\rho=2\tau$.

Here $\tau$ is nontrivial because $P_4=0$. Its $H^0$ and $H^3$
in degrees zero vanish, the latter by $H_1(\pm\tau)=0$.
The paired Euler identity then requires
\[
 h^2(\mathcal O_Z(\tau))+h^2(\mathcal O_Z(-\tau))
 \ge4-C=2,
\]
a contradiction. Thus neither of the inverse noncanonical
profiles can occur, and every actual torsion class has the
canonical profile.

If a nontrivial class existed, a multiple would have prime order
$p$. Its connected quasi-\'etale cyclic cover has
$V'=p/420$, $\chi'=2p$, $p_g'=P_2'=0$ and
$h^{2,0}\ge2p-1$. For $p\ge3$, the weak uniform bound
$V'\ge(h^{2,0}-3)/60$ proved from
\cite[Section~12, equation~(12.5)]{LBVolume} gives $13p\le28$,
which is impossible. For $p=2$, the same source's
$P_2'=0,h^{2,0}\ge3$ consequence gives $V'\ge1/14>1/210$.
This proves torsion triviality as well.
\end{proof}
\paragraph{Restriction to a strict visible component}
\label{sec:c2353-strict-bracket-restriction}

\begin{proposition}
\label{prop:c2353-restriction-injective}
Suppose a canonical realization of C2-02353
has the strict quintic alternative
\[
 \operatorname{div}(a_5)=T+E,
 \qquad (420K^2T,420KT^2)=(4,2),
\]
with the visible and invisible components described in the
Lemma~\ref{lem:2353-low46}.  Its nonzero section $u_9$ does not
vanish along $T$.  Restriction of the complete $12K_Z$
system to $T$ is injective and has dimension two.
\end{proposition}

\begin{proof}
If $l=\operatorname{ord}_T(u_9)\ge1$, then $5l-9\ne0$.
Lemma~\ref{lem:visible-invisible-bracket-division} gives
$D_{a_5}(u_9)/a_5\in R_{10}$.
Since $P_{10}=1$ and $a_5^2\ne0$, it equals $c a_5^2$,
with $c\ne0$. Thus $D_{a_5}(u_9)=c a_5^3$ has visible order
both $l$ and $3$. But ampleness and positive canonical degree
give $9\ge4l$, which rules out $l=3$.

If a nonzero $s\in R_{12}$ vanishes along $T$ to order
$l\ge1$, then $5l-12\ne0$. The same lemma produces a nonzero
$D_{a_5}(s)/a_5\in R_{13}$, contrary to $P_{13}=0$.
Consequently no nonzero section of $R_{12}$ vanishes on $T$.
Since $P_{12}=2$, restriction has dimension two.
\end{proof}

The complete adjunction class $\theta_6$ of the quintic is
nonzero on its sole visible component $T$.  Its square lifts
to an actual section of $R_{12}$ by the complete restriction sequence.  Thus the restricted pencil can be
written using $\theta_6^2$ and an independent twelfth section;
$u_9|_T$ is also an actual nonzero section.

\paragraph{The strict quartic component in C2-02353 is impossible}
\label{sec:c2353-strict-marked-exclusion}

\begin{proposition}
\label{prop:c2353-quintic-prime}
In a hypothetical canonical realization of
C2-02353, the unique fifth-canonical member
is reduced and prime on the canonical model.
\end{proposition}

\begin{proof}
It suffices to exclude the strict alternative in the Lemma~\ref{lem:2353-low46}.  Work on a terminal
$\mathbb Q$-factorial minimal model $Y\to Z$, with
$K=K_Y$ and $I=420$.  Write the \emph{full} divisor of its
actual fifth section as
\[
 \operatorname{div}_Y(a_5)=T+E\sim5K,
 \qquad IK^2T=4,\qquad IKT^2=2.
\]
Here $T$ is the strict visible prime, while $E$ is the full
effective invisible remainder, including any exceptional
components.  Its pushforward to $Z$ is the invisible
reduced degree-one prime; the intersections and local
characters below are all computed on $Y$.

By Proposition~\ref{prop:c2353-restriction-injective},
$u_9|_T\ne0$ and the restriction of $R_{12}$ has dimension
two.  The normalized adjunction square is an actual
twelfth section restricting to $\theta_6^2$, where
$\theta_6\ne0$ is the coherent polar section of the ambient
two-form on $T$.  Choose an independent twelfth section
$w_{12}$.  Resolve the resulting pencil on a resolution of
the normalization of $T$, and let $C$ be a normalized general
connected moving curve.  The surface carries a nonzero
two-form, so $g(C)\ge1$.  Put $\ell=IK\cdot C$ and let
$d$ be the Stein degree of the restricted pencil.  Then
\[
 0<d\ell\le12IK^2T=48.
\]

Let $(r,k)$ be the primitive canonical markings of $C$.
The possible nontrivial inertia orders divide the actual
basket indices; thus they belong to
$\{2,3,4,5,6,7,12\}$.  The canonical degree satisfies
\[
 \frac\ell{420}\equiv\sum\frac{k}{r}\pmod{\mathbb Z}.
\]
The actual sections $\theta_6$ and $u_9|_T$ have effective
rational divisors on $C$ of degrees $6\ell/420$ and
$9\ell/420$, respectively.

There is a positive polar cost at every marked branch.
If $6k\equiv0\pmod r$ and $\theta_6$ were a unit there,
then $\theta_6^2$ would be a unit.  The possible source
indices of such a branch divide twelve, so $12K$ is Cartier
at its actual terminal source.  The pencil ratio would then
be regular with a fixed value at that point, and a general
member could not pass through it.  This contradicts its
being an actual marked branch of the general curve.
Consequently, with residues chosen in $\{0,\ldots,r-1\}$,
\[
 \sum\frac{\overline{6k}^{\,+}}r\le\frac{6\ell}{420},
 \qquad
 \sum\frac{\overline{9k}}r\le\frac{9\ell}{420},
\]
where $\overline{6k}^{\,+}=\overline{6k}$ when nonzero
and equals $r$ otherwise.  In particular inertia orders
two, three and six are impossible, because their cost is
at least one whereas $6\ell/420\le24/35<1$.

For a general full member $F\in|12K|$, the curve $C$ is a
log canonical centre of the effective boundary $T+F$ at
its generic point: there the surfaces are smooth and meet
simply.  Applying Lemma~\ref{lem:full-adjunction46} to this boundary gives
\[
 2g(C)-2+\sum(1-1/r)
 \le(K+T+F)C=18KC-E C.
\tag{*}\label{eq:c2353-extra-invisible-cost}
\]
The general moving curve is not contained in $E$.
Hence $E|_C$ is effective.  Notice that the boundary used
here is $T+F$, not $T+E+F$; the extra term is supplied by
the actual equality $T+E\sim5K$.

There is at least one marking, since otherwise
$K C=\ell/420\in\mathbb Z$, impossible for
$1\le\ell\le48$.  The inequality shows $g(C)=1$:
if $g\ge2$, its left side is at least $2+1/2$, while
$18\ell/420\le72/35<5/2$.
Discarding the nonnegative term $E C$ for the initial
finite enumeration, the exact degree, orbifold and two
section budgets leave precisely
\[
\begin{array}{c|c|c|c}
 \ell&\text{markings}&\deg\operatorname{div}(\theta_6)
                    &\deg\operatorname{div}(u_9)\\ \hline
 35&(1/12)&1/2&3/4\\
 45&(1/4,6/7)&9/14&27/28\\
 48&(2/5,5/7)&24/35&36/35
\end{array}
\]
Both section budgets are saturated in every row, and
$d=1$ since $\ell>24$.

For $\ell=35$, the index-twelve source is the unique
ordinary cyclic quotient with basket entry $(5,12)$.
The lower coordinate orders on the lifted arc of canonical
mark $1/12$ are $(7,5,11)$.  The exact order of
$\theta_6^2$ is twelve.  The only positive monomial of
canonical degree twelve and lower arc weight at most
twelve is $xy$.  Constant terms vanish for both pencil
generators at this marked base point.  Thus their leading
terms are constant multiples of $xy$, with the first
multiple nonzero.  The pencil parameter would be their
fixed leading-coefficient ratio, contradicting a general
pencil member.  This is the same ordinary-source argument
used for the marked quartic case; the unique index-twelve
basket entry has axial weight one, so its index-one cover
is smooth.

For $\ell=45$, denote the marked points by $P_4,P_7$.
Saturation gives the exact rational divisors
\[
 \operatorname{div}(\theta_6)=\tfrac12P_4+\tfrac17P_7,
 \qquad
 \operatorname{div}(u_9)=\tfrac14P_4+\tfrac57P_7.
\]
The ratio of two sections of the same rational line class
is therefore a rational function satisfying
\[
 \operatorname{div}\!\left(\frac{u_9^2}{\theta_6^3}\right)
 =P_7-P_4.
\]
The two points are distinct.  This would be a nonconstant
degree-one function on a smooth curve of genus one, which
is impossible.

For $\ell=48$, retain the effective term in
\eqref{eq:c2353-extra-invisible-cost}.  The orbifold cost
is $58/35$, so
\[
 E C\le\frac{18\cdot48}{420}-\frac{58}{35}
       =\frac25.
\]
Only the basket entry $(3,7)$ contributes to the linear
Weil--Riemann--Roch congruence modulo seven.  Its coefficient
is $(420/7)3^2\equiv1\pmod7$, hence
$i_T\equiv IK^2T\equiv4\pmod7$.
The actual local divisor equality $T+E\sim5K$ gives
$i_E\equiv5-i_T\equiv1\pmod7$.
At the mark $5/7$, a nonzero local defining numerator of
$E$ consequently has order at least five.  Therefore
\[
 E C\ge\frac57>\frac25,
\]
the final contradiction.  Every strict-branch marking
configuration is excluded, proving the proposition.
\end{proof}

\paragraph{Six marked endpoints for the remaining prime quintic}
\label{sec:c2353-prime-marked-reduction}

\begin{proposition}
\label{prop:c2353-six-marked-endpoints}
If C2-02353 has a canonical realization, its
prime quintic carries a restricted $12K$ pencil with connected
genus-one general curve.  Its necessary degree and marking
data belong to the following list:
\[
\begin{array}{c|l}
 \ell=420K\cdot C&\text{primitive canonical markings}\\ \hline
 48&(2/5,5/7)\\
 55&(5/7,5/12)\\
 60&(1/7),\ (3/7,5/7),\ (4/7,4/7),\ (5/7,5/7,5/7)
\end{array}
\]
\end{proposition}

\begin{proof}
The preceding proposition makes the unique quintic prime
on the canonical model, and Lemma~\ref{lem:2353-low46} removes all actual
codimension-two-Cartier torsion.  On a terminal minimal model
write its full divisor as
\[
 \operatorname{div}(a_5)=T+E\sim5K,
 \qquad IK^2T=5,
\]
where $E$ is effective and canonical-exceptional.
Its components are invisible.  We do not assume they are
all point-centred, and we do not replace the terminal square
$IKT^2$ by $25$.

Normal division on the canonical model gives injectivity
of $R_{12}$ on the prime member because $P_7=0$.
Likewise $u_9|_T\ne0$ because $P_4=0$.
The complete adjunction class $\theta_6$ is nonzero on
$T$, and its square lifts to an actual twelfth section.
Thus $\theta_6^2$ and an independent $w_{12}$ give a
nonconstant pencil on a resolution of the normalization.
For its normalized general connected curve $C$ and Stein
degree $d$,
\[
 0<d\ell\le12IK^2T=60.
\]
The same terminal-model marked adjunction and positive
polar-cost argument as in the strict case give
\begin{align*}
 \frac\ell{420}&\equiv\sum\frac{k}{r}\pmod{\mathbb Z},\\
 2g(C)-2+\sum(1-1/r)&\le\frac{18\ell}{420},\\
 \sum\frac{\overline{6k}^{\,+}}r&\le\frac{6\ell}{420},
 \qquad
 \sum\frac{\overline{9k}}r\le\frac{9\ell}{420}.
\end{align*}
Here the allowed inertia orders divide the basket indices,
and $\overline{6k}^{\,+}$ is the positive residue convention
used above.  The effective exceptional term $EC$ has only
been discarded from an upper bound obtained from the
actual boundary $T+F$.

The polar degree is at most $6/7<1$, so orders two,
three and six are again excluded.  There is at least one
marking, since $0<\ell/420<1$.  Every remaining marking
costs at least $3/4$ in orbifold adjunction.  Therefore
$g(C)\ge2$ would require a left side at least $11/4$,
whereas the right side is at most $18/7<11/4$.
Consequently $g(C)=1$.

The exact finite enumeration gives eleven rows:
\[
\begin{array}{c|l}
35&(1/12)\\
45&(1/4,6/7)\\
48&(2/5,5/7)\\
49&(1/5,11/12)\\
55&(5/7,5/12)\\
60&(1/7),\ (2/7,6/7),\ (3/7,5/7),\ (4/7,4/7),\\
  &(4/7,5/7,6/7),\ (5/7,5/7,5/7).
\end{array}
\]
The degree-$35$ row is excluded by the same unique
index-twelve leading monomial $xy$ as before.
In every row the $\theta_6$ budget is saturated.
The $u_9$ budget is also saturated except for
$\ell=60$ with markings $(1/7)$ or $(4/7,4/7)$,
where one effective integral divisor of degree one is
left over; it is not assumed to lie away from the markings.

For the four configurations
\[
 (45;1/4,6/7),\quad(49;1/5,11/12),\quad
 (60;2/7,6/7),\quad(60;4/7,5/7,6/7),
\]
the saturated rational divisors make
$u_9^2/\theta_6^3$ a function with one simple zero and
one simple pole at distinct marked points.  Indeed its
coefficient at a mark is
\[
 \frac{2\overline{9k}-3\overline{6k}}r.
\]
The four coefficient lists are respectively
$(-1,1)$, $(1,-1)$, $(-1,1)$ and $(-1,0,1)$.
Such a degree-one function cannot exist on a genus-one
curve.  Removing these rows gives exactly the six listed
in the proposition.  Every surviving degree is at least
$48>60/2$, so $d=1$.

\end{proof}

\paragraph{The eighteenth system leaves two K3 endpoints}
\label{sec:c2353-eighteenth-fibre-rank}

The following refinement uses the dimension of a complete restricted
linear system, rather than just the existence of one section.

\begin{lemma}[A rank-one restriction and the moving degree]
\label{lem:elliptic-residue-rank-one34}
Let $A$ be a nef and big $\mathbb Q$-divisor on a smooth projective
surface resolving a prime canonical member, and let
$f:S\to\mathbb P^1$ have connected genus-one general fibre $C$.
Write $A C=\ell/I$.  Suppose a vector space $V$ of dimension $h$
of actual degree-$m$ canonical sections restricts injectively to
the prime member.  Suppose further that its restrictions to the
geometric generic fibre have the canonical fractional lower orders
$\{mk_x/r_x\}$ at the marked points and ordinary nonnegative orders
elsewhere.  Put
\[
 b_m=\frac{m\ell}{I}-\sum_x\left\{\frac{mk_x}{r_x}\right\}.
\]
If $0\le b_m\le1$, then
\[
 (h-1)\ell\le I m A^2.
\]
The statement does not require normality of the original prime
surface or a section of the genus-one fibration.
\end{lemma}

\begin{proof}
Pass to a common resolution for the surface and the two linear
systems.  The fractional-divisor construction of the marked-curve adjunction interface applies to every actual section
of the same degree, with the same fractional character at a mark.
Choose one nonzero restricted section $s_0$.  Its rational zero
divisor on the generic fibre has the form
\[
 \operatorname{div}(s_0)=
 \sum_x\left\{\frac{mk_x}{r_x}\right\}P_x+Z_0,
 \qquad Z_0\ge0,\quad\deg Z_0=b_m,
\]
where $Z_0$ is integral.  For every other restricted section $s$,
the ratio $s/s_0$ is an ordinary rational function on $C$ and its
pole divisor is bounded by $Z_0$.  Thus all restrictions inject
into $H^0(C,\mathcal O_C(Z_0))$.  On a smooth genus-one curve an
effective divisor of degree zero or one has at most one independent
section, including over the algebraic closure of the base field.
Consequently the ratios of all sections of $V$ are constant on
the connected generic fibre.  They belong to $\mathbb C(\mathbb P^1)$.

The induced $h$-dimensional rational linear series on
$\mathbb P^1$ has moving degree $e\ge h-1$.
On the common resolution its moving divisor is $eC$, and the
actual degree-$m$ section divisors give an effective
$\mathbb Q$-divisor $mA-eC$.  Intersecting with the nef divisor $A$
gives $e\ell/I\le mA^2$, as asserted.  In particular the proof
uses the full effective remainder; it does not replace a terminal
prime square by its Hodge upper bound.
\end{proof}

\begin{proposition}[Two remaining marked K3 configurations]
\label{prop:c2353-two-k3-endpoints}
In a hypothetical realization of C2-02353, the
restricted twelfth pencil on the prime quintic has connected
genus-one general curve of degree $\ell=60$ and has one of the
following two necessary marking configurations:
\[
 \begin{array}{c|c|c}
 \ell&\text{primitive canonical markings}&
       \text{a section of the pencil}\\ \hline
 60&(1/7)&\text{the unique marked point}\\
 60&(4/7,4/7)&\text{the residual degree-one divisor of }u_9
 \end{array}
\]
In both cases the smooth minimal normalization surface is K3.

\end{proposition}

\begin{proof}
Use the six configurations of
Proposition~\ref{prop:c2353-six-marked-endpoints}.
For the pullback $A$ of the ambient canonical divisor on the
resolved prime quintic,
\[
 I=420,\qquad A^2=\frac5{420},\qquad
 I\,18A^2=90.
\]
The actual Hilbert function has $P_{18}=3$ and $P_{13}=0$.
Normal division on the canonical threefold therefore makes the
restriction of $R_{18}$ to the prime quintic injective.
The residual integral degrees on the six generic curves are
\[
\begin{array}{c|c|c}
 \ell&\text{markings}&b_{18}\\ \hline
 48&(2/5,5/7)&1\\
 55&(5/7,5/12)&1\\
 60&(1/7)&2\\
 60&(3/7,5/7)&1\\
 60&(4/7,4/7)&2\\
 60&(5/7,5/7,5/7)&0
\end{array}
\]
In the four rows with $b_{18}\le1$, the lemma would give
$2\ell\le90$.  Each has $\ell\ge48$, which is impossible.
The two rows in the statement are exactly the survivors.

In the first row the unique marking is Galois-invariant over the
pencil field, and therefore gives a section.
In the second, let $P,Q$ denote the two geometric marked points.
The exact ninth-section budget says
\[
 \operatorname{div}(u_9|_C)=\frac17(P+Q)+R,
 \qquad R\ge0,\qquad\deg R=1,
\]
with $R$ an integral divisor.  Both the actual zero divisor and
the equal-coefficient fractional marking sum are Galois-invariant
over the pencil's function field.  Hence $R$ is defined over that
field.  An effective integral degree-one divisor is a rational
point and gives a section of the relatively minimal fibration.
This remains true if its support is one of the markings; in that
event the selected marking itself is rational.

Ordinary component trace and visibility give $p_g=1$ for the
smooth normalization surface.  The elliptic fibration has base
$\mathbb P^1$ and a section, hence no multiple fibres.  The
characteristic-zero elliptic-fibration formulas give $q=0$,
$\chi=2$, and trivial canonical bundle on the relatively minimal
surface.  It is K3 in both cases.
\end{proof}

\paragraph{Saturated moving degree removes the fifth-index points}
\label{sec:c2353-fifth-index-avoidance}

\begin{proposition}
\label{prop:c2353-fifth-index-avoidance}
For a hypothetical realization of C2-02353, its
canonical prime quintic $D$ avoids every point of canonical index
five.  On any crepant terminal model its strict transform $T$
avoids every actual index-five terminal point.  Consequently
\[
 (IK^2T,IKT^2)=(5,25),\qquad I=420,
\]
and every component of the full remainder
$\operatorname{div}(a_5)_Y-T$ is contracted to a point of the
canonical model.
\end{proposition}

\begin{proof}
Write $H=K_Z|_D$, interpreted on the finite normalization when
necessary, and let $A$ be its pullback to a common resolution.
The two remaining marked configurations give a connected
twelfth pencil with
\[
 I A^2=5,\qquad I A C=60.
\]
The moving part therefore consumes all the intersection
$I A(12A)=60$.  Every fixed divisor on the resolution has
$A$-degree zero.  As $K_Z$ is ample and the normalization is
finite, the twelfth pencil has no fixed curve on $D$.

Suppose $P\in D$ has canonical index five.  Take its local
canonical cover with point $P^\#$ and a local canonical frame.
The two actual degree-twelve numerators $U,W$ are regular
semi-invariant functions with nontrivial character.  Both vanish
at $P^\#$.  The pullback of the prime quintic is a reduced
Cartier surface: it is generically reduced because the cover is
quasi-etale, and Cohen--Macaulayness excludes embedded or purely
point-supported nilpotents.  Normalize its reduced surface
branches.  This is a finite map.  At every point above $P^\#$,
$U$ and $W$ still vanish.  They have no common curve, since the
finite image of such a curve would be a fixed curve of the
twelfth pencil on $D$.

It follows that a general combination $W-tU$ has a genuine
moving curve branch through a point above $P^\#$.  For precision,
on a normal local surface the nonzero nonunit $W-tU$ has a
height-one zero through that point; no such height-one zero can
be a common component of $U,W$.  Its finite image belongs to
the moving general twelfth curve.  The already proved connected
pencil has Stein degree one, so this branch occurs on the
normalized general curve $C$ used in the marking calculation.

Pull this branch to a common resolution of the canonical and
terminal models.  The numerator $U$ vanishes on the whole inverse
image of $P^\#$, since it is the pullback of a regular function
vanishing there.  Crepancy preserves the canonical frame.
The actual identity $u_{12}|_T=\theta_6^2$ therefore forces a
positive zero of $\theta_6$ at a point of $C$ above $P$.
Every terminal local index over $P$ divides five: a local
trivialization of $5K_Z$ pulls back to one of $5K_Y$.
Thus this zero has either ordinary positive integral order or
positive fifth-order fractional cost.

In both remaining configurations the divisor of $\theta_6$ on
$C$ already has its entire degree $6/7$ at the seventh-order
markings.  These points cannot lie over $P$, because their
terminal index seven does not divide five.  The additional
positive zero just found contradicts this exact budget.
Hence $D$ contains no canonical index-five point.

The canonical index-transfer list is
$\{1,2,3,4,5,6,7,12\}$.  The image of any actual terminal
index-five point has local index divisible by five, and must
therefore have index five.  Thus $T$ avoids all such terminal
points.  Its local class is zero on every virtual basket
constituent of index five.  In the weak Riemann--Roch square
congruence
\[
 IKT^2\equiv-
 \sum_{(b,r)}(I/r)\overline{bi_T}(r-\overline{bi_T})
 \pmod{2I},
\]
every other summand is divisible by five.  Hence $IKT^2\equiv0
\pmod5$.  Of the exhaustive visible states
$IKT^2\in\{-3,11,25\}$, only $25$ remains.

Finally, with $E=5K-T\ge0$ the full exceptional remainder,
\[
 IKE^2=25-10\cdot5+25=0.
\]
Restricting to a general pullback of a sufficiently ample
canonical hyperplane gives an effective exceptional curve
divisor with square zero.  Negative definiteness of the
exceptional intersection form forces it to vanish.  Every
component of $E$ therefore has point image, as asserted.
\end{proof}

\paragraph{Local polar length leaves only the two fourth-seventh markings}
\label{sec:c2353-seventh-local-length}

Proposition~\ref{prop:actual-source-transfer46} makes the canonical contraction an isomorphism near the actual ordinary source $Q_7$.

\begin{proposition}[The remaining marked configuration]
\label{prop:c2353-seventh-local-length}
The degree-$60$ configuration with a single mark $1/7$ is
impossible.  The restricted twelfth pencil on the prime quintic
has exactly the configuration
\[
 \ell=60,\qquad (4/7,4/7).
\]
The quintic is not yet assumed normal.
\end{proposition}

\begin{proof}
Work on the smooth local canonical cover at $Q_7$, with coordinates
$(x,y,z)$ of characters $(1,6,3)$.  Let $f$ be the fifth numerator,
let $U,W$ be the two twelfth numerators, and write
\[
 \alpha=\iota_\delta(dx\wedge dy\wedge dz),\qquad
 \theta=\delta f,\qquad U=\theta^2\pmod f.
\]
The characters of $f,U,W$ are six and that of $\theta$ is three.
The identity is the actual polar-square lift, in the common local
canonical frame.
Separately, our positive monomial filtration on $(x,y,z)$ will
have weights $(1,6,3)$; after eliminating $y$, it has weights
$(1,3)$ on $(x,z)$.  The canonical frame itself has group
character three.  Group characters are taken modulo seven,
whereas filtration weights are ordinary nonnegative integers.

We first justify the exact local colength
\[
 \operatorname{length}\frac{\mathbb C[[x,y,z]]}{(f,W,\theta)}=6.
 \tag{1}\label{eq:c2353-local-colength-six}
\]
The saturated moving degree removes all fixed curves of the
twelfth pencil on the canonical prime quintic.  The preceding
local-isomorphism lemma therefore removes the possible terminal
$A$-null fixed curve at this point as well.  Since $U=\theta^2$
modulo $f$, the displayed ideal is primary for the maximal ideal.
For a general parameter $t$, the curve
$(f,W-tU)$ is a reduced complete intersection; generic reducedness
and Cohen--Macaulayness also exclude embedded nilpotents.  The
element $\theta$ is a parameter on its one-dimensional local ring.
Its colength is the sum of its orders on the normalization
branches.  This follows from the finite-normalization exact
sequence, whose finite-length kernel and cokernel contributions
cancel after multiplication by $\theta$.

The already established two configurations give total polar
degree $6/7$, entirely at $Q_7$: either one full-inertia branch
has numerator order six, or two full-inertia branches have order
three each.  There is no additional ordinary-inertia branch,
which would cost a positive integral amount.  Thus the sum of
the cover-branch orders is six.  The ideal after adding $\theta$
is $(f,W,\theta)$, independent of $t$, proving (1).

At least one of $f,W$ has a nonzero linear $y$ term.  Indeed, this
is their only possible linear monomial, and their only quadratic
monomial is $z^2$.  The only possible linear term of $\theta$ is
$z$, and it has no quadratic monomial.  If neither $f$ nor $W$
has a linear term and $\theta$ has a linear $z$ term, eliminate
$z$ equivariantly.  The resulting series for $z$ in $(x,y)$ has
ordinary order at least three.  Both remaining equations then
have order at least three, so their colength is at least nine.
One may see this directly in degree-at-most-four jets: there are
fifteen monomials and at most six generator multiples.
If $\theta$ has no linear term, its ordinary order is at least
three.  The three-variable jets of degree at most three have
twenty monomials, while the two order-two generators and the
order-three generator span at most $4+4+1$ dimensions there.
The colength is at least eleven.  Both alternatives contradict (1).

We use one more finite-jet observation.  After a linear $y$
elimination, give $(x,z)$ weights $(1,3)$.  A character-three
series has initial part
\[
 h_3=az+bx^3,
\]
and a character-six series has initial part
\[
 g_6=cz^2+dx^3z+ex^6.
\]
There are twelve monomials of weight at most six.  If $h_3=0$,
the two generators span at most one dimension in those jets,
contradicting colength six.  If $h_3\ne0$, its multiples of
weight at most six span five dimensions, with multipliers
$1,x,x^2,x^3,z$.  Colength six therefore forces
\[
 g_6\notin h_3\langle x^3,z\rangle.
 \tag{2}\label{eq:c2353-weighted-coprime}
\]
Equivalently these weighted initial forms have no common factor.
No abstract weighted intersection formula is needed here.

If $f$ has a linear $y$ term, eliminate $f$ and apply this
observation to $h=\theta|_{f=0}$ and $g=W|_{f=0}$.
The initial equation of the general moving curve is
$g_6-t h_3^2$.  Its discriminant as a quadratic in $z$ is
\[
 d^2-4ce+4t(cb^2+ea^2-dab).
\]
The coefficient of $t$ is nonzero by (2), including the case
$a=0$.  For general $t$ the quadratic coefficient is nonzero
and the discriminant is nonzero.  It has two distinct roots
$z=c_i x^3$, with $a c_i+b\ne0$.

If $f$ has no linear term, $W$ does.  Eliminate the linear $y$
term of $W-tU$.  Its solution has weighted order at least six.
The weight-six part of $f$ is still
$cz^2+dx^3z+ex^6$, independent of this substitution.
The characters of the coefficients of
$\partial_x,\partial_y,\partial_z$ in $\delta$ are respectively
$5,3,0$, and their least weighted orders are $5,3,0$.
Since $f$ has no linear $y$ term, its polar therefore satisfies
\[
 \theta_3=\delta_z(0)(2cz+dx^3).
\]
The finite-jet observation forces $\delta_z(0)\ne0$,
$c\ne0$, and $d^2-4ce\ne0$.  Again there are two distinct
weighted roots and the polar has nonzero order-three initial
value on each.

In either case the implicit function theorem gives exactly two
smooth branches
\[
 z=x^3\psi_i(x^7),\qquad y=x^6\phi_i(x^7).
\]
All higher weighted terms differ by multiples of seven, explaining
the indicated argument $x^7$.  Each branch is invariant under
$\mu_7$, with uniformizer $x$ of character one, and its polar
has exact numerator order three.  With ambient canonical
character three its primitive canonical fraction is
$-3/7\equiv4/7$.  Thus the two branches give precisely
$(4/7,4/7)$, excluding the single $1/7$ endpoint.
\end{proof}

\paragraph{The global twelfth pencil is connected and has no fixed divisor}
\label{sec:c2353-global-twelfth-pencil}

\begin{corollary}
\label{cor:c2353-global-twelfth-pencil}
In a hypothetical realization of C2-02353, the complete
global twelfth canonical pencil has connected general fibre and no
fixed divisor on the canonical model.  If $F_Y$ is its primitive
general moving member on a terminal model, then
\[
 420K_Y^2F_Y=12.
\]
\end{corollary}

\begin{proof}
The restricted pencil on the prime quintic is nonconstant and has
Stein degree one.  A nontrivial finite factor in the Stein
factorization of the global pencil would remain a finite factor
after restriction to this horizontal surface.  Hence the global
pencil is also connected.  A general connected restricted fibre therefore lies in one connected global fibre and dominates the original parameter line with the same Stein factor.

On a common resolution $\pi:X\to Z$ of the pencil and of the
quintic, put $L=\pi^*K_Z$ and let $M$ be the basepoint-free moving
divisor of the twelfth system.  Both $L$ and $M$ are nef.
The actual section divisors give
\[
 12L-M\ge0,\qquad \pi^*D=\widetilde D+E\sim_{\mathbb Q}5L,
 \qquad E\ge0.
\]
Restriction to the smooth resolved quintic gives
$420LM\widetilde D=\ell=60$.  Positivity against the full
effective remainders therefore gives
\[
 \frac{60}{420}=LM\widetilde D
 \le5L^2M\le5\cdot12L^3=\frac{60}{420}.
\]
Every inequality is an equality.  Thus $420L^2M=12$.
On the canonical model the fixed divisor of the system has zero
intersection with $K_Z^2$.  Ampleness of $K_Z$ makes this divisor
zero.  Connectedness makes the general moving member primitive,
and projection gives the asserted terminal-model degree.
\end{proof}

\paragraph{The full curve and its nontrivial three-torsion line}\label{sec:2353-node46}
\begin{proposition}\label{prop:c2353-node39}
Suppose C2-02353 is realized.  Let $D_5$ be its prime canonical
quintic and let $D_{12}$ be a general member of its complete
twelfth pencil on the canonical model $Z$.
Their full scheme-theoretic intersection $\Gamma$ is an integral
curve with exactly one singularity: an ordinary node at the
unique index-seven point.  Its normalization $\nu:C\to\Gamma$
is a smooth genus-one curve.  Consequently
\[
 p_a(\Gamma)=2,\qquad
 \nu^*\omega_\Gamma\simeq\omega_C(P+Q),
\]
where $P,Q$ are the two points above the node.  The ninth section gives a point $R\in C$ with
$P+Q\sim2R$.

\end{proposition}

\begin{proof}
The preceding propositions give a connected genus-one general curve of scaled degree $60$, two canonical marks $4/7,4/7$, and a connected global twelfth pencil with no fixed divisor.
The saturated moving degree on the quintic leaves no fixed
curve there.  Since $K_Z$ is ample, there are no positive
residual components in the full intersection.  Its moving
component occurs with multiplicity one, by generic smoothness
on a resolution of the connected pencil.  Locally on the
canonical index cover the two equations form a regular
sequence.  Thus their intersection is Cohen--Macaulay; taking
invariants in characteristic zero preserves the required
one-dimensional unmixedness.  Generic reducedness therefore
makes the complete curve reduced.  It follows that $\Gamma$
is integral and that $C$ is its normalization.

Apply Lemma~\ref{lem:late-ci-trace46} with
$(a,b,d)=(5,12,18)$.  Its exact degree is
\[
 d\deg K_Z|_C=18\cdot\frac{60}{420}=\frac{18}{7}.
\]
At either seventh branch, the trace order obeys
\[
 c\equiv18\cdot4+1\equiv3\pmod7,
 \qquad c\ge0.
\]
The two minimal local contributions already exhaust the
entire degree:
\[
 2g(C)-2+2\cdot\frac{7-1+3}{7}=\frac{18}{7}.
\]
Hence $c=3$ at both seventh branches; every other branch has
trivial inertia and trace order zero.  The final assertion of the same lemma shows that the full local intersection is smooth and has one branch at every other point, including a point with trivial inertia. Thus the coarse curve is smooth away from the unique seventh point.

At that point the local colength calculation gives an equivariant
linear elimination and precisely two smooth branches
\[
 z=x^3\psi_i(x^7),\qquad i=1,2,
 \qquad\psi_1(0)\ne\psi_2(0),
\]
with $\mu_7$ acting on $(x,z)$ with characters $(1,3)$.
After subtracting their mean and dividing $z$ by an invariant
unit, the cover curve is equivariantly
\[
 A=\mathbb C[[x,z]]/(z^2-x^6).
\]
This coordinate change is legitimate because
$(\psi_1-\psi_2)/2$ is an invertible series in $x^7$.
Every element of $A$ is uniquely $a(x)+zb(x)$.  Its invariant
part is generated over $\mathbb C[[u]]$ by $1,v$, where
\[
 u=x^7,\qquad v=x^4z,\qquad v^2=u^2.
\]
Indeed, the invariant monomials in the first summand have
exponents divisible by seven, and in the second summand have
$x$-exponents congruent to four.  This proves the complete
local-ring identity
\[
 A^{\mu_7}=\mathbb C[[u,v]]/(v^2-u^2),
\]
the ordinary node.  It has delta invariant one.  The
normalization exact sequence gives
$p_a(\Gamma)=g(C)+1=2$; ordinary nodal adjunction gives the
displayed dualizing pullback.

For completeness, the proved polar and ninth-section
divisors on $C$ are
\[
 \operatorname{div}(\theta_6)=\frac37(P+Q),\qquad
 \operatorname{div}(a_9)=\frac17(P+Q)+R,
 \qquad\deg R=1.
\]
The ratio $a_9^2/\theta_6^3$ is an ordinary rational function,
with divisor $2R-P-Q$.  This recovers $P+Q\sim2R$ and explains
the degree-two eighteenth map as the degree-two
linear series associated with the node pair.
\end{proof}

\begin{proposition}\label{prop:c2353-three-torsion39}
In any realization of C2-02353, the elliptic K3 normalization above has a section of exact order three, with the proved
ninth-section point $R$ as zero section.
\end{proposition}
\begin{proof}
Work first on the geometric generic normalized curve $C$ of the
twelfth pencil. Put $\Sigma=P+Q$ for the two fourth-seventh markings.
The canonical orbifold line has degree $1/7$ and can be written
formally as
\[
 \mathcal A=\mathcal D+\frac47\Sigma,\qquad
 \deg\mathcal D=-1.
\]
Here additive notation denotes tensor product of line bundles and the
rational marked extension; $\mathcal D$ is its underlying integral
line bundle. The normalized integral line for degree $m$ is
\[
 \mathcal M_m=m\mathcal D+
       \left\lfloor\frac{4m}{7}\right\rfloor\Sigma.
 \tag{E353.1}
\]
The exact polar and ninth-section divisors proved above give
\[
 6\mathcal D+3\Sigma=0,\qquad
 9\mathcal D+5\Sigma=R,\qquad
 \Sigma\sim2R.
\]
Subtracting the first equation from the second yields
$3(\mathcal D+R)=0$. Thus
\[
 \varepsilon:=\mathcal D+R\in\operatorname{Pic}^0(C)
 \quad\hbox{satisfies}\quad 3\varepsilon=0.
 \tag{E353.2}
\]
It remains to prove that $\varepsilon$ is nonzero. This is where the
full complete-intersection module supplies more information than
the ordinary marked-degree inequality.

Let $\Gamma=D_5\cap D_{12}$ be the coarse curve and let
$\mathcal F_m$ be the pushforward of its degree-$m$ canonical twist
from the canonical covering stack. At the unique node the full
cover ring is
\[
 A=\mathbb C[[x,z]]/(z^2-x^6),\qquad
 \operatorname{wt}(x,z)=(1,3)\pmod7.
\]
The numerator of a degree-$m$ canonical section has character
$-3m\equiv4m\pmod7$. Writing $u=x^7$, its character-four part is
\[
 A_4=x^4\mathbb C[[u]]\ \oplus\ xz\mathbb C[[u]].
\]
On the two normalized branches $z=\pm x^3$, these summands give
$x^4(a(u)+b(u))$ and $x^4(a(u)-b(u))$. Their coefficients are
independent. Since $8\equiv1\pmod7$, this proves that $\mathcal F_8$
is the full normalization pushforward
\[
 \mathcal F_8=\nu_*\mathcal M_8,\qquad
 \mathcal M_8=8\mathcal D+4\Sigma=8\varepsilon.
 \tag{E353.3}
\]
There is no additional matching condition at the node. Away from it,
the exhausted full-trace budget shows that the complete curve is
smooth with trivial branch inertia, so the same identity extends
globally.

Lemma~\ref{lem:complete-ci-canonical-restriction}, with
$(a,b,m)=(5,12,8)$ and denominator degrees $(3,-4,-9)$,
makes $H^0(8K_Z)\to H^0(\Gamma,\mathcal F_8)$ surjective.
Since $P_8=0$, (E353.3) gives
\[
 H^0(C,8\varepsilon)=0.
\]
A trivial degree-zero line bundle on the genus-one curve has a
nonzero section. Therefore $\varepsilon\ne0$, proving exact order
three.

The point is defined over the pencil field $\mathbb C(t)$.
One can avoid any descent ambiguity in the fractional notation:
the ordinary line $\mathcal M_7$ is the pullback of the invertible
degree-seven canonical twist on the nodal curve, and
\[
 \mathcal M_7-R
 =7\mathcal D+4\Sigma-R=7\varepsilon=\varepsilon.
\]
Both $\mathcal M_7$ and $R$ are defined over $\mathbb C(t)$.
With $R$ as the origin this is a rational nonzero three-torsion
point of the elliptic generic fibre. Properness extends it to
a section of the relatively minimal elliptic K3 surface.
\end{proof}

\paragraph{A proper intersection from four sections}\label{sec:2353-frame46}
\begin{lemma}[Nontrivial normalization twist prevents a node base point]
\label{lem:r40-nodal-torsion}
Let $\Gamma$ be an integral projective curve with a single ordinary
node $q$, and let $\nu:C\to\Gamma$ be its smooth genus-one
normalization. Write $\nu^{-1}(q)=\{P,Q\}$ and assume
$P+Q\sim 2R$ for a point $R\in C$. Let $L$ be an invertible
sheaf on $\Gamma$ such that
\[
             \nu^*L\simeq\mathcal O_C(2R)\otimes\eta,
             \qquad 0\ne\eta\in\operatorname{Pic}^0(C).
\]
Then $h^0(\Gamma,L)=1$, and its nonzero section does not vanish at
$q$.
\end{lemma}
\begin{proof}
The degree of $L$ is two and $p_a(\Gamma)=2$, so its Euler
characteristic is one. If $\omega_\Gamma\otimes L^{-1}$ had a
nonzero section, its pullback would be a nonzero section of
$\mathcal O_C(P+Q-2R)\otimes\eta^{-1}\simeq\eta^{-1}$,
which is impossible for a nontrivial degree-zero line bundle on a
smooth genus-one curve. Serre duality gives $h^1(\Gamma,L)=0$,
and hence $h^0(\Gamma,L)=1$.
If a nonzero section vanished at $q$, its pullback would vanish at
both $P$ and $Q$, because $L$ is invertible at the node. These two
zeros exhaust its degree. It would follow that
$\nu^*L\simeq\mathcal O_C(P+Q)\simeq\mathcal O_C(2R)$,
contrary to $\eta\ne0$.
\end{proof}

\begin{proposition}[The degree-fourteen section avoids the seventh point]
\label{prop:r40-c2353-node-avoidance}
Assume C2-02353 is realized. Let $a\in R_5$ be its prime quintic
section, and choose a basis $u,b$ of $R_{12}$ with
$u|_{D_5}=\theta_6^2$. For a general choice of $b$, there is
$c\in R_{14}$ whose value at the unique index-seven point is
nonzero, in the local trivialization of $14K_Z$.
\end{proposition}
\begin{proof}
Use the preceding conclusions that the entire curve
$\Gamma=D_5\cap\operatorname{div}(b)$ is integral, smooth except
for one ordinary node $q$ at the unchanged index-seven point,
and has genus-one normalization $C$. Let $P,Q\in C$ lie over
$q$. The ninth section gives $R\in C$ with
$P+Q\sim2R$. Write the integral part of the canonical orbifold
line as $D$, so that the ordinary degree-$m$ normalization line is
\[
             M_m=mD+\left\lfloor\frac{4m}{7}\right\rfloor(P+Q).
\]
The preceding proposition proves that $\epsilon=D+R$ has exact order three in
$\operatorname{Pic}^0(C)$.

Let $\mathcal F_{14}$ denote the complete degree-fourteen
canonical twist on $\Gamma$, obtained from the canonical cover.
It is invertible. At $q$, degree fourteen has trivial character
modulo seven, so its module is the invariant local ring of the
node. Away from $q$, the exhausted full trace budget shows that
$\Gamma$ is smooth with trivial canonical inertia. Its pullback is
\[
 \nu^*\mathcal F_{14}=M_{14}
 =14D+8(P+Q)\sim2R+14\epsilon=2R+2\epsilon.
\]
The normalization twist $2\epsilon$ is nonzero. The lemma implies
that $H^0(\Gamma,\mathcal F_{14})$ has dimension one and its
nonzero section is nonzero at $q$.

Lemma~\ref{lem:complete-ci-canonical-restriction}, now with
$(a,b,m)=(5,12,14)$ and denominator degrees $(9,2,-3)$,
gives a surjective restriction map with kernel the image of
$R_9\oplus R_2$. Its image has dimension
$P_{14}-P_9-P_2=2-1-0=1$. A lift $c$ of the nonzero section
therefore has the asserted nonzero value at $q$.
\end{proof}

\begin{lemma}[A proper intersection for C2-02353]
\label{thm:r40-c2353-finite}
With $a,u,b,c$ as above, the four sections
$a^{84},u^{35},b^{35},c^{30}$ of $420K_Z$ have no common zero.
The triple intersection $(a_5,u_{12},c_{14})$ is zero dimensional,
has total rational intersection number $2$, and has local length
at most $24$ on the smooth index-twelve cover at $Q_{12}$.
\end{lemma}
\begin{proof}
On the complete curve $\Gamma=D_5\cap\operatorname{div}(b)$, the
polar identity $\operatorname{div}(\theta_6)=\frac37(P+Q)$ implies
that $u=\theta_6^2$ has no zero away from the node, after taking
a Cartier power. Proposition~\ref{prop:r40-c2353-node-avoidance}
gives $c$ nonzero at the node. Thus $u^{35}$ and $c^{30}$ have
no common zero on $\Gamma$, and the four displayed sections
have no common zero on $Z$.

If $(a,u,c)$ contained a curve, $b^{35}$ would be nowhere zero
on that curve, contrary to the ampleness of $420K_Z$.
The intersection is therefore proper, with
\[
 (5K_Z)(12K_Z)(14K_Z)=2.
\]
Its local contributions are nonnegative, as follows by taking
Cartier powers. On the smooth index-twelve cover the local
equations form a regular sequence, and a local length $\ell$
contributes $\ell/12$. Hence $\ell\le24$.
\end{proof}


\label{sec:2353-normality46}
Use the prime quintic $D=(a_5=0)$, the proper intersection just proved,
and the ordinary unchanged sources $Q_7$ and $Q_{12}$.
All following calculations take place on actual local canonical covers.

\begin{proposition}[The twelfth-index slice is smooth on its cover]
\label{prop:r40-c2353-Q12-smooth}
The degree-five numerator has a nonzero linear term at $Q_{12}$.
Consequently $D$ is normal and klt at $Q_{12}$, with surface quotient
type $\tfrac1{12}(1,5)$.
\end{proposition}
\begin{proof}
Take smooth cover coordinates $(x,y,z)$ of characters $(1,11,5)$
modulo twelve. Write the actual local numerators as $f,u,c,\theta$.
Their characters are respectively $11,0,2,6$. The square lift means
\[
                         u=\theta^2+fh,
\]
where $h$ is a regular semi-invariant of character one.

Suppose $f$ has no linear term. There is no degree-two monomial of
character eleven, so $\operatorname{ord}f\ge3$. Its degree-three part is
\[
                         f_3=x(Ay^2+Bz^2).
\]
The only degree-two monomial of character six is $xz$; there is no
degree-one monomial of that character. Thus $\operatorname{ord}\theta\ge2$.
The only degree-one monomial of character one is $x$, so
$\operatorname{ord}h\ge1$. It follows that $\operatorname{ord}u\ge4$,
and its degree-four part is divisible by $x^2$. Finally
$\operatorname{ord}c\ge2$, and $c_2$ is a scalar multiple of $x^2$.

Lemma~\ref{thm:r40-c2353-finite} gives
\[
 \ell=\operatorname{length}\mathbb C[[x,y,z]]/(f,u,c)\le24.
\]

For a Cohen--Macaulay local ring $A$ and a parameter $g\in\mathfrak m^s$,
the Hilbert--Samuel function comparison
\[
 \ell(A/(g,\mathfrak m^N))\ge
 \ell(A/\mathfrak m^N)-\ell(A/\mathfrak m^{N-s})
\]
gives $e(A/gA)\ge s e(A)$. Apply this successively to the regular
sequence in the regular local ring of the cover. It gives
$\ell\ge\operatorname{ord}f\operatorname{ord}u\operatorname{ord}c$.
If any of these orders exceeds $(3,4,2)$, this contradicts
$\ell\le24$. Otherwise $c_2$ is a nonzero multiple of $x^2$,
and $u_4=q_2c_2$ for a homogeneous quadratic polynomial $q_2$.
Replacing $u$ by $u-q_2c$ leaves the ideal unchanged and increases
its order to at least five. B\'ezout then gives
\[
                             \ell\ge3\cdot5\cdot2=30>24.
\]
The contradiction proves that $f$ has a linear term. Its only possible
linear monomial is $y$. Equivariant elimination of $y$ leaves the
smooth cover with coordinates $(x,z)$ of weights $(1,5)$, as claimed.
\end{proof}

\begin{theorem}[Global normality of the quintic]
\label{thm:r40-c2353-normality}
In every hypothetical realization of C2-02353, the canonical quintic
$D$ is normal.
\end{theorem}
\begin{proof}
The surface is $S_2$. Locally on a canonical cover it is the reduced
hypersurface cut out by the actual quintic numerator in a
Cohen--Macaulay Gorenstein threefold; its quotient is again
Cohen--Macaulay in characteristic zero. Thus any failure of normality
has a curve centre.

Suppose $\Gamma$ is such a conductor curve, and put
$w=420K_Z\Gamma$. On a resolution $S$ of the normalization, canonical
adjunction with its effective integral conductor $C$ gives
\[
          420A\bigl(\operatorname{div}(\alpha|_S)+C\bigr)=30,
          \qquad A=\pi^*(K_Z|_D).
\]
The generic transverse ring of $D$ along $\Gamma$ is a reduced
one-dimensional Gorenstein ring. Its conductor has total degree
$2\delta\ge2$ over the residue field. Consequently $2w\le30$.
The previously proved point-centred-remainder theorem says that
the canonical model is smooth at a generic point of every curve
contained in $D$, so the transverse germ here is a reduced plane curve.

Generic-centre marked adjunction for the nonnormal prime $D\sim5K_Z$
has coefficient at most six. Index-five avoidance leaves inertia
orders $2,3,4,6,7,12$. The primitive markings obey
\[
 \sum k_x/r_x\equiv w/420\pmod{\mathbb Z},\qquad
 -2+\sum(1-1/r_x)\le6w/420.
\]
For $1\le w\le15$, the exhaustive exact list is
\[
\begin{array}{c|c|c|c}
 w&\text{primitive markings}&\text{degree-nine cost}&9w/420\\\hline
 5&(3/7,7/12)&31/28&3/28\\
10&(1/2,2/3,6/7)&17/14&3/14\\
10&(1/6,6/7)&17/14&3/14\\
15&(1/2,1/4,2/7)&37/28&9/28\\
15&(3/4,2/7)&37/28&9/28
\end{array}
\]
Here the degree-nine cost is $\sum\overline{9k_x}/r_x$. Each cost
exceeds the full degree of the ninth section. Hence the actual
section $a_9$ vanishes identically on every conductor curve.
Every listed row has an index-seven mark, so every conductor curve
passes through the unique point $Q_7$.

Use the local coordinates $(x,y,z)$ at $Q_7$, of characters
$(1,6,3)$. On the two branches of the general twelfth complete
intersection, the proved divisor formula for $a_9$ gives numerator
order exactly one. Indeed its only additional zero is $R$, and
$R\ne P,Q$: otherwise $P+Q\sim2R$ on the elliptic normalization
would force $P=Q$. The only linear monomial of the character-one
ninth numerator is $x$. Thus $a_9$ has nonzero linear $x$ term and
can itself be chosen as the equivariant coordinate $x$.

If the quintic numerator $f$ had a linear $y$ term, $D$ would be
normal near $Q_7$, and no conductor could pass through it. We may
therefore take the no-linear-term case of the local-length
six argument. That argument proves that the coefficient of $z^2$
in $f$ is nonzero. Thus $f_z$ has a nonzero linear $z$ term.
Changing $x$ to the ninth numerator does not alter this assertion.

The lift of every conductor germ is contained both in $(x=0)$,
because $a_9$ vanishes on it, and in $(f_z=0)$, because it belongs
to the singular locus of the reduced hypersurface. The implicit
function theorem makes
\[
                             x=f_z=0
\]
a single smooth cover curve parametrized by $y$, of character six.
It is invariant under the full group. Its primitive canonical
fraction at $Q_7$ is therefore
\[
                   -3\cdot6^{-1}\pmod7=3,\qquad\text{i.e. }3/7.
\]
Only the first row of the displayed list has that mark.
Every conductor curve must therefore have $w=5$ and the second
mark $7/12$. In particular it must pass through $Q_{12}$.
Proposition~\ref{prop:r40-c2353-Q12-smooth} gives normality in an
entire neighbourhood of that point, a contradiction. Hence $D$
has no conductor curve, and $S_2$ now proves normality everywhere.
\end{proof}

\begin{corollary}[Rational klt quintic and no crepant remainder]
\label{cor:r40-c2353-klt}
The surface $D$ has rational and klt singularities. The pair $(Z,D)$
is purely log terminal. On every crepant terminalization the full
quintic divisor is its prime strict transform, with no exceptional
component.
\end{corollary}
\begin{proof}
The complete divisorial sequence and Serre duality give
\[
 H^1(D,\mathcal O_D)=0,\qquad
 \chi(\mathcal O_D)=\chi(\mathcal O_Z)+P_6=2.
\]
Here $H^2(Z,-5K_Z)$ is dual to $H^1(Z,6K_Z)=0$.
The normalization resolution has a K3 minimal model,
so its Euler characteristic is also two. Normality gives
$\pi_*\mathcal O_S=\mathcal O_D$, and Leray yields
\[
 \sum_{P\in D}\operatorname{length}(R^1\pi_*\mathcal O_S)_P
       =\chi(\mathcal O_D)-\chi(\mathcal O_S)=0.
\]
Thus all singularities of $D$ are rational.

Put $H=K_Z|_D$; adjunction gives $K_D=6H$.
The complete canonical-cover quotient sheaves are the reflexive
restrictions on the now normal surface. Exact restriction and
canonical vanishing give
\[
                    H^0(D,7H)=H^0(D,8H)=0.
\]
The Weil gap criterion \cite[Lemma~10.6(i)]{LBVolume} applies to these
reflexive restrictions: when $K_D\sim_{\mathbb Q}cH$, the vanishing
$H^0(D,jH)=0$ for $j>c$ implies that $D$ is klt wherever $jH$
is Cartier. Here $c=6$, and we use $j=7,8$.
Degree seven covers the canonical ambient index-seven points;
degree eight covers indices two and four. At ambient indices
$1,2,3,6$, the divisor $K_D=6H$ is Cartier, and a rational
Gorenstein surface singularity is Du Val. There are no index-five
points on $D$, and the remaining index-twelve point is klt by
Proposition~\ref{prop:r40-c2353-Q12-smooth}. The source-transfer index list
$\{1,2,3,4,5,6,7,12\}$ is exhaustive. Hence $D$ is klt everywhere.

The point-centred-remainder argument already showed that the ambient
model is regular at the generic point of every curve contained in $D$.
Hence the divisorial different is zero. Divisorial adjunction and inversion
of adjunction make $(Z,D)$ purely log terminal.
For a crepant terminalization write the actual full member as
$g^*D=T+\sum c_iE_i$. The coefficients $c_i$ are nonnegative
integers, and the corresponding pair discrepancies are $-c_i$.
Pure log terminality forces $c_i<1$, hence every $c_i=0$.
\end{proof}

\Needspace{5\baselineskip}
\paragraph{The seventh quotient and the complete discrepancy count}
\begin{proposition}\label{prop:2353-two-polar46}
The positive polar divisor $E=\operatorname{div}_D\theta_6$ is
reduced and has exactly one of the following forms:
\[
\begin{array}{c|c|l}
 \text{seventh cover}&420HC&
       \text{primitive canonical marks on its normalization}\\\hline
 \text{smooth}&30&4/7;\ 7/12,11/12\\
 D_7&25&1/7;\ 11/12\\
    &5&3/7;\ 7/12.
\end{array}
\]
In the first case there is one component; in the second there are
exactly two. The semicolon separates the seventh and twelfth sources.
All components are rational. At the twelfth source, the two lifted
polar branches on the smooth surface cover are smooth and transverse.
\end{proposition}
\begin{proof}
Let $\pi:W\to D$ be the minimal resolution and $\sigma:W\to S$
its contraction to the K3 minimal model. Write
$\Delta=\pi^*K_D-K_W$. Here $K_D=6H$ and $H^2=1/84$.
Away from the seventh and twelfth sources all singularities of $D$
are Du Val, even when the low-index sources coalesce. Indeed the
full terminal quintic has no exceptional remainder; at every other
terminal point its ambient index divides six, because it avoids
all index-five points. Its klt surface germ is therefore Gorenstein
and Du Val. Its crepant image on $D$ remains canonical.

At $Q_7$ the quintic numerator $f$ has character six in coordinates
$(x,y,z)$ of characters $(1,6,3)$. A linear term is a nonzero
$y$ term; then its surface quotient is $\frac17(1,3)$, with chain
$[3,2,2]$. If it has no linear term, the colength-six argument gives
\[
 f_6=c z^2+d x^3z+e x^6,
 \qquad c\ne0,\quad d^2-4ce\ne0
\]
for the positive filtration $(1,6,3)$. The index-one surface is
klt and Gorenstein, hence Du Val. Its quadratic part has rank one.
A rank-one quadratic Du Val germ has a nonzero cubic part in its
remaining variables, and the only allowed cubic monomial is $xy^2$.
Its coefficient is therefore nonzero. Complete the $z$ square.
For weights $(2,5,6)$ the resulting degree-twelve initial equation
is $z^2+xy^2+\lambda x^6$, with $\lambda\ne0$; every other term
has higher weight. The Jacobian algebra of this $D_7$ equation has
no class of weight greater than ten, so analytic finite determinacy
removes all higher terms, equivariantly. Thus the second cover is
exactly $D_7$, with its specified $\mu_7$ action.

The quotient $D_7/\mu_7$ has central vertex $U$ of square $-3$,
two short arms $A_+,A_-$ of square $-2$, and a long arm
$D_0(-2)-L(-3)$. This graph can be checked directly from the finite
group description: $D_7$ is the binary dihedral plane quotient of
order twenty, whose invariant generators have degrees $4,10,12$;
scalar multiplication by $\zeta_7^2$ induces characters $(1,6,3)$.
The enlarged group is $D_{12,5}$ in \cite[\S\S5,7]{WemyssMcKay},
and $12/5=[3,2,3]$ gives the central and long-arm weights, with the
two short $[2]$ arms. Equivalently, after blowing up the plane, the
projective stabilizers are $2,2,5$ and the last quotient is
$5/3=[2,3]$. The central square is then fixed by adjunction.

Solving the rational intersection equations gives
\[
\begin{array}{c|c|c}
 \text{germ}&\Delta\text{ in graph order}&-\Delta^2\\\hline
 \frac17(1,3)&(3,2,1)/7&3/7\\
 D_7/\mu_7&(6,3,3,5,4)/7&10/7\\
 \frac1{12}(1,5)&(1,1,1)/2&1.
\end{array}
\]
Write the last chain as $B(-3)-V(-2)-C(-3)$.
There is no further discrepancy cost, by the preceding Du Val
argument. Since $K_D^2=3/7$, we obtain
\[
 -K_W^2=1\quad\hbox{or}\quad2,
\]
respectively. These are exactly the numbers of point blowups from
$S$ to $W$.

We need to know where every positive canonical component goes.
Its image is an irreducible component of the actual intrinsic
polar $E$, of degree $1\le w\le30$. The square lift $u_{12}|_D=
\theta^2$ and a generic log-canonical threshold at that component
give marked adjunction coefficient at most twelve. The other
twelfth pencil generator is nonzero on that component, because
the restricted pencil has no fixed curve. Hence its marks satisfy
\[
 \sum k/r\equiv w/420\pmod{\mathbb Z},\qquad
 -2+\sum(1-1/r)\le12w/420,\qquad
 \sum\{12k/r\}\le12w/420.
\]
The allowed inertias are $2,3,4,6,7,12$. Genus one would allow
at most one nontrivial mark, whose degree already forces $w\ge35$;
higher genus is impossible. There are at most five marks. Enumerating
these explicitly bounded primitive fractions gives thirty-five
rows, all containing an index-seven mark. Thus every positive
component of $K_W$ meets the complete seventh exceptional fibre.
This conclusion concerns all actual intrinsic polar components;
it makes no claim about arbitrary curves of the same degree.

In the smooth-cover case the single blowup gives $K_W=E_0$, with
$E_0^2=-1$. Its intersection with each $(-3)$ curve is one and
with every $(-2)$ exceptional curve is zero. Thus it meets the
first seventh vertex and both twelfth ends, transversely at their
smooth open parts. Its image has degree thirty.

In the $D_7$ case, suppose the two blowups were nested. Then
$K_W=E_1+2E_2$, with $E_1^2=-2,E_2^2=-1$. At each of the four
$(-3)$ vertices $U,L,B,C$, adjunction gives $K_W R=1$.
Consequently $E_2$ misses all four and $E_1$ meets each once.
In particular $E_1$ is not $\pi$-exceptional, since it meets
fibres over two different points. Adjunction at the seventh
$(-2)$ vertices now makes $E_2$ miss them as well. But a final
$(-1)$ curve is not exceptional for the minimal resolution $\pi$,
so $E_2$ has positive $H$-degree. It would be a positive polar
component missing $Q_7$, a contradiction.

The two blowups are therefore disjoint and $K_W=E_1+E_2$.
Each $E_i$ meets one of $U,L$ once and no $(-2)$ vertex. Positivity
\[
 6HE_i=(K_W+\Delta)E_i=-1+\Delta E_i>0
\]
forces it to meet a twelfth end as well. The two curves exhaust
the unit canonical intersections at $B,C$, so each meets exactly
one end once. All stated intersections are transverse and avoid
the nodes of the exceptional graphs.

The actual integral polarization determines the pairing. Choose
an integral Weil representative of $H$ and write
$\pi^*H=M-\sum t_iR_i$, with $M$ integral and $0\le t_i<1$.
At $Q_7$ the class has order seven and $K_D=6H$, so
$t_i=\Delta_i$. At $Q_{12}$ orient the rays as
$(1,5),(3,3),(5,1)$ divided by twelve. The numerator character of
$H$ is seven; the representative $x^7$ gives fractional complements
$(5,3,1)/12$. Thus in
\[
 6M=K_W+\sum(\Delta_i+6t_i)R_i
\]
the twelfth coefficients are $(3,2,1)$, and the $D_7$ coefficients
are $(6,3,3,5,4)$. A curve joining $U$ to a twelfth end of
coefficient $c$ has $6ME_i=-1+6+c$. Integrality forces $c=1$;
the other value $c=3$ gives eight, not divisible by six.
Therefore $U$ joins $C$ and $L$ joins $B$. Their degrees are
\[
 420HE_1=70(-1+6/7+1/2)=25,
 \qquad
 420HE_2=70(-1+4/7+1/2)=5.
\]

The marks are the fractional coefficients of $\pi^*H$ at these
transverse intersections. At the seventh vertices they are
$4/7$ in the smooth-cover case, and $1/7,3/7$ at $U,L$.
At the two twelfth ends they are $7/12,11/12$ respectively.
No further mark occurs: each positive curve misses all other
exceptional fibres, since $K_W$ is their sum and a $(-2)$ fibre
has zero canonical intersection. Their normalizations are the
smooth rational $E_i$. At the twelfth ends the primitive arc orders
are $(1,5)$ and $(5,1)$ on the surface cover. Each lift is smooth,
and the two tangent directions are distinct. The asserted table
and transversality follow. Finally all positive canonical
coefficients are one, so the polar divisor is reduced.
\end{proof}


\label{sec:2353-exclusion46}
\begin{theorem}\label{thm:2353-exclusion46}
There is no canonical threefold with $\chi=2$, $K^3=1/420$ and
terminalization basket
$3(1,2)+4(1,3)+(1,5)+2(2,5)+(3,7)+(5,12)$.
\end{theorem}
We use the prime normal klt quintic $D=(a_5=0)$ and the complete
polar classification of Proposition~\ref{prop:2353-two-polar46}.
The final obstruction uses an identity on the whole quintic and
its transverse jet in the ambient threefold.

\paragraph{The normalized curvature cannot vanish on the quintic}
Write $a=a_5$ and, in a local canonical frame,
$\theta=\delta a$, $\gamma=\delta\theta$, and
$\mathcal D_a f=5a\delta f-j\theta f$ for a section of degree $j$.
The actual normalized square of Lemma~\ref{lem:polar-square46} is
\[
 u=\theta^2-\frac{10}{11}a\gamma+\frac1{11}a^2\rho\in R_{12},
 \qquad v=-\mathcal D_a u/12\in R_{18}.
\]
The coefficient $\rho$ is regular in the local frame.
Corollary~\ref{cor:polar-cubic-covariant} gives global sections
$w\in R_{26}$ and $s\in R_{14}$ satisfying
\begin{equation}\label{eq:2353-curvature-chain46}
 v^2-u^3=a^2w,\qquad
 \mathcal D_a v=-18u^2+a^2s,\qquad
 \mathcal D_a w=2vs.
\end{equation}
Its restriction formulas give
\begin{equation}\label{eq:2353-curvature-restriction46}
 w|_D=-\theta^2s|_D/13,
 \qquad s|_E=(975/121)\gamma^2|_E,
 \quad E=\operatorname{div}_D\theta.
\end{equation}

Let $b=b_9$ generate $R_9$. It restricts nontrivially to $D$,
since $R_4=0$, and is nonzero on some component of the reduced
polar divisor $E$. Otherwise $b/\theta$ would have no divisorial
pole on the normal surface and would extend to a nonzero section
of $3H$. But complete restriction gives
\[
 H^0(D,3H)=0,
 \qquad H^1(Z,-2K_Z)=H^2(Z,3K_Z)^*=0.
\]

Suppose $s|_D=0$. Since $R_9=\mathbb C b$, normal division gives
$s=\kappa ab$. The first identity of
\eqref{eq:2353-curvature-restriction46} gives $w=a w_{21}$, and
\eqref{eq:2353-curvature-chain46} then gives
\[
 \mathcal D_a w_{21}=2\kappa vb,
 \qquad w_{21}|_D=-\frac{2\kappa}{21}\theta^2b.
\]
The second identity of \eqref{eq:2353-curvature-restriction46}
makes $\gamma$ vanish on $E$. At the generic point of a component
where $b$ is nonzero, the ambient threefold is smooth and
$\mathfrak l=(a,\theta)$ is the maximal ideal of a regular
local ring of dimension two: $D$ is normal and $E$ is reduced.
Both $\delta a$ and $\delta\theta$ lie in $\mathfrak l$, so
$\delta$ preserves every power. The normalized formula now gives
$u\in\mathfrak l^2$, $v\in\mathfrak l^3$, and
\[
 a^3w_{21}=v^2-u^3\in\mathfrak l^6.
\]
Since $a$ is a regular parameter, $w_{21}\in\mathfrak l^3$.
Its restriction to $D$ has order at least three, whereas
$-2\kappa\theta^2b/21$ has order exactly two if $\kappa\ne0$.
Thus $\kappa=0$ and $s=0$.

Now $\mathcal D_a w=0$. A nonzero $w$ of $a$-order $k$ has
leading derivative coefficient $(5k-26)\theta$, never zero for
an integer $k$, because $D$ is visible. Hence $w=0$. Then
$v^2=u^3$, and the rational section $v/u$ has nonnegative valuation
at every prime, of degree six. It extends to a nonzero element of
$R_6=0$, the contradiction. We have proved $s|_D\ne0$.

The quotient $R_{14}/aR_9$ has dimension one. The finite-frame
section $c_{14}$ spans it and is a unit at $Q_7$; hence so is $s$.
In particular \eqref{eq:2353-curvature-restriction46} makes
$\gamma$ a unit on every polar germ at $Q_7$.

\paragraph{Sharp restriction degrees force three nonzero coefficients}
At $Q_{12}$ the quintic is smooth on its cover. Choose equivariant
coordinates $(x,a,z)$ of characters $(1,11,5)$, using its nonzero
linear term. The ambient derivation $\delta$ has character seven.
The two transverse reduced polar branches give
\begin{equation}\label{eq:2353-theta46}
 \theta|_D=A xz+O((x,z)^6),\qquad A\ne0.
\end{equation}
The two axis directions are distinguished by their characters.
The intrinsic section $\gamma$ has degree seven on $E$ and is
nonzero on each component, by its unit value at their seventh
source. In the degree-thirty case its whole degree is
\[
 7\cdot30/420=1/12+5/12.
\]
In the two-component case its degree on the small component is
$7\cdot5/420=1/12$. The seventh costs are zero in degree seven.
In both cases it therefore has exact numerator order one on the
$x$ branch at $Q_{12}$.

The coefficient $\delta x$ has character eight and no monomial
of ordinary degree below four. If $c=\delta z(0)$, differentiating
\eqref{eq:2353-theta46}, including normal-direction terms, gives
$\gamma|_D=Acx+O((x,z)^5)$. Thus $c\ne0$.

The section $b$ is nonzero on the whole degree-thirty component
in the first case. Its exact degree is exhausted by
\[
 9\cdot30/420=1/7+1/4+1/4.
\]
In the other case it vanishes on the degree-five component,
since that component's forced degree-nine cost exceeds its degree:
\[
 9\cdot5/420-(6/7+1/4)=-1.
\]
It is therefore nonzero on the degree-twenty-five component, where
$9\cdot25/420=2/7+1/4$. In either case the $z$ branch at $Q_{12}$
has exact numerator order three. Its character is three, so
\begin{equation}\label{eq:2353-b46}
 b|_D=B x^3+Cz^3+O((x,z)^7),\qquad C\ne0.
\end{equation}
No assumption about $B$ is needed.

\paragraph{The actual degree-seventeen identity and its ambient jet}
Define actual divided sections
\[
 d=\mathcal D_a b\in R_{15},\qquad
 e=\frac{\mathcal D_a d-135ub}{a}\in R_{16},\qquad
 j=\frac{ud+9vb}{a}\in R_{22}.
\]
Their numerators vanish along $D$; division at its generic point
and normality make them global sections. Direct polynomial
expansion before restriction gives
\[
 e=-95\theta\delta b+\frac{855}{11}\gamma b
       +a\left(25\delta^2b-\frac{135}{11}\rho b\right),
 \qquad j|_D=-\theta e|_D/19.
\]
Thus $g=(\mathcal D_a e-304j)/a\in R_{17}$ is also an actual
section. Since $P_{17}=P_{12}=2$, multiplication by $a$ gives
$R_{17}=aR_{12}$ and hence $g|_D=0$. Its exact expansion is
\begin{equation}\label{eq:2353-mixed46}
 550\theta\delta^2b-950\gamma\delta b
       +475(\delta\gamma)b-251\theta\rho b=0\quad\hbox{on }D.
\end{equation}

We record the full ambient character check to exclude a missing
normal derivative. In variables $(x,a,z)$ of characters $(1,11,5)$,
\[
\begin{array}{c|c|l}
 \text{coefficient}&\text{character}&\text{lowest possible terms}\\\hline
 \delta x&8&\text{ordinary degree at least four}\\
 \theta&6&Axz;\ \text{the normal term }ax^2z\text{ has degree four}\\
 \delta z&0&c;\ \text{the normal term }xa\text{ has degree two}\\
 b&3&Bx^3+Cz^3+D a^2z\\
 \rho&2&\text{ordinary degree at least two}\\
 \delta\gamma&8&\text{ordinary degree at least four}.
\end{array}
\]
In particular the normal cubic $Da^2z$ contributes
$2D\theta^2z$ to $\delta^2b|_D$, of order five. Terms of $b$
linear in $a$ first occur in ordinary degree five and contribute
no degree-one term after two derivatives. The allowed normal
terms of $\delta z$ and $\theta$ likewise contribute only above
the orders displayed below. Consequently
\[
\begin{aligned}
 \theta|_D&=Axz+O(6),& \gamma|_D&=Acx+O(5),\\
 \delta b|_D&=3cCz^2+O(6),&
 \delta^2b|_D&=6c^2Cz+O(5).
\end{aligned}
\]
Here $O(q)$ means total degree at least $q$ in $(x,z)$.
The other two terms of \eqref{eq:2353-mixed46} have order at least
seven. Its degree-three part is therefore
\[
 (550\cdot6-950\cdot3)Ac^2C\,xz^2
     =450Ac^2C\,xz^2=0,
\]
contradicting $A,c,C\ne0$. This proves the theorem.

The calculations use the whole surface identity and its ambient
normal derivatives. They are unchanged by equivariant coordinate
changes or by an equivariant canonical-frame unit: the extra terms
then have the higher orders already included in the character list.



\providecommand{\OcticRoot}{../octic}

\section{The octic basket and its conductor}
\label{sec:octic-foundation-revised}

We turn to the basket C2-02350.  Its first canonical section has degree
eight.  The geometry of this section is more restrictive than its
Hilbert function alone suggests: its normalization is birational to a
K3 surface, but the section itself is nonnormal.  The conductor and the
divisor of the restricted two-form then leave only eleven possibilities.

\begin{figure}[htbp]
\centering
\begin{tikzpicture}[>=Stealth,font=\small,
 box/.style={draw=inkblue,rounded corners=2pt,align=center,inner sep=5pt,
             text width=61mm,minimum height=15mm},
 wide/.style={box,text width=128mm},
 route/.style={->,draw=inkblue,line width=.5pt}]
\node[wide] (foundation) {Section~\ref{sec:octic-foundation-revised}: the prime octic, its normalization and conductor\\The full conductor and polar budget give eleven total states};
\node[box,below=10mm of foundation,xshift=-34mm] (census) {The exhaustive census\\A paired conductor counts as one total state};
\node[box,below=10mm of foundation,xshift=34mm] (atlas) {Section~\ref{sec:octic-atlas-revised}: common local geometry\\Q5, Q12 and all seven Q14 formats\\Canonical forests and K3 obstructions};
\node[box,below=15mm of census] (eight) {Section~\ref{sec:octic-eight-states-revised}: eight total states\\All required node/cusp and polar alternatives};
\node[box,below=15mm of atlas] (higher) {Section~\ref{sec:octic-higher-delta-revised}: $(7,2),(7,3),(7,4)$\\Conductor filtrations\\Cohomology and multiplication};
\node[wide,below=13mm of eight,xshift=34mm,fill=inkblue!7] (end) {Section~\ref{sec:octic-completion-revised}: all eleven states are excluded\\The basket \basket{C2-02350} is impossible};
\draw[route] (foundation)--(census);
\draw[route] (foundation)--(atlas);
\draw[route,dashed] (census)--(eight);
\draw[route,dashed] (census.south)--++(0,-4mm)-|(higher.north);
\draw[route] (atlas)--(higher);
\draw[route] (atlas.south west)--(eight.north east);
\draw[route] (eight)--(end);
\draw[route] (higher)--(end);
\end{tikzpicture}
\caption{The octic exclusion. The division into eight states and three
higher-delta states is a division of the exhaustive census, not an
assumption that the first eight conductor schemes are reduced. Solid
arrows indicate inputs; dashed arrows indicate the allocation of cases.}
\label{fig:octic-proof}
\end{figure}
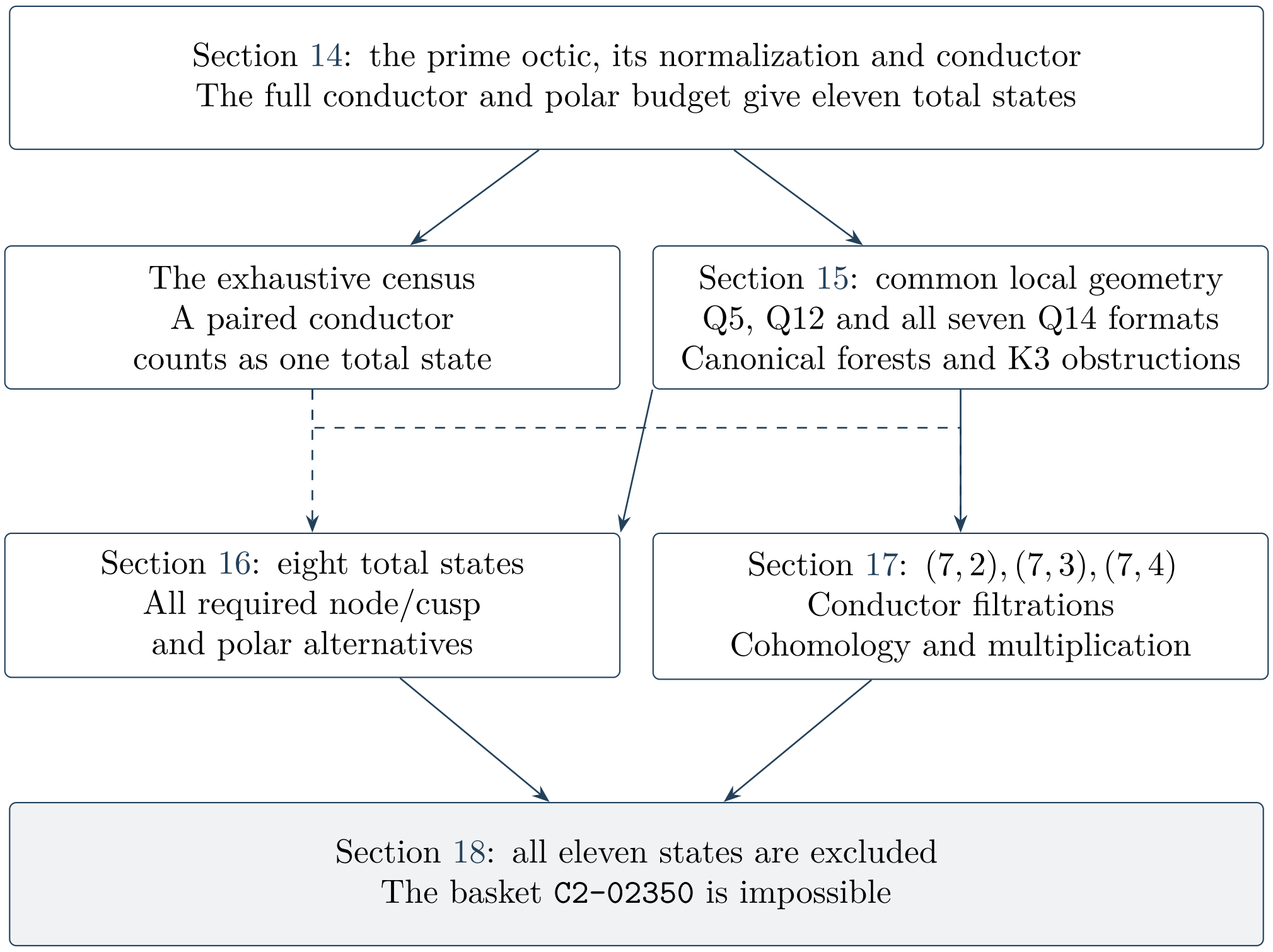

We first establish this reduction.  Section~\ref{sec:octic-atlas-revised}
provides the local surface calculations used throughout the exclusions.
We treat the eight states outside higher delta in degree seven in
Section~\ref{sec:octic-eight-states-revised}, and keep the full conductor
algebras for the remaining three states in
Section~\ref{sec:octic-higher-delta-revised}.  Figure~\ref{fig:octic-proof}
records these dependencies.

\subsection{C2-02350: the prime octic}
\label{sec:octic-prime46}

We first prove that the unique degree-eight member is reduced.  The
only alternative is twice a degree-four prime, and its associated
quadratic cover gives the following small auxiliary case.

\begin{lemma}\label{lem:octic-quadratic-cover46}
There is no canonical threefold with ample canonical divisor,
\[
 V=1/210,\quad\chi=3,\quad p_g=q=0,\quad
 (P_2,\ldots,P_7)=(0,0,1,0,0,1),\quad P_{13}=0.
\]
\end{lemma}
\begin{proof}
The low plurigenera give $\delta\ge8$ directly. The complete packing
calculation of Section~\ref{sec:equality-census}, with the same
intrinsic index bound, leaves the two explicit baskets
\eqref{eq:cover-baskets46}. Both have $I=420$, $v=2$ and $P_3=0$.
There are two generically independent holomorphic two-forms:
the alternative would contradict the $1/14$ bound in
\cite[Lemma~2.7]{ChenJiangForms}. The two-form argument of
Section~\ref{sec:chi3} supplies two small visible primes of scaled
degrees at most $\lfloor2v/3\rfloor=1$. Already the individual
degree-one necessary list is empty in each basket. Explicitly,
testing every tuple $i_Q\in\mathbb Z/r_Q$ in
\eqref{eq:degree-character46}--\eqref{eq:square-character46} with
$n=1$ and $-1\le M<1/2$ gives no solution. These are finite
congruences modulo $420$ and $840$, with the tuples grouped only
under permutations of identical constituents. Thus neither basket
can occur.
\end{proof}

\begin{proposition}\label{prop:octic-prime46}
In a hypothetical realization of C2-02350 the unique member
$D=(a_8=0)$ is a reduced prime, visible to the ambient two-form.
Its prime strict transform $T$ on a crepant terminal model has
\[
 (IK^2T,IKT^2)=(8,64),\qquad I=420.
\]
Any remainder in the full terminal divisor of $a_8$ is supported
on exceptional divisors with point centres.
\end{proposition}
\begin{proof}
We use the plurigenera listed in \eqref{eq:octic-small-P46}.
All four necessary opposite torsion states have $H_2=H_3=0$.
Squaring excludes every effective degree-one torsion translate.
The visible RR list through scaled degree eight is precisely
$(n,M)=(n,n^2)$, $1\le n\le8$. Null descent and the preceding
vanishings give a visible-prime floor of four. Visible saturation
(Lemma~\ref{lem:visible-saturation34}) exhausts the full divisor
on $Z$. Its complete dualizing space has dimension
$P_9-P_1+h^{2,0}=1$, so there is only one distinct visible prime.
Thus $D=T_8$ or $D=2T_4$.

In the latter alternative $T_4\sim4K_Z+\tau$ and the actual
divisor identity gives $2\tau=0$. Since $P_4=0$, its order is
exactly two. Of the four necessary states exactly two have
$H_4(\tau)=H_4(-\tau)>0$. For the first one, complete twisted
prime restriction has $G_5(\tau)=1$ and $G_{10}(0)=0$, excluding
it by squaring. The remaining state has
\[
 C(\tau)=2,\quad H_4(\tau)=H_7(\tau)=1,\quad
 H_2=H_3=H_5=H_6=H_{11}=H_{13}=0.
\]
The opposite Euler formula gives $\chi(\mathcal O_Z(\tau))=1$.
Also $H_1(\tau)=0$ by squaring, and nontrivial torsion has no
degree-zero section. Therefore $h^2(\tau)-h^1(\tau)=1$.
There is a canonical degree-ten section coprime to $T_4$, since
its restriction has dimension $P_{10}-H_6(\tau)=1$.
The two-twist Koszul inequality applied to this actual pair gives
\[
 h^2(\tau)\le P_{15}-H_{11}(\tau)-P_5+H_1(\tau)=1.
\]
Thus $h^2(\tau)=1$, $h^1(\tau)=0$.

The connected quasi-\'etale cyclic cover associated to $\tau$ is
canonical, with $\pi_*\mathcal O_{Z'}=\mathcal O_Z\oplus
\mathcal O_Z(\tau)$ and $K_{Z'}=\pi^*K_Z$. Its invariants are
exactly those of Lemma~\ref{lem:octic-quadratic-cover46}, by the
displayed section dimensions. This is a contradiction.
The complete packing census, rather than the truncated plurigenus
sequence, supplies the two possible cover baskets used here.

Hence $D$ is a reduced visible prime. Its terminal strict transform
has degree eight, and the individual RR list gives square $64$.
In the full equality $8K\sim T+E$, the effective exceptional
remainder therefore has $K^2E=KE^2=0$. Restricting to a general
canonical hyperplane and using the negative definite exceptional
curve matrix shows that no component of $E$ has a curve centre.
\end{proof}

\begin{proposition}\label{prop:octic-sources46}
For C2-02350 the ordinary terminal points of indices
$3,5,12,14$ have unchanged canonical neighbourhoods. Every other
canonical index divides four. In particular $420K_Z$ is Cartier.
\end{proposition}
\begin{proof}
This is the row 02350 application of
Proposition~\ref{prop:actual-source-transfer46}. Its fifth
constituent is the singleton $(2,5)$, its third constituent is a
singleton, and it contains no $(1,6)$. The mandatory-subbasket
test leaves only indices two and four for non-Gorenstein
point-centred crepant divisors. The index-twelve and index-fourteen
points are unchanged by the high-index test, and the small-index
argument then applies to the third and fifth points. Index
transfer accounts for every remaining point, giving the stated
list.
\end{proof}


\subsection{Complete restriction and the K3 normalization}
\label{sec:octic-foundation46}

Throughout this section let $Z$ be a canonical model with
\begin{equation}\label{eq:octic-basket46}
 K_Z^3=\frac1{420},\qquad \chi(\mathcal O_Z)=2,\qquad
 B=5(1,2)+(1,3)+2(1,4)+(2,5)+(5,12)+(5,14).
\end{equation}
Here $B$ is the basket of a crepant terminalization, in the convention
$(b,r)=\frac1r(1,-1,b)$.  The preceding reductions give
$p_g(Z)=q(Z)=0$ and $h^0(Z,\Omega_Z^{[2]})=1$.  We fix a nonzero
reflexive two-form $\alpha$.  We use the reduced-prime conclusion of
Proposition~\ref{prop:octic-prime46} and the actual-source conclusion of
Proposition~\ref{prop:octic-sources46}.  Thus the unique octic section
$a\in H^0(Z,8K_Z)$ has a reduced integral divisor $D$, and the
canonical contraction has the unchanged ordinary quotient neighborhoods
\begin{equation}\label{eq:octic-sources46}
 Q_3:\tfrac13(1,2,1),\quad Q_5:\tfrac15(1,4,2),\quad
 Q_{12}:\tfrac1{12}(1,11,5),\quad Q_{14}:\tfrac1{14}(1,13,5).
\end{equation}
Every other non-Gorenstein canonical-model point has index two or four.
The full exceptional remainder of the octic on a terminalization is
point-centred.  In particular $Z$ is smooth at the generic point of every
curve in the nonnormal locus of $D$.

Write $\nu:N\to D$ for the finite normalization.  The restriction of
$K_Z$ defines an integral Weil class $H$ on $N$, which is ample and
$\mathbb Q$-Cartier.  We also use $H$ for the corresponding canonical
restriction on $D$; the meaning of every sheaf on $D$ is the complete
character restriction, rather than the quotient of global section rings.
Put
\[
 T_j=H^0(D,\mathcal O_D(jH)),\qquad
 S_j=H^0(N,\mathcal O_N(jH)).
\]
On a canonical index cover, $\mathcal O_D(jH)$ is the appropriate
character summand of the hypersurface quotient by $a$.  These sheaves are
Cohen--Macaulay and have rank one on $D$.  Adjunction and extension across
finitely many points give
\begin{equation}\label{eq:octic-adjunction46}
 \omega_D\simeq\mathcal O_D(9H),\qquad H^2=\frac8{420}=\frac2{105}.
\end{equation}
In particular adjunction here is an isomorphism of actual sheaves.

\paragraph{Complete restriction and the twelfth pencil}

Substituting \eqref{eq:octic-basket46} in Reid's formula
\eqref{eq:RR46} gives
\begin{equation}\label{eq:octic-small-P46}
\begin{array}{c|rrrrrrrrrrrrrrr}
 m&2&3&4&5&6&7&8&9&10&11&12&13&14&15&18\\\hline
 P_m&0&0&0&0&0&0&1&0&1&1&2&0&2&1&2.
\end{array}
\end{equation}
The complete restriction sequence and canonical vanishing yield
\begin{equation}\label{eq:octic-complete-cohomology46}
 h^0(\omega_D)=1,\quad h^1(\mathcal O_D)=0,\quad
 h^2(\mathcal O_D)=1,\quad \dim T_{13}=0,\quad\dim T_{18}=1.
\end{equation}
For clarity, $\dim T_9=1$ even though $P_9-P_1=0$: the extra
section comes from $H^1(Z,K_Z)$, which has dimension one.  The actual
section is the restriction of the ambient two-form under adjunction.
Locally write $\alpha=\iota_\delta\eta$, where $\eta$ is a canonical
frame.  Its coefficient in the residue generator is
\[
                         \theta=\delta a\pmod a.
\]
It is a nonzero section of $\mathcal O_D(9H)$.  The frame-change terms
are multiples of $a$, so this local expression defines the stated global
section, including on canonical covers.

Let $W$ resolve $N$ and the rational map induced by the two twelfth
sections.  The pullback of $\alpha$ is a nonzero holomorphic two-form on
$W$: pullback to a resolution of the canonical threefold followed by
restriction proves regularity, and visibility proves nonvanishing.
Dualizing trace gives
\[
 H^0(W,K_W)\hookrightarrow H^0(N,\omega_N)
                  \hookrightarrow H^0(D,\omega_D).
\]
Consequently
\begin{equation}\label{eq:octic-pg46}
                              p_g(W)=1.
\end{equation}
This conclusion does not assume rational singularities of $N$.

\begin{proposition}\label{prop:octic-connected-pencil46}
The restricted twelfth pencil has connected general fibre.  If $C$ is its
smooth general fibre on $W$, then
\[
                 49\le \ell=420H\cdot C\le96,
                    \qquad g(C)\in\{1,2\}.
\]
\end{proposition}
\begin{proof}
Normal division by the prime $D$ shows that the twelfth restrictions are
independent: their kernel is $aH^0(Z,4K_Z)=0$.  Resolve their base locus
and take the Stein factorization.  If its connected general fibre is $C$
and its finite map to the parameter line has degree $d$, positivity gives
$0<d\ell\le96$.  The existence of a nonzero two-form excludes a ruled
surface, hence $g(C)\ge1$.

The canonical restrictions of nonzero sections in degrees
$10,11,12,14,15$ are nonzero; their potential division kernels are,
respectively, the zero spaces in degrees $2,3,4,6,7$.  Together with
$\theta$ they give nonzero sections $a_m$ of weights
$m=9,10,11,12,14,15$ on a general $C$.  Here $a_{12}$ is one of the two
pencil generators.  At a normalized branch over a non-Gorenstein source,
let $r$ be its inertia order and let $k/r$, $0<k<r$, $(k,r)=1$, be the
canonical marking: a weight-$m$ numerator has order congruent to $mk$
modulo $r$.  The possible $r$ are
$2,3,4,5,6,7,12,14$.  The marked adjunction lemma
\cite[Lemma~10.2]{LBVolume}, applied to $D$ and a general
twelfth member, and the effective divisors of these six sections give
\begin{equation}\label{eq:octic-pencil-adjunction46}
 2g(C)-2+\sum_x(1-r_x^{-1})\le\frac{21\ell}{420},\qquad
 \frac\ell{420}\equiv\sum_x\frac{k_x}{r_x}\pmod{\mathbb Z},
\end{equation}
\begin{equation}\label{eq:octic-pencil-sections46}
           \sum_x\frac{\overline{mk_x}}{r_x}
                      \le\frac{m\ell}{420}
       \quad(m=9,10,11,12,14,15).
\end{equation}
All residues in this last formula are nonnegative, including zero.
The hypotheses of marked adjunction hold at the generic point of a moving
intersection curve: there $D$ is smooth and a general pencil member meets
it simply.  Pulling back a Cartier multiple verifies that the section
divisors used in \eqref{eq:octic-pencil-sections46} are effective even
before rationality of $N$ has been established.

Suppose $\ell\le48$.  There is a mark, because $0<\ell/420<1$.
The first inequality makes genus at least two impossible and permits at
most four marks.  Exhausting the finite list of primitive multisets gives
exactly the following possibilities.  In the right column a monomial
ratio has weight zero and is an honest rational function on $C$.
\begin{equation}\label{eq:octic-pencil-finite-list46}
\begin{array}{c|l|l}
\ell&\text{markings}&\text{ratio or remaining local case}\\\hline
30&1/2,4/7&a_9^{11}/a_{11}^{9}\\
30&1/14&Q_{14}\\
35&1/2,7/12&a_9^{11}/a_{11}^{9}\\
35&1/3,3/4&a_{10}^{7}/a_{14}^{5}\\
35&2/3,5/12&a_{11}^{14}/a_{14}^{11}\\
35&1/4,5/6&a_9^{14}/a_{14}^{9}\\
35&1/6,11/12&a_9^{10}/a_{10}^{9}\\
35&1/12&Q_{12}\\
36&4/5,2/7&a_9^4/a_{12}^3\\
40&2/3,3/7&a_{11}^5/(a_9^5a_{10})\\
40&1/6,13/14&a_9^{10}/a_{10}^{9}\\
42&1/2,3/5&a_9^{11}/a_{11}^{9}\\
45&1/4,6/7&a_9^5/a_{15}^3\\
45&3/4,5/14&Q_4,Q_{14}\\
48&2/5,5/7&a_{11}^3a_{12}/a_9^5.
\end{array}
\end{equation}
For each displayed ratio the section budgets it uses are equalities;
there are no remaining integral zeros.  Substitution of the residues
shows that its divisor is the difference of the two marked points.
This is a degree-one rational function on a genus-one curve, a
contradiction.  The finite calculation is specified completely by
\eqref{eq:octic-pencil-adjunction46}--\eqref{eq:octic-pencil-sections46}:
choose one to four primitive marks in nondecreasing order, set
$\ell=420\sum k/r\bmod420$, and test the seven inequalities.  Thus the
table does not presuppose an earlier profile list.

For the other three cases use the actual ordinary source covers in
\eqref{eq:octic-sources46}.  At $\frac1r(1,-1,5)$ the positive coordinate
lower orders of a branch marked $k/r$ are
\[
             (-5^{-1}k,\,5^{-1}k,\,-k)\pmod r.
\]
Each twelfth generator has positive order at a mark of a general pencil
member.  When its source character is nonzero this is automatic.  When
it is zero, a unit numerator would make the quotient of the two
generators regular with a fixed value at the source, which a varying
general member cannot have.  In particular a zero-residue twelfth mark
costs at least one in its rational divisor.

For $(\ell;k/r)=(30;1/14)$ the coordinate lower orders are $(11,3,13)$;
the entire generator degree $6/7$ forces numerator order twelve.
The only character-ten monomial of positive weight at most twelve is
$y^4$.  For $(35;1/12)$ the lower orders are $(7,5,11)$ and the exact
numerator order is twelve; the sole positive character-zero monomial
of weight at most twelve is $xy$.  For $(45;3/4,5/14)$ the quarter mark
uses one whole unit of twelfth degree, leaving $2/7$ at $Q_{14}$.
The lower orders there are $(13,1,9)$; the only character-ten monomial
of weight at most four is $y^4$.

In each case the two generators have nonzero initial coefficients
$cM,c'M$ of the same sole monomial.  For a general parameter $\lambda$,
the member $t-\lambda s$ has nonzero initial coefficient
$(c'-\lambda c)M$ and cannot vanish along the alleged branch.  Higher
orders cannot cancel this term.  This excludes all three cases and
proves $\ell\ge49$, so $d=1$.

Finally $\ell\le96$ and \eqref{eq:octic-pencil-adjunction46} exclude
$g\ge3$: two marks cost at least one, giving a left side at least five
and a right side at most $24/5$.  With one mark its fraction is at most
$8/35$, so its inertia is at least five.  Equality would require inertia
five and fraction $8/35$, contrary to the degree congruence.  There cannot
be no mark.  Thus $g\le2$.
\end{proof}

\paragraph{A pluricanonical normalization argument}

\begin{lemma}[Full-conductor inclusion]\label{lem:octic-pluri-conductor46}
Let $D$ be an integral Cohen--Macaulay surface, Gorenstein in codimension
one, and let $\nu:N\to D$ be its finite normalization.  Let
$\mathcal L_m$ be a torsion-free $S_2$ sheaf agreeing with
$\omega_D^{\otimes m}$ outside finitely many points.  For $m\ge1$ there
is a natural injection
\[
                    \nu_*\omega_N^{[m]}\hookrightarrow\mathcal L_m.
\]
Consequently, for a resolution $W\to N$,
$H^0(W,mK_W)\hookrightarrow H^0(D,\mathcal L_m)$.
\end{lemma}
\begin{proof}
At a height-one point let $A$ be the one-dimensional Gorenstein local
ring of $D$ and let $B$ be its semilocal finite normalization.  In the
common total fraction field finite duality gives
\[
 \omega_B=\operatorname{Hom}_A(B,\omega_A)
          =\mathfrak c\omega_A,\qquad
 \mathfrak c=\{x:xB\subset A\}\subset A.
\]
Since $B$ is regular in dimension one and $\omega_A$ is invertible,
$\omega_B^{\otimes m}=\mathfrak c^m\omega_A^{\otimes m}
\subset\omega_A^{\otimes m}$.  This is multiplication of conductor
ideals, not a multiplicativity assertion about a trace map.
These inclusions agree in the function field and define a map off a
finite set.  Hartogs for the $S_2$ target extends the map: the images of
free generators in a finite presentation extend, and their relations
continue to vanish because the target is torsion-free.  The map remains
injective by its generic injectivity.  A regular pluricanonical form on
$W$ restricts on the smooth locus of $N$ and hence defines a section of
$\omega_N^{[m]}$, proving the last assertion.
\end{proof}

\begin{theorem}\label{thm:octic-K3-rational46}
The smooth minimal model of $W$ is a K3 surface.  The normal surface $N$
has rational singularities.  Its normalization quotient
$\mathcal Q=\nu_*\mathcal O_N/\mathcal O_D$ is acyclic:
\[
                         H^0(D,\mathcal Q)=H^1(D,\mathcal Q)=0.
\]
\end{theorem}
\begin{proof}
Apply Lemma~\ref{lem:octic-pluri-conductor46} with
$\mathcal L_2=\mathcal O_D(18H)$.  Equation
\eqref{eq:octic-complete-cohomology46} gives $P_2(W)\le1$; the square of
the nonzero canonical form gives equality.

By Proposition~\ref{prop:octic-connected-pencil46}, $W$ has a connected
fibration over $\mathbb P^1$ with genus-one or genus-two general fibre.
The relative irregularity bound gives $q(W)\le g(C)$.  In the genus-two
case equality would make $W$ birational to $\mathbb P^1\times C$, by
Beauville's lemma \cite[p.~345]{BeauvilleDebarre1982}; this has geometric
genus zero and is impossible.  Therefore $q(W)\le1$, and
$\chi(\mathcal O_W)=2-q(W)\ge1$.

Let $W_0$ be its smooth minimal model.  If $\kappa(W_0)=2$, surface
Riemann--Roch and vanishing imply
$P_2(W_0)=\chi(\mathcal O_{W_0})+K_{W_0}^2\ge2$, a contradiction.
If $\kappa(W_0)=1$, use its relatively minimal genus-one Iitaka
fibration $f:W_0\to B$.  Put $b=g(B)\le q(W_0)\le1$.  The characteristic-zero canonical
bundle formula, including multiple fibres, gives
\[
 f_*\omega_{W_0}^{\otimes2}
  \simeq (\omega_B\otimes L)^{\otimes2}
      \Bigl(\sum_i\bigl\lfloor2(1-1/m_i)\bigr\rfloor P_i\Bigr),
 \qquad \deg L=\chi(\mathcal O_{W_0}),
\]
where $m_i>1$ are the multiple fibre multiplicities.  If $b=1$, then $q=1$ and $\chi=1$, so
$\deg(\omega_B\otimes L)=1$ and $P_2\ge2$ by Riemann--Roch on the
elliptic base.  If $b=0$, the first-power formula gives
$1=p_g=h^0(\mathbb P^1,\mathcal O(\chi-2))$, forcing $\chi=2$
and $q=0$.  The ordinary base term then has degree zero.
Proper ellipticity requires a multiple fibre; each such fibre contributes
one to the displayed second-power divisor, so again $P_2\ge2$.  These are contradictions.
The canonical bundle formula and the classification are recalled in
\cite[\S13, Theorem 7; \S15, Theorem 5]{Peters1995}.

Positive geometric genus excludes Kodaira dimension $-\infty$.
In Kodaira dimension zero the possibilities with $p_g=1$ are K3 and
abelian surfaces; the latter has $q=2$ and is excluded.  Thus $W_0$ is
K3, and in particular $(p_g(W),q(W))=(1,0)$.

It remains to justify rationality, rather than assume it in the preceding
argument.  Both $D$ and $N$ are integral and projective, so their global
function spaces are $\mathbb C$.  The normalization sequence and
$H^1(D,\mathcal O_D)=0$ imply $H^0(D,\mathcal Q)=0$.
The quotient has dimension at most one, hence
$\chi(\mathcal Q)=-h^1(\mathcal Q)\le0$.
For a resolution $p:W\to N$, let
$e=\dim H^0(N,R^1p_*\mathcal O_W)\ge0$.
Leray gives $\chi(\mathcal O_N)=\chi(\mathcal O_W)+e=2+e$.
Since $\chi(\mathcal O_D)=2$, additivity gives
$\chi(\mathcal Q)=e\ge0$.  Therefore $e=0$ and $h^1(\mathcal Q)=0$.
The vanishing of $R^1p_*\mathcal O_W$ is precisely rationality for a
normal surface.  No earlier assertion about the Kodaira dimension or
irregularity of the normalization has entered the proof.
\end{proof}

\paragraph{The actual conductor and its cohomology}

Let $\mathfrak c$ be the full conductor ideal, viewed on both $D$ and $N$.
Set
\[
 C_D=\operatorname{Spec}_D(\mathcal O_D/\mathfrak c),\qquad
 C_N=\operatorname{Spec}_N(\mathcal O_N/\mathfrak c).
\]
The same symbol $C_N$ denotes its effective Weil divisor on $N$.
No reduction of either conductor scheme is made.  Let
$B_N=\operatorname{div}_N(\alpha_N)$, where $\alpha_N$ is the actual
restricted two-form.

\begin{proposition}\label{prop:octic-full-conductor46}
The conductor schemes are pure Cohen--Macaulay curves, and
\[
 H^1(C_D,\mathcal O_{C_D})=H^1(C_N,\mathcal O_{C_N})=0.
\]
Moreover $H^0(\mathcal O_{C_D})\simeq H^0(\mathcal O_{C_N})$.
Every reduced irreducible conductor image is a smooth rational curve;
every reduced union of them has arithmetic genus zero on each connected
component.  There is an actual sheaf isomorphism and divisor equality
\begin{equation}\label{eq:octic-full-polar46}
 \omega_N(C_N)\simeq\mathcal O_N(9H),\qquad
                  \operatorname{div}_N(\theta)=B_N+C_N.
\end{equation}
\end{proposition}
\begin{proof}
The conductor equals $\operatorname{Hom}_D(\nu_*\mathcal O_N,\mathcal O_D)$.
It has depth two.  Indeed if $x,y$ is a regular sequence on a local
Cohen--Macaulay surface ring $A$ and $y\phi=x\psi$ in such a Hom module,
then $\phi(b)$ is divisible by $x$ for every $b$, since $y$ is a
nonzerodivisor modulo $x$.  Division by $x$ is again a homomorphism.
This proves regularity of $y$ modulo $x$ on the Hom module.
The same parameter sequence proves depth two over the finite normalization.
Thus $\mathfrak c$ is rank-one reflexive on $N$, equal to
$\mathcal O_N(-C_N)$, and both quotient schemes are Cohen--Macaulay of
dimension one.  In particular there is no omitted zero-dimensional
saturation.

At the generic point of each divisor, finite Gorenstein duality gives
$\omega_N=\mathfrak c\nu^*\omega_D$.  Reflexivity extends this to
\eqref{eq:octic-full-polar46}, whose second equality follows from the
residue expression $\theta=\delta a$.  It also gives
\[
       \mathcal H om_N(\mathfrak c,\omega_N)
                              \simeq\mathcal O_N(9H).
\]
Applying $\mathcal H om_N(-,\omega_N)$ to the conductor sequence yields
\begin{equation}\label{eq:octic-conductor-adjunction46}
 0\longrightarrow\omega_N\longrightarrow\mathcal O_N(9H)
                    \longrightarrow\omega_{C_N}\longrightarrow0.
\end{equation}
Here $\omega_{C_N}=\mathcal E xt_N^1(\mathcal O_{C_N},\omega_N)$ is
the dualizing sheaf of the full, possibly nonreduced, curve.

For any $k$, multiplication by $\alpha_N$ followed by conductor trace
gives an injection $S_k\hookrightarrow T_{k+9}$.  It is an injection at
the generic point and extends by the $S_2$ property.  Thus $S_4=0$,
and taking powers gives $S_1=S_2=0$.  Also $\dim S_9\le\dim T_{18}=1$;
the nonzero section $\theta$ gives equality.  Rationality of $N$ and the
preceding theorem give $h^0(\omega_N)=1$ and $H^1(\omega_N)=0$.
The first map on global sections in
\eqref{eq:octic-conductor-adjunction46} is therefore an isomorphism.
It follows that $H^0(\omega_{C_N})=0$, hence $H^1(\mathcal O_{C_N})=0$
by curve duality.

The conductor square is the exact sequence
\[
 0\longrightarrow\mathcal O_{C_D}
      \longrightarrow\nu_*\mathcal O_{C_N}
                 \longrightarrow\mathcal Q\longrightarrow0.
\]
The acyclicity of $\mathcal Q$ gives the asserted $H^0$ isomorphism and
$H^1$ vanishing on $C_D$.  A reduced component or reduced union is a
quotient of $\mathcal O_{C_D}$; the kernel has dimension at most one,
so its $H^1$ vanishes.  For an integral component this says arithmetic
genus zero, which forces a smooth $\mathbb P^1$ by the normalization
sequence.  For a connected reduced union the same sequence says its
components are rational and its only gluings form a tree of ordinary
multibranch points (each local $\delta$ equals the number of branches
minus one).  This conclusion concerns the reduced union only and does
not remove any nilpotents from $C_D$ or $C_N$.
\end{proof}

\begin{proposition}[Twisted trace bounds]\label{prop:octic-twisted-trace46}
For a reduced conductor component $\Gamma$, put
$\mathcal F_{\Gamma,j}=(\mathcal O_D(jH)\otimes\mathcal O_\Gamma)/
\text{torsion}$.  Then
\begin{equation}\label{eq:octic-trace-lower46}
                 \deg\mathcal F_{\Gamma,j}\ge-1
               \quad\text{for }j=5,7\text{ and every }j\ge9.
\end{equation}
For $j=9$ the full canonical trace quotient
$\mathcal M_9=\omega_D/\nu_*\omega_N$ is acyclic.
\end{proposition}
\begin{proof}
For any $j$ conductor trace defines
\[
 \mathcal A_j=
 \nu_*\bigl(\omega_N\mathbin{\widehat\otimes}\mathcal O_N((j-9)H)\bigr)
       \hookrightarrow\mathcal O_D(jH),\qquad
 \mathcal M_j=\mathcal O_D(jH)/\mathcal A_j,
\]
where $\widehat\otimes$ means reflexive tensor product.  The inclusion is
first defined where the canonical class is Cartier and then extended
across finitely many points by Hartogs.  The quotient is not saturated;
all finite-length defects remain present.

The complete restriction sequence gives
$H^1(D,\mathcal O_D(jH))=0$ for the stated $j$.
For $j>9$, Serre duality identifies $H^2(D,\mathcal A_j)^*$ with
$H^0(N,\mathcal O_N((9-j)H))=0$, since $H$ is ample.
For $j=7,5$ the same group is respectively $S_2,S_4$, already zero.
Thus $H^1(\mathcal M_j)=0$.  For $j=9$, dualizing the normalization
sequence gives
$\mathcal M_9=\mathcal E xt_D^1(\mathcal Q,\omega_D)$.
The sheaf $\mathcal Q$ is Cohen--Macaulay of dimension one by the depth
lemma.  Proper duality and its acyclicity give
$H^0(\mathcal M_9)=H^1(\mathcal M_9)=0$.

At the generic point of $\Gamma$, the trace image is the conductor ideal
multiplied by $\omega_D$; it vanishes on the reduced conductor.  Therefore
the composite $\mathcal A_j\to\mathcal F_{\Gamma,j}$ is generically
zero and hence zero, since its target is torsion-free on the integral
curve.  Thus $\mathcal M_j$ surjects onto $\mathcal F_{\Gamma,j}$.
Its kernel has dimension at most one, so $H^1(\mathcal F_{\Gamma,j})=0$.
As $\Gamma\simeq\mathbb P^1$ and $\mathcal F_{\Gamma,j}$ is a line
bundle, this is precisely \eqref{eq:octic-trace-lower46}.
\end{proof}


\subsection{Nonnormality and the intrinsic polar divisor}
\label{sec:octic-normality-polar46}

\paragraph{A local defect under rational deformations}

For a rational normal surface germ $(V,P)$ and its minimal resolution
$\pi:\widetilde V\to V$, define $\Delta$ by
$\pi^*K_V=K_{\widetilde V}+\Delta$, and put
\[
 r(V,P)=\#\{\text{exceptional components}\},\qquad
 c(V,P)=-\Delta^2,\qquad d(V,P)=r(V,P)-c(V,P).
\]
The definition of $\pi^*K_V$ here is the numerical pullback; for a
$\mathbb Q$-Gorenstein germ it agrees with the usual pullback.
The defect $d$ need not be nonnegative for an arbitrary rational
singularity.  We shall justify nonnegativity only through a specified
deformation.

\begin{lemma}\label{lem:octic-defect-specialization46}
Consider a same-character quotient-hypersurface deformation at an
ordinary cyclic source.  Suppose its central surface germ is normal and
rational and its sufficiently small partial Milnor fibre $M_t$ is normal
with rational singularities.  Then
\begin{equation}\label{eq:octic-defect-specialization46}
       d(V_0,P)\ \ge\ \sum_{Q\in\operatorname{Sing}M_t}d(V_t,Q).
\end{equation}
The same inequality holds for a suitable displacement of any normal
rational Cartier surface germ in the actual canonical threefold $Z$;
the displacement can be chosen to have only Du Val singularities.
For these two kinds of families the inequality follows from
\begin{equation}\label{eq:octic-local-Noether46}
 d(V_0,P)-\sum_Qd(V_t,Q)=\chi_{\rm top}(M_t)-1.
\end{equation}
\end{lemma}
\begin{proof}
First establish the identity by relative Noether.  The indicated families
can be compactified so that, after resolving the singularities that
persist away from the distinguished germ, their complements outside a
Milnor neighborhood are smooth over the parameter disc.  One construction is to compactify
the ambient quotient, resolve its boundary away from the germ, and extend
the character divisor after twisting by a sufficiently positive divisor
supported on the boundary.  Prescribing a sufficiently high finite jet
retains the isolated hypersurface germ; equivariant finite determinacy
follows by averaging the infinitesimal coordinate changes over the finite
group.  General choices are smooth off the prescribed neighborhood.
For a displaced Cartier slice the ambient source need not be a
hypersurface.  Instead keep the actual projective threefold $Z$ and
prescribe the slice in a sufficiently high jet at the chosen point $P$
using a sufficiently positive Cartier line bundle.  Here this jet
prescription preserves the local analytic Cartier germ.  Indeed, with
$A=\mathcal O_{Z,P}^{\mathrm{an}}$ and equation $f$, the ideal
\[
             (f,\{v(f):v\in\operatorname{Der}_{\mathbb C}(A)\})
\]
is primary to the maximal ideal: off $P$, either $f$ is a unit, or
$Z$ and the Cartier surface are both smooth, so some derivative is a
unit.  To justify the latter assertion, a singular threefold local ring
has embedding dimension at least four, so quotienting by one element
cannot give a regular surface ring; normality of the surface therefore
makes its intersection with $\operatorname{Sing}Z$ isolated.
Here is the finite-determinacy argument explicitly.  Choose generators
$v_i$ of $\operatorname{Der}_{\mathbb C}(A)$ and write
$J=(f,v_1(f),\ldots,v_s(f))$.  Choose $k$ with $\mathfrak m^k\subset J$.
For $h\in\mathfrak m^{k+\ell}$, $\ell\ge2$, one can express
$h=bf+\sum a_i v_i(f)$ with $b,a_i\in\mathfrak m^\ell$.
Moreover $v_i(h)\in\mathfrak m^{k+\ell-1}\subset
\mathfrak m^{\ell-1}J$.  Hence, for $f_t=f+th$, the vectors of generators
$(f_t,v_1(f_t),\ldots,v_s(f_t))$ and $(f,v_1(f),\ldots,v_s(f))$
are related by a matrix $I+tM$ with entries of $M$ in
$\mathfrak m^{\ell-1}$.  This matrix is analytically invertible after
shrinking the representative, uniformly for $t$ on the unit interval.
It follows that
$h=b_t f_t+\sum a_{i,t}v_i(f_t)$ with analytic coefficients
$b_t,a_{i,t}\in\mathfrak m^\ell$.  Integrating the time-dependent
analytic derivation $-\sum a_{i,t}v_i$ removes $h$, while integrating
$b_t$ supplies an invertible multiplier.  The derivations preserve the
ambient analytic ideal, so the flow preserves $Z$.  This proves that the
prescribed finite jet determines the Cartier germ without an LCI
assumption on $Z$.
A general section satisfying the prescribed jet is transverse to the
smooth parts of the singular curves of $Z$ and avoids all other isolated
bad points.  Its remaining singularities are Du Val slices.  These vary
locally trivially outside the chosen Milnor neighborhood (and admit
simultaneous resolutions after a finite base change).
Choose the second pencil section to be a unit at $P$.  Its local
quotient with the central section is a defining equation of the same
Cartier germ, up to the analytic change of coordinates and invertible
multiplier just constructed.  We use the displacement of this defining
equation.  Only the central analytic germ and the existence of a
displaced normal rational slice are needed in the inequality; no
equivalence of arbitrary prescribed deformation directions is asserted.
The complement contributes the same terms on both sides.
In both constructions the projective family is flat and
$\mathbb Q$-Gorenstein, so $\chi(\mathcal O)$ and $K^2$ of its surface
fibres are constant.  Rational singularities do not change $\chi(\mathcal O)$
upon resolution.  Replacing a rational surface singularity by its
exceptional tree changes Euler characteristic by its number of vertices.
Subtracting Noether's formulas for the resolved fibres gives
\eqref{eq:octic-local-Noether46}.  This is the relative version of the
usual local Noether identity for surface smoothings; see also
\cite[\S\S1,5]{Wahl2015}.

Here the right side is nonnegative, including when $M_t$ is still
singular.  Let $\widetilde M_t$ resolve its rational singularities.
The Stein surface $M_t$ has the homotopy type of a real CW complex of
dimension at most two, including when it is singular
\cite[Satz 1]{Hamm1983}.  Resolving its isolated rational points replaces
points by trees of two-spheres and introduces no third homology.  Thus
$H_3(\widetilde M_t,\mathbb Q)=0$.
Its boundary is the link of the rational germ $V_0$, hence is a rational
homology sphere.  Poincare--Lefschetz duality and the homology sequence
of this boundary imply $H_1(\widetilde M_t,\mathbb Q)=0$.
Collapsing the exceptional trees preserves first homology, so
$H_1(M_t,\mathbb Q)=0$.  Consequently
$\chi_{\rm top}(M_t)=b_0(M_t)+b_2(M_t)\ge1$.
\end{proof}

\begin{lemma}\label{lem:octic-global-defect46}
If $D$ were normal, then
\begin{equation}\label{eq:octic-global-defect46}
                     K_D^2+\sum_Pd(D,P)\le19.
\end{equation}
\end{lemma}
\begin{proof}
Theorem~\ref{thm:octic-K3-rational46} gives rational singularities and
$(p_g,q)=(1,0)$ on a resolution $W$.  The $r=\sum_Pr(D,P)$ exceptional
classes together with $\pi^*K_D$ are independent, because the former
have negative-definite intersection matrix and the latter has positive
square and is orthogonal to them.  Thus $r+1\le h^{1,1}(W)$.
Noether and Hodge give $h^{1,1}(W)=20-K_W^2$, while
$K_W^2=K_D^2-\sum_Pc(D,P)$.  Substitution gives the inequality.
\end{proof}

\begin{proposition}\label{prop:octic-nonnormal46}
The prime octic $D$ is nonnormal.
\end{proposition}
\begin{proof}
Assume it is normal.  The four ordinary source neighborhoods of
\eqref{eq:octic-sources46} give local quotient-hypersurface equations
of characters $1,4,8,2$, respectively.  Adding general terms of the
same character gives the following rational normal germs:
\begin{equation}\label{eq:octic-normal-generic-graphs46}
\begin{array}{c|l|l|c|c|c}
Q&\text{covering equation}&\text{resolution graph}&r&c&d\\\hline
3&x&A_2&2&0&2\\
5&y&[3,2]&2&2/5&8/5\\
12&x^3z+y^4+b y^2z^2+c z^4&[4;A_2,A_2,A_1,A_5]&11&3&8\\
14&x^2+y^3z+z^6&[4;A_2,A_2,A_4]&9&15/7&48/7.
\end{array}
\end{equation}
Here $[b;\cdots]$ means a central curve of square $-b$ with the displayed
arms.  The parameters in the quartic equation are general.
We include the quotient calculation to specify the graphs, rather than
infer them from the basket alone.

At $Q_{12}$ the four displayed quartic monomials are all the lowest
ordinary-degree monomials of character eight.  A general member defines
a smooth plane quartic $E$ of genus three.  After the ordinary blowup,
$\mu_2\subset\mu_{12}$ is radial; the residual group $\mu_6$ acts on
$E$ with quotient $\mathbb P^1$ and stabilizers $3,3,2,6$.
In the $z$-chart write $(x,y,z)=(su,sv,s)$ and then $S=s^2$.
The residual characters of $(S,u,v)$ are $(5,4,3)$ modulo six.
At $u=0$ the two orbits have local weights $(2,1)$ modulo three;
at $v=0$ the weights are $(1,1)$ modulo two.  In the $x$-chart the
fixed point has weights $(1,5)$ modulo six.  These give the stated
$A_2,A_2,A_1,A_5$ arms.  The pre-resolution central square is
$-8/6=-4/3$; subtracting
$2/3+2/3+1/2+5/6=8/3$ makes it $-4$.

At $Q_{14}$ use weights $(9,5,3)$, with equation degree eighteen.
The ambient action is weighted radial scaling: its character vector is
$11(9,5,3)$ modulo fourteen.  The weighted curve has two points of
stabilizer three and one of stabilizer five.  In the $z$-chart, after
quotienting by the radial group, the residual local weights are $(2,1)$
modulo three, and in the $y$-chart they are $(4,1)$ modulo five.
Hence the arms are $A_2,A_2,A_4$.
The pre-resolution central square is
$-14\cdot18/(9\cdot5\cdot3)=-28/15$;
subtracting $2/3+2/3+4/5=32/15$ gives $-4$.

For each graph, let $M$ be the positive definite matrix having diagonal
entries the indicated positive self-intersection numbers and off-diagonal
entries $-1$ along edges.  Adjunction gives
$M\Delta=(b_i-2)_i$.  Solving gives the values of $c=\Delta^TM\Delta$
in the table.  The reduced exceptional cycle is anti-nef and has
arithmetic genus zero, so these singularities are rational by the
rationality criterion.  

Rationality is open in these deformations \cite[Theorem IV]{Elkik1978}.  Equivalently one can choose
a sufficiently small general parameter in the open rational locus
containing the given rational central germ.  Lemma
\ref{lem:octic-defect-specialization46} therefore bounds each actual
source defect below by its entry in the table.

At every other canonical-model point the index divides four, so the
octic is Cartier.  Choose the suitable local defining equation $f$
supplied by Lemma~\ref{lem:octic-defect-specialization46} and displace
it to $f+t$.
A general displaced slice avoids isolated exceptional ambient points and
meets only the generic Du Val singularities along the ambient singular
curves; thus it is smooth or has Du Val surface singularities.
This is a $\mathbb Q$-Gorenstein deformation.  The same lemma shows that
its original defect is nonnegative.  This treats all additional points
without asserting that their index-one covers are smooth.

Finally adjunction gives $K_D^2=81H^2=54/35$, and the table gives
\[
 K_D^2+\sum_Pd(D,P)\ge
 \frac{54}{35}+2+\frac85+8+\frac{48}{7}=20.
\]
This contradicts Lemma~\ref{lem:octic-global-defect46}.
\end{proof}

\paragraph{The strict global conductor budget}

For a reduced conductor image $\Gamma$, let
$w_\Gamma=420K_Z\cdot\Gamma$ and let $\delta_\Gamma$ be its transverse
plane-curve delta invariant, including residue-field degrees.  These
numbers are positive integers.  Finite duality at its generic point gives
\begin{equation}\label{eq:octic-conductor-degree46}
                   \nu_*C_N=\sum_\Gamma2\delta_\Gamma\Gamma.
\end{equation}
For example this is the familiar equality between twice the delta
invariant and conductor colength on the normalization of a reduced
Gorenstein curve; extending scalars to an algebraic closure proves the
statement with all branch degrees included.
The actual polar equality \eqref{eq:octic-full-polar46} gives
\begin{equation}\label{eq:octic-polar-mass46}
                 420H\cdot(B_N+C_N)=72.
\end{equation}

\begin{lemma}\label{lem:octic-polar-gap46}
Every irreducible component of the actual intrinsic polar divisor $B_N$
has positive integral scaled $H$-degree at least three.
\end{lemma}
\begin{proof}
Its finite image on $D$ is an actual positive polar curve $\Gamma$.
The complete restriction in weight eighteen is surjective, because
$H^1(Z,10K_Z)=0$.  Thus $\theta^2$ has an ambient lift
$u\in H^0(Z,18K_Z)$, and $D$ is not a component of $(u)$.
At a smooth generic point of $D$, if the polar has order $j\ge1$, this
actual square lift makes $\Gamma$ a generic lc centre of
$(Z,D+(u)/(2j))$.  The adjunction coefficient is
$1+8+18/(2j)\le18$.  At a singular generic point use the generic lc
threshold $c\le1$ of $D$, whose coefficient $1+8c$ is at most nine.  The point-centred remainder permits the application
on a terminalization.  The marked adjunction lemma
\cite[Lemmas~10.2 and~10.5]{LBVolume} gives
\[
 2g(\Gamma^\nu)-2+\sum(1-1/r)\le18w/420,\qquad
             w/420\equiv\sum k/r\pmod{\mathbb Z}.
\]
For $w=1,2$ there are at most four marks, since each costs at least $1/2$.
There are no primitive multisets with denominators
$2,3,4,5,6,7,12,14$ satisfying these two conditions even with $g=0$.
This is the complete finite test: set $w=420\sum k/r\bmod420$ for every
multiset of size at most four and compare the two rational numbers.
Thus $w\ge3$.  A component upstairs has degree equal to this degree
multiplied by its positive integral generic mapping degree.
\end{proof}

\begin{lemma}\label{lem:octic-polar-nonzero46}
The divisor $B_N$ is nonzero.  Consequently
\begin{equation}\label{eq:octic-strict-budget46}
          \sum_\Gamma\delta_\Gamma w_\Gamma\le34,
                    \qquad 420H\cdot B_N\ge4.
\end{equation}
\end{lemma}
\begin{proof}
If $B_N=0$, the actual nonzero canonical form trivializes the reflexive
canonical sheaf of $N$.  A rational Gorenstein surface singularity is
Du Val, so $N$ would be a Du Val K3 model with integral Weil polarization
$H$ of square $2/105$.
On its minimal resolution the exceptional ADE root lattice has rank at
most nineteen.  Numerical pullback of the integral Weil divisor gives
\[
                        H^2-\sum_Pq_P\in2\mathbb Z,
\]
where on an $A_{r-1}$ component the possible local square corrections are
$q_P=k(r-k)/r$ modulo $2\mathbb Z$.
Only $A_{7a-1}$ root systems carry seven-primary denominators.  Within
rank nineteen a single such component is $A_6$ or $A_{13}$.
The required seven-primary congruence would respectively say
$k^2\equiv5$ or $3\pmod7$, both impossible.  At least two such
components, hence at least twelve root ranks, are needed.
At most seven ranks remain for the five-primary denominator, allowing
only one $A_4$; its congruence would say $k^2\equiv3\pmod5$, again
impossible.  The other ADE systems have discriminants supported only
at two and three.  This contradiction proves $B_N\ne0$.

By Lemma~\ref{lem:octic-polar-gap46}, its positive integral degree is at
least three.  Equations \eqref{eq:octic-conductor-degree46} and
\eqref{eq:octic-polar-mass46} make that degree even, so it is at least
four.  The asserted inequalities follow.
\end{proof}


\subsection{The eleven possible conductor states}
\label{sec:octic-conductor-reduction46}

This section derives the finite list of total conductor states from the
geometry just proved.  A state records the reduced conductor components,
their canonical degrees and markings, and their generic delta invariants.
It is a necessary configuration, not an existence assertion.

\paragraph{Full-source orders and the initial census}

At an ordinary source $\frac1n(1,-1,b)$ a branch of inertia $r\mid n$
marked $k/r$ has positive coordinate lower orders
\[
 (v_x,v_y,v_z)=(e,-e,be)\pmod r,
                       \qquad e=-kb^{-1}\pmod r.
\]
For a canonical numerator of weight $j$ define
\begin{equation}\label{eq:octic-source-min46}
 \lambda_{k/r}^{(n,b)}(j)=
 \min\{p v_x+qv_y+sv_z:p,q,s\ge0,\ p-q+bs\equiv-jb\pmod n\}.
\end{equation}
The full source character modulo $n$ is retained when $r<n$.
The actual torsion-free restriction line has degree bounded above by
\begin{equation}\label{eq:octic-source-degree46}
 U_j(\Gamma)=\frac{jw_\Gamma}{420}
                     -\sum_{k/r}\frac{\lambda_{k/r}^{(n,b)}(j)}r.
\end{equation}
To see this, trivialize the canonical line on each normalized cover
branch.  Each local numerator has order at least the displayed minimum;
passing to the coarse line subtracts that order divided by the inertia.
At all other points the possible image ideals add only nonnegative
integral vanishing.  Thus the formula is an upper bound, including any
hidden contact or singular lifted branch.  It is integral by the degree
congruence.  At unspecified marks of inertia $2,3,4$ initially use only
the fractional-character bound $\overline{jk}/r$; it does not presume
an ordinary local cover there.

Let $I_\Gamma$ denote the proper first-polar intersection length in the
transverse plane.  Generically $I_\Gamma\ge2\delta_\Gamma$.
On a conductor component $a=\delta a=0$, so
$\gamma=\delta^2a$ is a genuine section of weight ten.  If it vanishes,
$\eta=\delta^3a$ is a genuine section of weight eleven.
The elementary directional-polar estimates give
\begin{equation}\label{eq:octic-low-jet46}
 I_\Gamma\ge2;\qquad
 \gamma=0\Rightarrow I_\Gamma\ge3;\qquad
 \gamma=\eta=0\Rightarrow I_\Gamma\ge4.
\end{equation}
For completeness, if the derivation vanishes at the transverse origin,
its derivative of an equation of multiplicity $m\ge2$ has multiplicity
at least $m$, and the intersection length is at least four.
For a nonzero derivation straighten its flow.  Multiplicity at least
three already gives length at least six.  In multiplicity two, the
Weierstrass quadratic in the transverse variable shows that vanishing
of the second, then third directional jets raises the minimal polar
length from two to three, then four.  Equivalently, in the nodal form
$uv$ the vanishing second jet makes the constant direction tangent to
one branch; vanishing of the third kills its first transverse derivative,
raising that branch's contribution from two to three.  In the cuspidal
quadratic form the same two jet vanishings give the stated, or stronger,
bound.  Higher terms can increase these lengths but cannot reduce them.

The strict budget gives $1\le w_\Gamma\le34$.  At the generic
point of $\Gamma$, the ambient threefold is smooth and the reduced
transverse equation is singular.  Its log canonical threshold $c$ is
at most one, and threshold adjunction to $(Z,cD)$ has coefficient
$1+8c\le9$.  Hence the marked adjunction inequality of
\cite[Lemmas~10.2 and~10.5]{LBVolume} gives
\begin{equation}\label{eq:octic-conductor-marks46}
 -2+\sum(1-1/r)\le9w_\Gamma/420,
       \qquad w_\Gamma/420\equiv\sum k/r\pmod{\mathbb Z}.
\end{equation}
Thus there are at most five marks.  Generate all primitive multisets of
one to five marks of denominators $2,3,4,5,6,7,12,14$, set
$w=420\sum k/r\bmod420$, and apply
\eqref{eq:octic-conductor-marks46} and \eqref{eq:octic-low-jet46}.
If $U_{10}<0$, impose polar length at least three; if also $U_{11}<0$,
impose length at least four.  This gives exactly forty-seven rows.
The trace inequalities $U_j\ge-1$ for $j=5,7$ and $j\ge9$ leave
exactly twenty-one.  It suffices to check $9\le j\le428$, because
$\lambda(j+420)=\lambda(j)$ and $U_{j+420}=U_j+w$.
Two rows have two distinct marks forced to the same unique source
$Q_5,Q_{12}$ or $Q_{14}$; smoothness of the coarse conductor excludes
these.  The nineteen remaining rows are as follows.  The final columns
are the permissive full-source upper degrees; more precise source
assignments can only lower them.
\begin{equation}\label{eq:octic-nineteen46}
\begin{array}{r|r|l|rrr}
\text{row}&w&\text{markings}&U_9&U_{10}&U_{11}\\\hline
1&7&3/5,5/12&-1&0&-1\\
2&10&2/3,5/14&0&-1&-1\\
3&12&1/2,3/5,13/14&-1&0&-1\\
4&14&1/5,5/6&-1&0&-1\\
5&18&2/5,9/14&-1&0&0\\
6&21&1/2,3/4,4/5&-1&0&-1\\
7&21&1/3,4/5,11/12&0&0&-1\\
8&21&1/4,4/5&0&0&-1\\
9&24&1/5,6/7&-1&0&-1\\
10&25&1/2,11/12,9/14&-1&0&0\\
11&25&5/12,9/14&-1&0&0\\
12&26&1/3,4/5,13/14&0&0&-1\\
13&28&1/3,1/3,2/5&0&0&-1\\
14&28&2/3,2/5&0&0&0\\
15&28&3/4,2/5,11/12&-1&0&0\\
16&30&1/2,4/7&-1&0&0\\
17&31&4/5,11/12,5/14&0&0&-1\\
18&32&3/5,5/6,9/14&-1&0&0\\
19&33&3/4,2/5,13/14&-1&0&0.
\end{array}
\end{equation}
The local minimum in \eqref{eq:octic-source-min46} is independently
certified by a shortest-path calculation on $\mathbb Z/n$: the edges
have characters $1,-1,b$ and costs $v_x,v_y,v_z$.  A monomial witnessing
each minimum together with all Bellman inequalities proves minimality,
so no search cutoff for monomial exponents is used.

Refine each mark by the actual possible source indices:
\[
\begin{array}{c|c}
r&n\\\hline
2&2,4,12,14\\3&3,12\\4&4,12\\5&5\\6&12\\7&14\\12&12\\14&14.
\end{array}
\]
There are thirty-six source assignments before the full-character trace
test and twenty-eight afterwards.  At index four this uses no ordinary
chart: for inertia four the bound is fractional, while for inertia two
a nontrivial full character whose inertia character is zero has positive
even order at least two.  These are the only additional costs needed
there.  A smooth coarse curve cannot have an inertia-one contact at any
unchanged ordinary source: every invariant function in its ambient
maximal ideal has ordinary order at least two, and could not restrict
to a coarse uniformizer of order one.

For total conductors enumerate positive $\delta$ with
$\sum\delta w\le34$ and $\sum w\max(2\delta,I_{\min})\le72$.
The incidence graph of the components and their assigned ordinary sources
must be a forest, by Proposition~\ref{prop:octic-full-conductor46}.
The initial coarse constraints leave fifty-one aggregate states.
The subsequent exclusions are proved below; the finite enumeration is
only bookkeeping for these geometric statements.

\needspace{5\baselineskip}
\paragraph{A discriminant bound and forced smooth lifts}

\begin{lemma}\label{lem:octic-conductor-double-cover46}
Let $\Gamma$ have transverse delta one and generic nodal type.  The
quotient algebra
\[
 \mathcal A_\Gamma=
 (\nu_*\mathcal O_{C_N}\otimes_{\mathcal O_{C_D}}\mathcal O_\Gamma)
                         /\mathcal O_\Gamma\text{-torsion}
\]
is locally free of rank two and has discriminant of degree at most two.
At an unshared ordinary source of full inertia $n$, if its lifted
conductor is smooth and has mark $k/n$, its discriminant order is at least
\begin{equation}\label{eq:octic-branch-bound46}
               \left\lceil\frac{2a}{n}\right\rceil,
                   \qquad a=\overline{-1-9k}\in\{0,\ldots,n-1\}.
\end{equation}
\end{lemma}
\begin{proof}
Torsion is an ideal, so $\mathcal A_\Gamma$ is an actual quotient algebra
of the full conductor upstairs.  Its $H^1$ vanishes, since the kernel
has dimension at most one.  It is torsion-free on the smooth $\Gamma$,
hence locally free; its generic rank is two.  Trace splits it as
$\mathcal O_\Gamma\oplus L$.  Multiplication is a nonzero map
$L^{\otimes2}\to\mathcal O_\Gamma$ because the generic node has a
reduced separable conductor algebra.  Thus $\deg L\le0$, while
$H^1(L)=0$ gives $\deg L\ge-1$.  Its discriminant, including multiple
zeros, has degree $-2\deg L\le2$.

On the cover let $u$ be the smooth conductor parameter of character $c$.
The octic residue generator has character $9b$, so curve duality makes
a generator $t$ of the conductor quotient have character $c-9b$.
The relation $k=-b/c\pmod n$ makes $q=u^at$ invariant for exactly
$a=\overline{-1-9k}$.  If $t^2=h(u)$, its character forces
$\operatorname{ord}h\equiv-2a\pmod n$; hence the order of $q^2$ in the
coarse uniformizer $u^n$ is at least $\lceil2a/n\rceil$.
Normalization and the full conductor commute with taking invariants
here: this is true off the closed point by etale base change, and the
conductor ideals are $S_2$, so the equality extends over that point.
Generic delta one makes the cover conductor generically reduced and
Cohen--Macaulay, hence reduced.  No analogous discriminant assertion is
made for a generic cusp, where the square multiplication is zero.
\end{proof}

A coordinate of lower order one is forced to have actual order one in
either of two ways.  If all such orders were raised to $1+r$, there
would be no full-source invariant of order $r$ to furnish a coarse
uniformizer; or one of the trace upper bounds would drop below $-1$.
This tests all higher branch orders at once.  The exact minima give
\begin{equation}\label{eq:octic-forced-axes46}
\begin{array}{c|c|c}
\text{source and mark}&\text{forced parameter}&\text{branch lower bound}\\\hline
Q_3,\ 1/3\text{ or }2/3&\text{a weight-one parameter direction}&2\\
Q_5,\ 2/5&y&1\\Q_5,\ 3/5&x&1\\
Q_{12},\ 5/12&y&1\\Q_{12},\ 11/12&z&2\\
Q_{14},\ 9/14&x&1\\Q_{14},\ 5/14&y&2\\
Q_{14},\ 13/14&z&2.
\end{array}
\end{equation}
Only the rows in which the raising test succeeds are used.  For example
raising $z$ from one to thirteen at the $11/12$ mark fails weight nine
when $U_9=-1$, and fails weight eleven in row 17; raising $z$ at the
$13/14$ mark fails a trace bound in every occurrence used below.

\begin{proposition}\label{prop:octic-unshared-branch46}
Rows 12,13,15,17,19 of \eqref{eq:octic-nineteen46} are impossible.
The sole-component occurrence of row 10 is impossible, as is the pair
of rows 3 and 7.  More generally applying
\eqref{eq:octic-branch-bound46} at all unshared forced sources deletes
eight of the forty aggregate states left after row 2 is removed.
\end{proposition}
\begin{proof}
The five rows have delta one, are necessarily nodal because $3w>72$,
and have respectively the following unshared branch contributions:
\[
 (Q_3,Q_{14}):(2,2),\quad (Q_3,Q_5):(2,1),\quad
 (Q_5,Q_{12}):(1,2),\quad (Q_{12},Q_{14}):(2,2),\quad
 (Q_5,Q_{14}):(1,2).
\]
Rows of degree at least twenty-eight are sole components.  Row 12 can
only have row 1 as companion, which shares neither of its charged
sources.  In row 13 the two third-inertia marks must use $Q_3$ and
$Q_{12}$, since the coarse curve is smooth and $Q_3$ is unique.
Each displayed sum exceeds two.
The sole row 10 has contributions two and one from $Q_{12},Q_{14}$.
For the pair of rows 3 and 7, the latter has unshared contribution four;
it would have to be cuspidal, but $2\cdot12+3\cdot21>72$.
The same test on each source assignment gives the asserted count;
any unassigned polar remainder must be zero or at least three, by
Lemma~\ref{lem:octic-polar-gap46}.
\end{proof}

\paragraph{The degree-ten component, including its nonreduced conductor}

\begin{proposition}\label{prop:octic-degree10-common-base46}
Row 2 of \eqref{eq:octic-nineteen46} cannot occur with any delta invariant
or any companion component.
\end{proposition}
\begin{proof}
Its third-inertia mark is $Q_3$: assignment to $Q_{12}$ gives
$U_{13}=-2$.  The fourteenth lift is the smooth $y$-axis by coarse
smoothness.  There are no further non-Cartier contacts.  Indeed
\[
 U_5=\frac5{42}-\frac13-\frac{11}{14}=-1,
\]
and an additional inertia-one contact of index two or four would add
an integral positive weight-five cost; ordinary sources admit no such
contact at all.  Therefore its actual weight-nine line has degree
\begin{equation}\label{eq:octic-ten-degreezero46}
            \frac{9\cdot10}{420}-\frac3{14}=0.
\end{equation}
Both $U_{10}$ and $U_{11}$ equal $-1$, so $\gamma=\eta=0$.

Here $\delta\le3$.  The reduced plane curve types of delta at most three
are $A_1,A_2$; $A_3,A_4$; and $A_5,A_6,D_4,D_5,E_6$.
This short list follows by multiplicity: multiplicity two gives
$x^2+y^{n+1}$, and multiplicity three with delta three uses all its
delta in the first blowup, leaving three transverse branches, a cusp
and a transverse line, or the $(3,4)$ branch.  Under the two jet
vanishings their polar lengths are at least
\[
\begin{array}{c|rrrrrrrrr}
\text{type}&A_1&A_2&A_3&A_4&A_5&A_6&D_4&D_5&E_6\\\hline
I&4&5&6&7&8&9&7&8&8.
\end{array}
\]
These follow by substituting the branch parametrizations into
$A f_x+B f_y$: for $A_{2s-1}$ the two branches each contribute at
least $s+1$ when the quadratic-direction coefficient vanishes;
for $A_{2s}$ the orders $(2s+1,2)$ give at least $2s+3$.
For $D_4=xy(x-y)$ a direction on one tangent line contributes
$3+2+2$; for $D_5=x(y^2+x^3)$ the line and cusp contribute at least
eight; for $E_6=x^3+y^4$ the derivative orders are eight and nine.
In the cusp $A_2$, the third jet additionally kills the tangential
constant direction, giving five.
Lengths at least eight exceed the total mass.  Length seven leaves
forbidden residual mass two.  Thus delta three is impossible, and
delta two leaves only $A_3$ with length six.  Such a delta-two
component must be sole, as its mass is sixty and every other conductor
costs at least fourteen.

For that sole delta-two case the full conductor $C_D$ is a ribbon on
$\Gamma$: its nilradical $J$ has $J^2=0$ generically and therefore
everywhere by purity, and $J$ is a line bundle on the smooth support.
The canonical trace quotient $\mathcal M_9$ is an invertible sheaf
on this ribbon.  Away from $Q_{14}$ this follows from the Cartier
canonical frame.  At $Q_{14}$ its equivariant conductor ring is
\[
 \mathbb C[[y,\xi]]/(\xi^2),\qquad
                     \operatorname{wt}(y)=13,\quad\operatorname{wt}(\xi)=c,
                      \quad c\in\{1,5\}.
\]
These are the two normal characters to the smooth $y$-axis.
The invariant ring has generators $s=y^{14}$ and
$\epsilon=\xi y^c$, with $\epsilon^2=0$.  The octic dualizing
character is three, so its coefficient module has character eleven:
\[
 y^3\mathbb C[[s]]\ \oplus\ \xi y^{c+3}\mathbb C[[s]]
                         =y^3\mathcal O_{C_D}.
\]
Since $c+3<14$, no extra power of $s$ is lost.  This proves local
invertibility with the full nilpotent structure.
Its restriction to $\Gamma$ is trivial by
\eqref{eq:octic-ten-degreezero46}.  The map
$\operatorname{Pic}(C_D)\to\operatorname{Pic}(\Gamma)$ is injective:
$H^1(\mathcal O_{C_D})=0$ and the lifted constants give $H^1(J)=0$,
and use the square-zero units sequence.  Thus $\mathcal M_9$ is
trivial, contrary to its $H^0$ vanishing.
For sole delta one, $C_D=\Gamma$ and the same contradiction follows
immediately from its actual degree-zero line.

Only delta one with a companion remains.  The aggregate constraints
leave row 1 with delta one or two, row 3 with delta one, or row 4
with delta one.  If row 2 is nodal, its unshared $Q_3,Q_{14}$ each
cost two discriminant units, so only the companion row 3, which
shares $Q_{14}$, can remain.  If row 2 is cuspidal, its length at
least five permits only row 1 with delta one; that companion must
be nodal, since two cusps have minimum mass $50+21=71$ and any
increase exceeds seventy-two.

For this cuspidal/nodal pair, the reduced conductor components cannot
be disjoint, since the degree-zero component would give a section of
$\mathcal M_9$.  They meet once at a Cartier point, their source sets
being disjoint and the preceding hidden-contact guard excluding any
other non-Cartier meeting.  The trace quotient is invertible with
multidegrees $(0,-1)$.  Curve duality makes the normalization quotient
an invertible line of multidegrees $(-1,0)$, and the conductor algebra
is a finite flat rank-two algebra, split by trace.  Its multiplication
section has multidegrees $(2,0)$.  It is identically zero on the cusp
component, but is a nonzero constant on the nodal component.  This
cannot glue at their common node.

For the remaining nodal pair, rows 2 and 3, the lifted conductor at
$Q_{14}$ is the transverse union of the $y$- and $z$-axes.
Its ring is $\mathbb C[[y,z]]/(yz)$ with weights $(13,5)$.
The weight-nine coefficient module is
\[
            y^3\mathbb C[[y^{14}]]\oplus
                         z^5\mathbb C[[z^{14}]].
\]
Thus it splits at the coarse node; there is no gluing of constant
coefficients.  This is their only intersection, by the rational-tree
property.  The degree-zero summand of the global trace quotient
therefore supplies a section, again a contradiction.
\end{proof}

\paragraph{Four exact local conductor algebras}

\begin{proposition}\label{prop:octic-exact-local-conductor46}
Rows 7,10,14,18 of \eqref{eq:octic-nineteen46} are impossible after the
preceding discriminant cuts.
\end{proposition}
\begin{proof}
For row 7 the only possible companion was row 3, already excluded.
The fifth and twelfth lifts are forced smooth $z$-axes (raising their
orders fails weights fourteen and eleven).  The actual weight-nine
restriction has degree
$9\cdot21/420-1/5-3/12=0$.  Its $U_5=-1$ excludes every extra
inertia-one non-Cartier contact, while ordinary hidden contacts are
impossible.  Since the sole conductor has delta one, the full
canonical trace quotient is exactly this degree-zero line on
$\mathbb P^1$, contradicting its acyclicity.  This includes a cusp.

Row 10 cannot be sole.  Its sole possible companion is row 1, sharing
$Q_{12}$.  The lifted branches are the $z$- and $y$-axes, so the full
reduced conductor ring is $\mathbb C[[y,z]]/(yz)$ with weights $(11,5)$.
Its dualizing character is zero and the octic surface character is
nine; a quotient generator $t$ has character three.  On the row-10
$z$-branch an invariant generator is $q=z^9t$.  Its square coefficient
has character six and order at least six, so $q^2$ has coarse order
at least $(18+6)/12=2$.  The unshared $Q_{14}$ contributes one more.
Row 10 is necessarily nodal, and its discriminant budget is exceeded.

For row 18, the proper-inertia-six branch at $Q_{12}$ has lower orders
$(5,1,1)$.  Weight-nine trace forces $z$ to have order one, and coarse
smoothness then forces $y$ to have order one.  The two conjugate
branches can be straightened exactly to $z=\pm y$; the $x$-coordinate
is an equivariant function of $y,z$.  Thus the conductor ring is
$\mathbb C[[y,z]]/(z^2-y^2)$ with weights $(11,5)$.
Its dualizing character is six, so a quotient generator has character
nine.  Put $s=y^5z$ and $q=y^2zt$.  Every character-six function on
this node is $y^6H(s)$, hence
$q^2=s^2H(s)$.  This source costs at least two, and $Q_5,Q_{14}$
cost at least one each.  The sole component is nodal, a contradiction.

Finally consider row 14.  If its third mark is $Q_3$, that source
costs two and $Q_5$ costs one.  For its $Q_{12}$ assignment the lower
orders are $(2,1,1)$.  A coarse uniformizer of order three forces
$x,y$ to have orders two and one; the weight-nine trace bound forces
$z$ to have order one.  Write one branch
$y=t$, $x=t^2A(t^3)$, $z=tC(t^3)$, with $A(0)C(0)\ne0$.
Its four conjugates have
$x_j=i^jy^2A(i^jy^3)$, $z_j=i^{2j}yC(i^jy^3)$.
The invariant $xy=uA(u)$ has a formal inverse $u=U(xy)$.
The equivariant changes $X=x/A(U(xy))$, $Z=z/C(U(xy))$ put the full
reduced conductor exactly in the complete-intersection form
\[
                         Z^2=y^2,\qquad X^2=y^3Z.
\]
Its equation characters are ten and two, so its curve dualizing
character is five and its quotient generator $t$ has character eight.
Put $S=Xy$.  The character-four coefficient module has basis
$y^8,Xy^9,Zy,XZy^2$ over $\mathbb C[[y^{12}]]$; multiplication by
$S$ cycles these starting with $yZ$.  It is therefore exactly
$(yZ)\mathbb C[[S]]$.  The invariant generator $q=yZt$ satisfies
$q^2=S^2H(S)$.  The source contributes at least two discriminant
units, in addition to one from $Q_5$, again impossible.
\end{proof}

\paragraph{Three further component exclusions}

\begin{proposition}\label{prop:octic-degree24-pole46}
Row 9 of \eqref{eq:octic-nineteen46}, of degree twenty-four, is impossible.
\end{proposition}
\begin{proof}
Its delta is one and its polar length at most three.  Since $U_{11}=-1$,
\eqref{eq:octic-low-jet46} forces $\gamma\ne0$.  As $U_{10}=0$, every
source order in that bound is exact; at $Q_{14}$ the order of $\gamma$
is four.  The only possible companion is row 1, which avoids $Q_{14}$.
The inertia-seven lower orders at that source are $(3,4,1)$.
Weight-nine trace forces $z$ to have order one.  The only invariant
monomials of order seven are $xy$ and $yz^3$, so coarse smoothness
forces $y$ to have order four.

Write the two conjugate branches
$x=z^3A(\epsilon z^7)$, $y=\epsilon z^4C(\epsilon z^7)$,
$\epsilon=\pm1$, $C(0)\ne0$.  The invariant $s=yz^3$ has a formal
inverse $u=U(s)$ for $uC(u)$.  The equivariant changes
$X=x-z^3A(U(s))$, $Y=y/C(U(s))$ put the full reduced conductor in the
exact form
\[
                      X=0,\qquad F=Y^2-z^8=0.
\]
This is a reduced complete intersection; its ordinary square equals its
second symbolic power, since its powers are unmixed.  The octic equation
therefore has the form $a=A_0X^2+B_0XF+C_0F^2$.
On its two branches the coefficient characters $0,3,6$ have lower orders
$0,9,4$, and $\delta X,\delta F$ have lower orders two and seven.
The exact order-four value
$\gamma=2A_0(\delta X)^2+2B_0\delta X\delta F+2C_0(\delta F)^2$
forces $A_0(0)\ne0$.  Equivariant Morse preparation gives
\[
                         a=X'^2+F^\nu h(Y,z),
                   \qquad \nu=2\text{ or }3.
\]
The series $h$ is coprime to $F$ and squarefree, since no other conductor
passes through this source.  With $V=X'/F$ the actual normalization is
$V^2+h=0$ for a node and $V^2+Fh=0$ for a cusp.  These are normal
reduced double planes, with characters $(V,Y,z)=(3,13,5)$.
The canonical form $\omega=dY\wedge dz/(2V)$ has character one, hence
$Y\omega$ is invariant and is a regular reflexive canonical section
on the coarse normalization.

In the nodal case $h$ has character six.  Every such monomial has
$(Y,z)$-weight at least eight for weights $(1,2)$.
Under $Y=s,z=s^2t$, if its actual weighted order is even $d$, the form
$Y\omega$ has exponent $3-d/2\le-1$ in $s$ at a generic exceptional
point.  If $d$ is odd, normalization uses $s=\tau^2$ and the exponent
is $7-d\le-2$.  In the cuspidal case $h$ has character eight and
ordinary order at least four, so $Fh$ has order $d\ge6$.
The ordinary base blowup gives exponent $2-d/2\le-1$ for even $d$,
and $5-d\le-2$ for odd $d$.  This computes a pole of the invariant
form in all cases, without requiring a reduced initial face.
At a rational surface singularity every reflexive canonical section
extends to a resolution \cite[Lemma 54.9.6, Tag 0BBU]{StacksProject}.
Pulling back such an extension to a common resolution upstairs cannot
have the displayed pole.  This contradiction excludes both transverse
types and the row.
\end{proof}

\begin{proposition}\label{prop:octic-degree30-jet46}
Row 16, of degree thirty, is impossible.
\end{proposition}
\begin{proof}
It is sole with delta one.  Its polar length is two, so it is nodal and
$\gamma\ne0$.  The half-mark cannot come from index four, since the
weight-ten numerator there has minimum inertia-two order two and would
give $U_{10}=-1$; the index-two assignment is the only remaining one.
At $Q_{14}$ the inertia-seven lower orders are $(2,5,3)$.
The sole invariant of order seven is $xy$, forcing the first two orders
to be exact.  Raising the $z$ order to ten raises the weight-nine
numerator minimum from eight to fifteen, violating the trace bound.
Thus all three orders are exact.

With weights $(2,5,3)$ the only character-two monomials below weight
eighteen are $x^2$ of weight four and $xz^3,x^3y$ of weight eleven.
All first derivatives of $a$ vanish on the conductor.  Their lowest
terms successively force the coefficients of $x^2$, $x^3y$, $xz^3$
to vanish: use $a_x$ at order two, $a_y$ at order six and $a_z$ at
order eight.  Hence $\operatorname{wt}(a)\ge18$.
The coefficient characters of $(\delta x,\delta y,\delta z)$ are
$(10,8,0)$, whose minimum weights are $(6,9,0)$.  Thus $\delta$
lowers weight by at most three and $\gamma$ has order at least twelve.
But $U_{10}=0$ required the exact order five at this source:
$10\cdot30/420-5/7=0$.  The actual order twelve lowers the bound
to $-1$, contradicting $\gamma\ne0$.
\end{proof}

\begin{proposition}\label{prop:octic-degree21-three-mark46}
Row 6 is impossible.  For row 8 only its actual index-four assignment
survives, its fifth lift is the cusp $(x,y)=(t^3,t^2)$ in an equivariant
Morse plane, and its generic transverse type is nodal with polar length
two.
\end{proposition}
\begin{proof}
Each row has delta one and polar length at most three; $U_{11}=-1$
forces $\gamma\ne0$.  Since $U_{10}=0$ and weight ten is Cartier at
$Q_5$, $\gamma$ is a unit there.  The only constant derivation
coefficient is $\delta z$, and the only character-four quadratic is
$z^2$.  Both coefficients are nonzero.  Hence equivariant Morse
preparation puts the octic in the form $z'^2+f(x,y)$; its conductor
is in the graph $z'=0$.

For row 6, every source assignment has $U_9=-1$.
The fifth lower orders are $(3,2,1)$; raising $z$ from one to six
raises the least weight-nine order from one to six and is impossible.
Thus $z$ has order one, which the graph $z=h(x,y)$, with no constant
or linear term, cannot contain.  The same argument excludes the
index-twelve assignment of the quarter-mark in row 8.
For its index-four assignment, the fifth graph and the invariant
uniformizer force exact orders $(x,y)=(3,2)$; equivariant preparation
gives the cusp $g=x^2-y^3$.

If row 8 were generically cuspidal, its equation would be
$z'^2+g^3h$, with $h$ of character three and ordinary order at least
two.  There is no other conductor in a sole occurrence, and the
possible paired occurrences are excluded below; equivalently apply
the following local argument whenever no companion meets this source.
Its actual normalization is $V^2+gh=0$, $V=z'/g$, with $V$ of
character zero.  The residue $dx\wedge dy/(2V)$ is invariant.
The plane branch function has order $d\ge4$, and the ordinary base
blowup gives exponent $1-d/2$ for even $d$ or $3-d$ after the double
base change for odd $d$.  Both are negative, contrary to rationality.
Thus the sole remaining occurrence is nodal; nonzero $\gamma$ then
gives polar length two.
\end{proof}

\paragraph{Gluing and the remaining pairs}

We shall repeatedly use the following elementary gluing test.
If two reduced conductor components meet, they meet at one ordinary
coarse node.  The actual torsion-free weight-nine restriction
$\mathcal G_9$ on their union is a quotient of $\mathcal M_9$, so
$H^1(\mathcal G_9)=0$.  If both component lines have degree $-1$ and
$\mathcal G_9$ is invertible at the node, its normalization sequence
instead gives $h^1=1$.  Consequently such a locally free gluing is
impossible.  This remains valid with a nonreduced full conductor,
because the reduced union is an actual quotient.

\begin{proposition}\label{prop:octic-pair-gluing46}
The only possible two-component conductor remaining from
\eqref{eq:octic-nineteen46} has degrees seven and eighteen, both delta
one.  Its degree-seven component is cuspidal of polar length three;
its degree-eighteen component is nodal of length two; both second jets
are nonzero.  The residual intrinsic polar mass is fifteen.
\end{proposition}
\begin{proof}
The finite aggregate list after the preceding cuts has only the pairs
listed below.  Their source calculations give the exclusions explicitly.
All uses of a curve's dualizing character below take the character of
its ambient plane volume minus its equation character; subtracting the
octic dualizing character then gives the conductor quotient character.

\emph{Degrees seven and twelve.}
The four possibilities for their deltas are $(1,1),(1,2),(2,1),(3,1)$.
They share $Q_5$ and both weight-nine component lines have degree $-1$.
At least one second jet is nonzero: otherwise the low-jet bounds give
mass at least $76$ in the first three cases and at least $90$ in the
last.  A nonzero second jet is a unit at $Q_5$ and puts both lifted
branches in the same $(x,y)$ Morse plane of weights $(1,4)$.
They are smooth branches tangent to the $x$-axis.  Straighten one to
$y=0$ and the other to $y=x^qu(x^5)$, $q\equiv4\pmod5$.
The invariant intersection length is $\lceil q/5\rceil$; the coarse
node forces $q=4$.  The character-two module consists of pairs
$(x^2h_1(x^5),x^2h_2(x^5))$ with equal constant terms.  It is free over
the coarse nodal ring, generated by $(x^2,x^2)$.
The gluing test excludes all four cases.

\emph{Degrees twelve and fourteen.}
Both deltas are one and at least one second jet is nonzero by the
polar budget.  In the fifth Morse plane the reduced union is exactly
$y(y^2-x^3)=0$.  With $u=x^5$, $v=xy$ its invariant ring is
$\mathbb C[[u,v]]/(v(v^2-u))$.  Its character-two module is generated
freely by $x^2$: over $\mathbb C[[x^5]]$ its basis is
$x^2,x^3y,x^4y^2=x^2(1,v,v^2)$.  The component degrees are $(-1,-1)$,
so the gluing test applies.

\emph{Degrees fourteen and eighteen.}
Both components must be nodes of length two, since replacing either
length by three already exceeds seventy-two.  In particular both
second jets are nonzero and supply the fifth Morse plane.
The reduced union is $x(y^2-x^3)=0$.  Its dualizing character is one;
the octic dualizing character is three, hence the quotient generator
$t$ has character three.  On the degree-fourteen cusp lift an invariant
generator is $q=x^2t$.  Its square coefficient restricts as $xH(s)$,
$s=xy$, and $q^2=x^5H(s)=s^2H(s)$.
The shared fifth point therefore costs two discriminant units; the
unshared proper-inertia-six twelfth point costs at least two by the
node algebra $z^2=y^2$ computed above.  This is impossible.

\emph{Degrees seven and twenty-five.}
At $Q_{12}$ both smooth lifted branches are graphs over the $y$-axis,
whose character is eleven.  After straightening the first, the first
nonzero differences in $x,z$ have orders eleven and seven modulo
twelve.  The coarse node makes the first contact $q$ equal to seven
or eleven.  If $q=7$, an equivariant elimination of $x$ gives the
exact plane curve $Z(Z-y^7u(y^{12}))=0$ with weights $(y,Z)=(11,5)$.
Its dualizing character is six; the quotient generator has character
nine.  On the degree-twenty-five branch, $y^9t$ is invariant and its
square has order at least twenty-four.  The shared source costs two,
and the unshared $Q_{14}$ costs one more.  This component is necessarily
nodal, contradicting its discriminant bound.
If $q=11$, eliminate $z$ instead.  The weight-nine character-three
module consists of pairs $y^9h_i(y^{12})$ with equal constant terms,
so it is free over the coarse node.  Both component degrees are $-1$,
and the gluing test applies.

\emph{Degrees twelve and twenty-one (row 8).}
Both are nodal, since their minimum polar mass is sixty-six.
The degree-twenty-one second jet is nonzero and only its index-four
assignment remains.  The fifth Morse plane has union
$y(x^2-y^3)=0$.  Its dualizing character is four and the quotient
generator has character one.  On the degree-twelve $x$-axis the
invariant generator $x^4t$ has square order at least ten, costing two.
The unshared $Q_{14}$ costs at least two more, a contradiction.

\emph{Degree seven with delta two and degree eighteen.}
The latter must have length two, hence is nodal and has nonzero second
jet.  In the fifth Morse plane the exact equation is
$z^2+x^2y^\nu h(x,y)$, $\nu=4$ or $5$, and $h$ is squarefree and
coprime to $xy$.  Normalize by $V=z/(xy^2)$, of character three.
The normalization is the normal double plane
$V^2+y^{\nu-4}h=0$; its full conductor image is $xy^2=0$.
On the degree-eighteen $y$-axis the invariant trace-free generator is
$y^3V$.  Its square coefficient has character one and therefore order
at least four in $y$.  Its coarse square order is at least two;
$Q_{14}$ costs another one.  This uses the full doubled conductor,
not a reduction of it, and contradicts the node budget.

\emph{Degrees seven and twenty-one (row 8).}
The fifth plane again has union $y(x^2-y^3)=0$ and quotient generator
of character one.  On the degree-twenty-one cusp lift write
$(x,y)=(u^3,u^2)$, with parameter character two.  The invariant generator
is $yt$.  Its square coefficient has parameter order one modulo five,
but order one is absent from the cusp semigroup, so its order is at
least six.  The shared fifth contribution is at least two.
The quarter-source is unshared: the degree-seven component has
$U_7=7/60-1/5-11/12=-1$, which excludes a hidden inertia-one contact
there.  On the fifth cusp the least weight-nine numerator has
parameter order six, so the refined ninth bound on the degree-twenty-one
component is
\[
                 9\cdot21/420-6/5-1/4=-1.
\]
It is sharp.  At the actual index-four source this forces a smooth
lifted parameter of character three.  The octic cover is Gorenstein,
with dualizing character one, so its quotient generator has character
two.  Multiplication by the square of the parameter makes it invariant;
its square has coarse order at least one.  This argument uses only
Gorenstein finite duality on the index-one cover and does not assume
that the ambient index-four source is ordinary or even cDV.
The total contribution is at least three, a contradiction.

The only pair left has degrees seven and eighteen and deltas $(1,1)$.
At least one second jet is nonzero, since otherwise the polar mass is
at least $4\cdot7+3\cdot18=82$.
The fifth Morse plane has axis union $xy=0$; its quotient generator
has character two.  On the degree-seven $x$-axis the invariant generator
$x^3t$ has coarse square order at least two.  Its unshared $Q_{12}$
costs one more, so the degree-seven component cannot be nodal.
It is cuspidal; both components cannot be cuspidal because
$3(7+18)>72$.  The degree-eighteen node must have length two.
If the degree-seven second jet vanished, its weight-eleven gap would
make its cusp length at least five; the minimum mass would be
$5\cdot7+2\cdot18=71$, leaving forbidden remainder one, and any
increase is too large.  Thus its second jet is nonzero and its length
is three.  The total is $21+36=57$, leaving fifteen.
\end{proof}

\paragraph{The last degree-twelve double conductor}

\begin{proposition}\label{prop:octic-degree12-double46}
The sole degree-twelve conductor with delta two is impossible.
\end{proposition}
\begin{proof}
Its fifth and fourteenth lifts are the smooth $x$- and $z$-axes.
If $\gamma\ne0$, its degree-ten line has degree zero, so the $Q_{14}$
order is exactly four.  Straighten its conductor support to $x=y=0$.
The octic is $A x^2+Bxy+Cy^2$, and on the $z$-axis the coefficient
orders are at least $0,6,12$; the derivation orders of $x,y$ are at
least two and ten.  The exact second-jet order forces $A$ to be a
unit.  Equivariant Morse preparation gives
$a=X^2+Y^\nu h$, $\nu=4$ or $5$.
Normalize with $V=X/Y^2$, of character three.  Its double-plane
branch function has character six and $(Y,z)$-weight at least eight
for weights $(1,2)$.  The invariant form $Y\,dY\wedge dz/(2V)$ has
exactly the pole computed in Proposition~\ref{prop:octic-degree24-pole46}.
This is impossible, so $\gamma=0$.

For $A_3$ the polar length is at least six; for $A_4$ it is at least
seven.  The total mass therefore leaves only $A_3$ with length exactly
six.  Its two normalized branches have conductor order two and actual
two-form order one each.  Writing $\Lambda$ for their reduced sum,
\[
            C_N=2\Lambda,\qquad B_N=\Lambda,
                   \qquad\operatorname{div}(\theta)=3\Lambda.
\]
Indeed the component already uses all mass $6\cdot12=72$, and
ampleness of $H$ excludes any further positive divisor.

Its half-mark must be at index two.  Assignment to $Q_{12}$ gives
$U_9=9/35-3/2-2/5-5/14=-2$; assignment to index four gives
$U_{18}=18/35-1-4/5-5/7=-2$.  Thus the conductor avoids $Q_{12}$,
including hidden contacts, and $D$ is normal there.  The octic itself
contains $Q_{12}$, because its numerator character is eight.
On its normal cover, the polar of weight nine has nontrivial
character three, hence is a nonzero nonunit.  A nonzero nonunit on a
normal two-dimensional local domain has a curve of zeros through the
origin.  The cover is normal: it is reduced and normal off the origin
by the free quotient action, is a hypersurface and hence $S_2$, and
satisfies $R_1$; a normal local ring is a domain.
The resulting polar curve cannot lie over the conductor which avoids
this point.  Its positive degree contradicts the exhausted mass.
\end{proof}

\begin{theorem}[Exhaustive conductor states]\label{thm:octic-eleven-states46}
Every hypothetical octic under consideration has one of the following
eleven total conductor states:
\begin{equation}\label{eq:octic-eleven46}
\begin{array}{c|l|l}
\text{degrees and deltas}&\text{markings}&\text{additional restrictions}\\\hline
(7,1),(7,2),(7,3),(7,4)&3/5,5/12&\text{one component in each state}\\
(12,1)&1/2,3/5,13/14&\text{cusp; half-mark at index two}\\
(14,1),(14,2)&1/5,5/6&\delta=2:\ A_3,\ \gamma\ne0,\ I=4\\
(18,1)&2/5,9/14&\text{node or cusp}\\
(21,1)&1/4,4/5&\text{node; actual quarter-source of index four}\\
(25,1)&5/12,9/14&\text{node}\\
(7,1)+(18,1)&(3/5,5/12);(2/5,9/14)&\text{cusp/node; lengths }3,2.
\end{array}
\end{equation}
The table has ten sole-component states and one two-component state.
\end{theorem}
\begin{proof}
The construction of \eqref{eq:octic-nineteen46}, all possible source
assignments and the total mass and forest constraints were specified
above.  They give fifty-one aggregate states.
Proposition~\ref{prop:octic-degree10-common-base46} removes eleven;
Proposition~\ref{prop:octic-unshared-branch46} removes eight more;
Proposition~\ref{prop:octic-exact-local-conductor46} removes four.
There remain twenty-eight.  The degree-twenty-four, degree-thirty and
three-mark degree-twenty-one exclusions remove two, one and three,
respectively.  The seven pair calculations remove ten further states,
leaving twelve, and Proposition~\ref{prop:octic-degree12-double46}
removes the last degree-twelve double state.  These are the eleven
states displayed.  Each enumeration uses at most four components,
since each component has degree at least seven and $\sum\delta w\le34$.
No larger state can be omitted.

For the listed type refinements, a sole degree-twelve node has
unshared discriminant contributions one and two from $Q_5,Q_{14}$,
so it is cuspidal.  Its half-mark is at index two: the same
weight-nine and weight-eighteen bounds displayed in
Proposition~\ref{prop:octic-degree12-double46} exclude the index-twelve
and index-four assignments independently of the delta invariant.
A degree-fourteen delta-two component with
$\gamma=0$ has length at least six, exceeding the mass; with
$\gamma\ne0$ the $A_4$ option has length five and leaves forbidden
mass two.  Thus it is $A_3$ of length four.  The degree-twenty-one
refinement follows from its fifth invariant-form obstruction, now
with all its companions excluded.  Degree twenty-five is too large
for a cusp.  The pair refinements were proved above.
\end{proof}


\section{Local surface geometry for the octic}
\label{sec:octic-atlas-revised}

The preceding reduction supplies rational normalization singularities
and a K3 minimal model.  We use these two facts in different ways.
Rationality constrains the local equations through canonical forms;
the K3 model constrains the exceptional curves through the divisor of
the actual restricted two-form.  We prove the common tools and local
charts here, with their hypotheses stated explicitly.  The subsequent
case arguments will supply the appropriate polar budget and incidences.

\subsection{Canonical forests and the degree-seven stack}
\label{sec:octic-forest-tools46}

Let $p:W\to N$ be a smooth resolution and let $\sigma:W\to S$ be its
morphism to the K3 minimal model.  The actual two-form is the pullback of
a nowhere-vanishing form on $S$, so
\begin{equation}\label{eq:octic-canonical-forest46}
 A:=\operatorname{div}_W(\alpha)=K_{W/S}.
\end{equation}
Thus every positive coefficient is an integer, and the support of $A$
is the exceptional forest of a sequence of point blowups.  A curve born
at a point on components of coefficients $a,b$ has coefficient $1+a+b$;
at a point on just one component it has coefficient $1+a$; at an ordinary
point it has coefficient one.  Every edge in the final graph joins a
curve to one of its parents.  In particular adjacent positive coefficients
are unequal, a vertex has at most two smaller neighbors, and if there are
two then its coefficient is one plus their sum.

\begin{lemma}\label{lem:octic-forest-minimum46}
Every connected subtree of the positive canonical forest has a unique
vertex of smallest coefficient.  For every adjacent pair of coefficients
$a,b$ one has $\gcd(a+1,b+1)=1$.
\end{lemma}
\begin{proof}
Induct on the blowups.  A blowup away from the existing forest adds an
isolated coefficient-one vertex.  A blowup on one component adds a leaf
whose coefficient is larger than its parent's.  A blowup at an
intersection subdivides an edge $a$--$b$ by $a+b+1$, again larger than
either endpoint.  For a connected subtree containing the new vertex,
remove a new leaf, or restore the subdivided old edge when both parents
are present.  This reduces its minimum to that of a connected old
subtree.  If only one parent is present, removing the new vertex has the
same effect.  The singleton case is immediate.  This proves uniqueness
of the minimum in every connected subtree, not just in the whole forest.
For shifted coprimality, a new leaf has shifted coefficients $a+1,a+2$.
At a corner blowup the new shifted coefficient is $(a+1)+(b+1)$, which
preserves the gcd on each of the two new edges.  Induction proves the
second assertion.
\end{proof}

\begin{lemma}\label{cor:octic-root-budget}
If a negative-definite configuration of smooth rational $(-2)$-curves on
$W$ has zero coefficient in $A$ on every component, its rank is at most
nineteen.  All its curves and intersections survive unchanged on $S$.
\end{lemma}
\begin{proof}
Adjunction gives $A E=K_WE=0$ for each such curve.  Since $A$ is
effective and does not contain $E$, the curve is disjoint from its entire
support.  It is therefore unaffected by the contractions to $S$.
The Neron--Severi group of a projective complex K3 surface has signature
$(1,\rho-1)$ and $\rho\le20$, proving the rank bound.
\end{proof}

We call these surviving curves \emph{protected curves}.  When a local
resolution is used below, zero order always means the order of the actual
ambient two-form, not of a general section of the canonical sheaf.
Every irreducible component of the intrinsic polar divisor on $N$ has
rational strict transform in the forest, by
\eqref{eq:octic-canonical-forest46}.  Its positive degree must be charged
to the global mass budget.

\begin{lemma}\label{lem:octic-degree7-stack46}
Suppose the reduced conductor is a sole degree-seven component $\Gamma$.
Its normalized canonical-index stack with its canonical restriction is
\[
 (\mathcal G,H)\simeq(\mathbb P(5,12),\mathcal O(1)),\qquad
 G:=\bigoplus_{m\ge0}H^0(\mathcal G,mH)=\mathbb C[X_5,Y_{12}],
 \qquad \omega_{\mathcal G}=\mathcal O(-17).
\]
There are no additional non-Cartier contacts of $\Gamma$ with $Z$.
In generic delta one its full normalization quotient is
$\mathcal Q=\mathcal O_{\mathcal G}(-26)$, and the full conductor
algebra upstairs is
\begin{equation}\label{eq:octic-degree7-boundary46}
 G[T_{26}]/(T_{26}^2-cX_5^8Y_{12}),
\end{equation}
with $c\ne0$ for the node and $c=0$ for the cusp.
\end{lemma}
\begin{proof}
The component is a smooth rational curve with precisely the markings
$3/5,5/12$.  The source bounds give
$U_7=7/60-1/5-11/12=-1$.  An additional inertia-one contact at an
index-two or index-four point would add a positive integral weight-seven
cost, violating the trace bound.  An inertia-one contact at an ordinary
source is already impossible by coarse smoothness.  Additional
Gorenstein ambient singularities are allowed and do not change the
canonical-index stack.

The effective stack over the smooth rational curve consequently has
two cyclic points of coprime orders five and twelve.  It is the root
stack $\mathbb P(5,12)$: the usual quotient construction has trivial
generic stabilizer, the two specified local stabilizers, and coarse space
$\mathbb P^1$.  Its Picard group is infinite cyclic, generated by
$\mathcal O(1)$ of degree $1/60$.  Since $\deg H=7/420=1/60$, $H$
is this generator.  Its local markings are indeed $3/5$ and $5/12$:
at the first chart the parameter has character twelve modulo five,
and at the second it has character five modulo twelve.  The graded
section ring and canonical class follow from this quotient description.

For delta one the full conductor image is reduced: it is generically
reduced and Cohen--Macaulay, so it has no nilpotent section supported at
finitely many points.  It is therefore the actual smooth stack curve
$\mathcal G$.  Finite Gorenstein duality on the canonical stack gives
\[
 \mathcal Q\simeq\omega_{\mathcal G}\otimes\mathcal O_{
 \mathcal G}(-9H)=\mathcal O_{\mathcal G}(-26).
\]
The conductor algebra is a locally free rank-two module.  Its trace
splits off the unit, and its multiplication is a section of
$\mathcal O_{
 \mathcal G}(52)$.  The only monomial of degree fifty-two is
$X_5^8Y_{12}$.  A generic node has a reduced separable algebra, so its
coefficient is nonzero; a generic cusp has square-zero trace-free
part, so the coefficient is zero.  This proves the full algebra, not
merely a generically defined double cover.
\end{proof}

\begin{lemma}\label{lem:octic-weighted-line-splitting}
Every vector bundle on $\mathcal G=\mathbb P(5,12)$ is a sum of line
bundles.  Writing $G=\mathbb C[X_5,Y_{12}]$ and $G_d=0$ for $d<0$, one has
\[
 h^0(\mathcal G,O(d))=\dim G_d,\qquad
 h^1(\mathcal G,O(d))=\dim G_{-17-d}.
\]
\end{lemma}
\begin{proof}
Pull a bundle back to $\mathbb A^2\setminus\{0\}$ and extend it
reflexively over the origin.  The resulting graded module over
$\mathbb C[X,Y]$ has depth two.  A graded minimal free resolution
and Auslander--Buchsbaum make it graded free.  Sheafification gives
the splitting.  The cohomology formulas follow from the section ring
and Serre duality with $\omega_{\mathcal G}=O(-17)$.
\end{proof}

\begin{lemma}[Intrinsic higher jets]\label{lem:octic-intrinsic-jets}
Let $a$ be a section of $mK_Z$ and let $\delta$ be the degree-one
twisted derivation defined by the two-form.  On a reduced conductor
curve $\Gamma$, suppose
$a|_\Gamma=\delta a|_\Gamma=\cdots=\delta^{k-1}a|_\Gamma=0$.
Then $\delta^ka|_\Gamma$ is a regular section of
$O_\Gamma((m+k)H)$ on the canonical-index stack.
\end{lemma}
\begin{proof}
Write $a=f\eta^m$ and $\alpha=\iota_\delta\eta$ in a local
canonical frame.  Replacing $\eta$ by $g\eta$ replaces $f$ by
$g^{-m}f$ and $\delta$ by $g^{-1}\delta$.  The leading term of
the $k$th iterate is $g^{-m-k}\delta^kf$; every other term contains
a lower iterate.  These vanish on $\Gamma$, proving the required
transition rule.  Regularity follows from the regular derivation
on each canonical cover.
\end{proof}


\subsection{Positive coefficients and protected root configurations}
\label{sec:octic-forest-chain}

The uniform K3 reduction gives a useful criterion stronger than merely
excluding equal coefficients on adjacent components of the canonical
exceptional forest.

\begin{proposition}\label{prop:octic-root-chain-monotonicity}
Let \(\sigma:W\to S\) be a birational morphism between smooth projective
surfaces, with \(S\) a K3 surface, and let
\(A=K_{W/S}=\operatorname{div}(\alpha)\).
Suppose \(E_1,\ldots,E_n\), \(n\geq3\), form a chain of smooth
rational \((-2)\)-curves on \(W\).  Then their coefficients in \(A\)
are either all zero or are positive and strictly monotone along the chain.
In the latter case, the endpoint of larger coefficient \(m\) meets a
component of \(A\) outside the chain whose coefficient is at least
\(m+1\).
\end{proposition}

\begin{proof}
Write \(m_i=\operatorname{mult}_{E_i}A\), allowing zero.
Adjunction gives \(AE_i=0\).  Hence if \(m_i=0\), every component
of \(A\) is disjoint from \(E_i\), and zero propagates along the
chain.  Otherwise all the \(m_i\) are positive.

The coefficient rule preceding
Lemma~\ref{lem:octic-forest-minimum46} shows that adjacent positive
coefficients are unequal.  A vertex with two smaller neighbours has
coefficient one plus their sum.

For the indicated chain, adjunction and effectivity of the remaining
components give
\[
  2m_i\geq m_{i-1}+m_{i+1},\qquad
  m_0=m_{n+1}=0.
\]
Thus the sequence, after adjoining zero endpoints, is discretely concave.
For an interior index \(2\leq i\leq n-1\), its chord inequalities give
\[
 m_{i-1}\geq\frac{i-1}{i}m_i,
 \qquad
 m_{i+1}\geq\frac{n-i}{n+1-i}m_i.
\]
The two fractions sum to at least one.  Therefore an interior strict
local maximum would satisfy
\(m_{i-1}+m_{i+1}\geq m_i\), contradicting the parent formula
\(m_i=1+m_{i-1}+m_{i+1}\).
Discrete concavity already excludes an interior strict local minimum,
and adjacent equality is impossible.  The sequence must therefore be
strictly monotone.

An exceptional curve is born with self-intersection \(-1\).
The larger endpoint has final self-intersection \(-2\), so exactly one
subsequent blowup had its centre on that curve.  The exceptional curve
created by that blowup is still adjacent to it and has coefficient at
least one larger.  It cannot be its sole neighbor inside the monotone
chain, whose coefficient is smaller.  This proves the last assertion.
\end{proof}

For the seven-curve arm at the compulsory twelfth source, a positive
coefficient vector must therefore be strictly monotone, not simply
nonconstant.  If its larger endpoint is the free end of the full local
resolution graph, the additional component in the proposition is the
strict transform of an intrinsic positive polar curve on the normalized
surface.  It is not another exceptional curve over that source.  Its
multiplicity must be charged to the global polar degree budget.  This
last conclusion is conditional on which endpoint is larger; the
proposition alone does not choose that orientation or exclude the
remaining conductor state.

\begin{corollary}[Forks are protected]\label{cor:octic-protected-forks}
Every connected configuration of smooth rational \((-2)\)-curves with
Dynkin diagram \(D_n\) or \(E_n\) on a surface birational to a K3 surface
has zero order for the holomorphic two-form on every component.  In
particular it survives unchanged on the K3 minimal model.
\end{corollary}

\begin{proof}
Zero propagates along a connected \((-2)\)-configuration by adjunction.
It is therefore enough to rule out the alternative in which all its
curves belong to the canonical exceptional forest.

In a graph obtained by point blowups of a smooth surface, two exceptional
\((-1)\)-curves are never adjacent.  This follows inductively because
a newly created \((-1)\)-curve meets only parent curves whose squares
have just dropped to at most \(-2\).  A curve of final square \(-2\)
was born with square \(-1\) and had exactly one subsequent blowup centred
on it.  At its birth at most one parent can remain a \((-2)\)-curve in
the final surface: if there were two parents, one already had square at
most \(-2\) and its square dropped to at most \(-3\).  Its unique later
blowup creates at most one further neighbor of final square \(-2\).
Thus the subgraph of final \((-2)\)-curves in the exceptional forest
has valency at most two and, being a forest, is a disjoint union of
chains.  It cannot contain a \(D\)- or \(E\)-fork.
\end{proof}

\begin{proposition}[A local forest coefficient bound]\label{prop:octic-local-coefficient-bound}
Suppose every curve in the complete resolution fibre over a normal
surface point has square at most \(-2\).  Let \(T\) be an upper
bound for the available scaled degree of the intrinsic polar divisor,
and suppose every intrinsic polar curve of coefficient at least two
has scaled degree at least \(d>0\).  Every positive canonical
coefficient on that resolution fibre is at most
\[
                         M=\left\lfloor T/d\right\rfloor-1.
\]
A positive \(A_n\)-chain attached to a further positive exceptional
curve is impossible when \(n\geq M\).
\end{proposition}

\begin{proof}
A positive exceptional curve with square at most \(-2\) has had a
subsequent blowup centred on it.  The curve born at the last such
blowup is still adjacent to it and has strictly larger canonical
coefficient.  Following such larger neighbors must eventually leave
the finite local resolution graph.  Its first exterior curve is an
intrinsic polar curve and has coefficient at least one larger than
that of the starting curve, hence at least two.  It alone consumes
scaled degree at least \(d(m+1)\).  This proves the bound.

For the chain assertion, the cases \(M\leq1\) are immediate from the bound
and the prohibition on adjacent equal positive coefficients.  Otherwise
the preceding proposition makes its positive coefficients strictly monotone.  If \(n>M\), there are already too
many distinct positive integers.  If \(n=M\), the coefficients must
be \(1,2,\ldots,M\), in one of the two orientations.  If the
additional positive curve meets the coefficient-one endpoint,
adjunction on that \((-2)\)-curve is impossible.  If it meets the
coefficient-\(M\) endpoint, its coefficient is a positive integer
less than \(M\): equality is forbidden and the bound excludes a
larger value.  The coefficient-\(M\) curve would then have two
smaller parents, one of coefficient \(M-1\), giving
\(M=1+(M-1)+c\) with \(c>0\), again impossible.
\end{proof}


\subsection{A discriminant bound for reduced conductor covers}
\label{sec:reduced-conductor-double-cover45}

\begin{lemma}
\label{lem:reduced-conductor-double-cover45}
Let $C_N$ be the full conductor scheme on the normalization of the
actual octic, and assume $H^1(C_N,\mathcal O_{C_N})=0$.
Let $\Gamma\simeq\mathbb P^1$ be a reduced conductor image and
let $R$ be the reduced union of the conductor components mapping
onto it.  If $R\to\Gamma$ has generic separable degree two,
then its actual discriminant divisor has degree at most two.
The full conductor may be nonreduced and may have arbitrary
generic delta invariant.
\end{lemma}

\begin{proof}
The quotient $\mathcal O_{C_N}\to\mathcal O_R$ gives
$H^1(\mathcal O_R)=0$, because its kernel has dimension at most one.
The finite pushforward $p_*\mathcal O_R$ is torsion-free over
the smooth curve $\Gamma$, hence locally free of rank two.
Trace splits it as $\mathcal O_\Gamma\oplus L$.
The trace-free square is a nonzero map $L^2\to\mathcal O_\Gamma$,
so $\deg L\leq0$.  Its $H^1$ vanishes, so $\deg L\geq-1$.
The discriminant is therefore a section of $L^{-2}$ of degree
zero or two.  Its zeros are counted with their full multiplicities.
\end{proof}

For the octic under consideration, the required vanishing for the
full conductor is Proposition~\ref{prop:octic-full-conductor46}.
Thus this discriminant bound applies in every delta invariant.


\subsection{Weighted canonical forms and forced source coefficients}
\label{sec:uniform-weighted-pole45}

The following criterion packages the local rationality tests into a single
inequality.  Its proof retains the entire exceptional scheme.  In particular,
neither a reduced initial equation nor a nonzero numerator on the reduced
exceptional curve is required.  This is essential when a weighted initial equation has repeated factors.

\begin{theorem}[Weighted canonical-form inequality]
\label{thm:weighted-canonical-inequality45}
Let $D^\#=(f=0)$ be a reduced hypersurface surface germ in a smooth
three-dimensional germ over $\mathbb C$, with a diagonal finite cyclic
action preserving it.  Let $N$ be the quotient of its finite normalization.
Assume that $N$ has rational singularities.  Write
\[
 \Omega=\operatorname{Res}\frac{dx\wedge dy\wedge dz}{f}.
\]
Let $q$ be a monomial for which $q\Omega$ is invariant and is a regular
canonical form on the normalization.  For every positive integer weight
vector $w=(w_x,w_y,w_z)$ one then has
\begin{equation}
 \boxed{\operatorname{ord}_w(f)
       <w_x+w_y+w_z+\operatorname{ord}_w(q).}
 \label{eq:octic46:octic_weighted_criterion46:45-W1}
\end{equation}
For a nonnormal hypersurface, it suffices that $q$ belong to its full
conductor ideal.  For a normal hypersurface, every regular invariant
canonical numerator is allowed.
\end{theorem}

\begin{proof}
Put $d=\operatorname{ord}_w(f)$, $e=\operatorname{ord}_w(q)$, and let
$f_d$ be the initial form.  Suppose $d\ge w_x+w_y+w_z+e$.
On a smooth covering chart of the weighted blowup, with exceptional
parameter $s$, write
\begin{equation}
 f=s^dF,\qquad q=s^eq',\qquad
 q\Omega=s^{w_x+w_y+w_z+e-d-1}q'\eta,
 \label{eq:octic46:octic_weighted_criterion46:45-W2}
\end{equation}
where $\eta$ generates the dualizing module of the hypersurface strict
transform.  The exponent in \eqref{eq:octic46:octic_weighted_criterion46:45-W2} is at most $-1$.

We explain why some chart detects a pole even if $f_d$ is nonreduced.
The weighted initial numerator has degree $e<d$ and is nonzero in
$\mathbb C[x,y,z]/(f_d)$.  This hypersurface ring is Cohen--Macaulay of
dimension two.  Its principal ideal is saturated at the irrelevant ideal,
so the numerator defines a nonzero section on the full exceptional
hypersurface $E=(f_d=0)\subset\mathbb P(w)$, interpreted here on the
weighted projective stack.  Its twist $O_E(e)$ is a line bundle there;
no invertibility on the coarse weighted projective space is asserted.
The stack curve $E$ is pure Cohen--Macaulay.  Its nonzero section is
therefore nonzero at some generic component, with that component's full
scheme multiplicity retained.  Equivalently, this is the unmixedness
argument in the displayed graded hypersurface ring, followed by a
smooth covering chart.

Pass to a covering chart at such a generic component and let $A$ be the
one-dimensional local ring of the strict transform there.  It is reduced
and Gorenstein, $s$ is a nonzerodivisor, and $q'\notin sA$.  Let $B$ be
its finite normalization, a product of discrete valuation rings.  Finite
birational duality identifies, inside the same rational canonical forms,
\begin{equation}
 \omega_B=\operatorname{Hom}_A(B,\omega_A)
          =\mathfrak c\,\omega_A\subseteq\omega_A=A\eta,
 \label{eq:octic46:octic_weighted_criterion46:45-W3}
\end{equation}
where $\mathfrak c=(A:B)$ is the full conductor.
If \eqref{eq:octic46:octic_weighted_criterion46:45-W2} were regular on every normalization branch, it would belong
to $\omega_B$, hence to $A\eta$.  Its negative power of $s$ would then
force $q'\in sA$, a contradiction.  Thus it has a pole on at least one
normalization branch.

The chart covers and the original cyclic group do not remove that pole.
For a finite map of smooth surface germs, canonical orders along a prime
with ramification index $r$ satisfy
$\operatorname{ord}_{\rm cover}=r\operatorname{ord}_{\rm base}+r-1$.
Equivalently, the pullback of a regular canonical form is regular.
The assumed invariant form descends to a regular reflexive canonical
section on $N$.  Rationality gives its extension to a resolution, and
then to every dominating resolved cover.  The pole just obtained is
impossible.
\end{proof}

The finite-duality identification is the birational specialization of
finite duality \cite[Section 48.11, Tag 0AWZ]{StacksProject}.
Extension of canonical forms over a rational surface singularity is
\cite[Lemma 54.9.6, Tag 0BBU]{StacksProject}.  The criterion is a necessary weighted condition; it
does not assert that satisfying all its inequalities proves rationality.

\paragraph{Three Q12 coefficients are forced in every delta-one case}

At Q12 the ordinary cover has characters $(1,11,5)$ and the octic
equation has character $8$.  For the sole degree-seven conductor of
generic delta one, its full ideal is $(x,z)$ and $f\in(x,z)^2$.
The residue has character $9$.  Each of the following numerators belongs
to the conductor and makes the residue invariant:
\begin{equation}
\begin{array}{c|c|c|c}
w&q&\sum w+\operatorname{ord}_w(q)&
 \text{all allowed terms of smaller weight}\\\hline
(1,1,3)&x^3&8&x^3z\\
(5,7,1)&z^3&16&z^4\\
(3,1,3)&y^2z&12&y^2z^2.
\end{array}
\label{eq:octic46:octic_weighted_criterion46:45-W4}
\end{equation}
Consequently
\begin{equation}
 [x^3z]f\ne0,\qquad[z^4]f\ne0,
 \qquad[y^2z^2]f\ne0.                         \label{eq:octic46:octic_weighted_criterion46:45-W5}
\end{equation}
These conclusions hold in both the nodal and cuspidal alternatives,
independently of the vanishing of $\gamma_{10}$.
Each coefficient is therefore forced by its own weighted test.
For example an $x^8$-only initial equation may have an empty exceptional
set in the $x$-chart, but the full exceptional scheme in the other charts
still violates \eqref{eq:octic46:octic_weighted_criterion46:45-W1}.  Replacing that scheme by its reduced support
would lose exactly the information used by \eqref{eq:octic46:octic_weighted_criterion46:45-W3}.

\paragraph{The Q14 clauses are direct instances of the same theorem}

At a normal ordinary Q14 point, the cover characters are $(1,13,5)$,
$f$ has character $2$, and the residue has character $3$.
The four rows
\begin{equation}
\begin{array}{c|c|c|l}
w&q&\sum w+\operatorname{ord}_w(q)&
\text{allowed lower terms}\\\hline
(6,1,2)&y^3&12&y^3z\\
(2,5,3)&xz^2&18&x^2,\ x^3y,\ xz^3\\
(3,4,1)&xz^2&13&x^2,\ xz^3,\ z^6\\
(17,11,1)&z^5&34&xz^3,\ z^6,\ yz^9,\ z^{20}
\end{array}
\label{eq:octic46:octic_weighted_criterion46:45-W6}
\end{equation}
give a nonzero coefficient in every listed lower-term set.  These clauses
apply to arbitrary degenerate initial equations and will make the local
case division exhaustive.  The displayed lists are finite: for a bound
$D$ and positive weights, enumerate $0\le i<D/w_x$,
$0\le j<D/w_y$, $0\le k<D/w_z$, and retain precisely the monomials
of the required full source character and weight below $D$.


\subsection{The normal fourteenth source: five isolated formats}
\label{sec:octic-q14-atlas}
\label{sec:octic-q14-five-strata}

We now describe the normal ordinary $Q_{14}$ source.  Throughout this
subsection the coarse surface germ is rational, the cover coordinates
$(x,y,z)$ have characters $(1,13,5)$, and the equation has character two.
The descriptions below allow all terms of larger weight.  They will be
used with several different conductor states, so we separate these local
calculations from the global polar and Picard-rank arguments.

\begin{proposition}\label{prop:octic-q14-atlas}
After an equivariant analytic coordinate change, the germ belongs to
one of the five formats $G,H,J,K,L$ in
\eqref{eq:octic46:octic_q14_forest46:45-QG}, or to one of the two
formats of weighted degree $38$ and $48$ described in
Subsections~\ref{sec:octic-q14-degree38} and
\ref{sec:octic-q14-degree48}.  In the five isolated formats the complete
quotient graphs and their orientations are those in the table.
\end{proposition}

\begin{proof}
The four weighted tests in
\eqref{eq:octic46:octic_weighted_criterion46:45-W6} force $y^3z$ and
at least one term in each of the other three indicated sets.
If $[x^2]f\ne0$, equivariant Morse splitting gives $x^2+g(y,z)$,
without changing the nonzero $y^3z$ coefficient.  Applying the last
weighted test to this prepared equation shows that the first nonzero
term among $z^6,yz^9,z^{20}$ gives respectively $G,H,J$.

Suppose $[x^2]f=0$.  If both $x^3y$ and $xz^3$ occur, we obtain $K$.
If only the former occurs, the third weighted test forces $z^6$,
giving $L$.  Otherwise $xz^3$ occurs and $x^3y$ does not.  The presence
or absence of $x^4y^2$ gives the degree-$38$ and degree-$48$ formats,
respectively.  This proves that the seven alternatives are exhaustive.

The following table gives each weighted initial equation, the weights,
and the full quotient resolution graph.  Each arm is listed from the
central curve outwards; all entries are positive numbers denoting the
negative of a self--intersection.  Repeated 2's are abbreviated by
$2^j$.
\begin{equation}
\begin{array}{c|c|c|c}
 \text{format}& f_{\rm in}&\text{weights}&
       \text{centre; arms}\\ \hline
 G&x^2+y^3z+z^6&(9,5,3)&4;\ [2,2],\ [2,2],\ [2^4]\\
 H&x^2+y^3z+yz^9&(13,8,2)&2;\ [2^3],\ [2],\ [4,2^4]\\
 J&x^2+y^3z+z^{20}&(30,19,3)&2;\ [2],\ [2^2],\ [2,2,4,2^4]\\
 K&x^3y+xz^3+y^3z&(1,1,1)&3;\ [2,4],\ [2^6],\ [3,2,2]\\
 L&x^3y+y^3z+z^6&(13,15,9)&2;\ [2,2,2,4],\ [3,2^6],\ [3].
\end{array}                                                   \label{eq:octic46:octic_q14_forest46:45-QG}
\end{equation}
Every displayed coefficient is nonzero.  Each initial hypersurface cone
in this table is isolated away from its vertex; the ordinary quartic in
row $K$ is a smooth plane quartic.  Thus the weighted strict transform
is smooth along its initial curve before the specified finite
quotients.  Arbitrary terms of larger weight do not change these
quotient graphs.  These are entire analytic strata, not merely their
generic members.

For completeness, the raw central squares and arm corrections are
\begin{equation}
\begin{array}{c|c|c}
 G&-28/15&2/3,\ 2/3,\ 4/5\\
 H&-7/16&3/4,\ 1/2,\ 5/16\\
 J&-7/57&1/2,\ 2/3,\ 27/38\\
 K&-8/7&4/7,\ 6/7,\ 3/7\\
 L&-28/65&10/13,\ 7/15,\ 1/3.
\end{array}                                                   \label{eq:octic46:octic_q14_forest46:45-QC}
\end{equation}
Subtracting the corrections gives the central entries in \eqref{eq:octic46:octic_q14_forest46:45-QG}.
In $G,L$ the original group is entirely radial, of order $14$.
In $H,J$ its radial subgroup has order $7$, with a residual involution
on the initial rational curve.  In $K$ the radial subgroup has order
$2$, and the residual $\mu_7$ quotient of the genus--three quartic is
rational; its three fixed coordinate points give the listed arms.
For the long arms of $H,J$, the quotient types are respectively
$\frac1{16}(1,13)$ and $\frac1{38}(1,31)$, with the orientations shown.
For $L$, the three quotients are
$\frac1{13}(1,4)$, $\frac1{15}(1,13)$, and $\frac13(1,1)$.
These quotient types and attachment orientations determine the displayed
continued fractions and intersection matrices.

\end{proof}

\paragraph{Orders of the actual two-form.}
Write
\[
 \alpha=\delta_x\,dy\wedge dz-\delta_y\,dx\wedge dz
             +\delta_z\,dx\wedge dy,
 \qquad \operatorname{char}(\delta_x,\delta_y,\delta_z)=(10,8,0).
\]
The following observations will be used repeatedly.  In $H$, the
smallest weights of actual summands give positive central order:
$\delta_x=z^2$ gives order one and $\delta_z=1$ gives order two
after the radial quotient of order seven.  The minimum central
order in $J$ is three.  In $L$, even $dx\wedge dy$ has total weight
$28$, hence order one after the radial quotient of order fourteen.
Thus the central curve is positive in all three formats.

In $K$, a unit $\delta_z$ gives zero order on the $A_6$ arm and
both curves of the $[2,4]$ arm.  Indeed, on the $x$-chart we have
$y/x\sim-(z/x)^3$ and $S=x^2$ after the radial quotient; the leading
form is a unit multiple of $w^2dS\wedge dw$ in $\frac17(1,2)$,
whose orders on both resolved curves are zero.  The $y$-chart gives
a unit form on $A_6$.  Higher-weight terms cannot cancel these
initial forms at the generic points of the exceptional curves.
If $\delta_z(0)=0$, every ambient summand has ordinary degree at
least four, so the main coefficient is positive.

In $G$, the weight-fourteen part of the actual form is
\begin{equation}
 A z^2dy\wedge dz+B dx\wedge dy.
 \label{eq:octic46:octic_q14_forest46:45-GF}
\end{equation}
If $A=B=0$, the central coefficient is positive.  Otherwise it is
zero, and the leading coefficients on the two $A_2$ arms are
proportional to $6B+2iA$ and $6B-2iA$; on the $A_4$ arm the
coefficient is proportional to $B$.  If $B=0$ and $A\ne0$, the
initial form in the $A_4$ chart is a unit times $u^5dS\wedge du$.
The four orders are positive, the first is one, and they are bounded
above by $1,2,3,4$.  The root-chain criterion therefore makes them
exactly $1,2,3,4$.

\begin{corollary}\label{cor:octic-q14-J}
The actual octic cannot have format $J$ at $Q_{14}$.
\end{corollary}
\begin{proof}
The centre and the first curve of each of its three arms form a $D_4$
configuration of $(-2)$-curves.  Corollary~\ref{cor:octic-protected-forks}
makes their actual two-form orders zero, whereas the central order
in format $J$ is at least three.
\end{proof}


\subsection{The normal fourteenth source: the degree-$38$ format}
\label{sec:octic-q14-degree38}

We retain the rationality and source-character assumptions of
Subsection~\ref{sec:octic-q14-atlas}.  Assume
$[x^2]f=[x^3y]f=0$ and
$[xz^3]f[y^3z]f[x^4y^2]f\ne0$.
We give the full resolution calculation, including repeated roots of
the auxiliary cubic.  Only its local conclusions will be used here;
the resulting global obstructions depend on the conductor state.

The original index--one chart has $\mu_{14}$--weights $(1,13,5)$ on
$(x,y,z)$, and $f$ has character $2$.  With weights $(5,9,11)$ the only
character--$2$ monomials of weight at most $38$ are
\[
 x^2,\quad x^3y,\quad x^4y^2,\quad xz^3,\quad y^3z.
\]
After scaling, the initial equation is
\begin{equation}
 f_{38}=x^4y^2+xz^3+y^3z.                       \label{eq:octic46:octic_q14_degree38_46:38-1}
\end{equation}
The original $\mu_{14}$ is the radial subgroup with parameter of character
$3$.  The exceptional curve of the weighted blowup is rational after this
quotient.  The $y$--axis and $z$--axis charts give respectively an $A_8$
chain and the chain
\begin{equation}
 (-2),\ (-3),\ (-2),\ (-2),                    \label{eq:octic46:octic_q14_degree38_46:38-2}
\end{equation}
listed from the main curve $E$ outwards.  We justify the remaining chart,
including its degenerations, below.

\paragraph{The uniform $x$--axis chart.}
Put $x=s^5$, $y=s^9u$, $z=s^{11}v$, and first divide out the radial group
by setting $S=s^{14}$.  The residual quotient is
\[
 \mu_5(S,u,v)=(1,4,1),
\]
and the strict equation begins
\[
 u^2+v^3+u^3v+S(bv^2+cu^3)+S^2dv+S^3e+\cdots.
\]
Equivariant analytic Morse splitting in $u$ changes this to
\begin{equation}
 u^2+g(v,S)=0,\qquad
 g=g_3+g_8+g_{13}+\cdots,\qquad
 g_3(v,S)=v^3+bSv^2+dS^2v+eS^3.               \label{eq:octic46:octic_q14_degree38_46:38-3}
\end{equation}
Here $g_j$ is homogeneous of degree $j$ in $(v,S)$; the nonzero coefficient
of $v^3$ is essential, but the binary cubic need not have distinct roots.
Normality of the original surface makes this double--cover equation normal.

Blow up this chart with weights $(u,v,S)=(3,2,2)$.  Its exceptional curve
$F$ has coarse normalization $\mathbb P^1$, with coordinate $v/S$.
In the $S$--chart write
\[
 S=t^2,\qquad v=t^2b,\qquad u=t^3a,\qquad T=t^5.
\]
After dividing by the radial $\mu_5$, the chart is
\begin{equation}
 \bigl\{a^2+g_3(b,1)+T^2g_8(b,1)+T^4g_{13}(b,1)+\cdots=0\bigr\}/\mu_2,
 \quad (T,a,b)\longmapsto(-T,-a,b).            \label{eq:octic46:octic_q14_degree38_46:38-4}
\end{equation}
The cover is normal, and the involution has no divisorial fixed locus on
it.  Its residue two--form has trivial character.  Hence \eqref{eq:octic46:octic_q14_degree38_46:38-4} is
Gorenstein.  It is a normal partial resolution of the rational surface
singularity $N$; every such partial resolution has rational singularities,
by the Leray spectral sequence applied to a resolution.  Consequently
all singularities in \eqref{eq:octic46:octic_q14_degree38_46:38-4} are rational Gorenstein surface singularities,
that is, Du Val singularities.  The proper old curve meets the other chart at a smooth point, as checked
below.  In particular every curve subsequently introduced
at $F$ is a $(-2)$--curve and these subsequent resolutions are crepant.

Before resolving those Du Val points, $F^2=-5/2$ and its discrepancy over
the $x$--axis quotient is $-4/5$.  Indeed the weights $(3,2,2)$, degree $6$,
give square $-1/2$ and discrepancy $0$ before the radial $\mu_5$ quotient;
the quotient multiplies the square by $5$ and changes the discrepancy to
$(0+1)/5-1$.  Thus
\[
 K\cdot F=(-4/5)(-5/2)=2.
\]
The coarse curve $F$ is smooth: in \eqref{eq:octic46:octic_q14_degree38_46:38-4}, at $T=0$, the quotient is the
$b$--line.  Its strict transform is still smooth after the Du Val
resolution.  Crepancy and adjunction therefore give
\begin{equation}
 F^2=-4                                                     \label{eq:octic46:octic_q14_degree38_46:38-5}
\end{equation}
on the chosen smooth resolution, uniformly in all higher coefficients.
For a cubic with distinct roots the attached Du Val curves are three
$A_1$'s; no such genericity is needed below.

The proper transform of the old main curve is $S=0$.  It meets $F$ at
$[u:v:S]=[\pm i:1:0]$, away from the roots of $g_3(v,S)$, since
$g_3(v,0)=v^3$.  In the $v$--chart put $v=t^2$, $u=t^3a$, $S=t^2b$.
The residual $\mu_2$ exchanges the two points $a=\pm i$, and after the
radial quotient $U=t^5$ one has
\[
 S^5=U^2b^5.
\]
Thus the self--intersection correction for the old main curve at this
chart is exactly $2/5$, independently of the Du Val resolutions away
from it.  The raw square of $E$ is $-532/495$; the other two corrections
are $8/9$ and $7/11$.  Therefore
\begin{equation}
 E^2=-\frac{532}{495}-\frac89-\frac7{11}-\frac25=-3.          \label{eq:octic46:octic_q14_degree38_46:38-6}
\end{equation}
The graph needed below consists of $E^2=-3$, its $A_8$ arm, the fixed arm
\eqref{eq:octic46:octic_q14_degree38_46:38-2}, and the adjacent $F^2=-4$ with some Du Val configuration attached.
This is a good resolution, and need not be made relatively minimal.

\paragraph{The actual two--form.}
Write the actual form as
\[
 \alpha=\delta_x\,dy\wedge dz-\delta_y\,dx\wedge dz
             +\delta_z\,dx\wedge dy,
\]
where the coefficient characters are $10,8,0$ modulo $14$ respectively.
If $\delta_z(0)\ne0$, the $y$--chart
\[
 y=s^9,\quad x=s^5u,\quad z=s^{11}v,\quad S=s^{14}
\]
has $v=-u^4+\cdots$, residual quotient
$\mu_9(S,u)=(5,4)\simeq\frac19(1,8)$, and
$dx\wedge dy$ is a nonzero constant multiple of $dS\wedge du$.
The entire $A_8$ arm is therefore protected.

If instead $\delta_z(0)=0$, the following order calculation applies.

The smallest remaining weighted degrees of the three summands of
$\alpha$ are respectively $42,42,28$.  Dividing by the radial group of
order $14$ shows that the coefficient $e$ of $E$ in the actual canonical
forest is a positive integer.  All eight coefficients $a_1,\ldots,a_8$
on its $A_8$ arm are then positive as well, since a zero $(-2)$--curve
cannot meet a positive canonical component.


\subsection{The normal fourteenth source: the degree-$48$ format}
\label{sec:octic-q14-degree48}

Suppose now that $[x^2]f=[x^3y]f=[x^4y^2]f=0$, whereas
$[xz^3]f[y^3z]f\ne0$.  The coarse surface germ is again rational.
We retain all coefficient specializations, including a reducible
initial exceptional divisor.

Use weights $(3,11,15)$.  The character--$2$ monomials of weight at
most $48$, under the original $\mu_{14}$ action $(1,13,5)$, are precisely
\[
 x^2,\ x^3y,\ x^4y^2,\ y^3z,\ xz^3,
 \ x^5y^3,\ x^6z^2,\ x^{11}z,\ x^{16}.
\]
The equivariant coordinate change $z\mapsto z+\lambda x^5$ eliminates
the $x^5y^3$ coefficient and does not introduce any of the three
forbidden lower terms.  After scaling, the entire weight--$48$ face is
\begin{equation}
 f_{48}=xz^3+y^3z+B x^6z^2+C x^{11}z+D x^{16}.              \label{eq:octic46:octic_q14_degree48_46:48-1}
\end{equation}
The original group is the radial $\mu_{14}$ with parameter of character
$5$.  We analyze its weighted blowup without imposing a nondegeneracy
condition on $(B,C,D)$ or on higher terms.

\paragraph{The three chart types.}
In the $x$--chart put
\[
 x=s^3,\quad y=s^{11}u,\quad z=s^{15}v,\quad S=s^{14}.
\]
The residual group and equation are
\begin{equation}
 \mu_3(S,u,v)=(2,1,0),\qquad
 vu^3+v^3+Bv^2+Cv+D+S G_1(u,v)+S^2G_2(u,v)+\cdots=0.      \label{eq:octic46:octic_q14_degree48_46:48-2}
\end{equation}
The equation is invariant.  The residue canonical character is
$2+1+0=0$ modulo $3$.  The group has no divisorial fixed locus on the
surface: its fixed line $S=u=0$ meets the surface at the finitely many
roots of $v^3+Bv^2+Cv+D$.  Thus this quotient is Gorenstein.

The weighted strict transform is normal.  It is Cohen--Macaulay and is
smooth at the generic points of all its exceptional components: when
$D\ne0$ the initial curve is reduced and generically smooth; when
$D=0$ it is the reduced union of the two distinct factors
$z$ and $y^3+xz^2+Bx^6z+Cx^{11}$.  The other points of codimension one
come from the original normal surface.  Serre's criterion applies.
Consequently \eqref{eq:octic46:octic_q14_degree48_46:48-2} is a normal partial resolution of the rational
surface $N$.  It is rational by Leray and therefore, being Gorenstein,
has only Du Val singularities.  Their resolutions are crepant and
introduce only $(-2)$--curves, even when the cubic in $v$ has repeated
roots or when the exceptional divisor is reducible.

The $y$--axis chart has $y=s^{11}$, $x=s^3u$, $z=s^{15}v$.
Its leading equation is linear in $v$, so after putting $S=s^{14}$ it
is the fixed quotient
\begin{equation}
 \mu_{11}(S,u)=(3,8)\simeq\tfrac1{11}(1,10),              \label{eq:octic46:octic_q14_degree48_46:48-3}
\end{equation}
with an $A_{10}$ resolution chain.  The $z$--axis chart is similarly
smooth before its residual quotient and gives
\begin{equation}
 \tfrac1{15}(1,11),\qquad
 (-2),(-2),(-3),(-2),(-2),                                \label{eq:octic46:octic_q14_degree48_46:48-4}
\end{equation}
listed from the main exceptional divisor.  These charts include every
possible nontrivial stabilizer locus of $\mathbb P(3,11,15)$: the
remaining index--$3$ points lie in \eqref{eq:octic46:octic_q14_degree48_46:48-2}.

\paragraph{The main exceptional curves have no $(-1)$ specialization.}
If $D\ne0$, the main curve $E$ is a coarse $\mathbb P^1$.  Indeed in
\eqref{eq:octic46:octic_q14_degree48_46:48-2} its quotient has coordinate $v$, with
$u^3=-(v^3+Bv^2+Cv+D)/v$; the weighted completion is rational and
smooth as a coarse curve.  Its raw square is
\[
 E^2_{\rm raw}=-14\frac{48}{3\cdot11\cdot15}
              =-\frac{224}{165}.
\]
The discrepancy of the weighted blowup is
\begin{equation}
 a=\frac{(3+11+15-48)}{14}-1=-\frac{33}{14}.               \label{eq:octic46:octic_q14_degree48_46:48-5}
\end{equation}
Thus $K_{\rm partial}\cdot E=16/5$.  The $A_{10}$ and all $x$--chart
resolutions are crepant.  For the chain \eqref{eq:octic46:octic_q14_degree48_46:48-4}, the coefficient of the
canonical correction at its first vertex is $1/5$.  Hence on the smooth
resolution
\begin{equation}
 K_W\cdot E=3,\qquad E^2=-5.                            \label{eq:octic46:octic_q14_degree48_46:48-6}
\end{equation}

If $D=0$, the main divisor is instead the reduced union
\[
 E_0=(z=0),\qquad
 E_1=(y^3+xz^2+Bx^6z+Cx^{11}=0).
\]
Their coarse curves are both smooth rational curves.  In the $x$--chart
the first has coordinate $u^3$ and the second coordinate $v$.
Let $E_{\rm tot}=E_0+E_1$.  Their weighted degrees are $15$ and $33$,
so
\[
 E_{\rm tot}\cdot E_0=-\frac{14}{33},\qquad
 E_{\rm tot}\cdot E_1=-\frac{14}{15}.
\]
Together with \eqref{eq:octic46:octic_q14_degree48_46:48-5} these give
$K_{\rm partial}\cdot E_0=1$ and
$K_{\rm partial}\cdot E_1=11/5$.
The $A_{10}$ chain is attached to $E_0$, while \eqref{eq:octic46:octic_q14_degree48_46:48-4} is attached to
$E_1$.  Resolving the latter subtracts $1/5$ from the second canonical
intersection.  All $x$--chart resolutions are again crepant.  Therefore
\begin{equation}
 E_0^2=-3,\qquad E_1^2=-4.                              \label{eq:octic46:octic_q14_degree48_46:48-7}
\end{equation}
These calculations apply also when the two curves meet in a Du Val
point of an arbitrarily special higher--order type.  Their strict
transforms stay smooth, and rationality makes the complete exceptional
divisor a tree with transverse intersections.  In particular every
exceptional curve in the resulting resolution has square at most $-2$.

\paragraph{Orders of the actual form.}
Write
$\alpha=\delta_xdy\wedge dz-\delta_ydx\wedge dz
 +\delta_zdx\wedge dy$, with coefficient characters $10,8,0$.
If $\delta_z(0)\ne0$, then in \eqref{eq:octic46:octic_q14_degree48_46:48-3} the term $dx\wedge dy$ is a unit
multiple of $dS\wedge du$.  Thus the $A_{10}$ arm is protected.

If $\delta_z(0)=0$, the minimum weighted degrees of the three actual
ambient two--form summands are respectively $56,42,28$.
After division by the radial group, every main component has positive
canonical coefficient.  This also holds in the reducible case: at a
generic point of either component, the pullback of the relevant ambient
form is regular, and every summand has radial order at least one.
A vanishing leading coefficient only increases that order.

The $A_{10}$ arm is therefore entirely positive and is attached to a
positive main curve.


\subsection{A common finite test for positive exceptional cores}
\label{sec:octic-q14-core-test}

After the exclusion of $J$, the remaining graphs admit one common
calculation.  We state it with an external-intersection bound so that
every later use can specify its own polar budget.

\begin{proposition}\label{prop:octic-q14-positive-cores}
Suppose that the main curve of one of the five graphs in
\eqref{eq:octic46:octic_q14_forest46:45-QG} has positive canonical
coefficient.  Let $v_i$ be the total intersection with intrinsic polar
components at the vertices of its connected positive exceptional core.
If $\sum v_i\le12$, the necessary canonical-forest conditions leave
only the following possibilities:
\begin{enumerate}
\item in $H$, the long arm is omitted, the coefficients are
$2$ at the centre, $(3,4,5)$ on $A_3$, and $1$ on $A_1$; a single
polar component of coefficient six meets the end of $A_3$;
\item in $K$, the final $(-4)$ and the $[3,2,2]$ arm are omitted;
the coefficients are $1$ at the centre, $2$ on the short arm and
$(2,3,4,5,6,7)$ on $A_6$, with external intersections three and eight.
The contribution three is a single coefficient-three polar component.
\end{enumerate}
In particular a bound $\sum v_i\le10$ leaves only the $H$ possibility.
\end{proposition}

\begin{proof}
Format $J$ is excluded by Corollary~\ref{cor:octic-q14-J},
so consider $G,H,K,L$.  Let $M$ be the negative intersection matrix
of a connected positive exceptional subtree at Q14.  Write $z_i>0$ for its canonical coefficients,
$k_i=-2-M_{ii}$, and $v_i\ge0$ for the intrinsic polar intersection at its
vertices.  An adjacent exceptional vertex not in the positive subtree
has coefficient zero, so adjunction reads
\begin{equation}
                     Mz+v=k.                              \label{eq:octic46:octic_q14_forest46:45-F1}
\end{equation}
If $j$ is such a zero boundary vertex, we also impose
\begin{equation}
             k_j-\sum_{i\text{ in the subtree}} M_{ji}z_i\ge0.
                                                               \label{eq:octic46:octic_q14_forest46:45-F2}
\end{equation}
Any further positive subtree beyond that boundary is ignored; doing so
weakens this necessary inequality and cannot cause a false exclusion.

The assumed bound is
\begin{equation}\label{eq:octic-core-bound}
                         \sum_i v_i\le12.
\end{equation}

We enumerate every nonnegative integral $v$ satisfying \eqref{eq:octic-core-bound}, compute
$z=M^{-1}(k-v)$ exactly, and retain only positive integral solutions of
\eqref{eq:octic46:octic_q14_forest46:45-F1}--\eqref{eq:octic46:octic_q14_forest46:45-F2}.  The following necessary canonical-forest conditions
are then imposed.

\begin{enumerate}
\item Every connected subtree has a unique minimum coefficient.  This
includes the requirement that adjacent positive coefficients differ.
It can be checked by taking, for each level $h$, connected components
of the vertices with coefficient at least $h$ and requiring at most
one vertex of coefficient exactly $h$ in each component.
\item Adjacent coefficients $a,b$ satisfy
$\gcd(a+1,b+1)=1$.  At a free blowup the shifted coefficients are
consecutive, and at a corner blowup the new shifted coefficient is
the sum of its two old shifted coefficients; induction proves this
condition.
\item A vertex has at most two smaller neighbors.  If it has two,
its coefficient is one plus their sum.  They are its surviving
creation parents; any neighbor created later has larger coefficient.
\item Every $v_i$ is partitioned in all possible ways into the positive
multiplicities of its attached intrinsic polar curves.  We require
that the enlarged tree, with those attached leaves, also satisfies the
unique-minimum rule and the shifted-coprimality rule.  Imposing the
parent rule only on the original exceptional subtree already suffices
for the exclusions below; ignoring further parent restrictions on
the new leaves makes the test weaker.
\end{enumerate}

To make this finite calculation explicit, clear $M^{-1}$ by
$|\det M|$ and evaluate $z=M^{-1}(k-v)$ for every nonnegative
integral vector in the stated simplex.  All tests are integral after
this denominator is cleared; no separate bound on $z$ is imposed.
For a subtree with $n$ vertices the number of $v$'s is exactly
$\binom{n+12}{n}$.  The corresponding exact computation is included in the source.
Omissions in the table below
refer to the corresponding vertex positions in each arm of \eqref{eq:octic46:octic_q14_forest46:45-QG}.
\begin{equation}
\begin{array}{c|c|c|c}
\text{format}&\text{possible main positive cores}&
 \text{numbers of enumerated }v&\text{survivors}\\\hline
G&\text{full}&293930&0\\
H&\text{two cores}&646646,\ 6188&0,\ 1\\
K&\text{four cores}
 &2704156,\ 1352078,\ 293930,\ 125970&0,0,0,1\\
L&\text{eight cores}
 &\text{listed below}&0\text{ in every case}.
\end{array}                                                   \label{eq:octic46:octic_q14_forest46:45-F4}
\end{equation}
The two $H$ cores are full and with its long arm omitted.  The four $K$ cores
are full, with its end $[4]$ omitted, with its $[3,2,2]$ arm omitted,
and with both omitted, in that order.  For $L$ the omitted branches
are independently the final $[4]$, the
entire $[3,2^6]$, and the final $[3]$.  Its eight enumeration counts
are
\[
 (5200300,2704156,18564,6188),\qquad
 (2704156,1352078,6188,1820).
\]
All possibilities for the main positive component are exhausted:
a zero $(-2)$ vertex cannot border a positive vertex, so an omission
can begin only at one of the specified non--2 vertices.

For the $H$ survivor, any partition of the external contribution six
into smaller positive parts repeats a minimum along a path to the
existing coefficients $1,2,3,4,5$.  It is therefore a single polar
component.  The same argument makes the contribution three in $K$
a single component.  Its total external contribution is eleven,
which proves the last assertion.
\end{proof}


\subsection{The protected boundary at a normal twelfth source}
\label{sec:q12-normal-boundary45}

Assume that the octic is normal at the ordinary Q12 source.
The source characters are $(x,y,z)=(1,11,5)$ and the equation
has character eight.

\begin{proposition}\label{prop:octic-normal-q12-boundary}
Suppose the coarse surface germ is rational.  The coefficients of
$x^3z$ and $z^4$ are nonzero.  After ordinary blowup and normalization,
the exceptional component through the $x$--axis quotient chart and
every other component of the complete minimal resolution have square
at most $-2$.  The $x$--axis chart supplies an $A_5$ chain.
If the actual ambient two-form has $\delta_z(0)=0$, its order on
the component attached to that chain is positive.  If
$\delta_z(0)\ne0$, the five $A_5$ curves and a disjoint $A_1$
curve have order zero.  Thus the fixed protected boundary is $A_5+A_1$.
\end{proposition}

\begin{proof}
For the residue $\Omega$, both $x^3\Omega$ and $z^3\Omega$ are
invariant regular canonical forms.  With weights $(1,2,2)$ the
first form has threshold eight, and the only equation monomial below
eight is $x^3z$.  With weights $(5,7,1)$ the second has threshold
sixteen, and the only equation monomial below sixteen is $z^4$.
The weighted canonical-form inequality therefore forces both terms.
The ordinary quartic initial equation is consequently
\[
 f_4=x^3z+c z^4+b y^2z^2+a y^4,\qquad c\ne0,
\]
after a nonzero scalar normalization.  This lists every ordinary
degree-four monomial of character eight; no lower-degree equation
term exists.

The original cyclic group has radial subgroup of order two on the
ordinary exceptional projective plane.  For an orbit of irreducible
components of total plane degree $d_i$, the intersection of the full
exceptional divisor with its quotient component is $-d_i/3$.
Indeed the ordinary blowup gives degree $-d_i$ before the source
quotient; the radial ramification contributes $2^2$ and the source
group has order twelve.  Intersections with other components and
the corrections from resolving quotient or higher singularities
are nonnegative.  Thus they can only decrease the self-intersection
of the final strict transform.

If $a\ne0$, the quartic is irreducible.  To see this, regard it as a
cubic in $x$ over $\mathbb C(y,z)$: the $z$--valuation of
$-(c z^4+b y^2z^2+a y^4)/z$ is $-1$, so that expression is not a
cube.  Its sole quotient component has raw square $-4/3$ and hence
integer final square at most $-2$.  This includes the specialization
$b^2=4ac$ and does not assume that the plane quartic is smooth.

If $a=0$, the exceptional equation is
$z(x^3+b y^2z+c z^3)$.  The line $z=0$ has raw full-divisor
intersection $-1/3$.  At $[1:0:0]$ the derivative in $z$ is a unit;
after removing the radial subgroup, the quotient chart is
$\frac16(1,5)$ and its $A_5$ correction is $5/6$.
Thus the line has final square at most $-1/3-5/6=-7/6$, hence at
most $-2$.

For $b\ne0$, the cubic is smooth and irreducible.  For $b=0$, its
three distinct lines are permuted transitively by the residual group,
so in either case their orbit has total degree three and raw
full-divisor intersection $-1$.  The three points with $y=0$ and
$x^3+c z^3=0$ are smooth points of the strict transform, permuted
transitively, with residual stabilizer of order two.  Its action on
the two local surface coordinates is $\frac12(1,1)$.
The resulting $A_1$ correction is $1/2$.  The quotient cubic-orbit
component therefore has final square at most $-3/2$, hence at most
$-2$.  All remaining singular charts are resolved minimally.
They are rational, as normal birational partial resolutions of the
rational original surface, and their minimal exceptional curves
have squares at most $-2$.

Finally the actual two-form has components of characters
$(\delta_x,\delta_y,\delta_z)=(8,6,0)$.
If $\delta_z(0)=0$, every component of the form has ordinary weighted
degree at least four.  On removing the radial subgroup by $S=s^2$,
its order on every main exceptional component is therefore at least
one.  If $\delta_z(0)\ne0$, the leading form at the $x$--axis chart
is a unit multiple of $dx\wedge dy$, which becomes a unit multiple
of $dS\wedge du$.  This is the canonical generator on
$\frac16(1,5)$, so all five $A_5$ orders are zero.

The three points $y=0$, $x^3+c z^3=0$ give the additional protected
$A_1$ for all $a,b$.  In the $x$--chart $x=s$, $y=su$, $z=sv$,
the derivative of the strict equation in $v$ at a root of
$1+c v^3$ is $-3$.  These points are smooth.  Their stabilizer
before removing the radial subgroup is $\mu_4$, with action
$(s,u)\mapsto(is,-u)$.  With $S=s^2$, its residual action is
$\mu_2(S,u)=(-S,-u)$.  The leading actual form is
\[
       \delta_z(0)\,dx\wedge dy
              =\tfrac12\delta_z(0)\,dS\wedge du.
\]
All other terms have larger ordinary weight.  This is a unit canonical
form on the $A_1$ quotient, so its exceptional curve has order zero.
The orbit is disjoint from the $x$--axis $A_5$ point and from every
degenerating chart in the $y$--direction.
\end{proof}

This is a local lemma only.  Its application to a global conductor state
still requires a justified polar coefficient bound or an independent
count of disjoint protected roots.


\subsection{C2-02350: the degree-seven conductor and its nodal charts}
\label{sec:octic-local45}

Assume that the conductor image is a sole degree-seven component
with generic delta one.  Lemma~\ref{lem:octic-degree7-stack46}
identifies its canonical-index stack as $\mathbb P(5,12)$ and gives
the full conductor algebra.  We compute the nodal normalization at Q12.

\paragraph{The forced nodal Q12 coefficients}

At the ordinary twelfth source straighten the smooth lifted conductor
as $(x,z)=0$ in coordinates of characters $(1,11,5)$.  The full
conductor ideal is exactly $(x,z)$, by purity and generic delta one.
The node algebra of Lemma~\ref{lem:octic-degree7-stack46} has
trace-free relation $T^2=cX_5^8Y_{12}$, so its determinant has order
exactly eight along the $y$-axis.  The normal quadratic part of $f$ has
coefficient orders at least six, ten, two for $x^2,xz,z^2$.
Consequently its determinant forces both $[y^6x^2]f$ and
$[y^2z^2]f$ to be nonzero.  Theorem
\ref{thm:weighted-canonical-inequality45} also forces $[x^3z]f$ and
$[z^4]f$ to be nonzero.  All four statements apply to the actual
local equation and allow arbitrary higher character terms.

\paragraph{The resulting local normalization is concrete}

A normalized representative of the forced leading format is
\begin{equation}
 f=y^6x^2+x^3z+y^2z^2+z^4.
 \label{eq:octic46:octic_local_46:45-O4}
\end{equation}
Set $A=y^6$, $B=x^2/2$, $C=y^2+z^2$, and adjoin a character-two variable
$t$.  The actual finite normalization is
\begin{equation}
 xt+Bx+Cz=0,\qquad zt-Ax-Bz=0,\qquad t^2-B^2+AC=0.
 \label{eq:octic46:octic_local_46:45-O5}
\end{equation}
These are the minors of
\[
 \begin{pmatrix}x&t-B&C\\z&A&-t-B\end{pmatrix}.
\]
The last equation proves finiteness; the first two identify the fraction
field.  The determinantal algebra is Cohen--Macaulay.  Away from the
origin it is the normalization of a smooth point or of the transverse
node, so it is regular in codimension one and hence normal.
On its conductor the equation is $t^2+y^8=0$.  The full cyclic action
exchanges the two branches, and its invariant ring is
$\mathbb C[[y^2t]]$, because $(y^2t)^2=-y^{12}$.  Thus the coarse
conductor is smooth and its double cover is simply branched at Q12.

The weighted blowup $(x,y,z,t)=(2,1,4,4)$ resolves, after taking the
source quotient, into the following graph.  There are two central
curves $E,F$ of square $-3$.  The curve $E$ has an $A_7$ arm, an $A_1$
arm and a single $-2$ bridge to $F$; the curve $F$ has further $A_1$
and $A_2$ arms.  There are fourteen vertices altogether.

Here are the local charts, included to make the graph verifiable.
At the x-chart the smooth coordinates have action $\frac1{24}(1,21)$;
removing the radial reflection of order three leaves
$\frac18(1,7)$, giving $A_7$.  The two finite points
$(u,\tau)=(0,\pm i)$ of $\tau^2=u^4/4-1$ are exchanged by the
source group and give one $A_1$.  At the z-chart, eliminating $\tau$
gives
\[
 v^2+u^3+u^2v^6+s^6=0.
\]
Removing the radial order-three reflection by $S=s^3$, and replacing
$v$ equivariantly by $v\sqrt{1+u^2v^4}$, gives an $A_2$ surface with
$\mu_{16}$ weights $(S,u,v)=(5,14,13)$.  Its universal plane cover is
\[
 S+iv=a^3,\quad S-iv=b^3,\quad u=-ab.
\]
The group of order 48 is generated by
$(a,b)\mapsto(\epsilon a,\epsilon^{-1}b)$ and
$(a,b)\mapsto(\zeta_{48}^5b,\zeta_{48}^{37}a)$.
Its scalar kernel has order eight and its projective quotient has order
six.  Blowing up gives a $-3$ centre with $A_2,A_1,A_1$ arms; the original
exceptional curve attaches at one $A_1$ arm.  This is the second centre
and bridge above.

The first centre has self-intersection $-15/16$ before these local
resolutions.  Subtracting $7/8$, $1/2$ and $11/16$ gives $-3$.
The reduced exceptional divisor is anti-nef with arithmetic genus zero;
its multiples have vanishing $H^1$, so the quotient normalization is
rational.  Exact inversion of its intersection matrix gives
\begin{equation}
 -\Delta^2=\frac{83}{30},\qquad
 r+\Delta^2=\frac{337}{30},\qquad
 q_{C}=\frac{23}{30},\quad \Delta_{C}=\frac{11}{15}.
 \label{eq:octic46:octic_local_46:45-O6}
\end{equation}
The conductor meets the finite $A_1$ arm.  The forced nonzero leading
coefficients keep these quotient charts and the $A_2$ chart unchanged
under the remaining higher character terms.

\paragraph{The actual Q12 form when the second jet is nonzero}

Assume in addition that $\gamma_{10}|_\Gamma\ne0$.  Its saturated order
two at Q12 forces $\delta_z(0)\ne0$.  The actual two-form then has the
same exceptional orders as $dx\wedge dy$; the other character components
have strictly higher relevant weight.  All fourteen exceptional orders
are zero.  On the first two quotient charts this follows from
$dx\wedge dy$ becoming a unit times $dS\wedge dv$ or $du\wedge dS$.
On the universal plane cover at the last chart its leading coefficient is
\[
 (a^6+4a^3b^3+b^6)\,da\wedge db.
\]
It is nonzero at the projective branch orbits $(a/b)^3=\pm1$ and $ab=0$,
so the orders on every arm are zero.  The twelve noncentral $-2$ curves
are therefore disjoint from the full canonical exceptional divisor and
survive unchanged on the K3 minimal model.

The nodal exclusion will show that its second polar is indeed
nonzero.  We have proved the local graphs and their actual-form
orders under this explicit hypothesis, without using that exclusion.


\subsection{C2-02350: the cuspidal twelfth source}
\label{sec:octic-cusp-q12-minimal45}

The following local calculation applies at the unshared Q12 source,
both for the sole degree-seven cusp and for the degree-seven component
of the two-component state.

\begin{theorem}
\label{thm:octic-cusp-q12-minimal45}
Suppose that the normalization of the octic has a K3 minimal model, and
that at an ordinary Q12 source its full conductor is the generic
delta--one cusp with ideal $(x,z)$, in source coordinates of characters
$(1,11,5)$.  Suppose that the actual second polar has order exactly two
on the $y$--axis.  Then, after an equivariant change straightening the
rank--one quadratic part, the four coefficients
\[
 [y^2z^2]f,\quad [x^3z]f,\quad [z^4]f,\quad [y^7x^3]f
\]
are nonzero.  The normalization has the sixteen-vertex resolution
described below.  The actual global two-form has order zero on all its
exceptional curves.  In particular the K3 minimal model contains the
protected root configuration
\[
                 A_8\oplus A_2\oplus A_1\oplus A_2,
                 \qquad\text{of rank }13.
\]
\end{theorem}

\paragraph{An exhaustive equivariant preparation}

The weighted canonical-form inequality already forces $y^2z^2,x^3z,z^4$.
Write the quadratic part along the conductor as
$A_0(y)x^2+2B_0(y)xz+C_0(y)z^2$.  The cusp condition is
$A_0C_0=B_0^2$.  Here $C_0=y^2\cdot\text{unit}$ and
$B_0\in y^{10}\mathbb C\{y^{12}\}$, so $B_0/C_0$ is regular.
The equivariant change $z\mapsto z-(B_0/C_0)x$ therefore removes both
the mixed and the pure-$x^2$ quadratic coefficients.  The latter change
does not alter any of the three already forced lowest coefficients.
Distribute the remaining terms as
\[
                  f=Ax^2+Cz^2,\qquad A\in(x,z),
                  \qquad C=y^2+\text{higher terms},
\]
allowing harmless nonzero scalar rescalings.  Adjoin $t$ of character two:
\begin{equation}
 xt+Cz=0,\qquad zt-Ax=0,\qquad t^2+AC=0.             \label{eq:octic46:octic_cusp_q12_46:45-C12-1}
\end{equation}
These finite birational determinantal equations give the actual
normalization.  They are Cohen--Macaulay and agree with the normalization
off the special source; their charts below have only isolated singularities,
which also verifies normality.

Blow up weights $(x,y,z,t)=(2,1,4,4)$.  On the $y$--chart put
\[
 y=s,\quad x=s^2u,\quad z=s^4w_0,\quad t=s^4\tau_0,
                  \qquad S=s^3.
\]
The transform of $C/s^2$ is a unit.  Thus the first equation in
\textup{\eqref{eq:octic46:octic_cusp_q12_46:45-C12-1}} eliminates $w_0$; at its lowest order,
$w_0=-u\tau_0$.  The remaining equation has lowest part
$\tau_0^2-u^2\tau_0$.  The radial subgroup of order three has been
removed by $S=s^3$, and the residual action is
\[
                 \mu_4(S,u,\tau_0)=(3,1,2).
\]
Equivariant Weierstrass preparation in $\tau_0$, completion of its square,
and preparation of the remaining factor in $u$ give
\begin{equation}
 F=\tau^2+u\{u^3+a(S)u^2+b(S)u+c(S)\}=0,             \label{eq:octic46:octic_cusp_q12_46:45-C12-2}
\end{equation}
where
\begin{equation}
 a(S)=S^3a_0(S^4),\qquad b(S)=S^2b_0(S^4),
 \qquad c(S)=S c_0(S^4).                              \label{eq:octic46:octic_cusp_q12_46:45-C12-3}
\end{equation}
These preparations use invariant units and keep $S,u$ fixed.  The
factor $u$ is retained because the proper conductor is $u=\tau=0$.
The lowest quartic coefficient is a nonzero scalar multiple of the square
of the forced $x^3z$ coefficient.  Thus the preparation is exhaustive,
not an assumption of a generic Newton face.  The hypersurface in
\textup{\eqref{eq:octic46:octic_cusp_q12_46:45-C12-2}} is normal: off $S=0$ it is an open set of the actual
normalization, and on $S=0$ it is smooth away from the origin.  In
particular $c$ cannot be identically zero, since otherwise it would be
singular along $u=\tau=0$.

The desired coefficient $[y^7x^3]f$ is, up to a nonzero scalar, $c_0(0)$.
We show that $c_0(0)=0$ is impossible.  There are precisely three
possibilities: $b_0(0)\ne0$; $b_0(0)=0$ but $[S^5]c\ne0$; or both
of these coefficients vanish.  The terms in \textup{\eqref{eq:octic46:octic_cusp_q12_46:45-C12-3}} cover
every higher specialization in these three cases.

\paragraph{The actual form, before any quotient resolution}

The exact order-two polar jet and the nonzero $y^2z^2$ coefficient give
$\delta_z(0)\ne0$.  Indeed all other quadratic contributions evaluated
on the normal derivation have order strictly larger than two.
The characters of $(\delta_x,\delta_y,\delta_z)$ are $(8,6,0)$.
Consequently, on the first $y$--chart,
\[
 \delta_z=A_1(S,u,w_0),\qquad
 \delta_x=s^4B_1(S,u,w_0),\qquad
 \delta_y=s^6C_1(S,u,w_0),\qquad A_1(0)\ne0,
\]
with regular coefficients.  Direct differentiation gives the following
identity for the actual global form, up to the common fixed sign convention:
\begin{align}
 \alpha={}&\frac{A_1}{3}\,du\wedge dS
   +\frac{S^2B_1}{3}\,dS\wedge dw_0\notag\\
 &+S^3C_1\left(\frac43w_0\,dS\wedge du
             -\frac23u\,dS\wedge dw_0\right)
   +S^4C_1\,dw_0\wedge du.                         \label{eq:octic46:octic_cusp_q12_46:45-C12-4}
\end{align}
After the preceding preparations,
$w_0\in(u\tau,u^3,S^3u^2)$, with analytic coefficients.
Formula \textup{\eqref{eq:octic46:octic_cusp_q12_46:45-C12-4}} will be used term by term.  In particular
we do not infer that an arbitrary reflexive two-form on a quotient
singularity is a K\"ahler two-form.

\paragraph{First stratum: a nonzero $S^2u^2$ coefficient}

Suppose $b_0(0)\ne0$ and write $\operatorname{ord}_S c=4k+1$, where
$k\ge1$.  Blow up weights $(S,u,\tau)=(1,1,2)$ and use
$S=\sigma,u=\sigma v,\tau=\sigma^2w$.  Removing the radial
subgroup of order two by $\rho=\sigma^2$ gives
\[
 w^2+v^4+\rho a_0(\rho^2)v^3+b_0(\rho^2)v^2
       +\rho^{2k-1}\widetilde c(\rho^2)v=0,
 \qquad b_0(0)\widetilde c(0)\ne0.
\]
The remaining involution is
$(\rho,v,w)\mapsto(-\rho,-v,w)$.  The parameter-dependent
equivariant Morse lemma in $v$, followed by multiplication by units,
puts this in the form
\begin{equation}
                 w^2+v'^2-\rho^{4k-2}=0.             \label{eq:octic46:octic_cusp_q12_46:45-C12-5}
\end{equation}
It is $A_{4k-3}$, and the remaining involution is free away from the
origin and preserves its canonical generator.  Its quotient is
$D_{2k+1}$, with the convention $D_3=A_3$.

Here $w_0=\sigma^3W(\rho,v,w)$.  Substitution in
\textup{\eqref{eq:octic46:octic_cusp_q12_46:45-C12-4}} gives
\begin{equation}
 \alpha=\text{unit}\,dv'\wedge d\rho
              +\rho^2\beta,                        \label{eq:octic46:octic_cusp_q12_46:45-C12-6}
\end{equation}
where $\beta$ is a regular K\"ahler two-form on the smooth
three-variable ambient germ restricted to \textup{\eqref{eq:octic46:octic_cusp_q12_46:45-C12-5}}.
For example $S^2dS\wedge dw_0$ becomes a multiple of
$\rho^2d\rho\wedge dW$; the other error terms have still larger
radial order.  All coefficients descend regularly, since they have even
$\sigma$-power after the radial involution.

Relative to the Du Val generator $\eta$ of \textup{\eqref{eq:octic46:octic_cusp_q12_46:45-C12-5}},
$dv'\wedge d\rho=2w\eta$.  Every regular K\"ahler two-form
has its coefficient in the Jacobian ideal.  Thus the canonical
coefficient of \textup{\eqref{eq:octic46:octic_cusp_q12_46:45-C12-6}} vanishes at the origin.  It is
invariant under the symplectic involution, so it vanishes at the
quotient point as well.  Every curve of the quotient Du Val resolution
therefore has positive actual $\alpha$-order.

For $k\ge2$ this is impossible: a $D$-configuration of $(-2)$-curves
on a surface birational to a K3 surface has order zero.  Indeed every
positive component is exceptional for the morphism to the K3 minimal
model, whereas the final $(-2)$-subgraph of a point-blowup forest is
a union of chains.

For $k=1$ the quotient is $A_3$.  On its smooth plane cover $(a,b)$,
the leading invariant $w$ is a nonzero multiple of $ab$, and the
quotient action is $\frac14(1,3)$.  The three toric valuations of
$ab$ are all one.  The error $\rho^2\beta$ has ordinary plane order
at least six, so its invariant divisorial orders on these curves are
at least two.  The nonzero leading coefficient cannot cancel.  Thus
the three adjacent $A_3$ curves all have actual order one, impossible
in a canonical exceptional forest.  This excludes the entire first
stratum.

\paragraph{Second stratum: the first remaining linear term is $S^5u$}

Suppose $b_0(0)=0$ and $[S^5]c\ne0$.  Blow up weights
$(S,u,\tau)=(3,5,10)$.  The initial equation is
\[
                  \tau^2+u^4+\lambda S^5u=0,
                  \qquad\lambda\ne0.
\]
The full source $\mu_4$ is radial.  In the $S$--chart put
$S=\sigma^3,u=\sigma^5v,\tau=\sigma^{10}w$ and
$\rho=\sigma^4$.  Near $v=w=0$ the equation solves
$v=w^2\cdot\text{unit}$; after this elimination the residual
chart action is $\frac13(1,2)$ on $(\rho,w)$.  Thus it gives an
$A_2$ chain attached to the main exceptional curve $\rho=0$.

The leading part of the actual form is
\[
 du\wedge dS=\text{unit}\,\sigma^7w\,d\sigma\wedge dw
            =\text{unit}\,\rho w\,d\rho\wedge dw.
\]
The two $A_2$ valuations of $\rho w$ are both one, as is the
coefficient on the main exceptional curve.  The terms of
\textup{\eqref{eq:octic46:octic_cusp_q12_46:45-C12-4}} cannot cancel this leading term: here
$w_0=\sigma^{15}W$, so already
$S^2dS\wedge dw_0$ has radial order at least five after
$\rho=\sigma^4$, and the other terms have higher order.  Thus
adjacent curves again have equal positive coefficients, a
contradiction.

\paragraph{Third stratum: rationality alone rules out the remaining jets}

Finally suppose $b_0(0)=[S]c=[S^5]c=0$.  Translate $u$ by
$a(S)/4$ and write the resulting equation as
\[
 F=\tau^2+u^4+A(S)u^2+B(S)u+C(S),
 \qquad
 \operatorname{ord}A\ge6,\quad
 \operatorname{ord}B\ge9,\quad
 \operatorname{ord}C\ge12.
\]
This translation preserves $du\wedge dS$ and the $\mu_4$ action.
The normal quotient of this $y$--chart has rational singularities:
it is a normal birational partial resolution of the rational
normalization of the octic.  The regular canonical form
\[
          S^2\operatorname{Res}\frac{dS\wedge du\wedge d\tau}{F}
\]
is invariant, since the residue and $S^2$ both have character two.
For weights $(S,u,\tau)=(1,3,6)$ its equation has order at least
twelve, while the canonical-form threshold is
$1+3+6+2=12$.  This contradicts the weighted canonical-form inequality.
The third stratum is therefore impossible too, proving $[S]c\ne0$.

\paragraph{The resulting protected graph}

With $[S]c\ne0$, the $y$--chart is equivariantly $A_1/\mu_4$, with
characters $(S,u,\tau)=(3,1,2)$.  Its quotient is
$\frac18(1,3)$, resolved by a chain $G(-3)-J(-3)$.
The two initial components $E_1,E_2$ meet $G$ at distinct points;
the proper reduced conductor meets $J$.  The $x$--chart is the
smooth cover of $\frac18(1,7)$ and gives an $A_7$ on $E_1$.
The $z$--chart is the same forced $A_2/\mu_{16}$ chart as in the
nodal calculation: it gives a centre $F(-3)$ with arms
$A_1,A_1,A_2$, one of the $A_1$ arms being the bridge $T$ to $E_2$.
Resolving the quotient corrections gives $E_1^2=E_2^2=-2$.
Thus the full graph is
\[
\begin{gathered}
 E_1(-2)\text{ with its }A_7\text{ arm, joined to }G(-3),\\
 G(-3)\text{ joined also to }J(-3)\text{ and }E_2(-2),\\
 E_2(-2)-T(-2)-F(-3),
 \quad F\text{ with its remaining }A_1,A_2\text{ arms}.
\end{gathered}
\]
The raw intersections before resolving the quotient points are
\[
 (E_1+E_2)E_1=-3/8,\quad (E_1+E_2)E_2=-9/16,
 \quad E_1E_2=3/8.
\]
The self-intersection corrections are $7/8+3/8$ for $E_1$ and
$3/8+11/16$ for $E_2$, giving the asserted final squares.

The actual form has zero orders on this graph.  In the $x$--chart its
leading term is a unit times the smooth symplectic generator; its
other terms have larger radial order.  In the $y$--chart the universal
plane cover may be written
$u=b^2$, $S'=a^2$, $\tau'=iab$, with $S=S'+u^3/4$ up to units.
Then $du\wedge dS$ is a unit times $ab\,da\wedge db$.
Its orders on the two rays of $\frac18(3,1)$ are zero.
Formula \textup{\eqref{eq:octic46:octic_cusp_q12_46:45-C12-4}} shows directly that its error terms have
strictly larger orders.  At the $z$--chart the same character estimate
and the previously computed plane cover give the leading coefficient
$a^6+4a^3b^3+b^6$, which is nonzero at every projective branch orbit;
all its resolution orders are zero as well.

The thirteen $(-2)$-curves consequently survive unchanged on the K3
minimal model.  This conclusion is local at Q12.  It does not assume
that the conductor is smooth in the ambient stack at another source,
does not rule out hidden Cartier ambient contacts, and does not by
itself determine the complete global conductor configuration.


\subsection{C2-02350: a multiplicity bound at the cuspidal fifth source}
\label{sec:octic-q5-cusp-multiplicity45}

The following local bound applies to either second-polar alternative.

\begin{theorem}
\label{thm:octic-q5-cusp-multiplicity45}
At an ordinary Q5 source let the cover coordinates have characters
$(x,y,z)=(1,4,2)$, let the octic equation have character four, and let its
full conductor be the generic delta--one cusp with ideal $(y,z)$.
Suppose the quotient normalization $N$ is rational.  Then
\[
                  \operatorname{mult}_P N\le10.
\]
If the $z^2$ coefficient is zero, the stronger bound
$\operatorname{mult}_P N\le9$ holds.  Every exceptional curve of the
minimal resolution of $P\in N$ consequently has self-intersection at
least $-10$.
\end{theorem}

\paragraph{The two quadratic alternatives are exhaustive}

Put $\Omega=\operatorname{Res}(dx\wedge dy\wedge dz/f)$.
Its character is three.  Since $z$ is in the full conductor,
$z\Omega$ is a regular invariant canonical form on the normalization.
The weighted canonical-form inequality with weights $(1,1,1)$ has
threshold four.  The only character-four monomials of degree less than
four in $(y,z)^2$ are $z^2$ and $xy^2$.  Thus
\begin{equation}
                     [z^2]f\ne0\quad\hbox{or}\quad[xy^2]f\ne0.
                                                        \label{eq:octic46:octic_q5_multiplicity46:45-Q5M-1}
\end{equation}

\paragraph{The double-plane alternative}

If $[z^2]f\ne0$, the equivariant Morse lemma with parameters $x,y$
puts $f$ in the form $z^2+g(x,y)$.  The rank--one quadratic part along
the cusp means $g\in(y^3)$.  Adjoining $t=z/y$ gives the actual
normalization cover
\begin{equation}
 B=\mathbb C\{x,y,t\}/(t^2+G(x,y)),\qquad G=g/y^2\in(y),
 \qquad\mu_5(x,y,t)=(1,4,3).                           \label{eq:octic46:octic_q5_multiplicity46:45-Q5M-2}
\end{equation}
It has multiplicity two.  This finite birational ring is normal: it is
the normalization at the general cusp point, agrees with the original
normal surface off the conductor, and has no nonregular codimension-one
points after shrinking the germ.

We record the elementary finite-quotient estimate used here.  If a finite
group of order $r$ acts faithfully on a two-dimensional local domain $B$
over $\mathbb C$, with $A=B^G$, then
\begin{equation}
                     e(A)\le r\,e(B).                 \label{eq:octic46:octic_q5_multiplicity46:45-Q5M-3}
\end{equation}
For every $b\in\mathfrak m_B$, the orbit polynomial
$\prod_{g\in G}(T-gb)$ gives $b^r\in\mathfrak m_A B$.
Polarization in characteristic zero gives
$\mathfrak m_B^r\subseteq\mathfrak m_A B$.  Monotonicity and the
finite-module rank formula for Hilbert--Samuel multiplicity give
\[
 r e(A)=e(\mathfrak m_A B,B)
       \le e(\mathfrak m_B^r,B)=r^2e(B).
\]
Taking $r=5$ and $e(B)=2$ proves the bound ten in this alternative.

\paragraph{The determinantal alternative}

Suppose $[z^2]f=0$.  By \textup{\eqref{eq:octic46:octic_q5_multiplicity46:45-Q5M-1}}, $[xy^2]f\ne0$.
On the conductor the quadratic coefficient of $y^2$ is $x$ times a
unit, and the mixed coefficient is divisible by $x^3$.  The rank--one
condition permits a regular equivariant $y$-shift which removes both
the mixed and the kernel quadratic terms along the conductor.
Every remaining term of transverse order at least three is divisible
by $y^2$ or by $z^2$.  Thus write $f=A y^2+Cz^2$ with
$\partial A/\partial x(0)\ne0$ and $C\in(y,z)$.
Replace $x$ by the character-one coordinate $A$.  We obtain
\[
                  f=xy^2+Cz^2,\qquad C\in(y,z),
                  \operatorname{char}(C)=0.
\]
No variable has character zero, so $\operatorname{ord}C\ge2$.
The actual normalization cover is the finite birational
Hilbert--Burch ring
\begin{equation}
 B=\frac{\mathbb C\{x,y,z,t\}}
   {(yt+Cz,\ zt-xy,\ t^2+xC)},\qquad
            \mu_5(x,y,z,t)=(1,4,2,3).                 \label{eq:octic46:octic_q5_multiplicity46:45-Q5M-4}
\end{equation}
The same codimension-one argument verifies its normality; the
determinantal presentation supplies the Cohen--Macaulay condition.

Its tangent cone is
\begin{equation}
 B_0=\frac{\mathbb C[x,y,z,t]}{(yt,\ zt-xy,\ t^2)},
               \qquad H_{B_0}(s)=\frac{1+2s}{(1-s)^2}.
                                                        \label{eq:octic46:octic_q5_multiplicity46:45-Q5M-5}
\end{equation}
For clarity, the three displayed initial quadrics form a standard
basis.  Their two Hilbert--Burch syzygies lift to the exact syzygies
\[
 t(yt+Cz)-C(zt-xy)-y(t^2+xC)=0,
\]
\[
 x(yt+Cz)+t(zt-xy)-z(t^2+xC)=0.
\]
Their leading syzygies generate the syzygy module of the initial
quadrics, proving \textup{\eqref{eq:octic46:octic_q5_multiplicity46:45-Q5M-5}}.  In particular $B_0$ is
Cohen--Macaulay of multiplicity three.

The following two invariant parameters give a sharper estimate than
the general quotient bound:
\begin{equation}
                  p=xz^2+y^2z,\qquad q=x^5+y^5+z^5.
                                                        \label{eq:octic46:octic_q5_multiplicity46:45-Q5M-6}
\end{equation}
Their ordinary degrees are three and five.  The reduced support of
$B_0$ is the union of the planes $(x,t)=0$ and $(y,t)=0$.
On the first plane, $p=y^2z$, $q=y^5+z^5$; on the second,
$p=xz^2$, $q=x^5+z^5$.  On either plane they vanish simultaneously
only at the origin.  Since $B_0$ is Cohen--Macaulay, their initial
forms are a homogeneous regular sequence.  The filtered regular
sequence criterion therefore gives
\begin{equation}
 \operatorname{length}B/(p,q)
 =\operatorname{length}B_0/(p,q)
 =3\cdot5\cdot e(B_0)=45.                            \label{eq:octic46:octic_q5_multiplicity46:45-Q5M-7}
\end{equation}
With $A=B^{\mu_5}$, the finite-module rank formula gives
\[
                  e((p,q),A)=45/5=9.
\]
Since $(p,q)\subseteq\mathfrak m_A$, monotonicity yields
$e(A)\le9$, as asserted.  This proof does not need a semicontinuity
argument or a genericity assumption on $C$.

\paragraph{From multiplicity to exceptional self-intersections}

Let $Z=\sum z_iE_i$ be the fundamental cycle on the minimal resolution
of the rational singularity $P\in N$, and put $E_i^2=-b_i$.
All $z_i\ge1$ and $b_i\ge2$.  Artin's multiplicity formula and
$p_a(Z)=0$ give
\begin{equation}
 e(N)=-Z^2=2+K_WZ=2+\sum_i z_i(b_i-2)\ge b_i.
                                                        \label{eq:octic46:octic_q5_multiplicity46:45-Q5M-8}
\end{equation}
This is the claimed uniform bound $b_i\le10$.
We use Artin's multiplicity formula in the form
\cite[Corollary~6]{ArtinRational1966}.


\subsection{C2-02350: the fifth source with a nonzero second polar}
\label{sec:octic-q5-cusp-nonzero45}

\begin{theorem}
\label{thm:octic-q5-cusp-nonzero45}
At the ordinary Q5 source of the degree--seven delta--one cusp,
suppose that the actual second polar is a unit and that the octic
normalization has a K3 minimal model.  Then its local normalization
is equivariantly
\begin{equation}
     (t^2+x^2y+y^4=0)/\mu_5(1,4,3),\qquad z=yt.
                                                        \label{eq:octic46:octic_cusp_q5_unit46:45-Q5P-1}
\end{equation}
The minimal resolution has a centre of square $-3$, two $A_1$ arms,
and one further $-3$ arm met by the reduced conductor.  Every exceptional
curve has actual two-form order zero.  In particular its two $A_1$
curves are protected on the K3 minimal model.
\end{theorem}

The unit polar forces $[z^2]f\ne0$ and $\delta_z(0)\ne0$.
Indeed the other normal derivation component has character two and
restricts to the $x$--axis with order at least two; all its quadratic
contributions have positive order.  The equivariant Morse preparation
therefore gives the double-plane normalization
\[
 F=t^2+G(x,y),\qquad G\in(y),\qquad
 \operatorname{char}(G)=1,\qquad \mu_5(x,y,t)=(1,4,3).
\]
Its residue has character two, so $x^3\operatorname{Res}(dx\wedge
dy\wedge dt/F)$ is a regular invariant canonical form.  For weights
$(x,y,t)=(1,4,8)$, the weighted canonical threshold is sixteen.
The only allowed monomials of smaller weight in $G$ are
\begin{equation}
                        x^2y,\qquad x^7y,\qquad x^3y^2.
                                                        \label{eq:octic46:octic_cusp_q5_unit46:45-Q5P-2}
\end{equation}
At least one is therefore nonzero.

\paragraph{The $x^3y^2$ alternative is incompatible with the actual form}

Suppose $[x^2y]G=0$ and $[x^3y^2]G\ne0$.  Blow up weights
$(x,y,t)=(2,3,6)$.  In the $x$--chart put
$x=s^2,y=s^3v,t=s^6w$ and remove the radial $\mu_5$ by $S=s^5$.
The remaining action is
\[
                         \mu_2(S,v,w)=(-S,-v,w).
\]
The equation has a nonzero $v^2$ coefficient.  Its coefficient linear
in $v$ is $S^{2\ell-1}$ times a unit for some $\ell\ge1$, because
the nonzero coefficient of $y$ in $G$ has smallest exponent
$x^{2+5\ell}$.  That coefficient cannot vanish identically: the
generic conductor is an ordinary cusp.  Equivariant Morse preparation
in $v$ consequently gives
\begin{equation}
                         w^2+v'^2-S^{4\ell-2}=0.
                                                        \label{eq:octic46:octic_cusp_q5_unit46:45-Q5P-3}
\end{equation}
Its quotient is $D_{2\ell+1}$, with $D_3=A_3$.

Here is the actual-form interface.  On $z=yt$, one has
\begin{equation}
 \alpha=(\delta_z-t\delta_y)\,dx\wedge dy
       +y\delta_x\,dy\wedge dt+y\delta_y\,dt\wedge dx.
                                                        \label{eq:octic46:octic_cusp_q5_unit46:45-Q5P-4}
\end{equation}
The source characters give
$\delta_x=s^3B$, $\delta_y=s^4C$ and $\delta_z=A$ with $A(0)\ne0$.
Direct substitution in \textup{\eqref{eq:octic46:octic_cusp_q5_unit46:45-Q5P-4}} gives, after the radial quotient,
\[
                 \alpha=\text{unit}\,dS\wedge dv'
                                  +S^2\beta,
\]
where $\beta$ is a regular K\"ahler two-form.  For example the
$y\delta_x\,dy\wedge dt$ terms have factors $s^{14}ds$ or
$s^{15}$, which become $S^2dS$ or $S^3$; the remaining errors have
at least the same order.  The leading coefficient relative to the
Du Val canonical generator is a unit times $w$.

Every exceptional order of the quotient Du Val configuration is
positive.  For $\ell\ge2$ its $D$-fork contradicts the canonical
forest theorem.  For $\ell=1$, its smooth plane cover has action
$\frac14(1,3)$ and the leading coefficient is $ab$; all three $A_3$
orders are one.  The $S^2\beta$ terms have strictly larger toric
orders, so adjacent equal positive coefficients give the contradiction.

\paragraph{The remaining $x^7y$ alternative gives an adjacent equal pair}

Suppose $[x^2y]G=[x^3y^2]G=0$.  Then \textup{\eqref{eq:octic46:octic_cusp_q5_unit46:45-Q5P-2}} forces
$[x^7y]G\ne0$.  Blow up weights $(3,7,14)$ and use
$x=s^3,y=s^7v,t=s^{14}w$ in the $x$--chart.  Again the full source
$\mu_5$ is radial; set $S=s^5$.  The equation has the form
\[
                         w^2+v+\lambda v^4+\cdots=0.
\]
The linear derivative in $v$ is nonzero, even when $\lambda=0$.
After eliminating $v$, the residual action is
$\frac13(2,1)$ on $(S,w)$.  It gives an $A_2$ chain attached to the
main exceptional curve $S=0$.

The leading term of the actual form is a unit times
\[
 dx\wedge dy=\frac35S\,dS\wedge dv
                         =\text{unit}\,S w\,dS\wedge dw.
\]
Both $A_2$ curves and the main exceptional curve thus have order one.
For the error estimate, the source characters give
$\delta_x=s^7B$, $\delta_y=s^6C$ and $\delta_z=A$, with $A$ a unit.
In \textup{\eqref{eq:octic46:octic_cusp_q5_unit46:45-Q5P-4}}, the term $t\delta_y dx\wedge dy$ has an
extra factor $S^4$, and the two remaining terms have radial orders
at least six and five.  They cannot cancel the displayed order-one
terms.  This contradicts the canonical forest.  We have proved
\[
                              [x^2y]G\ne0.
\]

\paragraph{The complete $D_{5m}/\mu_5$ family and its actual coefficients}

The equivariant parameter Morse lemma in $x$, keeping $y$ fixed, now
puts the normal cover in the form
\begin{equation}
                     t^2+x^2y+y^{5m-1}=0,
                     \qquad m\ge1.                    \label{eq:octic46:octic_cusp_q5_unit46:45-Q5P-5}
\end{equation}
The pure $y$ term is nonzero because the normalization is normal.
All its higher unit factors are absorbed equivariantly.  This is
$D_{5m}$ divided by the given $\mu_5$.

Set $h=5m-2$.  Its universal plane cover can be chosen so that
\[
       x=a^{2h}+b^{2h},\qquad y=(ab)^2,
       \qquad t=ab(a^{2h}-b^{2h}),
\]
up to nonzero constants.  The binary dihedral group is enlarged by
the scalar $\mu_5$.  The total scalar kernel has order ten and the
projective group has order $2h$.  Blowing up the plane origin gives a
central curve with preliminary square $-5/h$, two order-two branch
orbits and one order-$h$ orbit.  The first two arms are $A_1$.
For $m=1$, the last arm is a single $-3$ curve and the final central
square is $-3$.  For $m\ge2$, the last arm has continued fraction
\[
 \frac{5m-2}{5m-7}
                =[\underbrace{2,\ldots,2}_{m-2},3,3],
\]
and the central square is $-2$.  These statements also follow directly
by subtracting the two corrections $1/2$ and the last correction
$(h-5)/h$ from $-5/h$.

The reference form is
\[
       dx\wedge dy
          =\text{nonzero constant}\,ab(a^{2h}-b^{2h})\,da\wedge db.
\]
After the scalar quotient its order on the centre is $m-1$.
At either order-two orbit it has order $\lfloor m/2\rfloor$ on the
$A_1$ arm.  To see the parity explicitly, after the scalar quotient
write the local coordinates as $(N,u)$ with residual action
$\frac12(1,1)$.  The coefficient is a unit times $N^{m-1}$ for odd
$m$, and a unit times $N^{m-1}u$ for even $m$.

These are also the orders of the actual form.  In the universal
plane grading, $x,y,t,z$ have respective degrees
$10m-4,4,10m-2,10m+2$.  The terms of
\textup{\eqref{eq:octic46:octic_cusp_q5_unit46:45-Q5P-4}} other than its unit multiple of $dx\wedge dy$
have total degree at least $10m+10$.  After the scalar quotient their
coefficients are divisible by $N^m$.  For odd $m$ this gives strictly
larger arm valuation; for even $m$ a possible same-valuation term is
proportional to $N^m$, and cannot cancel the distinct leading monomial
$N^{m-1}u$ identically along the exceptional curve.

If $m=2$, the centre and both $A_1$ arms have coefficient one,
which is impossible.  If $m\ge3$, both arm coefficients are smaller
than the central coefficient.  They are two lower parents of that
central curve in the canonical forest, so the blowup rule would require
\[
                  m-1=1+2\lfloor m/2\rfloor.
\]
This equality is false for both parities.  Hence $m=1$.

In that case the central and the two $A_1$ orders are zero.  At the
remaining $\frac13(1,1)$ chart the leading form is
$v\,dN\wedge dv$; its exceptional order is zero as well.
The higher terms have larger radial order and cannot cancel this
leading monomial.  This proves \textup{\eqref{eq:octic46:octic_cusp_q5_unit46:45-Q5P-1}} and the full
actual-order statement.

The theorem makes no assertion that the reduced conductor has no other
singular normalization contacts.  In particular no global
self-intersection for that conductor is inferred here.


\section{The conductor states outside higher delta in degree seven}
\label{sec:octic-eight-states-revised}

We now exclude eight of the eleven total conductor states.  A state
records all conductor components, so the pair $(7,1)+(18,1)$ counts
as one state.  Each argument uses the local charts of
Section~\ref{sec:octic-atlas-revised}, together with the polar degree
and source incidences specific to that state.  We retain both the
nodal and cuspidal alternatives whenever the census permits them.
In particular the state $(14,2)$ still has a nonreduced conductor;
it is included here because its local square-lift argument is independent
of the degree-seven conductor algebras.

\subsection{C2-02350: the nodal state $(7,1)$}
\label{sec:octic-nodal-gamma-zero}

We first exclude a vanishing second polar for the degree-seven node,
using its forced Q12 normalization and the K3 model.

\begin{theorem}
In the sole degree-seven, generically nodal conductor state of C2-02350,
the second polar jet cannot vanish identically on the conductor:
\begin{equation}
                    \gamma_{10}|_\Gamma\ne0.
\label{eq:octic46:octic_node_gamma0_46:45-G0}
\end{equation}
\end{theorem}

\paragraph{The zero second jet would make the two-form vanish along the conductor}

The actual conductor stack is
\(\mathcal G=\mathbb P(5,12)\), with \(H=\mathcal O(1)\),
and its normalized conductor is a connected double cover.  Thus the
transverse binary quadratic form at the generic point of \(\Gamma\)
is anisotropic over \(\mathbb C(\Gamma)\).  If
\(\gamma_{10}|_\Gamma=0\), the normal component of the polar
vector field is isotropic for that form and must vanish.  The vector
field is tangent to \(\Gamma\).  Its induced tangential part belongs to
\[
 H^0(\mathcal G,T_{\mathcal G}\otimes H)
       =H^0(\mathbb P(5,12),\mathcal O(18))=0.
\]
Consequently it vanishes on the conductor.  This argument uses the
actual derivation and the ideal of the conductor; it remains valid at
hidden Cartier ambient singular contacts.

Write \(R=C_N\) for the reduced conductor on the normalization.
The restricted two-form then vanishes along \(R\), so
\begin{equation}
 B_N=rR+B',\qquad r\geq1,\qquad R\not\subseteq\operatorname{Supp}B'.
\label{eq:octic46:octic_node_gamma0_46:45-G1}
\end{equation}
The degree data are
\begin{equation}
 420H B_N=58,\qquad 420H R=14.
\label{eq:octic46:octic_node_gamma0_46:45-G2}
\end{equation}
Thus \(1\leq r\leq4\).  Every component of \(B'\) has positive
scaled degree at least three.

\paragraph{The entire twelfth-source resolution would lie in the canonical forest}

The complete forced resolution graph has two \((-3)\)-centres \(E,F\).
The centre \(E\) has a seven-curve \((-2)\)-arm
\(A_1,\ldots,A_7\), a one-curve \((-2)\)-arm \(A\), and a
\((-2)\)-bridge \(T\) to \(F\).  The centre \(F\) has a one-curve
\((-2)\)-arm \(U\) and a two-curve \((-2)\)-arm \(V_1,V_2\).
The strict conductor meets \(A\) transversely.  This is the full
fourteen-vertex resolution graph, not a chosen subgraph.

We first show that the actual two-form has order at least three along
\(E\).  At the ordinary source use cover coordinates of characters
\((1,11,5)\), with conductor \((x,z)=0\).  Write
\[
 \alpha=\delta_z\,dx\wedge dy+
         \delta_x\,dy\wedge dz+
         \delta_y\,dz\wedge dx.
\]
All three coefficients lie in \((x,z)\).  For the weighted blowup
with weights \((2,1,4)\), the only possible two-form terms of weight
at most nine are
\begin{equation}
       cxy\,dx\wedge dy,
       \qquad d x^2y^2\,dx\wedge dy.
\label{eq:octic46:octic_node_gamma0_46:45-G3}
\end{equation}
They have weights six and nine and would give orders one and two on
\(E\), respectively.  The other two coefficient characters are eight
and six, and their terms have total two-form weight at least twelve.
There is no weight-three term, since the vector field vanishes on the
conductor.

On the first arm chart set \(x=s^2,y=sv\), and remove the radial
order-three subgroup by \(S=s^3\).  The reference two-form
\(\alpha_0=dx\wedge dy\) is a unit multiple of
\(dS\wedge dv\) on the smooth cover of the \(A_7\) quotient.
The two terms in \eqref{eq:octic46:octic_node_gamma0_46:45-G3} give, respectively, leading coefficients
\((Sv)^k\), \(k=1,2\).  The remaining coefficient is divisible by
\(S^{k+1}\) and is regular on this smooth cover.  At the first
exceptional ray \((7/8,1/8)\), its valuation is at least
\[
       \left\lceil\frac{7(k+1)}8\right\rceil>k.
\]
Thus \(E\) and its adjacent arm curve would have the same positive
canonical coefficient \(k\).  This is impossible in a point-blowup
canonical forest: every edge joins a curve to a parent with strictly
smaller coefficient.  Hence
\begin{equation}
                         e:=\operatorname{ord}_E(\alpha)\geq3.
\label{eq:octic46:octic_node_gamma0_46:45-G4}
\end{equation}

A \((-2)\)-curve of zero two-form order is disjoint from the effective
canonical divisor by adjunction.  Therefore all the arm curves adjacent
to \(E\), and the bridge \(T\), have positive order.  If \(F\)
had order zero, adjunction at \(F\) would bound the bridge coefficient
by one, since \(K_WF=1\).  Adjunction at the bridge would instead give
\(2\operatorname{ord}_T(\alpha)\geq e\geq3\).  This contradiction
shows that \(F\), and consequently both its arms, also have positive
order.  All fourteen curves therefore belong to the canonical
exceptional forest of the morphism to the K3 minimal model.

\paragraph{A forest cannot fit the remaining polar degree}

Let the coefficients of \(E,F,A,T,U,V_1,V_2\) be
\(e,f,a,t,u,v_1,v\), respectively, and those of the seven-curve
arm be \(m_1,\ldots,m_7\).  Put \(\ell=m_7\).
All these numbers are positive integers.

An intrinsic polar curve other than \(R\) can meet this connected
fourteen-curve forest at most once.  Otherwise adjoining that curve
would create a cycle in the point-blowup exceptional graph.  Such
intersections are transverse because that graph has simple normal
crossings.  Hence the sum of its contributions to local adjunction is
at most its coefficient in \(B'\).  Since each such curve has scaled
\(H\)-degree at least three, the sum of all external contributions is
at most
\begin{equation}
                        \left\lfloor\frac{58-14r}{3}\right\rfloor.
\label{eq:octic46:octic_node_gamma0_46:45-G5}
\end{equation}

Sum the adjunction equations over the fourteen curves.  Internal edges
cancel the valency terms.  The two \((-3)\)-centres contribute two,
and the four terminal vertices contribute
\(\ell+a+u+v\).  The conductor contribution \(r\) has already been
separated.  Thus the external sum is exactly
\[
                         2+\ell+a+u+v-r.
\]
Combining with \eqref{eq:octic46:octic_node_gamma0_46:45-G5} gives
\begin{equation}
 \ell+a+u+v
 \leq\left\lfloor\frac{52-11r}{3}\right\rfloor\leq13.
\label{eq:octic46:octic_node_gamma0_46:45-G6}
\end{equation}

The adjunction equations also imply
\begin{equation}
 a\geq\frac{e+r}{2},\quad
 u\geq\frac f2,\quad
 v_1\geq\frac{2f}{3},\quad
 v\geq\frac f3,\quad
 t\geq\frac{e+f}{2},\quad
 t\leq1+\frac{11f}{6}.
\label{eq:octic46:octic_node_gamma0_46:45-G7}
\end{equation}
For the two-curve arm, the two individual inequalities are
\(2v_1\geq f+v\) and \(2v\geq v_1\); these give the displayed
bounds.  In particular,
\begin{equation}
                         f\geq\frac{3(e-2)}8.
\label{eq:octic46:octic_node_gamma0_46:45-G8}
\end{equation}

By the canonical-forest chain criterion, the positive sequence
\(m_1,\ldots,m_7\) is strictly monotone.
If it increases, adjunction at \(A_1\) gives
\(m_1\geq e+1\), and therefore \(\ell\geq e+7\).
Together with \eqref{eq:octic46:octic_node_gamma0_46:45-G4}, \eqref{eq:octic46:octic_node_gamma0_46:45-G7}, and \(u,v\geq1\), this gives
\[
 \ell+a+u+v\geq(e+7)+\left\lceil\frac{e+1}{2}\right\rceil+2
                  \geq14,
\]
contradicting \eqref{eq:octic46:octic_node_gamma0_46:45-G6}.

Thus the arm decreases.  If \(e<m_1\), then \(A_1\) has two
smaller adjacent parents, so
\(m_1=1+e+m_2\).  Discrete concavity along the \((-2)\)-chain
makes its successive positive drops nondecreasing.  Each of its six
drops is therefore at least \(e+1\), and
\(m_1\geq6e+7\).  But adjunction at \(E\) gives
\(m_1\leq3e+1-a-t\leq3e-1\), impossible.
Adjacent equal coefficients are also impossible.  Consequently
\begin{equation}
                         e>m_1\geq7,\qquad e\geq8.
\label{eq:octic46:octic_node_gamma0_46:45-G9}
\end{equation}

Suppose \(a>e\).  Both its adjacent curves \(E\) and \(R\)
have smaller coefficient, since \(r\leq4<e\).  The blowup-parent
rule then gives \(a=1+e+r\geq10\).  Equation \eqref{eq:octic46:octic_node_gamma0_46:45-G8} gives
\(f\geq3\), so \(u\geq2\) and \(v\geq1\).  Therefore
\(\ell+a+u+v\geq1+10+2+1=14\), again impossible.
Equality \(a=e\) is forbidden.  We must have \(a<e\).

Now \(E\) has the two smaller adjacent curves \(A_1\) and \(A\)
as its parents.  Thus
\[
 e=1+m_1+a\geq8+\frac{e+1}{2},
 \qquad\text{and hence }e\geq17.
\]
Finally \eqref{eq:octic46:octic_node_gamma0_46:45-G7}--\eqref{eq:octic46:octic_node_gamma0_46:45-G8} imply
\[
 \ell+a+u+v
 \geq1+\frac{e+1}{2}+\frac{5(e-2)}{16}
 =\frac{13e+14}{16}>13.
\]
This final contradiction proves the theorem.


\paragraph{The five isolated fourteenth-source formats.}
\label{sec:octic-q14-node-cores}

All results here are for the sole degree--seven conductor with generic
$\delta=1$ and nodal transverse type, after the zero-second-jet branch
has been excluded.  The Q12 resolution then has zero actual two--form
order on all fourteen curves, including its two $(-3)$ centres; its
other twelve curves are protected roots.  At Q5 both curves of the
$[-3,-2]$ chain also have zero order.  The actual intrinsic polar divisor
has scaled degree $58$.

\paragraph{The Q12--Q5 curves after contraction.}
A zero-order $(-2)$ curve is disjoint from the canonical forest.  Each
of the three zero-order $(-3)$ centres at Q12 and Q5 has canonical
intersection $1$, so it meets precisely one coefficient-one canonical
root and no other canonical component.  Contracting the forest changes
its square to $-2$.  Several of these three centres can meet the same
root.  Their image intersection matrix is obtained from the disjoint
source matrices by adding a matrix of ones on each block of a partition
of the three centres.  All five partitions have rank $16$; their
nonzero determinants are
\begin{equation}
 -120,\quad -1272,\quad -248,\quad -168,\quad -2600.
                                                        \label{eq:octic46:octic_q14_forest46:45-Q16}
\end{equation}
Thus the Q12--Q5 curves give sixteen independent classes on the K3
surface.  We compare them with the protected curves at Q14.

\paragraph{The polar bound for the common core test.}
An intrinsic polar component of multiplicity at least two avoids
$Q_5,Q_{12}$ and has scaled degree at least five.  A coefficient-one
component has degree at least three.  At most one such component
can meet a connected positive core, by the unique-minimum rule,
and no polar curve can meet the same core twice.  Thus
\begin{equation}
 5\sum_i v_i-2\le58,\qquad \sum_i v_i\le12.
 \label{eq:octic46:octic_q14_forest46:45-F3}
\end{equation}
Proposition~\ref{prop:octic-q14-positive-cores} now applies.

\paragraph{Formats H and L}

With $J$ already excluded by Corollary~\ref{cor:octic-q14-J},
the local atlas makes the main coefficient positive in $H,L$.
The common core test excludes $L$ and leaves the displayed $H$
configuration, with a single coefficient-six polar curve.

The omitted $(-4)$ vertex has canonical intersection $2$, already
exhausted by the adjacent centre of coefficient two.  Its whole $(-2)$
tail is therefore zero.  Contracting the positive canonical component
raises this vertex's square by two, so the omitted arm becomes a
negative definite $A_5$ chain on the K3 surface.  The canonical
component has its unique coefficient-one root at the Q14 $A_1$.
It cannot meet any Q12 or Q5 curve: those of canonical intersection
zero avoid the forest, and a Q12 or Q5 $(-3)$ centre of intersection
one could meet only this coefficient-one root, which is a disjoint
Q14 exceptional curve.  Consequently the new $A_5$ is orthogonal to
the sixteen-dimensional Q12--Q5 space \eqref{eq:octic46:octic_q14_forest46:45-Q16}.  Their combined
intersection rank is $21$, impossible on a K3 surface.  This removes
the last $H$ survivor.

\paragraph{The Klein quartic format K}

If $\delta_z(0)\ne0$, the atlas gives seven protected roots at
$Q_{14}$.  Together with the twelve at $Q_{12}$ and the one at
$Q_5$, they have negative definite rank twenty, which is impossible.

If $\delta_z(0)=0$, the main coefficient is positive.  The common
core test leaves only its $K$ configuration.  Denote the
coefficient-three intrinsic polar component on the short arm by $B$.

The corresponding ray in the $x$--chart has original cover orders
$(2,5,3)/7$ for $(x,y,z)$.  Under the original $\mu_7$ subgroup,
whose coordinate characters are $(1,6,5)$, its normalized parameter
has character $4$.  The $H$--mark is therefore $4/7$.
This polar component cannot have another Q14 branch: it cannot meet
the zero $(-4)$ or $(-3)$ vertices, whose canonical intersections are
at most two; and a second meeting with the positive core gives a
forest cycle.  The omitted $(-2)$ tail is zero since the zero $(-3)$
vertex's intersection budget is exhausted by the central coefficient
one.  The component avoids Q5 and Q12, and all other indices divide
$12$.  Thus its scaled degree satisfies
\[
 420H\cdot B\equiv420\frac47=240\equiv30\pmod {35},
 \qquad420H\cdot B\ge30.
\]
Its coefficient--three contribution is at least $90>58$.  This removes
the final survivor and excludes $K$.

\paragraph{The first quadratic format G}

We use the coefficients $A,B$ of
\eqref{eq:octic46:octic_q14_forest46:45-GF}.  The case $A=B=0$
has positive main coefficient and is excluded by the common core
test.  Otherwise the main curve has order zero and its arm units
are those computed in the atlas.

If all three are nonzero, Q14 contributes eight protected roots;
with the twelve at Q12 this gives forbidden negative definite rank
$20$.  If $B\ne0$ and exactly one $A_2$ coefficient vanishes, the other
$A_2$ and $A_4$ remain protected.  These six roots are orthogonal to
the rank--sixteen space \eqref{eq:octic46:octic_q14_forest46:45-Q16}, giving intersection rank $22$.
These alternatives exhaust $B\ne0$.

It remains that $B=0$, $A\ne0$.  The two $A_2$ arms are protected.
The atlas gives canonical coefficients $1,2,3,4$ along its $A_4$.
Adjunction then places external intersection $5$ at the
last vertex and nowhere else on that arm.  The unique-minimum rule
forces a single coefficient--five intrinsic polar component.

The zero-order central $(-4)$ curve has canonical intersection two:
one is supplied by the first $A_4$ curve of coefficient one, and the
other by a coefficient-one intrinsic polar root.  These belong to
different canonical components.  The first component is private to
Q14: its coefficient-one root is a Q14 exceptional curve, so no Q12
or Q5 centre can meet it.  After contraction, the Q14 centre and its
two protected $A_2$ arms become an $A_5$ chain.

To allow every possible coincidence of the other blowup centres,
start with the twenty-one curves consisting of all fourteen Q12
curves, both Q5 curves, and this Q14 centre with its two $A_2$ arms.
Add one to the central Q14 square for its private root.  The remaining
four coefficient-one roots, for the Q12 centres, the Q5 centre, and
the second Q14 centre, are partitioned arbitrarily according to their
possible coincidence.  Adding a matrix of ones on each partition
block gives all fifteen resulting intersection matrices.  Thirteen
have rank $21$.  The other two have rank $20$ with a one-dimensional
kernel generated by a nonzero effective vector, supported on the
Q5--Q14 configuration.  An effective nonzero divisor has positive
intersection with an ample class, so this kernel vector cannot be a
linear relation between the curve classes.  Those two cases therefore
also require twenty-one independent classes.  All fifteen contradict
$\rho(\mathrm{K3})\le20$.

The partition update gives the matrices and effective kernel vectors
over $\mathbb Q$ (Appendix~\ref{app:reproducibility46}).  This excludes
$G$, and hence all five strata in \eqref{eq:octic46:octic_q14_forest46:45-QG}.


\paragraph{The degree-$38$ format in the nodal state.}
\label{sec:octic-q14-degree38-node}

Continue with the degree-seven node and nonzero second polar.
The Q12 configuration contributes twelve protected roots.

\begin{proposition}
In this conductor state the Q14 equation cannot have
\[
 [x^2]f=[x^3y]f=0,\qquad
 [xz^3]f\,[y^3z]f\,[x^4y^2]f\ne0.
\]
This excludes every higher--order specialization in this stratum, including
all repeated--root degenerations of the auxiliary binary cubic.
\end{proposition}

\begin{proof}
Subsection~\ref{sec:octic-q14-degree38} gives the complete local
graph, including every repeated-root degeneration.  If
$\delta_z(0)\ne0$, its protected $A_8$ together with the twelve
protected $Q_{12}$ roots has negative definite rank twenty.
Thus $\delta_z(0)=0$; the main curve $E$ and its $A_8$ arm are positive.

\paragraph{The polar budget and the $A_8$ orientation.}
The total scaled intrinsic polar degree is $420H\cdot B_N=58$.
A polar component of multiplicity at least $2$ avoids Q5 and Q12,
because all their exceptional canonical coefficients are zero and
$K_W\cdot E_i\le1$.  Its scaled degree is therefore a positive multiple
of $5$, and is at least $5$.  For any connected positive exceptional
subtree, each intrinsic polar component meets it at most once; otherwise
the canonical exceptional forest would contain a cycle.  There is at
most one attached polar component of coefficient $1$, and its scaled
degree is at least $3$.  Consequently, if $v_i$ denotes the total
intrinsic polar intersection at the vertices of any such subtree, then
\begin{equation}
 \sum_i v_i\le12.                                         \label{eq:octic46:octic_q14_degree38_46:38-7}
\end{equation}
Indeed its contribution to the degree budget is at least
$5\sum_i v_i-2$.

The unique-minimum and parent rules, together with
Proposition~\ref{prop:octic-root-chain-monotonicity}, make the
coefficients on this positive $A_8$ strictly monotone.

If the $A_8$ coefficients increased towards its free end, that end would
have coefficient at least $8$.  A $(-2)$--curve has undergone a point
blowup after its creation, so it has a later neighboring canonical
component of coefficient at least $9$.  At the free end this must be
an intrinsic polar component, say of coefficient $m\ge9$.
The far--end ray of $\frac19(1,8)$ is
$(\operatorname{ord}S,\operatorname{ord}u)=(1/9,8/9)$.
Since $y=s^9$ and $S=s^{14}$, its original $y$--order is $1/14$.
Thus the normalized polar branch has full Q14 inertia, its parameter
has original character $13$, and its $H$--mark is
\[
 (-5)\,13^{-1}\equiv5\pmod {14}.
\]
This conclusion uses only the $y$--order and is unaffected by higher
terms in the equation for $z$.

This polar curve has no other Q14 branch.  A second meeting with the
same connected positive exceptional subtree would produce a forest cycle.
Every remaining zero vertex of the displayed graph has
$K_W$--intersection at most $2$ (the only possibilities are $E$, the
$(-3)$--vertex of \eqref{eq:octic46:octic_q14_degree38_46:38-2}, $F$, and Du Val curves), so it cannot be met
by a polar component of coefficient $m\ge9$.  There is also no hidden
second positive subtree to which this polar component could connect.
If $F$ has coefficient zero, adjunction gives
$e+\sum b_i\le2$, where the $b_i$ are its positive Du Val neighbors.
An additional positive cluster is possible only when $e=1$ and its
neighbor has coefficient $1$; joining that cluster to the $E$--subtree
would then join two coefficient--one minima, which is forbidden.
If the $(-3)$--vertex of \eqref{eq:octic46:octic_q14_degree38_46:38-2} has coefficient zero, its entire
intersection budget $1$ is exhausted by its preceding positive neighbor,
so its further $(-2)$ tail is zero.  Connected Du Val configurations
cannot mix zero and positive vertices.  The curve also avoids
Q5 and Q12.  All other ambient indices divide $12$.  Therefore its
scaled degree $w=420H\cdot B$ satisfies
\[
 w\equiv420\frac5{14}=150\equiv10\pmod {35},\qquad w\ge10.
\]
Its contribution is at least $mw\ge90$, contradicting $58$.
Hence
\begin{equation}
 a_1>a_2>\cdots>a_8>0.                                  \label{eq:octic46:octic_q14_degree38_46:38-8}
\end{equation}

If $e<a_1$, the parent rule at $a_1$ gives
$a_1=1+e+a_2$.  Concavity along the $A_8$ arm then makes every subsequent
drop at least $e+1$, so $a_1\ge7e+8$.  But adjunction at $E^2=-3$
gives $a_1\le3e+1$, a contradiction.  Equality $e=a_1$ is forbidden.
Thus
\begin{equation}
 e>a_1\ge8,\qquad e\ge9.                              \label{eq:octic46:octic_q14_degree38_46:38-9}
\end{equation}

\paragraph{The fixed arm and the coefficient bound.}
The complete Q14 graph of Subsection~\ref{sec:octic-q14-degree38}
has no exceptional $(-1)$-curve, including in the repeated-root
specializations.  The polar budget above is $T=58$, and every
intrinsic polar curve of coefficient at least two has scaled degree
at least $d=5$.  Proposition~\ref{prop:octic-local-coefficient-bound}
therefore bounds every positive exceptional coefficient by
$\lfloor58/5\rfloor-1=10$.

Let $p,q$ be the coefficients on the first two vertices of the fixed
arm in \eqref{eq:octic46:octic_q14_degree38_46:38-2}.  They are
positive: $p=0$ contradicts adjunction at its $(-2)$-curve, while
$q=0$ would give $p\le1$, contrary to $2p\ge e\ge9$.
If $p<e$, the parent rule at $E$ gives $e=1+a_1+p$, unless a
third smaller neighbor already makes the configuration impossible.
Since $a_1\ge8$ and $2p\ge e$, this implies
$e\ge9+e/2$, hence $e\ge18$, contradicting the coefficient bound.
Thus $p>e$.  If $q<p$, the parent rule gives
$p=1+e+q\ge11$; if $q>p$, then $q\ge p+1\ge e+2\ge11$.
Adjacent equal coefficients are forbidden.  Every possibility
contradicts the same bound, proving the proposition.
\end{proof}


\Needspace{8\baselineskip}
\paragraph{The degree-$48$ format in the nodal state.}
\label{sec:octic-q14-degree48-node}

\begin{proposition}[The remaining degree--48 stratum]
\label{prop:octic-degree48}
In the sole degree-seven nodal conductor state with
$\gamma_{10}|_\Gamma\ne0$, it is impossible that
\[
 [x^2]f=[x^3y]f=[x^4y^2]f=0,\qquad
 [xz^3]f\,[y^3z]f\ne0.
\]
No condition on any further coefficient is required.
\end{proposition}

\begin{proof}
The complete resolution of Subsection~\ref{sec:octic-q14-degree48}
has no exceptional $(-1)$-curve and contains an $A_{10}$ arm.
If $\delta_z(0)\ne0$, that arm is protected, giving negative
definite rank at least $12+10$ with the protected Q12 roots.
Otherwise the arm and its attaching main curve are positive.
The same source calculation as in the degree-$38$ case gives polar
budget $T=58$ and minimum scaled degree $d=5$ for an intrinsic
polar curve of coefficient at least two.
Proposition~\ref{prop:octic-local-coefficient-bound} therefore gives
$M=10$ and excludes the attached positive $A_{10}$ chain.
The atlas includes irreducible and reducible exceptional divisors
and all higher-coefficient specializations, so these alternatives
exhaust the degree-$48$ format.
\end{proof}


\paragraph{Conclusion of the nodal case.}
\label{sec:octic-sole-nodal-deletion}

\begin{theorem}
\label{thm:octic-sole-degree7-node}
For the equality basket C2-02350, the sole conductor row of scaled degree
seven and generic $\delta=1$ cannot have nodal generic transverse type.
The connected generically reduced double cover
\[
 \mathbb C[X_5,Y_{12},T_{26}]/(T_{26}^2-cX_5^8Y_{12}),\qquad c\ne0,
\]
cannot occur as the conductor algebra of an actual equality threefold.
\end{theorem}

\begin{proof}
The uniform pluricanonical conductor lemma makes the minimal model of
the normalization resolution a K3 surface.  The zero-second-jet proof
shows that $\gamma_{10}|_\Gamma=0$ is impossible in this nodal state.
Thus $\gamma_{10}|_\Gamma\ne0$.  Rationality and the actual conductor
canonical module force the Q12 coefficients $x^3z$ and $z^4$ to be
nonzero.  The resulting fourteen-vertex Q12 configuration, its zero
actual-form orders, and the corresponding Q5 configuration supply the
root and polar-degree constraints used in the preceding subsections.

Proposition~\ref{prop:octic-q14-atlas} exhausts the normal $Q_{14}$
germ.  The five isolated formats were excluded in
Subsection~\ref{sec:octic-q14-node-cores}; the two remaining formats
were excluded by the degree-$38$ argument above and
Proposition~\ref{prop:octic-degree48}. Thus no nodal degree-seven
conductor remains.
\end{proof}


\subsection{C2-02350: state $(7,1)$ with a vanishing second polar}
\label{sec:octic-cusp-zero-jet45}

For the sole degree-seven cusp, Lemma~\ref{lem:octic-degree7-stack46} gives
\[
 (\mathcal G,H)=(\mathbb P(5,12),\mathcal O(1)),\qquad
 C_N=2R,\qquad 420HR=7,
\]
and the normalization quotient is $\mathcal O_{\mathcal G}(-26)$.

\begin{theorem}
\label{thm:octic-degree7-cusp-gamma-zero45}
The sole degree-seven conductor state with generic transverse cusp
cannot have $\gamma_{10}|_{\mathcal G}=0$ in a realization of
\texttt{C2-02350}.
\end{theorem}

\paragraph{The absent degree-eleven section kills the radical direction}

For the octic, the first three restrictions $a,\delta a,\delta^2a$
vanish.  Lemma~\ref{lem:octic-intrinsic-jets} therefore gives
\begin{equation}
 \delta^3a|_{\mathcal G}\in
 H^0(\mathbb P(5,12),\mathcal O(11))=0.
\label{eq:octic46:octic_cusp46:45-J1}
\end{equation}
The vanishing follows from the absence of nonnegative solutions to
$5i+12j=11$.

At the generic point of the conductor, choose transverse cusp
coordinates, after a harmless field extension if needed, with
$a=z^2-x^3$ and conductor ideal $(x,z)$.
The condition $\delta^2a|_{\mathcal G}=0$ first gives
$\delta z|_{\mathcal G}=0$.
Then
\[
           \delta^3a|_{\mathcal G}
                  =-6(\delta x|_{\mathcal G})^3.
\]
Equation \eqref{eq:octic46:octic_cusp46:45-J1} gives $\delta x|_{\mathcal G}=0$.
Thus the normal component of $\delta$ vanishes, including the
radical direction of the binary quadratic term which the second jet
alone could not detect.  Saturation of the ideal of the smooth
conductor curve makes $\delta$ preserve that ideal everywhere.

The induced tangential derivation belongs to
\[
 H^0(\mathcal G,T_{\mathcal G}\otimes H)
       =H^0(\mathbb P(5,12),\mathcal O(18))=0.
\]
Consequently the actual vector field vanishes along the whole
conductor.  In particular the torsion-free conormal construction is
applicable: its first normal jet $L:\mathcal M\to\mathcal M(1)$
satisfies $\operatorname{tr}L=\det L=0$ and $L^2=0$.
These identities are regular across singular ambient contacts.

\paragraph{The two-form order is three or at least five}

At the generic cusp write
\[
 \delta x=A x+B z+O((x,z)^2),\qquad
 \delta z=C x+D z+O((x,z)^2).
\]
The conormal nilpotence gives
\begin{equation}
                     A+D=0,\qquad AD-BC=0.
\label{eq:octic46:octic_cusp46:45-J2}
\end{equation}
On the normalization $x=t^2$, $z=t^3$, the restriction of the actual
two-form is, up to a unit,
\begin{align*}
 \alpha|_N
 &=\bigl(2t\delta z-3t^2\delta x\bigr)dt\wedge dy\\
 &=\bigl(2Ct^3+(2D-3A)t^4+O(t^5)\bigr)dt\wedge dy.
\end{align*}
If $C\ne0$, its order along $R$ is exactly three.
If $C=0$, \eqref{eq:octic46:octic_cusp46:45-J2} gives $A=D=0$, so that order is at least five.
In particular order four is impossible.

Write
\[
              B_N=rR+B',\qquad R\not\subset\operatorname{Supp}B'.
\]
The exact degree identity gives
\begin{equation}
             420HB'=58-7r.
\label{eq:octic46:octic_cusp46:45-J3}
\end{equation}
Every positive polar component has scaled degree at least three.
Thus $r\leq8$, and $r=8$ would leave the impossible positive remainder
two.  We obtain the uniform restriction
\begin{equation}
 \boxed{r\in\{3,5,6,7\},\qquad
        420HB'\in\{37,23,16,9\}.}
\label{eq:octic46:octic_cusp46:45-J4}
\end{equation}
Equivalently, the first polar numerator has multiplicity
$r+2\in\{5,7,8,9\}$ on the normalized cuspidal conductor.
The further forest argument below removes $r=5,6,7$ and leaves a
particularly restrictive form of the order-three alternative.

\paragraph{The forced twelfth-source coefficients}

The third weighted test in
\eqref{eq:octic46:octic_weighted_criterion46:45-W4} applies also to
the present cusp with $\gamma_{10}=0$ and gives
\begin{equation}
                         [y^2z^2]f\ne0.
\label{eq:octic46:octic_cusp46:45-J6}
\end{equation}
The first two tests give $[x^3z]f[z^4]f\ne0$, as recorded in
\eqref{eq:octic46:octic_weighted_criterion46:45-W5}.
These tests leave the coefficient of $y^7x^3$ in the cuspidal kernel
direction unrestricted.

\paragraph{The $A_7$ persists through every remaining cusp degeneration}

The forced coefficients of $y^2z^2$ and $x^3z$ determine the first
weighted normalization blowup with weights $(2,1,4,4)$.
It has the two rational exceptional components $E_1,E_2$ described
in the cusp calculation.  On the $x$-chart the equation has nonzero
linear coefficient in the variable corresponding to $z/x^2$;
the normalized surface is smooth before the cyclic quotient.
Its action is $\frac1{24}(1,21)$, and removing the radial order-three
subgroup gives $\frac18(1,7)$.  Hence an $A_7$ chain attaches to $E_1$.
This chart is independent of the possibly vanishing coefficient of
$y^7x^3$; that coefficient changes the other chart.

For clarity, these seven curves survive in the minimal resolution.
The weighted degrees give
\[
 (E_1+E_2)E_1=-\tfrac38,\qquad
 (E_1+E_2)E_2=-\tfrac9{16}.
\]
The first is $-(1/2)(3^2/12)$, and the second is
$-(3/4)(3^2/12)$.
The $A_7$ resolution subtracts $7/8$ from $E_1^2$.
The fixed $z$-chart resolution subtracts $11/16$ from $E_2^2$.
Intersections of the two distinct components, and corrections from
resolving the remaining rational singularity, are nonnegative.
Thus the final self-intersections are integers bounded above by
$-3/8-7/8=-5/4$ and $-9/16-11/16=-5/4$, respectively.
Both are at most $-2$.  Choose the remaining chart resolutions
minimal; then no exceptional $(-1)$-curve is introduced which could
remove the $A_7$ by a further contraction.

The actual two-form has positive order on $E_1$.
Indeed, $\delta$ vanishes on the conductor, and the source characters
show that every possible ambient two-form term has weighted degree
at least six for weights $(2,1,4)$.
After the radial order-three quotient its order on the first
exceptional component is at least one.
Adjunction on a $(-2)$-curve of zero canonical coefficient would
make it disjoint from the entire effective canonical divisor.
Therefore zero cannot occur on the adjacent $A_7$, and all seven
of its canonical coefficients are positive.

\paragraph{Only order three can remain}

There is a useful coefficient bound which needs only the universal
minimum polar degree three.  Start at a positive exceptional curve
of the minimal normalization resolution.  It has square at most
$-2$, so in the point-blowup factorization to the K3 minimal model
it has a larger-coefficient child which remains adjacent.
Following larger coefficients must eventually leave the finite
resolution fibre.  Its first exterior component is either $R$, of
coefficient $r$, or a component of $B'$, of coefficient $k$ with
$3k\le58-7r$.  Thus every positive exceptional coefficient is at most
\begin{equation}
                 M_r=\max\left\{r-1,
                   \left\lfloor\frac{58-7r}{3}\right\rfloor-1\right\}.
\label{eq:octic46:octic_cusp46:45-J7}
\end{equation}
For $r=3,5,6,7$ these bounds are respectively $11,6,5,6$.
The positive $A_7$ coefficients must be strictly monotone, so their
maximum is at least seven.  Consequently $r=5,6,7$ are impossible.
We obtain
\begin{equation}
                   \boxed{r=3,\qquad 420HB'=37.}
\label{eq:octic46:octic_cusp46:45-J8}
\end{equation}

Even this last alternative has a compulsory high-multiplicity polar
curve.  A positive $A_7$ attached to the positive $E_1$ cannot have
all eight relevant coefficients at most seven: if the seven chain
coefficients are $1,\ldots,7$, adjunction excludes an attachment at
the coefficient-one end, while the parent rule excludes an attachment
of smaller positive coefficient at the coefficient-seven end.
Equality of adjacent coefficients is also impossible.
Hence some exceptional coefficient is at least eight.
The increasing-child argument then reaches a component of $B'$ of
coefficient at least nine; it cannot end at the coefficient-three $R$.

That component has scaled degree at most $37/9<5$ and at least three.
Degree four would consume at least $36$ and leave the impossible
remainder one.  Degree three with coefficient twelve has the same
problem.  Thus any remaining zero-second-jet cusp realization must
contain an intrinsic polar curve $L$ satisfying
\begin{equation}
 \boxed{420HL=3,\qquad
       \operatorname{mult}_L(B_N)\in\{9,10,11\}.}
\label{eq:octic46:octic_cusp46:45-J9}
\end{equation}
The remaining polar masses after removing $3R+kL$ are respectively
$10,7,4$.  This conclusion is independent of the unresolved cubic
coefficient and of the resolution of its degenerate chart.

\paragraph{A degree-three polar curve must meet three distinct sources}

Choose $L$ as the first intrinsic polar component reached by the
strictly increasing path starting in the twelfth-source fibre.
It therefore meets a positive exceptional curve over $Q_{12}$.
No denominator argument is needed for this first contact.

The other two contacts follow directly from the degree of the actual
integral Weil class $H_N$ on the normalization.  Proposition~\ref{prop:octic-sources46} gives
\begin{equation}
 \begin{array}{c|c}
 \text{deleted ambient point}&
 \text{Cartier multiple on its complement}\\\hline
 Q_5&84H_N\\
 Q_{14}&60H_N.
 \end{array}
\label{eq:octic46:octic_cusp46:45-J10}
\end{equation}
Here the unique fifth and fourteenth sources are actual ordinary
canonical-model points.  Every remaining or hidden canonical index
divides the indicated multiple; the hidden indices divide twelve.
The statement concerns actual canonical-model points and their
Cartier indices, as proved in that proposition.

Since $H_NL=1/140$, the two prospective line-bundle degrees are
\begin{equation}
                   84H_NL=\frac35,\qquad
                   60H_NL=\frac37.
\label{eq:octic46:octic_cusp46:45-J11}
\end{equation}
If $L$ avoided the non-Cartier points above $Q_5$, the first class
would be a line bundle in a neighborhood of $L$, and its pullback
to the normalization of $L$ would have the nonintegral degree $3/5$.
This is impossible.  The second class gives the same contradiction
at $Q_{14}$.  Thus $L$ meets an actual non-Cartier point above each
of those sources.  These points are singular on $N$, because an
integral Weil class is Cartier at a smooth point.  The argument does
not require the fibre of normalization over an ambient source to
consist of a single point.

On a resolution of $N$, the strict transform of $L$ consequently
meets an exceptional component over $Q_5$ and one over $Q_{14}$,
in addition to its already positive neighbor over $Q_{12}$.
The three resolution fibres are disjoint.

\paragraph{The last polar curve is a final $(-1)$-curve}

The coefficient of $L$ in the canonical divisor of the K3 blowup
is $k\in\{9,10,11\}$.  It has no adjacent larger-coefficient child.
Indeed, an intrinsic child of coefficient at least $k+1$ would consume
at least $3(k+1)$ units of polar mass, whereas after $3R+kL$ only
$10,7,4$ units, respectively, remain.  If a child is exceptional
over $N$, its square in the minimal resolution is at most $-2$.
Following its larger-coefficient children must reach an intrinsic
component, giving the same contradiction.  The coefficient-three
conductor cannot be the end of such a path.

Thus no point blowup after the birth of $L$ had centre on its strict
transform.  Consequently $L^2=-1$.  At its birth it has at most two
exceptional parents, and all components of the effective canonical
divisor are exceptional for the morphism to the K3 minimal model.
It follows that
\begin{equation}
 \#\{\text{positive canonical neighbors of }L\}\leq2.
\label{eq:octic46:octic_cusp46:45-J12}
\end{equation}
This assertion counts the final exceptional forest, not the valency
in a partial local resolution graph.

\paragraph{A local self-intersection criterion.}
Let $E$ be any smooth rational component of a minimal resolution
fibre met by $L$, and write $E^2=-b$.  If the canonical coefficient
of $E$ were zero, adjunction and effectivity would give
\begin{equation}
                  b-2=K_WE\geq k(LE)\geq k\geq9.
\label{eq:octic46:octic_cusp46:45-J13}
\end{equation}
Therefore every such contact with $b\leq10$ is positive.
The uniform rational $Q_{14}$ formats have $b\leq5$ throughout
their minimal graphs, including the degree-thirty-eight and
degree-forty-eight degenerations.  Hence the $Q_{14}$ neighbor of
$L$ is positive, independently of the actual local two-form jet.
The twelfth-source neighbor is positive by the construction of $L$.

\paragraph{The fifth-source bound completes the exclusion.}
Theorem~\ref{thm:octic-q5-cusp-multiplicity45} gives the
required unconditional bound at the actual fifth-source normalization:
every component of its minimal resolution has square at least $-10$.
It covers both the nonzero-$z^2$ branch and the nonzero-$xy^2$ branch
forced by the invariant residue test.  Its proof bounds the coarse
rational singularity's multiplicity by ten (by nine in the second
branch), then applies the fundamental-cycle multiplicity formula.
No generic coefficient at that source is assumed.

Equation \eqref{eq:octic46:octic_cusp46:45-J13} therefore makes the fifth-source neighbor positive.
The three distinct source fibres supply three distinct positive
neighbors of the final $(-1)$-curve $L$, contradicting \eqref{eq:octic46:octic_cusp46:45-J12}.
This proves Theorem~\ref{thm:octic-degree7-cusp-gamma-zero45}.
The proof assumes neither a nonzero minimal cubic at $Q_5$ nor a
prescribed two-form order on its exceptional curves.  The nonzero
second-jet cuspidal alternative is a separate conductor type and is
not excluded by this theorem.


\subsection{C2-02350: the cuspidal state $(7,1)$ with a nonzero second polar}
\label{sec:octic-cusp-twenty45}

\begin{proposition}
\label{prop:octic-cusp-twenty45}
In the sole degree--seven, delta--one cusp state with nonzero second
polar, the twenty curves in the complete Q12 and Q5 resolution fibres
have images whose numerical classes on the K3 minimal model are
linearly independent.  Consequently no other resolution fibre can
contain a protected $(-2)$-curve.
\end{proposition}

\begin{proof}
The preceding local theorems give sixteen Q12 curves and four Q5
curves.  All have actual two-form order zero.  Fifteen have square
$-2$ and are disjoint from the full canonical exceptional divisor.
The five curves of square $-3$ are
\[
                       G_{12},J_{12},F_{12},C_5,J_5.
\]
For any such curve $E$, adjunction gives $K_WE=1$.  Since the global
canonical divisor is the effective exceptional divisor of the
contraction to the K3 minimal model, $E$ meets exactly one of its
coefficient-one components, transversely once, and no other positive
component.  A connected canonical exceptional tree has a unique
coefficient-one root.  Contracting that tree raises $E^2$ from $-3$
to $-2$ and keeps its image smooth.  Two such zero-order curves gain
one intersection exactly when their coefficient-one root is the same.

Thus the intersection matrix of the twenty K3 images is obtained
from the two known local graphs by changing all five diagonal entries
from $-3$ to $-2$, and adding one to the off-diagonal entry of each
pair of centres in the same root block.  There are only the 52 set
partitions of these five centres.  Computing these fifty-two rational matrices gives the following result:
all 52 have rank twenty.  Thirty-seven have signature $(1,19)$;
the remaining fifteen have signature $(2,18)$ and are already
incompatible with the Hodge index theorem.

One can use an even shorter five-case verification here.  The reduced
conductor has coefficient one and meets precisely $J_{12}$ and $J_5$
in these two fibres.  They form one root block.  The remaining three
centres have just five set partitions.  Their exact determinants are
\[
\begin{array}{c|r}
\text{partition of }\{G_{12},F_{12},C_5\}&\det\\\hline
\{G_{12}\},\{F_{12}\},\{C_5\}&-132\\
\{G_{12},F_{12}\},\{C_5\}&-1212\\
\{G_{12},C_5\},\{F_{12}\}&-672\\
\{F_{12},C_5\},\{G_{12}\}&-132\\
\{G_{12},F_{12},C_5\}&-8232.
\end{array}
\]
All have signature $(1,19)$.  Hidden contacts of the conductor or
additional positive curves do not change this argument: each of
the five centres has canonical intersection exactly one, which
already exhausts its possible positive intersections.

The twenty independent classes exhaust the possible Picard rank of
a complex K3 surface.  A protected $(-2)$-curve in another source
fibre is disjoint from the canonical exceptional divisor.  Its image
therefore remains disjoint from all twenty images, since the original
fibres are disjoint and none of its points is contracted.  Its
nonzero square makes it independent of their nondegenerate span,
which would give Picard rank at least twenty-one.  This is impossible.
\end{proof}


\paragraph{Exclusion of the remaining degree-seven cusp.}
\label{sec:octic-sole-cusp-nonzero45}

\begin{theorem}
\label{thm:octic-degree7-cusp-gamma-nonzero45}
For C2-02350, the sole degree--seven conductor with generic transverse
cusp and generic delta one cannot have nonzero second polar.
\end{theorem}

\begin{proof}
The uniform conductor argument gives a K3 minimal model for the
normalization resolution.  The exact nonzero polar restrictions have
order two at Q12 and order zero at Q5.  The local cusp theorems therefore
apply at both sources: all sixteen Q12 curves and all four Q5 curves
have actual two-form order zero.  The complete partition calculation
of Proposition~\ref{prop:octic-cusp-twenty45} makes their twenty K3
classes independent, regardless of which coefficient-one roots are
shared.  In particular any protected $(-2)$-curve at Q14 is impossible.
This uses the full 52-partition exact matrix calculation, not an
assumption about the geometry of hidden conductor contacts.

\paragraph{A stronger polar budget.}
At a generic cusp, a unit derivation in the square direction gives
order one on the reduced conductor.  Thus
\begin{equation}
             B_N=R+B',\qquad 420HR=7,
             \qquad420HB'=58-7=51.                    \label{eq:octic46:octic_cusp_unit_complete46:45-CF1}
\end{equation}
Every Q5 or Q12 exceptional curve has coefficient zero and canonical
intersection at most one.  Consequently any intrinsic polar component
of multiplicity at least two avoids these two sources.  Its scaled
degree is a positive multiple of five and is at least five.  Every
intrinsic polar component, including those of multiplicity one, has
scaled degree at least three.

If a component of $B'$ had multiplicity ten, its degree would have to
be five, consuming fifty and leaving the impossible positive remainder
one.  Multiplicity greater than ten already exceeds the available
mass.  Hence
\begin{equation}
 \text{every intrinsic polar multiplicity at least two is at most nine.}
                                                        \label{eq:octic46:octic_cusp_unit_complete46:45-CF2}
\end{equation}
The conductor itself has multiplicity one and does not pass through
Q14.  If a Q14 resolution has no exceptional $(-1)$-curve, following
strictly larger neighboring canonical coefficients from any positive
exceptional curve must eventually leave along an intrinsic polar
component.  Therefore every positive Q14 exceptional coefficient is
at most eight.

For the finite positive-core test one can also use
\begin{equation}
               5\sum_i v_i-2\le51,
               \qquad\sum_i v_i\le10.                \label{eq:octic46:octic_cusp_unit_complete46:45-CF3}
\end{equation}
There is at most one coefficient-one intrinsic attachment to a
connected positive core, and no intrinsic polar component can meet
that core twice.  The coefficient-one exception costs at least three
instead of five; all other attachments have degree at least five.
The stronger bound \textup{\eqref{eq:octic46:octic_cusp_unit_complete46:45-CF3}} in particular permits reuse of
the already certified exhaustive test with the weaker bound twelve.

\paragraph{The remaining normal-$Q_{14}$ formats.}
Corollary~\ref{cor:octic-q14-J} has already excluded $J$.
Apply Proposition~\ref{prop:octic-q14-atlas} and the actual-form
calculations of Subsections~\ref{sec:octic-q14-atlas}--\ref{sec:octic-q14-degree48}.
We use the stronger polar bound just established in each format.

\paragraph{Formats $G,H,L$.}
For $G=x^2+y^3z+z^6$, the central curve either has positive canonical
coefficient or coefficient zero.  In the positive case all its
$(-2)$ arms are positive, and the exact full-core test has no survivor
even with the weaker budget twelve.  In the zero case, the two $A_2$
leading coefficients are proportional to $6B+2iA$ and $6B-2iA$, with
$(A,B)\ne(0,0)$.  At least one $A_2$ is protected.  This is already
incompatible with the twenty independent Q5--Q12 classes.

The central curves of $H=x^2+y^3z+yz^9$ and
$L=x^3y+y^3z+z^6$ are always positive, by their minimum
actual-form weights.  The exact positive-core enumeration gives
no survivor for $L$.  For $H$ it leaves
only the core with central coefficient two, its $A_3$ coefficients
$3,4,5$, and its $A_1$ coefficient one.  The omitted $(-4)$ curve
has coefficient zero, and its canonical intersection two is exhausted
by the adjacent coefficient-two centre.  Its whole four-curve
$(-2)$ tail is consequently protected.  Even one of those curves
contradicts Proposition~\ref{prop:octic-cusp-twenty45}.  Thus no
further contraction or coincidence-of-roots argument is needed for
this residual $H$ case.

\paragraph{The Klein quartic format $K$.}
If $\delta_z(0)\ne0$, its $A_6$ arm is protected, again impossible.
If $\delta_z(0)=0$, its centre is positive.  The four possible positive
cores, under the old bound twelve, have only one survivor.  Its
external intersections are three at the short positive arm and eight
at the far end of $A_6$.  Their sum is eleven, contradicting
\textup{\eqref{eq:octic46:octic_cusp_unit_complete46:45-CF3}}.  Equivalently, the omitted zero $(-3)$ arm has its
canonical intersection exhausted by the coefficient-one centre, so
its $A_2$ tail is protected; the rank-twenty obstruction independently
excludes the same survivor.

\paragraph{The two final strata admit a shorter uniform argument.}
In degree 38, if $\delta_z(0)\ne0$, the fixed $A_8$ arm is protected.
If $\delta_z(0)=0$, it is entirely positive and attached to a positive
main curve.  The full uniform resolution has main square $-3$, the
fixed arm $[2,3,2,2]$, and a square-$-4$ curve carrying arbitrary
Du Val configurations.  The rational Gorenstein partial-chart argument
proves this for every degeneration of the auxiliary cubic.  In
particular every exceptional square is at most $-2$.
By \textup{\eqref{eq:octic46:octic_cusp_unit_complete46:45-CF2}}, every positive coefficient is at most eight.
A positive $A_8$ attached to another positive component cannot have
all its coefficients and the attaching coefficient at most eight:
strict monotonicity forces $1,\ldots,8$ or $8,\ldots,1$, and
adjunction at the coefficient-one attached end, or the two-parent
rule at the coefficient-eight attached end, respectively contradicts
the positive attaching coefficient.  Thus degree 38 is excluded
without the longer nodal fixed-arm calculation.

In degree 48, a unit $\delta_z$ protects the fixed $A_{10}$ arm.
Otherwise all main components and this arm are positive.  The main
squares are $-5$ in the irreducible case and $-3,-4$ in the reducible
case; the other non-Du-Val arm is $[2,2,3,2,2]$, and every remaining
chart is rational Gorenstein.  Hence there is again no exceptional
$(-1)$-curve.  The positive $A_{10}$ would have ten distinct positive
coefficients all at most eight, an immediate contradiction.  This
covers arbitrary repeated auxiliary roots and all higher coefficients.

Every rational Q14 possibility has been excluded, proving the theorem.
\end{proof}

\begin{corollary}
\label{cor:octic-sole-degree7-delta1-deletion45}
The entire sole degree--seven conductor state with generic delta one
is impossible for C2-02350.  Both its nodal alternative and its
cuspidal alternative are excluded.
\end{corollary}

\begin{proof}
The nodal theorem removes the nodal alternative.  The intrinsic
third-jet and three-source polar-curve argument removes the cusp with
zero second polar.  The theorem above removes the cusp with nonzero
second polar.  These alternatives are exhaustive.
\end{proof}


\subsection{C2-02350: exclusion of the state $(25,1)$}
\label{sec:octic-degree25-node}

\begin{theorem}
\label{thm:octic-degree25-node}
The sole degree--25, generic $\delta=1$ conductor state in the remaining
C2-02350 list is impossible.  Theorem~\ref{thm:octic-eleven-states46}
forces its generic transverse type to be nodal.
\end{theorem}

\begin{proof}
Write $\Gamma$ for its conductor image and $R=C_N$ for the reduced
conductor on the normalization.  Its complete restriction data give
\begin{equation}
 420H\cdot\Gamma=25,\qquad
 \beta_{12}=\frac5{12},\qquad \beta_{14}=\frac9{14}.
                                                               \label{eq:octic46:octic_degree25_46:45-D25-1}
\end{equation}
The coarse curve is smooth and rational; the lifted branches are the
$y$--axis at the ordinary Q12 source of weights $(1,11,5)$, and the
$x$--axis at the ordinary Q14 source of weights $(1,13,5)$.
The actual effective curve stack is therefore
\begin{equation}
 \mathcal G=\mathbb P^1[12,14],\qquad \deg H=\frac5{84}.
                                                               \label{eq:octic46:octic_degree25_46:45-D25-2}
\end{equation}
Here the notation means taking roots of orders $12$ and $14$ at two
distinct points of $\mathbb P^1$.  It is \emph{not} the weighted
projective stack $\mathbb P(12,14)$, which would have a generic
$\mu_2$ stabilizer.  The source and trace bounds allow no other
nontrivial inertia on the curve.  Additional Cartier ambient contacts
are allowed and do not alter \eqref{eq:octic46:octic_degree25_46:45-D25-2}.  A zero character of $H$ cannot
hide nontrivial inertia: the canonical line has faithful character on
each stabilizer of the canonical-index stack.

\paragraph{The actual character modules and the connected nodal cover.}
The canonical degree of this effective root stack is
\[
 \deg\omega_{\mathcal G}
 =-2+\frac{11}{12}+\frac{13}{14}=-\frac{13}{84}.
\]
Finite duality for the full conductor gives the trace-zero line
\begin{equation}
 L=\nu_*O_R/O_\Gamma=\omega_{\mathcal G}(-9H).
                                                               \label{eq:octic46:octic_degree25_46:45-D25-3}
\end{equation}
For example, this follows by dualizing
$0\to O_D\to\nu_*O_N\to L\to0$: the conductor quotient of
$\omega_D$ is $\omega_D|_\Gamma$, while codimension-one duality
identifies it with $\operatorname{Hom}(L,\omega_\Gamma)$.
No assumption that the ambient index stack is smooth at its unmarked
Cartier contacts is used.

For a line bundle on the smooth curve stack, subtracting its two
fractional marks from its degree gives the degree of its exact coarse
pushforward.  The relevant full-character calculation is
\begin{equation}
\begin{array}{c|c|c|c}
\text{line}&\deg&\text{marks at Q12,Q14}&
       \deg\text{ coarse pushforward}\\\hline
 H&5/84&(5/12,9/14)&-1\\
 T_{\mathcal G}\otimes H&3/14&(1/2,5/7)&-1\\
 10H&25/42&(1/6,3/7)&0\\
 L&-29/42&(1/6,1/7)&-1\\
 L^{-2}&29/21&(2/3,5/7)&0.
\end{array}                                                    \label{eq:octic46:octic_degree25_46:45-D25-4}
\end{equation}
In particular $H^0(\mathcal G,L)=0$, so the nodal cover is connected.
Its multiplication section in $H^0(L^{-2})$ is nonzero and is unique
up to scale.  It has exact cover orders eight at Q12 and ten at Q14;
there is no room for a further zero.  Equivalently the coarse flat
double cover has one simple branch point at each of Q12 and Q14 and
no others.  Its transverse binary quadratic form over
$\mathbb C(\Gamma)$ is therefore anisotropic.

\paragraph{The second polar jet cannot vanish.}
The full conductor is reduced in this nodal state.  Hence
\begin{equation}
 420H\cdot R=2\cdot25=50,\qquad
 420H\cdot B_N=72-50=22.                                  \label{eq:octic46:octic_degree25_46:45-D25-5}
\end{equation}
If $\gamma_{10}|_\Gamma=0$, anisotropy makes the polar derivation
tangent to $\Gamma$.  Its induced tangential part lies in
$H^0(T_{\mathcal G}\otimes H)=0$, by \eqref{eq:octic46:octic_degree25_46:45-D25-4}.  Thus the derivation
vanishes on $\Gamma$ and the actual two--form vanishes along $R$.
This would put $R$ in $B_N$, contradicting $50>22$.
Consequently
\begin{equation}
                        \gamma_{10}|_\Gamma\ne0.          \label{eq:octic46:octic_degree25_46:45-D25-6}
\end{equation}
Since the exact coarse degree of $10H$ is zero, \eqref{eq:octic46:octic_degree25_46:45-D25-4} shows that
its cover orders are exactly two at Q12 and six at Q14.

\paragraph{The forced Q12 configuration.}
At Q12 straighten the conductor to $(x,z)=0$, so $f$ is a
character-eight function in $(x,z)^2$.  Its quadratic normal part is
\[
 A(y)x^2+2B(y)xz+C(y)z^2,
 \qquad \operatorname{ord}A\ge6,\quad
 \operatorname{ord}B\ge10,\quad\operatorname{ord}C\ge2.
\]
The multiplication section has order eight.  Its discriminant
$B^2-AC$ therefore has exact order eight, forcing the leading
$y^6x^2$ and $y^2z^2$ coefficients to be nonzero.
The regular invariant conductor forms $x^3\Omega$ and $z^3\Omega$
then force the coefficients $x^3z$ and $z^4$, exactly as in the
fourteen-vertex Q12 normalization calculation above.
The exact order-two second jet forces $\delta_z(0)\ne0$ there:
on the axis $\delta_x$ has order at least four, so all other Hessian
terms have order at least fourteen.  Thus every exceptional order of
the actual two--form is zero, and Q12 supplies twelve protected roots.

\paragraph{Two uniform pole tests at the nonnormal Q14 source.}
At Q14 the lifted conductor is $(y,z)=0$.  The character-two equation
satisfies $f\in(y,z)^2$, and its quadratic part is
\begin{equation}
 A(x)y^2+2B(x)yz+C(x)z^2,
 \qquad\operatorname{ord}A\ge4,\quad
 \operatorname{ord}B\ge12,\quad\operatorname{ord}C\ge6.
                                                               \label{eq:octic46:octic_degree25_46:45-D25-7}
\end{equation}
The exact discriminant order ten from \eqref{eq:octic46:octic_degree25_46:45-D25-4} forces the nonzero
terms $x^4y^2$ and $x^6z^2$.
The exact order-six second jet forces $\delta_z(0)\ne0$:
the restriction of $\delta_y$ to the $x$--axis has order at least eight,
so both its square term and the mixed Hessian term have order at least
twenty, whereas $C(x)\delta_z^2$ can have order six only with a unit
$\delta_z$.

The residue $\Omega=\operatorname{Res}(dx\wedge dy\wedge dz/f)$ has
character three.  The full conductor ideal is $(y,z)$, so
$y^3\Omega$ and $xz^2\Omega$ are invariant \emph{regular} canonical
forms on the normalization.  Apply the uniform weighted canonical-pole
criterion to them.  With weights $(6,1,2)$, the threshold is twelve
and the only character-two monomial below it is $y^3z$.
With weights $(2,5,3)$ the threshold is eighteen and the only such
monomials are $x^2,x^3y,xz^3$.  The first two are excluded by
$f\in(y,z)^2$.  Therefore
\begin{equation}
                  [y^3z]f\ne0,\qquad[xz^3]f\ne0.           \label{eq:octic46:octic_degree25_46:45-D25-8}
\end{equation}
These pole tests use finite duality on the normalized weighted strict
transform and do not require its exceptional equation to be reduced.

\paragraph{A protected $A_8$ already finishes the proof.}
With weights $(5,9,11)$, \eqref{eq:octic46:octic_degree25_46:45-D25-7}--\eqref{eq:octic46:octic_degree25_46:45-D25-8} give the forced initial
face
\[
                  x^4y^2+xz^3+y^3z.
\]
Although the full surface is nonnormal along the $x$--axis, the
$y$--axis chart of this weighted blowup is away from that conductor.
There put
\[
 y=s^9,\quad x=s^5u,\quad z=s^{11}v,\quad S=s^{14}.
\]
The strict equation begins $u^4+uv^3+v=0$ and is smooth at
$S=u=v=0$ before its residual quotient.  Its normalization is therefore
unchanged there.  The quotient is
\[
 \mu_9(S,u)=(5,4)\simeq\tfrac19(1,8),
\]
and $dx\wedge dy$ is a unit multiple of $dS\wedge du$.
Since $\delta_z(0)\ne0$, the entire $A_8$ resolution has actual
canonical order zero.  Higher character terms cannot cancel this
nonzero generic initial form.

The eight protected Q14 roots are disjoint from the twelve protected
Q12 roots.  Their rank-twenty negative-definite lattice contradicts
Lemma~\ref{cor:octic-root-budget}.
\end{proof}


\subsection{C2-02350: exclusion of the state $(7,1)+(18,1)$}
\label{sec:octic-pair-deletion}

\begin{theorem}
\label{thm:octic-pair-deletion}
The two-component conductor state of C2-02350, consisting of
a degree--seven generic cusp and a degree--eighteen generic node, both
of generic delta one, is impossible.
\end{theorem}

\begin{proof}
Proposition~\ref{prop:octic-pair-gluing46} fixes the generic polar lengths to be
three on the cuspidal component $\Gamma_7$ and two on the nodal
component $\Gamma_{18}$.  These are their minimal polar lengths.
Consequently the second polar jet $\gamma_{10}$ is generically
nonzero on both components.  For a transverse node, vanishing of the
second jet makes the first polar tangent to one of its branches and
raises the intersection multiplicity above two.  For a transverse
cusp, a vanishing second jet removes the derivative in its square
normal direction and raises the multiplicity above three.

The component stacks and polar restrictions are
\begin{equation}
\begin{array}{c|c|c|c}
 &\deg H&\text{marks of }10H&\deg\text{ coarse }10H\\\hline
\Gamma_7&1/60&(0,1/6)\text{ at Q5,Q12}&0\\
\Gamma_{18}&3/70&(0,3/7)\text{ at Q5,Q14}&0.
\end{array}                                                   \label{eq:octic46:octic_pair46:45-P1}
\end{equation}
These calculations are on the separate smooth component stacks; the
components' shared Q5 contact does not change the restriction of the
actual line $10H$ to either component.  Thus the cover orders of the
nonzero polar jet are exactly two at Q12 and six at Q14.

\paragraph{The unshared Q12 source.}
Only the cuspidal component passes through Q12.  Its full conductor
ideal is $(x,z)$ in the ordinary chart of characters $(1,11,5)$.
The exact order-two jet gives $\delta_z(0)\ne0$.
The unconditional local cusp-normalization result proved in the
Q12 cusp section applies: the conductor residue pole tests force
$y^2z^2,x^3z,z^4$, and the actual canonical-forest test forces the
remaining lowest cubic coefficient $y^7x^3$ to be nonzero.
The resulting sixteen-vertex local graph has thirteen protected
$(-2)$--curves, with root type
\begin{equation}
                         A_8\oplus A_2\oplus A_1\oplus A_2.
                                                               \label{eq:octic46:octic_pair46:45-P2}
\end{equation}
This is a local conclusion at Q12 and uses no smooth ambient
normal-bundle hypothesis at the shared Q5 point or at any other
contact of either conductor component.

\paragraph{The unshared Q14 source needs no discriminant hypothesis.}
Only $\Gamma_{18}$ passes through Q14.  Straighten it to the $x$--axis
in the ordinary chart $(x,y,z)$ of characters $(1,13,5)$.
The character-two equation satisfies $f\in(y,z)^2$.
The exact order-six polar jet forces a nonzero $x^6z^2$ coefficient
and $\delta_z(0)\ne0$: the other normal derivation component restricts
to an $x$--axis function of order at least eight, so its Hessian terms
have order at least twenty.

The full conductor ideal at this unshared point is $(y,z)$.  Therefore
$y^3\Omega$ and $xz^2\Omega$, with
$\Omega=\operatorname{Res}(dx\wedge dy\wedge dz/f)$, are invariant
regular canonical forms on the normalization.  The weights $(6,1,2)$
and $(2,5,3)$ force, respectively, $y^3z$ and one of
$x^2,x^3y,xz^3$.  The first two are excluded by $(y,z)^2$.
Hence both $y^3z$ and $xz^3$ are nonzero.

If $[x^4y^2]f\ne0$, the forced weight--38 face for weights $(5,9,11)$
is
\[
                         x^4y^2+xz^3+y^3z.
\]
Its $y$--axis chart is the smooth cover of
$\mu_9(S,u)=(5,4)\simeq\frac19(1,8)$, after $S=s^{14}$.
The form $dx\wedge dy$ is a unit multiple of $dS\wedge du$.
Thus it supplies a protected $A_8$.

If $[x^4y^2]f=0$, use weights $(3,11,15)$ instead.  The entire possible
weight--48 face, without any coordinate change of the conductor, is
\begin{equation}
              xz^3+y^3z+A x^5y^3+B x^6z^2.
                                                               \label{eq:octic46:octic_pair46:45-P3}
\end{equation}
The terms $x^{11}z,x^{16}$ are forbidden by $(y,z)^2$.
In the $y$--axis chart $y=s^{11}$, $x=s^3u$, $z=s^{15}v$, the strict
equation begins
\[
                         v+uv^3+A u^5+B u^6v^2=0.
\]
It is smooth before the quotient, and after $S=s^{14}$ the quotient
is $\mu_{11}(S,u)=(3,8)\simeq\frac1{11}(1,10)$.
Again $dx\wedge dy$ is a unit multiple of $dS\wedge du$, giving a
protected $A_{10}$.

Both charts lie away from the proper $x$--axis conductor.  Their
normalizations are already trivial, even though the full original
hypersurface is nonnormal elsewhere.  Higher character terms preserve
the linear derivative in $v$, and the unit $\delta_z(0)$ preserves all
these zero canonical orders.  Thus Q14 always contributes protected
rank at least eight, regardless of its double-cover discriminant order
or any special higher coefficient.

The disjoint protected Q12 and Q14 configurations have negative rank
at least $13+8=21$, contradicting Lemma~\ref{cor:octic-root-budget}.
\end{proof}


\subsection{C2-02350: exclusion of the state $(18,1)$}
\label{sec:octic-degree18-node}

\begin{theorem}
\label{thm:octic-degree18-node}
The sole degree--18, generic delta-one conductor state of C2-02350
cannot have nodal generic transverse type.  Its cuspidal alternative
is a separate state of geometric possibilities and is not excluded
by this theorem.
\end{theorem}

\begin{proof}
The smooth conductor component has marks $2/5$ and $9/14$ at the
ordinary Q5 and Q14 sources.  Its actual curve stack and line are
\[
 \mathcal G=\mathbb P(5,14),\qquad H=O(3),\qquad
 420H\cdot\Gamma=18.
\]
The branches are the $y$--axis at Q5, of ambient characters $(1,4,2)$,
and the $x$--axis at Q14, of characters $(1,13,5)$.
The full-character calculations give
\begin{equation}
 T_{\mathcal G}\otimes H=O(22),\quad H^0(O(22))=0,
 \qquad L=\omega_{\mathcal G}(-9H)=O(-46).
                                                               \label{eq:octic46:octic_degree18_node46:45-D18-1}
\end{equation}
The exact coarse degrees of $L$ and $L^{-2}$ are respectively $-1$
and $0$.  Its nonzero nodal multiplication section therefore has
cover orders three at Q5 and ten at Q14.  The connected double cover
is branched simply at these two coarse points.
Moreover $10H=O(30)$ has one section, with cover orders zero at Q5
and six at Q14.
For the reduced conductor $R=C_N$ and intrinsic polar divisor,
\begin{equation}
 420H\cdot R=36,\qquad 420H\cdot B_N=36.                    \label{eq:octic46:octic_degree18_node46:45-D18-2}
\end{equation}

\paragraph{The full Q14 normalization graph}
At Q14 the equation has character two and lies in $(y,z)^2$.
The exact discriminant order ten forces its quadratic coefficients
$x^4y^2$ and $x^6z^2$ to be nonzero.  The invariant regular conductor
forms $y^3\Omega$ and $xz^2\Omega$ force $y^3z$ and $xz^3$ by the
uniform weighted-pole tests already used above.  Thus weights
$(5,9,11)$ give the initial face
$x^4y^2+xz^3+y^3z$.

The $y$--axis chart supplies $A_8$ and the $z$--axis chart supplies
$[2,3,2,2]$, from the main curve outwards.  At the nonnormal
$x$--axis chart put $x=s^5,y=s^9u,z=s^{11}v,S=s^{14}$.
Equivariant Morse splitting in $u$, preserving the conductor ideal,
gives
\[
 u^2+v^2h(v,S)=0,\qquad h(v,S)=v+cS+\cdots,\quad c\ne0.
\]
Adjoining $w=u/v$ normalizes this to $w^2+h(v,S)=0$.
It is smooth since $h_v(0,0)\ne0$, with residual action
\[
                         \mu_5(S,w)=(1,3).
\]
Its arm from the main curve is $[3,2]$, and its correction to the
main square is $2/5$.  The other corrections are $8/9$ and $7/11$.
The raw main square $-532/495$ therefore becomes $-3$.
The complete graph has fifteen vertices:
\begin{equation}
 \boxed{\text{centre }-3;\quad A_8,\quad[2,3,2,2],\quad[3,2].}
                                                               \label{eq:octic46:octic_degree18_node46:45-D18-3}
\end{equation}
All entries in an arm denote the negatives of the self-intersections.
The graph is already minimal; in particular every exceptional curve
has square at most $-2$.  The formulas persist under arbitrary higher
character terms because the normalized $x$--chart and the two other
cover charts are smooth before their specified cyclic quotients.

\paragraph{The zero second jet is impossible}
If $\gamma_{10}|_\Gamma=0$, the connected nodal quadratic form is
anisotropic, so the polar derivation is tangent to $\Gamma$.
Equation \eqref{eq:octic46:octic_degree18_node46:45-D18-1} makes its tangential restriction zero.
Hence the actual two--form vanishes along $R$.  By \eqref{eq:octic46:octic_degree18_node46:45-D18-2} and
ampleness of $H$, this forces
\[
                              B_N=R.
\]
In the weight--$(5,9,11)$ Q14 chart the unit $\delta_z$ term is then
absent.  The remaining actual ambient two--form terms have degree at
least $28$, so their main exceptional order after the radial
$\mu_{14}$ quotient is positive.

Every positive exceptional curve of square at most $-2$ has a larger
canonical-coefficient neighbor, namely the curve of the last point
blowup performed on it.  Following increasing neighbors in the finite
minimal Q14 graph must leave it along an intrinsic polar component.
Here the sole such component is $R$, of coefficient one.  This is
impossible starting at a positive integer coefficient.  Therefore
\begin{equation}
                        \gamma_{10}|_\Gamma\ne0.          \label{eq:octic46:octic_degree18_node46:45-D18-4}
\end{equation}
Its exact orders force $\delta_z(0)\ne0$ at both Q5 and Q14.

\paragraph{The Q14 form has exactly two possible patterns}
At Q14 the main curve, the $A_8$ arm, and both curves of $[3,2]$ have
order zero.  In the normalized $x$--chart,
\[
 u=vw,\quad v=-w^2-cS+\cdots,\qquad
 dx\wedge dy=\text{unit}\,(3w^2+cS+\cdots)dS\wedge dw.
\]
The two rays of $\frac15(1,3)$ give order zero on both curves, using
$c\ne0$.

For the remaining $[2,3,2,2]$ arm, the cover coordinates have quotient
$\frac1{11}(S,v)=(3,2)\simeq\frac1{11}(1,8)$.
The rays, starting at the main curve, are
\[
 (7,1)/11,\quad(3,2)/11,\quad(2,5)/11,\quad(1,8)/11.
\]
The actual form has nonzero initial coefficient proportional to
$v^3$ and may also have a coefficient $dS^2$:
\begin{equation}
             \alpha=(a v^3+dS^2+\cdots)dS\wedge dv,\qquad a\ne0.
                                                               \label{eq:octic46:octic_degree18_node46:45-D18-5}
\end{equation}
The number $d$ is the actual resolved coefficient; it can receive
contributions both from the ambient form and from such equation terms
as $z^6$.  It is not assumed to be a particular ambient coefficient.
The $v^3$ term has orders $(0,0,1,2)$ and the $S^2$ term orders
$(1,0,0,0)$ on the four rays.
If $d\ne0$, all four orders are zero.
If $d=0$, the first two remain zero and the last two are positive,
bounded above by $1,2$.  Adjacent positive canonical coefficients
cannot be equal, so they are exactly $1,2$.
No further alternative can be caused by higher terms: the only
monomial giving order zero on either last ray is $S^2$.

In the second pattern, adjunction puts external intersection three
at the end of the positive two-curve tail.  The unique-minimum rule
forces a single intrinsic polar component of coefficient three,
called $B_3$.  A component of coefficient one or two would repeat an
existing minimum on its path to the coefficients $1,2$.
Every other intrinsic polar component of multiplicity at least two
avoids Q14: the zero vertices have canonical intersection at most one,
and the positive tail has no other external intersection.

\paragraph{A local K3-forest criterion forces the minimal Q5 node}
At Q5 the unit quadratic coefficient of $z^2$ permits equivariant
Morse splitting.  The exact nodal discriminant order three gives
\[
 f=z^2+x^2g(x,y),\qquad g(0,y)=a y^3+\cdots,\quad a\ne0.
\]
The normalization is obtained by $z=xt$:
\begin{equation}
 t^2+g(x,y)=0,\qquad\mu_5(x,y,t)=(1,4,1),
                                                               \label{eq:octic46:octic_degree18_node46:45-D18-6}
\end{equation}
with equation character two.  The actual form has a nonzero leading
$dx\wedge dy$ term.

We claim that $[x^2]g\ne0$.  If it were zero, apply the canonical-pole
criterion to $x\operatorname{Res}(dx\wedge dy\wedge dt/(t^2+g))$.
The residue has character four, so this form is invariant and regular.
Weights $(x,y,t)=(1,4,6)$ give threshold twelve.  The only lower
character-two monomials after Morse splitting are $x^2,x^3y,x^7$.
Thus $x^3y$ or $x^7$ must occur.

If $x^3y$ occurs, the isolated initial cone is
$t^2+y^3+x^3y$, with weights $(4,6,9)$ and degree eighteen.
After the radial $\mu_5$ quotient the graph is
\begin{equation}
                  \text{centre }-2;\quad A_3,\quad A_1,\quad[3].
                                                               \label{eq:octic46:octic_degree18_node46:45-D18-7}
\end{equation}
Its raw main square is $-5/12$ and the three corrections are
$3/4,1/2,1/3$.  The actual $dx\wedge dy$ has central order
$(4+6)/5-1=1$.  Each adjacent $(-2)$ arm is positive, and its first
coefficient must differ from one, hence is at least two.  Already
these two arms contradict adjunction at the central $(-2)$ curve:
its coefficient twice is two, whereas their sum is at least four.

If $x^3y$ is absent, $x^7$ must occur.  The isolated initial cone
$t^2+y^3+x^7$ has weights $(6,14,21)$ and degree forty-two.
Its quotient graph is
\begin{equation}
                  \text{centre }-2;\quad A_1,\quad A_2,\quad[2,2,3].
                                                               \label{eq:octic46:octic_degree18_node46:45-D18-8}
\end{equation}
The raw square is $-5/42$ and the corrections are
$1/2,2/3,5/7$.  The actual central order is
$(6+14)/5-1=3$.  On its $A_1$ chart the leading form is a unit times
$T^3u\,dT\wedge du$, giving order two.  On its $A_2$ chart it is a
unit times $T^3dT\wedge dv$, giving orders $(2,1)$.
Higher terms can only lower these upper bounds.  Positivity and
adjunction force the exact coefficients $2$ and $(2,1)$ respectively.
The path consisting of the first $A_1$ vertex, the centre, and the
first $A_2$ vertex then has coefficients $2,3,2$, contradicting the
unique-minimum rule.  Both alternatives are impossible, proving the
claim.  All these weighted faces are isolated, so the displayed graphs
are unchanged by higher terms.

With $[x^2]g\ne0$, the Q5 normalization is the quotient of
$t^2+x^2+y^3$ by $\mu_5(1,4,1)$.  Its complete graph is
\begin{equation}
                         [2,2,3,2,2]=\tfrac1{15}(1,11).
                                                               \label{eq:octic46:octic_degree18_node46:45-D18-9}
\end{equation}
On the universal plane cover, $t+ix=u^3$, $t-ix=v^3$, $y=-uv$,
and $dx\wedge dy$ is a nonzero constant times
$(u^3+v^3)du\wedge dv$.  Its order is zero on all five rays.
Thus Q5 supplies four protected roots and every polar component of
multiplicity at least two avoids Q5.

\paragraph{The normal Q12 source must supply five more protected roots}
If the Q14 form is in its zero-order pattern, every multiple polar
component would avoid both Q5 and Q14.  If it is in the second pattern,
only the one component $B_3$ can fail to avoid Q14.
For a curve avoiding Q5 and Q14, all remaining ambient indices divide
$12$.  Its scaled degree is consequently a positive multiple of $35$.
A multiplicity-two contribution would be at least $70>36$.
Therefore the only possible multiple intrinsic polar component is
$B_3$ of coefficient three.  Every intrinsic polar coefficient is at
most three.

The normal-Q12 boundary lemma of
\S\ref{sec:q12-normal-boundary45} now applies.
It forces $x^3z,z^4$, and its complete resolution has no exceptional
$(-1)$ curve.  The ordinary $x$--axis chart supplies an $A_5$.
If $\delta_z(0)=0$, its attached main component is positive, so all
five $A_5$ vertices are positive.  The increasing-child argument,
using the intrinsic coefficient bound three, bounds every positive
minimal exceptional coefficient by two.  This is impossible on a
strictly monotone positive $A_5$ chain.  Hence $\delta_z(0)\ne0$
and the five Q12 $A_5$ roots are protected.

\paragraph{The final rank contradiction covers both Q14 patterns}
If all Q14 exceptional orders are zero, its twelve $(-2)$ curves are
protected.  Q5 and Q12 add four and five disjoint protected roots,
giving negative rank $12+4+5=21>19$.

In the remaining pattern, contract the Q14 positive tail of coefficients
$1,2$ together with its canonical component.  Its coefficient-one root
is the first tail curve, private to Q14.  The adjacent zero $(-3)$
curve gains one in self-intersection and becomes $-2$.
The main and $x$--arm zero $(-3)$ centres each meet a coefficient-one
intrinsic root, different from that private component.  These two
roots may be distinct or equal.  After all contractions, the thirteen
surviving Q14 curves therefore have one of two explicit intersection
matrices.  Their determinants are respectively
\[
                              18,\qquad 180.
\]
Both have rank thirteen.  The four protected Q5 roots and the five
protected Q12 roots are orthogonal to this space, because protected
curves are disjoint from the entire canonical forest.  Their combined
intersection rank is $13+4+5=22>20$, impossible on a K3 surface.
This excludes the nodal alternative.
\end{proof}


\paragraph{The cuspidal degree-eighteen alternative.}
\label{sec:octic-degree18-cusp}

\begin{theorem}
\label{thm:octic-degree18-cusp}
The sole degree--18, generic delta-one conductor state of C2-02350
cannot have cuspidal generic transverse type.  Together with
Theorem~\ref{thm:octic-degree18-node}, this deletes the entire
numerical degree--18 state.
\end{theorem}

\begin{proof}
The conductor image again has smooth component stack
$\mathcal G=\mathbb P(5,14)$ with $H=O(3)$, meeting the ordinary Q5
source along its $y$--axis and Q14 along its $x$--axis.
Its reduced normalized conductor $R$ has
\begin{equation}
 C_N=2R,\qquad420H\cdot R=18,\qquad420H\cdot B_N=36.
                                                               \label{eq:octic46:octic_degree18_cusp46:45-C18-1}
\end{equation}
The source-support guard in Section~\ref{sec:octic-conductor-reduction46}
shows that this conductor avoids
Q12: its complete source marks involve only Q5 and Q14, and a smooth
coarse conductor cannot have an invisible inertia-one contact with an
ordinary source.  Indeed every nonconstant invariant ambient function
has ordinary order at least two and cannot provide a coarse order-one
uniformizer on such a branch.

If the second jet $\gamma_{10}$ vanished on $\Gamma$, the normalized
ambient form of a generic $A_2$ cusp would have order at least two
along $R$.  Equation \eqref{eq:octic46:octic_degree18_cusp46:45-C18-1} would force $B_N=2R$.
But the octic is normal at Q12, where its actual polar numerator
$\theta_9$ has nontrivial character $-9\cdot5\equiv3\pmod{12}$.
The nontrivial polar character forces a divisorial zero through
that ordinary source.
Neither $R$ nor its conductor contribution meets Q12.  This contradicts
the exhausted polar degree.  Hence $\gamma_{10}|_\Gamma\ne0$.

The exact coarse degree of $10H=O(30)$ is zero, with marks
$0$ and $3/7$ at Q5 and Q14.  The second jet is therefore a unit at
Q5 and has exact order six on the Q14 cover.  At both sources it forces
$\delta_z(0)\ne0$.  At a generic cusp the corresponding ambient form
has order exactly one on the normalized conductor, so
\begin{equation}
             B_N=R+B',\qquad420H\cdot B'=18,\qquad R\not\subset B'.
                                                               \label{eq:octic46:octic_degree18_cusp46:45-C18-2}
\end{equation}
Every irreducible intrinsic polar component has positive scaled
$H$--degree at least three.
Thus every coefficient of $B'$ is at most six.

\paragraph{Five protected roots at the normal Q12 source}
The complete normal-Q12 boundary lemma applies.  Its ordinary $x$--axis
chart contains an $A_5$, and all components of the complete resolution
have square at most $-2$.  If $\delta_z(0)=0$ there, the component
attached to the $A_5$ has positive canonical order, and all five chain
coefficients are positive.  Following increasing blowup-child
coefficients from any minimal exceptional curve must leave its finite
fiber along an intrinsic polar component.  The coefficient bound six
from \eqref{eq:octic46:octic_degree18_cusp46:45-C18-2} consequently bounds all these exceptional coefficients
by five.

A positive $A_5$ attached to a positive component cannot satisfy that
bound.  Its coefficients are strictly monotone, hence would be either
$1,2,3,4,5$ or $5,4,3,2,1$.  In the first case the first adjunction
inequality forces the attaching coefficient to be zero.  In the second,
the attaching coefficient cannot equal five, and the two-smaller-parent
rule at the first vertex again forces it to be zero.  Thus
$\delta_z(0)\ne0$ and the Q12 $A_5$ is protected.

\paragraph{Ten protected roots at the cuspidal Q14 source}
At Q14, $f\in(y,z)^2$ has character two.  Write its quadratic normal
part as $A(x)y^2+2B(x)yz+C(x)z^2$.
The exact order-six second jet forces
$C(x)=c x^6+\cdots$ with $c\ne0$ and $\delta_z(0)\ne0$.
Generic cuspidality makes the quadratic rank one, so
$A(x)C(x)=B(x)^2$.  Since $\operatorname{ord}B\ge12$,
\[
                         \operatorname{ord}A\ge18.
\]
In particular the $x^4y^2$ coefficient vanishes.
The invariant regular conductor forms $y^3\Omega$ and $xz^2\Omega$
force $y^3z$ and $xz^3$.  With weights $(3,11,15)$ the complete possible
initial equation is therefore
\begin{equation}
                 xz^3+y^3z+a x^5y^3+b x^6z^2,\qquad b\ne0.
                                                               \label{eq:octic46:octic_degree18_cusp46:45-C18-3}
\end{equation}
The coefficient $a$ is retained; it is not eliminated by a coordinate
change that would move the conductor.
Its $y$--axis chart is smooth before the quotient and has residual
$\frac1{11}(1,10)$.  The unit $\delta_z$ term makes
$dx\wedge dy$ a unit canonical form on that chart.  Its $A_{10}$ is
therefore protected, exactly as in the two-component-state proof.

\paragraph{The Q5 cusp must have a nonzero quadratic term after normalization}
At Q5 the unit second jet forces a nonzero quadratic $z^2$ term.
Equivariant analytic Morse splitting, preserving the conductor ideal
$(x,z)$, gives
\[
                   f=z^2+x^3h(x,y).
\]
Generic delta-one cuspidality means $h(0,y)$ is not identically zero.
Adjoining $t=z/x$ normalizes the surface to
\begin{equation}
 F=t^2+G(x,y)=0,\qquad G=xh(x,y)\in(x),\qquad
 \mu_5(x,y,t)=(1,4,1),\quad\operatorname{char}F=2.
                                                               \label{eq:octic46:octic_degree18_cusp46:45-C18-4}
\end{equation}
The leading actual ambient form is a unit multiple of $dx\wedge dy$.

The residue of \eqref{eq:octic46:octic_degree18_cusp46:45-C18-4} has character four, so $x\Omega_N$ is an
invariant regular canonical form.  With weights $(x,y,t)=(2,3,7)$
its threshold is fourteen.  After Morse splitting the only
character-two equation monomials below fourteen are
\[
                         x^2,\qquad x^3y,\qquad y^3.
\]
The last is excluded by $G\in(x)$.  Rationality therefore forces
$x^2$ or $x^3y$.

Suppose $x^2$ were absent.  Then $x^3y$ is nonzero.  Use weights
$(x,y,t)=(2,8,7)$, of degree fourteen.  The entire initial equation is
\[
                         t^2+x^3y+b x^7.
\]
All further terms have weight at least twenty-four.  In the $x$--chart
put $x=s^2$, $y=s^8v$, $t=s^7u$, and $S=s^5$.
The source $\mu_5$ is radial, and the residual action is
$\mu_2(S,u)=(-S,-u)$.  The strict equation gives
\[
                         v=-u^2-b+O(S^2).
\]
This chart is already normal and its quotient is $A_1$.
The actual form has expansion
\begin{equation}
                 \alpha=(aSu+b' S^2+\text{higher terms})
                              dS\wedge du,\qquad a\ne0.
                                                               \label{eq:octic46:octic_degree18_cusp46:45-C18-5}
\end{equation}
Indeed $dx\wedge dy=(2/5)S\,dS\wedge dv$ supplies $Su$;
the possible $S^2$ term includes the ambient $\delta_y$ contribution.
All other actual terms have central order at least two.
Thus the main exceptional curve has order one, and on the $A_1$ ray
$(1/2,1/2)$ both displayed monomials have order one as well.
They cannot cancel identically along that exceptional curve because
$u/S$ is a nonconstant coordinate there.  The two adjacent exceptional
curves consequently have equal positive coefficient one, impossible
in the canonical forest of point blowups of a K3 surface.
Resolving any other chart does not change this local adjacency.
This proves
\begin{equation}
                              [x^2]G\ne0.                  \label{eq:octic46:octic_degree18_cusp46:45-C18-6}
\end{equation}

\paragraph{Six protected Q5 roots for every higher cusp coefficient}
The character of $h(0,y)$ is one, so its finite order is
\[
                            n=4+5k,\qquad k\ge0.
\]
A parameter-dependent equivariant Morse change in $x$ changes \eqref{eq:octic46:octic_degree18_cusp46:45-C18-4}
to
\[
                         t^2+X^2+y^{2n}=0
\]
up to nonzero constants and units; the critical value has order exactly
$2n$ because $G=xh$ has two roots whose separation has order $n$.
The changes have $dx\wedge dy$ equal to a unit times $dX\wedge dy$.
The universal plane description is
\[
 t+iX=u^{2n},\qquad t-iX=v^{2n},\qquad y=uv.
\]
The full group is cyclic of type
\begin{equation}
                  \frac1{10n}(1,8n-1).
                                                               \label{eq:octic46:octic_degree18_cusp46:45-C18-7}
\end{equation}
Indeed the two $2n$th powers have source character one and $uv$ has
source character four.  Since $n\equiv4\pmod5$, the action is free
away from the origin and $(8n-1)^2\equiv1\pmod{10n}$.

The first three and last three resolution curves are $(-2)$--curves.
The first three fan rays are
\begin{equation}
 \frac{(8n-1,1)}{10n},\qquad
 \frac{(6n-2,2)}{10n},\qquad
 \frac{(4n-3,3)}{10n},                                   \label{eq:octic46:octic_degree18_cusp46:45-C18-8}
\end{equation}
and the last three are their coordinate exchanges.  They are six
distinct rays for every $n\ge4$.  The actual leading two--form is a
nonzero constant times
\[
                      (u^{2n}+v^{2n})\,du\wedge dv.
\]
It has order zero on all six rays: for the $j$th ray on the left,
its numerator coordinates satisfy
$A_j+(2n+1)B_j=10n$.
Here the assertion about the actual form, rather than an arbitrarily
chosen canonical numerator, follows directly from normalization.
In the Morse coordinates the ambient form is
\[
 (\delta_z-t\delta_x)\,dx\wedge dy
       +x\delta_x\,dy\wedge dt-x\delta_y\,dx\wedge dt,
 \qquad z=xt.
\]
The coefficient $\delta_z$ is a unit, whereas $\delta_x$ and
$\delta_y$ have characters four and two and vanish at the origin.
On the universal plane, $x$ and $t$ have order at least $2n$,
$y$ has order two, and the change from $x$ to $X$ has unit
$X$--derivative.  The numerator of $dx\wedge dy$ has degree $2n$;
the numerators of $dy\wedge dt$ and $dx\wedge dt$ have degrees at
least $2n$ and $4n-2$.  Thus the terms involving $t\delta_x$,
$x\delta_x$ or $x\delta_y$, and every nonconstant change in the unit
coefficient, have numerator degree strictly greater than $2n$.
This proves that every higher actual form term has larger total degree
on this plane.
Since $A_j\ge B_j$, such a term has strictly positive order on that
ray and cannot cancel the displayed zero-order initial form.
The same argument applies on the right.  Thus all six outer roots are
protected, without assuming $k=0$ or any minimal cubic coefficient.

Finally Q5, Q12 and Q14 supply mutually disjoint protected root ranks
six, five and ten.  Their negative definite direct sum has rank
\[
                             6+5+10=21>19,
\]
contradicting the K3 minimal model.  The cuspidal alternative is
therefore impossible, and the state $(18,1)$ is excluded.
\end{proof}


\subsection{C2-02350: exclusion of the state $(21,1)$}
\label{sec:octic-degree21-node45}

Consider the sole degree-twenty-one conductor of
Theorem~\ref{thm:octic-eleven-states46}, with marks $(1/4,4/5)$
and its actual index-four assignment.  The theorem proves that its
generic transverse type is nodal, its second polar jet is nonzero,
and its first-polar length is exactly two.
Thus the full conductor consumes forty-two polar units, and
\begin{equation}
                        420H B_N=30.
\label{eq:octic46:octic_degree21_46:45-T1}
\end{equation}
The octic is normal at the ordinary sources $Q_{12}$ and $Q_{14}$.
The uniform K3 theorem applies to the normalization resolution.

\begin{theorem}
\label{thm:octic-degree21-node45}
The sole degree-twenty-one conductor state cannot occur in a
volume-$1/420$ realization of \texttt{C2-02350}.
\end{theorem}

\paragraph{The fifth source necessarily supplies two protected $A_4$ chains}

At the fifth source the nonzero second jet gives a nonzero quadratic
$z^2$ coefficient and $\delta_z(0)\ne0$.
The exact conductor in the equivariant Morse plane is the cusp
\[
 g=x^2-y^3=0,\qquad
 \operatorname{char}(x,y,z)=(1,4,2)\pmod5.
\]
Because the generic transverse type is nodal, the actual local octic is
\begin{equation}
                  f=z^2+g^2h(x,y),\qquad\operatorname{char}h=0.
\label{eq:octic46:octic_degree21_46:45-T2}
\end{equation}
The absence of further conductor components makes $h$ reduced and
coprime to $g$ whenever it is not a unit.  Its finite normalization is
\begin{equation}
                  V^2+h(x,y)=0,\qquad V=z/g,
                  \operatorname{char}V=0.
\label{eq:octic46:octic_degree21_46:45-T3}
\end{equation}
In particular its residue form
$\Omega=dx\wedge dy/(2V)$ is invariant.

Suppose $h(0)=0$.  If $[xy]h=0$, the auxiliary weights
$(x,y,V)=(1,1,2)$ give weighted order at least four, equal to the
canonical-form threshold $1+1+2=4$ for the invariant numerator one.
Indeed $xy$ is the only nonconstant invariant monomial of weight
less than four.  The uniform weighted canonical-form inequality
therefore contradicts rationality of the coarse normalization.
Hence $[xy]h\ne0$.

Equivariant Morse splitting now makes \eqref{eq:octic46:octic_degree21_46:45-T3} an $A_1$ cover with
the symplectic action $\mu_5(1,4,0)$.  Its quotient is $A_9$:
with $X=x^5$, $Y=y^5$, its invariant ring is
\[
                 \mathbb C[[X,Y,V]]/(XY+cV^{10}),\qquad c\ne0.
\]
The actual ambient form on the normalization is $q\Omega$, where
\begin{equation}
                   q=2V\delta_z+2h\,\delta g+g\,\delta h.
\label{eq:octic46:octic_degree21_46:45-T4}
\end{equation}
Since $\delta_x(0)=\delta_y(0)=0$ and $\delta_z(0)\ne0$,
the coefficient of $V$ in $q$ is nonzero.  The other terms lie in
$(V^2,X,Y)$ in the invariant ring.  On the nine $A_9$ rays,
the orders of $V,X,Y$ are $1,i,10-i$, respectively.
Consequently the actual form has order one on every interior ray
$2\leq i\leq8$.  This gives adjacent positive equal coefficients,
impossible in the point-blowup canonical forest over a K3 surface.
The argument allows possible endpoint cancellation and does not
replace the actual two-form by a chosen residue form.

Thus $h(0)\ne0$.  The normalization cover has two smooth branches
$V=\pm\sqrt{-h}$, preserved individually by $\mu_5$.
Each quotient is $\frac15(1,4)$, hence $A_4$.
Equation \eqref{eq:octic46:octic_degree21_46:45-T4} is a unit on both branches, so all eight exceptional
curves have zero actual two-form order.  They are protected roots.

\paragraph{A global degree-five gap and the normal twelfth source}

Every intrinsic positive polar curve avoids those two $A_4$ fibres.
Indeed a zero-order $(-2)$-curve has canonical intersection zero
and cannot meet any positive canonical component.
Proposition~\ref{prop:octic-sources46} makes $84H$ Cartier away from
the unique fifth source.  It follows that every such polar curve
$L$ satisfies
\begin{equation}
                         420HL\in5\mathbb Z_{>0}.
\label{eq:octic46:octic_degree21_46:45-T5}
\end{equation}
In particular every intrinsic polar coefficient is at most six.
Following a larger-coefficient child out of any minimal resolution
fibre, as in the canonical-forest lemma, bounds every positive
exceptional coefficient by five.

The normal-$Q_{12}$ boundary lemma supplies an $A_5$ attached to
a main curve.  If $\delta_z(0)=0$ there, the main curve and all
five arm curves are positive.  A positive $A_5$ attached to a
further positive curve cannot have all coefficients at most five.
Thus
\begin{equation}
                         \delta_z(Q_{12})\ne0.
\label{eq:octic46:octic_degree21_46:45-T6}
\end{equation}
The $A_5$ is protected.

Proposition~\ref{prop:octic-normal-q12-boundary} also supplies
a protected $A_1$, disjoint from the $A_5$, for every specialization
of the remaining normal-$Q_{12}$ coefficients.

We have proved the existence of fourteen disjoint protected roots:
\begin{equation}
                       2A_4+A_5+A_1.
\label{eq:octic46:octic_degree21_46:45-T7}
\end{equation}

\paragraph{A ramification bound at $Q_{12}$.}
Every intrinsic polar component through $Q_{12}$ has coefficient at
most three.  Here is a local proof including the coordinate and
different assertions.  Divide the actual vector field by its
invariant unit coefficient $\delta z$; then $\delta' z=1$.
Solve the analytic flow with initial transverse plane $z=0$.
The first two coordinates of its inverse flow are holomorphic first
integrals $X,Y$, and $Z=z$ is its time coordinate.  Uniqueness of
this flow makes the construction equivariant: $X,Y,Z$ retain
characters $1,11,5$.  Also $X=x+O(\mathfrak m^2)$ and
$Y=y+O(\mathfrak m^2)$, because the other vector-field coefficients
vanish at the origin.  Thus $\delta'=\partial_Z$ and the actual
ambient form, obtained by contracting a local volume form, is a
unit times $dX\wedge dY$.

The degree-four initial surface equation is unchanged by this
coordinate change.  In particular its restriction to $X=Y=0$
begins with the forced nonzero term $cZ^4$.
Weierstrass preparation gives a monic degree-four equation in $Z$,
so the normal cover is finite of degree four over the smooth
$(X,Y)$-plane.  At a generic ramification curve take a transverse
base parameter $t$ and a parameter along its image.
After a separable residue-field extension, the corresponding
upstairs transverse parameter $u$ satisfies $t=u^e$ times a unit.
Pullback of the base two-form therefore has order exactly $e-1$.
Characteristic zero makes this extension tame, and its degree
four gives $e\leq4$.  The source quotient is etale in codimension
one, so the same coefficient bound holds on $N$.

At least one intrinsic polar curve passes through $Q_{12}$: the
normal-source support lemma applies to the actual degree-nine polar,
whose source character is nontrivial.  The local coefficient bound
will therefore exclude any proposed polar allocation consisting of
only one curve of coefficient greater than three.

\paragraph{The complete fourteenth-source list becomes impossible}

The local atlas and Corollary~\ref{cor:octic-q14-J} leave the
formats $G,H,K,L$ and the degree-thirty-eight and degree-forty-eight
formats.  Their resolution descriptions allow all higher-order terms.

For any connected positive exceptional core, the sum of the intrinsic
polar contributions at its vertices is at most six, by \eqref{eq:octic46:octic_degree21_46:45-T1} and
\eqref{eq:octic46:octic_degree21_46:45-T5}.  A positive polar component meets such a core at most once,
since otherwise the canonical forest would contain a cycle.
Proposition~\ref{prop:octic-q14-positive-cores} applies with
this stronger bound.  Its $K$ possibility has external sum eleven,
so only the $H$ core remains.
The results are:
\begin{equation}
\begin{array}{c|c}
 \text{positive main format}&\text{surviving necessary core}\\\hline
 G&\text{none}\\
 H&\text{centre }2,\ A_3=(3,4,5),\ A_1=1;\ \text{one external }6\\
 K&\text{none}\\
 L&\text{none}.
\end{array}
\label{eq:octic46:octic_degree21_46:45-T8}
\end{equation}
The main curves in $H,L$ always have positive actual order.
The one $H$ survivor requires a single coefficient-six intrinsic
polar curve.  By \eqref{eq:octic46:octic_degree21_46:45-T5} this consumes all thirty polar units and
leaves no curve which can meet $Q_{12}$ with coefficient at most
three.  It is impossible.

In $K$, a zero main order requires $\delta_z(0)\ne0$.
The fixed $A_6$ and the first $(-2)$ in its $[2,4]$ arm are then
protected.  With \eqref{eq:octic46:octic_degree21_46:45-T7} this gives twenty-one negative root ranks,
impossible on a K3 surface.  Its positive-main case is excluded by
\eqref{eq:octic46:octic_degree21_46:45-T8}, or directly by its positive $A_6$ and the coefficient bound.

For the degree-thirty-eight and degree-forty-eight strata, a nonzero
$\delta_z(0)$ protects respectively $A_8$ and $A_{10}$, contradicting
\eqref{eq:octic46:octic_degree21_46:45-T7}.  If $\delta_z(0)=0$, each of these long arms attaches to a
positive main curve.  All curves in the complete minimal resolution
have square at most $-2$, even in the degenerate auxiliary charts.
Their positive coefficients would all be at most five, again
impossible.  Thus both complete strata are excluded.

It remains to discuss $G$ with zero main coefficient.  Its weight-fourteen
actual form is
\begin{equation}
                      A z^2dy\wedge dz+B dx\wedge dy.
\label{eq:octic46:octic_degree21_46:45-T9}
\end{equation}
If $B\ne0$, its $A_4$ arm and at least one $A_2$ arm are protected:
the two $A_2$ unit coefficients are proportional to
$6B+2iA$ and $6B-2iA$, which cannot both vanish.
These six roots together with \eqref{eq:octic46:octic_degree21_46:45-T7} give forbidden negative rank
twenty.  If $A=B=0$, the main coefficient is positive and \eqref{eq:octic46:octic_degree21_46:45-T8}
already applies.  Hence only $B=0$, $A\ne0$ needs consideration.

\paragraph{The last $G$ allocation violates the actual source congruence}

In this last case the two $A_2$ arms are protected.  At the $A_4$
point use the $y$-chart of weights $(9,5,3)$:
\[
 y=s^5,\qquad x=s^9u,\qquad z=s^3v,\qquad S=s^{14},
 \qquad v=-u^2+\cdots.
\]
The residual action is $\mu_5(S,u)=(4,1)$ and the main curve is
$S=0$.  The actual form coefficient on its smooth cover restricts
at $S=0$ to a nonzero multiple of $u^5$ times a unit.
Every $A_4$ coefficient is positive and its $i$th coefficient,
starting from the main curve, is at most $i$.
The strict-monotonicity forest lemma therefore forces the exact
sequence
\begin{equation}
                            (1,2,3,4).
\label{eq:octic46:octic_degree21_46:45-T10}
\end{equation}
At the free end adjunction gives total external contribution five.
Its required larger child has coefficient at least five, so this
contribution is one intrinsic polar curve $L$ of coefficient five.
The global mass and \eqref{eq:octic46:octic_degree21_46:45-T5} force $420HL=5$.

The outer-end curvette has the ordinary source $H$-mark $5/14$.
One can read this directly from the chart: a lifted curvette has
orders $(S,u)=(1,4)$, so its branch on the original source cover
has a smooth $y$ parameter.  That parameter has character thirteen,
whereas the canonical polarization has character $-5=9$ modulo
fourteen.  The normalized character is
\begin{equation}
                      9\cdot13^{-1}=5\pmod {14}.
\label{eq:octic46:octic_degree21_46:45-T11}
\end{equation}
Higher terms in the weighted equation do not change this parameter.

The curve $L$ avoids $Q_5$.  It also avoids $Q_{12}$, because its
coefficient five exceeds the ramification bound three.  It has no
second branch at $Q_{14}$: a second intersection with the positive
$A_4$ would create a forest cycle, while the zero main curve has
canonical intersection two and the protected $A_2$ curves have
canonical intersection zero, both less than five.
All remaining actual canonical indices divide twelve.  The degree
congruence for the integral Weil polarization on the normalization
of $L$ consequently gives
\begin{equation}
                  420HL\equiv420\frac5{14}
                        \equiv10\pmod {35}.
\label{eq:octic46:octic_degree21_46:45-T12}
\end{equation}
This contradicts $420HL=5$.
Every rational fourteenth-source format has been excluded, proving
Theorem~\ref{thm:octic-degree21-node45}.


\subsection{C2-02350: exclusion of the state $(14,2)$}
\label{sec:octic-degree14-delta2-45}

\begin{theorem}
\label{thm:octic-degree14-delta2-45}
The sole conductor state of scaled degree fourteen and generic
transverse delta invariant two cannot occur in a volume-$1/420$
realization of \texttt{C2-02350}.
\end{theorem}

\begin{proof}
Theorem~\ref{thm:octic-eleven-states46} gives generic transverse
type $A_3$, a nonzero second polar and first-polar length four.
The conductor curve is a smooth coarse $\mathbb P^1$ with marks
$(1/5,5/6)$ at $Q_5,Q_{12}$.  We obtain a contradiction at the
twelfth source by constructing its full normalization and then applying
the weighted canonical-form inequality.

At $Q_{12}$ the ambient characters are $(1,11,5)$.
The mark $5/6$ has coordinate lower orders $(5,1,1)$.
The actual degree-nine trace bound is saturated at degree $-1$.
If the order of $z$ rose from one to seven, the smallest
character-three numerator order would rise from three to nine,
violating that bound.  Thus $z$ has order one.
Coarse smoothness forces $y$ to have order one as well: without it
there is no full-source invariant of order six.
The two conjugate smooth branches can consequently be straightened
equivariantly to
\begin{equation}
                  I_\Gamma=(x,g_0),\qquad g_0=z^2-y^2.
\label{eq:octic46:octic_degree14_delta2_46:45-U7}
\end{equation}
The group interchanges the two branches and a single branch has
stabilizer $\mu_6$.

The degree-ten second jet has exact order two on this lifted
conductor: its global underlying degree is zero, and its only
nontrivial local fractional contribution is $2/6$ here.
Among the ordinary degree-four equation terms, membership in
$I_\Gamma^2$ forces the restriction at $x=0$ to be a multiple
of $(z^2-y^2)^2$.  The exact second-jet order therefore forces
its coefficient to be nonzero and $\delta_z(0)\ne0$.
Thus
\begin{equation}
                      f|_{x=0}=g_0^2\cdot\text{unit}.
\label{eq:octic46:octic_degree14_delta2_46:45-U8}
\end{equation}

There is an elementary square-lifting argument which works even
though $g_0$ is not an ambient coordinate at the origin.
Write $f=\sum_{j\geq0}x^jF_j(y,z)$ and absorb the invariant unit
in \eqref{eq:octic46:octic_degree14_delta2_46:45-U8} into $G_0=g_0\sqrt{\text{unit}}$.
Membership in $I_\Gamma^2$ makes $F_1$ divisible by $G_0$.
Set $G_1=F_1/(2G_0)$.
The rank-one transverse quadratic at the generic points of both
branches makes $F_2-G_1^2$ vanish there; since $G_0$ is reduced,
it is divisible by $G_0$ in $\mathbb C[[y,z]]$.
Set $G_2=(F_2-G_1^2)/(2G_0)$.
The generic $A_3$ type similarly kills the cubic term after these
square completions, so $F_3-2G_1G_2$ is divisible by $G_0$.
Taking $G_3=(F_3-2G_1G_2)/(2G_0)$ gives
\begin{equation}
       G=G_0+xG_1+x^2G_2+x^3G_3,
       \qquad f=G^2+x^4h,
       \qquad \operatorname{char}G=10,
       \quad \operatorname{char}h=4.
\label{eq:octic46:octic_degree14_delta2_46:45-U9}
\end{equation}
All divisions are by a reduced equation of both actual branches;
no meromorphic coefficient along the conductor is used.

Adjoin $V=G/x^2$, of character eight.  The finite birational algebra
is the complete intersection
\begin{equation}
                 x^2V=G,\qquad V^2+h=0.
\label{eq:octic46:octic_degree14_delta2_46:45-U10}
\end{equation}
It has no component contained in $x=0$, since the latter quotient
has dimension one.  It agrees with the original surface for
$x\ne0$, and at the generic points over $x=0$ it is precisely the
smooth normalization of a tacnode.  Its codimension-one points are
therefore normal, and the complete-intersection property gives
$S_2$.  Thus \eqref{eq:octic46:octic_degree14_delta2_46:45-U10} is the actual normalization.
The curve $x=0$ is Cohen--Macaulay and generically reduced, hence is
the actual reduced conductor $R$ at this source.

The normalization algebra is generated over the original surface
ring by $V$, with $V^2=-h$ and $x^2V=G$.
Thus $x^2$ is in the full conductor, and
\begin{equation}
           x^3\operatorname{Res}
                   \frac{dx\wedge dy\wedge dz}{f}
\label{eq:octic46:octic_degree14_delta2_46:45-U12}
\end{equation}
is a regular form on the normalization.  It is invariant: its total
character is $3+(1+11+5)-8=0\pmod {12}$.

Use auxiliary weights $(x,y,z)=(1,2,2)$.
The canonical-form threshold is $1+2+2+3=8$, and the only
character-eight equation monomial of smaller weight is $x^3z$.
Rationality would therefore force $[x^3z]f\ne0$.
But in $f=G^2+x^4h$, a character-ten $G$ has no linear term;
its quadratic part involves only $y^2,z^2$.
The character-four function $h$ has no constant term.
Consequently the ordinary degree-four part of $f$ is a square in
$y,z$ and contains no $x^3z$, a contradiction.
\end{proof}


\subsection{C2-02350: exclusion of the state $(12,1)$}
\label{sec:octic-degree12-cusp-reduction45}

The remaining degree--twelve conductor is necessarily a generic
delta--one cusp.  Its known marks are
$(1/2,3/5,13/14)$ at the actual sources Q2,Q5,Q14, and
$H\Gamma=1/35$.  Put $C_N=2R$; then
\[
                420HR=12,\qquad420HB_N=48.
\]
The statements below do not assume that the full canonical-cover
conductor is smooth at the possibly nonordinary second source.

\begin{proposition}
\label{prop:octic-degree12-cusp-reduction45}
In a realization of this state the second polar is nonzero.  In the
forced ordinary Q14 normalization chart described below, the actual
tangential derivation component $\delta_z$ must be a unit.  Consequently
all normalization-resolution exceptional curves have actual two-form
order zero, and the intrinsic canonical divisor is one of
\begin{equation}
                 B_N=R_{12}+L_{36},
       \qquad\text{or}\qquad
                 B_N=R_{12}+L_{21}+L_{15},              \label{eq:octic46:octic_degree12_reduction46:45-D12-1}
\end{equation}
where subscripts denote scaled $H$-degrees and all coefficients are
one.  The smooth normalization resolution is therefore a blowup of
its K3 minimal model at respectively two or three distinct points.
This proposition is a necessary reduction; by itself it does not
exclude the two alternatives in \textup{\eqref{eq:octic46:octic_degree12_reduction46:45-D12-1}}.
\end{proposition}

\paragraph{The Q14 normalization is forced for both polar alternatives}

The ordinary source has characters $(x,y,z)=(1,13,5)$; the conductor
is the $z$--axis, with full ideal $(x,y)$, and $f$ has character two.
The invariant conductor numerator $xz^2$ with weights $(3,4,1)$
has canonical threshold thirteen.  The only lower character-two
terms are $x^2,xz^3,z^6$, and the last two do not belong to $(x,y)^2$.
Thus $[x^2]f\ne0$.  Similarly $y^3$ with weights $(6,1,2)$ forces
$[y^3z]f\ne0$.

Equivariant Morse splitting in $x$ and the cusp condition give
$f=x^2+y^3h(y,z)$, with $\partial h/\partial z(0)\ne0$.
Replacing $z$ by $h$ is an equivariant analytic coordinate change
preserving the conductor.  Hence the actual normalization is
\[
              t^2+yz=0,\qquad t=x/y,
              \mu_{14}(t,y,z)=(2,13,5).
\]
On its universal plane cover use
$y=a^2,z=b^2,t=iab,x=ia^3b$, with harmless constants suppressed
below.  The quotient is $\frac1{28}(1,23)$.  Its minimal resolution
has chain
\begin{equation}
                      [2,2,2,2,3,3],                  \label{eq:octic46:octic_degree12_reduction46:45-D12-2}
\end{equation}
whose rays, from the free $A_4$ end towards the conductor, are
\[
 \frac1{28}(1,23),\quad\frac1{28}(2,18),\quad
 \frac1{28}(3,13),\quad\frac1{28}(4,8),\quad
 \frac1{28}(5,3),\quad\frac1{28}(11,1).
\]
The reduced conductor meets the final $(-3)$ curve.
The actual form has plane coefficient
\begin{equation}
 \alpha=\{4ab\,\delta_x-6a^2b^2\,\delta_y
                       -2a^4\,\delta_z\}\,da\wedge db,
                                                        \label{eq:octic46:octic_degree12_reduction46:45-D12-3}
\end{equation}
where the derivation components are evaluated at
$(x,y,z)=(a^3b,a^2,b^2)$ and have source characters $(10,8,0)$.

\paragraph{The zero-second-jet alternative is impossible}

On the normalized marked conductor, the degree of the coarse
degree-eleven restriction is at most
\[
 \frac{11}{35}-\left(\frac12+\frac35+\frac3{14}\right)=-1.
\]
Thus, if the second polar vanishes, the intrinsic third jet vanishes
as well.  At a generic transverse cusp the second jet kills the
square normal component of the derivation, and the third jet kills
the radical normal component.  The derivation is therefore tangent
to the conductor.

For completeness, tangential derivations remain regular on the
normalization even if a source lift of the curve is singular.
Trivialize the twisting line on the source cover, and write the
fraction-field extension of a derivation as
$D=h(t)d/dt$ on the normalized branch.  If $h$ had order $-k<0$,
the value semigroup of the original branch would be closed under
subtracting $k+1$ from every positive element.  In characteristic
zero its leading derivative coefficient never vanishes.  Iteration
would make every semigroup value divisible by $k+1$, contradicting
that the normalization is birational.  Hence $h$ has no pole.
This argument is compatible with the finite source action and local
frames of the twisting line.

The orbifold canonical degree of the normalized marked curve is at
least
\[
 -2+(1-1/2)+(1-1/5)+(1-1/14)=8/35.
\]
Consequently $\deg(T_{\mathcal G}\otimes H)\le-1/5<0$.
The induced tangential derivation is zero.  Thus $\delta$ vanishes
along the conductor, including at Q14.  In \textup{\eqref{eq:octic46:octic_degree12_reduction46:45-D12-3}} each
derivation component is then divisible by $a^2$.

On each of the first four rays of \textup{\eqref{eq:octic46:octic_degree12_reduction46:45-D12-2}}, every allowed
actual-form monomial has strictly positive order.  This can be
checked without a genericity assertion.  For the first three rays
the only plane monomial of canonical order zero that can occur in
\textup{\eqref{eq:octic46:octic_degree12_reduction46:45-D12-3}} is $a^4$; for the fourth ray the two other formal
possibilities $a^2b,b^2$ cannot occur from the source-coordinate
monomials.  The $a^4$ term would require a constant $\delta_z$,
which has just been excluded.  Therefore the whole $A_4$ is positive.

At a generic cusp, a derivation vanishing along the conductor makes
the normalized form vanish to order at least three on $R$.
Write $B_N=rR+B'$, so $r\ge3$ and
$420HB'=48-12r$.  Hence $r=3$ or $4$.  Every intrinsic polar
coefficient is at most four: for $r=3$ the remaining mass is twelve
and every component has degree at least three; for $r=4$ there is
no remaining mass.  Following increasing neighboring coefficients
from a positive exceptional curve on the minimal rational resolution
therefore bounds its coefficient by three.  The positive $A_4$
would have four strictly monotone positive coefficients all at most
three, a contradiction.  This argument needs no conormal nilpotence
or smooth normal-bundle hypothesis at Q2.

\paragraph{A nonzero polar with vanishing tangential Q14 value is impossible}

The degree-ten coarse restriction has degree zero.  A nonzero polar
therefore has exact cover order four at Q14 and is a unit at Q5.
At Q5 the preceding local theorem forces the $D_5/\mu_5$ format,
with all four actual exceptional orders zero.
At Q14 the exact order four gives
$\delta_x|_\Gamma=z^2\cdot\text{unit}$.  It does not, by itself,
kill the tangential component $\delta_z$.

Suppose first that $\delta_z(0)=0$.
The same ray enumeration as above makes the $A_4$ orders positive.
The forced $ab^5$ term in \textup{\eqref{eq:octic46:octic_degree12_reduction46:45-D12-3}} gives the upper bounds
$(4,3,2,1)$ on that chain and order zero on both $(-3)$ curves.
Other terms can lower individual bounds, so a leading-term claim
alone would not suffice.  Strict monotonicity of a positive root
chain, together with its last upper bound one, forces the exact
orders
\[
                             (4,3,2,1,0,0).
\]
Adjunction at its free end gives total external polar intersection
five and gives zero external intersection elsewhere on the chain.
The unique-minimum rule makes that external intersection a single
intrinsic component of multiplicity five.

Here the generic conductor has coefficient one, so after removing
$R$ the polar mass is $48-12=36$.  This multiplicity-five component
avoids Q5, whose zero-order exceptional curves have canonical
intersection at most one.  Its scaled degree is a positive multiple
of five and at most $36/5$, so it must be five.
However the free $A_4$ ray has original $y$-order $1/14$.
Its branch parameter has original inertia character thirteen, giving
$H$-mark $5/14$.  It has no second Q14 branch: another contact with
the positive chain gives a forest cycle, while a multiplicity-five
curve cannot meet either zero $(-3)$ curve.  All other indices divide
twelve.  Thus its degree must satisfy
\[
                        w\equiv420\frac5{14}
                                  \equiv10\pmod {35},
\]
contradicting $w=5$.

\paragraph{The remaining unit tangent branch has only two or three roots}

The only remaining alternative is $\delta_z(0)\ne0$.
The two nonzero terms $ab^5$ and $a^4$ in
\textup{\eqref{eq:octic46:octic_degree12_reduction46:45-D12-3}} together give order zero on all six Q14 curves.
All four Q5 curves also have order zero.  A multiple intrinsic polar
component consequently avoids both Q5 and Q14, so its scaled degree
is a positive multiple of thirty-five.  Its contribution would be
at least seventy, exceeding the remaining mass thirty-six.
Therefore every intrinsic canonical coefficient is one.

Every positive exceptional coefficient on a minimal normalization
resolution would have a strictly larger child path ending at an
intrinsic component.  Such an endpoint cannot have coefficient one.
There are thus no positive exceptional coefficients anywhere, including
at Q12 and at hidden contacts.  All positive canonical curves are
intrinsic coefficient-one roots, with no later blowup on them.  They
are disjoint $(-1)$-curves.

There is exactly one additional canonical branch through the Q5
central $(-3)$ curve, and one through the first Q14 $(-3)$ curve.
At Q5 the $D_5$ plane cover gives source leading orders
$(\operatorname{ord}x,\operatorname{ord}y)=(3,2)$; the branch
parameter has inertia character two, so its $H$-mark is $4/5$.
At Q14 the strict canonical branch has leading equation
$a^3+\lambda b^5=0$, $\lambda\ne0$.  On the original index cover
its source orders are $(x,y,z)=(9,5,3)$, giving inertia character
eleven and $H$-mark $11/14$.

If these are different curves, neither can have another contact at
the other source: the other zero $(-3)$ intersection is already
exhausted by $R$, and all $(-2)$ curves have canonical intersection
zero.  Their respective degree congruences are
\[
                    w_5\equiv21\pmod {35},\qquad
                    w_{14}\equiv15\pmod {35}.
\]
The total mass thirty-six forces $w_5=21,w_{14}=15$, with no further
component.  If they are the same curve, its degree is
$w\equiv1\pmod {35}$.  Degrees one and two are impossible, so the
same budget forces $w=36$.  This proves
\textup{\eqref{eq:octic46:octic_degree12_reduction46:45-D12-1}} and the claimed two- or three-point blowup
description.  The remaining global configuration still requires
an exclusion or construction argument.


\paragraph{The normal twelfth source completes the degree-twelve case.}
\label{sec:octic-degree12-q12-completion45}

This section completes Proposition~\ref{prop:octic-degree12-cusp-reduction45}.
In its surviving alternative every exceptional curve on the \emph{minimal}
normalization resolution has actual two-form order zero.  The canonical
divisor consists of either two or three disjoint $(-1)$ curves, with
scaled degrees
\begin{equation}
                       (12,36)\quad\hbox{or}\quad(12,21,15).
                                                        \label{eq:octic46:octic_degree12_complete46:45-D12Q-1}
\end{equation}
The source Q12 is outside the conductor, so its octic germ is normal.

\begin{theorem}
\label{thm:octic-degree12-cusp-excluded45}
The sole degree--twelve, delta--one cusp state of C2--02350 is impossible.
\end{theorem}

\paragraph{The fixed charts and the actual form}

The normal Q12 boundary lemma forces the ordinary quartic part
\begin{equation}
              f_4=x^3z+c z^4+b y^2z^2+a y^4,
                         \qquad c\ne0.                 \label{eq:octic46:octic_degree12_complete46:45-D12Q-2}
\end{equation}
The source action is $\mu_{12}(x,y,z)=(1,11,5)$.
The same lemma and the all-zero minimal-resolution condition force
$\delta_z(0)\ne0$.  There are always protected $A_5$ and $A_1$
charts, at the $x$--axis and the orbit $y=0$, respectively.

When $a=0$ and $b\ne0$, there is one further protected $A_2$ orbit.
It consists of the two points $x=0$, $y/z=\pm\sqrt{-c/b}$ on the
cubic component.  In the $z$--chart, after the radial $\mu_2$ quotient,
their stabilizer acts on the two smooth coordinates with characters
$(2,1)$ modulo three.  The leading form is a unit canonical form.
This $A_2$ is present in every degeneration considered below; it must
not be omitted from the resolution graph.

Only the $y$--chart can require further analysis.  Put
\[
 y=s,\qquad x=su,\qquad z=sv,\qquad S=s^2.
\]
After removing the radial $\mu_2$, its remaining action and equation are
\begin{equation}
 \mu_6(S,u,v)=(5,1,3),\qquad
 F=u^3v+c v^4+b v^2+a+S(\cdots),\qquad \operatorname{char}F=0.
                                                        \label{eq:octic46:octic_degree12_complete46:45-D12Q-3}
\end{equation}
Every normal \emph{quotient} partial chart is rational, since it is a
normal birational partial resolution of the rational coarse surface germ.
No rationality of an arbitrary canonical-cover chart is asserted.
Normality here can also
be checked directly: its exceptional initial curve is generically
reduced, and any singular locus away from that curve would give a
nonisolated singular locus on the original normal surface.

The source characters of the actual derivation give
\[
             \delta_z=A,\qquad\delta_x=s^4B,
             \delta_y=s^2C,
\]
where $A$ is a unit and $A,B,C$ are analytic after the radial quotient.
The lowest possible terms of $C$ are $uv$ and $S^2$.
Direct substitution into
$\delta_xdy\wedge dz-\delta_ydx\wedge dz+\delta_zdx\wedge dy$
gives the useful exact expression
\begin{align}
 \alpha={}&\frac A2du\wedge dS+\frac{S^2B}{2}dS\wedge dv
       +S^2C\,dv\wedge du \notag\\
       &+\frac{SC}{2}\{v\,dS\wedge du-u\,dS\wedge dv\}.
                                                        \label{eq:octic46:octic_degree12_complete46:45-D12Q-4}
\end{align}
All subsequent coordinate preparations are equivariant analytic
preparations.  Their effect on this formula is included below.

We will use the following minimality guard repeatedly.  A rational curve
with actual two-form order zero cannot be a $(-1)$ curve: its intersection
with the effective divisor of the form would be nonnegative, whereas
adjunction gives $K\cdot E=-1$.  Thus the original zero-order components
cannot disappear by a minimal-resolution contraction.  On the other
hand, a positive-order weighted exceptional curve \emph{can} be an
extraneous $(-1)$ curve.  Such a curve occurs in one case below, and is
explicitly contracted before using the all-zero hypothesis.

\paragraph{An irreducible quartic: the ordinary and double-root cases}

Suppose $a\ne0$.  If $b^2\ne4ac$, the minimal Q12 graph is the
usual centre of square $-4$ with arms $A_5, A_2, A_2, A_1$.
It contains eleven zero-order curves.  Together with the four Q5 and
six Q14 curves this gives twenty-one curves on the smooth resolution.

Let $M$ be their intersection matrix before contracting the canonical
roots.  A root meeting these curves with incidence vector $v$ adds
$vv^{\mathsf t}$ to their intersection matrix on the K3 surface.
The root of degree twelve meets the conductor ends at Q5 and Q14.
If the other root has degree thirty-six, it meets the Q5 and Q14
remaining $(-3)$ curves and meets the Q12 centre with total
intersection two.  These two intersections are allowed to lie on the
\emph{same} root.  If the other roots have degrees twenty-one and
fifteen, their Q12 incidences are $j$ and $2-j$, $j=0,1,2$.
The resulting four matrices have determinants
\begin{equation}
                    63936,\quad9792,\quad10224,\quad14976.
                                                        \label{eq:octic46:octic_degree12_complete46:45-D12Q-5}
\end{equation}
They all have rank twenty-one, impossible on a complex projective K3.

If $b^2=4ac$, the two double roots $v=\pm v_0\ne0$ in the
$y$--chart are exchanged by $\mu_6$.  At either root their stabilizer
is $\mu_3$ with characters
\[
                        (S,u,w=v-v_0)=(2,1,0).
\]
The equation starts with $w^2+\text{unit}\,u^3$ and possibly $Su$;
it has no linear term.  The quotient has invariant canonical generator
and is rational and Gorenstein, hence is a non-smooth Du Val germ.
Formula~\textup{\eqref{eq:octic46:octic_degree12_complete46:45-D12Q-4}} is a regular K\"ahler two-form on the
cover.  Its coefficient relative to the residue belongs to the
Jacobian ideal and vanishes at the origin.  Since that residue is
invariant, its descended coefficient vanishes at the quotient origin
as well.  Every exceptional curve of this Du Val resolution therefore
has positive actual order.  The new curves have square $-2$; the old
zero-order rational curves cannot have square $-1$ by the guard above.
Thus these positive curves survive on the minimal resolution, a
contradiction.

\paragraph{A line and a smooth cubic: complete equivariant preparation}

Suppose $a=0$, $b\ne0$.  Equivariant Morse splitting in $v$, keeping
$S,u$ fixed, writes the equation as
\begin{equation}
                         v^2+h(S,u)=0,
                                                        \label{eq:octic46:octic_degree12_complete46:45-D12Q-6}
\end{equation}
where $h$ is invariant, and $h(0,u)$ is a nonzero unit times $u^6$.
The last assertion follows from the nonzero $u^3v$ coefficient before
Morse splitting.  A unit normalization makes this coefficient one.
The complete list of possible terms of degree at most six is
\begin{equation}
        h=\varepsilon Su+\lambda S^2u^2+u^6
                 +\mu S^3u^3+\nu S^6+\text{higher terms}.
                                                        \label{eq:octic46:octic_degree12_complete46:45-D12Q-7}
\end{equation}
This list is exhaustive, since $S^iu^j$ has character $j-i$ modulo six.

\paragraph{The case $\varepsilon\ne0$.}
The cover is an ordinary node and its quotient is
$\frac1{12}(1,5)$, with chain $[3,2,3]$.
On a universal plane cover, the original two exceptional branches
have leading equations $a=\pm\text{unit}\,b^5$.
They meet the same end $(-3)$ curve at distinct points.
This follows either by an equivariant node preparation, or from
$S'=S+\text{unit}\,u^5+\cdots$, $u'=u+O(S^5)$ and
$S'=a^2$, $u'=b^2$, $v=ab$.
The actual form has leading coefficient $ab$ and has order zero on
all three quotient curves.  Higher terms cannot cancel that monomial
identically on any exceptional curve.

The complete Q12 graph has thirteen curves.  Denote the strict line
by $E_l$ and the strict cubic by $E_c$.  The line has square $-2$
and an $A_5$ arm.  The cubic has square $-3$ and arms $A_2,A_1$.
Both meet the same end of the chain $G(-3)-M(-2)-J(-3)$.
For clarity, the intersection corrections are
\[
 E_lE_c=\frac5{12},\qquad
 \operatorname{corr}_Y(E_l^2)=\operatorname{corr}_Y(E_c^2)=\frac5{12},
\]
and the fixed corrections are $5/6$ on the line and $2/3+1/2$ on
the cubic.  Using $E_{\rm tot}E_l=-1/3$ and
$E_{\rm tot}E_c=-1$ gives their asserted final squares.

There are three Q12 $(-3)$ centres, namely $E_c,G,J$.
Each must meet exactly one coefficient-one canonical root.
Including Q5 and Q14 gives twenty-three curves.
The degree-thirty-six alternative gives one incidence assignment;
the degree-twenty-one/fifteen alternative gives all $2^3=8$
assignments of these three centres.  Direct evaluation gives rank
twenty-three for all nine matrices, contradicting the K3 rank bound.
The matrices can also be recovered from the computation sources
described in Appendix~\ref{app:reproducibility46}.

\paragraph{The case $\varepsilon=\lambda=0$.}
Blow up weights $(S,u,v)=(1,1,3)$.  The order-two source subgroup
is radial; after removing it write $T=s^2$ in either axis chart.
The leading actual form becomes a unit multiple of $dT\wedge dU$.
All errors in \textup{\eqref{eq:octic46:octic_degree12_complete46:45-D12Q-4}} are regular K\"ahler forms with
coefficients divisible by $T^2$.  This follows directly from the
weighted orders: the $B$ term has total order at least six, and the
$C$ terms have still larger order because $C\in(uv,S^2)+\cdots$.

If $\nu=0$, the $S$--axis chart has action
$\mu_3(T,U,V)=(2,1,0)$ and equation
\[
       V^2+U^6+\mu U^3+T G_8(U)+\cdots=0,
                         \qquad G_8(0)=0.
\]
The last equality is forced by characters: a pure $S^8$ term is
not invariant.  The quotient is therefore a non-smooth rational
Gorenstein germ, hence Du Val.  The actual form is K\"ahler there,
so its residue coefficient vanishes and every Du Val exceptional
curve has positive order.  The weighted main component has order
zero, as do the original components.  None of these rational curves
can be $(-1)$; the added Du Val curves are $(-2)$.  Thus the positive
curves are on the minimal resolution, a contradiction.

If $\nu\ne0$, each of the two axis charts has two smooth fixed
points, at $V=\pm\sqrt{\nu}$ or at the two roots of the normalized
$u^6$ coefficient.  The residual action on their smooth coordinates
is $\frac13(1,2)$, and the actual form is a unit canonical form.
This supplies four disjoint protected $A_2$ chains.  With the fixed
Q12 $A_5+A_2+A_1$, the Q5 two $A_1$ curves, and the Q14 $A_4$,
there are twenty-two mutually orthogonal negative root classes.
They are untouched by all canonical contractions.  A projective K3
has negative definite N\'eron--Severi rank at most nineteen, so this
case is impossible, regardless of any additional singular points of
the weighted initial curve.

\paragraph{The case $\varepsilon=0$, $\lambda\ne0$.}
Blow up weights $(S,u,v)=(1,1,2)$.
The two initial components $v=\pm\sqrt{-\lambda}\,Su$ are exchanged
by $\mu_6$, giving one weighted exceptional curve $E$.
It has preliminary square $-1/3$ and actual order one.
This positive curve must not be used directly as an obstruction.

In the $u$--axis chart, write
$u=s$, $S=sU$, $v=s^2V$.
Parameter Morse preparation in $U$ gives
\[
                  V^2+U^2+s^2=0,
                \qquad\mu_6(s,U,V)=(1,4,1).
\]
Its quotient is $\frac1{12}(1,7)$, with chain $[2,4,2]$.
The curve $E$ meets its middle $(-4)$ curve; the old line and cubic
meet its two outside $(-2)$ curves.  In the universal plane,
the actual leading form is
\[
                 s\,ds\wedge dU
                      =\text{unit}\,(a^4-b^4)\,da\wedge db.
\]
All three quotient curves consequently have actual order zero.
The other terms of \textup{\eqref{eq:octic46:octic_degree12_complete46:45-D12Q-4}} have strictly higher
plane degree and cannot cancel either displayed leading monomial.

If $\nu\ne0$, the $S$--axis chart is another identical
$[2,4,2]$ quotient, again with $E$ meeting its middle curve and
all orders zero.  Each chart contributes $1/3$ to the correction of
$E^2$.  Thus
\[
                           E^2=-\frac13-\frac13-\frac13=-1.
\]
Contract this \emph{extraneous} order-one $(-1)$ curve.  The resulting
minimal $Y$ graph has two joined $(-3)$ centres, each with two
$(-2)$ arms.  All six curves have actual order zero.  The old line
and cubic meet the two arms of the first centre.

The inverse intersection matrix of this six-vertex graph gives
corrections $2/3$ for either old branch and $1/6$ for their mutual
intersection.  After the fixed $A_5,A_2,A_1$ corrections, the old
line and cubic have squares $-2$ and $-3$ respectively.  The complete
Q12 graph thus has sixteen curves and precisely three $(-3)$ centres.
With Q5 and Q14 this gives twenty-six curves.  The same one plus
eight possible assignments of the three centres in
\textup{\eqref{eq:octic46:octic_degree12_complete46:45-D12Q-1}} all give intersection matrices of rank
twenty-six.  Their determinants are
\begin{equation}
 -648960;\quad -111360,61760,32000,110400,
                   133440,20480,38720,-168960.
                                                        \label{eq:octic46:octic_degree12_complete46:45-D12Q-8}
\end{equation}
This is impossible on the K3 surface.

It remains to check $\nu=0$ in this weighted format.
In the $S$--axis chart, Morse preparation in $U$ gives
\begin{equation}
                 V^2+U^2+s^{6k+2}=0,\qquad k\ge1.
                                                        \label{eq:octic46:octic_degree12_complete46:45-D12Q-9}
\end{equation}
The exponent follows from characters; its coefficient cannot vanish
identically, since that would make the strict transform reducible
and nonnormal away from the exceptional divisor.  The quotient is a
dihedral quotient with scalar kernel of order six and projective
dihedral parameter $h=3k+1$.  Its minimal graph has a centre of square
$-2$, two $A_1$ arms, and a third arm
\[
                        [\underbrace{2,\ldots,2}_{k-1},4].
\]
The curve $E$ meets the final $(-4)$ curve.  The correction at this
end is again $1/3$, so $E$ is again $(-1)$ and must be contracted.
After its contraction the last $(-4)$ and the middle $(-4)$ of the
other axis become $(-3)$ and acquire an edge.  No remaining
exceptional curve is $(-1)$.

The actual order on the dihedral centre is $k>0$.
Indeed on its universal plane put
\[
       s=ab,\qquad U=\text{unit}\,(a^{2h}-b^{2h}),
                    \qquad V=\text{unit}\,(a^{2h}+b^{2h}).
\]
The form $s\,ds\wedge dU$ has total plane degree
$2h+4=6k+6$, giving order $k$ after the scalar quotient of order six.
The actual errors have the form
$s^4$ times regular multiples of $ds\wedge dU$ or $ds\wedge dV$,
or $s^5$ times a regular multiple of $dU\wedge dV$.
They have larger central order.  Parameter Morse changes may add
functions of $s$ to $U,V$, but their differentials wedge to zero
with $ds$, so the same estimate remains valid.
This positive centre survives the contraction of $E$ and lies on
the minimal resolution.  It contradicts the all-zero hypothesis.

\paragraph{A line and a three-line orbit}

Finally suppose $a=b=0$.  The $y$--chart has
\[
                       F=u^3v+c v^4+\varepsilon Su+\cdots.
\]
Its residue has character three, so $v\operatorname{Res}(dS\wedge
du\wedge dv/F)$ is invariant and regular canonical.
With ordinary weights its canonical threshold is four.
The only invariant equation monomials of degree below four are
$v^2$ and $Su$; the former is absent.  The weighted canonical-form
inequality therefore forces $\varepsilon\ne0$.

Equivariant parameter Morse preparation in $S,u$, keeping $v$ fixed,
gives $S'u'+c v^4=0$.  The quotient is
$\frac1{24}(1,11)$ with chain
\[
                              [3,2,2,2,2,3].
\]
The original $S,u$ remain in the ideal $(S',u')$: their characters
five and one cannot be supplied by a pure power of $v$, of character
three.  On the universal plane, $S'=a^4$, $u'=b^4$, $v=ab$.
Formula~\textup{\eqref{eq:octic46:octic_degree12_complete46:45-D12Q-4}} therefore has coefficient of ordinary
plane degree at least six.

The four middle $(-2)$ rays are
\[
            \frac1{24}(3,9),\quad\frac1{24}(5,7),\quad
            \frac1{24}(7,5),\quad\frac1{24}(9,3).
\]
A plane coefficient monomial of canonical order zero on one of them
has degree at most four; the only possibilities are among
$a^4,ab,b^4$.  None can occur in the actual form.  All four middle
orders are thus positive.  The old components have order zero and
cannot be $(-1)$, while every new quotient curve has square at most
$-2$.  The positive $A_4$ survives on the minimal resolution, giving
the final contradiction.

\paragraph{Conclusion.}
The four cases use twenty-two contraction matrices: four of size
twenty-one, nine of size twenty-three and nine of size twenty-six.
Their graphs, canonical-root incidences and determinants were specified
above.  The weighted two-component case is evaluated after contracting
its extraneous $(-1)$-curve, as required by the actual-form calculation.

The four quartic cases above exhaust \textup{\eqref{eq:octic46:octic_degree12_complete46:45-D12Q-2}}.
Each contradicts either the actual orders on the minimal resolution
or the K3 intersection lattice.  Thus neither alternative in
\textup{\eqref{eq:octic46:octic_degree12_complete46:45-D12Q-1}} occurs, completing the degree--twelve cusp exclusion.


\subsection{C2-02350: exclusion of the state $(14,1)$}
\label{sec:octic-degree14-delta1-node45}

\begin{proposition}
\label{prop:octic-degree14-delta1-node45}
In the sole degree-fourteen conductor state with generic delta one,
the generic transverse singularity cannot be nodal.  This conclusion
does not assume that the second polar jet is nonzero.
\end{proposition}

\begin{proof}
The actual conductor has marks $(1/5,5/6)$ at $Q_5,Q_{12}$.
Its coarse curve is smooth $\mathbb P^1$.
Generic delta one and purity make the full image conductor reduced:
the full conductor quotient is Cohen--Macaulay of dimension one
and generically reduced, hence reduced everywhere.
The same assertion applies on the ordinary canonical covers.
We will use the exact coherent conductor rings at the two sources.

At $Q_5$ the coordinate lower orders are $(2,3,4)$.
Coarse smoothness forces $x,y$ to have the exact orders two and
three, since $xy$ must supply the invariant uniformizer of order five.
The lifted branch is the plane cusp, and every possible $z$ term
belongs to its semigroup $\langle2,3\rangle$.
An equivariant analytic coordinate change therefore gives the
exact ideal
\begin{equation}
              I=(z,g),\qquad g=y^2-x^3,
              \qquad \operatorname{char}(x,y,z)=(1,4,2).
\label{eq:octic46:octic_degree14_node46:45-V1}
\end{equation}
This is the full conductor ideal, not just its radical support.

The equation belongs to $I^2$.  For example, its constant term in
$z$ is divisible by $g^2$, and its coefficient of $z$ is divisible
by $g$, by generic double vanishing along the reduced cusp.
The form $z\operatorname{Res}(dx\wedge dy\wedge dz/f)$ is invariant
and regular on the normalization, since $z$ lies in the full conductor.
If the $z^2$ coefficient of $f$ vanished, every term of $f\in I^2$
would have ordinary order at least four.
Indeed the cross term $zg$ requires a coefficient of character four,
of order at least one, and $g^2$ requires a coefficient of character
three, of order at least two.
The uniform canonical-form inequality with weights $(1,1,1)$ has
threshold $1+1+1+1=4$, and gives a contradiction.
Hence the quadratic $z^2$ coefficient is nonzero.

Equivariant Morse splitting, now justified independently of the
second polar jet, yields
\begin{equation}
                f=z'^2+g^2h(x,y),\qquad
                \operatorname{char}h=3.
\label{eq:octic46:octic_degree14_node46:45-V2}
\end{equation}
The exponent is exactly two under the nodal hypothesis.
The actual finite normalization is $V^2+h=0$, with $V=z'/g$
of character four.  Its reduced full conductor has coherent ring
\[
       \mathbb C[[u^2,u^3,V]]/(V^2+h(u)),
       \qquad x=u^2,\ y=u^3,\quad\operatorname{char}u=3.
\]
The invariant base is $\mathbb C[[s]]$, $s=u^5=xy$.
An invariant coefficient multiplying $V$ requires
$3n+4=0\pmod5$ in the cusp semigroup.
Its least exponent is two, and all exponents are $2+5j$.
Thus $q=xV$ is the primitive generator of the full trace-free module.
For the square coefficient, $3n=3\pmod5$ requires $n=1\pmod5$;
one is missing from the cusp semigroup, so the least exponent is six.
Consequently
\begin{equation}
                         q^2=-s^2H(s)
\label{eq:octic46:octic_degree14_node46:45-V3}
\end{equation}
and $Q_5$ contributes at least two discriminant units.

At $Q_{12}$ the actual trace bound and coarse smoothness force the
two transverse conjugate branches with full coherent ring
\begin{equation}
       A=\mathbb C[[y,z]]/(z^2-y^2),\qquad
       \operatorname{char}(y,z)=(11,5)\pmod {12}.
\label{eq:octic46:octic_degree14_node46:45-V4}
\end{equation}
This is the same exact reduced-curve straightening as in
\eqref{eq:octic46:octic_degree14_delta2_46:45-U7}.  It uses the degree-nine trace gap if $z$ rose from order
one to seven, and the invariant uniformizer of order six to force
$y$ of order one.  It does not use the generic transverse delta
beyond its already stated value one.

The coherent node \eqref{eq:octic46:octic_degree14_node46:45-V4} has dualizing character six, whereas the
surface residue has character nine.  Full conductor duality therefore
makes the trace-free cover module free over $A$, on a generator $t$
of character $6-9=9$.  Its square $t^2$ has character six.
Taking full $\mu_{12}$ invariants is essential here.
Put
\[
                          s=y^5z,\qquad q=y^2zt.
\]
The complete character-three coefficient module on $A$ is
$(y^2z)\mathbb C[[s]]$: multiplication by $s$ sends $y^2z$ to
$y^9$, and then repeats these two coherent normal forms.
Thus $q$ is a primitive trace-free generator.
Every character-six function on $A$ is $y^6H(s)$, so
\begin{equation}
                          q^2=s^2H(s).
\label{eq:octic46:octic_degree14_node46:45-V5}
\end{equation}
This source also contributes at least two discriminant units.
Passing only to one normalized branch before taking these invariants
would lose the nodal gluing condition and give an incorrect smaller
bound; \eqref{eq:octic46:octic_degree14_node46:45-V5} is a statement about the actual coherent node.

The reduced conductor cover has generic separable degree two under
the nodal hypothesis.  Lemma~\ref{lem:reduced-conductor-double-cover45}
bounds its total discriminant degree by two, whereas \eqref{eq:octic46:octic_degree14_node46:45-V3}--\eqref{eq:octic46:octic_degree14_node46:45-V5}
give at least four at two distinct points.
This contradiction proves the proposition.
\end{proof}

The surviving degree-fourteen delta-one alternative is therefore
cuspidal.  Its zero-second-jet case is also impossible, by the
\emph{coherent} degree-eleven bound, not by a gap in the normalized
weighted-line section ring.
If $\gamma|_\Gamma=0$, the third jet $\eta=\delta^3a|_\Gamma$
is intrinsic, because every change-of-frame error contains a lower
vanishing jet.  At the fifth cusp an actual character-three
coefficient has order at least six: the potential normalized order
one is absent from the semigroup $\langle2,3\rangle$.
At the proper-inertia-six source its order is at least one.
Consequently its actual coherent restriction line has degree at most
\[
              \frac{11\cdot14}{420}-\frac65-\frac16=-1,
\]
so $\eta=0$.
At the generic transverse cusp $z^2-x^3$, the vanishing of the second
jet first gives $\delta z|_\Gamma=0$; the third then gives
$-6(\delta x|_\Gamma)^3=0$.
Thus both normal components lie in $(x,z)$.
No vanishing of the induced tangential derivation is asserted or
needed.  On $x=t^2,z=t^3$, the first polar numerator has order at
least five, so the polar length is at least five.
This spends at least seventy units; length five leaves the forbidden
positive residual mass two, and a larger length exceeds seventy-two.
Thus any surviving realization of this ledger state has generic
cusp, nonzero second jet, and exact first-polar length three.


\paragraph{The cuspidal degree-fourteen alternative.}
\label{sec:octic-degree14-delta1-cusp45}

\begin{theorem}
\label{thm:octic-degree14-delta1-cusp45}
The sole degree-fourteen conductor state with generic delta one
cannot occur in a volume-$1/420$ realization of \texttt{C2-02350}.
\end{theorem}

The nodal type has already been excluded.  The coherent third-jet
argument in the preceding proposition forces the remaining cuspidal
type to have nonzero second jet and exact first-polar length three.
Let $R$ be its reduced conductor on the normalization.  Then
\begin{equation}
 C_N=2R,\qquad B_N=R+B',\qquad
                  420HR=14,\qquad420HB'=30.
\label{eq:octic46:octic_degree14_cusp46:45-W1}
\end{equation}
The actual second jet is a unit at $Q_5$ and has exact cover order
two at the proper-inertia-six source $Q_{12}$.

\paragraph{The complete fifth-source normalization is forced}

The exact fifth-source conductor cusp is $g=y^2-x^3$ in characters
$(x,y,z)=(1,4,2)$.  Equivariant Morse splitting and the generic cusp
type give
\[
             f=z^2+g^3h(x,y),\qquad\operatorname{char}h=0.
\]
Its normalization is
\begin{equation}
             V^2+gh=0,\qquad V=z/g,
             \operatorname{char}(x,y,V)=(1,4,4).
\label{eq:octic46:octic_degree14_cusp46:45-W2}
\end{equation}
This is the actual finite normalization, not just its generic model:
it is a hypersurface, is normal away from the displayed conductor,
and is smooth at its generic points because the transverse type is
a cusp.  Serre's criterion therefore applies.  The same argument is
used for the complete intersection below.
The residue $\Omega=dx\wedge dy/(2V)$ has character one.

We claim that $h$ is a unit.  If $h(0)=0$ and $[xy]h=0$,
the invariant form $y\Omega$ and weights $(1,1,3)$ have threshold
six, while the equation has weighted order at least six.
This contradicts rationality.  Thus a nonunit $h$ must have
$[xy]h\ne0$.
The weights $(x,y,V)=(4,6,11)$ then give the fixed isolated
initial equation
\begin{equation}
                         V^2+xy^3-x^4y.
\label{eq:octic46:octic_degree14_cusp46:45-W3}
\end{equation}
All subsequent terms have strictly larger weight.  Its weighted
curve is rational: since $x,y$ have even weights and $V$ odd
weight, its degree-zero function field is generated by $y^2/x^3$.
The source group is the radial $\mu_5$.

The $x$-axis, $y$-axis and remaining $V=0$ point give respectively
the quotients $\frac14(1,1)$, $\frac16(1,5)$ and
$\frac12(1,1)$.
The main raw square and corrections are
\[
 -\frac{5\cdot22}{4\cdot6\cdot11}=-\frac5{12},
 \qquad\frac14,\frac56,\frac12.
\]
Thus its final square is $-2$, with arms $[-4],A_5,A_1$.
Quasismoothness of \eqref{eq:octic46:octic_degree14_cusp46:45-W3} away from the origin makes this the
complete graph for all higher terms in this branch.

The leading actual form is a unit multiple of $dx\wedge dy$.
Its total weight is ten, and the radial order is five, so the main
canonical coefficient is exactly one.
Both adjacent $(-2)$ arms must be positive by adjunction.
Their coefficients cannot equal that of their parent, hence are at
least two.  Adjunction on the main $(-2)$-curve would then require
$2\geq2+2$, a contradiction.  This proves that $h$ is a unit.

After absorbing this unit, \eqref{eq:octic46:octic_degree14_cusp46:45-W2} is an $A_2$ cover with source
action giving $\frac1{15}(1,11)$ downstairs.
Its graph is
\begin{equation}
                            [2,2,3,2,2].
\label{eq:octic46:octic_degree14_cusp46:45-W4}
\end{equation}
All five actual-form orders are zero.
For a direct check, uniformize the $A_2$ by $x=st$,
$y\sim s^3+t^3$, $V\sim s^3-t^3$.
The cyclic rays have numerators
\[
                 (1,11),(2,7),(3,3),(7,2),(11,1)
\]
over fifteen, and the order of $V\,ds\wedge dt$ is
$\bigl(3\min(a,b)+a+b\bigr)/15-1=0$ on each.
The actual form differs by a unit: its coefficient relative to
$\Omega$ is $2V\delta_z+3gh\delta g+g^2\delta h$, and
$g=-V^2/h$ makes the other terms divisible by $V^2$.
The coefficient-one curve $R$ meets the unique $(-3)$ curve once;
the four $(-2)$ curves are protected and the $(-3)$ canonical
intersection is exhausted by $R$.

\paragraph{The complete twelfth-source graph is forced as well}

Straighten its actual conductor to $(x,g_0)$ with $g_0=z^2-y^2$,
in source characters $(1,11,5)$.
The exact second-jet order forces the coefficient of $g_0^2$ and
$\delta_z(0)$ to be nonzero.
The same square-lifting argument as for \eqref{eq:octic46:octic_degree14_delta2_46:45-U9}, now only through
the quadratic normal term, gives
\begin{equation}
          f=G^2+x^3h,\qquad G|_{x=0}=g_0\cdot\text{unit},
             \operatorname{char}G=10,\quad\operatorname{char}h=5.
\label{eq:octic46:octic_degree14_cusp46:45-W5}
\end{equation}
The finite birational normalization, with $V=G/x$ of character nine,
is the complete intersection
\begin{equation}
                         xV=G,\qquad V^2+xh=0.
\label{eq:octic46:octic_degree14_cusp46:45-W6}
\end{equation}
Its normality follows from $S_2$ and regularity in codimension one:
for $x\ne0$ it is the original normal locus, and at a generic
point over $x=0$ it is the smooth normalization of the transverse cusp.
The quotient by $x$ has dimension one, so there is no extra component
supported there.

Since $xV=G$ and $V^2=-xh$, the conductor contains $x$.
In particular $x^3\operatorname{Res}(dx\wedge dy\wedge dz/f)$ is
an invariant regular normalization form.
Weights $(1,2,2)$ and threshold eight force $[x^3z]f\ne0$.
This term cannot come from $G^2$, so $h_z(0)\ne0$.
After scaling, the ordinary tangent cone of \eqref{eq:octic46:octic_degree14_cusp46:45-W6} is
\begin{equation}
                      xV=z^2-y^2,\qquad V^2+xz=0.
\label{eq:octic46:octic_degree14_cusp46:45-W7}
\end{equation}
It is a smooth complete intersection of two quadrics in $\mathbb P^3$,
an elliptic quartic.

The ordinary blowup is therefore smooth along that elliptic curve
before the source quotient.  The radial subgroup has order two;
the remaining $\mu_6$ acts on the elliptic curve with quotient
$\mathbb P^1$ and stabilizer orders $6,2,3$.
The full quotient graph is a main curve $E$ of square $-3$, with
arms $A_5,A_1$, and one further $(-3)$ curve $F$.
Indeed the raw main square is $-4/3$ and the three corrections are
$5/6,1/2,1/3$, giving $-3$.

Here are the source charts, including the actual form.
At the $x$-axis the residual quotient is $\frac16(1,5)$ and
$dx\wedge dy$ becomes a unit multiple of $dS\wedge du$.
At the three nonzero $y=0$ points it gives $\frac12(1,1)$,
again with a unit form.
At the two points $x=V=0$, $z=\pm y$, the stabilizer-three chart
has coordinates $(S,w)$ with action $\frac13(1,1)$.
The second equation of \eqref{eq:octic46:octic_degree14_cusp46:45-W6} makes $x/y$ a unit times $w^2$,
so the actual form is a unit times $w\,dS\wedge dw$.
Its order on the single $(-3)$ curve is zero.
Its axes $S=0$ and $w=0$ are respectively the main curve and the
strict conductor; they meet that $(-3)$ curve at different points.
Thus $R$ meets $F$ once and is disjoint from $E$.
All eight exceptional curves have zero actual-form order, and the
$A_5+A_1$ are protected.

\paragraph{Thirteen fixed zero-order classes and a degree-five budget}

The fifth and twelfth graphs have thirteen zero-order curves.
The curve $R$ of coefficient one exhausts the canonical intersection
of the fifth $(-3)$ centre and of $F$.
The other $(-3)$ curve $E$ at $Q_{12}$ needs one additional
coefficient-one intrinsic polar curve, denoted $L$.
It cannot be $R$, by the incidence just established.
No positive polar component can meet a protected root.
Consequently every component of $B'$ avoids $Q_5$, and every
component of coefficient at least two also avoids $Q_{12}$.
The actual source-index theorem therefore gives
\begin{equation}
                         420H B_i\in5\mathbb Z_{>0}
                     \quad\text{for every component of }B'.
\label{eq:octic46:octic_degree14_cusp46:45-W8}
\end{equation}
Every positive exceptional coefficient at $Q_{14}$ is at most five.

After contracting the canonical forest, the intersection matrix of
the thirteen zero-order curves is their disjoint source matrix
plus $rr^t+ee^t$, where $r$ marks the fifth centre and $F$, and
$e$ marks $E$.
Children of the coefficient-one roots cannot meet these curves,
whose canonical intersections are at most one; contracting those
children first therefore gives exactly this update.
The matrix has determinant $48$, rank thirteen and signature $(1,12)$.
Thus eight further protected $Q_{14}$ roots would already contradict
the K3 Picard-rank bound.  The final matrices and their polarizations are computed below.

\paragraph{Only two fourteenth-source allocations survive the local tests}

After Corollary~\ref{cor:octic-q14-J}, the remaining rational
$Q_{14}$ formats and the external-contribution bound six give
the same positive-core test as for degree twenty-one.
Positive main curves in $G,K,L$ give no survivor;
the only $H$ survivor requires a coefficient-six polar curve, which
consumes all thirty units and leaves no coefficient-one curve through
$E$ at $Q_{12}$.
In the degree-thirty-eight and degree-forty-eight formats, a unit
$\delta_z$ protects respectively $A_8,A_{10}$, while a zero
$\delta_z$ makes the long arm positive, contradicting the coefficient
bound five.  These entire formats are excluded.

In format $G$, write the leading form as
$A z^2dy\wedge dz+B dx\wedge dy$.
If all three arm units are nonzero there are eight protected roots.
If $B=0$, $A\ne0$, its $A_4$ has coefficients $(1,2,3,4)$ and
forces a coefficient-five polar curve of degree five with mark
$5/14$.  It avoids $Q_5,Q_{12}$, so its degree is congruent to
$10$ modulo $35$, a contradiction.
The only remaining $G$ case is $B\ne0$ with one $A_2$ unit
cancelled.  Its other $A_2$ and its $A_4$ are protected.
At the cancelled point the form coefficient restricts along the main
curve to a nonzero cubic in the quotient coordinate.
Its two positive orders are bounded above by $(1,2)$ and cannot
be equal, so they are exactly $(1,2)$.
The free end needs one coefficient-three polar curve $B_3$;
the main $(-4)$ needs one coefficient-one curve.
The outer curvette of $B_3$ has a smooth $z$ parameter and
$H$-mark $13/14$.
It cannot meet another $Q_{14}$ component, or $Q_5,Q_{12}$.
Thus its scaled degree is $5$ modulo $35$, and the budget forces
it to be five.
The coefficient-one curve meets the main weighted curve once,
with $H$-mark $11/14$, and no other $Q_{14}$ component.
Its degree is therefore $15$ modulo $35$.
The remaining budget forces this curve to be precisely $L$, of
degree fifteen, meeting $E$ at $Q_{12}$.
These two components exhaust $B'$.

In format $K$, a unit $\delta_z$ protects $A_6$ and the first
$(-2)$ of the $[2,4]$ arm; both curves of that arm have order zero.
The third arm is $[3,2,2]$.
Its $(-3)$ always has order zero.  In its cyclic chart the actual
form coefficient is $c v^3+S\cdot(\text{higher coefficient})$,
$c\ne0$.
The two last curves are either both protected, or have orders
$(1,2)$: a nonzero $S$ term gives zero orders, while in its absence
positivity and the upper bounds $(1,2)$ force this pair.
The protected alternative gives nine roots and is impossible.
The positive alternative again forces $B_3$ of coefficient three,
degree five and mark $13/14$ at the free end.

The main $(-3)$ in $K$ requires a coefficient-one curve $L$.
Its generic main-curve contact has mark $1/2$.
To have positive degree at most fifteen, it must acquire a
seven-primary marking at the same source, since it avoids $Q_5$
and all other indices divide twelve.
The only available further contact is the zero $(-4)$ in the
$[2,4]$ arm, whose outer curvette has mark $9/14$.
Its total intersection with that curve is at most two.
One intersection would give degree $25$ modulo $35$, exceeding
the budget.  Two give degree $15$ modulo $35$.
Hence $L$ has degree fifteen, meets that $(-4)$ with total
intersection multiplicity two, and also meets $E$ at $Q_{12}$.
The two intersections may occur at one tangential branch; both
the degree congruence and the matrix update use their total
multiplicity and allow this possibility.
Again $B_3$ and $L$ exhaust $B'$.

\paragraph{The two final allocations violate the K3 intersection form}

In both allocations the full canonical forest is determined.
The coefficient-one curves $R,L$ are isolated final $(-1)$ roots;
the remaining component is the chain with coefficients $1,2,3$
formed by the positive $A_2$ and $B_3$.
There are no additional positive resolution curves: a minimal
resolution curve of square at most $-2$ has a larger child, and
following such children must end at an intrinsic component.
The three displayed intrinsic components exhaust the budget;
their coefficients and already fixed parents leave no other branch.

\paragraph{The polarization and the two remaining configurations.}
\label{sec:degree14-hodge45}

This calculation applies after the geometric reduction of the sole
degree-fourteen cusp to the last normal-Q14 format $G$.
It uses the actual polarization as well as the exceptional curves.
We keep track of the polarization during the contractions.

There are twenty zero-order curves on the normalization resolution:
the Q5 chain $[2,2,3,2,2]$, the Q12 star with central $(-3)$-curve
$E$, arms $A_5,A_1$, and one outer $(-3)$-curve $F$, and the Q14
central $(-4)$-curve $G$ with protected arms $A_2,A_4$.
Write $T$ for the central $(-3)$-curve at Q5.
The positive canonical divisor consists of the following components:
\[
 R,\quad L,\quad A_1+2A_2+3B,
 \qquad
 HR=\frac1{30},\quad HL=\frac1{28},\quad HB=\frac1{84}.
\]
Here $R$ and $L$ are disjoint final $(-1)$-curves.
The other three curves form a chain with squares $(-2,-2,-1)$.
The curve $R$ meets $T,F$ once, $L$ meets $E,G$ once, and the
three-curve chain meets $G$ at its coefficient-one end.
There are no other positive components or intersections.

After contracting this canonical forest to the K3 surface $S$,
the twenty images all have square $-2$.
Their intersection matrix $M$ has determinant $-924$.
Eliminating the negative-definite arms gives the four-by-four Schur
complement, in the order $(T,F,E,G)$,
\begin{equation}
 M_0=
 \begin{pmatrix}
 -2/3&1&0&0\\
 1&-2&1&0\\
 0&1&-2/3&1\\
 0&0&1&-8/15
 \end{pmatrix}.
 \label{eq:octic46:octic_degree14_hodge46:45-H1}
\end{equation}

Let $H_S$ be the rational divisor obtained by pushing the actual
polarization forward to $S$.  The blowdown intersection formula gives
\begin{equation}
 H_S^2=\frac2{105}+\frac1{900}+\frac1{784}
                  +\frac3{7056}=\frac{241}{11025}.
 \label{eq:octic46:octic_degree14_hodge46:45-H2}
\end{equation}
The factor three in the last term is essential: each of the successive
contractions of $B,A_2,A_1$ has current polarization degree $1/84$.
The intersections with the four central curves are
\begin{equation}
 h_0=\left(\frac1{30},\frac1{30},\frac1{28},\frac1{21}\right)^{\!t}.
 \label{eq:octic46:octic_degree14_hodge46:45-H3}
\end{equation}
All other sixteen intersections are zero.  In particular the last
entry is $1/28+1/84$, since both $L$ and the three-curve chain change
the intersection with $G$.

Let $Q$ be the unique rational combination of the twenty curves having
the same intersections with them as $H_S$.  Exact elimination yields
\begin{equation}
 Q^2=h_0^tM_0^{-1}h_0=\frac{1627}{242550}>0,
 \qquad
 (H_S-Q)^2=\frac1{66}>0,
 \qquad Q(H_S-Q)=0.
 \label{eq:octic46:octic_degree14_hodge46:45-H4}
\end{equation}
Thus two orthogonal rational divisor classes have positive square on
the smooth projective surface $S$, contrary to the Hodge index theorem.
The last $G$ configuration is impossible.

\paragraph{The companion $K$ configuration.}
For completeness, the corresponding reduction in the Klein-quartic
format $K$ leaves twenty-three zero-order curves.
Its Q14 part is the central $(-3)$-curve with arms
$[2,4]$, $A_6$, and a single $(-3)$-curve; the positive terminal
$A_2$ has already been removed.
Contracting $R$ adds the outer product of the incidence vector
$e_T+e_F$.  Contracting $L$ adds the outer product of
$e_E+e_{\mathrm{main},14}+2e_{-4,14}$.
The terminal positive chain adds one to the remaining Q14
$(-3)$-curve's square.  The resulting matrix has determinant
$1656$, signature $(1,22)$ and rank twenty-three.  These twenty-three independent
divisor classes cannot lie in the Neron--Severi group of a complex
projective K3 surface, whose rank is at most twenty.


Both final allocations are impossible.  This proves
Theorem~\ref{thm:octic-degree14-delta1-cusp45}.


\section{Conductor filtrations in degree seven}
\label{sec:octic-higher-delta-revised}

It remains to treat $(7,2),(7,3),(7,4)$, whose reduced marked
curve is $\mathbb P(5,12)$. We retain the full conductor and its
actual line-bundle characters. After the ribbon and the preliminary
delta-three reductions, cohomology of the conductor filtration forces
a degree-three section on the normalization while forbidding its
square. The delta-four double branch is excluded by the degrees of
the successive multiplication maps. These arguments use the saturated
layers without requiring a global splitting of the conductor algebra.

\paragraph{Duality on the full conductor.}

\begin{lemma}\label{lem:octic-full-conductor-duality}
Let $\mathcal B=C_D$ be the full image conductor on the canonical-index
stack and put $Q=\nu_*O_N/O_D$.  Then
\[
 Q\simeq\omega_{\mathcal B}(-9H),\qquad
 H^1(\mathcal B,O_{\mathcal B})=0,
 \qquad H^i(\mathcal B,O_{\mathcal B}(9H))=0\quad(i=0,1).
\]
\end{lemma}
\begin{proof}
Finite duality applied to the normalization sequence gives
$\mathcal Ext_D^1(Q,\omega_D)=\omega_D|_{\mathcal B}$.
Cohen--Macaulay biduality and $\omega_D=O_D(9H)$ give the displayed
isomorphism.  The first vanishing is Proposition~\ref{prop:octic-full-conductor46};
curve duality translates the acyclicity of $Q$ from
Theorem~\ref{thm:octic-K3-rational46} into the last two vanishings.
The coarse pushforward is exact, so these statements include the
entire nonreduced conductor.
\end{proof}

\begin{lemma}[A square obstruction on the normalization]
\label{lem:octic-conductor-square-obstruction}
For the full image conductor $\mathcal B=C_D$ of the octic, one cannot
have
\[
 h^1(\mathcal B,O_{\mathcal B}(3H))=0,
 \qquad h^1(\mathcal B,O_{\mathcal B}(6H))>0.
\]
\end{lemma}
\begin{proof}
Put $Q=\nu_*O_N/O_D$. Lemma~\ref{lem:octic-full-conductor-duality}
and curve duality give
\[
 H^0(Q(jH))\simeq
 H^1(\mathcal B,O_{\mathcal B}((9-j)H))^*.
\]
The complete restriction formula \eqref{eq:complete-restriction46}
and \eqref{eq:octic-small-P46} give $T_3=T_6=0$.
Moreover $H^1(D,O_D(3H))=0$: in the octic restriction sequence,
positive canonical cohomology vanishes and
\[
 H^2(Z,O_Z(-5K_Z))^*=H^1(Z,O_Z(6K_Z))=0.
\]

Twist the normalization sequence by the actual line bundle $jH$ on
the canonical-index stack and take coarse pushforward. Both operations
are exact; the normalized character sheaf is $O_N(jH)$, since it is
$S_2$ and agrees with that reflexive sheaf in codimension one.
Thus the sequences of complete sections give
\[
 S_3\simeq H^0(Q(3H)),\qquad S_6\hookrightarrow H^0(Q(6H)).
\]
The asserted cohomology would imply $S_3\ne0$ and $S_6=0$.
Since $N$ is integral, the square of a nonzero section of $O_N(3H)$
is a nonzero section of $O_N(6H)$, a contradiction.
\end{proof}

We will apply this lemma directly to filtrations of $B=O_{\mathcal B}$.
Every twist is by the actual $H$; on a layer supported on
$\mathcal G=\mathbb P(5,12)$ it is the usual $O_{\mathcal G}(1)$.
For $G=\mathbb C[X_5,Y_{12}]$ we use
$h^0(O_{\mathcal G}(d))=\dim G_d$ and
$h^1(O_{\mathcal G}(d))=\dim G_{-17-d}$.
In particular $O_{\mathcal G}(3)$ and $O_{\mathcal G}(6)$ are acyclic.

\subsection{C2-02350: the full conductor algebra of the state $(7,2)$}
\label{sec:octic-degree7-delta2-ribbon45}

\begin{proposition}
\label{prop:octic-degree7-delta2-ribbon45}
Suppose the sole conductor image has scaled degree seven and generic
delta two.  Its full image conductor is a ribbon $\mathcal B$ over
$\mathcal G=\mathbb P(5,12)$, and its nilpotent line is necessarily
\begin{equation}
               \mathcal L=\mathcal O_{\mathcal G}(17).
\label{eq:octic46:octic_delta2_ribbon46:45-RB1}
\end{equation}
Its nilpotent direction is $y$ at the fifth source and $x$ at the
twelfth source.  If $\mathcal J$ is the restriction to $\mathcal G$
of the trace-free line of the full normalization-conductor double
algebra, then
\begin{equation}
                   \mathcal J=\mathcal O_{\mathcal G}(-43).
\label{eq:octic46:octic_delta2_ribbon46:45-RB2}
\end{equation}
For transverse type $A_3$ its reduced trace-zero square is a nonzero
multiple of $X_5^{10}Y_{12}^3$.  For transverse type $A_4$ its square
lands in the nilpotent line and is a nonzero multiple of
$X_5^{11}Y_{12}^4$.
\end{proposition}

\begin{proof}
The degree-seven source and trace argument identifies
$(\mathcal G,H)=(\mathbb P(5,12),\mathcal O(1))$, with
\[
 G=\mathbb C[X_5,Y_{12}],\qquad \omega_{\mathcal G}=\mathcal O(-17).
\]
Lemma~\ref{lem:octic-degree7-stack46} proves this statement using only
the reduced component, so it applies to the present nonreduced conductor.
The actual lifts are the $x$-axis at $Q_5$ and the $y$-axis at $Q_{12}$;
additional Gorenstein contacts create no further stack structure.

Write $B=\mathcal O_{\mathcal B}$.  It is Cohen--Macaulay and has
generic length two over its reduced curve.  Its nilradical has square
zero: the square vanishes generically, and a Cohen--Macaulay curve has
no nonzero subsheaf supported at finitely many points.  Its nilradical
is torsion-free of rank one on the smooth stack $\mathcal G$, hence is
a line $\mathcal L=\mathcal O(m)$.  Locally $B$ is the ordinary
ribbon algebra $A[\epsilon]/(\epsilon^2)$ over a smooth curve ring
$A$; in particular it is Gorenstein.  A global splitting is not
assumed.

Let $Q=\nu_*\mathcal O_N/\mathcal O_D$.  Full conductor duality gives
\begin{equation}
                        Q=\omega_B(-9H).
\label{eq:octic46:octic_delta2_ribbon46:45-RB3}
\end{equation}
This is Lemma~\ref{lem:octic-full-conductor-duality}.
Thus $Q$ is invertible over $B$, and ribbon adjunction yields
\begin{equation}
              Q|_{\mathcal G}=\mathcal J
              =\omega_{\mathcal G}\otimes\mathcal L^{-1}(-9H)
              =\mathcal O(-26-m).
\label{eq:octic46:octic_delta2_ribbon46:45-RB4}
\end{equation}

The same lemma gives $H^1(B)=0$ and the acyclicity of $B(9H)$.
Since $\mathcal O_{\mathcal G}(9)$ is acyclic, the ribbon sequence
implies that $\mathcal L(9)$ is acyclic.  The unit of $B$ makes
$H^0(B)\to H^0(\mathcal O_{\mathcal G})$ surjective, so $H^1(B)=0$
also gives $H^1(\mathcal L)=0$.
Writing $G_d=0$ for $d<0$, the exact necessary conditions are
\begin{equation}
              G_{m+9}=G_{-26-m}=G_{-17-m}=0.
\label{eq:octic46:octic_delta2_ribbon46:45-RB5}
\end{equation}

Here the signs of the local characters matter.  The fibre character
of the line $H$ is $+2$ at $Q_5$ and $+5$ at $Q_{12}$.
The opposite characters $-2j,-5j$ describe section numerators and
must not be substituted for these line-module characters.
At $Q_5$ the nilpotent fibre is generated by either $y$ or $z$,
of characters four or two.  At $Q_{12}$ it is generated by either
$x$ or $z$, of characters one or five.  The distinct characters make
these the only possibilities.  Thus the four congruence classes and
their surviving cohomology representatives are
\begin{equation}
\begin{array}{c|c|c|c}
 Q_{12}\text{ direction}&Q_5\text{ direction}&m\pmod {60}
                                      &\text{solutions of \eqref{eq:octic46:octic_delta2_ribbon46:45-RB5}}\\\hline
 x&y&17&17\\
 z&y&37&-23\\
 x&z&41&-19\\
 z&z&1&\text{none}.
\end{array}
\label{eq:octic46:octic_delta2_ribbon46:45-RB6}
\end{equation}
This is finite without a search assumption: every integer at least
44 belongs to the semigroup $\langle5,12\rangle$, so \eqref{eq:octic46:octic_delta2_ribbon46:45-RB5}
forces $-60\le m\le34$.
It is essential to retain the last vanishing in \eqref{eq:octic46:octic_delta2_ribbon46:45-RB5}.
The two conditions coming from $Q$ alone would leave $m=-59$ in
the last congruence class; $G_{42}\ne0$ violates $H^1(B)=0$ there.

Now let $C=\nu_*\mathcal O_{C_N}$.  The full conductor square is
\[
                         0\longrightarrow B\longrightarrow C
                              \longrightarrow Q\longrightarrow0.
\]
Since $Q$ is invertible, $C$ is locally free of rank two over $B$.
The characteristic-zero algebra trace splits the unit summand:
$C=B\oplus\mathcal N$, with $\mathcal N\simeq Q$ invertible and
$\operatorname{Tr}\mathcal N=0$.  The Cayley--Hamilton identity
then gives a multiplication map
\begin{equation}
                         \mu:\mathcal N^{\otimes2}\longrightarrow B.
\label{eq:octic46:octic_delta2_ribbon46:45-RB7}
\end{equation}
This statement holds over the nonreduced ribbon itself.

Generic delta two has precisely the transverse alternatives $A_3$
and $A_4$.  For $A_3$, after a separable splitting extension if needed,
the full conductor algebra is
$k[\epsilon]/(\epsilon^2)\oplus k[\epsilon]/(\epsilon^2)$ over its
diagonal ribbon.  The trace-zero square is a unit modulo $\epsilon$.
Thus \eqref{eq:octic46:octic_delta2_ribbon46:45-RB7} reduces to a nonzero map
$\mathcal J^2\to\mathcal O_{\mathcal G}$, and necessarily
\begin{equation}
                              G_{52+2m}\ne0.
\label{eq:octic46:octic_delta2_ribbon46:45-RB8}
\end{equation}
For $A_4$, the local full conductor algebra is
$k[t]/(t^4)$ over $k[t^2]/(t^4)$.  Its trace-zero generator has
square $t^2$, a nonzero nilpotent.  Consequently $\mu$ is zero
modulo $\mathcal L$ generically, hence everywhere on the integral
reduced curve.  It lands in $\mathcal L$ and factors through
$\mathcal J^2$, giving a nonzero map
$\mathcal J^2\to\mathcal L$.  Necessarily
\begin{equation}
                              G_{52+3m}\ne0.
\label{eq:octic46:octic_delta2_ribbon46:45-RB9}
\end{equation}
Neither argument uses a global splitting of $B$.

For $m=-23$ the two spaces in \eqref{eq:octic46:octic_delta2_ribbon46:45-RB8}--\eqref{eq:octic46:octic_delta2_ribbon46:45-RB9} are $G_6,G_{-17}$;
for $m=-19$ they are $G_{14},G_{-5}$.  All four vanish.
These two entire ribbon classes are therefore impossible, leaving
only $m=17$ and \eqref{eq:octic46:octic_delta2_ribbon46:45-RB1}--\eqref{eq:octic46:octic_delta2_ribbon46:45-RB2}.
Finally
\[
 G_{86}=\mathbb C\,X_5^{10}Y_{12}^{3},\qquad
 G_{103}=\mathbb C\,X_5^{11}Y_{12}^{4},
\]
which gives the last assertions of the proposition.
\end{proof}

\paragraph{The resulting algebra normal form.}
For the surviving $m=17$, the ribbon splits because
$H^1(T_{\mathcal G}\otimes\mathcal L)=H^1(\mathcal O(34))=0$.
Also $H^1(\mathcal L)=0$, so an invertible line on the ribbon with
restriction $\mathcal J$ has a unique isomorphism class of lifts.
In these splittings $\mu$ has components in $G_{86}\oplus G_{103}$.
In the $A_3$ case the nilpotent component can be removed by rescaling
the trace-free generator by a unit $1+s$, $s\in H^0(\mathcal L)$:
$G_{17}=\mathbb C X_5Y_{12}$ and
$(X_5Y_{12})(X_5^{10}Y_{12}^3)=X_5^{11}Y_{12}^4$.
Thus there is a single normalized multiplication in each of the
$A_3,A_4$ types at the level of this full conductor algebra.
This determines the full conductor algebra in the two transverse types.


\subsection{C2-02350: the twelfth and fifth sources in the state $(7,2)$}
\label{sec:degree7-delta2-q12-45}

The full ribbon theorem fixes the nilpotent line to $\mathcal O(17)$.
At Q12 its full ideal is $(z,x^2)$, in characters
$(x,y,z)=(1,11,5)$, and the trace-free line on the reduced conductor
is $\mathcal O(-43)$.  The following calculation does not assume a
polar alternative until the actual-form paragraph.

\begin{proposition}
\label{prop:degree7-delta2-q12-45}
For generic type $A_3$, the complete Q12 graph has a strict line of
square $-2$ with an $A_5$ arm, and a strict cubic of square $-3$
with arms $A_2,A_1$.  They meet opposite ends of
$[2,2,2,3,2,2,2]$.
For generic type $A_4$, replace this joining chain by
$[2,2,2,2,4,2,2,2,2]$.
If the second polar is nonzero, all actual-form orders are zero,
and the protected root ranks are fifteen and seventeen respectively.
If the derivation vanishes along the reduced conductor, these graphs
contain respectively a positive $A_9$ and a positive $A_{10}$.
\end{proposition}

\begin{proof}
The regular conductor numerators $x^3,z^3,y^2z$ force the
$x^3z,z^4,y^2z^2$ coefficients to be nonzero.  Write
\begin{equation}
 f=A(x,y)x^4+2B(x,y,z)x^2z+C(x,y,z)z^2,\qquad
 \operatorname{char}(A,B,C)=(4,1,10).
 \label{eq:octic46:octic_delta2_q12_46:45-D2T1}
\end{equation}
We choose $A$ independent of $z$.  The initial coefficients are
\begin{equation}
 B=\beta x+\cdots,\qquad C=c y^2+d z^2+\cdots,\qquad\beta cd\ne0.
 \label{eq:octic46:octic_delta2_q12_46:45-D2T2}
\end{equation}
The orders of $B(0,y,0)$ and $C(0,y,0)$ are at least eleven and
exactly two.  The corresponding square completion in $z$ is regular
and preserves $(z,x^2)$; perform it to arrange $B(0,y,0)=0$.
The finite determinantal normalization adjoins $t$ of character one:
\begin{equation}
 x^2t+Bx^2+Cz=0,\quad zt-Ax^2-Bz=0,\quad t^2-B^2+AC=0.
 \label{eq:octic46:octic_delta2_q12_46:45-D2T3}
\end{equation}
This algebra is Cohen--Macaulay, is normal off the conductor, and is
its smooth normalization at every generic conductor point.
It is therefore normal by $R_1$ and $S_2$.

For $A_3$, the full boundary multiplication is the nonzero monomial
$X^{10}Y^3$.  At Q12, $X$ has order one and $Y$ is a unit.
Thus \eqref{eq:octic46:octic_delta2_q12_46:45-D2T3} gives $A(0,y)=y^8\cdot\text{unit}$.
For $A_4$, the square lands in the nilpotent line and its coefficient
is $X^{11}Y^4$.  It gives
\begin{equation}
             A(0,y)=0,\qquad A_x(0,y)=y^9\cdot\text{unit}.
 \label{eq:octic46:octic_delta2_q12_46:45-D2T4}
\end{equation}
These are exact boundary orders, and cannot be raised by a local
specialization.

\paragraph{The normalized $y$-chart.}
Put $y=s,x=su,z=sv$ and $S=s^2$.
The strict equation is
$\widetilde A u^4+2\widetilde B u^2v+\widetilde C v^2$,
where $\widetilde A=A(su,s)$,
$\widetilde B=B(su,s,sv)/s$ and
$\widetilde C=C(su,s,sv)/s^2$.
Near the origin $\widetilde C$ is a unit.
Equivariant Morse preparation in $v$ gives
$v'^2+u^4h(S,u)=0$, whose actual normalization is
\begin{equation}
 T^2+h(S,u)=0,\qquad T=v'/u^2,\qquad
                      \mu_6(S,u,T)=(5,1,1).
 \label{eq:octic46:octic_delta2_q12_46:45-D2T5}
\end{equation}
The coefficient of $u^2$ in $h(0,u)$ is nonzero.
For $A_3$, the preceding boundary order supplies a nonzero $S^4$
term, so parameter Morse preparation in $u$ gives
$T^2+U^2+S^4=0$.
For $A_4$, it gives
$h(S,0)=0$ and $h_u(S,0)=S^5\cdot\text{unit}$.
The critical point has order five in $S$, and its critical value has
order ten with nonzero coefficient.
There is no independent pure-$S$ term which could cancel it,
because $h(S,0)$ is identically zero.
Thus the second equation becomes $T^2+U^2+S^{10}=0$.

Write $n=4$ or $10$.  On its universal plane,
$S=pq$, $U\sim p^n-q^n$, $T\sim p^n+q^n$.
The full quotient is
\begin{equation}
                       \frac1{6n}(1,5n-1).
 \label{eq:octic46:octic_delta2_q12_46:45-D2T6}
\end{equation}
These weights induce the stated $(5,1,1)$ action; their continued
fractions give the joining chains in the proposition.

\paragraph{The two attachments.}
The ordinary quartic is
$z(\beta'x^3+c y^2z+d z^3)$, with all three coefficients nonzero.
Its line and smooth cubic meet opposite ends of \eqref{eq:octic46:octic_delta2_q12_46:45-D2T6}, since
their strict transforms are the two branches of $S=0$.
Their local mutual intersection before resolving that chart is
$1/(6n)$, and the end correction is $(5n-1)/(6n)$ on each.
The fixed charts are $A_5$ on the line and $A_2,A_1$ on the cubic.
The total ordinary exceptional divisor meets the quotient line
and cubic in $-1/3$ and $-1$, respectively.
Their final squares are consequently
\begin{equation}
 -\frac13-\frac1{6n}-\frac{5n-1}{6n}-\frac56=-2,\qquad
 -1-\frac1{6n}-\frac{5n-1}{6n}-\frac23-\frac12=-3.
 \label{eq:octic46:octic_delta2_q12_46:45-D2T7}
\end{equation}
Every exceptional square is at most $-2$, so this is the minimal
resolution; there is no extraneous $(-1)$ component.

\paragraph{The actual form.}
If the second polar is nonzero, its exact order two at Q12 gives
$\delta_z(0)\ne0$.  The main components and fixed $A_5,A_2,A_1$
charts have order zero by ordinary blowup.
In the $y$-chart the exact actual form is
\begin{equation}
 \alpha=\frac A2du\wedge dS+\frac{S^2B}{2}dS\wedge dv
 +S^2C\,dv\wedge du+
 \frac{SC}{2}(v\,dS\wedge du-u\,dS\wedge dv).
 \label{eq:octic46:octic_delta2_q12_46:45-D2T8}
\end{equation}
Here $A$ is a unit and the lowest terms of $C$ are $uv,S^2$;
these letters now denote derivation coefficients.
Parameter Morse changes keep $S$ fixed.
On the universal plane the original $u,T$ have degree at least $n$.
The original $v$ has degree at least $3n$: it equals
$u^2\Phi(S,u,T)$ with $\Phi(S,0,0)=0$, by the earlier elimination
of $B(0,y,0)$.
Possible pure-$S$ translations of $u,T$ have order at least five
by characters, and therefore plane degree at least $n$.

The first term of \eqref{eq:octic46:octic_delta2_q12_46:45-D2T8} is a unit multiple of
$dU\wedge dS$, with leading numerator a nonzero combination of
$p^n,q^n$.
Every other term has larger plane degree.
On every ray in \eqref{eq:octic46:octic_delta2_q12_46:45-D2T6}, the leading numerator has canonical
order zero.  A higher-degree monomial has larger valuation than
its minimum-degree leading monomial, since both ray coordinates
are positive.  Thus all asserted zero orders hold for the actual
form with all higher terms present.
The protected ranks are
\begin{equation}
 (5+1+3)+3+2+1=15,\qquad (5+1+4)+4+2+1=17.
 \label{eq:octic46:octic_delta2_q12_46:45-D2T9}
\end{equation}

If instead $\delta$ vanishes along the conductor, every ambient
form summand has ordinary total weight at least four.
The old line therefore has positive order after the radial
involution.  Adjunction propagates positivity along its connected
$(-2)$ chain.  That chain comprises the $A_5$, the old line,
and the first three or four vertices in the joining chain.
It is a positive $A_9$ or $A_{10}$, respectively.
\end{proof}


\paragraph{The fifth source when the second polar is nonzero.}
\label{sec:octic-degree7-delta2-q5-nonzero45}

We apply Proposition~\ref{prop:octic-degree7-delta2-ribbon45} at Q5:
the nilradical is $L=O(17)$, the trace-free restriction is $J=O(-43)$,
and their multiplication is given by $X^{10}Y^3$ for $A_3$ and
$X^{11}Y^4$ for $A_4$, up to nonzero constants.

\begin{proposition}
\label{prop:octic-degree7-delta2-q5-nonzero45}
Suppose the second polar is nonzero in this state.  At Q5 the generic
$A_3$ case has precisely the normalization
\[
            (w^2+y^2+x^3=0)/\mu_5(1,4,4).
\]
Its minimal graph is $[2,2,3,2,2]$, with all actual orders zero,
and hence four protected root curves.
In the generic $A_4$ case the only local possibility not already
contradicting rationality or the canonical forest is
\[
            (w^2+y^2+x^8=0)/\mu_5(1,4,4).
\]
Its graph is $[2,2,2,4,2,2,2]$, all actual orders are zero,
and it supplies six protected root curves.
The assertions are invariant under equivariant analytic preparation;
the actual form is controlled as in the proof, not replaced by an
unrelated model form.
\end{proposition}

\begin{proof}
The nonzero degree-ten polar is a unit at Q5.  Equivariant Morse
splitting in $z$ gives $\delta_z(0)\ne0$ and
\begin{equation}
 f=z^2+y^4h(x,y),\qquad
 N^\#=(w^2+h(x,y)=0),\qquad w=z/y^2,
 \quad \mu_5(x,y,w)=(1,4,4).
                                                        \label{eq:octic46:octic_delta2_q5_46:45-D2Q5-1}
\end{equation}
Here $h$ has character three.  The full conductor ideal is $(z,y^2)$;
thus its normalized algebra modulo $y^2$ is exactly the rank-two
ribbon cover, with trace-free generator $w$.
The global boundary multiplication therefore gives
\begin{align}
 A_3:&\quad h(x,0)=x^3\cdot\text{unit},                 \label{eq:octic46:octic_delta2_q5_46:45-D2Q5-2}\\
 A_4:&\quad h=yh_1(x,y),\qquad
                  h_1(x,0)=x^4\cdot\text{unit}.       \label{eq:octic46:octic_delta2_q5_46:45-D2Q5-3}
\end{align}
These are exact first boundary orders, not generic-coefficient choices.
The cover in \eqref{eq:octic46:octic_delta2_q5_46:45-D2Q5-1} is normal: a repeated branch divisor away
from $y=0$ would produce an additional conductor component, and the
generic $A_3$ or $A_4$ condition excludes a repeated factor $y$.

The actual form on this normalization is
\begin{equation}
 \alpha=(\delta_z-2yw\delta_y)\,dx\wedge dy
       +y^2\delta_x\,dy\wedge dw-y^2\delta_y\,dx\wedge dw.
                                                        \label{eq:octic46:octic_delta2_q5_46:45-D2Q5-4}
\end{equation}
The source characters of $\delta_x,\delta_y,\delta_z$ are
$(4,2,0)$, and $\delta_z$ is a unit.  This formula will control
all possible higher terms.

\paragraph{The cases with a nonzero $y^2$ coefficient}
If $[y^2]h\ne0$, equivariant parameter Morse splitting in $y$,
keeping $x$ fixed, gives $w^2+y'^2+x^n=0$.
In \eqref{eq:octic46:octic_delta2_q5_46:45-D2Q5-2} one has $n=3$.  In \eqref{eq:octic46:octic_delta2_q5_46:45-D2Q5-3} one has $n=8$:
the critical value is a nonzero constant times $-x^8$, since the
linear $y$ coefficient has exact order four and no pure $x$ term
was present before completing the square.

The quotient is
\[
                 \frac1{5n}(1,4n-1).
\]
On the universal plane, $x=ab$, $y'=a^n-b^n$ and
$w=a^n+b^n$, up to nonzero constants.  The leading coefficient
of $dx\wedge dy$ is a nonzero multiple of $a^n+b^n$.
The continued fractions are
\[
 \frac{15}{11}=[2,2,3,2,2],\qquad
 \frac{40}{31}=[2,2,2,4,2,2,2].
\]
Direct ray valuation gives order zero on every curve in both chains.
The remaining terms in \eqref{eq:octic46:octic_delta2_q5_46:45-D2Q5-4} have strictly larger ordinary
degree on the universal plane.  They cannot cancel either of the
two displayed monomials identically on an exceptional curve.
The same holds for the nonconstant part of the unit $\delta_z$.
Thus these are the actual zero orders, and every $(-2)$ curve is
protected under the canonical contractions to the K3 model.

\paragraph{Generic $A_3$ with $[y^2]h=0$}
The coefficient of $x^3$ is nonzero by \eqref{eq:octic46:octic_delta2_q5_46:45-D2Q5-2}.
If $[xy^3]h\ne0$, use weights
\[
                       (x,y,w)=(6,4,9).
\]
The initial equation is $w^2+x^3+cxy^3$, $c\ne0$, and the
source $\mu_5$ is radial.  In the $y$--chart put
$y=s^4$, $x=s^6u$, $w=s^9v$ and $S=s^5$.
At either nonzero root of $u^2+c=0$, $v=0$, the remaining
stabilizer acts by $\frac12(1,1)$ on $(S,v)$.
The form $dx\wedge dy$ is a unit times
\[
                             S v\,dS\wedge dv.
\]
The main exceptional curve and the adjacent $A_1$ therefore both
have actual order one.  In these weights $\delta_x$ has weight
at least four and $\delta_y$ at least twelve.  All corrections
in \eqref{eq:octic46:octic_delta2_q5_46:45-D2Q5-4} have weighted degree at least twenty-five, whereas
$dx\wedge dy$ has degree ten.  After the radial quotient they
cannot cancel either order-one coefficient.  Adjacent equal
positive coefficients are impossible in a canonical blowup forest.

Suppose now that $[xy^3]h=0$ as well.  The residue of the normal
cover has character one, so $y\Omega$ is an invariant regular
canonical form.  With weights $(x,y,w)=(4,1,6)$ its canonical
threshold is twelve.  The only equation monomials below twelve
are $y^2,xy^3,y^7$.  The first two are absent, so $[y^7]h\ne0$.
Use weights $(14,6,21)$.  The entire initial equation is
\begin{equation}
                            w^2+x^3+c y^7,
                           \qquad c\ne0.               \label{eq:octic46:octic_delta2_q5_46:45-D2Q5-5}
\end{equation}
Again the source $\mu_5$ is radial.  The rational main exceptional
curve has preliminary square $-5/42$.  Its three cyclic charts have
arms $[2]$, $[2,2]$ and $[2,2,3]$, listed from the centre.
The corrections at the ends
attached to the main curve are $1/2,2/3,5/7$, respectively, so its
final square is $-2$.

The actual order of the main curve is three, since
$\operatorname{wt}(dx\wedge dy)=20$ and the radial quotient has
order five.  Its three adjacent arm orders are all two.
These orders can be seen directly.  At the $w=0$ orbit of order
two the form is $S^3v\,dS\wedge dv$, giving order two on $A_1$.
At the $x=0$ orbit of order three it is $S^3dS\wedge du$;
the $A_2$ orders from the centre outwards are $(2,1)$.
At the $y=0$ orbit of order seven the same expression has orders
$(2,1,0)$ from the centre outwards on $[2,2,3]$.

In these weights the source lower bounds are
$\operatorname{wt}\delta_x\ge6$ and
$\operatorname{wt}\delta_y\ge18$.
The other terms in \eqref{eq:octic46:octic_delta2_q5_46:45-D2Q5-4} have weighted degree at least
forty-five, against the leading degree twenty.  They acquire at
least five additional powers of $S$ and have larger orders on
the three indicated adjacent rays.  Thus the order-three centre
has three distinct smaller positive neighbors, all of coefficient
two.  A curve in a canonical blowup forest has at most two such
parents.  This excludes the final $A_3$ alternative without
assuming a generic higher jet.

\paragraph{Generic $A_4$ with $[y^2]h=0$}
The coefficient of $x^4y$ is nonzero by \eqref{eq:octic46:octic_delta2_q5_46:45-D2Q5-3}.
If $[xy^3]h\ne0$, weights $(x,y,w)=(4,6,11)$ give the isolated
initial equation
\[
                            w^2+x^4y+cxy^3.
\]
In the $x$--chart, $x=s^4$, $y=s^6u$, $w=s^{11}v$, the nonzero
roots of $1+cu^2=0$, $v=0$, give a residual
$\frac12(1,1)$ chart after $S=s^5$.
The leading form is again a unit times $Sv\,dS\wedge dv$.
The main curve and its adjacent $A_1$ have order one, a forbidden
equal adjacent pair.  Here the source lower bounds are six and
eight for $\delta_x,\delta_y$, and all errors in
\eqref{eq:octic46:octic_delta2_q5_46:45-D2Q5-4} have weighted degree at least thirty-five, against
the leading ten.  The actual-order assertion therefore survives
all higher terms.

If $[xy^3]h=0$, the same invariant $y\Omega$ test forces
$[y^7]h\ne0$.  Weights $(x,y,w)=(3,2,7)$ give the complete
initial equation
\begin{equation}
            w^2+y\{x^4+c x^2y^3+d y^6\}=0,
                         \qquad d\ne0.                 \label{eq:octic46:octic_delta2_q5_46:45-D2Q5-6}
\end{equation}
The middle coefficient $c$ is unrestricted.  In particular we do
not assume that the quadratic discriminant $c^2-4d$ is nonzero.
The reduced weighted exceptional curve $E\subset\mathbb P(3,2,7)$
is integral, and complete-intersection duality gives
\[
                   \omega_E=O_E(2),\qquad
                   h^1(E,O_E)=h^0(E,O_E(2))=1.
\]
The sole degree-two monomial is $y$.
This arithmetic-genus calculation also applies when
$c^2=4d$, when the curve is nodal after the relevant quotient.
The source $\mu_5$ acts radially and therefore does not change
the exceptional curve.

This contradicts rationality.  More explicitly, the normal weighted
partial resolution of a rational surface germ has vanishing first
structure-sheaf cohomology over a small affine neighborhood.
The exact sequence defining its reduced exceptional curve and the
vanishing of second cohomology for a proper map with fibres of
dimension at most one imply $H^1(O_E)=0$.  Alternatively one may
resolve the partial surface and use the corresponding effective
exceptional cycle on a resolution of the rational singularity.
Thus even the discriminant-zero degeneration of \eqref{eq:octic46:octic_delta2_q5_46:45-D2Q5-6} is
excluded.
\end{proof}

\paragraph{Connection to the whole nonzero-polar alternative.}
The exact Q12 normalization calculation supplies fifteen protected
roots for generic $A_3$ and seventeen for generic $A_4$.
The latter together with the six Q5 roots already exceeds the
maximal negative rank nineteen on the K3 surface.
In the $A_3$ case the nineteen Q5--Q12 roots leave no room for
any protected Q14 root.  The previously verified Q14 seven-stratum
analysis still applies, with the stronger polar mass forty-four.
All multiple polar curves avoid the zero-order $(-3)$ Q5--Q12
centres, so they have scaled degree at least five.  The positive
degree-thirty-eight and degree-forty-eight root chains are excluded
by the bound seven on their exceptional coefficients; the remaining
positive-core alternatives either fail the forest test proved above or
leave a protected root.  The exceptional $K$-core would require
external coefficient sum eleven and hence polar mass at least
$5\cdot11-2=53>44$.  Thus the $A_3$ alternative is excluded as well,
once these already established Q12 and Q14 interfaces are included.


\subsection{C2-02350: exclusion of the state $(7,2)$}
\label{sec:octic-degree7-delta2-zero-polar45}

We use the full ribbon-algebra reduction and the exact Q12
normalization charts.  The sole remaining ribbon has nilradical
$O_{\mathbb P(5,12)}(17)$ and is generated by $x$ at Q12 and by
$y$ at Q5.  The generic transverse type is $A_3$ or $A_4$.
At Q12 the normalization has an old line of square $-2$, an
attached $A_5$, and the following chain joining the old line to
the old cubic:
\begin{equation}
 A_3:\ [2,2,2,3,2,2,2],\qquad
 A_4:\ [2,2,2,2,4,2,2,2,2].                            \label{eq:octic46:octic_delta2_zero46:45-D2Z-1}
\end{equation}
All these are curves of the minimal normalization resolution.

\begin{theorem}
\label{thm:octic-degree7-delta2-zero-polar45}
The second polar cannot vanish in the sole degree--seven,
generic-delta-two state.
\end{theorem}

\paragraph{A higher intrinsic-jet gap forces tangency}

Suppose $\gamma_{10}|_\Gamma=0$.
In the full Q12 conductor ideal $(z,x^2)$, the quadratic transverse
part is a nonzero multiple of $y^2z^2$.  Consequently
$\delta_z|_\Gamma=0$.
At the Q12 origin the other derivation components, of characters
eight and six, vanish automatically.  Hence $\delta$ preserves
the maximal ideal at that point.  The source equation has ordinary
multiplicity four, so
\begin{equation}
                         \delta^n f\in\mathfrak m^4
                              \quad\hbox{for every }n. \label{eq:octic46:octic_delta2_zero46:45-D2Z-2}
\end{equation}

By Lemma~\ref{lem:octic-intrinsic-jets}, whenever the lower jets
vanish, $\delta^nf|_\Gamma$ is a section of $O_\Gamma(8+n)$.
At Q12, \eqref{eq:octic46:octic_delta2_zero46:45-D2Z-2} makes this section vanish to cover order at
least four.  In the complete section ring
\[
                         \mathbb C[X_5,Y_{12}],
\]
that means divisibility by $X^4$, and a nonzero such section has
degree at least twenty.  Induction therefore gives
\begin{equation}
                    \delta^nf|_\Gamma=0
                              \qquad(0\le n\le11).
                                                        \label{eq:octic46:octic_delta2_zero46:45-D2Z-3}
\end{equation}

If the normal value of $\delta$ were nonzero at the generic point
of $\Gamma$, a formal flow box would make $\delta=\partial_t$,
with $\Gamma=(t=z=0)$.  Its first integral $z$ is the square
direction of the transverse quadratic form, because the second
jet vanishes.  Weierstrass preparation, which does not change the
polar intersection on $f=0$, writes the transverse equation as
\[
                    F=z^2+a(t)z+b(t),\qquad
                    \operatorname{ord}b\ge12.
\]
The generic type $A_4$ would require discriminant order five,
which is impossible since $a^2$ has even order and $b$ has order
at least twelve.  The $A_3$ type requires $\operatorname{ord}a=2$.
If $N=\operatorname{ord}b<\infty$, the two roots have orders
$N-2$ and two, and evaluation of
$F_t=a'z+b'$ gives branch orders $N-1$ and three.
Thus
\[
                         I(F,F_t)=N+2\ge14.
\]
If $b=0$, the polar has a common branch, contrary to its properness.
Either way the nonzero-normal-value alternative is impossible:
its polar contribution would be at least $7\cdot14=98>72$.
This is also the usual polar identity
$I(F,F_t)=\mu(F)+I(F,z)-1$, but the preceding calculation proves
the required case directly.

The derivation is therefore tangent to $\Gamma$.
Its induced tangent section lies in
\[
 H^0(\mathbb P(5,12),T_\Gamma\otimes H)=H^0(O(18))=0.
\]
Consequently $\delta$ vanishes along the conductor everywhere.
No smooth ambient normal-bundle hypothesis is used.

\paragraph{Positive root chains at Q12}

All derivation coefficients vanish at the Q12 origin.  The actual
ambient two-form has ordinary degree at least four.  On the old
line of the ordinary exceptional quartic, after the radial
$\mu_2$ quotient, it therefore has strictly positive order.
This line has square $-2$ on the minimal normalization resolution.
Zero propagation through adjacent $(-2)$ curves now implies that
the connected chains
\begin{equation}
                     A_5+E_l+A_3=A_9\quad(A_3\hbox{ type}),
 \qquad
                     A_5+E_l+A_4=A_{10}\quad(A_4\hbox{ type})
                                                        \label{eq:octic46:octic_delta2_zero46:45-D2Z-4}
\end{equation}
are entirely positive.  The notation on the left describes the
chains geometrically; the type in parentheses is the generic
transverse type, not an additional curve.

Recall two canonical-forest facts.  A positive $(-2)$ chain has
strictly monotone positive integral coefficients.  Also every
positive exceptional curve on a minimal rational resolution has
a larger-coefficient child path ending at an intrinsic polar
curve: a terminal curve in such a path would otherwise be a
$(-1)$ curve.  Thus the largest exceptional coefficient is
strictly smaller than the largest intrinsic coefficient.
Here the intrinsic canonical polar mass is
\begin{equation}
                              420HB_N=44.              \label{eq:octic46:octic_delta2_zero46:45-D2Z-5}
\end{equation}

\paragraph{The cuspidal transverse type $A_4$ is impossible}

The reduced conductor normalization is one curve $R$ of scaled
degree seven.  Since $\delta$ vanishes along $\Gamma$, the generic
parametrization $x=t^2$, $z=t^5\cdot\text{unit}$ shows that the
actual form vanishes to order at least three on $R$.
Indeed $\delta_z\in(x,z)$ has order at least two, and its factor
$dx$ adds one; the other two terms have no smaller order.
Write
\[
                            B_N=rR+B',\qquad r\ge3.
\]
Then $r\le6$ and $420HB'=44-7r$.  Since every other positive
curve has scaled degree at least three, every intrinsic
coefficient is at most seven (the maximum occurs at $r=3$).
Every exceptional coefficient is therefore at most six.
The positive $A_{10}$ in \eqref{eq:octic46:octic_delta2_zero46:45-D2Z-4} is impossible.

\paragraph{The two-branch transverse type $A_3$ is impossible}

The unique reduced conductor multiplication is
\[
                            T_{43}^2=cX^{10}Y^3,
                                  \qquad c\ne0.
\]
Its normalized graded ring is $\mathbb C[X_5,S_6]$, with
$Y=S^2$ and $T=X^5S^3$, after a scalar rescaling.
The coarse reduced cover is a connected double cover of
$\Gamma$, branched at its two marked points.  Thus it is one
curve $R$ of scaled degree fourteen.  The vanishing of $\delta$
along $\Gamma$ makes its actual canonical coefficient at least
one.  Write $B_N=rR+B'$.

If $r\ge2$, the remaining mass is at most sixteen, so all
intrinsic coefficients are at most five and the positive
$A_9$ is already impossible.  Hence $r=1$.
The remaining mass is thirty, so every intrinsic coefficient is
at most ten and every exceptional coefficient is at most nine.
The nine positive coefficients on the $A_9$ must consequently
be exactly $1,2,\ldots,9$ in one order.
Its coefficient-nine endpoint has a larger child, necessarily
an intrinsic curve of coefficient ten.  The mass bound forces
the complete intrinsic divisor
\begin{equation}
                           B_N=R_{14}+10L_3.            \label{eq:octic46:octic_delta2_zero46:45-D2Z-6}
\end{equation}
In particular $L$ meets a positive Q12 exceptional curve.

The scaled degree three of $L$ forces actual non-Cartier contacts
at both Q5 and Q14.  Indeed $84H_N$ is Cartier away from the
Q5 fibre and at every smooth point of that fibre; if $L$ avoided
its non-Cartier points, its degree would be
$84\cdot3/420=3/5$, impossible for a line bundle on a complete
curve.  Similarly $60H_N$ gives $3/7$ and forces a Q14 contact.
The strict transform of $L$ therefore meets exceptional curves
over both sources.

At Q5 every minimal exceptional self-intersection is at least
$-3$.  This assertion is geometric and does not require a unit
actual form.  First the full ribbon orientation gives
$\mathfrak c=(z,y^2)$ and $f\in\mathfrak c^2$; the latter follows
from the generic $A_3/A_4$ equation and extends over the parameter
by torsion-freeness.  The invariant conductor form $z\Omega$ with
ordinary weights has canonical threshold four.  The only
character-four equation monomial below four in $\mathfrak c^2$
is $z^2$.  Thus its coefficient is a unit even when the actual
second polar vanishes, and the Morse normalization
$w^2+h(x,y)$ is valid here as well.
The boundary coefficient $x^3$ is fixed, so the
three possible normalization formats are precisely:
\[
 \begin{array}{c|c}
 {}[y^2]h\ne0 & [2,2,3,2,2]\\
 {}[y^2]h=0,\ [xy^3]h\ne0
   & \text{centre }-2\text{ with arms }[2],[3],[2,2,2]\\
 {}[y^2]h=[xy^3]h=0
   & \text{centre }-2\text{ with arms }[2],[2,2],[2,2,3].
 \end{array}
\]
The last case uses the same invariant-pole forcing of $y^7$ as
in Proposition~\ref{prop:octic-degree7-delta2-q5-nonzero45}.
For the middle case, the raw centre square is $-5/12$ and the
three corrections are $1/2,1/3,3/4$, giving square $-2$.
These isolated weighted initial equations give the displayed
charts for every higher-jet specialization.
At the normal Q14 source, the seven-stratum resolution proved above
analysis gives the uniform bound $-E^2\le5$.

If either Q5 or Q14 contact curve $E$ had canonical coefficient
zero, adjunction would give
\[
                       K_W E=-E^2-2\le3,
\]
whereas the term $10L$ in the effective canonical divisor would
contribute at least ten.  This is impossible.  Hence both
contacts are positive, as is the already known Q12 contact.
The coefficient-ten curve $L$ has no larger child: every other
intrinsic coefficient is one, and every exceptional coefficient
is at most nine.  It is therefore a final $(-1)$ curve with at
most two positive parents.  Its three distinct positive
neighbors at Q5,Q12,Q14 contradict that rule.
This excludes the final zero-polar alternative.


\paragraph{Conclusion of the delta-two case.}
\label{sec:degree7-delta2-complete45}

\begin{theorem}
\label{thm:degree7-delta2-complete45}
The C2-02350 conductor cannot consist of a single component of
scaled degree seven and generic delta invariant two.
\end{theorem}

\begin{proof}
The full ribbon theorem, Proposition~\ref{prop:octic-degree7-delta2-ribbon45},
fixes the nilpotent line to $O_{\mathbb P(5,12)}(17)$ and the
normalization-conductor trace-free restriction to $O(-43)$.
It also fixes the exact boundary multiplications in the two possible
generic transverse types $A_3,A_4$.
This removes the other ribbon orientations before any local
generic-coefficient assumption is made.

If the second polar vanishes, the higher intrinsic-jet argument of
Theorem~\ref{thm:octic-degree7-delta2-zero-polar45} forces the actual
derivation to vanish along the reduced conductor.  The exact Q12
normalization then contains a positive $A_9$ in type $A_3$ and a
positive $A_{10}$ in type $A_4$.
The polar mass is forty-four.  The latter type violates the
coefficient bound immediately; the former forces exactly
$B_N=R_{14}+10L_3$, whose three positive contacts at Q5,Q12,Q14
contradict the two-parent rule for the final $(-1)$ curve $L$.
This is the zero-polar theorem, including its uniform Q5 and Q14
self-intersection bounds.

Suppose instead that the second polar is nonzero.
Proposition~\ref{prop:degree7-delta2-q12-45} gives protected negative
root ranks fifteen and seventeen at Q12 in the two types.
Proposition~\ref{prop:octic-degree7-delta2-q5-nonzero45} gives ranks
four and six at Q5, respectively.  In type $A_4$ their sum is
twenty-three, already greater than the maximal negative rank
nineteen on a complex projective K3 surface.

For type $A_3$ the sum is nineteen.  It remains to show that the
normal Q14 source contributes a protected root.  We record why
the exhaustive Q14 analysis proved above applies with room to spare.
Every Q5 and Q12 exceptional curve has actual coefficient zero
and self-intersection at least $-3$.  An intrinsic polar curve
of coefficient at least two therefore avoids both sources.
Its scaled degree is at least five, while every positive intrinsic
curve has scaled degree at least three.  Since the total mass is
forty-four, an intrinsic coefficient is at most eight and a
positive coefficient on a minimal Q14 exceptional curve is at
most seven.  The increasing-child path proves the latter bound.
For a connected positive Q14 core, the external intersection
vector $v$ also satisfies
\[
                         5\sum_i v_i-2\le44,
                         \qquad \sum_i v_i\le9.
\]
Here at most one coefficient-one attachment is possible by the
unique-minimum property of a connected canonical forest.

After Corollary~\ref{cor:octic-q14-J}, six Q14 strata remain
from the uniform weighted residue analysis.  Their consequences
in the present state are:
\[
\begin{array}{c|l}
G&\text{a positive centre fails the core test; a zero centre protects an }A_2,\\
H&\text{the sole permitted positive core leaves a protected }A_4,\\
L&\text{no positive core survives even the weaker budget twelve},\\
K&\text{a unit polar protects }A_6;\text{ the sole other core requires }
       \sum v_i=11,\\
38&\text{a unit polar protects }A_8;\text{ otherwise a positive }A_8
       \text{ violates the bound seven},\\
48&\text{a unit polar protects }A_{10};\text{ otherwise a positive }A_{10}
       \text{ violates the bound seven}.
\end{array}
\]
The graphs and the arbitrary higher-jet cases are those in the
normal Q14 proof; no smoothness or genericity of an auxiliary
cubic is added here.  The $K$ core can also be rejected by its
remaining protected $A_2$, independently of the mass bound.
Thus a surviving normal Q14 germ would contribute a protected
root disjoint from the nineteen Q5--Q12 roots.  Their negative
definite rank would be at least twenty, again impossible.

Both polar alternatives and both generic transverse types have
now been excluded.
\end{proof}


\subsection{C2-02350: the second polar in the state $(7,3)$}
\label{sec:octic-degree7-delta3-nonzero}

\begin{theorem}
\label{thm:octic-degree7-delta3-nonzero}
For the sole degree--seven conductor state of C2-02350 with generic
$\delta=3$, the restriction $\gamma_{10}|_\Gamma$ must vanish.
\end{theorem}

\begin{proof}
Suppose $\gamma_{10}|_\Gamma\ne0$.  The smooth marked conductor is
$\mathbb P(5,12)$ with $H=O(1)$.  Its degree-ten section is a nonzero
multiple of $X_5^2$, hence is a unit on the Q5 cover and has exact
order two on the Q12 cover.  We use the uniform K3 reduction throughout.

\paragraph{A flat conductor algebra at Q12}
Put $A=\mathbb C[[y]]$ and use ordinary source coordinates
$\mu_{12}(x,y,z)=(1,11,5)$, with $f$ of character eight and conductor
support $I=(x,z)$.  Write $\mathfrak c$ for the full conductor.
The algebra
\[
                     B=O_D/\mathfrak c
\]
is Cohen--Macaulay of dimension one and therefore finite free of rank
three over $A$.  Indeed the full conductor is a maximal Cohen--Macaulay
module by finite duality, and its quotient has no zero-dimensional
torsion.  Its generic length is the delta invariant.  Generic length
three also gives $I^3B=0$: this holds over the fraction field and extends
by $A$--torsion-freeness.

The derivatives $f_x,f_z$ belong to $\mathfrak c$.  For example
$f_z\operatorname{Res}(dx\wedge dy\wedge dz/f)=dx\wedge dy$ is a regular
canonical form on the normalization, and finite duality identifies such
forms with $\mathfrak c$ times the residue.
The exact order-two second polar forces both
$[y^2z^2]f\ne0$ and $\delta_z(0)\ne0$.
The characters of the other normal derivation component give
$\delta_x|_\Gamma\in y^4 A$; their contributions to the quadratic
polar have order at least fourteen.

Reduce $f_z=0$ in $B$.  All terms of normal degree at least three vanish.
The full character list for its remaining terms gives
\begin{equation}
 y^2 u(y)z+y^{10}a(y)x+y^{11}b(y)x^2
                 +y^3c(y)xz+y^7d(y)z^2=0,\qquad u(0)\ne0.
                                                               \label{eq:octic46:octic_delta3_nonzero46:45-D3-1}
\end{equation}
Dividing by $y^2$ is valid in the free $A$--module $B$.
Solving this relation in the nilpotent ideal $IB$, whose cube is zero,
gives
\[
             z=a_8(y)x+b_9(y)x^2,\qquad
             \operatorname{ord}a_8\ge8,\quad
             \operatorname{ord}b_9\ge9.
\]
Consequently $B$ is generated by $x$, and the surjection
$A[x]/(x^3)\longrightarrow B$ is an isomorphism between free modules
of the same rank.  An equivariant regular change of $z$ therefore gives
\begin{equation}
                         \mathfrak c=(z,x^3).             \label{eq:octic46:octic_delta3_nonzero46:45-D3-2}
\end{equation}
The high orders in that change preserve the three coefficients used
below.

The invariant conductor forms $x^3\Omega$ and $z^3\Omega$ and the
uniform weighted canonical-form inequality force $[x^3z]f\ne0$ and
$[z^4]f\ne0$.  Since $f_z\in\mathfrak c$, mixed terms $x^jz$ with
$j<3$ are absent.  Since $f_x\in\mathfrak c$, pure $x$ terms of degree
at most three are absent.  Pure degrees four or five would give generic
delta two after the ordinary Morse splitting in $z$ over
$\mathbb C((y))$, and are also absent.  Thus
\[
             f=A_0x^6+2B_0x^3z+C_0z^2,\qquad B_0(0)\ne0,
\]
where $A_0$ has character two and
$C_0=c y^2+d z^2+\text{higher terms}$ with $cd\ne0$.
The discriminant $B_0^2-A_0C_0$ is a unit.  Put
$q=B_0+\sqrt{B_0^2-A_0C_0}$, choosing the square root for which $q$ is a
unit.  There is an actual analytic factorization
\begin{equation}
            f=\left(x^3+\frac{C_0}{q}z\right)
                         \left(A_0x^3+qz\right).          \label{eq:octic46:octic_delta3_nonzero46:45-D3-3}
\end{equation}
Both factors are separately semi-invariant.  This also proves that the
generic transverse type is nodal $A_5$, rather than cuspidal $A_6$.

\paragraph{Nineteen protected roots at Q12}
The second factor in \eqref{eq:octic46:octic_delta3_nonzero46:45-D3-3} is smooth, with coordinates $(x,y)$.
Its quotient is $\frac1{12}(1,11)$, and the actual form has unit leading
term $dx\wedge dy$.  It supplies a protected $A_{11}$.

The first factor has isolated homogeneous cubic initial equation
\begin{equation}
                         x^3+c' y^2z+d' z^3=0,
                         \qquad c'd'\ne0.                \label{eq:octic46:octic_delta3_nonzero46:45-D3-4}
\end{equation}
The projective cubic is smooth.  Ordinary blowup has a smooth elliptic
exceptional curve, the scalar source subgroup has order two, and the
residual group has order six.  Its quotient is rational and has three
orbifold orbits of stabilizers six, three and two.  The resolution
graph has a $(-3)$ centre with arms
\begin{equation}
                             A_5,\qquad A_2,\qquad A_1.
                                                               \label{eq:octic46:octic_delta3_nonzero46:45-D3-5}
\end{equation}
For completeness the raw centre square is $-1$; the three endpoint
corrections are $5/6,2/3,1/2$, summing to two.

The actual form has order zero on these eight roots.  In the $y$--chart,
$y=s,x=su,z=sv$, the strict equation eliminates $v$, and after
$S=s^2$ the residual action is $\frac16(1,5)$.
Here $dx\wedge dy$ is a unit times $dS\wedge du$.
At the orbit $(0:\pm\sqrt{-d'/c'}:1)$ the residual action is
$\frac13(1,2)$, and at the three points on $y=0$ it is
$\frac12(1,1)$; the same differential calculation gives a unit
canonical form at both orbits.  The coefficients
$\delta_x,\delta_y$ have ordinary orders at least four and two,
respectively.  Their terms acquire positive powers of $S$, and cannot
cancel these units.  The isolated smooth cubic initial cone makes all
these calculations unchanged by the higher terms of \eqref{eq:octic46:octic_delta3_nonzero46:45-D3-4}.
The two disjoint normalization points over Q12 therefore supply
\begin{equation}
                           11+5+2+1=19
                                                               \label{eq:octic46:octic_delta3_nonzero46:45-D3-6}
\end{equation}
protected negative roots.

\paragraph{The Q5 alternatives both contradict the K3 model}
At Q5 use source coordinates of characters $(1,4,2)$ and conductor
the $x$--axis.  The unit second polar forces a unit $z^2$ coefficient
and $\delta_z(0)\ne0$.  An equivariant Morse change preserving the
conductor puts the generic $A_5$ equation in the form
\[
             f=z^2+y^6h(x,y),\qquad h(x,0)\not\equiv0.
\]
Normalization adjoins $w=z/y^3$, giving
\begin{equation}
              w^2+h(x,y)=0,\qquad
              \mu_5(x,y,w)=(1,4,0),\qquad \operatorname{char}h=0.
                                                               \label{eq:octic46:octic_delta3_nonzero46:45-D3-7}
\end{equation}
It is normal: a repeated divisor of $h$ away from $y=0$ would give a
second conductor component, and $h(x,0)\not\equiv0$ excludes a repeated
factor $y$.  As a hypersurface, it is Cohen--Macaulay.

If $h(0,0)\ne0$, there are two smooth normalization branches, each with
quotient $\frac15(1,4)$ and unit actual form $dx\wedge dy$.
Their protected roots together with \eqref{eq:octic46:octic_delta3_nonzero46:45-D3-6} contradict the maximal
negative rank nineteen of a projective K3 surface.

Suppose $h(0,0)=0$.  The quotient in \eqref{eq:octic46:octic_delta3_nonzero46:45-D3-7} is Gorenstein because
its canonical character is $1+4+0=0$.  It is rational, hence Du Val.
Its source cover is finite and etale away from the origin, so is itself
Du Val.  The only invariant quadratic monomial of $h$ is $xy$, and
there are no invariant cubic monomials.  If $[xy]h=0$, the hypersurface
would have rank-one quadratic part $w^2$ and zero cubic part on its
quadratic kernel.  This is incompatible with every Du Val normal form.
Thus $[xy]h\ne0$.

Equivariant Morse splitting gives $w^2+XY=0$.  Its quotient is
$A_9=\frac1{10}(1,9)$, with universal variables
$X=u^2,Y=v^2,w=uv$.  The actual form has leading numerator a nonzero
multiple of $uv$.  Indeed, substituting $z=y^3w$ gives
\[
 \alpha=(\delta_z-3y^2w\delta_y)\,dx\wedge dy
       +y^3\delta_x\,dy\wedge dw-y^3\delta_y\,dx\wedge dw.
\]
Here $\delta_z$ is a unit and the other two coefficients vanish.
On the universal plane the first term is a unit multiple of $uv$;
every correction to that unit has positive degree.  The remaining
displayed terms have numerator degree at least ten.
The quotient action is symplectic, so each nonconstant invariant
numerator has positive integral order on every $A_9$ exceptional ray.
The term $uv$ has order exactly one on all nine rays.  No higher
monomial can cancel it identically along any such divisor, since
distinct characters on its open torus are linearly independent.
All nine canonical coefficients are therefore one.  Two adjacent
positive coefficients of a K3 blowup forest cannot be equal.  This is
the final contradiction.
\end{proof}

\begin{proposition}
\label{prop:octic-degree7-delta3-curvilinear}
In the sole degree--seven, delta-three state, the special fibre of the
full Q12 conductor algebra cannot be curvilinear.  Thus any surviving
Q12 conductor algebra has special Hilbert function $(1,2)$ and
$\gamma_{10}|_\Gamma=0$.
\end{proposition}

\begin{proof}
If the special fibre is curvilinear, equivariant Nakayama and the
distinct characters of $x,z$ show that it is generated by either $x$
or $z$.  Since $I^3B=0$ and $B$ is free of rank three, the corresponding
regular equivariant straightening gives exactly one of
\[
                \mathfrak c=(z,x^3),\qquad
                \mathfrak c=(x,z^3).
\]
The shifts in the first case have orders eight and nine as above;
in the second case they are
$x\mapsto x-a_4(y)z-b_9(y)z^2$, with orders at least four and nine.

A plane-curve conductor containing a regular parameter has multiplicity
at most two.  This is the first-blowup adjoint condition
$\mathfrak c\subset\mathfrak m^{m-1}$ for multiplicity $m$.
Equivalently it follows by pulling a conductor residue through that
ordinary blowup.  Consequently the generic type here is $A_5$ or $A_6$.
In the second orientation its equation belongs to $(x,z^3)^2$:
complete the square in $x$ over $\mathbb C((y))$ and use delta three,
then extend the resulting coefficient vanishings over $\mathbb C[[y]]$.
This contradicts the forced $z^4$ coefficient, which is unchanged by
the indicated high-order shift.

In the first orientation the regular conductor numerator $y^2z$ and
weights $(3,1,3)$ force $[y^2z^2]f\ne0$.  The factorization
\eqref{eq:octic46:octic_delta3_nonzero46:45-D3-3} and its smooth $A_{11}$ branch therefore apply without a
condition on the second polar.  If $\gamma_{10}\ne0$, the theorem
already gives a contradiction.  Otherwise the rank-one quadratic
part is $C_0(y)z^2$, and its evaluation on the normal derivation shows
$\delta_z|_\Gamma=0$.  The actual form on the smooth normalization
branch consequently vanishes at its Q12 fixed point: the other terms
contain derivatives of $z=O(x^3)$ and cannot contribute a unit.
All eleven $A_{11}$ coefficients are positive.

Here $420H\cdot B_N=72-2\cdot3\cdot7=30$.  Every intrinsic polar
curve has degree at least three, so every intrinsic coefficient is at
most ten.  The $A_{11}$ resolution contains no $(-1)$--curve.
Following the larger-coefficient blowup child from any of its positive
vertices must therefore leave this finite fibre along an intrinsic
polar curve.  Every chain coefficient is at most nine.  But a positive
$(-2)$ chain of length eleven in a canonical blowup forest has strictly
monotone positive integer coefficients and hence maximum at least
eleven.  This is impossible.
\end{proof}


\subsection{C2-02350: multiplicity three in the state $(7,3)$}
\label{sec:degree7-delta3-triple-algebra45}

The two nilpotent layers already give the square obstruction.

\begin{theorem}
\label{thm:degree7-delta3-triple-algebra45}
In the sole degree-seven conductor state of C2-02350 with generic
$\delta=3$, the generic transverse curve cannot have multiplicity three.
\end{theorem}

\begin{proof}
The actual reduced index stack is
\[
 \mathcal G=\mathbb P(5,12),\qquad H=\mathcal O(1),\qquad
 G=\mathbb C[X_5,Y_{12}],\qquad\omega_{\mathcal G}=\mathcal O(-17).
\]
The reduced-source identification and the exclusion of extra
non-Cartier contacts use only the saturated degree-seven trace
bound, and hence apply to every generic delta.
At Q5 the lifted curve is the $x$-axis in characters $(1,4,2)$;
at Q12 it is the $y$-axis in characters $(1,11,5)$.
All schemes and line bundles below are on this actual index stack.

Let $B=\mathcal O_{C_D}$ be the full image conductor algebra and
$I$ its nilradical.  A reduced plane curve of multiplicity three
and delta three has conductor exactly the square of its maximal
ideal.  One may see this without choosing its analytic type:
the first-blowup adjoint condition puts the conductor inside
$\mathfrak m^2$, and both quotients have length three.
Equivalently the three possible types $D_4,D_5,E_6$ give this
ideal directly.
Thus $I^2=0$ generically.  Since $B$ is Cohen--Macaulay of dimension
one, it has no finite-support subsheaf, and $I^2=0$ everywhere.
The nilradical is a rank-two vector bundle on the smooth stack
$\mathcal G$.  Its fibres are generated by $(y,z)$ at Q5 and by
$(x,z)$ at Q12.

\paragraph{The two nilpotent line bundles are forced.}
By Lemma~\ref{lem:octic-weighted-line-splitting}, write
$I=\mathcal O(m_1)\oplus\mathcal O(m_2)$.
Lemma~\ref{lem:octic-full-conductor-duality} gives $H^1(B)=0$
and the acyclicity of $B(9H)$.

Constants lift to $B$, and $\mathcal O_{\mathcal G}(9)$ is
acyclic.  Thus each summand satisfies
\begin{equation}
 G_{m+9}=G_{-26-m}=G_{-17-m}=0.                 \label{eq:octic46:octic_delta3_triple46:45-F3-1}
\end{equation}
The two fibre matchings have classes $(17,1)$ or $(41,37)$ modulo
$60$.  The exhaustive calculation in
\eqref{eq:octic46:octic_delta2_ribbon46:45-RB6}, applied to
\eqref{eq:octic46:octic_delta3_triple46:45-F3-1}, excludes the first
and gives $-19,-23$ for the second.  Therefore
\begin{equation}
             I=U\oplus V,\qquad
             U=\mathcal O(-19),\quad V=\mathcal O(-23).
                                                        \label{eq:octic46:octic_delta3_triple46:45-F3-2}
\end{equation}
Twisting $0\to I\to B\to O_{\mathcal G}\to0$ by $jH$,
for $j=3,6$, gives $H^1(B(jH))\simeq H^1(I(j))$ because
$O_{\mathcal G}(j)$ is acyclic. The two layers in degree three
are $O(-16),O(-20)$ and have no first cohomology. In degree six
they are $O(-13),O(-17)$, whose first cohomology has dimensions
zero and one. Hence
\[
                         h^1(B(3H))=0,\qquad h^1(B(6H))=1.
\]
Lemma~\ref{lem:octic-conductor-square-obstruction} excludes this case.
\end{proof}


\subsection{C2-02350: the remaining double branch in the state $(7,3)$}
\label{sec:octic-degree7-delta3-double-algebra45}

\begin{theorem}
\label{thm:octic-degree7-delta3-double-algebra45}
In the sole degree-seven conductor state, generic delta three cannot
have transverse multiplicity two.  Together with
Theorem~\ref{thm:degree7-delta3-triple-algebra45}, this excludes the
entire degree-seven, delta-three numerical state.
\end{theorem}

\begin{proof}
The two generic types are $A_5,A_6$.
The nonzero-second-polar alternative and every curvilinear special
fibre at $Q_{12}$ were excluded in Theorem~\ref{thm:octic-degree7-delta3-nonzero} and
Proposition~\ref{prop:octic-degree7-delta3-curvilinear}.
We may therefore assume
\begin{equation}
 \gamma_{10}|_\Gamma=0,\qquad
 B|_{Q_{12}}\text{ has special Hilbert function }(1,2),
 \qquad420H B_N=30.
\label{eq:octic46:octic_delta3_double46:45-FD1}
\end{equation}
Here $B=\mathcal O_{C_D}$ is the full image-conductor algebra on
$\mathcal G=\mathbb P(5,12)$, with $H=\mathcal O(1)$ and
$G=\mathbb C[X_5,Y_{12}]$.  We retain the full-scheme statements
$H^1(B)=0$ and $H^0(B(9H))=H^1(B(9H))=0$.
These are the full-conductor vanishings of
Lemma~\ref{lem:octic-full-conductor-duality}, including all Cartier contacts.

\paragraph{The saturated socle filtration.}
Let $I$ be the nilradical of $B$ and put $K=\operatorname{Ann}_B(I)$.
Generically $B=k[e]/(e^3)$.  Thus $I^3=0$ everywhere by purity,
$K\subset I$ is a line on $\mathcal G$, and $I^2\subset K$.
The line $K$ is saturated: if a local curve parameter $y$ satisfies
$yh\in K$, then $yIh=0$ and Cohen--Macaulayness gives $Ih=0$.
Hence $L=I/K$ is also a line.  Write
\[
                    L=\mathcal O(l),\qquad K=\mathcal O(k).
\]
Locally over a smooth parameter ring the algebra is
\begin{equation}
 A[e,f]/(e^2-a(y)f,ef,f^2),\qquad a(y)\not\equiv0.
\label{eq:octic46:octic_delta3_double46:45-FD2}
\end{equation}
The intrinsic multiplication $L^2\to K$ is therefore a nonzero
global section of $\mathcal O(k-2l)$.

The quotient $B/K$ is a ribbon.  Its first cohomology vanishes,
so the unit shows $H^1(L)=0$.
The exact sequences obtained by twisting the filtration by $9H$
give
\begin{equation}
 H^1(L(9))=0,\quad H^0(K(9))=0,\quad
 h^0(L(9))=h^1(K(9)),\quad h^1(K)\le h^0(L).
\label{eq:octic46:octic_delta3_double46:45-FD3}
\end{equation}
No splitting of $B$ is used for these inequalities.
In terms of the two-variable section spaces they say
\begin{equation}
 \begin{gathered}
 G_{-17-l}=G_{-26-l}=G_{k+9}=0,\
 \dim G_{l+9}=\dim G_{-26-k},\
 \dim G_{-17-k}\le\dim G_l,\qquad G_{k-2l}\ne0.
 \end{gathered}
\label{eq:octic46:octic_delta3_double46:45-FD4}
\end{equation}
Every integer at least 44 is in $\langle5,12\rangle$.
Consequently \eqref{eq:octic46:octic_delta3_double46:45-FD4} puts $-60\le l\le17$ and
$2l\le k\le34$; the enumeration is finite with proved bounds.

At $Q_{12}$ the two special nilpotent characters are one and five.
Thus the characters of $(L,K)$ are $(1,5)$ or $(5,1)$.
At $Q_5$ the possible pairs are
\begin{equation}
                         (4,3),\quad(2,4),\quad(4,2).
\label{eq:octic46:octic_delta3_double46:45-FD5}
\end{equation}
The first is the curvilinear fibre generated by $y$, with socle
$y^2$, and its multiplication is a unit at $Q_5$.
The second includes the curvilinear fibre generated by $z$ as well
as its noncurvilinear degeneration.  The last is the other
noncurvilinear orientation.  These exhaust the possibilities because
the actual source coordinates generate the conductor algebra and
their characters are distinct.
Line-module characters of $H$ are $+2,+5$ at $Q_5,Q_{12}$.
Applying \eqref{eq:octic46:octic_delta3_double46:45-FD4}--\eqref{eq:octic46:octic_delta3_double46:45-FD5} leaves exactly
\begin{equation}
\begin{array}{c|c|c}
 l&k&(L,K)\text{ at }Q_5\\\hline
 -19&-23&(2,4)\\
 -23&-19&(4,2)\\
 -23&29&(4,3).
\end{array}
\label{eq:octic46:octic_delta3_double46:45-FD6}
\end{equation}
The proved degree bounds and the source characters determine this
necessary list; the argument below uses only the two line-bundle layers.

Put $B_0=B/K$. The actual twisted sequences
\[
 0\to K(j)\to B(jH)\to B_0(jH)\to0,\qquad
 0\to L(j)\to B_0(jH)\to O_{\mathcal G}(j)\to0
\]
compute the required cohomology without a retraction.
For $j=3$, both $K(j)$ and $L(j)$ have vanishing first cohomology
in all three rows, so $h^1(B(3H))=0$.
In the first row, $K(6)=O(-17)$ has one-dimensional first
cohomology, while $L(6)=O(-13)$ is acyclic.
In the other two rows, $L(6)=O(-17)$ has one-dimensional first
cohomology and $H^1(K(6))=0$. Since $O_{\mathcal G}(6)$ is acyclic,
every row gives $h^1(B(6H))=1$.
Lemma~\ref{lem:octic-conductor-square-obstruction} excludes all three
rows and completes the delta-three case.
\end{proof}


\subsection{C2-02350: multiplicity two in the state $(7,4)$}
\label{sec:octic-degree7-delta4-double}

\begin{theorem}
\label{thm:octic-degree7-delta4-double}
In the sole degree-seven conductor state with generic $\delta=4$,
the generic transverse curve cannot have multiplicity two.
\end{theorem}

\begin{proof}
The possible transverse types are $A_7,A_8$, and the intrinsic polar
degree is
\[
                  420H\cdot B_N=72-2\cdot4\cdot7=16.
\]
Every intrinsic polar component has scaled degree at least three.
First suppose $\gamma_{10}|_\Gamma=0$.
For $A_7$, transverse Morse coordinates give $v^2-u^8$ and branches
$v=\pm u^4$.  The $v$-component of the normal derivation vanishes on
the conductor.  If its coefficient of $u$ is nonzero, the actual form
has order one on both branches; otherwise both orders are at least two.
Their total polar mass is therefore fourteen or at least twenty-eight.
The first leaves the impossible residual degree two, and the second
exceeds sixteen.  This argument also applies when the branches are
exchanged over the function field of $\Gamma$.
For $A_8$, write $u=t^2,v=t^9$.  The transverse actual form
$\delta_v\,du-\delta_u\,dv$ has order at least three, since
$\delta_v\in(u,v)$.  Its polar mass is at least $3\cdot7=21$.
Thus $\gamma_{10}|_\Gamma\ne0$.

Work on the actual reduced index stack
$\mathcal G=\mathbb P(5,12)$, with
$G=\mathbb C[X_5,Y_{12}]$ and $H|_{\mathcal G}=O(1)$.
The nonzero second polar is a multiple of $X_5^2$; it is a unit at
Q5 and has order two on the ordinary cover at Q12.
Let $B=O_{C_D}$, let $I$ be its nilradical, and set
\[
\begin{gathered}
 J_2=\operatorname{sat}_B(I^2),\qquad
 J_3=\operatorname{sat}_B(I^3),\\
 L=I/J_2=O(l),\quad M=J_2/J_3=O(m),\quad K=J_3=O(k).
\end{gathered}
\]
The generic algebra is $k[e]/(e^4)$.  Purity and saturation therefore
give
\[
 I^4=0,\qquad IJ_2\subset J_3,\qquad IJ_3=J_2^2=0.
\]
Consequently $L,M,K$ are line bundles on $\mathcal G$, and $B/J_2$
is a ribbon with nilpotent line $L$.

We determine $L$ using only this ribbon quotient.
At Q12 the coordinates $(x,y,z)$ have characters $(1,11,5)$ and the
equation has character eight.  Write its normal quadratic part as
$A(y)x^2+2B_1(y)xz+C(y)z^2$.
The character bounds are
$\operatorname{ord}_y A\ge6$, $\operatorname{ord}_y B_1\ge10$,
and $\operatorname{ord}_y\delta_x|_\Gamma\ge4$.
In the second polar, the terms $A\delta_x^2$ and
$2B_1\delta_x\delta_z$ thus have order at least fourteen.
Its exact order two forces $C=y^2u(y)$ with $u(0)\ne0$ and
$\delta_z(0)\ne0$.
The Jacobian element $f_z$ belongs to the generic conductor of
$A_7$ or $A_8$, hence to the full conductor by torsion-freeness of $B$.
In $B/J_2$ all normal quadratic terms vanish, so $f_z=0$ gives
\[
                  y^2u(y)z+y^{10}a(y)x=0.
\]
Saturation makes this quotient torsion-free over $\mathbb C[[y]]$.
Dividing by $y^2$ shows that $L$ is generated by $x$, with fibre
character one.  At Q5 the unit second polar similarly forces a unit
$z^2$ coefficient; the $yz$ coefficient is divisible by $x^3$.
The same Jacobian relation makes $L$ generated by $y$, of character
four.  Since the line-module characters of $H$ are $+5,+2$ at Q12,Q5,
\begin{equation}
             5l\equiv1\pmod {12},\qquad2l\equiv4\pmod5,
             \qquad l\equiv17\pmod {60}.
\label{eq:octic-delta4-ribbon-character}
\end{equation}

Lemma~\ref{lem:octic-full-conductor-duality} gives
$H^1(B)=H^1(B(9))=0$ and $H^0(B(9))=0$.
The first two vanishings pass to $B/J_2$.
Its ribbon sequence, the lifting of constants by the unit, and the
acyclicity of $O_{\mathcal G}(9)$ yield
\[
             H^1(L)=H^1(L(9))=0,
             \qquad G_{-17-l}=G_{-26-l}=0.
\]
Every integer at least forty-four lies in $\langle5,12\rangle$,
so $l\ge-60$.  By \eqref{eq:octic-delta4-ribbon-character}, the only
possibility below seventeen is $l=-43$, excluded by $G_{17}\ne0$.
Thus $l\ge17$.

The displayed saturated-product relations give well-defined maps
\[
                    L^{\otimes2}\longrightarrow M,
                    \qquad L\otimes M\longrightarrow K.
\]
Both are nonzero at the generic point, where they are the products of
$e,e^2,e^3$.  Hence $m\ge2l$ and $k\ge l+m\ge51$.
But $K(9)\subset B(9)$ gives $G_{k+9}=0$, whence $k\le34$.
This contradiction excludes both transverse types.
\end{proof}

\subsection{C2-02350: multiplicity three and completion of the state $(7,4)$}
\label{sec:octic-degree7-delta4-triple}

\begin{theorem}
\label{thm:octic-degree7-delta4-final}
The sole degree-seven conductor state with generic $\delta=4$ is
impossible.
\end{theorem}

\begin{proof}
By Theorem~\ref{thm:octic-degree7-delta4-double}, it remains to treat
transverse multiplicity three; multiplicity at least four has
$\delta\ge6$.  The formula
$\delta=\sum\delta_i+\sum_{i<j}I(C_i,C_j)$ gives the four types below.
For three branches the intersection numbers are $(2,1,1)$; for two
branches one has a $(2,5)$ cusp and a transverse line, or a $(2,3)$
cusp and a line of intersection three.  The irreducible case has
semigroup $\langle3,5\rangle$.
\[
\begin{array}{c|c|c}
 \text{type}&\text{equation}&\text{image conductor ideal}\\\hline
 D_6&u(v^2+u^4)&(u^3,uv,v^2)\\
 D_7&u(v^2+u^5)&(u^3,uv,v^2)\\
 E_7&u(u^2+v^3)&(u^2,uv,v^3)\\
 E_8&u^3+v^5&(u^2,uv,v^3).
\end{array}
\]
Indeed the conductor orders on the normalization branches are,
respectively, $(2,3,3)$, $(2,6)$, $(3,5)$ and $8$.
They put each displayed ideal inside the conductor; its colength is
four, so equality follows.  Every generic conductor algebra therefore
has Hilbert function $(1,2,1)$ and quadratic multiplication of rank one.

Keep $B=O_{C_D}$ on the actual index stack $\mathcal G=\mathbb P(5,12)$,
and let $I$ be its nilradical.  Purity extends $I^3=0$ from the
generic point.  Put
\[
 K=\operatorname{sat}_B(I^2)=O(k),\qquad F=I/K,
                  \qquad B_0=B/K.
\]
Then $IK=0$, $F$ is a rank-two vector bundle on $\mathcal G$, and
$B_0$ is a square-zero thickening with nilradical $F$.
At Q5 and Q12 the special fibre has length four, at most two normal
generators, and cube-zero maximal ideal.  Its Hilbert function must
therefore be $(1,2,1)$, so $I^2\to K$ is nonzero on both fibres.
The quadratic multiplication retains rank one.  At Q12 the three
quadratic characters of $x^2,xz,z^2$ are distinct.  A sole nonzero
mixed term would have rank two; thus the nonzero square is $x^2$ or
$z^2$.  At Q5 it is similarly $y^2$ or $z^2$.

By Lemma~\ref{lem:octic-weighted-line-splitting}, write
$F=O(l_x)\oplus O(l_z)$, labelled by the Q12 characters one and five.
The source characters give
\begin{equation}
 5l_x\equiv1,\quad5l_z\equiv5\pmod {12},\qquad
             \{2l_x,2l_z\}\equiv\{2,4\}\pmod5.
\label{eq:octic-delta4-filtration-characters}
\end{equation}
Moreover $5k\equiv2$ or $10\pmod {12}$ for the Q12 squares
$x^2$ or $z^2$, and $2k\equiv3$ or $4\pmod5$ for the Q5 squares
$y^2$ or $z^2$.
At either source, a nonzero square of the summand $O(l_i)$ gives
the necessary inequality $k\ge2l_i$.

The sequences $0\to K\to B\to B_0\to0$ and
$0\to F\to B_0\to O_{\mathcal G}\to0$, together with
Lemma~\ref{lem:octic-full-conductor-duality}, give
\begin{equation}
\begin{gathered}
 H^1(F)=H^1(F(9))=0,\qquad H^0(K(9))=0,\\
 \chi(K(9))+\chi(F(9))=0,\qquad h^1(K)\le h^0(F).
\end{gathered}
\label{eq:octic-delta4-filtration-cohomology}
\end{equation}
Constants lift by the unit of $B$, which gives the last inequality.
These statements use the actual twists of the filtration.

For completeness the finite bounds follow directly from these conditions.
Write $g(d)=\dim G_d$, with $g(d)=0$ for $d<0$.
Since every integer at least forty-four belongs to $\langle5,12\rangle$,
\eqref{eq:octic-delta4-filtration-cohomology} gives $l_x,l_z\ge-60$
and $k\le34$.  A nonzero source square gives $k\ge-120$.
In this range $g(-26-k)\le2$, and the Euler equality gives
$g(l_x+9)+g(l_z+9)=g(-26-k)$.
Every degree at least $164$ has at least three monomials: multiply a
monomial of degree $d-120\ge44$ by the three degree-$120$ monomials.
Thus the complete search is confined to
\[
                  -60\le l_x,l_z\le154,\qquad-120\le k\le34.
\]
Using $h^0(O(d))=g(d)$, $h^1(O(d))=g(-17-d)$, the source characters
and the two nonzero-square inequalities leaves exactly
\begin{equation}
\begin{array}{c|c|rrr}
 \text{Q12 square}&\text{Q5 square}&l_x&l_z&k\\\hline
 x^2&y^2&-19&-23&34\\
 x^2&z^2&-19&-23&22\\
 x^2&z^2&-19&37&-38\\
 z^2&y^2&-19&-23&14\\
 z^2&y^2&41&-23&-46\\
 z^2&z^2&-19&-23&2.
\end{array}
\label{eq:octic-delta4-filtration-candidates}
\end{equation}
Here the inequalities are only necessary conditions; no realization of
these filtrations is assumed.

The same two sequences now compute the cohomology needed for
Lemma~\ref{lem:octic-conductor-square-obstruction}.
For $j=3,6$, the sheaf $O_{\mathcal G}(j)$ is acyclic, so
$H^i(B_0(j))=H^i(F(j))$ for $i=0,1$.
In the four rows with $F=O(-19)\oplus O(-23)$, both $K(3),K(6)$
have vanishing $H^1$, while
\[
                         h^1(F(3))=0,\qquad h^1(F(6))=1.
\]
For $(-19,37,-38)$, one has $h^1(K(3))=g(18)=0$,
$h^1(F(3))=0$, and
\[
 h^1(K(6))=g(15)=1,\qquad
                H^0(F(6))=H^1(F(6))=0
\]
because $g(43)=0$.
For $(41,-23,-46)$, the relevant groups satisfy
$h^1(K(3))=g(26)=0$, $h^1(K(6))=g(23)=0$ and
$(h^1(F(3)),h^1(F(6)))=(0,1)$.
Thus every row gives
\[
                         h^1(B(3))=0,\qquad h^1(B(6))=1.
\]
The square obstruction excludes all six candidates, and hence every
remaining transverse type.
\end{proof}

\section{Completion of the basket exclusions}
\label{sec:octic-completion-revised}

The local and global arguments are now complete.  We assemble them
against the exhaustive conductor list, and hence exclude the last
numerical basket under consideration.

\subsection{C2-02350: exclusion of the octic basket}
\label{sec:octic-complete-exclusion46}

\begin{theorem}\label{thm:octic-complete-exclusion46}
There is no minimal complex projective threefold $X$ of general type
with $\mathbb Q$-factorial terminal singularities, $\chi(\mathcal O_X)=2$,
$K_X^3=1/420$, and Reid basket
\[
       5(1,2)+(1,3)+2(1,4)+(2,5)+(5,12)+(5,14).
\]
\end{theorem}
\begin{proof}
Propositions~\ref{prop:octic-prime46} and \ref{prop:octic-sources46}
give the reduced visible prime octic on its canonical model and the
actual source neighborhoods used throughout the preceding proof.
Theorem~\ref{thm:octic-K3-rational46} proves that its normalization has
rational singularities and K3 minimal model.  Proposition
\ref{prop:octic-nonnormal46} shows that its conductor is nonempty.
Full conductor duality, the actual polar degree bound and the complete
source-character reduction give precisely the eleven necessary states
of Theorem~\ref{thm:octic-eleven-states46}.

Every state has been excluded.  The following table records the exact
coverage; a pair of components is one total state.
\begin{center}
\begin{tabular}{c|p{103mm}}
State $(w,\delta)$ & Exclusion \\\hline
$(7,1)$ & The nodal case is Theorem~\ref{thm:octic-sole-degree7-node};
the two cuspidal polar alternatives are
Theorems~\ref{thm:octic-degree7-cusp-gamma-zero45} and
\ref{thm:octic-degree7-cusp-gamma-nonzero45}.\\
$(7,2)$ & Theorem~\ref{thm:degree7-delta2-complete45}.\\
$(7,3)$ & The multiplicity-two and multiplicity-three cases are
Theorems~\ref{thm:octic-degree7-delta3-double-algebra45} and
\ref{thm:degree7-delta3-triple-algebra45}.\\
$(7,4)$ & Theorem~\ref{thm:octic-degree7-delta4-final}, using the
saturated conductor filtration and the square obstruction.\\
$(12,1)$ & Theorem~\ref{thm:octic-degree12-cusp-excluded45}.\\
$(14,1)$ & Proposition~\ref{prop:octic-degree14-delta1-node45} and
Theorem~\ref{thm:octic-degree14-delta1-cusp45}.\\
$(14,2)$ & Theorem~\ref{thm:octic-degree14-delta2-45}.\\
$(18,1)$ & Theorems~\ref{thm:octic-degree18-node} and
\ref{thm:octic-degree18-cusp}.\\
$(21,1)$ & Theorem~\ref{thm:octic-degree21-node45}, with its
actual, possibly nonordinary, index-four source.\\
$(25,1)$ & Theorem~\ref{thm:octic-degree25-node}.\\
$(7,1)+(18,1)$ & Theorem~\ref{thm:octic-pair-deletion}.
\end{tabular}
\end{center}
This exhausts all necessary conductors and contradicts the existence of
the first octic.  Hence the basket itself is impossible.
\end{proof}

\paragraph{The finite calculations.}
\label{sec:octic-arithmetic46}

The finite calculations use integers and exact rational numbers.
Their inputs and bounds are specified in the proofs: source-character
minimum orders give the conductor list; quotient graphs give the
intersection matrices; the polar mass bounds the exceptional forests;
and conductor cohomology bounds the line-bundle degrees.
Appendix~\ref{app:reproducibility46} records the corresponding
computation sources, including the matrices, monomial lists and
bounded enumerations.  The local geometry, rationality and full
conductor identities needed for these calculations were proved above.



\clearpage
\part{Vega and Altair: construction, recognition and geometry}

\section{A common reconstruction framework}
\label{sec:exceptional-common-framework}

We now turn to the two baskets that survive the exclusions. The
quartic basket $\mathcal B_4$ gives Vega, and the quintic basket
$\mathcal B_5$ gives Altair. We retain the notation
$Z_\sharp=\mathrm{Vega}$, $Z_*=\mathrm{Altair}$ for the explicit
models and $R_\sharp,R_*$ for their canonical rings. An arbitrary
realization will be denoted by $Z$.

There are three separate tasks. We must construct each model and
verify its canonical polarization, recover its entire canonical
ring from an arbitrary realization of the basket, and prove ordinary
Artinian rigidity. The first task does not depend on the recognition
theorem. The last task concerns all ordinary deformations, without
requiring a chosen involution or algebra presentation to remain fixed.
We begin with the common surface and algebraic arguments, and then
carry out the two model-specific reconstructions.

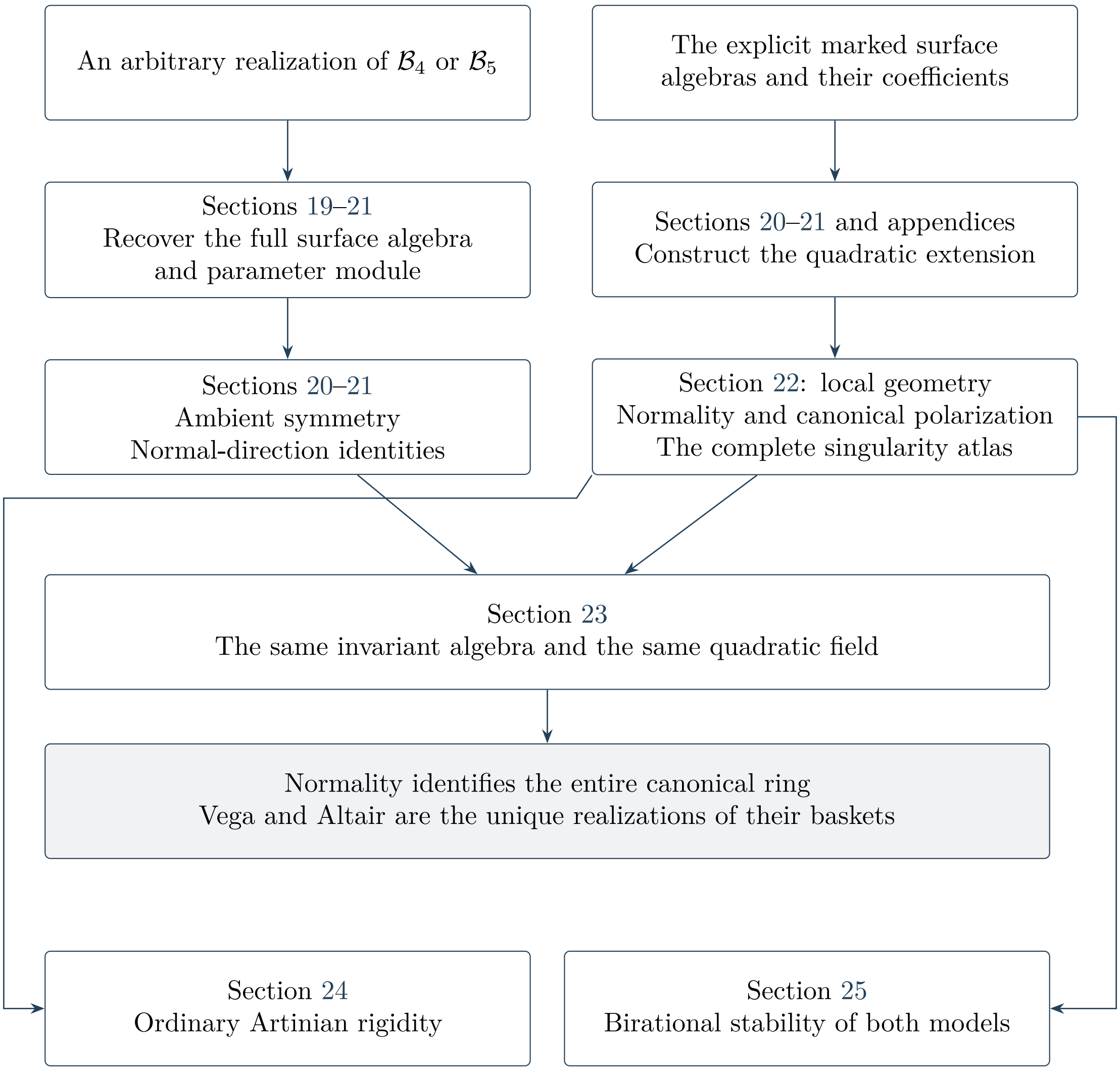
\begin{figure}[htbp]
\centering
\begin{tikzpicture}[>=Stealth,font=\small,
 box/.style={draw=inkblue,rounded corners=2pt,align=center,inner sep=5pt,
             text width=60mm,minimum height=15mm},
 wide/.style={box,text width=128mm},
 route/.style={->,draw=inkblue,line width=.5pt}]
\node[box] (arbitrary) {An arbitrary realization of $\mathcal B_4$ or $\mathcal B_5$};
\node[box,right=8mm of arbitrary] (explicit) {The explicit marked surface\\algebras and their coefficients};
\node[box,below=8mm of arbitrary] (surface) {Sections~\ref{sec:exceptional-common-framework}--\ref{sec:altair-model}\\Recover the full surface algebra\\and parameter module};
\node[box,below=8mm of explicit] (build) {Sections~\ref{sec:vega-model}--\ref{sec:altair-model} and appendices\\Construct the quadratic extension};
\node[box,below=8mm of surface] (symmetry) {Sections~\ref{sec:vega-model}--\ref{sec:altair-model}\\Ambient symmetry\\Normal-direction identities};
\node[box,below=8mm of build] (geometry) {Section~\ref{sec:exceptional-canonical-geometry}: local geometry\\Normality and canonical polarization\\The complete singularity atlas};
\node[wide,below=13mm of symmetry,xshift=34mm] (same) {Section~\ref{sec:exceptional-final-recognition}\\The same invariant algebra and the same quadratic field};
\node[wide,below=7mm of same,fill=inkblue!7] (normal) {Normality identifies the entire canonical ring\\Vega and Altair are the unique realizations of their baskets};
\node[box,below=12mm of normal,xshift=-34mm] (rigid) {Section~\ref{sec:exceptional-ordinary-rigidity}\\Ordinary Artinian rigidity};
\node[box,below=12mm of normal,xshift=34mm] (maps) {Section~\ref{sec:stability-completion}\\Birational stability of both models};
\draw[route] (arbitrary)--(surface);
\draw[route] (surface)--(symmetry);
\draw[route] (explicit)--(build);
\draw[route] (build)--(geometry);
\draw[route] (symmetry)--(same);
\draw[route] (geometry)--(same);
\draw[route] (same)--(normal);
\draw[route] (geometry.east)--++(5mm,0)|-(maps.east);
\draw[route] (geometry.south west)--++(-2mm,-3mm)--++(-75mm,0)|-(rigid.west);
\end{tikzpicture}
\caption{Construction and recognition for Vega and Altair. The outer
arrows indicate properties proved directly from the explicit rings.
Recognition then transfers those properties to an arbitrary realization.}
\label{fig:exceptional-proof}
\end{figure}


\subsection{Small adjunction and rational surface singularities}
\label{sec:quintic-surface-tools46}

We first isolate a surface argument used for both models. The point is
to prove rationality before identifying the minimal surface as K3.
A bicanonical gap will then distinguish the K3 case from the
properly elliptic case.

\begin{lemma}[A small square prevents contracted irregularity]
\label{lem:small-square-irregularity46}
\label{lem:small-adjunction-irregularity34}
Let $D$ be a normal projective surface with $h^1(\mathcal O_D)=0$,
let $H$ be an ample $\mathbb Q$-Cartier Weil divisor, and let
$\pi:W\to D$ be a resolution with $p_g(W)=1$. Put
$A=\pi^*H$, $d=A^2$ and $u=AK_W$. If $u\ge0$ and $u^2<d$, then
$q(W)=0$, provided the minimal model has nonnegative Kodaira dimension.
If also $h^2(\mathcal O_D)=1$, then $D$ has rational singularities.

\end{lemma}
\begin{proof}
Contract $W$ to its smooth minimal model $W_0$, and use the orthogonal
basis of total transforms of the successive exceptional curves:
\[
 A=\sigma^*A_0-\sum_i a_iE_i^*,\quad
 d=A_0^2-\sum_i a_i^2,\quad
 u=A_0K_{W_0}+\sum_i a_i,
 \qquad a_i\ge0.
\]
The class $A_0$ is nef and big. If $W_0$ were of general type,
Hodge index would give
$1\le K_{W_0}^2\le (A_0K_{W_0})^2/A_0^2\le u^2/d<1$.
Thus its Kodaira dimension is zero or one.

Suppose first that it is one and $q(W)>0$. Its elliptic base has
positive genus. Indeed for a base $\mathbb P^1$, relative duality and
Leray give $R^1f_*\mathcal O_{W_0}=\mathcal O(q-2)$ and
$q=h^0(\mathcal O(q-2))=\max(q-1,0)$, which is impossible.
We show that every $A$-null curve is vertical for this elliptic map.
If its image $C_0$ were horizontal, its geometric genus would be at
least one. Put $s=K_{W_0}C_0\ge1$, $P=A_0K_{W_0}$ and $b=P/s$.
If $m_i$ are its successive multiplicities at the blowup centres,
adjunction and the genus drop formula give
\[
 C_0^2\ge\sum_i m_i(m_i-1)-s,
 \qquad A_0C_0=\sum_i a_im_i.
\]
The class $A_0-bC_0$ is perpendicular to the nonzero nef isotropic
class $K_{W_0}$, so its square is nonpositive. Hence
\[
\begin{split}
 d&\le b^2\Bigl(s+\sum_i m_i\Bigr)
                       -\sum_i(a_i-bm_i)^2\\
  &\le b^2s+2b\sum_i a_i
   =2bu-b^2s\le u^2/s\le u^2.
\end{split}
\]
The second inequality uses $m_i\ge1$ at the centres on the curve:
$b^2m_i-(a_i-bm_i)^2\le2ba_i$; unused centres only decrease the
left side. This contradicts $u^2<d$. The positive-genus fibration
therefore descends through $\pi$, since it is constant on its connected
exceptional fibres. Its Stein factor gives a positive-genus curve
quotient of $D$, contradicting $h^1(\mathcal O_D)=0$.

If the Kodaira dimension is zero and $q>0$, the minimal model with
$p_g=1$ is abelian. Now $u=\sum a_i$. An image $C_0$ of an $A$-null
curve with geometric genus at least two satisfies
$C_0^2\ge\sum m_i(m_i-1)+2$. Set $b=A_0C_0/C_0^2$ and apply Hodge
index to $A_0-bC_0$. The preceding calculation becomes
\[
 d\le2bu-2b^2\le u^2/2,
\]
again impossible. Every nonconstant image is therefore a smooth
elliptic translate. If two images have different directions, write
$h=C_0D_0\ge1$, $p=A_0C_0$, $q'=A_0D_0$, and let $c$ be the sum of
the blowup weights at their common centres. Distinct elliptic directions
meet transversely, so there are at most $h$ such centres. Hence
\[
 p+q'\le u+c,\qquad A_0^2\le2pq'/h,\qquad
 \sum_i a_i^2\ge c^2/h,
\]
and
\[
 d\le(2pq'-c^2)/h\le u^2/h\le u^2.
\]
Thus all elliptic images are parallel. The quotient of the abelian
surface by this elliptic subgroup descends to $D$, again contradicting
$h^1(\mathcal O_D)=0$. If there are no such images, the birational map
to the abelian surface itself descends and gives the same contradiction.
This proves $q(W)=0$.

Finally $p_g(W)=1$ and $q(W)=0$ give $\chi(\mathcal O_W)=2$.
If $h^1(\mathcal O_D)=0$ and $h^2(\mathcal O_D)=1$, then
$\chi(\mathcal O_D)=2$ as well. The Leray sequence gives
$\operatorname{length}R^1\pi_*\mathcal O_W
=\chi(\mathcal O_D)-\chi(\mathcal O_W)=0$.
Thus the singularities are rational.

\end{proof}

\begin{corollary}[A bicanonical gap forces the K3 minimal model]
\label{cor:exceptional-bicanonical-k3}
Let $D,H,W,A,d,u$ be as in
Lemma~\ref{lem:small-square-irregularity46}, and assume
$u\ge0$, $u^2<d$ and $h^0(D,\omega_D^{[2]})\le1$.
Then the smooth minimal model of $D$ is K3.
\end{corollary}
\begin{proof}
The lemma gives $q(W)=0$ and excludes general type. A regular
bicanonical form on $W$ restricts to the regular locus of $D$ and
extends as a reflexive form. Hence
$P_2(W)\le h^0(D,\omega_D^{[2]})\le1$.
If the minimal surface were properly elliptic, its base would be
$\PP^1$ and $\chi(\mathcal O)=2$. The elliptic canonical-bundle
formula would then give
\[
 P_2=1+\#\{\text{multiple fibres}\}\ge2:
\]
at least one multiple fibre is necessary for positive Kodaira
dimension, and every multiplicity $m\ge2$ contributes
$\lfloor2(m-1)/m\rfloor=1$. This is a contradiction.
Thus the minimal model has Kodaira dimension zero; $p_g=1,q=0$
identify it as a K3 surface.
\end{proof}


\subsection{Parameter sections and the underlying canonical module}
\label{sec:exceptional-common-module}

Let $R$ be a positively graded domain with $R_0=\CC$, and let
$a\in R_e$ be nonzero. Suppose that the complete surface algebra $B$
containing $R/(a)$ is finite over three actual sections
$x_1,x_2,x_3\in R/(a)$. Then $R/(a)$ is finite over
$A_0=\CC[x_1,x_2,x_3]$. Lift these sections and a finite set of
module generators to $R$. Induction on degree shows that the lifts
generate $R$ over $\CC[a,x_1,x_2,x_3]$. If $\dim R=4$, the four
parameters are therefore algebraically independent. We shall use
this argument only after constructing the required parameter
sections on the complete surface.

Write $\Omega_A^i(\CC)$ for the $i$th syzygy in the minimal Koszul
resolution of the residue field of a weighted polynomial ring $A$.

\begin{lemma}[The two intermediate cohomology classes]
\label{lem:exceptional-common-module}
Let $A=\CC[x_1,x_2,x_3,x_4]$, with weights $w_i\ge2$, and put
$\mathfrak m=A_+$ and $W=\sum_iw_i$. If a finite graded $A$-module
$M$ has dimension four, depth at least two, and
\[
 H^2_{\mathfrak m}(M)=\CC(-1),\qquad
 H^3_{\mathfrak m}(M)=\CC,
\]
then, for a graded free module $F$,
\begin{equation}
 M\simeq F\oplus\Omega_A^2(\CC)(-1)\oplus\Omega_A^3(\CC).
 \tag{QR.1}\label{eq:exceptional-common-module}
\end{equation}
\end{lemma}
\begin{proof}
Local duality gives
$\Ext_A^2(M,A)=\CC(W+1)$ and $\Ext_A^1(M,A)=\CC(W)$.
In a minimal free resolution the last term is consequently
$A(-W-1)$. After homogeneous changes of bases and parameters, the
last dual differential is the row
$(x_1,x_2,x_3,x_4,0,\ldots,0)$: its cokernel is exactly
$\CC(W+1)$. Its kernel is
$\Omega_A^2(\CC)(W+1)\oplus G$, with $G$ free.

No minimal generator of the Koszul summand has degree $-W$, since
that would require $w_i+w_j=1$. Thus the quotient $\CC(W)$ comes
from $G$. Minimality makes the contributing summand $A(W)$, and
the preceding image is
\[
 \Omega_A^2(\CC)(W+1)\oplus\mathfrak m(W).
\]
Its minimal free cover is the sum of the two Koszul covers. Extra
free generators have zero image after a homogeneous basis change.
Dualizing and using Koszul self-duality proves the decomposition.
\end{proof}

For an arbitrary canonical realization of either exceptional basket,
the local-cohomology hypothesis follows from canonical vanishing,
duality and $p_g=q=0$, $h^2(\mathcal O)=1$. Once the parameter frame
has been proved, the Hilbert series determines the free summand.
The two possibilities are
\[
\begin{array}{c|c|l}
 &\text{parameter weights}&\text{degrees of the free generators}\\\hline
\text{Vega}&(4,10,12,14)&0,15,18,20,20,21,21,23,26,41\\
\text{Altair}&(5,8,12,14)&0,15,16,18,19,21,22,24,25,40.
\end{array}
\]
The weight products are both $6720$. At volume $1/420$ the full
module therefore has rank sixteen; the two Koszul summands have
rank three each. In particular the full canonical ring has depth
two, although its invariant algebra will be Cohen--Macaulay.

\begin{corollary}[Rigidity of the underlying module]
\label{cor:exceptional-common-module-rigidity}
Suppose that a finite graded $A$-algebra $R$ satisfies the preceding
lemma, $R_1=0$, and $\omega_R\simeq R(1)$. Then
\begin{equation}\label{eq:exceptional-module-Ext46}
 \Ext_A^1(R,R)_0=0.
\end{equation}
\end{corollary}
\begin{proof}
Write $F=\bigoplus_d A(-d)^{f_d}$. Since
$\omega_R=\Hom_A(R,A(-W))$, Koszul self-duality gives
\[
 \Hom_A(\Omega_A^2(\CC)(-1),A(-W))
       \simeq\Omega_A^3(\CC)(1),\qquad
 \Hom_A(\Omega_A^3(\CC),A(-W))\simeq\Omega_A^2(\CC).
\]
Comparison with $R(1)$ cancels the matching Koszul Hilbert series
and yields $\sum f_dt^{W+1-d}=\sum f_dt^d$.
Thus $f_W=f_1=0$: $R_1=0$, and the Koszul generator degrees
are greater than one.

Dimension shifting and the self-dual Koszul complex identify
\[
 \Ext_A^1(R,R)_0\simeq
 \Tor_1^A(\CC,R)_{W+1}\oplus(R/\mathfrak mR)_W.
\]
The first relations of $\Omega_A^2(\CC)(-1)$ have degrees
$w_i+w_j+w_k+1<W+1$, and the first relation of
$\Omega_A^3(\CC)$ has degree $W$. Hence the first term vanishes.
The generator degrees of the two Koszul summands are
$w_i+w_j+1$ and $w_i+w_j+w_k$; none equals $W$, since $w_i\ge2$.
The second term is therefore $\CC^{f_W}=0$.
\end{proof}

This corollary concerns the module only. Multiplication on that module
may still deform. The ordinary rigidity argument below includes all
multiplication directions and all changes of the parameter frame.

\subsection{Finite graded extensions and normality}
\label{sec:exceptional-finite-algebra-laws}

We shall describe the invariant rings as finite free algebras. The
positive grading makes both the extension calculations and their
changes of coordinates finite.

\begin{lemma}\label{lem:exceptional-finite-laws}
Let $A=\CC[a,x_1,\ldots,x_s]$ have positive weights, with
$\deg a=e$, and let $M=\bigoplus_{i=0}^{r-1}A(-d_i)$ have basis
$b_0=1,b_1,\ldots,b_{r-1}$, where $0=d_0\le d_i\le D$.
Fix a commutative unital multiplication on $M/aM$.
Its homogeneous commutative unital extensions are defined by
finitely many polynomial equations in the coefficients of
\[
 b_ib_j=\sum_h c_{ij}^h b_h,\qquad
 \deg c_{ij}^h=d_i+d_j-d_h.
\]
No product has $a$-order greater than $\lfloor2D/e\rfloor$;
checking basis associators through order $\lfloor3D/e\rfloor$
proves associativity in every degree. A homogeneous change of a
degree-$d_i$ generator has $a$-order at most $\lfloor d_i/e\rfloor$.
\end{lemma}
\begin{proof}
Every homogeneous piece of a positively weighted polynomial ring is
finite dimensional. The unit, symmetry and central-fibre conditions
are linear conditions on the structure constants. Associativity is
$A$-trilinear, so it suffices to test basis triples. The coefficient
on $b_h$ in such a triple has degree
$d_i+d_j+d_k-d_h\le3D$; the product bound is obtained with two
inputs. The bound on coordinate changes follows in the same way.
\end{proof}

For Vega, $(e,D)=(4,41)$ gives product and associator bounds
$20,30$; for Altair, $(e,D)=(5,40)$ gives $16,24$.
The same finite-degree argument applies to a presented module,
provided every module relation is imposed on the multiplication.
It also explains why the homogeneous formal coordinate changes
used in recognition are polynomial: each generator has a fixed
positive degree.

\begin{lemma}[Lifting a normal central fibre]
\label{lem:exceptional-normality-lift}
Let $S$ be a finitely generated positively graded Cohen--Macaulay
$\CC$-algebra with $S_0=\CC$,
and let $a$ be a homogeneous nonzerodivisor of positive degree.
If $S/(a)$ is a normal domain, then $S$ is a normal domain.
\end{lemma}
\begin{proof}
The separated $a$-adic filtration has associated graded ring
$(S/(a))[t]$, so $S$ is a domain. Work at the homogeneous vertex.
A height-one prime containing $a$ is $(a)$ and is regular. For a
height-one prime $\mathfrak p$ not containing $a$, take a prime
$\mathfrak q$ minimal over $(\mathfrak p,a)$; it has height two.
The quotient $S_{\mathfrak q}/(a)$ is regular, being a height-one
localization of the normal central fibre. Regularity lifts through
a nonzerodivisor, so $S_{\mathfrak q}$ and then $S_{\mathfrak p}$
are regular. Serre's criterion proves normality at the vertex.
The nonnormal locus is homogeneous and closed; if nonempty it
would contain the vertex. Thus it is empty.
\end{proof}

\subsection{Multiplication on a presented graded module}\label{ag:finite-multiplication-schemes-without-a-cohenmacaulay-assumption}

\begin{proposition}\label{ag:prop:multiplication-scheme}
Let \(k\) be a computable field of
characteristic zero, \[
A=k[x_1,\ldots,x_s],\qquad \deg x_i>0,
\] and let \(M\) be an explicitly presented finite graded \(A\)-module.
Fix \(1_M\in M_0\). The commutative, associative, unital graded
\(A\)-algebra structures on \(M\), with this unit and the given
\(A\)-module structure, form an explicitly computable affine scheme of
finite type over \(k\).
\end{proposition}
\begin{proof}
The space \[
V=\operatorname{Hom}_A(M\otimes_AM,M)_0
\] is finite dimensional. Images of finitely many homogeneous generators
lie in specified finite-dimensional graded pieces, and compatibility
with the module relations imposes finite linear equations. A
multiplication is an element \(\mu\in V\). Commutativity and
\(\mu(1_M,m)=m\) are linear conditions. Associativity is the vanishing
of \[
\mu(\mu(-,-),-)-\mu(-,\mu(-,-))
\in\operatorname{Hom}_A(M^{\otimes_A3},M)_0,
\] so gives finitely many quadratic equations. The relevant Hom spaces
and equations are computable by graded linear algebra and syzygies.
These constructions commute with scalar extension from the field \(k\).
\end{proof}

No freeness assumption is needed. In the free case
\(M=\bigoplus A(-d_i)\) with \(0\leq d_i\leq D\), and for a parameter
\(a\) of degree \(e\), product coefficients have \(a\)-order at most
\(\lfloor2D/e\rfloor\) and associator coefficients at most
\(\lfloor3D/e\rfloor\). These are the bounds of Lemma~\ref{lem:exceptional-finite-laws}. For a nonfree module, all relations remain part of the calculation.

Graded module automorphisms fixing the unit act on this scheme. The
quotient groupoid classifies structures with the chosen parameter
algebra and underlying module. To obtain moduli of canonical models, one
must also allow changes of parameter frame, identify redundant
presentations and impose the geometric conditions. A tangent calculation
in one format does not prove ordinary rigidity.

\subsection{The full quadratic normalization}
\label{sec:exceptional-quadratic-normalization}

\begin{lemma}\label{lem:exceptional-quadratic-normalization}
Let $S$ be a normal noetherian domain containing $\CC$. Suppose
$0\ne F,G\in S$ and $H\in S$ satisfy $H^2=FG$, no height-one prime
contains both $F$ and $G$, and $F$ is not a square in $K=\Frac S$.
For $v^2=F$ and $I=(F):_S H$, the full integral closure of $S$ in
$K(v)$ is
\[
                  R=S\oplus I/v.
\]
This finite quadratic cover is etale in codimension one. If the
grading is compatible and $\omega_S\simeq S(1)$, then the
reflexive canonical module is $\omega_R\simeq R(1)$.
\end{lemma}
\begin{proof}
At a height-one valuation positive on $F$, the element $G$ is a
unit; hence $\operatorname{ord}F=2\operatorname{ord}H$.
For $x\in I$ the colon condition gives
$2\operatorname{ord}x\ge\operatorname{ord}F$ at every height-one
prime. Normality implies $x^2/F\in S$, and similarly $xy/F\in S$
for $x,y\in I$. Thus $S\oplus I/v$ is a finite integral algebra.

Conversely, trace splits an integral element into its invariant part
in $S$ and an integral anti-invariant part $\xi$. Both $v\xi$ and
$(H/v)\xi$ are integral invariants, so they belong to $S$.
Therefore $v\xi\in I$, which proves completeness of the displayed
normalization. The even valuations of $F$ give an unramified
quadratic extension at every height-one prime. Pullback of canonical
forms there, followed by reflexive extension, proves the final
assertion.
\end{proof}

In the two constructions the hypotheses will be checked on the
independently constructed surface algebra. Once the invariant ring
and nontrivial quadratic field have been identified for an arbitrary
realization, this lemma identifies the entire canonical ring.
The canonicity of the projective model is a separate geometric
assertion, proved from its local charts in
Section~\ref{sec:exceptional-canonical-geometry}.

\subsection{Recognizing the canonical
polarization}\label{ag:recognizing-the-canonical-polarization}

The following criterion makes explicit the polarization check required
after a normal graded algebra has been constructed.

\begin{proposition}\label{ag:prop:canonical-recognition}
Let \(R\) be a finitely generated normal
positively graded domain, \(R_0=k\), of dimension at least three, and
\(Z=\operatorname{Proj}R\). Choose \(r>0\) with \(R^{(r)}\) generated in
degree one. Suppose that:

\begin{enumerate}
\item
  \(Z\) has canonical singularities;
\item
  \(\mathcal O_Z(r)\simeq\omega_Z^{[r]}\);
\item
  for \(0\leq j<r\), there are isomorphisms
  \(\widetilde{R(j)}\simeq\omega_Z^{[j]}\);
\item
  for every \(0\leq i,j<r\), these identifications respect the
  multiplication map
  \(\widetilde{R(i)}\otimes\widetilde{R(j)}\to\widetilde{R(i+j)}\),
  using \(\mathcal O_Z(r)\simeq\omega_Z^{[r]}\) when \(i+j\geq r\).
\end{enumerate}

Then \(K_Z\) is ample and \(R\simeq R(Z,K_Z)\) as graded algebras.
\end{proposition}
\begin{proof}
Write \(m=qr+j\). Twisting over
\(\operatorname{Proj}R^{(r)}\) gives \[
\widetilde{R(m)}
\simeq\widetilde{R(j)}\otimes\mathcal O_Z(r)^{\otimes q}
\simeq\omega_Z^{[m]},
\] compatibly with products. Normality gives depth at least two at the
irrelevant maximal ideal, so the Serre--Grothendieck sequence identifies
\(R_m\) with the global sections of its associated twist. Finally
\(\mathcal O_Z(r)\) is ample.
\end{proof}

These are finite residue-degree data after \(r\) is chosen, but
canonical singularities and the sheaf identifications are substantive
hypotheses. A shift such as \(\omega_R\simeq R(1)\), especially for a
non-Cohen--Macaulay ring, should not be used as a substitute for the
whole recognition argument.


\section{Vega: the quartic surface and its canonical-ring candidate}
\label{sec:vega-model}

The first quartic carries more information than its Hilbert series.
Its polar and its integral fifth root determine a marked elliptic
K3 surface. We shall recover this polarized surface from every
realization of $\mathcal B_4$, construct its invariant algebra,
and then construct the ambient normal cover. Finally we extend the
surface involution to an arbitrary ambient realization. The local
geometry and final recognition are proved in
Sections~\ref{sec:exceptional-canonical-geometry} and
\ref{sec:exceptional-final-recognition}.

\subsection{The first quartic of the basket $\mathcal B_4$}
\label{sec:quartic-reduction46}

\begin{lemma}[Small visible surfaces]\label{lem:small-visible46}
Let $T$ be a visible prime component of a member of $|eK_Y|$,
and $W$ a smooth model of its normalization. If
$e(e+1)^2V<1$, then $\kappa(W)\in\{0,1\}$. If the complete
dualizing space of the member has dimension one, then
$p_g(W)=1$ and $T$ is its only distinct visible prime.
\end{lemma}
\begin{proof}
Put $d=K^2T$ and $A=K|_W$. Normalization adjunction, including
its effective conductor, and Hodge index give
\[
 A^2=d\le eV,\qquad AK_W\le d+KT^2\le d+d^2/V.
\]
Visibility gives a nonzero holomorphic two-form on $W$. If its
minimal model were of general type, its nef canonical pullback
$B$ would have integral square at least one, whereas
\[
 1\le B^2\le(AB)^2/d\le d(1+d/V)^2
                                  \le e(e+1)^2V<1.
\]
This is impossible. The last assertion follows from ordinary
dualizing trace and the independence of its contributions at the
distinct generic points, including for repeated components.
\end{proof}

\begin{proposition}\label{prop:quartic-prime46}
For every realization of C2-02347, the unique quartic canonical
member $D$ is a reduced visible prime. Its normalization has
$p_g=1$ and Kodaira dimension zero or one. The complete $12K_Z$
system is a connected pencil with no divisorial fixed part.
\end{proposition}
\begin{proof}
The row has $I=420$, $v=1$, and
\[
 P_4=P_8=P_{10}=1,\quad P_5=P_6=P_{13}=0,\quad P_{12}=2.
\]
Every one of its $32$ necessary opposite torsion states has
$H_{12}(\pm\tau)\ge1$. Multiplication into $R_{13}=0$ therefore
excludes any effective $K_Z+\tau$. The visible RR list through
degree three is $(1,1),(2,4),(3,9)$. Null descent removes degree
one, giving a visible-prime floor of two.

For $a\in R_4\setminus\{0\}$ the complete dualizing space has
dimension $P_5+1=1$, so there is a unique distinct visible prime.
It is not repeated, by Lemma~\ref{lem:polar-square46} and $P_6=0$.
Its Hodge mass is four. A visible degree two or three would be
a Hodge equality with that smaller mass, a contradiction. Thus
it has degree four and exhausts the entire divisor on $Z$.
Lemma~\ref{lem:small-visible46}, with $4\cdot5^2/420=5/21$,
gives the assertions about its normalization.

Choose $p\in R_{12}$ independent of $a^3$. If $p$ vanished along
the prime $D$, normal division would give $p/a\in R_8=\mathbb Ca^2$.
Thus $a^3,p$ have no common divisor and span the complete system.
If its Stein degree were at least two,
\eqref{eq:pencil-area46} would give
$1\le k\le6\cdot7^2/420=7/10$, impossible.
\end{proof}


\subsection{Normality of the canonical quartic}
\label{sec:c2347-canonical-normality-continuation}

The argument in this subsection takes place on the canonical model.  It
therefore does not replace the canonical conductor by a terminal strict
transform's conductor.  Primeness is proved in Proposition~\ref{prop:quartic-prime46}.

\begin{lemma}[Canonical index transfer in the present basket]
\label{lem:c2347-canonical-index-continuation}
Let $g:Y\longrightarrow Z$ be a crepant terminalization of a canonical
realization of \texttt{C2-02347}.  Every local canonical index on
$Z$ belongs to $\{1,2,3,4,5,6,7\}$.  In particular $420K_Z$ is Cartier.
\end{lemma}

\begin{proof}
The terminal indices are among $1,2,3,4,5,7$. Proposition~\ref{prop:actual-source-transfer46}
shows that every additional canonical index is at most six. All the resulting
indices divide $420$.
\end{proof}

\begin{proposition}[Normality of the actual canonical quartic]
\label{prop:c2347-canonical-normality-continuation}
In a hypothetical canonical realization of \texttt{C2-02347},
the unique member $D=(a_4=0)$ is normal.  If $T$ is its strict transform
on a crepant terminal model, then
\[
 420K_Y^2T=4,\qquad 420K_YT^2=16.
\]
In particular the strict terminal state $(4,2)$ is impossible.
\end{proposition}

\begin{proof}
Proposition~\ref{prop:quartic-prime46} gives a reduced prime $D\sim4K_Z$.
It is Cohen--Macaulay: locally a Cartier member on the canonical cover,
its divisorial exact sequence descends by the character decomposition.
In particular it is $S_2$.  The nonzero ambient holomorphic two-form
restricts nontrivially to the normalization of $D$, by Proposition~\ref{prop:quartic-prime46}.

Let $S$ be a smooth resolution of the normalization $N$ of $D$ and
let $A$ be the pullback of $K_Z|_N$.  Write $C$ for the effective
conductor divisor, with its strict transform on $S$ understood.
All terms exceptional over points have zero intersection with $A$.
Canonical threefolds are Gorenstein in codimension two, so ordinary
Cartier adjunction and normalization give the exact numerical identity
\[
 A^2=\frac4{420}=\frac1{105},\qquad
 A(K_S+C)=5A^2=\frac{20}{420}.                 \tag{1}
\]
The restricted two-form gives $K_S$ an effective representative, hence
$AK_S\ge0$.

Suppose $D$ is not normal.  Since it is $S_2$, there is a curve
$\Gamma$ along which it is not normal.  At its generic point $Z$ is
Gorenstein and the transverse surface singularity is Du Val.  The
one-dimensional local ring of $D$ is a reduced Gorenstein ring, since
$D$ is Cartier in that surface.  For such a ring the conductor degree
on the normalization is twice its delta invariant.  This statement
includes branch residue degrees and does not require a plane-curve
embedding.  Indeed, for the finite normalization $R\subset\overline R$
and conductor $\mathfrak c$, one-dimensional Gorenstein duality applied
to $\overline R/R$ gives
$\ell(R/\mathfrak c)=\ell(\overline R/R)=\delta$; consequently
$\ell(\overline R/\mathfrak c)=2\delta$.
Non-normality therefore contributes at least
$2\Gamma$ to the pushforward of $C$.

Put $w=420K_Z\Gamma$.  Lemma~\ref{lem:c2347-canonical-index-continuation}
makes $w$ a positive integer.  From (1), ampleness, and the preceding
conductor lower bound, $2w\le20$, so $1\le w\le10$.

Take the generic log canonical threshold $t$ of $D$ along $\Gamma$.
Then $0<t\le1$, and $\Gamma$ is a log canonical centre of
$(Z,tD)$ at its generic point.  If $t<1$ this follows from a divisor
computing the threshold.  If $t=1$, a purely log terminal boundary
would be normal; non-normality gives a log canonical centre at the
generic curve.  Thus the canonical-model version of normalized
marked adjunction applies.  Lemma~\ref{lem:full-adjunction46} applies on the canonical
model; the Cartier multiple has just been verified.
For the canonical inertia fractions $k_x/r_x$ on $\Gamma^\nu$ it gives
\[
 2g(\Gamma^\nu)-2+\sum_x(1-1/r_x)
 \le (K_Z+tD)\Gamma
 \le\frac{5w}{420},\qquad
 \sum_x\frac{k_x}{r_x}\equiv\frac w{420}\pmod{\mathbb Z}. \tag{2}
\]
Every nontrivial $r_x$ is one of $2,3,4,5,6,7$ and
$\gcd(k_x,r_x)=1$.

There are at most four marked branches: each costs at least $1/2$,
whereas their total cost in the genus-zero relaxation of (2) is at
most $2+50/420<5/2$.  The exact enumeration of (2) for $1\le w\le10$
has only the following two rows, both with $w=10$:
\[
 \begin{array}{c|c|c}
 w&\text{marks}&-2+\sum(1-1/r_x)\\\hline
 10&(1/6,6/7)&-13/42\\
 10&(1/2,2/3,6/7)&1/42.
 \end{array}
\]
The accompanying standalone audit checks all multisets of at most
four primitive fractions, including composite inertia six.  It does
not use a prime-power denominator estimate.

Consequently the conductor contribution from this one curve is at
least $20/420$.  Equation (1) forces equality and $AK_S=0$.
The zero-adjunction minimalization argument now applies directly:
write $K_S=\sigma^*K_{S_0}+E$ on a smooth minimal model.  The
nonzero two-form implies nonnegative Kodaira dimension, and both
terms have nonnegative intersection with $A$.  Hodge index forces
$K_{S_0}$ to be numerically zero; positive geometric genus leaves
the K3 and abelian cases.  Every component of $E$ is $A$-null, so
the morphism to $N$ factors through $S_0$.  An abelian surface has
no birational contraction to a normal surface, and an integral Weil
polarization on a smooth abelian surface has integral square, which
contradicts $A^2=1/105$.  In the K3 case the contracted configurations
are ADE configurations, so $N$ is a Du Val K3.

The restriction of $K_Z$ defines an integral ample Weil divisor on
$N$: it is Cartier at every generic curve point, since the
non-Gorenstein locus of $Z$ is finite.  Its square is $1/105$.
The Du Val K3 square obstruction
(Lemma~\ref{lem:k3-square105-extension}) excludes this.  This
contradiction proves normality of $D$.

Finally a normal $D$ cannot contain a generic singular curve of
$Z$.  Transversely, a singular Du Val surface has embedding
dimension three; quotient by the one Cartier equation of $D$
has dimension one and embedding dimension at least two.  It is
not a discrete valuation ring, contrary to normality of $D$.
Hence the exceptional divisorial part $E_Y$ in
$g^*D=T+E_Y\sim4K_Y$ is supported over points.  Intersecting with a
general pullback of a sufficiently divisible ample canonical divisor
therefore gives $K_Y^2E_Y=K_YE_Y^2=0$.  Thus
\[
 420K_Y^2T=420\cdot4K_Y^3=4,\qquad
 420K_YT^2=420\cdot16K_Y^3=16,
\]
as claimed.
\end{proof}


\subsection{Rationality, klt singularities and the full quartic member}
\label{sec:c2347-small-adjunction-rationality}

We apply the common adjunction lemma to the normal quartic, and then
use the complete restriction gaps to control its discrepancies.

\begin{proposition}[Rational and klt quartic; no crepant remainder]
\label{prop:c2347-rational-klt-continuation}
In a realization of \texttt{C2-02347}, the normal
canonical quartic $D$ has rational and klt singularities.  The pair
$(Z,D)$ is purely log terminal.  For every crepant terminalization
$g:Y\to Z$, the full divisor of the quartic section is its prime
strict transform $T$, with no exceptional component, and $T$ is
normal and klt.
\end{proposition}

\begin{proof}
Normality was proved in
Proposition~\ref{prop:c2347-canonical-normality-continuation}.
Put $H=K_Z|_D$ and let $\pi:S\to D$ be a resolution.  The
common visible-surface reduction gives $p_g(S)=1$ and
$\kappa(S)\in\{0,1\}$.  Since $D$ is normal, adjunction gives
$K_D=5H$, and all resolution discrepancies are exceptional over
points.  Thus, for $A=\pi^*H$,
\[
 d=A^2=\frac1{105},\qquad u=AK_S=\frac1{21},
 \qquad \frac{u^2}{d}=\frac5{21}<\frac13.
\]
In particular a general-type minimal model would give
$1\le K_{S_0}^2\le u^2/d=5/21$, which is impossible.
The full divisorial sequence and Serre duality give
\[
 H^1(D,\mathcal O_D)=0,\qquad
 \chi(\mathcal O_D)=\chi(\mathcal O_Z)+P_5=2.
\]
Here $H^2(Z,\mathcal O_Z(-4K_Z))$ is dual to
$H^1(Z,\mathcal O_Z(5K_Z))=0$ by vanishing.
Lemma~\ref{lem:small-adjunction-irregularity34} gives $q(S)=0$,
so $\chi(\mathcal O_S)=2$.  Leray therefore gives
\[
 \sum_{P\in D}\operatorname{length}
       (R^1\pi_*\mathcal O_S)_P
   =\chi(\mathcal O_D)-\chi(\mathcal O_S)=0.
\]
Every surface singularity of $D$ is rational.

We use the Weil-gap criterion of \cite[Lemma~10.6(i)]{LBVolume}:
when $K_D\sim_{\mathbb Q}cH$, the vanishing $H^0(D,jH)=0$ for
$j>c$ implies that $D$ is klt wherever $jH$ is Cartier. Here $c=5$.

The Cohen--Macaulay canonical-cover quotient agrees with the
reflexive restriction sheaf on the normal quartic.  Its complete
restriction dimensions give
\[
 H^0(D,6H)=H^0(D,7H)=H^0(D,8H)=0.
\]
These three degrees cover all canonical ambient indices
$1,2,3,4,6,7$, and each exceeds the canonical weight five.
At the remaining index five, $K_D=5H$ is Cartier, so rationality
makes the surface singularity Du Val.  The established canonical
index list is exhaustive, hence $D$ is klt everywhere.

The different in divisorial adjunction is zero: at a generic curve
point the normal quartic is regular and Cartier in the canonical
Gorenstein threefold, which is consequently regular there.
Inversion of adjunction gives that $(Z,D)$ is purely log terminal.
In
\[
 g^*D=T+\sum_i c_iE_i
\]
every $c_i$ is a nonnegative integer, since this is the actual
divisor of the pulled-back quartic canonical section and $g$ is
crepant.  Its pair discrepancy is $-c_i$, whereas pure log
terminality requires $c_i<1$.  Thus all $c_i=0$.
The crepant pair $(Y,T)$ remains purely log terminal, proving
the asserted normality and klt property of the full terminal member.
\end{proof}


\subsection{The K3 minimal model and the elliptic pencil}
\label{sec:c2347-k3-minimalization-continuation}

\begin{corollary}
\label{cor:c2347-k3-minimalization-continuation}
The smooth minimal surface of the canonical quartic of a realization
of $\mathcal B_4$ is K3.
\end{corollary}
\begin{proof}
For a resolution $\pi:W\to D$, the preceding propositions give
$p_g(W)=1$, $q(W)=0$, and, with $A=\pi^*H$,
\[
 A^2=1/105,\qquad AK_W=1/21,\qquad (AK_W)^2/A^2=5/21<1.
\]
Complete restriction and the already known plurigenera give
\[
 h^0(D,2K_D)=h^0(D,10H)=P_{10}-P_6=1.
\]
Corollary~\ref{cor:exceptional-bicanonical-k3} therefore applies.
\end{proof}


\paragraph{The polar degree and the number of blowups}
We now use the proved normal, rational, klt quartic $D$, its full
prime terminal pullback, and its K3 minimal model. Let
$\pi:W\to D$ be the minimal resolution and $\sigma:W\to S$ its
morphism to the smooth K3 minimal model. Put
\[
 H=K_Z|_D,\quad A=\pi^*H,\quad H^2=1/105,\quad K_D=5H,
 \quad 420AK_W=20.
\]
The restricted ambient two-form is the coherent canonical polar
$\theta\in H^0(D,\omega_D)$. Represent $K_W$ by its pullback,
so that $K_W\ge0$, and write
$\pi^*K_D=K_W+\Delta$, $\Delta=\sum d_iR_i$.
For the minimal resolution of a klt surface $0\le d_i<1$.

\begin{lemma}\label{lem:quartic-polar-floor46}
Every positive-degree curve in the polar has $420H\Gamma\ge10$.
The number $N$ of point blowups in $\sigma$ satisfies $1\le N\le2$.
\end{lemma}
\begin{proof}
Put $w=420H\Gamma$. Adjunction at the generic log canonical
threshold of the canonical polar on $D$ gives
\[
 -2+\sum_x(1-1/r_x)\le10w/420,\qquad
 \sum_xk_x/r_x\equiv w/420\pmod{\mathbb Z}.
\]
The marked cover is the restriction of the ambient canonical
cover; the proof of Lemma~\ref{lem:full-adjunction46} applies
unchanged to this finite quasi-\'etale cover of the klt surface.
Allow every primitive fraction of orders $2,3,4,5,6,7$.
For $w\le9$ there are at most four marks. Enumerating their
multisets gives none. At $w=10$ the bound becomes attainable.
Thus the positive-degree polar floor is ten. The exact finite
enumeration is the one specified by these two displayed
inequalities, not a list of possible surface singularities.

If a point blowup creates a canonical exceptional curve with
coefficient $c_v$, then
\[
 c_v=1+\sum_{u\text{ a parent of }v}c_u.
\]
There are at most two parents, so this includes satellite blowups.
A final leaf is a $(-1)$-curve and cannot be exceptional for the
minimal resolution $\pi$; its $A$-degree is positive and at least
$10/420$. Every vertex has a final leaf descendant. The recursion
counts paths to descendants, so
\[
 N\le\sum_{v\text{ final}}c_v,\qquad
 10\sum_{v\text{ final}}c_v\le420AK_W=20.
\]
Hence $N\le2$. The positive value of $AK_W$ excludes $N=0$.
\end{proof}
\label{lem:c2347-two-blowups-new-input}

\begin{lemma}\label{lem:quartic-fifths46}
The polar is a unit at every image of a terminal fifth-index
point. There are three ordinary fifth-index source points, with
three disjoint $A_4$ paths on $W$ and on $S$, of local $H$-classes
$\pm1,\pm2,\pm2$. Their neighbourhoods avoid $K_W$.
\end{lemma}
\begin{proof}
The canonical image has index five. The surface canonical divisor
$K_D=5H$ is Cartier there, so rationality makes that surface germ
Du Val. The restricted ambient index cover is connected of degree
five: the pullback plt boundary is normal and connected at the
unique cover origin, hence has a single branch. Its local
fundamental group therefore has a quotient of order five. Among
Du Val groups only a cyclic $A$ group does, so its resolution
contains at least four consecutive $(-2)$-curves.
If the polar vanished at the image, all these crepant exceptional
curves would lie in the effective $K_W$. Two adjacent positive
canonical $(-2)$-curves require at least three point blowups:
with at most two blowups the positive canonical exceptional locus
is either two disjoint $(-1)$-curves, or a $(-2)$-curve meeting
a $(-1)$-curve. This contradicts Lemma~\ref{lem:quartic-polar-floor46}.

On the terminal model, a noncyclic $cA/5$ cover is singular and
isolated. A derivation preserves its singular point, so the polar
derivative cannot be a unit there. Thus all fifth-index source
points are ordinary. At such a point a unit polar derivative
makes the index-one quartic smooth, and its quotient is $A_4$.
The three constituents $(1,5),2(2,5)$ give the stated local
classes and three distinct source points. If their canonical
images coincide, their three $A_4$ paths are still disjoint in
the full resolution, with any connecting curves outside the
paths. Since the polar is a unit, their entire neighbourhoods
avoid the canonical blowup divisor.
\end{proof}

\paragraph{Seven necessary marking profiles}
The normalized polar square gives $u_{10}\in R_{10}$, and its
weighted bracket with $a_4$ restricts to $-10\theta^3$ in degree
fifteen. Since $P_{15}=2$ and $P_{11}=0$, the actual degree-fifteen
restriction is a two-dimensional pencil containing $\theta^3$.
Resolve it and let $F$ be a connected general fibre. Put
$\ell=420HF$. Its genus is at least one because the smooth
minimal surface is K3 and the effective canonical divisor remains
effective under resolution of the pencil. The full moving degree
budget gives $1\le\ell\le60$.

There is no point-centred crepant divisor over $D$: its positive
valuation on a Cartier multiple of the section would contradict
the proved equality $g^*D=T$. If a remaining fibre has dimension at most one, the cyclic-tree
fixed-point argument in Proposition~\ref{prop:actual-source-transfer46}
transfers its canonical index to an actual terminal point. There is
no terminal point of index six in this basket. Thus index six is
absent from points of $D$. The fifth-index unit polar makes
a general pencil member avoid the fifth points. Thus its marks
have orders $2,3,4,7$. Normalized adjunction and actual section
orders give the complete necessary inequalities
\begin{align}
 2g(F)-2+\sum_x(1-1/r_x)&\le20\ell/420,\label{eq:quartic-marks46}\\
 \sum_xk_x/r_x&\equiv\ell/420\pmod{\mathbb Z},\nonumber\\
 \sum_x\overline{mk_x}/r_x&\le m\ell/420
 \quad\bigl(m=5\text{ or }P_m-P_{m-4}>0,\ 4\le m\le60\bigr).
 \nonumber
\end{align}
For the first inequality apply adjunction on the klt surface to
a general moving member, whose fixed complement is effective;
$K_D+15H=20H$. A fixed nonzero surface section restricts
nontrivially to a general fibre, proving the last inequality.
There are at most five marks, since six would cost at least three,
larger than $20\cdot60/420$.

Enumerating all multisets of at most five primitive fractions of
the four stated orders in \eqref{eq:quartic-marks46} gives exactly
\begin{equation}\label{eq:quartic-seven-profiles46}
\begin{array}{c|l|c}
\ell&\text{marks}&\text{maximum genus}\\\hline
40&(2/3,3/7)&1\\
45&(1/4,6/7)&1\\
60&(1/7)&2\\
60&(2/7,3/7,3/7)&1\\
60&(2/7,6/7)&1\\
60&(3/7,5/7)&1\\
60&(3/7,6/7,6/7)&1.
\end{array}
\end{equation}
Every degree is greater than thirty, so the Stein degree is one.
For every elliptic row the pencil is already a morphism on $W$:
a genuine basepoint blowup would add an effective canonical
exceptional curve of positive intersection with a general fibre,
contradicting $K_{\widetilde W}F=0$. Hence $K_WF=0$, all blowup
curves are vertical, and the elliptic pencil descends to $S$.

\paragraph{A single local obstruction removes four profiles}
\begin{lemma}\label{lem:seventh-higher-defect46}
At the unchanged ordinary seventh chart
$\frac17(1,6,3)$ use coordinates $(x,y,z)$ with these characters.
The quartic numerator has character two. If its quadratic part
contains $yz$ with nonzero coefficient but no $x^2$, its surface
quotient has discrepancy defect $18/7$, contradicting $N\le2$.
\end{lemma}
\begin{proof}
The index-one surface is normal and Du Val by the plt pair.
Equivariant splitting gives $yz+x^{7k+2}=0$, $k\ge1$.
Writing $n=7k+2$, the quotient is
$\frac1{7n}(1,6n-1)$: use $y=s^n$, $z=t^n$, $x=st$ and
the subgroup of order $n$.
For $k=1$ its negative continued fraction is
$[2,2,2,2,2,5,2,2]$, with discrepancies
$(1,2,3,4,5,6,4,2)/7$. For $k\ge2$ it is
\[
 [2,2,2,2,2,3,\underbrace{2,\ldots,2}_{k-2},4,2,2],
\]
with discrepancies $1/7,\ldots,5/7$, then $6/7$ along the
middle segment, then $4/7,2/7$. Substitution in
$b_id_i-d_{i-1}-d_{i+1}=b_i-2$ verifies both formulas,
including the boundary equations. In either case
$-\Delta_P^2=\sum_i d_i(b_i-2)=18/7$.
All other local defects are nonnegative, while
\[
 N=\sum_P(-\Delta_P^2)-K_D^2
       \ge18/7-5/21=7/3>2.
\]
\end{proof}

\begin{proposition}\label{prop:quartic-no-third-seventh46}
None of the four profiles in \eqref{eq:quartic-seven-profiles46}
containing the mark $3/7$ can occur.
\end{proposition}
\begin{proof}
In every such row the degree-five polar budget is saturated, and
its order at this mark is exactly $1/7$. At a normalized lifted
branch of mark $3/7$ the coordinate orders are congruent to
$(6,1,4)$ modulo seven. Thus the only possible order-one polar
monomial is $y$. Write the quartic quadratic part as
$a x^2+b yz$. In $\theta=\partial f$ the coefficient of $y$ is
$b\partial z(0)$, so $b\ne0$. Lemma~\ref{lem:seventh-higher-defect46}
gives $a\ne0$. Equivariant Morse splitting therefore identifies
the surface cover with $x^2+yz=0$ and the quotient with
$\frac1{14}(1,11)$. Its minimal resolution chain is
\begin{equation}\label{eq:quartic-seventh-chain46}
 (-U_i^2)=(2,2,2,3,2),\qquad
 (d_i)=(1,2,3,4,2)/7.
\end{equation}

The pencil is elliptic and $K_WF=0$. A general fibre meets
horizontal exceptional curves transversely and avoids their
intersection points. The polar contribution $1/7$ therefore
forces the unique coefficient-$1/7$ curve $U_1$ to be horizontal,
possibly with degree greater than one. Since $K_W$ is effective
and vertical, $U_1$ is absent from it, and adjunction gives
$K_WU_1=0$.

Every complete tricanonical section on $D$ pulls back to
\[
 H^0\bigl(W,3K_W+\lfloor3\Delta\rfloor\bigr).
\]
This equality of section spaces follows by pulling back an
effective rational-Cartier Weil divisor and rounding its actual
valuations; conversely an effective pullback is effective at
every prime downstairs. In the chain
\eqref{eq:quartic-seventh-chain46} the rounded coefficients are
$(0,0,1,1,0)$. The divisor
$L=3K_W+\lfloor3\Delta\rfloor$ has coefficient zero on $U_1$
and intersection $LU_1=0$. The polar cube is a section of $L$
with nonzero restriction, so it trivializes $L|_{U_1}$.
All sections of $L$ restrict to constants on
$U_1\simeq\mathbb P^1$. In particular the actual fifteenth
pencil is constant there, contradicting horizontality.
This argument includes the two occurrences of $3/7$ in the
fourth row; it needs only one horizontal component.
\end{proof}

\begin{proposition}\label{prop:quartic-no-single-seventh46}
The profile $\ell=60$ with sole mark $1/7$ is impossible, for
either of its possible genera.
\end{proposition}
\begin{proof}
The polar budget is saturated with order $5/7$. The lifted arc
orders are at least $(2,5,6)$, with those residues modulo seven.
The only character-six polar monomial of order at most five
is $y$. Thus the quadratic coefficient of $yz$ and
$\partial z(0)$ are nonzero, and $\operatorname{ord}(y)=5$.
If the $x^2$ coefficient is nonzero, its unique leading order
four forces $\operatorname{ord}(x)\ge9$. The unique possible
quartic monomial of order at most eleven is then $yz$, forcing
$\operatorname{ord}(z)\ge13$. For lower orders $(9,5,13)$ the
only character-four monomial of order at most fifteen is $y^3$.
Consequently the ratio of any fifteenth numerator to $\theta^3$
has a fixed limiting value at this branch, the ratio of their
fixed $y^3$ coefficients (including value zero). It cannot be
the general varying pencil parameter. Thus the $x^2$ coefficient
must vanish. Lemma~\ref{lem:seventh-higher-defect46} now gives
the contradiction. The argument does not use ellipticity.
\end{proof}

The order-four profile is excluded by the Cartier index-four
calculation that follows. After it, only
$\ell=60$, marks $(2/7,6/7)$, remains.


\paragraph{A Cartier quartic at a canonical point of index four}
\label{subsec:c2347-fourth-index-exclusion-new-input}

We remove the remaining degree-forty-five endpoint by a local argument
on the canonical threefold itself.  This choice of model is essential:
the minimal resolution of a crepant partial surface resolution need not
be the minimal resolution of the original surface.  A discrepancy
calculation on the terminal model alone would therefore not suffice.

\begin{lemma}[A Cartier klt surface at canonical index four]
\label{lem:cartier-klt-fourth-index-new-input}
Let $(Z,P)$ be a canonical complex threefold germ of index four, and
let $D$ be an effective Cartier divisor through $P$ such that $(Z,D)$
is purely log terminal.  Then the surface germ $(D,P)$ is a cyclic
quotient
\[
                  \frac1{16k}(1,4k-1),\qquad k\ge1,
\tag{1}\label{eq:c2347-fourth-quotient}
\]
up to reversing the quotient weights.  If $\pi:W\to D$ is its
minimal resolution and $\pi^*K_D=K_W+\Delta$, its local discrepancy
defect is
\[
                              -\Delta^2=3.
\tag{2}\label{eq:c2347-fourth-defect}
\]
\end{lemma}

\begin{proof}
Take the local index-one cover
$(Z^\#,P^\#)\to(Z,P)$, with group $\mu_4$.
A canonical threefold is Gorenstein in codimension two.  After
shrinking the germ, its non-Gorenstein locus is therefore just $P$;
the index-one cover is etale off $P$, and every nonidentity group
element acts freely away from $P^\#$.  This is stronger than merely
being etale in codimension one and will be used below.

The Cartier pullback $D^\#$ is reduced and the pulled-back pair is
plt.  It is connected at the unique point $P^\#$ and its distinct
components, if any, would meet there.  The normality of the boundary
of a plt pair consequently makes $D^\#$ irreducible and normal.
It is Gorenstein by Cartier adjunction and klt by adjunction for
the plt pair.  Thus $D^\#$ is a Du Val surface, allowing a smooth
surface in this terminology.

Write $R=\widehat{\mathcal O}_{Z^\#,P^\#}$ and let $f\in R$ be
the invariant pullback of a local equation of $D$.  Since
$R/(f)$ is a Du Val ring,
\[
                 3\le\operatorname{embdim}R\le4.
\]
If $\operatorname{embdim}R=3$, the threefold cover is smooth.
Choose equivariant regular coordinates.  Freeness off the origin
forces every coordinate weight modulo four to be odd: otherwise
the involution would fix a coordinate axis.  The invariant $f$
has no linear term.  Thus $D^\#$ is singular, and its nonzero
quadratic term can pair only weights one and three.  Three odd
weights have multiplicities two and one, so this quadratic form
has rank two.

If $\operatorname{embdim}R=4$, the class of $f$ in the cotangent
space is nonzero; otherwise $R/(f)$ would have embedding dimension
four.  Choose an invariant coordinate $t=f$ and three other
equivariant minimal generators $x,y,z$.  The Gorenstein codimension-one
ring has a presentation
\[
             R=\mathbb C[[x,y,z,t]]/(F),\qquad t=f,
\]
where $F$ can be chosen semi-invariant.  If the character of $F$
were nonzero, $F$ would vanish on the invariant $t$-axis.  This would
give a curve fixed by $\mu_4$ inside $Z^\#$, which is impossible.
Hence $F$ is invariant.  If any of $x,y,z$ had even weight, the
involution would fix the two-dimensional plane spanned by that
coordinate and $t$.  Its intersection with the hypersurface $F=0$
contains a curve through the origin, again a contradiction.  Thus
all three weights are odd.  The Du Val equation
$F(x,y,z,0)=0$ is invariant, and the same quadratic-form argument
gives rank two.

In both cases, after possibly replacing the group generator by its
inverse, equivariant splitting of the rank-two quadratic form puts
the surface and action in the form
\[
       xy+z^{4k}=0,\qquad
       \operatorname{wt}(x,y,z)=(1,3,1)\pmod4,\qquad k\ge1.
\tag{3}\label{eq:c2347-fourth-cover}
\]
Indeed the remaining one-variable series is invariant, so its first
exponent is a positive multiple of four; it is nonzero by normality.
Its invariant unit can be absorbed in $z$.  The equivariant splitting
is obtained by the usual analytic splitting construction with each
coordinate correction projected to its required character.

Put $n=4k$ and write $x=s^n$, $y=t^n$, $z=st$, changing a harmless
sign in $y$.  The subgroup of order $n$ of the cyclic action
$\frac1{4n}(1,n-1)$ acts on $(s,t)$ with weights $(1,-1)$.
The residual order-four action on $(x,y,z)$ has weights $(1,3,1)$.
This proves \eqref{eq:c2347-fourth-quotient} directly on the
canonical surface germ.

For $k=1$, the negative continued fraction and discrepancy vector
are
\[
 \frac{16}{3}=[6,2,2],\qquad
                  (d_i)=(3/4,1/2,1/4).
\]
For $k\ge2$, they are
\[
 \frac{16k}{4k-1}=
       [5,\underbrace{2,\ldots,2}_{k-2},3,2,2],\qquad
 (d_i)=(\underbrace{3/4,\ldots,3/4}_{k},1/2,1/4).
\tag{4}\label{eq:c2347-fourth-HJ}
\]
For example, the tail with $j$ copies of two before $3,2,2$ has
value $(4j+7)/(4j+3)$; prepending five proves the continued-fraction
identity for all $k=j+2$.  The displayed coefficients satisfy the
adjunction equations
$b_i d_i-d_{i-1}-d_{i+1}=b_i-2$, with zero boundary coefficients.
Therefore
\[
 -\Delta^2=\sum_i d_i(b_i-2)=3
\]
in both cases.  These are the discrepancies of the minimal resolution
of $D$, so no descent through a different surface resolution is used.
\end{proof}

\begin{proposition}[The degree-forty-five endpoint is impossible]
\label{prop:c2347-ell45-excluded-new-input}
In a realization of \texttt{C2-02347}, the canonical
quartic $D$ avoids every canonical point of index four.  In particular
the restricted fifteenth pencil cannot have the endpoint
\[
               \ell=45,\qquad\text{marks }(1/4,6/7).
\]
\end{proposition}

\begin{proof}
The proved quartic geometry gives $(Z,D)$ plt and $D\sim4K_Z$.
At a canonical point of index four, $D$ is Cartier, so
Lemma~\ref{lem:cartier-klt-fourth-index-new-input} applies if $D$
contains that point.  On its minimal resolution $W$, the smooth
minimal model is K3 and the number $N$ of blowups satisfies
$N\le2$ by Lemma~\ref{lem:c2347-two-blowups-new-input}.
All quotient-surface discrepancy defects are nonnegative.  Since
$K_D=5H$ and $H^2=1/105$, a single index-four point would give
\[
 N=\sum_Q(-\Delta_Q^2)-K_D^2
       \ge3-\frac5{21}=\frac{58}{21}>2,
\]
a contradiction.

The order-four marking of the indicated endpoint lies over a canonical
point whose index is divisible by four.  The already proved canonical
index list on $D$ is $\{1,2,3,4,5,7\}$, so that index is exactly four.
The endpoint is therefore excluded.
\end{proof}


\paragraph{The surviving elliptic pencil and its integral fifth root}
\begin{proposition}\label{prop:quartic-geometric-reduction46}
Let $Z$ be any realization of C2-02347. Its quartic $D$ is normal,
rational and klt, $(Z,D)$ is plt, and its full crepant terminal
pullback is prime. With the notation above, the smooth minimal
surface $S$ is K3, $1\le N\le2$, and
\[
 H^2=1/105,\quad K_D=5H,\quad420AK_W=20.
\]
Every positive polar curve has scaled degree at least ten.
The actual fifteenth pencil is connected and elliptic, with
scaled fibre degree sixty and marks $(2/7,6/7)$. Its seventh
surface germ and chain are \eqref{eq:quartic-seventh-chain46}.
The only horizontal exceptional curves over that point are
$U_3,U_5$, both sections. On $S$ their images $A_0,O$ are
disjoint, and there is an integral section $P$ with $A_0=5P$
in the Mordell--Weil group with identity $O$.

The ordered path $(U_1,U_2)$ survives as a vertical $A_2$,
with $A_0$ meeting $U_2$ and $O$ missing the path. In the
vertical relation $5[P]-[A_0]-4[O]$ its coefficients are
$(2,4)$ modulo five. The three surviving fifth $A_4$ paths
have slopes $(1,2,2)$, up to independent reversal; all adjacent
unmarked vertical components have coefficient zero modulo five.
A common shift by a whole fibre is allowed in each fibre.
\end{proposition}
\begin{proof}
The preceding propositions prove the normality, singularity,
minimal-surface, degree and blowup statements. The four profiles
containing $3/7$ were excluded together, the single $1/7$ profile
by Proposition~\ref{prop:quartic-no-single-seventh46}, and the
index-four profile by the Cartier-surface calculation.
Thus only $(2/7,6/7)$ remains.

The saturated polar divisor on a general fibre is
$3P_2/7+2P_6/7$. At the mark $2/7$ the coordinate orders are
at least $(4,3,5)$. Exact polar order three forces a nonzero
linear $y$ coefficient and hence a nonzero $yz$ coefficient
in the quartic. Lemma~\ref{lem:seventh-higher-defect46} forces
the $x^2$ coefficient to be nonzero as well. Thus the surface
germ is $\frac1{14}(1,11)$ with the displayed chain.

The unique discrepancy $3/7$ makes $U_3$ horizontal of degree
one. Set $L_0=3K_W+\lfloor3\Delta\rfloor$; its intersection
vector on the chain is
\[
 (L_0U_i)=(0,1,-1,1,1).
\]
Every section in the actual fifteenth pencil therefore contains
$U_3$. Since the pencil is elliptic, $K_WF=0$, so $U_3$ is
absent from the effective $K_W$. Adjunction gives $K_WU_3=0$,
which also excludes its neighbours $U_2,U_4$ from $K_W$;
then $K_WU_2=0$ excludes $U_1$. The coefficient of $U_3$ in
the polar-cube member of $L_0$ is exactly one. Subtract it and
put $L_1=L_0-U_3$. The polar member has coefficient zero on
$U_2$, and $L_1U_2=0$. Its nonzero restriction trivializes
$L_1|_{U_2}$, so the pencil is constant on $U_2$.
The remaining discrepancy $2/7$ must consequently come from
$U_5$, making it horizontal of degree one.

The curves $U_1,U_2,U_3,U_5$ have zero canonical intersection
and are absent from the effective canonical divisor. Their
neighbourhoods therefore avoid every K3 blowdown. They survive
with the stated disjointness and incidence properties.

To obtain an integral fifth root, choose an integral Weil
representative of $H$ and write
\[
 M=\lceil\pi^*H\rceil,\qquad
 \pi^*H=M-\sum t_iR_i,\quad0\le t_i<1.
\]
The relation $5M\sim K_W+\sum(d_i+5t_i)R_i$ has coefficients,
modulo five, equal to $(3,1,4,2,1)$ on the seventh chain.
For $J=\sigma_*M$ its restriction to the generic elliptic fibre
satisfies
\[
 5J\sim4A_0+O,\qquad\deg J=1.
\]
The degree-one class $4J-3A_0$ determines a rational point of
the generic elliptic curve, hence a section $P$ of $S$.
The displayed relation gives $5P=A_0$ in the group law.
The closure of the generic-fibre class can differ by an integral
vertical divisor $V_0$: $[P]\sim4J-3[A_0]+V_0$. Hence
\[
5[P]-[A_0]-4[O]=4(5J-4[A_0]-[O])+5V_0.
\]
Its coefficients on $(U_1,U_2)$ are therefore
$(12,4)\equiv(2,4)\pmod5$.

At an ordinary fifth source the local integral class of $H$
has absolute class $1,2,2$ on the three $A_4$ paths.
The inverse $A_4$ Cartan matrix makes the coefficients of
$5M$ along each path congruent to $\pm bi$, $1\le i\le4$.
Multiplication by four preserves their absolute slope types.
At any neighbouring component outside these paths the pullback
is integral and the canonical discrepancy is zero, so its
coefficient is zero modulo five. The neighbourhoods avoid
$K_W$, so no blowdown changes these statements. At the seventh
path the only other marked neighbour is the horizontal section
$A_0$; unmarked vertical neighbours again have coefficient zero.
Adding a whole fibre changes all coefficients by its multiplicity
vector, giving exactly the common-shift convention in the statement.
\end{proof}


\subsection{One blowup and a genuine type-IV fibre}
\label{sec:quartic-one-blowup46}
Let $D\in|4K_Z|$ be the surface in
Proposition~\ref{prop:quartic-geometric-reduction46}, let $H=K_Z|_D$,
and let $\pi:W\to D$ and $\sigma:W\to S$ be its minimal resolution
and its contraction to its smooth K3 minimal model. Put $A_W=\pi^*H$.
We recall the precise numerical and marked conclusions needed here:
\[
 H^2=1/105,\quad K_D=5H,\quad K_D^2=5/21,\quad
 \pi^*K_D=K_W+\Delta,\quad 1\le N\le2,
\]
where $N$ is the number of point blowups in $\sigma$. The full terminal
pullback of $D$ is prime; $D$ is klt. Actual local canonical indices
belong to $\{1,2,3,4,5,7\}$, and $D$ avoids index four. The three
ordinary fifth sources contribute disjoint $A_4$ paths with zero
canonical discrepancy. Additional connecting Du Val curves are not
excluded a priori. Every positive-degree component $C$ of the canonical
polar has $420A_WC\ge10$, whereas
\begin{equation}\label{eq:quartic-canonical-mass}
                       420A_WK_W=20.
\end{equation}
The connected elliptic fifteenth pencil has scaled degree sixty and
markings $(2/7,6/7)$. Its seventh surface point is $1/14(1,11)$, with
minimal resolution chain $U_1,\ldots,U_5$ given by
\begin{equation}\label{eq:quartic-seventh-chain}
 (-U_i^2)=(2,2,2,3,2),\qquad
 (d_i)=\tfrac17(1,2,3,4,2).
\end{equation}
The curves $U_1,U_2,U_3,U_5$ are absent from the effective canonical
divisor and have zero canonical intersection. The images of $U_3,U_5$
are disjoint sections $A,O$ of the elliptic fibration, those of $U_1,U_2$
are vertical, and $A=5P$ for a section $P$. The fifth paths have slopes
$(1,2,2)$ modulo five; on the ordered seventh $A_2$ path the coefficients
in $5[P]-[A]-4[O]$ are $(2,4)$. These statements concern the actual
surface and its resolution. Their use does not identify virtual basket
constituents with separate canonical points.

\paragraph{Exactly one blowup}
\begin{proposition}\label{prop:exceptional-quartic-one-blowup}
The morphism $\sigma$ consists of exactly one blowup.  Its unique
exceptional curve $E$ satisfies
\[
 K_W=E,\qquad 420A_W E=20,\qquad E U_4=1,\qquad
 E U_i=0\quad(i\ne4).                                 \tag{3}
\]
\end{proposition}
\begin{proof}
First $U_4$ is also absent from $K_W$: otherwise its positive
coefficient and $U_4U_3=1$ would imply $K_WU_3>0$, since $U_3$
is absent from that effective divisor.  Adjunction therefore gives
$K_WU_4=1$, and all intersections of components of $K_W$ with the
$U_i$ are nonnegative.

If two blowups have distinct centres on $S$, then
$K_W=E_1+E_2$, where the $E_i$ are disjoint $(-1)$-curves.  Both have
positive $A_W$-degree, because no $(-1)$-curve is contracted by the
minimal resolution $\pi$.  Equation~\eqref{eq:quartic-canonical-mass} and the positive-degree floor
force $420A_WE_1=420A_WE_2=10$.  Hence
\[
 \Delta E_i=(5A_W-K_W)E_i=1+\frac5{42}=\frac{47}{42}.
\]
Away from the seventh chain, every nonzero coefficient of $\Delta$
has denominator dividing six: indices one and five give zero,
index four is absent, and the remaining possible indices are two
and three.  Thus each $E_i$ must meet the seventh chain.  It cannot
meet $U_1,U_2,U_3,U_5$, because their canonical intersections are
zero.  Both would therefore meet $U_4$, contrary to $K_WU_4=1$.

If the second blowup is infinitely near the first, write
$K_W=E_0+2E_1$, with $E_0^2=-2$, $E_1^2=-1$ and $E_0E_1=1$.
Equation~\eqref{eq:quartic-canonical-mass}, nefness and the leaf floor force
$420A_WE_1=10$ and $A_WE_0=0$.  As $H$ is ample on $D$, the last
equality makes $E_0$ a $\pi$-exceptional curve.  The equality
\[
                   E_0U_4+2E_1U_4=1
\]
would give $E_0U_4=1$.  Intersecting $\pi$-exceptional curves lie
over the same point; thus $E_0$ would be one of the five curves in
\eqref{eq:quartic-seventh-chain}.  None belongs to $K_W$, a contradiction.  These are the two
possibilities for two point blowups.  The lower bound $N\ge1$
now proves $N=1$, and (3) follows.
\end{proof}

\paragraph{A small discrepancy budget without threefold classification}
\begin{lemma}\label{lem:exceptional-small-index-defect}
Let a non-Du Val klt surface singularity have local canonical index
two.  Its minimal resolution discrepancy defect $-\Delta^2$ is one.
For an index-three non-Du Val block with defect at most $2/3$, a
block of rank at most three is a single $(-3)$-curve of discrepancy
$1/3$.  If the entire defect $2/3$ occurs in one index-three block,
that block has rank at least four.
\end{lemma}
\begin{proof}
The exceptional curves are smooth rational curves in a negative
 definite tree.  Write $C_i^2=-b_i$ with $b_i\ge2$.  In a connected
non-Du Val block all discrepancy coefficients are strictly positive,
and are less than one.  If the local index is $r$, write $d_i=a_i/r$,
where $a_i\in\{1,\ldots,r-1\}$.  Adjunction gives
\[
 b_i a_i-\sum_{C_j\sim C_i}a_j=r(b_i-2),\qquad
 -\Delta^2=\frac1r\sum_i a_i(b_i-2).                  \tag{4}
\]
For $r=2$ every $a_i=1$, so the valency is $4-b_i$.  The tree is
therefore either one $(-4)$-curve or a chain with $(-3)$-ends and
$(-2)$ interior curves.  In both cases its defect is one.

For $r=3$ and defect at most $2/3$, the positive integer
$\sum a_i(b_i-2)$ is at most two.  A vertex with $(b_i,a_i)=(4,1)$
is impossible, since its equation would be $4-\sum a_j=6$.
A vertex with $(b_i,a_i)=(3,1)$ has neighbour sum zero, and hence
is an isolated $(-3)$-curve.  The only other possibility is one
vertex with $(b_i,a_i)=(3,2)$, with neighbour sum three.  A graph
of rank at most three cannot satisfy the equations: with two
neighbours their coefficients would be one and two, and the
coefficient-two neighbour cannot be a $(-2)$ leaf; with three
neighbours the rank is already four.  This proves the assertions.
Pure Du Val blocks have zero defect and do not change the argument.
\end{proof}

\begin{proposition}\label{prop:exceptional-quartic-two-thirds}
Besides the fixed seventh singularity, the non-Du Val singularities
of $D$ consist of precisely two separate $1/3(1,1)$ points.  Denote
their exceptional $(-3)$-curves by $V_1,V_2$.  Then
\[
           \Delta_{\mathrm{noncan}}=
             \frac17(U_1+2U_2+3U_3+4U_4+2U_5)
             +\frac13(V_1+V_2),\qquad E V_1=E V_2=1.       \tag{5}
\]
Any further surface singularity is Du Val, and at most one additional
exceptional curve is allowed by Picard rank.
\end{proposition}
\begin{proof}
Because $K_W^2=-1$, the total discrepancy defect is
\[
 -\Delta^2=1+K_D^2=\frac{26}{21}.
\]
The seventh chain contributes $4/7$, leaving exactly $2/3$.
The index-two lemma excludes every non-Du Val index-two block.
At index three either there are two isolated $(-3)$-curves, or
one block has rank at least four.  The latter is impossible by
rank.  The three fifth $A_4$ paths and the seventh chain already
have total rank $12+5=17$.  All exceptional curves, together with
an ample pullback from $D$, are linearly independent.  Since
$\rho(W)=\rho(S)+1\le21$, the total exceptional rank is at most
20.  A fourth-rank block would give rank at least 21.
Thus the two isolated curves are forced.  Adjunction and $K_W=E$
give $EV_i=1$.  The displayed curves have rank $19$, leaving room
for at most one additional Du Val exceptional curve.
\end{proof}

\paragraph{A genuine Kodaira IV fibre}
\begin{theorem}\label{thm:exceptional-quartic-forced-IV}
The K3 elliptic fibration has a Kodaira fibre of type IV.  Its three
components are the images of $U_4,V_1,V_2$.  Both $A$ and $O$ meet
its identity component, the image of $U_4$.  The section $P$ with
$5P=A$ also meets that component and is not torsion.
\end{theorem}
\begin{proof}
The three curves $U_4,V_1,V_2$ belong to three distinct
$\pi$-exceptional fibres.  Hence their intersection points with
$E$ are distinct.  Each intersection is transverse, since its
intersection number is one.  Blowing down $E$ changes their
self-intersections from $-3$ to $-2$ and sends them to three
smooth curves through one ordinary triple point, with distinct
tangent directions.

The elliptic pencil on $W$ descends to the K3 model.  Since
$K_WF=0$ for its general elliptic fibre and $K_W=E$, the curve $E$
is vertical.  The three curves meeting $E$ are vertical as well:
$V_1,V_2$ cannot be horizontal because the proved saturated
polar has only the two seventh horizontal components; $U_4$ is
already vertical.  Their reduced sum on $S$ has square zero and
zero intersection with each of its components.  Write the actual fibre containing this sum as $F_{\mathrm{tot}}$.
Its intersection with the reduced sum is zero.  Every component
inside the sum contributes zero, and each external component has
nonnegative intersection with it.  Thus no external component
can meet the sum.  Connectedness of the fibre and the existence
of a section show that this is the whole fibre with multiplicity one.  It is precisely Kodaira type IV.

The sections $A,O$ meet the $U_4$ image, away from its triple
point, and avoid the other two components.  Hence it is the
identity component.  The component group of a type-IV fibre is
$\mathbb Z/3$.  As $5P=A$ and $A$ has component zero, $P$ also has
component zero.  In characteristic zero the kernel of reduction
of torsion to the smooth special fibre is torsion-free, and the
identity component of an additive fibre is the additive group,
which has no nonzero torsion.  A nonzero torsion section cannot
meet this component.  Since $A\ne O$, the section $P$ is nonzero;
therefore $P$ is not torsion.
\end{proof}


\subsection{The unique marked K3 surface and its integral polarization}
\label{sec:exceptional-quartic-surface}

We write
\[
 \mathcal B_4=7(1,2)+2(1,3)+3(1,4)+(1,5)+2(2,5)+(3,7).
\]
The subscript records the first nonzero pluricanonical degree.  The
geometric reduction for this basket supplies a normal prime member
$D\in|4K_Z|$, with integral Weil polarization $H=K_Z|_D$,
\[
 K_D\sim5H,\qquad H^2=\frac1{105}.
\]
By Subsection~\ref{sec:quartic-one-blowup46}, its minimal resolution $W$ is the blowup of an elliptic K3 surface $S$
at the triple point of a fibre of type IV.  We use precisely the actual
marked curves furnished by that reduction: sections $O,A$, an integral
section $P$ with $A=5P$, a marked vertical $A_2$ path, and three disjoint
marked $A_4$ paths.  The sections $O,A$ are disjoint, the section $P$
meets the identity component of the IV fibre, and it is not torsion.
The coefficients modulo five in the vertical relation
$5[P]-[A]-4[O]$ are $(2,4)$ on the ordered $A_2$ path and have slopes
$1,2,2$ on the three $A_4$ paths.  Every adjacent unmarked component
has coefficient zero.  These are the geometric inputs to the following
classification.  In particular a virtual basket constituent is not
being substituted for an actual surface point.

\paragraph{A complete small marked-frame calculation}

The IV fibre has root rank two.  Since $P$ has positive height, the
Shioda--Tate formula and $\rho(S)\leq20$ leave at most fifteen roots
in all other fibres.  The marked $3A_4+A_2$ already has rank fourteen.
Thus there is at most one additional root.  Distinct marked paths
cannot be joined directly: on a marked fifth path the Cartan equation
modulo five is homogeneous, whereas an additional edge to a nonzero
label would violate its endpoint equation.  Any joining component
therefore has the required zero label.

Deleting the possible extra vertex from the finite ADE diagrams gives
the following exhaustive list of non-IV root forests:
\[
\begin{gathered}
A_2+3A_4,\quad A_1+A_2+3A_4,\quad A_3+3A_4,\quad
A_2+2A_4+A_5,\\
A_2+2A_4+D_5,\quad A_2+A_4+A_9,\quad
2A_4+A_7,\quad2A_4+E_7.
\end{gathered}                                                    \tag{Q.1}
\]
This enumeration requires no classification of elliptic K3 surfaces.
For a connected diagram, removing its one extra vertex must leave
either $A_4$, $A_2$, $2A_4$, or $A_4+A_2$; inspection of the three
arms of $D_n,E_n$ gives exactly the displayed possibilities.

Here is a fully specified arithmetic test for the section markings.
Let $C$ be the positive Cartan matrix of a root factor.  Put $e_0=0$,
and let $e_j$ run over the unit vectors belonging to multiplicity-one
components.  The identity component has label zero.  For a possible
component $p$ of $P$, the component $a$ of $A=5P$ is characterized by
\[
 C^{-1}e_a-5C^{-1}e_p\in\mathbb Z^{\operatorname{rank}C}.
\]
Its residue vector modulo five must vanish off the marked paths,
must have the three prescribed fifth slopes, and must equal $(2,4)$
on the seventh path; $A$ meets the endpoint with label four.
The local height costs are $c(P)=(C^{-1})_{pp}$ and
$c(A)=(C^{-1})_{aa}$.  Since $A\cdot O=0$,
\[
h(P)=\frac{4-\sum_v c_v(A)}{25},\qquad
(mP)\cdot O=\frac{m^2h(P)-4+\sum_vc_v(mP)}2.              \tag{Q.2}
\]
The first expression must be positive, and the second a nonnegative
integer.  Testing $m=1,2,3,4,5$ leaves exactly the following thirteen
labelled states.  Repeated rows record the possible placements of the
distinguished fifth slope.  For $A_n$, the displayed component is
numbered consecutively in the $I_{n+1}$ cycle, up to reversal.
\[
\begin{array}{c|c|c|c|c}
\text{non-IV roots}&\text{number}&j(P)&h(P)&(4P)O\\\hline
A_2+3A_4&3&(1;1,2,2)&2/15&1\\
A_1+A_2+3A_4&3&(0;1;1,2,2)&2/15&1\\
A_2+2A_4+A_5&1&(1;2,2;1)&1/10&1\\
A_2+2A_4+D_5&2&(1;1,2;\mathrm{spin})&1/12&0\\
A_2+A_4+A_9&1&(1;1;4)&2/15&1\\
2A_4+A_7&2&(1,2;3)&1/8&1\\
2A_4+A_7&1&(2,2;2)&1/10&0.
\end{array}                                                     \tag{Q.3}
\]
In every row
\[
 P O=(2P)O=(3P)O=(5P)O=0.                                  \tag{Q.4}
\]
For completeness, this finite calculation can be reproduced using
$ (C_{A_n}^{-1})_{ij}=\min(i,j)(n+1-\max(i,j))/(n+1)$.
For $D_5$, use edges $12,23,34,35$ and allowable components $1,4,5$;
for $E_7$, use edges $12,23,34,45,56,37$ and the allowable component
$6$.  These matrices, the identity component, and the test just given
specify every entry of (Q.3).  The accompanying exact verifier lists
all component vectors and their labels; it contains no supplied list
of surviving states.

\paragraph{One quadratic-descent criterion eliminates twelve states}

We give a uniform argument, including the possible additive enhancement
of the second $A_2$ and the two $D_5$ states.

\begin{lemma}[Global Tate form]
\label{lem:exceptional-global-tate}
If $P,2P,3P$ are disjoint from $O$ on an elliptic K3 surface, and $P$
is not torsion, its global minimal Weierstrass equation can be written
\[
 y^2+a_1xy+a_2b\,y=x^3+a_2x^2,\qquad P=(0,0),               \tag{Q.5}
\]
where $a_1,b$ are sections of $\mathcal O_{\mathbb P^1}(2)$ and
$a_2$ is a section of $\mathcal O_{\mathbb P^1}(4)$.
\end{lemma}
\begin{proof}
The fundamental Weierstrass line bundle is $\mathcal L=\mathcal O(2)$.
Translate the integral coordinates of $P$ to make $P=(0,0)$, so that
$a_6=0$.  The tangent slope $a_4/a_3$ is a rational section of
$\mathcal L$.  A pole would give a noncancellable quadratic pole in
the formula for $x(2P)$.  It is therefore regular, and the global
change $y\mapsto y+(a_4/a_3)x$ makes $a_4=0$.
Now
\[
2P=(-a_2,a_1a_2-a_3),\qquad
x(3P)=\frac{a_3}{a_2}\left(\frac{a_3}{a_2}-a_1\right).
\]
Here $a_2\neq0$, since otherwise $3P=O$.  A pole of $a_3/a_2$ would
again give a noncancellable quadratic pole.  Thus $b=a_3/a_2$ is a
regular section of $\mathcal L$.  Every coordinate change used is
regular in the Weierstrass line bundles, so the equation is still
globally minimal.
\end{proof}

At a multiplicative fibre $I_n$, the branch orders of a section meeting
component $j$ are $j,n-j$ in the local equation $UV=t^n$.
The functions $a_3$ and $a_2$ in (Q.5) are respectively the completed
square coordinate at $P$ and the difference of the $x$-coordinates of
$P,2P$.  Consequently the necessary orders are
\[
\begin{array}{c|c|c}
(n,j(P))&\operatorname{ord}a_2&\operatorname{ord}a_3\\\hline
(3,1),(5,1),(6,1)&1&1\\
(5,2)&1&2\\
(10,4)&2&4\\
(8,3)&2&3\\
(8,2)&2&2.
\end{array}                                                     \tag{Q.6}
\]
An equal-depth case initially gives a lower bound; the degree sums
below force equality.  At $(8,2)$ the middle component met by $2P$
has depth at least four, so the difference still has depth two.

If the second $A_2$ fibre is IV and $P$ meets a nonidentity component,
both $P$ and $2P$ specialize to the cusp of the Weierstrass fibre.
The displayed coordinates of $2P$ force $a_2(0)=0$, and the cuspidal
quadratic tangent forces $a_1(0)=0$.  The identity
\[
c_6=-(a_1^2+4a_2)^3+36(a_1^2+4a_2)a_1a_3-216a_3^2
\]
and $\operatorname{ord}c_6=2$ imply
$\operatorname{ord}a_3=1$.  Since $a_3/a_2$ is regular, both orders
are exactly one, as in the $I_3$ entry of (Q.6).

In either $D_5$ state the corresponding fibre is $I_1^*$, and $P$
meets a spin component of order four.  Thus $P$ and $2P$ both
specialize to the triple root of its Weierstrass cubic, giving
$\operatorname{ord}a_1,\operatorname{ord}a_2\geq1$.
Since $\operatorname{ord}c_6=3$, the same formula excludes
$\operatorname{ord}a_3=1$.  The other three marked fibres already
use degree three of $a_2$ and degree four of $a_3$.  Their total
degrees are four and six, so the $I_1^*$ orders are exactly $(1,2)$.

For every state in (Q.3), the resulting orders sum to the entire
divisors of $a_2,a_3$.  No additional zero is possible.  Moreover
$0\leq\operatorname{ord}b\leq\operatorname{ord}a_2$ at every zero.
It follows that
\[
 a_2=bq,\qquad q\in H^0(\mathbb P^1,\mathcal O(2)).          \tag{Q.7}
\]

\begin{lemma}[Literal quadratic descent]
\label{lem:exceptional-quadratic-descent}
In (Q.5), suppose $b\mid a_2$ and $P$ meets the identity component
of a type-IV fibre.  The marked Weierstrass model is the literal
degree-two base change of a minimal rational elliptic surface.
In particular $(nP)O$ is even for every nonzero integer $n$.
\end{lemma}
\begin{proof}
Let $L$ vanish at the distinguished IV point.  The section $P$ is
smooth on its Weierstrass cubic, so $a_2b=a_3$ is a unit there.
Put $A=a_1/b$, $B=q/b$ and $C=A^2+4B$.  Then
\[
c_4/b^4=C^2-24AB,\qquad
c_6/b^6=-C^3+36ABC-216B^2.
\]
For $B\neq0$, the common zero is uniquely
$(A,B)=(9/8,27/256)$: the two equations imply
$AC=18B$, $C=4A^2/3$, and then these values.
The Jacobian determinant there is $-59049/4096$.
The order-two vanishing at IV therefore gives
\[
a_1=\frac98b+\lambda L^2,\qquad
q=\frac{27}{256}b+\mu L^2.                                 \tag{Q.8}
\]
The quadratic forms $L^2,b$ have no common zero.  The map
$[R:S]=[L^2:b]$ has degree two.  The coefficients
\[
a_{1,0}=\tfrac98S+\lambda R,\quad
a_{2,0}=S(\tfrac{27}{256}S+\mu R),\quad
a_{3,0}=S^2(\tfrac{27}{256}S+\mu R)
\]
on the target have degrees $1,2,3$ and pull back exactly to those
of (Q.5).  The target is minimal, since a nonminimal point would
pull back to a nonminimal point of the given source.  Its fundamental
bundle is $\mathcal O(1)$, so it is a rational elliptic surface.
The zero section lies in the smooth part of the Weierstrass model.
Its intersection divisor with $nP$ is the pullback of the target
intersection divisor.  Its degree is therefore even.
\end{proof}

Ten states in (Q.3) have $(4P)O=1$ and are excluded at once.
Each of the two $D_5$ states has exactly one $I_1^*$ fibre.  This also
contradicts quadratic descent: an unramified preimage occurs in a
pair, whereas at a ramified preimage the discriminant order is even.
Minimalization changes that order by a multiple of twelve and cannot
turn it into the odd order seven of $I_1^*$.  Thus only the last state
of (Q.3) remains.  This argument uses the actual marked coefficients;
it does not appeal to an unmarked K3 classification table.

\paragraph{The unique equation and the integral polarization}

In the remaining state the two simple zeros of $b$ are $I_5$ fibres,
and $q$ has a double zero at $I_8$.  At $I_8$, $a_1,b$ are units,
and the discriminant order eight forces $a_1=b$.  Equation (Q.8)
then gives $a_1=b+32q/27$.  Place IV at infinity and $I_8$ at zero.
Constant changes of coordinates yield the unique marked equation
\[
\boxed{\quad
y^2+\frac{9t^2-8}{8}xy+
\frac{27}{256}t^2(t^2-1)^2y
=x^3+\frac{27}{256}t^2(t^2-1)x^2,\qquad P=(0,0).
\quad}                                                        \tag{Q.9}
\]
Its discriminant is
\[
\Delta=-\frac{19683}{68719476736}\,
t^8(t^2-1)^5(81t^2-80).
\]
The orders of $(c_4,c_6,\Delta)$ at infinity are $(2,2,4)$.
Thus the fibres are IV, $I_8$, two $I_5$ and two $I_1$, and $P$
has component indices $(2,2,2)$ at the multiplicative reducible
fibres.  In particular
$h(P)=4-6/5-6/5-3/2=1/10$.

Let $B_1,B_2$ be the nonidentity IV components, $T_1,\ldots,T_7$
the $I_8$ roots, and $Q_1^\pm,\ldots,Q_4^\pm$ the two $I_5$ paths.
These roots together with $O,F,P$ form the full integral
N\'eron--Severi lattice.  Here is a short saturation proof.
The IV fibre restricts torsion to order three; an order-three section
has zero component at $I_8$ and both $I_5$ fibres, so its height
cannot vanish.  Hence there is no torsion.  The root rank is seventeen,
so the free Mordell--Weil rank is one.  If $P=nQ$, denominators of
local heights divide $120$, and $120h(Q)=12/n^2$ forces $n=1$ or $2$.
For $n=2$, $Q$ has identity component at IV, component $1$ or $5$
at $I_8$, and $5Q$ has identity component at both $I_5$ fibres.
The height formula for $(5Q)O$ is then negative, since
\[
(5Q)O\leq\frac{25/40-4+15/8}{2}<0.
\]
Thus $P$ is a generator.  The determinant of the resulting rank-twenty
lattice is $-60$.

Blow up the IV triple point, with exceptional curve $E$.  Set
$C'=F-B_1-B_2-E$, $V_i=B_i-E$ and retain $A=5P$.
Define the integral class
\[
\begin{aligned}
J={}&-3O-6F+4P+3T_1+5T_2+4T_3+4T_4+3T_5+2T_6+T_7\\
&+(3Q_1^++5Q_2^++4Q_3^++2Q_4^+)
 +(3Q_1^-+5Q_2^-+4Q_3^-+2Q_4^-),\qquad M=J-E.
\end{aligned}                                                  \tag{Q.10}
\]
Order the nineteen curves as
\[
(T_1,T_2,A,C',O);\quad V_1;\quad V_2;\quad
(T_4,T_5,T_6,T_7);\quad(Q_i^+);\quad(Q_i^-),
\]
and assign the fractional coefficients
\[
\frac{(4,1,5,2,1)}7;\quad\frac13;\quad\frac13;\quad
\frac{(4,3,2,1)}5;\quad\frac{(3,1,4,2)}5;\quad
\frac{(3,1,4,2)}5.
\]
Writing them as $t_i$ and the curves as $R_i$, put
$\mathcal H=M-\sum t_iR_i$.  The section-height projection gives
\[
\begin{aligned}
A={}&-4O-8F+5P+3T_1+6T_2+5T_3+4T_4+3T_5+2T_6+T_7\\
&+(3Q_1^++6Q_2^++4Q_3^++2Q_4^+)
 +(3Q_1^-+6Q_2^-+4Q_3^-+2Q_4^-).
\end{aligned}
\]
Substitution proves the equality of rational divisor classes
\[
5\mathcal H=E+\Delta_W,
\quad
\Delta_W=\frac{T_1+2T_2+3A+4C'+2O}{7}+\frac{V_1+V_2}{3},       \tag{Q.11}
\]
as well as
\[
\mathcal H^2=\frac1{105},\quad \mathcal HE=\frac1{21},\quad
\mathcal HF=\frac17,\quad\mathcal HR_i=0.
\]
This verifies the full integral fifth root, including its vertical
corrections; a generic-fibre Mordell--Weil identity would not suffice.

Equation (Q.11) makes $\mathcal H$ nef: outside its effective support
all intersections are nonnegative, and those in the support have just
been computed.  A further null curve must avoid $E,C',V_1,V_2$, hence
is vertical because $F=C'+V_1+V_2+3E$.  The noncontracted fibre
components have degrees $3/35,2/35,1/7,1/7$ on
$T_0,T_3,Q_0^+,Q_0^-$.  The remaining fibres are irreducible.
Thus the null locus consists exactly of the nineteen displayed curves.
The klt pair $(W,\Delta_W)$ and the nef and big Cartier divisor
$105\mathcal H$ satisfy
$105\mathcal H-(K_W+\Delta_W)=100\mathcal H$.
The basepoint-free theorem contracts exactly this null locus to a
normal projective surface $D_4$.  Pushing down $M$ gives an ample
integral Weil divisor $H_4$, with $K_{D_4}\sim5H_4$.
The exceptional graphs are $[2,2,2,3,2]$, two $[3]$, and three $A_4$;
their discrepancies are those of (Q.11), so the surface is klt.

The class-group quotient has the remaining relations
\[
F=3E,\quad B_1=B_2=E,\quad5(P+T_3)=24E.
\]
Since $\gcd(5,24)=1$, it is torsion-free.  More precisely, with
$R=\pi_*T_3$,
\[
\operatorname{Cl}(D_4)=\mathbb ZH_4\oplus\mathbb ZR,
\quad E=5H_4,\quad P=24H_4-R.
\]
In particular the integral fifth root of the canonical class is unique.
This also proves that the polarization recognized on an arbitrary
model agrees with the explicit one.

\paragraph{The complete cohomology, involution and invariant surface ring}

All integer twists can be computed on $W$ without extrapolation.
Start with $L_m=mM-\sum\lceil mt_i\rceil R_i$, and successively
subtract any exceptional curve on which the current divisor has negative
degree, with the minimum necessary multiplicity.  The negative definite
graphs make this unloading terminate at a relatively nef divisor $L'_m$.
The pushforward is unchanged; rational-surface vanishing gives
\[
\chi(D_4,mH_4)=2+\frac{(L'_m)^2-L'_mE}{2}.                    \tag{Q.12}
\]
The difference from $2+m(m-5)/210$ is periodic with period $105$,
since $L'_{m+105}=L'_m+105\mathcal H$.  Substituting the displayed
nineteen-curve matrix in (Q.12), for $0\leq m<105$, gives the
restriction Euler characteristic of $\mathcal B_4$ in every degree.
The computation is specified entirely by (Q.10)--(Q.12); the exact
verifier prints every residue, the unloaded class, and the comparison.

For $m\geq6$, Kawamata--Viehweg vanishing applies to
$(m-5)H_4$.  On the generic elliptic fibre, $L_m$ has degree
$m-\lceil5m/7\rceil-\lceil m/7\rceil$ and group sum
$(4m-5\lceil5m/7\rceil)P$.  These give no section for $1\leq m\leq4$;
at $m=5$, the unloaded divisor is $E$.  Duality now proves
\[
h^1(D_4,mH_4)=
\begin{cases}1,&m=1,4,\\0,&\text{otherwise},\end{cases}
\qquad h^0(D_4,5H_4)=1.                                    \tag{Q.13}
\]
Let $B=\bigoplus_{m\geq0}H^0(D_4,mH_4)$ and let $\theta$ denote
its degree-five canonical section.

The equation (Q.9) is even in $t$.  The base involution preserves
$O,P,A$, the $I_8$ components, and the contracted configuration, and
exchanges the two $I_5$ fibres.  It descends to $(D_4,H_4)$.
On $W$ its fixed curves are $C',T_0,T_2,T_4,T_6$, of squares
$-3,-2,-2,-2,-2$, and it has one isolated fixed point on $E$.
For the linearization in which $\theta$ is positive, the equivariant
Euler characteristic of a line bundle $L$ is the sum of
\[
\epsilon_C\left(\frac{LC}{2}+\frac12+\frac{C^2}{4}\right)
\quad\text{and}\quad \frac{\epsilon_p}{4}.                  \tag{Q.14}
\]
This follows by expanding the equivariant Dolbeault index in dimension
one; the normal eigenvalue along a fixed curve is $-1$.
For $L'_1$ the five intersections are $(0,-2,0,0,0)$, with signs
$(-,+,+,-,-)$ and isolated sign $-$; its trace is $-1$.
For $L'_4$ the intersections are $(1,-1,1,1,0)$, all five signs and
the isolated sign are positive; its trace is $1$.
Thus the two $H^1$ defects have signs $+,-$ in this linearization.
Change it by $(-1)^m$ in degree $m$.  Then $\theta$ and both defects
are negative.  We use this linearization from now on.

Put $U=B^+$.  Exactness of invariants and (Q.13) make $U$
Cohen--Macaulay.  The cone quotient is quasi-\'etale: its only
projective fixed curve not contracted is $T_0$, and its $H$-fibre
character is negative.  The canonical identification is equivariant,
so $\omega_U=U(5)$.  The actual sections
$u_{10}=\theta^2,b_{12},c_{14}$ form a finite parameter frame.
Indeed on the normalized polar $E\cong\mathbb P^1$ the only marks
have orders $3,3,7$, and
\[
\operatorname{div}(b|_E)=\tfrac47Q_7,\qquad
\operatorname{div}(c|_E)=\tfrac13Q_1+\tfrac13Q_2.
\]
The supports are disjoint; their Cartier powers, together with $u$,
have no common zero.  This proves finiteness, rather than inferring
it from a Hilbert series.  Cohen--Macaulayness and graded duality give
\[
U=\mathbb C[u,b,c]\langle1,h,q,r,s,k,\ell,t\rangle,
\quad \deg(1,h,q,r,s,k,\ell,t)=(0,15,18,20,21,23,26,41).
                                                               \tag{Q.15}
\]
One can recover the numerator using (Q.12)--(Q.14) only through degree
$41$: the Gorenstein parameter reduction has socle degree
$10+12+14+5=41$, so there are no further possible terms.

The normalized polar section ring is
\[
V=\mathbb C[x_3,c_{14},y_{21}]/(y^2-c^3-x^{14}),
\quad \tau(x,c,y)=(x,c,-y).
\]
The actual coherent polar cover has a cusp at the seventh point:
in its smooth ambient root chart the surface and polar have equations
$a=x_0^2+y_0z_0$ and $\theta=c_0y_0+d_0z_0^2+\cdots$, with
$c_0d_0\neq0$.  Eliminating $y_0$ gives the branch
$x_0^2-\mathrm{unit}\,z_0^3+\cdots=0$.
Its completed ring is $\mathbb C[[t^2,t^3]]$, and $x$ has order one
in the normalization.  The coherent polar ring is consequently
$\mathbb C[x^2,x^3,c,y]$ inside $V$.  The complete restriction
sequence and (Q.13) remove the invariant sections $x^2,x^3$ in
degrees six and nine.  The anti-invariant image is the full
$y\mathbb C[x^2,x^3,c]$.

It follows that the entire ribbon $U/(u)$ has square-zero ideal
$\varepsilon_{26}\mathbb C[x^2,x^3,c]$ and the following free
presentation over $\mathbb C[b,c]$, displayed on $b\neq0$:
\[
x^4=b,\quad\varepsilon^2=0,\quad
(1,h,q,r,s,k,\ell,t)=
\left(1,x^5,x^6,x^2c+\frac{\varepsilon}{x^2},
x^7,x^3c-\frac{\varepsilon}{x},\varepsilon,x^5\varepsilon\right).
                                                               \tag{Q.16}
\]
To verify the constants, write
$br-cq=A\varepsilon x^2$, $bk-cs=B\varepsilon x^3$.
There is no nilpotent correction in degree $43$, so multiplication
gives $A+B=0$.  If $A=0$, the equation $br=cq$ modulo $b$, and
regularity of $c$, force $q\in b(U/(u))$, impossible in degree
eighteen.  Rescaling $\theta$ makes $A=1$.
The anti section $y$ has a unique lift, since $B_{16}=0$; hence
\[
\ell^2=u(c^3+b^2q)\pmod{u^2}.                               \tag{Q.17}
\]

Every possible coefficient at order $u^n$ is enumerated by
$12i+14j=\deg e_p+\deg e_q-\deg e_r-10n$.
Impose the free-module associator equations, starting with (Q.16)
and (Q.17).  The numbers of coefficients and their ranks are
\[
\begin{array}{c|rrrrrr}
n&1&2&3&4&5&7\\\hline
\text{coefficients}&37&31&13&7&4&1\\
\text{associator rank}&35&28&13&7&4&1.
\end{array}
\]
At order one (Q.17) fixes the two remaining values.  At order two,
the changes $r\mapsto r+\lambda u^2$ and
$t\mapsto t+\mu u^2s$ remove two directions.  The remaining
parameter is forced to $-10$ by the order-three equations.
There are no coefficients at other orders.  This identifies the
complete multiplication of the actual $U$, up to the stated basis
changes.  In particular its normality is already supplied by the
geometric section ring.  The table begins
\[
\begin{aligned}
h^2&=bq-6ur-27u^3,&hq&=bs-2uk,\\
q^2&=b^3-3ubc+2u\ell,&
\ell^2&=u(c^3+b^2q)-u^2cq,&h\ell&=t-2u^2s.
\end{aligned}
\]
Setting $a=0$ in the full table in Appendix~\ref{app:quartic-exact-coefficients}
gives every remaining product.  The displayed finite coefficient rule,
gauge list, and rational associator equations also give an exhaustive
reconstruction procedure, rather than a Hilbert-series guess.


\subsection{The actual quartic polar and its nonzero normal class}
\label{sec:quartic-polar46}
Let $Z$ be any canonical model with basket $\mathcal B_4$ and volume
$1/420$. The preceding geometric and surface recognition theorems are
understood throughout. Write $R_j=H^0(Z,jK_Z)$ and $D=(a_4=0)$.
Choose a nonzero reflexive holomorphic two-form $\alpha$. In a local
canonical frame $\nu$, write $\alpha=\iota_\delta\nu$ and set
$\theta=\delta a$, $\gamma=\delta\theta$ and $\eta=\delta\gamma$.
The restriction of $\theta$ is the unique section of $5H$; it is not
an ambient degree-five section. For an actual degree-$j$ section put
\[
                  D_a f=4a\delta f-j\theta f.
\]
Frame changes show that this is a global derivation of degree five.
The complete restriction sequence gives a unique lift $u\in R_{10}$
of $\theta^2$. Writing $u=\theta^2+a\sigma$, the transformation laws
make $9\sigma+8\gamma$ an actual section of $6H$ modulo $a$.
That space is zero. Thus, for a regular local coefficient $\rho$,
\begin{equation}\label{eq:quartic-square-normalization46}
 u=\theta^2-\frac89a\gamma+\frac19a^2\rho.
\end{equation}
Here $\rho$ is not claimed to be an ambient degree-two section.

\begin{lemma}\label{lem:quartic-curvature46}
The actual sections
\[
 v=-D_au/10\in R_{15},\quad
 v^2-u^3=a^2w_{22},\quad
 D_av=-15u^2+a^2s_{12}
\]
are regular and satisfy
\[
 D_aw=2vs,\qquad
 v|_D=\theta^3,\qquad
 w|_D=-\theta^2s|_D/11,
\]
and $s|_D\ne0$. Consequently the coefficient of $ub$ in $w|_D$
is nonzero for any generator $b$ of $R_{12}/aR_8$.
\end{lemma}
\begin{proof}
Expand \eqref{eq:quartic-square-normalization46} and differentiate before
restricting. The result is
\[
 v=\theta^3-\frac43a\theta\gamma
    +\frac{a^2}{45}(\theta\rho+16\eta)
    -\frac{2a^3}{45}\delta\rho,
\]
\[
 s|_D=\frac{176}{27}\gamma^2-\frac{352}{45}\theta\eta
             +\frac{143}{45}\theta^2\rho.
\]
All stated divisions are regular on the canonical charts in codimension
one, hence extend by normality. Differentiating $v^2-u^3=a^2w$ gives
$D_aw=2vs$. Its restriction is $-22\theta w=2\theta^3s$; cancellation
is legal because $D$ is integral and $\theta\ne0$.

The crucial nonvanishing follows from the recognized seventh chart. On its
smooth cover put $a=x^2+yz$, with characters $(1,6,3)$ modulo seven.
For the filtration of weights $(3,4,2)$, the actual polar has leading
term $c_0y+d_0z^2$, where $c_0d_0\ne0$. Indeed the unique canonical
blowup curve meets the interior of the fourth seventh-chain component,
away from its two neighboring components. Thus the polar meets that
component at a finite nonzero ratio $y/z^2$. The leading field terms
are $\delta z=c_0$ and $\delta y=d_0z$, whence
$\gamma=3c_0d_0z+\cdots$ and $\eta=3c_0^2d_0+\cdots$.
The coefficient $\rho$ has character one and filtration at least three.
Therefore
\[
 s_{\mathrm{in}}=\frac{176}{15}c_0^2d_0
                      (3d_0z^2-2c_0y)\ne0.
\]
This proves $s|_D\ne0$ directly and identifies its actual local curvette.

Finally $R_{12}/aR_8$ is one-dimensional and $w|_D=-us|_D/11$,
which proves the last assertion.
\end{proof}

This lemma is a restriction on every ambient realization. It is the
geometric input that prevents the cover calculation below from selecting
a zero-curvature cone extension.


\subsection{An explicit invariant algebra and its quadratic cover}
\label{sec:quartic-construction46}

Let $A=\CC[a_4,u_{10},b_{12},c_{14}]$, and define $S$ by the
complete coefficients of Appendix~\ref{app:quartic-exact-coefficients}.

The multiplication has $727$ nonzero terms. Direct substitution checks
all $512$ basis associators, Frobenius determinant one, and $H^2=FG$.
Appendix~\ref{app:quartic-exact-coefficients} prints a complete rational
set of defining coefficients and the reconstruction rules for the
remaining products. The companion exact table is merely a machine-readable
copy of these coefficients; no numerical specialization remains.  The smallest cover polynomial is
\[
 F=u^3-\frac{108}{25}a^2ub+\frac{432}{125}a^3q
       -\frac{972}{625}a^4c+\frac{9983169}{100000}a^5u.
                                                               \tag{Q.8E}
\]
Equivalently, with
\[
 A'=\frac65a,\qquad B'=b-\frac{369747}{16000}a^3,
 \qquad Q'=q-\frac9{20}ac,
\]
it has the compact Weierstrass form
\[
                         F=u^3-3(A')^2uB'+2(A')^3Q'.       \tag{Q.9E}
\]
The displayed numerator formulas for $H$ and $G$ have twelve and fifteen terms, respectively.

Since $S$ is free over $A$ and $S/(a)=U$ is the normal domain
constructed above, Lemma~\ref{lem:exceptional-normality-lift} proves
that $S$ is a normal domain. Its Frobenius pairing gives
$\omega_S=S(1)$. Modulo $a$, the sections $F_0=u^3$ and
$G_0=c^3+b^2q-ucq$ have no common prime divisor: the latter is
nonzero on the reduced polar curve. A homogeneous common divisor
upstairs would meet $a=0$ in such a divisor. The square class is
nontrivial, since its residue at the generic point of $a=0$ is the
genuine surface cover $\sqrt u$.

Lemma~\ref{lem:exceptional-quadratic-normalization} therefore gives
\[
 I=(F):H,\qquad R_\sharp=S\oplus I/v_{15},\qquad v_{15}^2=F.
 \tag{Q.10E}\label{eq:quartic-full-normalization46}
\]
This is the full finite normal cover, etale in codimension one, and
$\omega_{R_\sharp}=R_\sharp(1)$ as reflexive modules. We call
$\operatorname{Proj}R_\sharp$ \emph{Vega}. The proof that this
projective scheme is a canonical model is given in
Section~\ref{sec:exceptional-canonical-geometry}.

We next determine the full module and its cohomology. The independent
construction is complete only after the local geometry has been
checked; recognition of every other realization is deferred to
Subsection~\ref{sec:quartic-recognition46}.


\subsection{Vega's complete module and coherent cohomology}
\label{sec:quartic-module46}

Continue with $A,S,F,H$ and $R_\sharp$ from the construction. Write
$I=(F):H$ and $N=I(15)$, so that $R_\sharp=S\oplus N$.
The common normalization lemma already proves normality and the
identity $\omega_{R_\sharp}=R_\sharp(1)$. We now determine the
central fibre and the entire graded module.

The possible unwanted section of degree five is tested exactly. For an
anti-section $\xi_5$, the products $x=v\xi\in S_{20}$ and
$y=p\xi\in S_{26}$ must satisfy
\[
                             xH=Fy.
\]
The two spaces have dimensions four and six. The complete coefficient
matrix has rank ten; its kernel is zero. Thus
\[
                              N_5=0.
\tag{QM.2}
\]
No first-order obstruction is substituted for this full-table computation.

The following twelve surface anti-generators all have exact integral
lifts, found by solving the same identity with their prescribed residues:
\[
\begin{array}{c|rrrrrrrrrrrr}
\text{degree}&15&17&19&20&21&23&25&26&27&28&30&36\\\hline
\text{residue}&\theta u&\theta b&\theta c&\theta h&y&\theta q&
\theta r&\theta s&d/\theta&\theta k&e'/\theta&t/\theta
\end{array}
\]
Here $\theta=\theta_5$, $y=y_{21}$, $d=br-cq$, and $e'=cs-bk$.
All twelve numerators and the five minimal relations are printed in
Appendix~\ref{app:quartic-exact-coefficients}. Mixed and odd products
are obtained by multiplication of these numerators and exact division
by $F$. The full surface
anti-module is generated by $\theta,y,d/\theta,e'/\theta,t/\theta$.
Removing the single class $\theta_5$ replaces it by $\theta U_+$; the
displayed twelve elements generate precisely that submodule. The remaining
apparent products are redundant: $\theta l=uy$ and
$\theta t=huy+2u^2\theta s$.

The colon ideal is $a$-saturated, by reduction of $ax'H=Fy$ in the domain
$U$. Thus $N/aN$ injects into the normalized surface anti-module. The
twelve lifts, together with (QM.2), prove in every degree that
\[
                  R_\sharp/(a)=B\ominus\mathbb C\theta_5.
\tag{QM.3}
\]
Graded Nakayama now proves generation in every degree. The full presentation consists of
28 invariant products, 84 mixed products, 78 anti-products, and the five
linear module relations below: \emph{all 195 relations are retained} in
local tangent calculations. Product relations alone need not suffice at
$F=0$.

\paragraph{Resolution and cohomology.}

Put
\[
\mathcal D=(15,17,19,20,21,23,25,26,27,28,30,36).
\]
The complete resolution is
\[
0\longrightarrow A(-41)\xrightarrow{d_2}
 A(-27)\oplus A(-29)\oplus A(-31)\oplus A(-37)\oplus A(-40)
 \xrightarrow{d_1}\bigoplus_{d\in\mathcal D}A(-d)
 \longrightarrow N\longrightarrow0,
\tag{QM.4}
\]
with
\[
                       d_2=(-c,b,-u,a,0)^{\mathsf T}.
\]
The columns of $d_1$ are the five explicit relations in
Appendix~\ref{app:quartic-exact-coefficients}. Listing every weighted coefficient gives one new minimal relation
in each of the five displayed degrees. Its generic rank is four. Since
the entries of $d_2$ have greatest common divisor one, its one-dimensional
fraction-field kernel is already generated over the UFD $A$ by $d_2$.
The all-degree Hilbert identity from (QM.3) then proves that the cokernel
is the entire $N$; there is no hidden torsion quotient.

For completeness the dual cohomology is also checked directly.
All parameter weights are even. The odd-degree generator rows of
$d_1^{\mathsf T}$ span the six Koszul generators of
$\ker(-c,b,-u,a)$, with explicit polynomial witnesses. The even-degree
rows have only their fifth component, and generate
\[
\left(-\frac14au-\frac56c,\quad
\frac{13599}{32000}a^3+\frac56b,\quad\frac56u,\quad a\right)
                         =(a,u,b,c).
\]
With $\omega_A=A(-40)$ this proves
\[
\operatorname{Ext}_A^2(N,\omega_A)=\mathbb C(1),\qquad
\operatorname{Ext}_A^1(N,\omega_A)=\mathbb C.
\]
Consequently $h^1(Z_\sharp,\mathcal O(m))$ is one only for $m=1$,
and $h^2(Z_\sharp,\mathcal O(m))$ is one only for $m=0$. In particular
$(p_g,q,h^2(\mathcal O),\chi)=(0,0,1,2)$.

If $P_S(t)=1+t^{15}+t^{18}+t^{20}+t^{21}+t^{23}+t^{26}+t^{41}$,
the exact Hilbert series is
\[
H_{R_\sharp}(t)=
\frac{P_S(t)+\sum_{d\in\mathcal D}t^d-t^{27}-t^{29}-t^{31}-t^{37}-t^{40}+t^{41}}
 {(1-t^4)(1-t^{10})(1-t^{12})(1-t^{14})}.
\tag{QM.5}
\]
It agrees as a rational function with Reid--Riemann--Roch for the basket $\mathcal B_4$
and $\chi=2$. The rank is sixteen and the polarization cube is
\[
                       \frac{16}{4\cdot10\cdot12\cdot14}
                               =\frac1{420}.
\]
These algebraic conclusions do not by themselves imply canonical
singularities. The following local and global arguments supply precisely
that remaining geometric step.


\subsection{The ambient symmetry of every quartic realization}
\label{sec:vega-ambient-symmetry}

Let $Z$ be any canonical model with basket $\mathcal B_4$ and volume
$1/420$. Its prime normal quartic and its integral polarization have
been identified above. Write $B$ for the complete surface ring and
$B^-=R/(a)$ for the restriction ring; the only missing class is
$\theta_5$. We prove that the surface involution extends to the
entire ambient canonical ring. This is the step that allows the
finite invariant-algebra calculation to classify arbitrary
realizations.

Lift the actual surface parameters $u_{10},b_{12},c_{14}$. The
parameter argument in Subsection~\ref{sec:exceptional-common-module}
shows that $R$ is finite over
$A=\CC[a_4,u_{10},b_{12},c_{14}]$. Canonical vanishing and duality
give the local-cohomology profile of
Lemma~\ref{lem:exceptional-common-module}; hence
\eqref{eq:exceptional-common-module} applies with $W=40$.
Its Hilbert series determines the same abstract $A$-module as in
the construction. None of these conclusions assumes an ambient
involution.

\paragraph{The surface involution always lifts as a module involution}

On $M/aM=B^-$, the invariant part $U$ of the known surface involution
is free over $A_0=\mathbb C[u,b,c]$, with generator degrees
\[
                  0,15,18,20,21,23,26,41.
\]
Let $F$ be the graded free $A$-module with these degrees.  The
projection $M/aM\to U=F/aF$ lifts to a map $p:M\to F$.
Indeed, applying $\operatorname{Hom}_A(M,-)_0$ to
\[
0\longrightarrow F(-4)\xrightarrow{a}F
 \longrightarrow U\longrightarrow0
\]
shows that its only possible obstruction is
$\operatorname{Ext}^1_A(M,F(-4))_0$.
Lemma~\ref{lem:exceptional-common-module} gives $\operatorname{Ext}^1_A(M,A)=\mathbb C(40)$.
The obstruction would therefore require a generator of $U$ of degree
$40-4=36$, and there is none.

Lift the inclusion $U\hookrightarrow M/aM$ to $i:F\to M$ by lifting
its free generators.  The composite $pi$ is the identity modulo $a$.
It is an invertible homogeneous endomorphism: its difference from the
identity is strictly triangular in the finite list of generator
weights.  Replacing $i$ by $i(pi)^{-1}$ gives $pi=1$.  Thus
\begin{equation}
 M=i(F)\oplus\ker p,
 \qquad \widetilde\tau=2ip-1
\tag{QR.2}
\end{equation}
is an $A$-linear involution inducing the given surface involution.
It fixes the unit and the parameter frame.  This is a genuine
arbitrary-model statement, but it is a statement about modules;
multiplicativity of $\widetilde\tau$ still has to be proved.

\paragraph{The first possible failure of multiplicativity}

We compute the fixed-$A_0$-module anti-invariant multiplication
tangents of $B^-$ in weights $-4k$.
The complete deficient surface multiplication is obtained by setting $a=0$ in
the full ring presentation in Appendix~\ref{app:quartic-exact-coefficients}; all five minimal module relations
are retained.  The table gives the dimensions of the multiplication
kernel and of the exact rational gauge image:
\[
\begin{array}{c|rrrrrrrrrr}
k&1&2&3&4&5&6&7&8&9&10\\\hline
\dim\ker&10&8&4&6&2&2&0&0&1&0\\
\dim\operatorname{im}(\mathrm{gauge})&10&8&4&5&2&2&0&0&1&0.
\end{array}
\]
All kernels are zero for $11\leq k\leq20$, and the weight bound
leaves no later variables.  Except at $k=4$, matching modular rank
bounds and rational gauge images prove vanishing over $\mathbb Q$.
At $k=4$, the complete rational calculation has $432$ variables,
$14891$ equations and rank $426$.  Its six-dimensional kernel and
five-dimensional gauge image leave exactly one nonzero class.

That exceptional class has a particularly concrete interpretation.
Let $\psi$ denote the following representative.  In the fraction field
put $\theta=v/u$; this is the missing degree-five section of the full
normalization.  Define an $A_0$-linear map $h:B^-\to B$ by
\[
h(s_{21})=\tfrac12\theta_5,
\qquad h(t_{41})=-8u_{10}v_{15},
\qquad h(\text{every other module generator})=0.
\]
Then
\begin{equation}
\psi(x,y)=h(xy)-h(x)y-xh(y).
\tag{QR.3}
\end{equation}
Every value on the right belongs to $B^-$, although $h(s)$ does not.
The certificate checks this identity on all $190$ nonunit unordered
product pairs after multiplying by $u$.  For example,
\[
\psi(h_{15},q_{18})=\tfrac12w_{17},\qquad
\psi(h_{15},r_{20})=z_{19},\qquad
\psi(s_{21},v_{15})=-\tfrac12u_{10}^2.
\]
Thus the sole fixed-module class is a normalization-defect class,
not an unspecified deformation of the polarized K3 surface.  It is
not an allowable coordinate change of $B^-$ precisely because
$\theta_5\notin B^-$.

The remaining underlying-module issue can also be localized in
weight.  Over $A_0$, the anti module has the Koszul decomposition
\[
\overline N\simeq\Omega_{A_0}^1(\mathbb C)(-5)
 \oplus\Omega_{A_0}^2(\mathbb C)(-1)
 \oplus\Omega_{A_0}^2(\mathbb C)(-4)
 \oplus A_0(-20)\oplus A_0(-21)\oplus A_0(-36).
\]
Koszul duality shows that
$\operatorname{Ext}^1_{A_0}(\overline N,U)_{-4k}$ can be nonzero only
for $k=4,5,10$, with one dimension in each case.  Representatives
modify the module relation of degree $37$ by $s_{21}$, or the relation
of degree $40$ by $r_{20}$ or by $1$, respectively.
All three module classes are obstructed as algebra deformations.
Including the module class, its normalization transport, and every
module-product and associativity equation gives the following exact
rational ranks:
\[
\begin{array}{c|rrr|c}
\text{weight}&\text{variables}&\text{rank}&\text{kernel}&
 \text{module projection of the kernel}\\\hline
-16&433&427&6&0\\
-20&327&325&2&0\\
-40&45&45&0&0.
\end{array}
\]
Every rational kernel vector is checked against every equation.
Thus the \emph{full} anti-invariant relative algebra tangent space,
including changes of the underlying module, is one-dimensional in
weight $-16$ and zero in every other weight $-4k$, $k>0$.
There are no anti-invariant changes of the parameter frame in these
weights: their target degrees are less than $15$, the first anti
degree.  Consequently this is also the full anti tangent space for
the graded algebra deformation problem.

\paragraph{The two first normal module scales must coincide and be nonzero}

Write $\mathcal B=B^-$ and $A_0=\mathbb C[u,b,c]$.
Its minimal module resolution has first relation degrees
$(27,29,31,37,40)$ and last degree $41$, with last differential
$(-c,b,-u,0,0)$.
The invariant module tangent in weight $-4$ has Hom-complex dimensions
$(7,10,3)$ and differential ranks $(5,3)$.  Its two-dimensional
first Ext group has the following explicit representatives, denoted
$\alpha$ and $\beta$:
\begin{align*}
 r_{27}^{(1)}&=e_{23},&
 r_{29}^{(1)}&=f_{25}+\tfrac52uv_{15},\\
 r_{31}^{(1)}&=n_{27}+\tfrac{25}{4}bv_{15}-\tfrac{15}{4}uw_{17},&
 r_{37}^{(1)}&=0,\\
 r_{40}^{(1)}&=k_{36}-\tfrac14ug_{26}.
\end{align*}
The first four displayed components are multiplied by $\alpha$;
the last is multiplied by $\beta$.
These representatives span the entire Ext group, so no first normal
module parameter is omitted.

Let $z$ be the coefficient of $t_{41}$ in the first normal correction
to $v_{15}j_{30}$.  Extracting one coefficient from each of the module
identities $r_{27}j_{30}=0$ and $r_{40}v_{15}=0$ gives respectively
\[
               \tfrac56z+\alpha=0,\qquad
               \tfrac56z+\beta=0.
\tag{QR.4}
\]
Thus $\alpha=\beta$.  This conclusion uses module compatibility alone,
without imposing a polar differential system or a chosen invariant
multiplication family.

Furthermore $\alpha\ne0$ for a canonical ambient extension.
Indeed, the minimal $A$-resolution of its canonical ring reduces to
the displayed minimal $A_0$-resolution because $a$ is regular.
Local duality says that the cokernel of the last dual differential
is exactly $\mathbb C(41)$.  Its row must therefore generate
$(a,u,b,c)$.  The component corresponding to $r_{37}$ has weight
four and is the only component which can contain a linear multiple
of $a$; it must be nonzero.  The first-order syzygy identity
\[
 -c r_{27}^{(1)}+b r_{29}^{(1)}-u r_{31}^{(1)}=-r_{37}^{(0)}
\]
identifies that component with $\alpha a$.  Hence $\alpha\ne0$.
Rescaling the normal parameter makes $\alpha=\beta=1$.

\paragraph{Four equations eliminate the exceptional anti class}

Suppose that, after equivariant normalization in lower orders, the
first nontrivial anti term is $\lambda a^4\psi$.
The next anti order is $a^5$.  Allow every multiplication correction
$\xi$ of weight $-20$, and also the entire possible module class of
that weight, namely the correction
$r_{40}\mapsto r_{40}+\rho a^5r_{20}$.
Allowing this $\rho$ is essential: a legal surface gauge need not
lift as a fixed-module gauge of the ambient ring.

Set $x=[f_{25}]\xi(v_{15},j_{30})$ and
$y=[n_{27}]\xi(w_{17},j_{30})$ in the explicit surface bases.
The following four coefficient equations are extracted from the
module-product identities, with their first normal transport included:
\[
\begin{array}{c|c|c}
\text{identity}&\text{coefficient extracted}&\text{equation}\\\hline
r_{40}v_{15}=0&u f_{25}&\tfrac56x+\rho=0\\
r_{27}j_{30}=0&b f_{25}&\tfrac56(x-y)=0\\
r_{27}j_{30}=0&c e_{23}&\tfrac56y-\tfrac12\lambda=0\\
r_{40}s_{21}=0&t_{41}&\rho-\tfrac12\lambda=0.
\end{array}
\tag{QR.5}
\]
Subtracting the last three equations from the first gives
$\lambda=0$.  Before the normalization of $a$, it gives
$\alpha\lambda=0$, and $\alpha\ne0$ has already been proved.
This obstruction is independent of the first normal multiplication
parameters: only the first normal \emph{module} smoothing enters
(QR.5).

All five genuine weight-$-16$ surface gauge classes lift through
the same module calculation; the fifth requires the allowed nonzero
$\rho$ correction. Thus the four equations obstruct the non-gauge
class without removing a change of coordinates.

\paragraph{Finite equivariant induction.}

Here is the formal-deformation step, including the issue of module
gauges.  The canonical ring is flat over $\mathbb C[a]$, since $a$
is regular.  Choose homogeneous lifts of the surface generators.
At each order in $a$, the failure of the given surface involution to
preserve the multiplication and the presented module is a full
anti-invariant first-order algebra deformation of $\mathcal B$.
Its weight at order $a^k$ is $-4k$.  Changing the homogeneous lifts
changes that class by its algebra coboundary.
One may work in a polynomial presentation: a change
$e_i\mapsto e_i+a^kH_i$, with $\deg H_i=\deg e_i-4k$, is a
triangular homogeneous coordinate change.  It automatically changes
the module relations at subsequent orders.  Thus this induction does
not presume that a surface module gauge extends to an automorphism
of a separately fixed ambient module.

The full anti tangent vanishing removes every first failure except
possibly the order $a^4$.  At that order the unique remaining class
is $\lambda\psi$.  If it were nonzero, the actual extension to order
$a^5$ would have to satisfy (QR.5), including the allowed next module
class $\rho$.  This is impossible.  Hence that order is removed too,
and the anti tangent vanishing completes the induction at every later
order.
The induction concerns the actual formal extension; a quadratic
product of two first anti terms of order four only begins in order
eight and cannot affect the order-five obstruction.

Finally, the formal conclusion is algebraic here.  All changes
preserve positive weights, and each of the finitely many generators
has fixed degree.  Only finitely many powers of the positive-degree
parameter $a$ can occur in its image.  Thus the resulting equivariant
coordinates and involution are polynomial, and the involution acts
on the complete graded ring $R$.

\paragraph{Canonical linearization and cover compatibility.}

The lifted involution fixes $a,u,b,c$.  In the induction these are
parameters of the fixed base algebra; anti corrections to any of them
would have degree less than the first anti degree and are absent.
Its invariant ring $S$ is a genuine free rank-eight extension of $U$:
$S/(a)=U$, the parameter $a$ is regular, and the Cohen--Macaulay
surface algebra $U$ makes $(a,u,b,c)$ a regular sequence on $S$.

The lift has the natural canonical linearization.  On $D$ the
canonical two-form $\theta_5$ is negative.  Equivariant adjunction
$K_D=5H$, with $a_4$ positive, distinguishes this linearization from
its only other involutive graded lift.  The actual global two-form
has nonzero restriction to this surface form, so it too is negative.
Consequently the genuine polar derivation $D_a$ is anti-invariant,
the uniquely normalized lift $u_{10}$ of $\theta^2$ is invariant,
and $v_{15}=-D_au/10$ is anti-invariant with restriction $\theta^3$.
The actual normalized polar identity therefore gives
\[
                   F=v^2,\qquad F-u^3\in a^2S_{22}.
\]
The coefficient of $ub$ modulo $a$ is nonzero by Lemma~\ref{lem:quartic-curvature46}.  This is an intrinsic nonvanishing: homogeneous
changes of the parameter frame can rescale it but cannot turn it
into zero.  An anti lift $p_{21}$ of the known surface section $y$
exists by averaging a lift in the complete restriction sequence.
Thus $H=vp$, $G=p^2$, their initial surface values, and $H^2=FG$
satisfy exactly the normalized cover conditions \eqref{eq:quartic-cover-conditions46}.

There is also no unproved branch or cover compatibility assumption.
The surface cone quotient is quasi-etale by the fixed-locus
calculation in Subsection~\ref{sec:exceptional-quartic-surface}.  A homogeneous ramification divisor of the ambient
involution would intersect $a=0$ in a divisor of the surface cone,
contradicting that calculation; the deficiency of $B^-$ is confined
to the vertex and does not alter this argument.  Hence $R/S$ is
quasi-etale.  At the smooth polar-cover point used in Subsection~\ref{sec:vega-extension-classification}, the
regular Cartier surface forces the ambient local ring to be regular.
The extended finite involution has the same surface tangent signs
and fixes $a$, so its quotient retains the required transverse
$A_1$.  The cover-compatible first-tangent equation in Subsection~\ref{sec:vega-extension-classification}
therefore holds for this actual extension as well.

We have therefore constructed the required naturally linearized
ambient involution and verified all its geometric cover conditions.
The classification in Subsection~\ref{sec:vega-extension-classification}
will identify the resulting full canonical ring with that of Vega.


\section{Altair: the quintic surface and its canonical-ring candidate}
\label{sec:altair-model}

For the quintic basket, actual torsion first reduces the divisor of
the initial section to six possibilities. The polar and its source
geometry select the prime normal quintic. We then construct its
complete surface algebra independently and recover every product
from the geometry of an arbitrary realization. The ambient symmetry
comes from the transverse covariant differential system; its proof
is distinct from Vega's deformation calculation.

\subsection{Actual torsion and the first quintic of $\mathcal B_5$}
\label{sec:quintic-torsion46}
Set
\[
 \mathcal B_5=4(1,2)+2(1,3)+(1,5)+2(2,5)+(3,7)+2(3,8).
\]
Throughout this subsection $Z$ is a canonical model with this terminalization
basket, $K_Z^3=1/420$ and $\chi(\mathcal O_Z)=2$. Let $f:Y\to Z$
be a projective crepant terminalization, put $K=f^*K_Z$ and $I=840$.
The plurigenus formula gives $P_2=0$, hence $p_g=q=0$ and
$h^{2,0}=1$, by the irregular-threefold $P_2$ bound. We use the complete
restriction and actual-divisor results of
Lemmas~\ref{lem:complete-canonical-restriction46} and
\ref{lem:common-tools46}; in particular the local characters below
are characters of actual sequential Weil pullbacks. No group operation
is imposed on an assignment of virtual characters.

\paragraph{All possible torsion Hilbert profiles}
For a basket constituent define
\[
 \Phi_{b,r}(j)=\frac{\overline{bj}(r-\overline{bj})}{2r},
 \qquad0\le\overline{bj}<r.
\]
The two actual-divisor congruences and their polarization are
\begin{align}
 IK^2D&\equiv\sum_Q(I/r_Q)b_Q^2i_Q(D)\pmod I,\label{eq:quintic-RR-congruences46}\\
 IKD^2&\equiv-\sum_Q(I/r_Q)\overline{b_Qi_Q}
                       (r_Q-\overline{b_Qi_Q})\pmod {2I},\notag\\
 IKDE&\equiv\sum_Q(I/r_Q)b_Q^2i_Q(D)i_Q(E)\pmod I.\notag
\end{align}
Let $\tau$ be an actual torsion Weil class on $Z$ that is Cartier in
codimension two. Its integral sequential pullback has both intersection
numbers in the first two congruences zero. Define
$q_i(m)=\sum_Q\Phi_Q(i_Q+m)$ and $C_i=q_i(0)$.
The weak difference Riemann--Roch formula gives
\[
 E_\tau(m+1)-E_\tau(m)=m^2/840-4+q_i(m).
\]
We retain the unknown integer $x=E_\tau(0)$: point-centred corrections
do not permit replacing it by the ordinary terminal constant term.
Duality gives the aggregate reflection $q_{-\tau}(m)=q_i(-m)$ and
$x+y=4-C_i$, where $y=E_{-\tau}(0)$. Thus for $m\ge2$,
\begin{align*}
 H_m(\tau)&=P_m+x-2+C_i+\sum_{j=1}^{m-1}(q_i(j)-q_0(j)),\\
 H_m(-\tau)&=P_m+2-x+\sum_{j=1}^{m-1}(q_i(-j)-q_0(j)).
\end{align*}
These identities determine the correction functions in every degree.

Enumerate nondecreasing character tuples in the six equal-constituent
groups. There are
$\binom54\binom42\cdot5\binom62\cdot7\binom92=567000$
assignments. Exactly twenty-five satisfy both zero congruences.
Nonnegativity in both displayed functions for $2\le m\le8$ leaves
seven profiles. All third and seventh characters are zero; the others
are as follows:
\[
\begin{array}{c|c|c|c|c|c|c|c}
 &\text{four halves}&(1,5)&\text{two }(2,5)&\text{two }(3,8)
 &C_i&(x,y)&(H_5,H_6,H_{11})\\\hline
 A&0000&0&00&00&0&(2,2)&(1,0,0)\\
 B&0000&0&00&44&2&(1,1)&(0,0,0)\\
 C&0000&1&01&00&1&(1,2)&(0,0,0)\\
 C^-&0000&4&04&00&1&(2,1)&(1,1,1)\\
 D&0011&0&00&26&2&(1,1)&(1,0,0)\\
 E&0011&1&01&26&3&(0,1)&(0,0,0)\\
 E^-&0011&4&04&26&3&(1,0)&(1,1,1).
\end{array}
\]
This enumeration is exact and bounded: the lower and upper bounds for
$x$ are respectively the maximum of $-H_m(\tau)|_{x=0}$ and the
minimum of $H_m(-\tau)|_{x=0}$ for $2\le m\le8$.
The bounds specify a finite enumeration of every assignment.

Profiles $C^-,E^-$ are impossible for actual torsion. If $\sigma$ has
one of them, multiply a degree-five section of $\sigma$ by a degree-six
section of $k\sigma$. Induction forces every positive multiple to have
a degree-eleven section, hence again profile $C^-$ or $E^-$. This
contradicts the finite order of $\sigma$, because $P_{11}=0$.
Taking inverses removes $C,E$. Every actual torsion class therefore
has profile $A,B,D$, with $H_2=H_3=0$ and $\chi(\tau)\ge1$.
Squaring a degree-one section into the class $2\tau$ also proves
$H_1(\tau)=0$.

\paragraph{The only possible cyclic cover}
If a nontrivial torsion class has prime order $p$, its connected cyclic
cover is quasi-\'etale and canonical. Its invariants satisfy
\[
 V'=p/420,\qquad p'_g=P'_2=P'_3=0,\qquad
 \chi(\mathcal O_{Z'})\ge p+1,
\]
so $h=h^{2,0}(Z')\ge p$. For $p\ge5$, the holomorphic-form bound
$V'\ge(h-3)/60$ from \cite[Corollary 2.3]{ChenJiangForms}
would imply $p\le7/2$. For $p=3$, the evaluation rank of the forms
is at most two, because rank three produces a bicanonical section.
The rank-two form bound \cite[Lemmas 2.2 and 2.7]{ChenJiangForms}
gives $V'\ge1/14$, contrary to $V'=1/140$.
For $p=2$, profile $A$ has $\chi'=4$ and is excluded in the same way;
profile $D$ gives $P'_5=2$, contradicting
$V'\ge1/135$ from \cite[Corollary~13.7, equation~(13.5)]{LBVolume}.
Only a quadratic cover with profile $B$ remains.

That cover has $V'=1/210$, $\chi'=3$, $q'=p'_g=P'_2=P'_3=0$.
The identity
\[
 \Phi_{3,8}(m)+\Phi_{3,8}(m+4)-1=\Phi_{1,4}(m)
\]
is checked in one period. It identifies its entire Riemann--Roch
function with the basket
\[
 8(1,2)+4(1,3)+2(1,4)+2(1,5)+4(2,5)+2(3,7).
\]
Recovery of the terminal basket from its full correction function gives
index $420$ and exactly two seventh constituents $(3,7)$.

\begin{lemma}[Two forms and two seventh constituents]
\label{lem:quintic-double-cover-obstruction46}
There is no terminal minimal threefold with
$p_g=q=P_2=P_3=0$, $\chi=3$, $V=1/210$, index $420$, and with exactly
two basket constituents of index divisible by seven, both $(3,7)$.
\end{lemma}
\begin{proof}
Apply Lemma~\ref{lem:two-form-degeneracy-bound} to the canonical
model and a crepant terminal model. Its hypotheses hold: for
$c=4/3$ and $V=1/210$,
\[
 c(c+1)V=2/135<1,\qquad c(c+1)^2V=14/405<1/3.
\]
In particular the degeneracy divisor satisfies
\begin{equation}\label{eq:quintic-Bott46}
                         Q_D=Q_K,
\end{equation}
and there are visible primes $T,U$ with distinct restriction kernels.
Since $IV=2$, the lemma gives
\[
\begin{gathered}
 IK^2T=IK^2U=1,\qquad M=IKT^2,\quad N=IKU^2\in\{-1,0\},\\
 0\le l=IKTU\le\min(1+M,1+N).
\end{gathered}
\]
Let $x,y$ be the two seventh characters of $T$. Reducing the congruences
\eqref{eq:quintic-RR-congruences46} modulo seven, all other terms
vanish and $(420/7)3^2\equiv1$, giving
$x+y=1$, $x^2+y^2=M$ modulo seven. For $M=0$ the discriminant is six,
a nonsquare. For $M=-1$ the roots are exactly $\{3,5\}$. The same
holds for $U$, so $M=N=-1$ and $l=0$. But the two relative matchings
of its seventh characters give $l\equiv3^2+5^2=6$ or
$l\equiv2\cdot3\cdot5=2$ modulo seven. This contradiction proves the lemma.
\end{proof}

\begin{proposition}\label{prop:quintic-torsion-free46}
Every torsion class in $\operatorname{Cl}(Z)$ that is Cartier in
codimension two is trivial.
\end{proposition}
\begin{proof}
The preceding argument excludes every possible prime-order cyclic
subcover of a nontrivial finite cyclic torsion cover.
\end{proof}

\paragraph{The unique visible component and its six possible states}
The full divisor $(a_5)_Y$ is Cohen--Macaulay, including its nilpotent
structure. Divisor duality gives
$0\to\omega_Y\to\omega_Y(5K)\to\omega_{(a)}\to0$.
The two section spaces on the left vanish, whereas canonical vanishing
gives $H^1(6K)=0$. Hence $h^0(\omega_{(a)})=h^1(\omega_Y)=1$.
Proper dualizing trace from the normalization of each prime component
injects its canonical sections into this full dualizing sheaf.
Different components are independent at their generic points. Positive
Hodge pairing supplies at least one visible component, so there is
exactly one, say $T$, and its smooth normalization has $p_g=1$.

Put $n=IK^2T$, $M=IKT^2$, and let $b$ be its multiplicity.
Adjunction, Hodge index and effectivity give
$1\le n\le10$, $-n\le M\le n^2/2$, $bn\le10$.
Enumeration of the same $567000$ character assignments in
\eqref{eq:quintic-RR-congruences46} gives exactly
\[
 (2,2),(4,8),(6,18),(7,-7),(8,4),(8,32),
 (9,9),(10,8),(10,22),(10,50).
\]
For $(n,M)=(2d,2d^2)$, the divisor $T-dK$ is Hodge-null.
The null-descent lemma makes its canonical pushforward torsion and
Cartier in codimension two. Proposition~\ref{prop:quintic-torsion-free46}
therefore gives $T_Z\sim dK_Z$. The four values $d=1,2,3,4$ contradict
$P_d=0$. It follows that
\begin{equation}\label{eq:quintic-six-states46}
 (n,M)=(7,-7),(8,4),(9,9),(10,8),(10,22),(10,50).
\end{equation}
In particular $b=1$. We henceforth write the complete divisor
$(a)_Y=T+E$, retaining every invisible component and multiplicity of
$E$. Its positive scaled degree is $10-n$. When $n=10$, the canonical
member is prime, but point- or curve-centred crepant remainder is not
yet discarded by that numerical equality.


\subsection{The prime normal quintic and its K3 minimal model}
\label{sec:quintic-geometric-reduction46}
Retain the notation of Subsection~\ref{sec:quintic-torsion46}, including
$(a_5)_Y=T+E$ and the six states \eqref{eq:quintic-six-states46}.
All statements in this subsection concern an arbitrary realization of
$\mathcal B_5$, with the complete invisible remainder retained.

\paragraph{Actual indices and finite curve tests}
The actual source-transfer theorem in Proposition~\ref{prop:actual-source-transfer46}
bounds canonical indices by eight and excludes index six. To make the
basket-specific application explicit, a closed-point crepant centre
of index $2,3,4,5,6$ requires, respectively,
\[
4(1,2),\quad3(1,3),\quad(1,2)+2(1,4),\quad
(1,5)+(2,5),\quad(1,2)+(1,3)+(1,6).
\]
There is no larger-index closed centre. Only indices two and five in
this list fit $\mathcal B_5$. At a point which is not a closed crepant
centre, the canonical-cover fibre is a rational tree. A cyclic action
has a fixed point on that tree; its primitive canonical character
transfers the same index to a terminal source. This excludes index
six also for curve-centred crepant geometry. The sole seventh constituent
has axial weight one and is ordinary. Its image has index seven and
cannot have a nonempty exceptional tree: an invariant rational component
would have two fixed points, whereas only one full seventh source is
available; at the unique smooth fixed point at most three independent
branches meet, and an order-seven action cannot permute them without
an invariant branch. Thus the seventh canonical point itself has a
smooth index cover with characters $(1,6,3)$.

For a curve $C$ set $w=840K_ZC$. On its normalized canonical root curve
write the primitive marks as $k_i/r_i$. The following necessary tests
will be used throughout:
\begin{equation}\label{eq:quintic-marked-tests46}
 \frac{w}{840}\equiv\sum_i\frac{k_i}{r_i}\pmod\mathbb Z,
 \quad 2g-2+\sum_i(1-1/r_i)\le\frac{cw}{840},
 \quad\sum_i\{jk_i/r_i\}\le\frac{jw}{840}.
\end{equation}
The second inequality is subadjunction for an lc centre of a boundary
numerically $(c-1)K_Z$; the third applies when a degree-$j$ section or
intrinsic jet is nonzero on the curve. A generic lc centre suffices;
closed non-lc points do not invalidate this normalized-centre inequality.
We conservatively allow every inertia $2,3,4,5,6,7,8$ in the finite
scans, and remove inertia six only when using its geometric absence.
There are at most $\lfloor2(2+cw_{max}/840)\rfloor$ marks. Enumerate
all their multisets, with $\gcd(k_i,r_i)=1$, and impose
\eqref{eq:quintic-marked-tests46}. This is a complete bounded procedure.
It gives
\[
\begin{array}{c|c|l}
c&w\text{ range}&\text{possible }w\\\hline
12&1\le w\le18&12,15\\
6&1\le w\le35&15,20,21,24,28,30,35\\
7/2&1\le w\le12&\varnothing.
\end{array}
\]
For $c=12,w\le12$ the sole row is
$w=12$, marks $(1/2,4/5,5/7)$. Genus one gives no marked row in these
ranges, and an unmarked positive degree is integral, so all these
small curves are rational. The exact program
\texttt{quintic\_marked\_scans.py} constructs all rows directly from
these three inequalities.

\paragraph{The actual polar and its curvature criterion}
Write a nonzero two-form as $\alpha=\iota_\delta\eta$ in a canonical
frame, and put $\theta=\delta a$, $\gamma=\delta\theta$.
An invisible prime is generically invariant under $\delta$, so every
power of its ideal is preserved. Complete restriction, including the
nilpotent structure of $(a)$, gives a unique lift $u\in R_{12}$ of
$\theta^2$. Write $u=\theta^2+a\sigma$ locally. If the frame changes
by a unit $h$ and $t=\delta h/h$, the coefficients, after removing the
appropriate powers of $h$, transform as
\[
 \theta'=\theta-5ta,\quad
 \gamma'=\gamma-11t\theta+(30t^2-5\delta t)a,\quad
 \sigma'=\sigma+10t\theta-25t^2a.
\]
The complete degree-seven restriction gap therefore gives
\[
                    11\sigma+10\gamma=a\rho.
\]
For an actual degree-$j$ section set $D_af=5a\delta f-j\theta f$.
The frame-change terms cancel; this is a global degree-six derivation,
with $D_aa=0$. Polynomial expansion gives global sections
\begin{equation}\label{eq:quintic-polar-identities46}
 v=-D_au/12,\quad v^2-u^3=a^2w_{26},\quad
 D_av=-18u^2+a^2s_{14},\quad D_aw_{26}=2vs,
\end{equation}
with restrictions
\begin{equation}\label{eq:quintic-curvature-restriction46}
 v_T=\theta_T^3,\quad (w_{26})_T=-\theta_T^2s_T/13,
 \quad s_D=\frac{975}{121}\gamma^2
       -\frac{325}{33}\theta\delta\gamma+\frac{104}{33}\theta^2\rho.
\end{equation}
These normal divisions are valid on the complete canonical-cover
quotients before normality of $D$ is known. On an invisible component
of multiplicity $d$, tangency gives
\[
             \operatorname{ord}(u,v,s,w_{26})\ge(2d,3d,2d,4d).
\]

Let $N$ be the finite normalization of $T_Z$ and still write $H$ for
its canonical root. Proper dualizing trace of the nonzero restricted
two-form, followed by multiplication, injects
$H^0(N,4H)$ into the complete degree-ten restriction on $(a)_Z$.
That target is zero because $P_{10}=P_5=1$ and $H^1(5K_Z)=0$.
Trace annihilates the other components with their full multiplicities,
so this remains true when the invisible remainder is nonreduced.
Taking powers on the integral normalization gives
\begin{equation}\label{eq:quintic-normalization-gaps46}
             H^0(N,H)=H^0(N,2H)=H^0(N,4H)=0.
\end{equation}
The section $s$ is nonzero. Otherwise $D_aw_{26}=0$; at the visible
prime the leading derivative of an order-$k$ nonzero section has
coefficient $(5k-26)\theta_T$, never zero. Hence $w_{26}=0$ and
$v^2=u^3$, making the integral rational degree-six section $v/u$
regular, contrary to $R_6=0$.
Consequently
\begin{equation}\label{eq:quintic-curvature-criterion46}
 s_T/\theta_T^2\text{ cannot be regular in codimension one on }N.
\end{equation}
Indeed regularity would give a degree-two section, hence $s_T=0$ by
\eqref{eq:quintic-normalization-gaps46}; the invisible orders would
then give $s/a\in R_9=0$.

\paragraph{Eliminating the five strict numerical states}
On a smooth normalization resolution $S$ of $T$, let $A=K|_S$ and let
$C_T$ be the intrinsic conductor. Normalization adjunction gives the
full equality
\[
 I A\bigl(\operatorname{div}(\alpha|_S)+C_T\bigr)=n+M.
\]
For a reduced full quintic, every positive-degree image in this divisor
has subadjunction coefficient $c\le12$: at a singular generic point
use the full quintic threshold, with $c\le6$; at a smooth polar of
order $j$, use $D+(2j)^{-1}(u=0)$, with $c=6+6/j$.
The states $(9,9)$ and $(10,8)$ have intrinsic mass eighteen and
reduced full canonical member. They would express eighteen as a
nonnegative integral combination of twelve and fifteen, impossible.

For $(10,22)$ the canonical member is prime and the full crepant
remainder $N_0$ has $IKTN_0=5n-M=28$. At a curve centre, the transverse
Du Val resolution has integral exceptional cycle $N_\Gamma$, and
\[
 TN_\Gamma=-N_\Gamma^2\in2\mathbb Z_{>0}.
\]
This follows by intersecting its complete Cartier pullback with the
exceptional curves and using the even ADE lattice. Splitting monodromy
and retaining residue degrees preserves evenness. Its singular image
curve has coefficient at most six and degree at least fifteen, so any
nonzero contribution is at least thirty, a contradiction.

For $(7,-7)$ the intrinsic mass is zero. At every positive-degree curve
of the visible normalization, $T$ is regular and its restricted form
is a unit. On the regular terminal model write $a=te$, where $t$ cuts
out $T$ and $e$ is the full invisible remainder. Then
$\theta_T=e\delta t$, with $\delta t$ a unit, while $s$ has order at
least twice that of $e$. Thus $s_T/\theta_T^2$ is regular, contradicting
\eqref{eq:quintic-curvature-criterion46}.

We give the additional argument for $(8,4)$, including a potentially
nonreduced remainder. If $F$ is an invisible prime of positive scaled
degree $s_0=IK^2F\in\{1,2\}$, set $N_F=IKF^2$.
On its smooth normalization $A^2=s_0/I$ and
$IAK_S\le s_0+N_F$. The following elementary lower bound makes the
character search finite. If $IAK_S\le-L<0$, put
$W=\lfloor3Is_0/L\rfloor$. Then
\begin{equation}\label{eq:quintic-ruled-bound46}
                    W<I/8,\quad2W<L
                    \quad\Longrightarrow\quad\text{impossible}.
\end{equation}
Run the smooth surface extremal program for $K_S+tA$, with rational
$t<L/s_0$ sufficiently close to it. The final ruled fibre has
$A$-degree below $2/t$, or a line on the plane endpoint has degree
below $3/t$; resolving a general line pencil handles the latter.
Thus there is a rational ruling of positive scaled fibre degree at
most $W$. Since its degree is below $1/8$, a general fibre has at least
two distinct fractional canonical marks. After finite base change,
split two into sections, send them birationally to the constant zero
and infinity sections on $\mathbb P^1\times B'$, and resolve.
Their images have the same $A_0$-degree $x$; since both original
sections are $A$-null, the blowup weights over each sum to $x$.
Therefore $A'K_{S'/B'}=-2x+\sum a_i\ge0$.
The vertical ramification index is at most the base ramification index,
and all resolution corrections are $A'$-null. Pushing down gives
$AK_S\ge(2g(B)-2)e\ge-2e$, contradicting $IAK_S\le-L$ and $2W<L$.
This proves \eqref{eq:quintic-ruled-bound46}.

It excludes $N_F\le-72$ for $s_0=1$ and $N_F\le-103$ for $s_0=2$.
Let $l=IKTF$. The simultaneous congruences, Hodge index and effectivity
now give the bounded test
\[
 0\le l\le36,\quad(8s_0-2l)^2\le56(s_0^2-2N_F),
 \quad N_F\le s_0^2/2.
\]
The null equality $(s_0,N_F)=(2,2)$ is incompatible with invisibility
by the null Hodge-pairing lemma. All remaining possibilities are
\[
              (s_0,N_F,l)=(1,-59,18),(2,-26,36),
                              (2,-40,8),(2,-68,36).
\]
For a completely explicit enumeration, separate the basket into the
$(2,8)$-, $3$-, $5$- and $7$-blocks; use congruence moduli $(8,16)$
for the first pair of residues and $(p,p)$ for the other blocks.
Their visible-character counts are $3,1,2,1$. For each visible tuple
retain every ordered character tuple for $F$, so all relative matchings
are present. Test the displayed finite ranges for $N_F,l$.
The program \texttt{quintic\_companion\_scan.py} performs exactly
this enumeration from the basket alone.

The positive remainder of total degree two is one degree-two prime,
two degree-one primes, or twice a degree-one prime. Except for the
row $(2,-40,8)$, the table gives $IKTN_0=0$. In that row the member
is reduced, with contacts $IKTN_0=28$ and $IKFN_0=42$.
Half their sum is thirty-five, which must be a sum of coefficient-six
degrees. The only possibilities are $35$ and $15+20$.
At each curve the separate contacts are nonnegative integers whose
sum is the even ADE contribution. The first partition permits the
$T$-contact only $0,35,70$; the second only $15a+20b$ with
$a,b\in\{0,1,2\}$. Neither gives twenty-eight. Thus this row is also
impossible, and the ambient canonical model is regular at every
generic curve point of $T_Z$.

The intrinsic mass is twelve. Outside the positive remainder, any
intrinsic curve is the single reduced degree-twelve row
$(1/2,4/5,5/7)$; a conductor would have coefficient at least two.
Its degree-seven section cost is $11/10>7\cdot12/840=1/10$, so
$\gamma$ vanishes there. At its generic point the parameter ideal
$J=(a,\theta)$ is preserved by $\delta$. The polar identities give
$u\in J^2$, $v\in J^3$, $a^2s\in J^4$; since
$(J^4:a^2)=J^2$, the curvature ratio is regular.
If the positive remainder is reduced, its intersection with $T_Z$
is singular on the reduced full member; coefficient six excludes
all degrees at most twelve. If the remainder is $2F_Z$, use the
boundary $F_Z+T_Z/2$. Unless it is generically plt, coefficient at
most $7/2$ gives the same contradiction. In the plt case the transverse
branches are smooth and transverse. Write $a=tz^2$ and
$\delta z=zh$. If the intrinsic polar order is $j$, the ideal
$J=(t,z^j)$ is preserved, and
\[
 \theta=z^2p,\quad\gamma=z^2g,\quad\sigma=z^2h',\qquad p,g,h'\in J.
\]
Thus $u=z^4U$ with $U\in J^2$, $v=z^6V$ with $V\in J^3$, and
$s=z^4S$. The identity for $D_av$ gives $t^2S\in J^4$, hence
$S\in J^2$. This is precisely regularity of $s_T/\theta_T^2$.
Every case contradicts the curvature criterion. Only $(n,M)=(10,50)$
remains; the canonical quintic is prime and $IKTN_0=0$, so no
curve-centred crepant remainder has centre in it. Point centres remain
allowed at this stage.

\paragraph{A K3 restriction at square $1/84$.}
\begin{lemma}[The degree-four section at square $1/84$]
\label{lem:quintic-k3-gap46}
Let $D_0$ be a K3 surface with only Du Val singularities and let $H$
be a nef and big integral Weil divisor with $H^2=1/84$. Then
$h^0(D_0,4H)=1$.
\end{lemma}
\begin{proof}
The exceptional ADE root lattice in the smooth K3 resolution has rank
at most nineteen. For a local class of $mH$, let $q_P(m)$ be the
minimum squared norm in its discriminant coset, using the positive
Cartan form. Surface Riemann--Roch and vanishing give
\[
 h^0(mH)=2+\frac12\left(\frac{m^2}{84}-\sum_Pq_P(m)\right).
\]
In particular $\sum q_P(1)\equiv1/84\pmod{2\mathbb Z}$.
A seven-primary discriminant component is necessary. In rank at most
nineteen, it must occur in $A_6(i)$, $1\le i\le3$, or in
$A_{13}(i)$, $1\le i\le6$, up to reversal. The remaining rank is
at most thirteen or six, respectively.

Here is a complete finite calculation, including nonprimitive local
classes. Form the Cartan matrices of every $A_n,D_n,E_6,E_7$ within
that remaining rank. For each inverse Cartan matrix, take every
fundamental-weight column modulo $\mathbb Z^n$. Equal columns modulo
$\mathbb Z^n$ represent the same coset; keep its smallest diagonal
value. Adjoin the zero coset with norm zero. These are all minimum
norms: the nonzero minimum representatives are the minuscule weights
in ADE type, and the listed fundamental weights include them all.
The number of resulting cosets is the Cartan determinant, which is
also checked. Multiplying a column by $m$ gives $q_P(m)$.
Enumerate arbitrary multisets of the remaining nonzero classes with
total rank at most the stated bound. Classes zero on $H$, including
$E_8$, may be omitted: doing so enlarges the available rank and leaves
every correction unchanged. No primitive embedding or gluing
restriction is imposed, so the test includes every geometric case.

This inverse-Cartan calculation performs $9480$ labelled tests. Exactly
seventeen satisfy the congruence for $m=1$; four already have a negative
value of the displayed dimension and are impossible geometrically.
All seventeen, including those four, have
\[
                         \sum_Pq_P(4)=46/21.
\]
For clarity, the enumeration is bounded solely by rank nineteen:
for $A_n$ use all $1\le i\le\lfloor(n+1)/2\rfloor$; for $D_n$
use its vector and spin classes; for $E_6,E_7$ use the nonzero class.
The two $D_4$ labels are harmlessly retained even though triality
identifies them. The attached standard-library program generates
these matrices and multisets anew and lists every surviving profile;
it reads no precomputed candidate list. Substitution gives
\[
 2+\frac12\left(16/84-46/21\right)=1,
\]
as required.
\end{proof}


\paragraph{Normality}
Suppose a curve $C$ is nonnormal on the prime quintic. The ambient
model is regular at its generic point. In the transverse regular
plane put $q_C=i(a,\delta a)$ and $w_C=840K_ZC$.
The proper polar intersection has total degree sixty, while singular
threshold adjunction gives $c\le6$. Hence $q_Cw_C\le60$ and
$w_C\ge15$. If $\delta$ preserves the transverse maximal ideal,
$a,\theta\in\mathfrak m^2$ force $q_C\ge4$. Equality consumes the
entire budget: $q_C=4,w_C=15$. Coprime quadratic initial forms then
force a node, with polar order two on each branch. The normalized
square gives $\sigma\in\mathfrak m^2$, $u\in\mathfrak m^4$,
$v\in\mathfrak m^6$, $a^2s\in\mathfrak m^8$, hence
$s\in\mathfrak m^4$. The curvature ratio is regular on both branches,
a contradiction.

Otherwise straighten $\delta=\partial_x$. For
$k=\operatorname{ord}_x a(x,0)$ the plane polar identity is
$q_C=\mu+k-1\ge k\ge2$. The first nonzero intrinsic jet has degree
$5+k$. It suffices to test $k=2,3,4$. The complete marked tests give
\[
\begin{array}{c|l}
k& (w_C;\text{marks})\\\hline
4&\varnothing\\
3&(15;1/7,7/8)\\
2&(15;1/7,7/8),(20;1/6,6/7),
  (24;3/5,3/7),(30;3/4,2/7).
\end{array}
\]
At the ordinary seventh source, the $k=3$ branch has coordinate orders
at least $(2,5,6)$ in characters $(1,6,3)$. A degree-eight numerator
has character four and minimum monomial order eight, whereas the
exhausted section budget requires one. The degree-twenty row uses
inertia six and is absent.
For the degree-fifteen row with $k=2$, lift its full eighth inertia
to a terminal endpoint. The source, allowing its full axial weight,
has coordinates $(x,y,z,t)$ of characters $(1,7,3,0)$ and equation
$xy+G(z^8,t)=0$. The branch has lower orders $(3,5,1,8)$.
The derivation coefficients have characters $(6,4,0,5)$ and lower
orders $(2,4,0,7)$, so differentiation lowers this filtration by at
most one. Therefore
\[
 \theta=cz^2+O_{\rm wt}(10),\qquad
 \gamma=2c\delta z(0)z+O_{\rm wt}(9).
\]
The required order-one jet forces $c\delta z(0)\ne0$, contrary to
$\theta|_C=0$. The endpoint transfer, with curve-centred sources
included, is proved in Subsection~\ref{sec:quintic-polar-sources46}.

For either remaining row $q_C=2$ and $\mu=1$, so the transverse curve
is a node and its conductor has multiplicity two. If $w_C=30$, the
conductor consumes all sixty units. The actual form on the finite
normalization has no divisorial zero, so $\omega_N\simeq\mathcal O_N$.
Its pullback is holomorphic, making $N$ Gorenstein canonical and hence
Du Val. The smooth minimal resolution has trivial canonical class
and $p_g=1$. It cannot be abelian, since then there are no rational
exceptional curves and the integral Cartier polarization has positive
even square, contrary to $H^2=1/84$. It is K3, and
Lemma~\ref{lem:quintic-k3-gap46} contradicts
\eqref{eq:quintic-normalization-gaps46}.

If $w_C=24$, the residual polar has degree twelve and marks
$(1/2,4/5,5/7)$. On $C$, the degree-seven jet has degree $1/5$, all
consumed at its fifth mark; it is a unit at the seventh source.
Equivariant flow and Morse splitting there give
$a=z^2+g(x,y)$, with $g$ of character six. The marked $3/7$ branch
has normalized parameter character six and $y$-order $B\equiv1$
modulo seven. Its monic branch equation in $x$ has degree $B$ and
character one. The nodal factorization $g=h^2l$ therefore has $l$ of
character four, which cannot be a unit. Conjugate unmarked branches
come in orbits of seven and have invariant norm; they cannot absorb
that character. The residual degree-twelve polar must consequently
have a seventh branch on which the jet is a unit. But its degree-seven
cost is $11/10>1/10$, impossible. The prime quintic is regular in
codimension one; its full Cohen--Macaulay quotient is $S_2$, hence it
is normal.

\paragraph{Rational singularities and the K3 minimal model}
Complete restriction gives $h^1(\mathcal O_D)=0$ and
$h^2(\mathcal O_D)=1$. The restricted actual two-form and dualizing
trace give $p_g(W)=1$ on the minimal resolution. With $A=\pi^*H$,
\[
            A^2=1/84,\qquad AK_W=6H^2=1/14,
            \qquad (AK_W)^2/A^2=3/7<1.
\]
The general-type case is excluded by Hodge index as in
Lemma~\ref{lem:small-square-irregularity46}; that lemma then gives
$q(W)=0$ and rational singularities of $D$.
Complete restriction also gives
\[
 h^0(D,2K_D)=h^0(D,12H)=P_{12}-P_7=1.
\]
Corollary~\ref{cor:exceptional-bicanonical-k3} now identifies the
smooth minimal model as K3.


\subsection{The polar curve and its actual sources}
\label{sec:quintic-polar-sources46}
We now identify the surface geometry needed for the algebraic recognition.
In particular the two eighth basket constituents will be shown to be
separate ordinary sources. This is proved below and is not assumed
from the basket.

\paragraph{The complete polar has only three possible profiles}
A useful gap criterion is the following. If $D$ is normal and rational,
$M$ is integral Weil, $M-K_D$ is ample and $H^0(D,M)=0$, then $D$ is
klt wherever $M$ is Cartier. On a resolution the sheaf
\[
 \mathcal J_M=\pi_*\mathcal O_W
          (K_W+\lceil\pi^*(M-K_D)\rceil)\subset\mathcal O_D(M)
\]
has $H^1=0$ by vanishing and the Leray edge injection. The quotient
has finite support and zero global sections, hence is zero. Where
$M$ is Cartier this is the unit multiplier ideal, exactly the klt
condition. Apply this to $M=7H$; the quintic is klt at the seventh source.

The full polar has scaled degree sixty. At a generic regular point of
a polar component of multiplicity $j$, if $\delta$ is tangent, the
ideal $J=(a,x^j)$ is preserved. The normalized polar identities give
$u\in J^2$, $v\in J^3$ and $a^2s\in J^4$, so the curvature ratio
is regular there. Not every component can be tangent, by
\eqref{eq:quintic-curvature-criterion46}. For a transverse component,
the first nonzero jet has degree $6+j$ and the actual square gives
subadjunction coefficient $c=6+6/j$, with $jw\le60$.
The marked tests have no row for $j=3,4,5$ and the degree floor twelve
excludes larger $j$. For $j=2$ the five rows are
\[
 (15;1/7,7/8),\quad(21;1/2,2/5,1/8),\quad(21;2/5,5/8),
 \quad(30;1/2,1/4,2/7),\quad(30;3/4,2/7).
\]
At a seventh mark $1/7$ or $2/7$, a degree-eight numerator has minimum
order eight or nine, whereas its saturated cost requires one or two.
The remaining degree-twenty-one rows consume forty-two units, leaving
eighteen, which cannot be a sum of the available degrees twelve and
fifteen. Thus a reduced component has a nonzero degree-seven first jet.

There are exactly thirty-two first-jet rows from
\eqref{eq:quintic-marked-tests46}, all with exhausted section budget.
The full eighth-source calculation in the normality proof removes six
rows containing $7/8$; absence of inertia six removes six, with one
overlap. All twenty-one remaining rows have a full seventh mark, so
$\gamma$ is a unit there and every polar branch through that source
occurs in this same list.
Equivariant flow and Morse splitting give $a=z^2+g(x,y)$ of character
six. A linear $y$ term makes the cover smooth and gives one branch
marked $4/7$. Otherwise the klt cover is Du Val. There is no quadratic
term of character six and its only cubic is $xy^2$, so the equivariant
Du Val equation is
\[
               a=z^2+x(y^2+x^{7k-2}),\qquad D_{7k}.
\]
For odd $k$ its branches have marks $(3/7,1/7)$; for even $k$ they
have $(3/7,4/7,4/7)$. The twenty-one global rows have degree at least
twenty-four, so at most two can occur. Matching their complete seventh
branch multisets and the total degree bound sixty leaves exactly
\begin{equation}\label{eq:quintic-three-polars46}
\begin{array}{c|c|c}
 w&\text{complete marks}&\text{seventh slice}\\\hline
60&(1/2,4/7)&\text{smooth}\\
60&(1/2,1/7,3/7)&D_{7k},\ k\text{ odd}\\
60&(3/4,3/4,4/7)&\text{smooth}.
\end{array}
\end{equation}
The bounded enumeration in \texttt{quintic\_marked\_scans.py} includes
all these branch partitions. Each row is one reduced rational curve
and exhausts the entire polar divisor.

\paragraph{Index eight, with curve-centred geometry retained}
At index eight there is no closed crepant centre. The canonical
index-one cover is therefore cDV. A general Du Val hyperplane through
it, avoiding the generic points of all curve centres, has no exceptional
divisor in its pullback. Inversion of adjunction and crepancy make its
strict transform normal, and adjunction gives a partial crepant
resolution of a Du Val surface. The reduced fibre is consequently
obtained from an ADE tree by contracting connected subtrees and has
at most three ends. The hyperplane need not be invariant: it describes
the intrinsic fibre on which the group already acts.

A cyclic two-power action on a tree with at most three ends has component
orbits of size at most two. If two components were exchanged by
$\mu_8$, their stabilizer $\mu_4$ would act faithfully on either
rational component; a kernel would give nontrivial terminal index at
its generic point. At most one of its two fixed points is also on the
other exchanged component. The other has exactly inertia four, absent
from the terminal basket. Thus all components are invariant and carry
faithful $\mu_8$ actions. A tree of $n$ such components has at least
$n+1$ distinct full eighth fixed points. Only two eighth constituents
are available, so the fibre is a point or one rational curve. In the
latter case both endpoints are ordinary terminal eighth points, since
their axial weights are one. A branch of proper inertia two or four
must lift to one of these endpoints: a faithful cyclic action on a
rational curve has no other nontrivial stabilizers. This also justifies
the full-inertia endpoint transfer used in the normality argument.

At such an endpoint allow the full terminal cover
$xy+G(z^8,t)=0$ with characters $(1,7,3,0)$.
The numerator characters of $(a,\theta,\gamma)$ are $(1,6,3)$ and
those of the derivation coefficients are $(6,4,0,5)$.
At a half or quarter marking, exhausted order one of $\gamma$ gives
\[
              \gamma=cz+O(\mathfrak m^2),\quad c\ne0,
              \qquad\delta z(0)\ne0.
\]
Indeed the only linear character-three eigenfunction is $z$. If
$\delta z(0)=0$, all derivation coefficients vanish; the only possible
linear term of $a$ is $x$, and $\delta x$ has no linear eigenfunction
of its character six. Thus $\delta^2a\in\mathfrak m^2$, a contradiction.
A regular vector field vanishes at an isolated singular point, so the
terminal endpoint is smooth. Moreover $\delta\gamma$ is a unit there.
Using the canonical-cover top frame and its crepant pullback, contraction
with the actual two-form commutes with the morphism. All derivatives of
$a$ pull back compatibly, so $\delta\gamma$ is already a unit on the
canonical germ. Its pushed vector field is nonzero there.

Equivariant flow splits that germ as $S\times\Delta_z$, with
$\delta=\partial_z$ and $z$ of character three. The surface $S$ is
smooth or Du Val and its action is faithful symplectic: its canonical
character is zero, and a kernel would fix a threefold divisor.
On the minimal equivariant resolution of $S$, a nontrivial symplectic
automorphism has no fixed curve. A nontrivial diagram action produces
proper inertia four by the same two-fixed-point argument; a trivalent
vertex would have three fixed points and hence be fixed pointwise.
Thus the diagram is $A_n$, with $n+1$ full eighth fixed points, forcing
$n=1$ in the singular case. The only coalesced alternative is therefore
\begin{equation}\label{eq:quintic-coalesced-source46}
 \{xy-t^2=0\}\times\Delta_z,\qquad
 \mu_8:(x,y,z,t)\mapsto(\zeta x,\zeta^7y,\zeta^3z,t).
\end{equation}
To check its characters, lift the symplectic cyclic action on $A_1$
to its smooth double cover. It has order sixteen and characters
$(c,-c,0)$ on $(x,y,t)$ for odd $c$. Both terminal endpoint types
being $(3,8)$ forces $c=\pm1$; exchange $x,y$ if necessary.

\paragraph{Half and quarter markings distinguish the two source types}
A half marking cannot be at an actual index-two point. The order-one
odd-character jet makes the normalized cover branch smooth. The complete
polar and its exhausted jet budget exclude any additional branch.
The reduced polar is the Gorenstein complete intersection $(a,\theta)$;
its dualizing character is $1+5+6=12$, even, whereas a smooth full
inertia-two parameter has odd differential character. Actual indices
four and six are absent, so the source must have index eight.

If its cover is smooth, a half-inertia branch has odd orders in $x,y,z$.
The exact jet gives $z$ order one, and cancellation in $\theta=0$
forces $y$ order one. Unless $a$ has a linear $x$ term, its initial forms
are
\[
 \theta_2=3Az^2+By^2,\qquad
 a_3=Az^3+By^2z+Cx^2y,\qquad AB\ne0.
\]
A linear $x$ term instead gives two smooth conjugate polar branches,
with inertia four, not two. Without it, $C$ must be nonzero, since
substituting $\theta_2=0$ otherwise leaves $-2Az^3$.
For $C\ne0$ the two initial forms are a regular sequence with four
nonaxis directions and an additional $x$-axis direction. Their complete
intersection is the tangent cone, so the axis supplies an extra actual
branch, contrary to the exhausted polar and jet budget.
Thus a half marking forces \eqref{eq:quintic-coalesced-source46}.

There the lower orders are $(3,1,1,2)$ on $(x,y,z,t)$.
The weight-three equation has only $cx+Az^3+By^2z$; the same argument
forces $cAB\ne0$. Eliminate $x$ to obtain
\[
 t^2+Ay z^3+By^3z+\text{terms of weight }>4=0,
 \qquad\operatorname{wt}(y,z,t)=(1,1,2).
\]
Its weighted blowup has a smooth elliptic exceptional curve of square
$-2$. The radial kernel is $\mu_2$, and the residual order-four action
has orbifold orders $(4,4,2)$. On the $z$ chart, with
$\rho=z,u=y/z,v=t/z^2$, the action is
$(\rho,u,v)\mapsto(\zeta^3\rho,-u,iv)$.
After dividing the radial kernel, the two order-four points have
opposite tangent and normal characters, giving $A_3$ arms, and the
order-two orbit gives an $A_1$ arm. The central rational curve has
square $-1$ before resolving these points; subtracting
$3/4+3/4+1/2$ gives the minimal graph
\[
              (-3)\text{ with arms }A_3,A_3,A_1.
\]
Its discrepancy coefficients are $1$ centrally,
$(3/4,1/2,1/4)$ on each long arm and $1/2$ on the short arm; the
square cost is one. The argument applies to every higher-order
perturbation of the displayed initial form.

For a quarter mark $3/4$ the lower orders are $(3,1,1,4)$.
The same initial equations force a nonzero linear $x$ coefficient.
In the coalesced model this would give $\operatorname{ord}(xy)=4$
while $\operatorname{ord}(t^2)\ge8$, impossible. The source is
therefore an ordinary terminal eighth point, with smooth linear
quintic cover and quotient slice $1/8(7,3)\simeq1/8(1,5)$.
One source has two smooth conjugate branches and thus exactly one
coarse quarter mark. Two quarter marks come from two distinct sources.

\paragraph{One canonical blowup and the two remaining K3 frames}
On the minimal resolution $W\to D$, the actual form has exactly one
positive-$H$-degree zero component, with coefficient one. Every
$(-1)$ curve has positive degree, by minimality of the resolution.
A point-blowup sequence from the K3 minimal model with two terminal
leaves would give two such components; with only one leaf and at least
two blowups, that leaf has coefficient at least two. Thus there is at
most one blowup.

The ordinary seventh slice has chain $[3,2,2]$, coefficients
$(3,2,1)/7$, and cost $3/7$. An ordinary eighth slice has chain
$[2,3,2]$, coefficients $(1/4,1/2,1/4)$, and cost $1/2$.
For odd $D_{7k}/\mu_7$ the cost is always $10/7$.
Here is a uniform graph computation. Put $s_0=7k-2$; the
binary-dihedral cover followed by the scalar action gives a central
rational curve with orbifold orders $(2,2,s_0)$, two short $(-2)$ arms,
and long Hirzebruch--Jung fraction $s_0/q$, $q\equiv-7\pmod{s_0}$.
For $k=1$ the centre is $-3$ and the long arm is $[2,3]$.
For odd $k\ge3$ the centre is $-2$ and the long arm is
$[2^{k-2},3,2,3]$. Solving the discrepancy equations gives central
and initial coefficients $6/7$, final three $6/7,5/7,4/7$, and short
coefficients $3/7$. The two $(-3)$ vertices alone contribute
$6/7+4/7=10/7$.

Since $K_D^2=3/7$, the odd-$D$ half-mark row of
\eqref{eq:quintic-three-polars46} would give
$K_W^2=3/7-10/7-1=-2$, contradicting the one-blowup bound.
Each other row has $K_W^2=-1$, so $W$ is exactly one K3 blowup,
with canonical exceptional curve $E$. For every surface resolution
component $C$,
\[
                       EC=K_WC=-2-C^2.
\]
In the ordinary row, $E$ meets the three $(-3)$ centres transversely
once and misses every $(-2)$ curve. Their images are three $(-2)$
curves meeting at an ordinary triple point. Their sum is nef and
isotropic; an adjacent arm meets it once, so it is primitive and is
an entire type-IV elliptic fibre. The adjacent seventh curve and four
eighth arms are sections, with the remaining seventh tail vertical.

In the half-mark row, blowing down $E$ gives an $I_0^*$ configuration:
the eighth central curve has multiplicity two, and its two first long
arm components, its short arm and the seventh centre are the four ends.
The next long-arm curves and next seventh curve are three sections,
each with an additional vertical $A_1$ tail.
We must include the low-index chains even when their canonical images
coalesce. Outside the polar $\theta$ is a unit; equivariant flow gives
$S\times\Delta_a$ with $S$ Du Val and faithful symplectic action.
The surface canonical character is six, while the normal fifth
character is primitive, so the actual index is two or three.
Its minimal equivariant resolution has only isolated stabilizers,
with opposite surface tangent characters; adjoining the normal
character gives ordinary terminal quotients on a crepant terminalization.
Basket invariance counts the remaining four halves and two thirds.
On the central surface these give disjoint marked $4A_1+2A_2$ paths.
Their further crepant surface resolution contains no $(-1)$ curve,
by adjunction, and hence is the minimal resolution of $D$. The marked
paths are not lost even if additional curves connect them.
They avoid the polar graph, survive the one blowdown and are vertical.
Together with its three tails they give an induced $7A_1+2A_2$
away from the distinguished $I_0^*$ fibre.

\begin{lemma}[An elliptic K3 packing obstruction]
\label{lem:quintic-K3-packing46}
An elliptic K3 surface with a section cannot have a distinguished
$I_0^*$ fibre and an induced disjoint vertical configuration
$7A_1+2A_2$ away from it. Unselected components may connect the marked
paths; they need not be complete fibres.
\end{lemma}
\begin{proof}
For a fibre containing $a$ selected singles and $b$ selected double
paths, assign mass $2a+3b$ and saving
$\epsilon=2a+3b-e(\text{fibre})$. If its affine graph is a tree and
$U$ is the unselected vertex set, edge counting gives
\[
 \epsilon=(a+b)-|U|-1
       =\sum_{v\in U}(\deg v-2)-e(U).
\]
A multiplicative cycle has no positive saving, since gaps separate
selected blocks. The small additive types II, III, IV likewise do not.
For $I_n^*$, positive saving requires both branching vertices unselected
and all four ends selected; otherwise a selected branching component
forces a missing end whose negative contribution removes any saving.
Thus it consumes at least four singles and has saving at most two.
Saving two means $U$ is independent. If the interior selected blocks
comprise $k$ singles and $l$ doubles, then
\[
                       a=4+k,\quad b=l,\quad n=2k+3l.
\]
This includes $I_0^*$ with $k=l=0$.
For IV$^*$, III$^*$, II$^*$, saving at most one follows from the single
branching vertex. Equality requires, respectively, $3A_2$, $5A_1$,
or $4A_1+A_2$ or $A_1+3A_2$: on an arm of length $L$, independent
unselected vertices require $L+1=2a_L+3b_L$, and the arm lengths are
$(2,2,2)$, $(3,3,1)$ and $(5,2,1)$.
Only two doubles are available. Every possible positive-saving fibre
therefore consumes at least four of the seven singles, so there is
at most one.

The distinguished fibre uses Euler number six, leaving eighteen.
The total assigned mass is twenty, so the savings must sum to two.
There is exactly one further $I_n^*$ with saving two, and every other
fibre has saving zero. Quadratically twist at this fibre and at the
distinguished $I_0^*$. They become $I_n$ and $I_0$, all other fibres
remain unchanged, and the Euler number drops by twelve. The resulting
elliptic surface has section, base $\mathbb P^1$ and $\chi=1$; it is
rational. The unchanged paths give $(3-k)A_1+(2-l)A_2$. The new
$I_n$ chain $A_{n-1}$ supplies $kA_1+lA_2$, since the selected blocks
and their separating vertices use $2k+3l-1=n-1$ vertices.
Thus a rational elliptic surface contains $M=3A_1+2A_2$ in its
vertical lattice.

Subtracting $(C\cdot O)F$ from a vertical class preserves all its
pairings and makes it perpendicular to $O,F$. Their orthogonal
complement is the even unimodular rank-eight lattice $E_8(-1)$:
the rational surface Néron--Severi lattice is unimodular,
$\det\langle O,F\rangle=-1$, and parity follows from $K=-F$.
But the two-primary discriminant form of $M$ is $(\mathbb Z/2)^3$,
with nonzero norms $-s/2$ for $s=1,2,3$ modulo $2\mathbb Z$.
None is isotropic. Every even overlattice of $M$ has odd index, so its
primitive closure still has two-primary discriminant length three.
A primitive sublattice of a unimodular lattice has the discriminant
group of its orthogonal complement, while that complement here has
rank one and hence cyclic discriminant. This is impossible.
\end{proof}

\begin{proposition}\label{prop:quintic-geometric-reduction46}
For every canonical model with basket $\mathcal B_5$ and volume $1/420$,
the first quintic is prime, normal and rational, with $K_D=6H$ and
$H^2=1/84$. Its complete polar is one reduced rational curve of degree
$HE=1/14$ and marks $(4/7,3/4,3/4)$, with nonzero first jet.
The seventh source is an unchanged ordinary quotient with smooth
linear quintic cover; the two eighth sources are distinct ordinary
quotients with smooth linear quintic covers. Its minimal resolution
is one blowup of a K3 surface at the triple point of a type-IV fibre.
\end{proposition}
\begin{proof}
The packing lemma excludes the coalesced half-mark frame. All other
assertions were proved above, retaining point-centred remainder until
it ceased to affect the actual surface. No classification step has
replaced virtual constituents by actual points without a source proof.
\end{proof}


\subsection{Independent construction of the complete quintic surface algebra}
\label{sec:quintic-explicit-surface46}
We construct the surface algebra directly, prove its normality, and
then recover it from every geometric realization in
Subsection~\ref{sec:quintic-surface-recognition46}.

Put $A_0=\mathbb C[\theta_6,b_8,z_{14}]$ and
\begin{equation}\label{eq:quintic-even-surface46}
\begin{split}
T&=A_0[P_{16},q_{18}]/(f,g),\\
f&=P^2-zq-\tfrac14b^4+\tfrac3{44}\theta^4b,\\
g&=q^2-bz^2+\tfrac8{11}\theta zP
                  +\tfrac7{11}\theta^2b^3+\tfrac4{121}\theta^6.
\end{split}
\end{equation}
The monic relations have relatively prime leading monomials $P^2,q^2$
for the ordinary degree order in $P,q$ over $A_0$. Thus $T$ is free
with basis $1,P,q,Pq$. It is a three-dimensional complete intersection.
Taking the coefficient of $Pq$ gives the Frobenius Gram matrix
\[
 \begin{pmatrix}
 0&0&0&1\\0&0&1&0\\0&1&0&0\\1&0&0&-8\theta z^2/11
 \end{pmatrix},
\]
of determinant one. Since the weights of $A_0$ sum to $28$ and
$\deg(Pq)=34$, finite duality gives $\omega_T=T(6)$.

\paragraph{The odd module and a complete multiplication rule}
Let $N$ have the six generators
\[
        e_{15},J_{21},K_{23},f_{19},L_{25},M_{27}
\]
and precisely the two $A_0$-linear relations
\begin{equation}\label{eq:quintic-surface-Koszul46}
                  ze-bJ+\theta K=0,
             \qquad zf-bL+\theta M=0.
\end{equation}
It is a direct sum of two Koszul cokernels, hence has rank four,
is torsion free, has depth two at the vertex, and is free elsewhere
over $A_0$. On $b\ne0$ use the basis $(e,K,f,M)$ and define the
$P,q$ actions by the following matrices, acting on columns:
\begin{equation}\label{eq:quintic-surface-action46}
\mathsf P=\begin{pmatrix}
0&b^3/4&7\theta z/11&3b^2\theta^2/44-z^2\\
b&0&3\theta^2/11&-7\theta z/11\\
-\theta z/b&-b\theta^2/4+z^2/b&0&-b^3/4\\
-\theta^2/b&\theta z/b&-b&0
\end{pmatrix},
\end{equation}
\[
\mathsf Q=\begin{pmatrix}
\frac{4\theta^2z}{11b}&-\frac{\theta(b^2\theta^2+4z^2)}{11b}
                         &bz&-7b^3\theta/22\\
\frac{4\theta^3}{11b}&-\frac{4\theta^2z}{11b}&14b\theta/11&-bz\\
z&-b^2\theta/2&-\frac{4\theta^2z}{11b}
                         &\frac{\theta(b^2\theta^2+4z^2)}{11b}\\
2\theta&-z&-\frac{4\theta^3}{11b}&\frac{4\theta^2z}{11b}
\end{pmatrix}.
\]
These rational expressions preserve the global module. For a direct
polynomial check, lifts on $(e,J,K,f,L,M)$ are
\[
\widetilde{\mathsf P}=\begin{pmatrix}
0&b^2\theta/4&b^3/4&4\theta z/11&3b\theta^3/44&3b^2\theta^2/44-4z^2/11\\
0&0&0&3b\theta/11&-4\theta z/11&-7bz/11\\
b&z&0&0&0&0\\
0&-\theta^3/4&-b\theta^2/4&0&-b^2\theta/4&-b^3/4\\
-\theta&0&z&0&0&0\\
0&0&0&-b&-z&0
\end{pmatrix}
\]
and
\[
\widetilde{\mathsf Q}=\begin{pmatrix}
0&-\theta^4/11&-b\theta^3/11&-3bz/11&-7b^2\theta^2/22+8z^2/11&-7b^3\theta/22\\
4\theta^2/11&0&-4\theta z/11&14b^2/11&3bz/11&0\\
0&0&0&0&0&-bz\\
-z&-b\theta^2/2&-b^2\theta/2&0&\theta^4/11&b\theta^3/11\\
2b&z&0&-4\theta^2/11&0&4\theta z/11\\
0&0&-z&0&0&0
\end{pmatrix}.
\]
Multiplying either relation vector by these matrices again gives a
combination of the two relation vectors. Their localizations are
\eqref{eq:quintic-surface-action46} and $\mathsf Q$.

Set
\[
 \mathsf H=\begin{pmatrix}0&0&0&-b\\0&0&b&0\\0&b&0&0\\-b&0&0&0\end{pmatrix},
 \qquad h(v,w)=v^{\mathsf t}\mathsf Hw.
\]
The product of two odd vectors is defined by
\begin{equation}\label{eq:quintic-odd-product46}
\begin{split}
vw={}&h(\mathsf P\mathsf Qv,w)+\tfrac8{11}\theta z^2h(v,w)
   +h(\mathsf Qv,w)P+h(\mathsf Pv,w)q+h(v,w)Pq.
\end{split}
\end{equation}
Substitute the six global generators, with
$J=(ze+\theta K)/b$ and $L=(zf+\theta M)/b$.
Every coefficient in the resulting twenty-one products is polynomial
in $A_0$; for example
\begin{equation}\label{eq:quintic-surface-odd-square46}
                         e^2=b^2z-2b\theta P+\theta^2q.
\end{equation}
The complete verification consists of
\[
[\mathsf P,\mathsf Q]=0,\quad
\mathsf P^{\mathsf t}\mathsf H=\mathsf H\mathsf P,\quad
\mathsf Q^{\mathsf t}\mathsf H=\mathsf H\mathsf Q,\quad
f(\mathsf P,\mathsf Q)=g(\mathsf P,\mathsf Q)=0
\]
and, for the $4^3$ basis triples,
\[
 \mathsf M(vw)x=\mathsf M(wx)v,
 \qquad
 \mathsf M(c_0+c_1P+c_2q+c_3Pq)
        =c_0I+c_1\mathsf P+c_2\mathsf Q+c_3\mathsf P\mathsf Q.
\]
Each is an identity in the displayed rational functions; torsion
freeness extends it across $b=0$. Thus these formulas define an
associative graded algebra $B=T\oplus N$ in all degrees.

\paragraph{Normality and smoothness away from the vertex}
For clarity we give an exhaustive Jacobian calculation for $T$.
On $\theta=1$ its singular ideal has triangular basis
\begin{equation}\label{eq:quintic-surface-singular-ideal46}
\begin{split}
q&=(99b^5+33b^2)z/10,\qquad
P=(99b^4+43b)z/20,\\
z^2&=(99b^5+53b^2)/20,\qquad
0=b^6+17b^3/33+16/1089.
\end{split}
\end{equation}
The last polynomial is $(b^3+1/33)(b^3+16/33)$; it has six distinct
nonzero roots and each gives two distinct $z$ values. On $\theta=0$,
$b\ne0$, the singular orbits are $z=q=0$, $P=\pm b^2/2$.
On $\theta=b=0$, $z\ne0$, the derivatives $f_q=-z$ and $g_b=-z^2$
are independent. These strata exhaust the punctured cone, so $T$ is
regular in codimension one. Its complete-intersection property and
Serre's criterion prove normality. A normal positively graded algebra
with degree-zero part $\mathbb C$ is connected and hence a domain.

The duality behind the odd product is valid globally. Put
\[
r=(z,-b,\theta)^{\mathsf t},\qquad
J_0=\begin{pmatrix}0&\theta&b\\-\theta&0&z\\-b&-z&0\end{pmatrix}.
\]
The Koszul complex identifies $\operatorname{coker}r$ with
$\ker r^{\mathsf t}$ with the appropriate graded shift. Consequently
the matrix $\left(\begin{smallmatrix}0&-J_0\\J_0&0\end{smallmatrix}\right)$
on the two summands of $N$ is a perfect $A_0$-duality. Combining it
with the Frobenius pairing of $T$ gives
\[
                  N\xrightarrow{\sim}\operatorname{Hom}_T(N,T),
           \qquad v\longmapsto(w\longmapsto vw).
\]
The $T$ action is faithful: every nonzero ideal of the finite domain
$T$ meets $A_0$ nontrivially and $N$ is $A_0$-torsion free.
The generic rank over $T$ is one. The displayed depth calculation
shows $N$ is reflexive, so at codimension one its square pairing is
a unit. Hence $B/T$ is etale there. Both summands are $S_2$;
therefore $B$ is normal and, by connected grading, a domain.

In fact $\operatorname{Spec}B\setminus\{0\}$ is smooth.
Over the regular locus of $T$, the reflexive rank-one module is free
and its square is a unit, giving an etale double cover of a regular
ring. At the twelve points of \eqref{eq:quintic-surface-singular-ideal46},
the ideal of the ten products in the basis $(e,K,f,M)$ has the same
triangular basis as the even singular ideal. Thus the odd products
generate the even maximal ideal in the transverse surface slice.
The odd fibre is the cokernel of
$[\mathsf P-PI\mid\mathsf Q-qI]$. Its $224$ three-by-three minors
vanish, whereas a two-by-two minor reduces to
$-9b^4/10-71b/220$, a unit on that scheme. The odd fibre therefore
has dimension two. If $\mathfrak n$ denotes the cover maximal ideal,
the even maximal ideal lies in $\mathfrak n^2$; consequently
$\dim\mathfrak n/\mathfrak n^2=2$, proving regularity of the slice.
At $\theta=0,b=1,P=\pm1/2$, the action matrix again has rank two,
and the linear terms of the odd products span the three even
cotangent coordinates $(\theta,z,q)$. The same argument applies,
with odd tangent coordinates $e,f$. Restoring the grading orbit
proves the assertion about the entire punctured cone.

All operations just used are polynomial reduction, determinant and
rank identities from the displayed matrices. The accompanying
\texttt{quintic\_surface\_algebra.py} repeats the identities and the
identification with the central fibre of the ambient table. Geometric
normality follows from the preceding depth and codimension arguments,
not from the output of that calculation.

\paragraph{The invariant algebra and the single missing class}
There is an involution fixing $b,z,q,e,f,J,L$ and negating
$\theta,P,K,M$; the product formulas verify it. Put
\[
             u=\theta^2,\quad k=\theta P,\quad U=B^{\iota}.
\]
The invariant even part is free over
$A_1=\mathbb C[b_8,u_{12},z_{14}]$ on $1,q,k,kq$.
In the odd Koszul summands the invariant relations eliminate
$\theta K$ and $\theta M$. Hence $U$ is free on
\begin{equation}\label{eq:quintic-invariant-surface-basis46}
                      1,q,k,kq,e,f,J,L.
\end{equation}
It is normal as an invariant subring of the normal domain $B$.
The anti-invariant module is generated over $U$ by $\theta,P,K,M$.
Remove only the degree-six class and put $C=B\ominus\mathbb C\theta$.
This is a subring: every positive-degree element has degree at least
six, so no product can recreate the omitted class. Its anti-part is
generated over $U$ by
\[
 -b\theta,P,\theta^3,\theta z,\theta e,K,\theta q,
                         \theta f,\theta J,M.
\]
Over $A_1$ add $\theta L,Pq$. The even-degree part of this module is
$\theta(b,u,z)$ together with the three free generators $P,\theta q,Pq$;
the odd-degree part has exactly the two Koszul relations
$z\theta e-b\theta J+uK=0$ and
$z\theta f-b\theta L+uM=0$.
Its minimal resolution has twelve generators, relations in weights
$26,28,32,35,39$ and one second relation of weight $40$.
This single missing class accounts for the two intermediate
cohomology groups when $\chi(\mathcal O)=2$.


\subsection{Recognition of the complete quintic surface algebra}
\label{sec:quintic-surface-recognition46}
Let $Z$ now be an arbitrary canonical model with the quintic basket.
Proposition~\ref{prop:quintic-geometric-reduction46} supplies the prime
normal rational divisor $D=(a_5)$ and its complete reduced polar.
We show that the full section algebra $B=\bigoplus_mH^0(D,mH)$ is
exactly the algebra constructed in Subsection~\ref{sec:quintic-explicit-surface46}.
The calculation retains every surface parameter until the actual
ambient polar derivatives force its value.

\paragraph{The actual parameter ring and the two Koszul summands}
The full divisorial restriction sequence, not just restriction in
codimension one, gives
\begin{equation}\label{eq:quintic-complete-surface-cohom46}
 h^1(D,mH)=\mathbf1_{m=1}+\mathbf1_{m=5},\qquad
 h^0(D,mH)=P_m-P_{m-5}+\mathbf1_{m=6}\quad(m\ge0).
\end{equation}
Here $P_0=1$ and $P_j=0$ for $j<0$. This follows from
Lemma~\ref{lem:complete-canonical-restriction46}, vanishing and duality.
In particular $\theta_6$ is precisely the missing class in $R/(a)$,
and $B_4=B_7=0$.
The actual normalized polar construction gives
\[
 u=\theta^2-10a\gamma/11+a^2\rho/11,\qquad
 v=-D_a u/12,\qquad s=(D_av+18u^2)/a^2,
\]
where $D_af_j=5a\delta f_j-j\theta f_j$ and
$D_bf_j=8b\delta f_j-j(\delta b)f_j$. Their frame-change terms cancel.
Their restrictions are
\begin{equation}\label{eq:quintic-curvature-surface46}
 u_D=\theta^2,\quad v_D=\theta^3,\quad
 s_D=975\gamma^2/121-325\theta\delta\gamma/33
                                  +104\theta^2\rho/33.
\end{equation}
The section $b$ does not vanish identically on the polar $E$; otherwise
$b/\theta\in B_2=0$. Its degree $8HE=4/7$ is exhausted at the
seventh mark, while the degree of $7H|_E$ is exhausted by the two
$1/4$ zeros of $\gamma$ at the eighth marks. Since
$s|_E=(975/121)\gamma^2$, the three sections $\theta,b,s$ have no
common projective zero. They form an actual polynomial parameter ring.

Temporarily halve the even grading. Then
$A'=\mathbb C[\theta_3,b_4,s_7]\subset T=B^{(2)}$,
and \eqref{eq:quintic-complete-surface-cohom46} makes $T$
Cohen--Macaulay with $\omega_T=T(3)$. Its exact Hilbert numerator is
$1+t^8+t^9+t^{17}$. It is free over $A'$ with those generator degrees.
The Artinian Gorenstein quotient has perfect pairing between degrees
eight and nine, so choose $p,q$ such that $1,p,q,pq$ is a free basis.
The two squares have no $pq$ coefficient, since the required parameter
weights are $-1$ and $1$. Their two monic equations reduce every
monomial to this basis, proving that $T$ is the resulting complete
intersection; no generation assumption is concealed in the series.

Write $N_n=B_{2n+1}$. Its sole local cohomology is
$H^2_{A'_+}(N)=\mathbb C\oplus\mathbb C(-2)$, annihilated by $A'_+$.
Duality gives $\operatorname{Ext}^1_{A'}(N,A')=
\mathbb C(14)\oplus\mathbb C(16)$. Its minimal resolution has length
one. Dualizing it gives a free presentation of these two residue
fields; the kernel of their minimal free cover is
$A'_+(14)\oplus A'_+(16)$. Its minimal presentation is the pair of
Koszul presentations. Any extra free summand splits off, and the two
Koszul cokernels already have the required rank four. Hence
\begin{equation}\label{eq:quintic-forced-Koszul46}
 N=E_{A'}\oplus E_{A'}(-2),\qquad E_{A'}=\Omega^2_{A'}\mathbb C.
\end{equation}
In original degrees the two relations are
$se-bj+\theta k=0$ and $sf-b\ell+\theta m=0$, with generator
weights $(15,21,23;19,25,27)$.

\paragraph{The complete polar ring fixes the boundary multiplication}
Place the three polar points at $\infty,0,1$. An integral translate
of its rational divisor is
$\Delta=-10[\infty]/7+3[0]/4+3[1]/4$.
The entire normalized polar ring is
\begin{equation}\label{eq:quintic-polar-ring46}
 V=\mathbb C[\eta_4,\gamma_7,p_{16}]/(p^2-\eta^4p-\eta\gamma^4).
\end{equation}
Indeed, in $\mathbb C(x)[t]$ take
$\eta=t^4/(x^3(x-1)^3)$, $\gamma=t^7/(x^5(x-1)^5)$ and
$p=x\eta^4$. The dimension in weight $d$ is
$\max(0,1-\lceil10d/7\rceil+2\lfloor3d/4\rfloor)$.
Writing each monomial with $p$-exponent at most one verifies equality
with the hypersurface series, using the period $28$ of the floor
remainders. Thus the displayed generators give the complete ring.

The reduced integral divisor and the depth-two restriction identify
$B/(\theta)$ with its image in $V$. The two distinct ordinary eighth
sources give two distinct roots of the quadratic $p$ equation at
$\gamma=0$. Translation and nonzero rescaling therefore give
\begin{equation}\label{eq:quintic-polar-even46}
 T/(\theta)=\mathbb C[b,s,p,q]/(p^2-b^2p-sq, q^2-bs^2),
 \quad b=\eta^2, s=\gamma^2, q=\eta\gamma^2.
\end{equation}
The coefficient of $sq$ cannot vanish, because the first equation
would factor in an integral polar ring; the coefficient of $bs^2$
cannot vanish because $q$ is not nilpotent. These justify the rescalings.
The odd difference from $V$ is exactly $t^7+t^{11}$. Hence its odd
part is generated by
\[
 e=\eta^2\gamma,\ j=s\gamma,\ k=p\gamma,\quad
 f=\eta^3\gamma,\ \ell=\eta s\gamma,\
 m=\eta(b^2-p)\gamma,
\]
with the two reduced Koszul relations.

These polar bases and the global Koszul bases can be normalized
simultaneously. The automorphism
$j\mapsto j+c\theta e$, $k\mapsto k+cbe$ removes the unwanted
$b^2\gamma$ term in $k$; the second summand has the analogous
operation. If the remaining leading coefficients of $k,m$ are
$A,C\ne0$, the coefficient-of-$pq$ pairing gives
$h(k,f)=Ab$, $h(e,m)=-Cb$, $h(j,f)=A\theta$, $h(e,\ell)=-C\theta$.
Pairing the first Koszul relation with $\ell$ and comparing
$\theta s$ forces $A=C$. One scalar change of $\theta$ normalizes
both summands. The degree-zero cross homomorphisms of the two summands
vanish. Degree-zero endomorphisms of one are the identity and the
nilpotent $e\mapsto0,j\mapsto\theta e,k\mapsto be$; restriction to
$\theta=0$ is injective on their span. The polar products thus fix
the full Frobenius pairing, which on $(e,k,f,m)$ is the matrix
$\mathsf H$ of Subsection~\ref{sec:quintic-explicit-surface46}.

\paragraph{Exhaustive self-adjoint actions: a four-parameter graph}
Here we give the entire coefficient space, including its degenerate
loci. A degree-$d$ endomorphism of $E_{A'}$ is a matrix $L$ with
$Lr=\lambda r$, $r=(s,-b,\theta)^{\mathsf t}$, modulo $r\varphi$;
its trace is $\operatorname{tr}L-\lambda$.
The coefficient weights are $d+g_j-g_i$, $g=(7,10,11)$.
In degrees $(6,7,8,9,10,11)$ the endomorphism dimensions are
$(4,6,4,4,7,7)$, and the traceless dimensions are $(3,4,3,3,5,5)$.
This follows by listing the nonnegative solutions of
$3i+4j+7k=d+g_j-g_i$ and solving $Lr=\lambda r$.
Consequently self-adjoint $p,q$ actions have $25$ coefficients;
the polar values leave the following $15$.
Write
\[
\mathsf P=\begin{pmatrix}U_p&V_p\\W_p&\operatorname{adj}U_p\end{pmatrix},
\qquad
\mathsf Q=\begin{pmatrix}U_q&V_q\\W_q&\operatorname{adj}U_q\end{pmatrix},
\]
where the two-by-two blocks are
\begin{align*}
U_p&=\begin{pmatrix}0&0\\b&b^2\end{pmatrix},\\
V_p&=\begin{pmatrix}
-bb_4\theta^2-b_3s\theta&b^2b_1\theta^2+b(b_0+b_4)s\theta-s^2\\
-(2b_3+1)\theta^2&bb_4\theta^2+b_3s\theta
\end{pmatrix},\\
W_p&=\begin{pmatrix}
-(bc_2\theta^2+s\theta)/b&(b^2c_0\theta^2+2bc_2s\theta+s^2)/b\\
-\theta^2/b&(bc_2\theta^2+s\theta)/b
\end{pmatrix},\\
U_q&=\begin{pmatrix}
-a_0\theta^3-a_2s\theta^2/b&a_1b\theta^3+(a_0+\alpha)s\theta^2+a_2s^2\theta/b\\
-a_2\theta^3/b&\alpha\theta^3+a_2s\theta^2/b
\end{pmatrix},\\
V_q&=\begin{pmatrix}b(s-\chi b\theta)&b^2(s+b_qb\theta)\\
b(\beta+1)\theta&-b(s-\chi b\theta)\end{pmatrix},\\
W_q&=\begin{pmatrix}s-db\theta&b(s+c_qb\theta)\\2\theta&-s+db\theta\end{pmatrix}.
\end{align*}
For a two-by-two matrix, $\operatorname{adj}U=(\operatorname{tr}U)I-U$.
The coefficients of $b^2[\mathsf P,\mathsf Q]$ give $22$ polynomials.
Their ideal over $\mathbb Q$ equals the following eleven-equation
ideal, without localizing at any parameter:
\begin{equation}\label{eq:quintic-surface-graph46}
\begin{aligned}
a_0&=-\alpha-\beta d+\chi,&
a_1&=(\beta d^2-2\chi d+\chi+d^2-d)/2,&a_2&=(\beta-1)/2,\\
b_0&=(-\beta d+2\chi+d-1)/2,&b_1&=\chi(1-d),&b_3&=-(\beta+1)/2,\\
b_4&=(\beta d+d-1)/2,&b_q&=-\chi,&c_0&=d(d-1),\\
c_2&=1/2-d,&c_q&=-d.
\end{aligned}
\end{equation}
This equality is verified by taking reduced lexicographic Groebner
bases in the eleven left variables followed by $(\alpha,\beta,\chi,d)$.
Both reduced bases are exactly the eleven displayed differences.
Thus no exceptional parameter component is lost.
Frobenius reconstruction gives
\[
 e^2=\theta^2q-2\theta bp+b^2s+(2-d)\theta b^3
       +[-\alpha+(\beta-1)(1-2d)/4]\theta^5.
\]
The actual square map, established next, forces
$d=2$, $\alpha=3(1-\beta)/4$. The even relations then become
\begin{equation}\label{eq:quintic-two-parameter-surface46}
\begin{split}
F_s={}&p^2-b^2p-(s-\theta b)q+\beta\theta^4b/4,\\
G_s={}&q^2-bs^2+(1-\beta)\theta sp+\chi\theta b^2s
 +2(\beta-\chi+1)\theta^2bp\\
 &+(\chi-(\beta+3)/2)\theta^3q+(\beta-1)^2\theta^6/16.
\end{split}
\end{equation}
The pairing and these actions determine every odd product.

\paragraph{Actual covariants, including exclusion of a degenerate generator}
Return to original degrees; write $b_i=\delta^ib$ and $\eta=\delta\gamma$.
The following successive divisions define genuine ambient sections:
\begin{equation}\label{eq:quintic-actual-covariants46}
\begin{gathered}
d_{14}=D_ab,\quad e_{15}=(D_ad-112ub)/a,\quad
j_{21}=(ud+8vb)/a,\quad g_{16}=(D_ae-255j)/a,\\
r_{27}=(vd+8u^2b)/a,\quad q^{\rm aux}_{22}=(r+ue/17)/a,\\
h_{17}=(D_ag-7344q^{\rm aux}+39168sb/13)/a.
\end{gathered}
\end{equation}
Substitute the normalized expression for $u$ in independent local jet
variables: every numerator has the displayed factor $a$. Frame-change
cancellation for $D_a$ and normality then give global sections.
Since $R_{17}=aR_{12}$, $h=\kappa au$. Define also
\begin{equation}\label{eq:quintic-actual-Q46}
Q_{18}=\frac{704b^2u-11(D_ba)^2+10aD_b^2a}{64a^2}
       =b^2\rho-25bb_2/4+425b_1^2/64.
\end{equation}
Universal expansion gives
\begin{equation}\label{eq:quintic-actual-square46}
 e_D^2/1088=\theta^2Q_D-2\theta b(g_D/216)+b^2(17s_D/39)
\end{equation}
and
\[
 e^2-34dg/27-1088uQ-18496b^2s/39+136abh/513=a^2J_{20},
\]
where, with $\rho_i=\delta^i\rho$,
\begin{equation}\label{eq:quintic-J-polynomial46}
\begin{split}
J={}&625b_2^2-21250b_1b_3/27+85000bb_4/513
 -13600\rho_2b^2/513+8500\rho_1bb_1/513\\
 &+10625\rho b_1^2/27-191000\rho bb_2/513+576\rho^2b^2/19.
\end{split}
\end{equation}
The complete degree-twenty space is
$\langle a^4,ae,bu,D_as\rangle$. At the ordinary seventh source,
$D_as$ has a nonzero linear term. The character $-3j$ of a weight-$j$
numerator gives minimum ordinary orders
$(2,1,2,2,1,1,2,2)$ for
$(b,b_1,b_2,b_3,b_4,\rho,\rho_1,\rho_2)$.
Every term in \eqref{eq:quintic-J-polynomial46} has order at least
three, as do $a^4,ae,bu$. Therefore
\begin{equation}\label{eq:quintic-J-span46}
                       J=A_0bu+B_0ae+C_0a^4.
\end{equation}
On the seventh polar branch write
$b=c_0x^4+\cdots$, $b_1=c_1x+\cdots$, $c_0\ne0$.
If $c_1\ne0$, the exact polar orders of $g,Q$ are one and two.
Thus they are the required independent even generators, and
\[
                    \lim\frac{Q(17s/39)}{(g/216)^2}=1.
\]
We now prove $c_1\ne0$.
If $c_1=0$, the exhausted polar degrees and the gaps $B_{10}=R_{11}=0$
give $g=\lambda b^2$ and $Q=\mu v+\nu a^2b$.
The latter follows because a nonzero polar restriction of $Q$ has
order exactly two, while its putative order is at least nine;
then $Q_D\in\mathbb C\theta^3$ and division by $a$ uses
$R_{13}=\mathbb Cab$.

On a canonical-root chart $b=1$, all positive $b$ derivatives vanish,
$Q=\rho$, and direct expansion gives
\begin{align*}
g&=(-5100\eta+2820\theta\rho+290a\delta\rho)/11,\\
v&=\theta^3-15a\theta\gamma/11+25a^2\eta/66
                +a^2\theta\rho/66-5a^3\delta\rho/132,\\
J&=-13600\delta^2\rho/513+576\rho^2/19,\\
h&=(484500\delta^4a-596220\gamma\rho-94810\theta\delta\rho
                    -62750a\delta^2\rho+51840a\rho^2)/143.
\end{align*}
Differentiate the full ambient expression before restricting.
Substitution of $\rho=\mu v+\nu a^2$ and $g=\lambda$ gives
\begin{align*}
J_D&=-95200\mu\theta\gamma^2/1881+2464\mu^2\theta^6/171
                  +(272\mu\lambda/4617-27200\nu/513)\theta^2,\\
s_D&=975\gamma^2/121-39\mu\theta^5/17+13\lambda\theta/612.
\end{align*}
Freeness over $\mathbb C[\theta,s]$ on $b=1$ makes these parameters
independent on every component, and then makes $\theta,\gamma$
independent. Comparison with \eqref{eq:quintic-J-span46} forces
$\mu=0$. Now $\rho=\nu a^2$ and $g=\lambda$ give
$\eta=2\nu a^2\theta/3-11\lambda/5100$.
The identity $h=\kappa au$ first gives $\kappa=330880\nu/143$ and then
\[
            0=119680\nu a^2(2\nu a^3-9\gamma)/1573.
\]
Since $\gamma|_E\ne0$, $\nu=0$ and $Q=0$.
But the universal identities
\[
D_ba=-d,\quad D_bd=-88be/85-896aQ/85,
\quad bD_be=-16b^2g/15-24dQ-20aD_bQ/3
\]
would give $D_b^4a=0$, while $D_ba|_D=8\theta b\ne0$.
This contradicts Lemma~\ref{ag:lem:no-time} on the normalized
canonical-root cover of $b$, as explained after that lemma. Its
tautological top form and the pulled-back two-form are holomorphic
on a projective resolution; the nonconstant finite derivative chain
would give a rational time. No extension of the rational vector field
across its poles is needed.
We have proved $c_1\ne0$.

\paragraph{The two remaining parameters from first polar derivatives}
Use actual covariants
$p_0=g/216$, $z_0=17s/39$, $q_0=Q$, $e_0=e/\sqrt{1088}$.
Their square is \eqref{eq:quintic-actual-square46}. Comparing its two
leading coefficients with the preceding Koszul normalization fixes
its possible scalar $\theta\mapsto c\theta$: $c^{-2}=c^{-1}=1$.
Thus \eqref{eq:quintic-two-parameter-surface46} holds in a normalization
coupled to the actual square map.
On $b_0=1,a=\theta_0=0$, differentiation before restriction gives
\begin{equation}\label{eq:quintic-first-polar-derivatives46}
 p_0=-425\eta/198,\quad z_0=425\gamma^2/121,\quad q_0=\rho,
 \quad\delta p_0=-16\gamma q_0/11,\quad
 \delta z_0=-14\gamma p_0/11.
\end{equation}
Write the raw even equations to first order as
\begin{align*}
F_0&=p_0^2+cb_0^2p_0+d_0b_0^4-z_0q_0+f_0\theta_0b_0q_0+O(\theta_0^3),\\
G_0&=q_0^2+ib_0z_0^2+j_0\theta_0z_0p_0
                                  +k_0\theta_0b_0^2z_0+O(\theta_0^2).
\end{align*}
Their first ambient derivatives vanish on $E$, because $\delta a=0$
there. Independence of $p_0$ on each component of the root polar,
with $q_0^2=-iz_0^2$, gives
$j_0=-8i/11$, $k_0=2if_0-32ic/11$.
For completeness the coordinate change is explicit: choose
$\lambda^2+c\lambda+d_0=0$, put $\Delta=2\lambda+c\ne0$, choose
$\alpha,\beta_0\ne0$ with $\beta_0^3=-i\alpha^4$, and set
$t_0=-\Delta\beta_0^2$. Then
\[
\theta_0=\alpha\theta,\quad b_0=\beta_0b,\quad
p_0=t_0p+\lambda\beta_0^2b^2,\quad
z_0=\alpha t_0s/\beta_0+2\alpha\lambda\beta_0\theta b,
\quad q_0=\beta_0t_0q/\alpha.
\]
It preserves the actual square after scaling $e_0$.
Substitution into \eqref{eq:quintic-two-parameter-surface46} gives
\[
 \beta=3/11,\qquad f_0=-c,\qquad
 \chi=(36\lambda-54c)/(11(2\lambda+c)).
\]
There is a unique degree-nineteen section, which lifts to $R_{19}$.
Its reconstructed odd product satisfies
\[
 ef=b^2q-(\beta+1)\theta^2bs/2
                 +(3\beta+1)\theta^3b^2/4-(\beta-1)\theta^3p/2.
\]
Hence $ef-b^2q$ is divisible by $\theta^2$, also in the coupled
raw normalization. On each component of the root polar put
$\xi=q_0/z_0$, so $\xi^2=-i\ne0$, and scale $f|_E=\xi e_0$.
The same first derivatives give
\[
 e_0=5\sqrt{17}\gamma/11,\quad
 \delta e_0=-18p_0/(5\sqrt{17}),\quad
 \delta q_0=\gamma\xi(-18p_0-27c)/11,
\]
and thus $\delta f|_E=-27c\xi/(5\sqrt{17})$.
On the other hand the complete degree-28 restriction space has basis
$z_0^2,\theta_0b_0z_0,\theta_0^2b_0^2,\theta_0^2p_0$.
The polar restriction of the actual section $D_bf$ must be a scalar
multiple of $z_0^2$. The derivative just obtained is constant on each
component, while $z_0$ is nonconstant there by finiteness over
$\mathbb C[b_0,z_0]$. Therefore $c=0$, and
\[
                         \boxed{\beta=3/11,\qquad\chi=18/11.}
\]
Center $P=p-b^2/2$, $z=s-\theta b$. The even equations become
\eqref{eq:quintic-even-surface46}; the raw actual covariants
$p_0,z_0,q_0$ differ from $P,z,q$ only by nonzero scalars, with no
remaining translations. The self-adjoint actions and Frobenius pairing
then give exactly \eqref{eq:quintic-surface-action46} and
\eqref{eq:quintic-odd-product46}. This proves recognition of the
complete surface algebra, including both non-Cohen--Macaulay
Koszul summands and every product.


\subsection{The ambient symmetry of every quintic realization}
\label{sec:quintic-involution46}
A unique extension inside an invariant presentation would not by itself
classify arbitrary threefolds. We now prove that every realization has
that invariant presentation. The proof uses a transverse differential
system for actual ambient covariants and the already recognized initial
surface.

The centered surface equations admit the symmetry fixing $\theta,b$
and negating $P,z,q$. Its actual odd square
$\epsilon^2=\theta^2q_0-2\theta bp_0+b^2z_0$ changes sign.
The full normal surface algebra is the integral closure of the even
ring in this quadratic field, so the symmetry extends by
$\Sigma\epsilon=i\epsilon$, with $\Sigma^2$ equal to parity.
The absence of translations in the final covariant normalization
means these are the signs of the actual $p_0,z_0,q_0$.

One further derivative character is needed. Normalize the first
Koszul generator $k_{23}$ by its polar value $p_0\gamma_0$, where
$\gamma_0=\epsilon/b$. The actual identity is
\[
           k_{23}=-85(D_bz_0)|_D/(112\sqrt{17}).
\]
Its polar value has character $-i$; the restriction kernel is
$\theta B_{17}=0$, so $\Sigma k_{23}=-ik_{23}$ on all of $D$.
For the actual ambient section $K_{27}=D_bQ$, the universal jet
identities give
\begin{equation}\label{eq:quintic-K-character46}
\sqrt{17}\left(b^2k_{23}-b\epsilon p_0-\epsilon\theta q_0
             +\frac{85\theta^2K_D}{144\sqrt{17}}\right)
                 =\frac{17\theta b^2h_D}{2736}=0.
\end{equation}
All terms have weight $39$; the first three have character $-i$.
As the ring is integral, $\Sigma K_D=-iK_D$.
This argument concerns the actual derivative, not a choice of sign
for an abstract odd generator.

\paragraph{The covariants generate the ambient field}
Put
\[
 R_0=\mathbb C[a,b,u,s,d,e,g,Q,v,K]\subset R.
\]
It contains the finite parameter ring $\mathbb C[a,b,u,s]$, hence
$R$ is finite over it. Let $R'$ be its normalization in its own field.
The divisor of $a$ on the normal affine cone of $R$ is one prime of
multiplicity one. In the finite extension from $R'$ there is one
prime above its image with ramification index one. The residue
extension is trivial: the deficient fibre has the same field as $B$
because $\theta=v/u$, and the residues of $R_0$ contain
$\theta,b,s,p_0,q_0,\epsilon$, which generate that field.
The discrete-valuation degree formula in characteristic zero gives
$[\operatorname{Frac}R:\operatorname{Frac}R_0]=1$.
Thus $R$ is precisely the normalization of $R_0$.
At the generic initial divisor, its normalized DVR $C$ and the local
ring $A$ of $R_0$ satisfy $C=A+aC$ by residue-field equality and
valuation one of $a$. Nakayama on $C/A$ gives $C=A$.
The covariant image is therefore smooth on a dense open neighborhood
of that divisor.

\paragraph{A transverse polynomial differential system}
On the canonical-root chart $b=1$ use coordinates
$(a,\theta,\gamma,p,q,r)$, where $p=p_0$, $q=q_0$, $r=\delta q$.
Put
\[
 E=(680\gamma-112aq)/11,\quad
 U=\theta^2-10a\gamma/11+a^2q/11,
\quad J=A_JU+B_JaE+C_Ja^4,\quad h=\kappa aU.
\]
The identities already proved are equivalent to
\begin{equation}\label{eq:quintic-transverse-ODE46}
\begin{aligned}
\delta a&=\theta,&\delta\theta&=\gamma,\\
\delta\gamma&=-198p/425+47\theta q/85+29ar/510,\\
\delta p&=-16\gamma q/11+8\theta r/9+48aq^2/935+aJ/170-13h/20520,\\
\delta q&=r,&\delta r&=486q^2/425-513J/13600.
\end{aligned}
\end{equation}
All other covariants are recovered polynomially. For example
$d=-8\theta$, $K=8r$, $u=U$, $e=E$ and
\[
 z_0=E^2/1088+2\theta p-Uq+ah/4104-a^2J/1088,
\]
\[
 v=\theta^3-15a\theta\gamma/11+25a^2\delta\gamma/66
                     +a^2\theta q/66-5a^3r/132.
\]
The closure $Y\subset\mathbb A^6$ of the actual image is an
irreducible threefold and has the same field as the root chart.
The polynomial vector field in \eqref{eq:quintic-transverse-ODE46}
preserves its defining ideal, since it is the actual derivation on
a dense open set.

Consider the linear transformation
\[
 \mathcal S(a,\theta,\gamma,p,q,r)=(-ia,\theta,i\gamma,-p,-q,-ir).
\]
Term by term $\mathcal S\delta=i\delta\mathcal S$.
Its action on $a=0$ is the recognized symmetry, including the
character of $r$ from \eqref{eq:quintic-K-character46}.
At a general point of that smooth surface, $\theta\ne0$ and
$\delta a=\theta$ is transverse. The holomorphic flow gives a
unique invariant three-dimensional germ through the given smooth
transverse surface germ. The two varieties $Y$ and $\mathcal S(Y)$
share that germ; their irreducibility gives $\mathcal S(Y)=Y$.
Thus the differential system preserves the actual algebraic image,
not merely an abstract family of possible integral varieties.

For the irreducibility just used, a reducible eighth-root cover over
$\mathbb C$ would make its defining coefficient a square in the
function field. A rational fourth-canonical section with square $b$
is integral over the full normal ring and would lie in $R_4=0$.
The same argument with $B_4=0$ treats the initial surface. Alternatively
one can apply the transverse germ argument on each component and
descend its induced field symmetry.

\paragraph{Canonical linearization, with no fixed divisor}
The homogeneous characters on the covariants are
\[
\begin{array}{c|rrrrrrrrrr}
 &a&b&u&s&d&e&g&Q&v&K\\\hline
\Sigma&-i&1&1&-1&1&i&-1&-1&1&-i.
\end{array}
\]
They commute with the grading and root-deck action, so the relations
homogenize to $R_0$ and extend uniquely to its normalization $R$.
The square is the degree-$m$ scalar $(-1)^m$, inducing an involution
$\iota$ of the canonical model.
Its fixed locus on $D$ is finite. In halved weights $(3,4,7,8,9)$,
if $\theta b\ne0$ a projective fixed point forces the scalar to have
$t^3=t^4=1$, hence $P=z=q=0$. The two equations give
$b^3=3\theta^4/11$ and $25\theta^6/121=0$, a contradiction.
For $\theta\ne0,b=0$ they similarly give $P=q=0$ and
$4\theta^6/121=0$. For $\theta=0,b\ne0$ they force $P=0$ and
leave finitely many points. The remaining locus is the seventh point.
Every projective fixed divisor would meet the ample $D$ in a curve,
so the quotient has no divisorial ramification.

To identify the natural canonical action, use the fixed polar point
other than the seventh point. Here $b,e$ are units. On the root chart
the grading correction by $i$ gives the local action
$(a,\theta,p)\mapsto(a,-\theta,-p)$.
These are parameters: the recognized surface Jacobian proves the
independence of $\theta,p$ and $a$ is a reduced transverse parameter.
The ambient index is one, since unit sections have coprime weights
$8,15$, and the smooth Cartier slice implies ambient smoothness.
The determinant is $+1$, so natural pullback fixes a local canonical
frame. Two graded lifts differ by $c^m$: their quotient acts by global
units on the reflexive canonical sheaves, and multiplication gives a
character of the degree semigroup, which contains all sufficiently
large integers. The units imply $c^8=1$, $c^{15}i=1$, hence $c=i$.
Therefore
\begin{equation}\label{eq:quintic-natural-involution46}
                 \iota^*|_{R_m}=i^m\Sigma|_{R_m}.
\end{equation}
Its signs on the displayed covariants are
$(+,+,+,+,-,+,-,+,-,-)$.
The unique two-form is anti-invariant: at this fixed polar point its
kernel vector has $a$ component $\theta=0$ and nonzero $\theta$
component $\gamma$. Its nonzero value lies in the span of
$dp\wedge da,da\wedge d\theta$, the negative eigenspace.
This is the character that annihilates the intermediate cohomology
of the invariant ring. It asserts nothing about the singularities
of the geometric quotient beyond the quasi-etale statement proved above.


\subsection{Altair's invariant algebra and full canonical-ring candidate}
\label{sec:quintic-construction46}
Let $A=\mathbb C[a_5,b_8,u_{12},z_{14}]$. Define $S$ as the free
$A$-module on
\[
 1,q_{18},k_{22},t_{\flat,40},e_{15},f_{19},j_{21},\ell_{25},
\]
with the products of Appendix~\ref{app:quintic-exact-coefficients}.
The auxiliary basis element $t_\flat$ differs from $qk$ by the explicit
triangular expression in that appendix. This convention is retained
in all module and deformation calculations. It has no geometric effect.
By construction $S$ is free, hence flat, over $A$.

All products are homogeneous and all associators vanish. A complete
check collects the coefficient of each of the eight basis elements in
$(E_iE_j)E_k-E_i(E_jE_k)$ for the $8^3$ triples. The Frobenius functional
extracting the top coefficient has degree $-40$ and determinant one:
modulo the parameter ideal, the complementary pairs are $1,t_\flat$,
$q,k$, $e,\ell$, $f,j$, with coefficients $1,1,-1,1$.
Thus $S$ is Cohen--Macaulay and $\omega_S=S(1)$.
By the central-fibre identification just proved, $S/(a)=U$.

Put
\begin{equation}\label{eq:quintic-ambient-cover46}
\begin{split}
F={}&b^2u+abe+\tfrac14a^2q+\tfrac{2207}{4400}a^4b,\\
H={}&-bk+\tfrac12a\ell-\tfrac{71}{440}a^2bu
                      -\tfrac1{220}a^3e-\tfrac1{48400}a^6,\\
G={}&zq+\tfrac14b^4-\tfrac3{44}u^2b
   +a(\tfrac1{10}bf-\tfrac2{11}ue)
   +\tfrac1{55}a^2k+\tfrac{9123}{4400}a^2bz
   +\tfrac{203}{24200}a^4u.
\end{split}
\end{equation}
Direct multiplication gives $H^2=FG$. Together with all basis associators, this is an exact identity in the printed coefficients.

\paragraph{Normality and the quadratic field.}
The independent surface construction gives $S/(a)=U$ normal.
Since $S$ is free over $A$, Lemma~\ref{lem:exceptional-normality-lift}
shows that it is a normal domain.

Modulo $a$, $F=b^2\theta^2$ and $G=P^2$. In the explicit even surface
algebra, both $b=P=0$ and $\theta=P=0$ have codimension two. Hence
$F,G$ have no common prime divisor on the central fibre, and none
upstairs: a homogeneous common divisor would meet the central fibre
in a common divisor. At every divisor where $F$ vanishes, $G$ is a unit,
so $\operatorname{ord}F=2\operatorname{ord}H$. The square class of $F$
is nontrivial, since its residue square root $-b\theta$ defines the
nontrivial field extension $\operatorname{Frac}B/\operatorname{Frac}U$.

Let $w^2=F$, put $p=H/w$, and set
\begin{equation}\label{eq:quintic-normalization46}
               I=(F):H,\qquad R_*=S\oplus I/w.
\end{equation}
Lemma~\ref{lem:exceptional-quadratic-normalization} shows that
this is the full integral closure in the nontrivial quadratic field,
etale in codimension one, and that $\omega_{R_*}=R_*(1)$ as
reflexive modules. We call $\operatorname{Proj}R_*$ \emph{Altair};
its canonical polarization and terminality will follow from the
local analysis in Section~\ref{sec:exceptional-canonical-geometry}.

\paragraph{The deficient fibre and the full module.}
The twelve numerator expressions in the appendix define generators
with the following weights and central residues:
\[
\begin{array}{c|rrrrrrrrrrrr}
& w&p&v&y&r_{21}&K_{23}&r_{24}&r_{25}&r_{27}&M_{27}&r_{31}&r_{34}\\\hline
\deg&14&16&18&20&21&23&24&25&27&27&31&34\\
a=0&-b\theta&P&\theta^3&\theta z&\theta e&K&\theta q&\theta f&\theta J&M&\theta L&Pq.
\end{array}
\]
For the last two one may use $r_{31}=bK_{23}-pe$ and $r_{34}=pq$.
Membership of every numerator in $I$ is verified by multiplication
with $H$ and division by $F$ in the free invariant algebra.
The ideal is $a$-saturated: from $axH=Fy$, reduction modulo $a$ in
$U$ gives $a\mid y$, and cancellation gives $x\in I$.
Therefore its residue embeds in the quadratic residue field and is
integral over $U$, hence lies in $B$.

The displayed residues generate the deficient algebra $C$ of
Subsection~\ref{sec:quintic-explicit-surface46}. The only possible additional residue is $\theta_6$.
It does not lift, since the map
\[
 S_{20}\oplus S_{22}\longrightarrow S_{50},\qquad
                      (x,y)\longmapsto Hx-Fy
\]
is injective. Its source bases are $(bu,a^4,ae)$ and $(bz,a^2u,k)$;
The target coefficients of
$b^2uk,b^3uz,abu\ell,ab^2ze,a^2b^2u^2,a^4bk$ give the minor
\[
\begin{pmatrix}
-1&0&0&0&0&-1\\
0&0&0&-1&0&0\\
1/2&0&1&0&0&1\\
0&0&1&-1&0&1\\
-71/440&0&-213/440&0&-1&-163/440\\
0&-1&1/11&0&0&9/320
\end{pmatrix},
\qquad\det=\tfrac12.
\]
This calculation uses only the printed products. Thus $R_* /(a)=C=B\ominus\mathbb C\theta_6$, and graded
Nakayama proves generation in every degree.

Write $N_*=I/w$ and
$\mathcal D=(14,16,18,20,21,23,24,25,27,27,31,34)$, retaining the
repeated degree. The full minimal resolution is
\begin{equation}\label{eq:quintic-full-module-resolution46}
0\longrightarrow A(-40)
 \xrightarrow{(z,u,-b,a,0)^{\mathsf t}}
 \bigoplus_{d=26,28,32,35,39}A(-d)
 \longrightarrow\bigoplus_{d\in\mathcal D}A(-d)
 \longrightarrow N_*\longrightarrow0.
\end{equation}
Its five first relations are printed in the appendix. The second
relation is a literal identity. Reduction modulo $a$ gives the
Koszul resolutions of $\theta(b,u,z)$ and the two odd surface
cokernels, together with free summands. The reduced complex is exact;
the homology sequence for reduction and graded Nakayama lift exactness
from left to right. This proves the resolution itself, not just an
alternating Hilbert numerator.

Dualizing with $\omega_A=A(-39)$ gives
\begin{equation}\label{eq:quintic-module-defects46}
 \operatorname{Ext}^2_A(N_*,\omega_A)=\mathbb C(1),
 \qquad \operatorname{Ext}^1_A(N_*,\omega_A)=\mathbb C.
\end{equation}
For an explicit check, the last transpose has image $(z,u,b,a)$.
In the preceding dual term the first four coordinates give the six
Koszul syzygies of $(z,u,-b,a)$; the fifth-coordinate image is exactly
$(a,b,u,z)$. Its surviving constant has degree zero.
The full ring consequently has depth two. Only its invariant
summand $S$ is Cohen--Macaulay.

Its Hilbert series is
\begin{equation}\label{eq:quintic-full-Hilbert46}
H_{R_*}(t)=\frac{
 \sum_{d\in\{0,15,18,19,21,22,25,40\}}t^d
 +\sum_{d\in\mathcal D}t^d
 -t^{26}-t^{28}-t^{32}-t^{35}-t^{39}+t^{40}}
 {(1-t^5)(1-t^8)(1-t^{12})(1-t^{14})}.
\end{equation}
The numerator at one is sixteen, giving volume $1/420$.


\section{Canonical geometry and the local singularity atlases}
\label{sec:exceptional-canonical-geometry}

The preceding sections have defined two finite graded normal rings.
We now prove that their projective spectra are canonical models.
For Vega there is one nonterminal canonical point; both the
crepant extraction and the actual completed germ will be needed.
For Altair the punctured canonical cone is smooth, and the complete
quotient atlas proves terminality. These arguments establish
existence independently of whole-basket recognition.

\subsection{Vega: the complete local atlas and crepant extraction}
\label{sec:quartic-geometry46}

For the quartic algebra, put
\[
 A=\mathbb C[a_4,u_{10},b_{12},c_{14}],\qquad
 S=A\langle1,h_{15},q_{18},r_{20},s_{21},k_{23},l_{26},t_{41}\rangle.
\]
Use the twelve anti-invariant generators
\[
(n_{15},n_{17},n_{19},n_{20},n_{21},n_{23},n_{25},n_{26},
 n_{27},n_{28},n_{30},n_{36})
\]
in the certified module gauges, and abbreviate $v=n_{15}$, $w=n_{17}$,
$z=n_{19}$. Thus $v^2=F$, and $n_j=X_j/v$ in the quadratic fraction-field
extension. The following calculation constructs an actual derivation of the
complete algebra, including the anti module.

\begin{proposition}\label{prop:quartic-global-D}
There is a unique deck-anti-invariant degree-five derivation with the prescribed polar
normalizations
\[
 Da=0,\qquad Du=-10v,\qquad Db=-12w,\qquad Dh=-15n_{20},
\]
and
\begin{align}
 Dc&=-14z-\frac{12}{5}av,\label{eq:quartic-Dc}\\
 Dv&=-15u^2+\frac{1188}{25}a^2b-
                 \frac{55945647}{100000}a^5.
 \label{eq:quartic-Dv}
\end{align}
All its values on the twelve anti generators are regular invariant sections.
\end{proposition}

\begin{proof}
Solve first for $\Delta=vD:S\to S$. For each invariant basis element, allow
its derivative to be an arbitrary anti section of the required weight,
using a linearly independent basis of the exact numerator span. For $Dc$
allow $-14z+\alpha av$, and for $Dv$ allow the right side of
\eqref{eq:quartic-Dv} with an undetermined coefficient of $a^5$.
There are $31$ scalar unknowns. The $28$ nontrivial invariant-product
Leibniz identities and
\[
 \Delta F=2F Dv
\]
give $1,670$ coefficient equations of rank $31$. Exact rational elimination
gives $\alpha=-12/5$ and the coefficient $-55945647/100000$ in
\eqref{eq:quartic-Dv}.

For every anti generator define
\[
 Dn_j=\frac{\Delta X_j-X_jDv}{F}.
\]
All twelve numerators divide exactly by $F$ in $S$. The division is checked
coefficient by coefficient in powers of $a$, starting with the monic term
$F\bmod a=u^3$. The resulting formulas satisfy all $84$ mixed and all $78$
anti-product Leibniz identities. Thus the fraction-field derivation preserves
the entire algebra. In particular, the construction leaves no denominator
on the locus $F=0$.

This is a specified finite rational calculation: list every monomial
in the stated degrees, expand the Leibniz rule in the free invariant
basis, and perform Gaussian elimination. Appendix~\ref{app:quartic-exact-coefficients}
gives all its inputs. The accompanying program reconstructs the system
and verifies every resulting identity, including the mixed and odd products.
\end{proof}

The identity $Da=0$ is essential here: $D$ itself preserves the chart $a=1$;
no Euler correction is required. The singular locus is preserved by a
derivation, since the relevant Fitting ideal of the module of differentials
is preserved, as is its radical in characteristic zero. Consequently a
finite reduced singular locus consists of zeroes of $D$. The independent
boundary calculation proves this finiteness for the projective singular
locus away from $a=0$.

\begin{proposition}\label{prop:quartic-polar-zero}
On the chart $a=1$, every zero of $D$ lies over one of two invariant points,
interchanged by the fourth-root chart action. Each has a unique reduced
lift to the complete cover algebra, with all twelve anti coordinates zero.
Writing $u=54\sqrt6/25$, the positive point is
\begin{align*}
 b&=\frac{329751}{16000},&c&=\frac{24u}{5},&
 q&=-\frac{29007u}{3200},\\
 r&=\frac{35108554221}{200704000},&
 l&=-\frac{8016563871u}{57344000},&h=s=k=t&=0.
\end{align*}
The other point is its Galois conjugate.
\end{proposition}

\begin{proof}
The equations $D(n_j)=0$ and the invariant multiplication equations allow
seven linear eliminations. In particular,
\begin{align*}
 b&=\frac{500000u^2+18648549}{1584000},\\
 q&=\frac{24377760c-850000u^3-136294083u}{4752000},\\
 r&=-\frac{47667200000cu+228577020000u^2-28075922096823}
                  {87306240000},\\
 s&=-\frac{50hu}{93},\qquad
 k=-\frac{h(1225000u^2+283065309)}{11457600},\\
 l&=\frac{549140711760c-49073800000u^3-2344359211353u}
                  {7741440000}.
\end{align*}
After the remaining elimination of $t$, the exact lexicographic Groebner
basis is
\[
 \boxed{h^2,\quad 5c-24u,\quad 625u^2-17496.}
\]
This is an invariant necessary zero scheme; its nonreduced structure is not
being identified with the full zero scheme without the anti equations.
Its reduced support is precisely the two displayed points.

At each point every anti square is zero, so every anti coordinate is zero
in a reduced residue field. All anti cross products are also zero, and the
full presentation vanishes there. These are therefore actual points of the
cover, each with a unique reduced lift. The equations $D(S)=0$ hold because
their values are anti sections. This proves both necessity and existence
on the reduced support. All eliminations take place in the rational free-module table; the
companion verifier retains their exact coefficient identities.
\end{proof}

It remains to examine this single projective orbit. The complete local
presentation has $195$ relations: $28$ invariant products, $84$ mixed
products, $78$ anti products, and five module relations. Its Jacobian has
rank $13$ in $22$ chart variables, hence embedding dimension nine. Thus the
actual local calculation leads to a cone-type canonical singularity, not
a hypersurface calculation. We now identify the quadratic
initial algebra and construct the crepant extraction.

\paragraph{Regularity of the boundary and finite singular locus}

The punctured special fibre is the actual quartic surface cone.  Its
index cover is smooth except at the single seventh-index orbit, where
the surface cover has an $A_1$ singularity.  At a smooth point of this
fibre, regularity lifts across the nonzerodivisor $a$.  At the remaining
orbit use $c=1$, with every other invariant coordinate zero.  The square
of the degree-twenty-one anti-invariant is one, so the cover is locally
\'etale over $S$.  In the ten chart variables
$(a,u,b,h,q,r,s,k,l,t)$, the invariant product equations have rank seven.
An exact seven-by-seven minor is $270$; after imposing $a=0$ the rank is
six.  In particular the $k_{23}^2$ product contains
\[
                         \frac{99}{20}ac^3.
\]
Thus the ambient cover smooths the surface $A_1$, and $R$ is regular
along the entire punctured boundary.

The actual boundary quotient charts are as follows.  The stated
stabilizers are those after passing to the full quadratic cover.
\[
\begin{array}{c|c|c|c}
\text{parameter chart}&\text{stabilizer}&\text{tangent characters}
 &\text{number of projective points}\\\hline
b=1&\mu_3&(1,2,2)&2\\
u=1&\mu_5&(1,4,4)&1\\
u=1&\mu_5&(2,3,4)&2\\
c=1&\mu_7&(2,5,6)&1
\end{array}
\tag{QG.1}
\]
For example, the third-index chart has $h^4=1$, $q=h^2$, $s=h^3$ and
$p_{21}^2=q$; the two residual orbits are distinguished by
$p_{21}/h^3=\pm1$.  The fifth-index chart has
$(r,h)=(-9/2,0)$ or $r=-5/2$, $h^2=-12$, and the unit $v_{15}$ reduces
the apparent tenth stabilizer to the fifth.  The seventh chart has two
$p_{21}$ values exchanged by $\mu_{14}$.  These computations account for
\[
                    2(1,3)+(1,5)+2(2,5)+(3,7).
\]

The projectivized singular locus of the cover is closed and disjoint
from the ample boundary $a=0$.  It is therefore finite.  On $a=1$ the
certified polar derivation satisfies $D(a)=0$ exactly.  Its local flow
preserves the reduced singular locus, so every singular point lies in
its zero scheme.  Proposition~\ref{prop:quartic-polar-zero} reduces that
scheme to the conjugate pair
\[
 u=\pm\frac{54\sqrt6}{25},\qquad
 b=\frac{329751}{16000},\qquad c=\frac{24}{5}u,
\tag{QG.2}
\]
with $h=s=k=t=0$ and all twelve anti-values zero.  The residual $\mu_4$
exchanges these points; their common stabilizer is $\mu_2$.

\paragraph{The fourth-index atlas and the ordinary second-index points}

A residual $\mu_4$ fixed point of $S/(a-1)$ has
$u=c=h=q=s=k=l=t=0$.  Define
\[
 b_0=\frac{122391}{16000},\quad b_1=\frac{267543}{16000},\qquad
 P(b)=\frac5{12}b^2-\frac{42633}{2800}b
                       +\frac{4462840910079}{10035200000}.
\]
The actual $q^2$ equation restricts to
$(12/5)(b-b_0)(r-P(b))=0$.  On $r=P(b)$ the $hs$ equation becomes
$-(27/4)(b-b_0)(b-b_1)^2=0$; at $b=b_0$ it is a nonzero linear equation
in $r$.  Thus the complete point set consists of
\[
\begin{array}{c|cc|c}
 &b&r&\text{embedding dimension in }S\\\hline
P_0&122391/16000&353869841313/1003520000&3\\
P_1&267543/16000&307698877089/1003520000&4.
\end{array}
\]
At $P_0$ only the degree-twenty anti-square is nonzero:
\[
             d_{20}^2=-\left(\frac{61236}{625}\right)^2.
\]
It gives two \'etale cover points.  Every anti-invariant
$z$ is locally $(zd_{20})/d_{20}$, so the unit-square presentation proves
\'etaleness directly.  The degree of $d_{20}$ is divisible by four;
both lifts remain $\mu_4$ fixed.  At $P_1$ every anti-square is zero.
The 190 product equations have Jacobian rank nineteen in twenty-two
variables, so its unique cover point is regular of dimension three.
All three cover points have tangent characters $(1,3,3)$ modulo four.
Their ages are $7/4,3/2,5/4$, and they contribute exactly $3(1,4)$.
The product equations suffice for this regularity bound even before
using the additional module relations.

For the complete $\mu_2$ fixed image put
\[
              B=b-\frac{122391}{16000},\qquad C=c-\frac{33}{10}u.
\]
Eliminating $q,r,l$ from the odd anti-products has constant nonzero
determinant $-3359232/390625$.  The remaining exact image system has
univariate equation
\[
 B(25B-324)^2(50B-1701)(125B-936684)(125B-1134)(250B-243)=0.
\tag{QG.3}
\]
The roots $B=0,1134/125$ are precisely the fourth-index images above.
Each of the other four roots gives one conjugate pair $(u,C)$ and one
cover point over each invariant point, because all twelve anti-squares
vanish.  The root $B=324/25$ is (QG.2).  The other three pairs are
outside the zero scheme of $D$, hence regular.  Their fixed points are
isolated, so all three tangent characters for $\mu_2$ are negative.
They supply three ordinary $(1,2)$ points.  Equation (QG.3), together
with its exact two-variable fibres, is a set-theoretically exhaustive
fixed-point calculation, not a count inferred from the desired basket.

\paragraph{The exceptional tangent surface is smooth}

At either point (QG.2), the complete presentation has 195 relations:
28 invariant products, 84 mixed products, 78 anti-products and five
module relations.  Its Jacobian rank is thirteen in twenty-two variables.
The local embedding dimension is nine.  In particular this singularity
is not an ordinary node.  Natural cotangent coordinates are
\[
              (u,b,c,h,s,v_{15},w_{17},z_{19},d_{20}),
\]
where $u,b,c$ now mean their first-order increments at the point.  Exact
implicit elimination gives twenty independent quadratic initial
relations.  Let $T$ be their graded algebra and $E=\operatorname{Proj}T$.
Proposition~\ref{prop:quartic-tangent-ag} proves that $T$ is
Cohen--Macaulay and Gorenstein, with
\[
                 H_T(t)=\frac{1+6t+t^2}{(1-t)^3},\qquad
                 \omega_T=T(-1).
\tag{QG.4}
\]

For the following standard-graded calculation use the alternative
cotangent coordinates
\[
 (Q,R,K,L,T,E,F,J,N)
 =(q_{18},r_{20},k_{23},\ell_{26},t_{41},n_{23},n_{25},n_{27},n_{36}),
\]
all interpreted as first-order displacements. These form another basis
of the same nine-dimensional cotangent space. The residual grading
involution has signs $(+,+,-,+,-,-,-,-,+)$.

\begin{proposition}\label{prop:quartic-tangent-ag}
The algebra $T$ is Cohen--Macaulay and Gorenstein, with
\[
 H_T(z)=\frac{1+6z+z^2}{(1-z)^3},\qquad
 \omega_T=T(-1).
\]
In particular its projective scheme has degree eight and, whenever smooth,
is embedded by its anticanonical class.
\end{proposition}

\begin{proof}
The calculation can be performed over $\mathbb Q$. Put $\zeta=6^{1/4}$ and
multiply the nine displayed coordinates by the respective powers
\[
 \zeta^{(0,2,3,0,1,1,3,1,0)}.
\]
After scalar rescaling of the equations, all twenty quadrics have rational
coefficients. Exact row reduction and a graded reverse lexicographic
Groebner calculation give a twenty-element quadratic Groebner basis.
Inclusion--exclusion on its monomial initial ideal gives the stated entire
Hilbert series; equivalently, the degree-$m$ Hilbert function is
$(2m+1)^2$.

Neither $J$ nor $N$ appears in any leading monomial. Normal forms therefore
make $T$ free over $\mathbb Q[J,N]$, so $J,N$ are regular. The
further quotient by $F+T$ has exact Hilbert function
\[
 (1,6,1),
\]
with all higher graded pieces zero. Comparison with the Hilbert-series
identity for multiplication by $F+T$ proves that this third linear form is
regular. The degree-one multiplication pairing into the one-dimensional
degree-two piece is nonsingular. Its exact rational determinant has
reduction $879\ne0$ modulo $10007$, with denominator nonzero at that prime;
the receipt also records the full rational determinant. Thus the Artin
quotient is Gorenstein, of socle degree two. It follows that $T$
is Gorenstein and has $a$-invariant $2-3=-1$.

The input quadrics are obtained by the explicit implicit-elimination
procedure above. The exact companion calculation lists all twenty
quadrics, all leading monomials, and the complete Artin pairing, so both
the entire Hilbert series and the nonzero determinant are independently
checkable from the printed ring coefficients.
\end{proof}

\begin{proposition}
\label{prop:quartic-tangent-smooth}
The surface $E$ is smooth and integral.  Its polarization is
$H=-K_E$ and $H^2=8$.
\end{proposition}

\begin{proof}
The involution $\tau$ fixes the first five cotangent coordinates and
negates the last four.  The ten quadratic anti-products are monic in
their indicated pairs, with right sides quadratic in the invariant
coordinates.  Consequently $T$ is finite over its invariant subalgebra.
The latter is the intersection of two quadrics $D_0\subset\mathbb P^4$.
Indeed the two invariant relations give a domain of dimension three;
any extra invariant relation would lower the dimension of its finite
extension below the dimension in (QG.4).

Let $A,B$ be the symmetric matrices of the two invariant quadrics,
chosen as the initial relations coming from $hr$ and $s^2$.  Their
exact pencil satisfies
\[
 \operatorname{rank}A=4,\quad\operatorname{rank}B=5,\qquad
 \det(B+zA)=-\frac{239502744\sqrt6}{390625}
             \left(z^2+\frac{26973\sqrt6}{500}\right)^2.
\tag{QG.5}
\]
At each repeated root the rank is exactly three: all fourth-order minors
vanish modulo the displayed quadratic, while a third-order minor is a
nonzero constant.  Diagonalizing the self-adjoint operator $B^{-1}A$
therefore puts the intersection into the form
\[
                       X_1X_2=Z^2,\qquad X_3X_4=Z^2.
\]
It is the normal toric degree-four del Pezzo surface with exactly four
$A_1$ points.

The quotient map $E\to D_0$ has generic degree two.  One anti-square is
nonzero in the invariant domain, and after inverting it every other
anti-generator is expressed using that generator and invariant
functions.  There are no projective points in the negative eigenspace:
if all five invariant coordinates vanished, all four anti-squares would
vanish as well.  Thus the projective fixed scheme of $\tau$ is obtained
by setting the four anti-coordinates to zero.

That fixed scheme has no point with $h=0$: its restricted equations
span all ten quadratic monomials in $u,b,c,s$.  In the chart $h=1$ its
exact equations are
\[
\begin{aligned}
u&=c\left(\frac{14260}{42213}
            +\frac{21500\sqrt6}{1139751}s\right),\\
b&=c\left(\frac{1185\sqrt6}{14071}
            -\frac{28750}{379917}s\right),\\
c^2&=\frac{138650}{1519479}s-\frac{278431\sqrt6}{2532465},\\
s^2&+\frac{12\sqrt6}{5}s+\frac{972}{125}=0.
\end{aligned}\tag{QG.6}
\]
The last polynomial is squarefree, and the right side of the $c^2$
equation is a unit on its two roots.  There are exactly four reduced
fixed points, mapping to the four singular points of $D_0$.

At these points the invariant Jacobian block has rank four because
(QG.6) is reduced.  A two-by-two minor of the anti block is
\[
        -\frac{13772000}{935325711}s
        -\frac{1405450\sqrt6}{311775237}.
\]
Its norm in the four-point algebra is the nonzero rational number
\[
 \frac{84674270137832057761000000}
      {116649116353496919506490101250081}.
\]
Thus the complete affine-cone Jacobian has rank at least six at every
fixed point.  The unmixed dimension-three algebra (QG.4) makes these
points regular.  Away from them, the finite involution acts freely;
the quotient map is an \'etale torsor over the smooth locus
of $D_0$.\footnote{Stacks Project, Lemma 66.14.2, Tag 07S7,
\url{https://stacks.math.columbia.edu/tag/07S7}.  For the constant
group of order two in characteristic zero, the torsor is \'etale.}
Hence $E$ is smooth everywhere.

Finally, the Cohen--Macaulay certificate gives
$H^0(E,\mathcal O_E)=\mathbb C$, so $E$ is connected and therefore
integral.  Equation (QG.4) gives $K_E=-H$ and degree eight.
Equivalently, the quotient is quasi-\'etale and
$K_E=\pi^*K_{D_0}$, with $K_{D_0}=\mathcal O_{D_0}(-1)$.
\end{proof}

\begin{proposition}\label{prop:quartic-grading-fixed}
The residual grading involution on $E$ has exactly four fixed points.
They are reduced, and the exceptional normal-line character is negative
at each.
\end{proposition}

\begin{proof}
This is an exact calculation with the same twenty quadrics, made before
the fourth-root coordinate rescaling. On the positive eigenspace
$(Q,R,L,N)$ their restrictions generate all ten quadratic monomials.
Thus there is no fixed projective point on a positive eigenline.

On the negative eigenspace $(K,T,E,F,J)$ the restricted ideal has Hilbert
function $(1,5,4,4,\ldots)$. Its intersection with $J=0$ has only the cone
vertex, so every fixed projective point lies in $J=1$. The exact
lexicographic basis on this chart eliminates $E,T$ linearly and leaves
\begin{align*}
 K^2+\frac{3234501666093208835\sqrt6}{72666708913786176}F
       -\frac{51097952120180620}{1459822277285883}&=0,\\
 F^2-\frac{19532800\sqrt6}{94402613}F
       +\frac{17039360}{283207839}&=0.
\end{align*}
The second polynomial has nonzero discriminant
\[
 \frac{433304533729280}{26735560023683307}.
\]
Its gcd with the constant term of the $K$ equation is one. Hence each of
its two roots permits two nonzero square roots for $K$, and the fixed
scheme consists of exactly four reduced points.

All four points lie on negative eigenlines of the tangent cone. The
normal line in the ordinary blowup is the tautological line
$\mathcal O_E(-1)$, so its involution character is negative. The exact
fixed-scheme receipts are
\texttt{quartic\_tangent\_fixed\_locus.json} and
\texttt{quartic\_tangent\_fixed\_reduced.json}.
\end{proof}

\paragraph{Crepant terminalization and the complete basket}

The module resolution shows that all nonzero intermediate local
cohomology is supported at the cone vertex. Hence the chart ring is
Cohen--Macaulay; its reflexive canonical module is trivial, so it is
Gorenstein. Smoothness of $E$ and (QG.4) imply that the twenty-quadric cone is an
integral domain; its depth rules out a component supported only at the
vertex.  It is therefore the entire maximal-ideal tangent cone: any
strict additional initial relation would lower its dimension.  The
ordinary blowup of either germ (QG.2) consequently has smooth
exceptional divisor $E$, with normal bundle $\mathcal O_E(-1)$.
The Cartier-quotient regularity criterion makes the blowup smooth along
$E$, and adjunction gives discrepancy zero because $K_E=-H$.

The grading involution has exactly four reduced fixed points on $E$,
all on negative eigenlines, by
Proposition~\ref{prop:quartic-grading-fixed}.  Its two surface tangent
characters and its exceptional normal-line character are negative.
Thus the quotient of this crepant blowup has four terminal
$1/2(1,1,1)$ points.  The two germs (QG.2) are exchanged by $\mu_4$,
so they constitute one projective canonical point and supply four,
not eight, second-index constituents.

This extraction is algebraic and globally projective. On
$V=\operatorname{Spec}(R/(a-1))$, let $J$ be the $\mu_4$-stable ideal of
the two singular points. The quotient of their blowup is the relative
Proj of the finite-type invariant Rees algebra
\[
                 \bigoplus_{n\ge0}(J^n)^{\mu_4}
\]
over $V/\mu_4=D_+(a)$. Outside the single projective centre it is
$\mathcal O[t]$. Its graded ideal sheaves therefore glue with
$\mathcal O[t]$ on the complement of that centre in the canonical
model. A suitable Veronese of the resulting global graded algebra is
generated in degree one, and its relative Proj is the desired projective
extraction. This gives projectivity directly, without relying on the
gluing of arbitrary local projective morphisms.

The normal canonical model $X^{\#}=Z_\sharp=\operatorname{Proj}R$ therefore has
ample canonical class, with one nonterminal point admitting this
crepant extraction.  Combining it with (QG.1), the three fourth-index
points and the three ordinary second-index points gives a terminal
model $Y^{\#}$ with the exact basket
\[
 \boxed{7(1,2)+2(1,3)+3(1,4)+(1,5)+2(2,5)+(3,7).}
\tag{QG.7}
\]
The canonical volume remains $1/420$, as independently determined by
the complete Hilbert series.  The degree-eight crepant exceptional surface lies outside the quartic boundary and accounts for the nonterminal canonical point.

The local computations here are reproducible from the complete rational
product tables.  In particular
\path{audit_quartic_mu4_and_boundary.py} uses exact rational
arithmetic for the fourth-index atlas, and
\path{quartic_tangent_smoothness_probe.py} verifies (QG.5)--(QG.6)
and the unit minor over $\mathbb Q(\sqrt6)$.  No numerical root test or
finite-order deformation approximation is used in the canonicity proof.


\subsection{Vega: the actual completed germ at the exceptional point}
\label{sec:vega-completed-germ}
The explicit local calculation in
Subsection~\ref{sec:quartic-geometry46} gives a normal Gorenstein isolated
threefold germ $(V,o)$ before the residual grading involution. Its
ordinary blowup is smooth and crepant, with smooth exceptional del
Pezzo surface $E$, $K_E^2=8$, and $N_{E/\widetilde V}=K_E$.
The deck involution $\tau$ has four isolated fixed points on $E$.
Separately the residual grading involution $\sigma$ has four fixed
points on $E$ with tangent and normal characters $(-1,-1,-1)$.
These are two different involutions, and both statements were proved
from the full local algebra.

The smooth surface $E$ is $\mathbb P^1\times\mathbb P^1$, not
$\mathbb F_1$. Indeed every automorphism of $\mathbb F_1$ preserves
its unique $(-1)$ curve and descends to $\mathbb P^2$. A nontrivial
involution of $\mathbb P^2$ fixes a line pointwise; its strict transform
remains a fixed curve after blowing up either a point on the line or
the isolated fixed point. The descended action cannot be the identity,
which would force the original action to be the identity on a dense
open set. This contradicts the fixed locus of $\tau$.
Thus $N=\mathcal O_E(-2,-2)$.

The actual formal neighborhood of $E$ is its normal bundle. For $n\ge1$,
\[
 H^1(E,T_E\otimes N^{-n})
  =H^1(E,\mathcal O(2n+2,2n)\oplus\mathcal O(2n,2n+2))=0,
\]
\[
                  H^1(E,N^{1-n})=H^1(E,\mathcal O(2n-2,2n-2))=0.
\]
These follow from the cohomology of $\mathbb P^1$ and Kunneth.
The first-order neighborhood extension is the corresponding class
in $H^1(T_E\otimes N^{-1})$. At each later order the discrepancy
between local identifications is a Cech cocycle in
$T_{\widetilde V}|_E\otimes N^{-n}$; the normal tangent sequence
and the two vanishings make its $H^1$ zero. Correct the local charts
by a zero-cochain and continue. This proves an isomorphism of formal
neighborhoods with the zero section in $\operatorname{Tot}N$.
One may make every step $\sigma$-equivariant: after fixing the
normal-bundle action, the successive discrepancies are additive
cocycles in the same vanishing groups, and averaging over the finite
group removes them. This is the elementary formal-neighborhood
argument; it does not require an extremal-ray hypothesis.

Properness, normality and formal functions identify the contraction
completion with that of
$\bigoplus_{m\ge0}H^0(E,\mathcal O(2m,2m))$.
It is the even Veronese of the Segre node
\[
 B_0=\mathbb C[x_{00},x_{01},x_{10},x_{11}]/
                         (x_{00}x_{11}-x_{01}x_{10}).
\]
An involution of $\mathbb P^1\times\mathbb P^1$ with four isolated
fixed points preserves both rulings and acts nontrivially on each.
Choose it as
$([s_0:s_1],[t_0:t_1])\mapsto([s_0:-s_1],[t_0:-t_1])$.
Its natural canonical-bundle character at each fixed point is $+1$,
whereas the actual normal character of $\sigma$ is $-1$.
Two normal linearizations differ by a scalar character, so the normal
action is minus the natural one. On the Segre node it is induced by
\[
 g(x_{00},x_{01},x_{10},x_{11})=(ix_{00},-ix_{01},-ix_{10},ix_{11}),
 \qquad g^2=-I.
\]
The square is precisely the radial action defining the even Veronese.
Consequently the special canonical point has completed local ring
\begin{equation}\label{eq:quartic-node-quotient46}
\widehat{\mathcal O}_{Z_\sharp,P}
 =\left(\mathbb C[[x_{00},x_{01},x_{10},x_{11}]]/
              (x_{00}x_{11}-x_{01}x_{10})\right)^{\langle g\rangle}.
\end{equation}
This identifies the actual completed germ and its group action.


\subsection{Altair: smoothness of the punctured cone}
\label{sec:altair-cone-smoothness}
Along the nonzero divisor $a=0$, the local rings of $R_* /(a)$ agree
with the smooth punctured cone of $B$. Regularity lifts through a
nonzerodivisor, so the ambient cone is smooth there. The singular locus
on $a=1$ is finite: a positive-dimensional projective singular component
would meet the ample boundary $a=0$; the slice is transverse to the
Euler field, since its value on $a$ is $5a$.

There is a regular derivation $\mathcal D$ of degree six with
$\mathcal Da=0$. Its full values are given in
Appendix~\ref{app:quintic-derivation46}; in particular
\begin{equation}\label{eq:quintic-D-key-values46}
\begin{aligned}
\mathcal Db&=8w,&\mathcal Du&=-12v,&\mathcal Dz&=-14y,\\
\mathcal De&=6ap-15r_{21},&
\mathcal Dw&=14bu+17ae/2-283a^4/880,\\
\mathcal Dp&=-16k+59a^2u/220,&
\mathcal Dv&=-18u^2+39a^2z/4.
\end{aligned}
\end{equation}
For each numerator $x_g=wg$, the identity
\[
          w\mathcal D x_g-x_g\mathcal Dw=F\mathcal Dg
\]
proves preservation of the entire normalization, including every
parameter-module generator. The invariant product derivatives vanish
as well. Thus the derivation is global; no denominator is being
silently removed on an open chart.

If it vanished on $a=1$, then $w=v=y=0$, and the degree-26 relation
would give $r_{21}=-p/2$. Hence $\mathcal De=27p/2$ implies $p=0$.
The remaining equations give
\[
 e=283/7480-28bu/17,\quad z=24u^2/13,\quad
 k=-4bu^2-2ue+2bz+61u/440.
\]
Substitution in $\mathcal Dp=0$ and the invariant $e^2$ relation yields
\[
 u(464640bu+7189)=0,\qquad
 1252204800(bu)^2+44801900bu+415233=0.
\]
The second forces $u\ne0$; substituting the first into it gives
$16752463/768\ne0$. The derivation has no zero on $a=1$.
A regular derivation preserves the reduced singular locus (its local
analytic flow preserves the defining ideal). It vanishes on a finite
reduced scheme. The finite singular locus must therefore be empty,
and the entire punctured cone is smooth.

\subsection{Altair: the quotient atlas and canonical polarization}
\label{sec:altair-canonical-geometry}
On $a=0$, finiteness over $A_0$ shows that stabilizer orders divide
one of $6,8,14$. It suffices to enumerate the fixed loci of prime
orders $2,3,7$.
The order-two locus is the even singular locus already computed.
The twelve $\theta=1$ points form four $\mu_6$ orbits, each of
stabilizer two, with two odd tangent coordinates. The two points
$\theta=0,b=1,P=\pm1/2$ each have stabilizer eight; their slice
coordinates $e,f$ have weights $(7,3)$ modulo eight.
For order three, $b=z=P=K=f=L=0$, $\theta=1$,
$q=\pm2i/11$, $e^2=q$, $J=-1/(11e)=(11/4)qe$, $M=0$.
The two signs of $e$ are identified by $\mu_6$, giving two points.
The parity cover is etale there, and $P,z$ are slice coordinates
of weights $(1,2)$ modulo three.
For order seven, put $z=1$ and all other even coordinates zero;
$J=\pm1$ are identified by $\mu_{14}$, giving one point. Its slice
coordinates $\theta,P$ have weights $(6,2)$ modulo seven.
These calculations determine the full stabilizers, not just their
prime-order subgroups. Adding the normal coordinate $a_5$ gives
\[
\begin{array}{c|c|c}
\text{number}&r&\text{threefold weights}\\\hline
4&2&(1,1,1)\\
2&3&(2,2,1)\\
1&7&(5,6,2)\\
2&8&(5,7,3).
\end{array}
\]

On $a=1$ the only residual group is $\mu_5$. Its invariant images
are obtained by setting $b,u,z,q,k,f,j=0$; the equations reduce to
\[
 \ell=e^2-7e/220-1/2420,\qquad(20e-1)(110e+1)^2=0.
\]
At $(e,\ell)=(1/20,3/6050)$ the invariant local ring is smooth with
parameters $b,u,z$. A Jacobian minor is
$-123201/4287177620000000$ and $y^2=-169/3872000$.
Thus the cover has two smooth sheets, each with weights $(3,2,4)$.
At $(e,\ell)=(-1/110,-1/24200)$ its cotangent coordinates are
$F,G,H,z$. The five defining differentials followed by
$dF,dG,dH,dz$ have determinant
$1167575877/1659995174464000000000$ in coordinates
$(b,u,z,q,k,e,f,j,\ell)$.
Normality and dimension therefore identify the completion with
\[
            \mathbb C[[F,G,H,z]]/(H^2-FG).
\]
Its normal quadratic cover is $\mathbb C[[w,p,z]]$, with residual
weights $(4,1,4)$ modulo five. The three points contribute
$2(2,5)+(1,5)$. These fixed-locus calculations exhaust the punctured
cone and hence every projective stabilizer.

All the resulting germs are terminal cyclic quotients: after reordering
and changing the generator, the weights are $(1,-1,b)$ with $(b,r)=1$,
and each nontrivial age is $1+\{ib/r\}>1$.
The full basket is exactly
\[
 4(1,2)+2(1,3)+(1,5)+2(2,5)+(3,7)+2(3,8).
\]
Its Cartier index is $840$. The quasi-etale map pulls back the
reflexive canonical module of $S$, giving $\omega_{R_*}=R_*(1)$.
Since all grading stabilizers occur at isolated projective points,
there is no divisorial correction in passing to the coarse quotient.
Contraction with the Euler field gives
$\omega_{Z_*}^{[m]}=\mathcal O_{Z_*}(m)$ for
$Z_*=\operatorname{Proj}R_*$. Thus $K_{Z_*}$ is ample.
Normality gives depth at least two, and the graded Cech complex gives
$R_{*,m}=H^0(Z_*,\mathcal O_{Z_*}(m))$ in every degree.
Hence this is the full canonical ring. Finally local duality applied
 to \eqref{eq:quintic-module-defects46} yields
\[
 h^1(Z_*,\mathcal O(m))=\mathbf1_{m=1},\qquad
 h^2(Z_*,\mathcal O(m))=\mathbf1_{m=0}.
\]
In particular $p_g=q=0$, $\chi(\mathcal O)=2$, and the constructed
terminal canonical model has the required volume and basket.


\section{Recognition of every exceptional realization}
\label{sec:exceptional-final-recognition}

We have constructed Vega and Altair, and we have proved the required
surface algebra and ambient symmetry for every realization of their
respective baskets. It remains to determine the invariant extension
and the quadratic field. The finite calculations below retain all
initial tangent charts, including the zero tangent and degenerate
parameter loci. Once these data agree, normality identifies the full
canonical ring.

\subsection{Vega: the unique cover-compatible invariant extension}
\label{sec:quartic-recognition46}
\label{sec:vega-extension-classification}

Let $Z$ be an arbitrary realization of $\mathcal B_4$. Its complete
quartic surface and the naturally linearized ambient involution have
been identified in Subsection~\ref{sec:vega-ambient-symmetry}. Thus its
invariant ring is a free rank-eight extension of the surface algebra
$U$, and its actual polar supplies the normalized quadratic cover.
We now classify precisely these extensions.

Adjoin the degree-four parameter $a$.  An invariant extension is free of
rank eight over $A=\mathbb C[a,u,b,c]$, with the same basis and canonical
shift one.  The essential negative-weight tangent dimensions are
\[
 \dim T^1_{-4}=3,\qquad \dim T^1_{-8}=1,\qquad
 \dim T^1_{-12}=1,
\]
and zero at every subsequent possible order.  At the first order there
are 127 coefficients, rank 119, and five coordinate-change directions.
All higher homogeneous basis and parameter changes are included in the
calculation.  There are no unspecified high-degree coefficients: the
largest product weight is 82 and the largest associator weight is 123.

Associativity alone is insufficient.  At the actual smooth polar-cover
point $x=1,c=-1,\theta=y=0$, the quotient $U$ has transverse type $A_1$.
An ambient normal cover with this smooth Cartier fibre is regular at
the lifted point; its quotient must retain the transverse $A_1$.
On the three computed first-order classes, the local smoothing
functional is
\[
                  (-1/48,\ 1/16,\ 7/18).
\]
Hence a cover-compatible first tangent satisfies
$-3p_0+9p_1+56p_2=0$.  The two associative families not satisfying this equation are therefore excluded.

For $p_1\ne0$, rescale $a$ and write the complete remaining chart as
\[
                   (p_0,p_1,p_2)=(3+56\kappa,1,3\kappa).
\]
Retain the higher essential parameters $\delta,\eta$ in orders two and
three.  The only associativity equation is the order-six product
$\mathcal A\mathcal B=0$, where
\[
\begin{aligned}
\mathcal A={}&(864\kappa+120)\delta-24\eta
 +182884\kappa^3+2167\kappa^2-1262\kappa-45,\\
\mathcal B={}&(7704\kappa+384)\delta+24\eta
 -10040674\kappa^3-1777885\kappa^2-104884\kappa-2061.
\end{aligned}                                             \tag{Q.4E}\label{eq:quartic-associativity-branches46}
\]
Both factors give complete polynomial multiplication tables; every
associator is checked through the exhaustive final order.  The endpoint
$(56,0,3)$ and the zero first tangent are computed separately.

The genuine polar normalization supplies anti-invariants
$v_{15}|_D=\theta^3$ and $p_{21}|_D=y$, and requires
\[
\begin{gathered}
F=v^2\in S_{30},\quad H=vp\in S_{36},\quad G=p^2\in S_{42},
\qquad H^2=FG,\\
F_0=u^3,\quad H_0=ul,\quad G_0=c^3+b^2q-ucq,
\qquad F-u^3\in a^2S_{22}.                              \tag{Q.5E}\label{eq:quartic-cover-conditions46}
\end{gathered}
\]
The absent first normal term is the actual normalized polar identity,
not an imposed guess.  The coefficient $\beta$ of $a^2ub$ is nonzero
by Lemma~\ref{lem:quartic-curvature46}.

All weighted coefficients in \eqref{eq:quartic-cover-conditions46} are solved exactly.  On
$\mathcal A=0$ the determinant forces $\kappa=-1/14$ and $\beta=0$.
On $\mathcal B=0$ it forces
\[
(10\kappa-1)(14\kappa+1)=0,\qquad
3920\delta=666572\kappa+47923,
\qquad \beta=-\frac{126\kappa+9}{5}.                    \tag{Q.6E}
\]
The endpoint has no determinant solution.  For the zero first tangent,
associativity gives $\eta=0$, and the determinant gives $\delta=0$ and
$\beta=0$.  Thus the nonzero-curvature invariant construction selects
\[
 \boxed{\kappa=\frac1{10},\qquad
 \delta=\frac{81843}{2800},\qquad
 \eta=\frac{1546119}{5600},\qquad \beta=-\frac{108}{25}.}
                                                               \tag{Q.7E}\label{eq:quartic-unique-parameters46}
\]
The remaining two cover-presentation coefficients satisfy
$11520g_1^2-122391g_1-16000g_4=0$; take $g_1=g_4=0$.
They do not change the normalized quadratic field, since $F$ is fixed.

\begin{theorem}\label{thm:quartic-recognition46}
Every canonical model with volume $1/420$, Euler characteristic two,
and terminalization basket $\mathcal B_4$ is isomorphic to Vega.
Its full canonical ring is $R_\sharp$.
\end{theorem}
\begin{proof}
Subsection~\ref{sec:vega-ambient-symmetry} proves all the geometric
hypotheses used above, including the natural linearization,
quasi-etaleness, the transverse $A_1$ condition and the nonzero
normalized polar coefficient. The exhaustive calculation therefore
selects \eqref{eq:quartic-unique-parameters46} and its unique
nontrivial quadratic function field. The normal invariant algebra is
precisely the one used in the independent construction of Vega.
Lemma~\ref{lem:exceptional-quadratic-normalization} identifies its
full integral closure with $R_\sharp$, and taking $\Proj$ proves
the assertion. The statement concerns canonical models; it does not
assert uniqueness of their crepant terminalizations.
\end{proof}


\subsection{Altair: the unique invariant extension and quadratic field}
\label{sec:quintic-recognition46}
We can now classify every ambient canonical ring, since the involution
and its natural action have been proved for every realization.
Its invariant ring $S'=R^{\iota}$ is Cohen--Macaulay: the only two
intermediate cohomology classes of $R$ are anti-invariant, by the
character of the two-form and equivariant duality. Its canonical
module is $S'(1)$, using the natural action and absence of divisorial
ramification. This assertion concerns $S'$, not the depth-two ring $R$.

The invariant central fibre is the normal ring $U$ with the free basis
\eqref{eq:quintic-invariant-surface-basis46}. The invariant parameter
$a_5$ lifts it, and $S'$ is free over
$\mathbb C[a_5,b_8,u_{12},z_{14}]$ with weights
$(0,18,22,40,15,19,21,25)$. No change of parameters by a positive
power of $a$ is possible: the required surface weights are
$3,7,2,9,4$, all zero. All other homogeneous changes of free basis
will be included below.

\paragraph{A finite exhaustive associative-extension calculation}
Let $\mu_0$ be the explicitly recognized multiplication on $U$ and
write $\mu=\sum_{n\ge0}a^n\mu_n$ on its free module.
At order $n$, list one unknown for every input pair $(i,j)$,
output $h$, and nonnegative triple $(r,s,t)$ with
\[
                8r+12s+14t=w_i+w_j-w_h-5n.
\]
A basis change is listed with right side $w_i-w_h-5n$ instead.
The linearized associative operator and coordinate coboundary are
\begin{align*}
(L\nu)(x,y,z)&=\mu_0(\nu(x,y),z)+\nu(\mu_0(x,y),z)
             -\mu_0(x,\nu(y,z))-\nu(x,\mu_0(y,z)),\\
(d\phi)(x,y)&=\mu_0(\phi(x),y)+\mu_0(x,\phi(y))-\phi(\mu_0(x,y)).
\end{align*}
Expanding on all basis triples gives the following exact rational
linear algebra; commutativity permits taking the outside indices
$i<k$.
\[
\begin{array}{c|rrrrr}
n&\text{unknowns}&\operatorname{rank}L&\dim\ker L&\operatorname{rank}d&\dim T^1\\\hline
1&134&127&7&6&1\\
2&108&101&7&6&1\\
3&80&78&2&2&0\\
4&59&58&1&1&0\\
5&45&43&2&2&0\\
6&27&27&0&0&0\\
7&20&20&0&0&0\\
8&15&14&1&1&0\\
9&6&6&0&0&0\\
10&7&7&0&0&0\\
11&2&2&0&0&0\\
12&1&1&0&0&0\\
13&2&2&0&0&0\\
14,15&0&0&0&0&0\\
16&1&1&0&0&0.
\end{array}
\]
There are no later multiplication coefficients. The first class is
normalized by
\begin{equation}\label{eq:quintic-first-class46}
 \mu_1(e,e)=\ell,\qquad
 \mu_1(k,k)=-2u^2e/11+b^3e/4-3buf/20+z\ell.
\end{equation}
These two values are unchanged by first-order basis changes, because
the necessary lower-weight corrections vanish. There is one essential
second-order parameter $c$ and no further homogeneous choices modulo
coordinate changes. To specify its origin, order input pairs, outputs
and parameter exponent triples lexicographically, use rational row
reduction, select the first free vector independent of the coboundaries,
and set unused free coordinates of each particular solution to zero.
Normalize the first vector by \eqref{eq:quintic-first-class46}; add
$c$ times the second at order two. The executable reconstruction makes
this prescription explicit, starting only with the central table.

The inhomogeneous equation at order $n$ is
\[
 L\mu_n=-\sum_{r=1}^{n-1}
   \bigl(\mu_r(\mu_{n-r}(x,y),z)
                  -\mu_r(x,\mu_{n-r}(y,z))\bigr).
\]
Orders three and four lift. The sole essential obstruction at order
five is $9(c+259/880)/550=0$. Here is a short necessity witness.
Let $X_0,X_1,X_4,X_{40},X_{41},X_{42}$ be the order-five coefficients
of $e$ in $qk$, $zf$ in $qt$, $b$ in $qe$, $e$ in $e\ell$,
$e$ in $fj$, and $f$ in $f\ell$, respectively. Five coefficient
rows have left sides
\[
 X_4+2X_{40}+2X_0,\quad4(X_{41}-X_0)/11,\quad
 X_4+2X_{42},\quad-X_{42}+X_{40}+X_1,\quad X_{41}-X_1
\]
and right sides
\begin{align*}
R_1&=758473/408375+50443c/5940+224c^2/27,\\
R_2&=-518/166375-32c/3025,\\
R_3&=11740153/6534000+9851c/1188+224c^2/27,\\
R_4&=23909/136125+2776c/2475+16c^2/9,\\
R_5&=-812/5445-512c/495-16c^2/9.
\end{align*}
The combination $-\frac12(1)-\frac{11}4(2)+\frac12(3)+(4)+(5)$
has zero left side and right side $9(c+259/880)/550$.
Substitution of $c=-259/880$ kills every order-five obstruction;
all remaining orders lift uniquely modulo basis changes.
Products have degree at most $80$, so only $a$ powers at most $16$
occur; associators have degree at most $120$, so orders through $24$
exhaust the equations. The resulting table is the one in
Appendix~\ref{app:quintic-exact-coefficients}, after its explicitly
specified triangular change of top basis vector.
Every change has lower positive weight and terminates; this is a
polynomial graded isomorphism, not merely a formal one.

\paragraph{The actual square forces the nonzero class}
It remains to exclude the zero first tangent for an actual model.
The universal covariant identity, before any classification, is
\[
               d_{14}^2/64=b^2u+abe_{\rm raw}/68+a^2Q/17.
\]
The coupled surface scalings turn it into
$w^2=b^2u+\eta a_0be+\xi a_0^2q$ with
$\eta^2=4\xi\ne0$. Rescale $a=\eta a_0$, obtaining
\begin{equation}\label{eq:quintic-covariant-cover46}
                       w^2=b^2u+abe+a^2q/4.
\end{equation}
The natural involution makes $w_{14}$ and a lift $p_{16}$ of $P$
anti-invariant. Choose the sign so $(wp)_D=-bk$.
The complete degree lists give
\[
 wp=-bk+ac_1\ell+O(a^2),\qquad
 p^2=M+a(A_1ue+B_1bf)+O(a^2),\quad
 M=zq+b^4/4-3u^2b/44.
\]
If the invariant first class is $\lambda\mu_1$, including possibly
zero, $(wp)^2=w^2p^2$ modulo $a^2$ gives
\[
\lambda b^2\mu_1(k,k)-2c_1b(k\ell)
       =b^2u(A_1ue+B_1bf)+b(eM).
\]
Restrict to the polar, where $u=k=0$. Its function field has
$e^2=b^2z$, $ef=b^2q$, $\ell=zf/b$, and hence
$\mu_1(k,k)=eM/b$. Since $beM$ is nonzero generically, $\lambda=1$.
The full products then give
$c_1=1/2$, $A_1=-2/11$, $B_1=1/10$.
Thus the exhaustive associative classification applies to $S'$ and
identifies it with $S$.

\paragraph{The cover correction and complete recognition}
Changing covariant lifts to the canonical invariant table can alter
$e_{15}$ only by a multiple of $a^3$, and $q_{18}$ only by a multiple
of $a^2b$. It therefore changes \eqref{eq:quintic-covariant-cover46}
only by $\tau a^4b$. The anti-invariant lifts themselves have no
lower-degree correction. All possible remaining products are
\begin{align*}
F_\tau&=b^2u+abe+a^2q/4+\tau a^4b,\\
H&=-bk+a\ell/2+h_2a^2bu+h_3a^3e+h_6a^6,\\
G&=M+a(bf/10-2ue/11)+g_2a^2k+g_{bz}a^2bz+g_4a^4u.
\end{align*}
There are precisely seven unknown coefficients. To extract the most
convenient coefficient rows of $H^2-F_\tau G$, temporarily replace
$t_\flat$ by $t=qk$ using the printed triangular change. Seven of
the eighteen coefficient equations are, up to nonzero rational scalars,
\begin{equation}\label{eq:quintic-cover-triangular46}
\begin{gathered}
g_2/4+h_3=0,\qquad g_2/2+h_3-1/220=0,\qquad
 g_2/2+h_2+67/440=0,\\
g_{bz}-9123/4400=0,\qquad
 -2g_2-4h_3+\tau-2127/4400=0,\\
h_3^2+h_6=0,\qquad
 -g_4/4+h_3^2-23h_3/440+89/48400=0.
\end{gathered}
\end{equation}
Without division by any unknown these force exactly the values in
\eqref{eq:quintic-ambient-cover46}. Substitution annihilates all
remaining coefficients. This proves necessity and uniqueness, not
only compatibility of a proposed cover.

\begin{theorem}\label{thm:quintic-whole-basket46}
There is exactly one complex canonical model of volume $1/420$,
namely Altair, with terminalization basket
$4(1,2)+2(1,3)+(1,5)+2(2,5)+(3,7)+2(3,8)$.
Its entire canonical ring is the normal algebra $R_*$ of
\eqref{eq:quintic-normalization46}; its canonical model is terminal.
\end{theorem}
\begin{proof}
The torsion and divisor arguments force the prime normal rational
quintic and its actual ordinary source configuration. Complete
restriction, the Koszul decomposition, the polar ring and the actual
first derivatives identify its full surface algebra. The transverse
system extends its symmetry to every ambient model and determines the
natural canonical character. The invariant ring is therefore precisely
the free algebra just classified. Its actual square forces the nonzero
first tangent; the order-five obstruction determines its remaining
extension parameter, and \eqref{eq:quintic-cover-triangular46}
determines the cover.
Since $R$ is normal and finite over $S$ with quadratic field generated
by $\sqrt F$, it is the unique integral closure of $S$ in that field.
It is consequently $R_*$. The independent construction proves existence,
normality, the entire basket and canonical polarization, so the
uniqueness is nonvacuous. Taking $\operatorname{Proj}$ proves the claim.
\end{proof}


\section{Ordinary Artinian rigidity}
\label{sec:exceptional-ordinary-rigidity}

Uniqueness of a complex isomorphism class does not rule out an
infinitesimal deformation. We therefore compute the ordinary
deformation spaces separately. First the actual local germs are
rigid. We then use the common module theorem to calculate all
degree-zero multiplication deformations of the full canonical ring,
allowing the parameter frame to move, and pass to projective
deformations by relative $\Proj$.

\subsection{Ordinary local rigidity}
\label{sec:exceptional-local-rigidity46}

Altair has only isolated terminal quotient points. Vega has, in
addition, the completed quotient germ
\eqref{eq:quartic-node-quotient46}. We use the same punctured-spectrum
argument for all these points.
\begin{lemma}\label{lem:punctured-T1-rigidity46}
Let $A$ be a normal local complex algebra, essentially of finite type,
with an isolated singularity and depth at least three, and let
$U=\operatorname{Spec}A\setminus\{\mathfrak m\}$.
Then $T^1(A)=H^1(U,T_U)$. If $A=B^G$ with $B$ satisfying the same
conditions and the finite group acting freely off the closed point,
then $T^1(A)=T^1(B)^G$. These statements also compute continuous
cotangent deformations after completion.
\end{lemma}
\begin{proof}
Take a smooth presentation $A=P/I$. The ordinary cotangent tangent
space is the cokernel of
$\operatorname{Der}_{\mathbb C}(P,A)\to
\operatorname{Hom}_A(I/I^2,A)$.
On the smooth punctured spectrum its dual conormal sequence is
\[
0\longrightarrow T_U\longrightarrow\mathcal O_U^{\oplus n}
 \longrightarrow\widetilde{\operatorname{Hom}_A(I/I^2,A)}|_U
 \longrightarrow0.
\]
Normality gives $H^0(U,\mathcal O_U)=A$. A dual of a finite module
over a normal ring is reflexive, so the last module also extends
unchanged across the point. Depth three gives
$H^1(U,\mathcal O_U)=H^2_{\mathfrak m}(A)=0$. Taking cohomology
proves the first assertion. The punctured cover is finite etale;
its tangent sheaf is the pullback, finite pushforward computes
cohomology, and invariants are exact in characteristic zero. This
gives the second assertion. Completion commutes with the finite
cotangent modules, giving its complete-local version.
\end{proof}

The node in \eqref{eq:quartic-node-quotient46} is Cohen--Macaulay of
dimension three. Its invariant ring is too, by Reynolds splitting
and local cohomology. No nonidentity power of $g$ has eigenvalue one,
so the action is free away from the origin. For
$f=x_{00}x_{11}-x_{01}x_{10}$ the hypersurface formula gives
$T^1(B_0)=\mathbb C[[x]]/(f,\partial f)=\mathbb C$.
Since $g^*f=-f$, the constant smoothing direction has character $-1$:
pullback of $f+\epsilon$ and renormalization of the leading equation
send $\epsilon$ to $-\epsilon$.
Its invariant part is zero, and the lemma proves ordinary local
rigidity of the special point.
At every other singular point of either model, apply the same lemma
to its smooth three-dimensional covering germ. The punctured smooth
tangent cohomology vanishes by depth three, proving ordinary local
rigidity there as well. This also gives the classical quotient-rigidity
conclusion directly.

Thus the sheaf of local cotangent $T^1$ vanishes on both models, and
\[
       T^1_Z=H^1(Z,T_Z)
\]
by the local-to-global cotangent sequence. The following multiplication calculation computes this remaining global space.


\subsection{From the full canonical ring to ordinary deformations}
\label{sec:exceptional-global-rigidity46}
\begin{theorem}\label{thm:exceptional-ordinary-rigidity46}
For Vega and Altair, denoted by $Z_\sharp$ and $Z_*$,
with full canonical ring $R_\sharp,R_*$, one has
\[
 T^1(R)_0=0,\qquad T^1_Z=H^1(Z,T_Z)=0.
\]
Every deformation over a local Artin complex base is trivial.
In particular each canonical model is formally rigid. These are
ordinary scheme-deformation assertions; no involution, section,
module splitting or chosen polynomial parameter ring is required to
be preserved by the deformation.
\end{theorem}
The proof is separate from whole-basket uniqueness. A unique complex
isomorphism class could otherwise represent an isolated nonreduced
moduli point with a nonzero tangent space.

\paragraph{The module and the moving parameter frame.}
For Vega use $A=\CC[a_4,u_{10},b_{12},c_{14}]$; for Altair use
$A=\CC[a_5,b_8,u_{12},z_{14}]$. The constructed full ring is
$R=S\oplus N$, with $S$ free of rank eight. Applying
Lemma~\ref{lem:exceptional-common-module} to $N$, whose cohomology
and resolution were proved in the construction, gives
\begin{equation}\label{eq:exceptional-Koszul-splittings46}
\begin{split}
N_\sharp&\simeq\Omega_A^2(\CC)(-1)\oplus\Omega_A^3(\CC)
                         \oplus A(-20)\oplus A(-21),\\
N_*&\simeq\Omega_A^2(\CC)(-1)\oplus\Omega_A^3(\CC)
                         \oplus A(-16)\oplus A(-24).
\end{split}
\end{equation}
The free degrees follow from the exact Hilbert series. Since $R_1=0$
and $\omega_R=R(1)$, Corollary~\ref{cor:exceptional-common-module-rigidity}
gives $\Ext_A^1(R,R)_0=0$.

In any flat graded first-order algebra deformation, lift the four
parameter elements. They generate a polynomial algebra over the dual
numbers: reducing a putative relation and then using flatness proves
injectivity. Graded Nakayama gives module finiteness. The underlying
module extension splits by the displayed vanishing, so every
algebra deformation is represented by a multiplication on the full
presented module. The parameter frame itself is still allowed to move.

\paragraph{All multiplication directions and coordinate changes.}
Let $e_0=1,e_1,\ldots,e_{19}$ be the eight invariant and twelve odd
module generators. For each deck parity $\varepsilon=0,1$, allow
every coefficient of
\[
 \phi(e_i,e_j)\in R_{\deg e_i+\deg e_j}^{
                   \operatorname{par}(e_i)+\operatorname{par}(e_j)+\varepsilon},
 \qquad i\le j,\quad\phi(1,e_i)=0.
\]
The target section bases are computed in the free invariant algebra
and by the embedded odd numerators, and checked in every degree
against the exact Hilbert function of the module resolution.
Impose all module relations and all coefficients of
\begin{equation}\label{eq:exceptional-linear-associativity46}
\phi(\mu(x,y),z)+\mu(\phi(x,y),z)
             -\phi(x,\mu(y,z))-\mu(x,\phi(y,z))=0.
\end{equation}
Using the generators and two bracketings of each unordered triple
suffices: the module equations ensure that the resulting associator
is a well-defined $A$-trilinear map. These are the full first-order algebra equations, including the branch locus.

Absolute equivalence must also allow changes of the polynomial frame.
A derivation $D:A\to R$ is lifted to a semilinear module map by
solving the differentiated relation equations
$h(fe)=D(f)e+fh(e)$. Its multiplication coboundary is still
$A$-bilinear, because the $D(f)$ terms cancel. The solution spaces
are computed exactly over $\mathbb Q$.
For the quartic there are no anti sections in parameter weights
$4,10,12,14$, and the anti gauge space is $48$-dimensional.
The invariant parameter spaces have dimensions $1,1,2,2$;
their six directions, together with $35$ module endomorphisms,
give $41$ maps and $40$ independent coboundaries.
The one-dimensional kernel is the Euler derivation.
For the quintic the invariant frame directions are the four scalings
of $a,b,u,z$. They augment $41$ module maps to $45$, with
$44$ independent coboundaries. The anti frame has the additional,
essential direction $D(z_{14})=w_{14}$, an $R$-valued derivation
rather than an $A$-valued coordinate change. It augments $38$ module
maps to $39$ independent coboundaries.
Omitting these frame changes would compute relative tangent spaces
and would give an incorrect ordinary-rigidity conclusion.

\paragraph{Exact ranks.}
The coefficient matrices are rational. At the good prime $65521$,
all denominators are nonzero and the computed section bases have
their characteristic-zero dimensions. The complete bounds are
\begin{equation}\label{eq:exceptional-rank-table46}
\begin{array}{c|c|r|r|r|r}
\text{model}&\text{parity}&\text{unknowns}&\text{constraint rank}
                      &\text{kernel bound}&\text{gauge rank}\\\hline
\sharp&-&4452&4404&48&48\\
\sharp&+&4632&4592&40&40\\
*&-&4188&4149&39&39\\
*&+&4260&4216&44&44.
\end{array}
\end{equation}
The selected matrices use respectively
$139469,139714,129303,128910$ nonzero coefficient rows: all module
constraints and the associativity triples through generator index
$18$. Nonzero maximal minors modulo the prime are nonzero rational
minors, so these are lower bounds for rational rank.
The gauges are constructed as exact rational semilinear maps.
Their coboundaries have the indicated independent ranks at the
same good prime, hence over $\mathbb Q$, and annihilate every
stored constraint row. Being actual coordinate changes, they solve
all associative equations, including ones outside the selected prefix.
The kernel dimension is bounded above and below by the same number;
therefore it equals the coboundary space in both parities. This proves
$T^1(R)_0=0$ over $\mathbb Q$ and then over $\mathbb C$.

The section bases, module equations, cochains and associators are
determined by the printed multiplication data. The accompanying exact
verification supplies the rational gauges and nonzero pivot minors.
The selected rows suffice because the genuine gauge subspace gives
the matching lower bound for the full kernel.

\paragraph{The projective deformation.}
The local theorem shows that every ordinary first-order deformation
of either canonical model is locally trivial at its singularities.
Its relative reflexive canonical powers are consequently well defined
and flat. For $m\ge2$, canonical vanishing gives
$H^1(Z,\omega_Z^{[m]})=0$, so global sections form flat modules
over the dual numbers and commute with restriction to the special
fibre. The degree-one space stays zero because $p_g=0$; the degree-zero
space is the scalar base ring because $H^1(\mathcal O_Z)=0$.
Multiplication therefore produces a flat graded deformation of the
entire canonical ring, not just a selected Veronese or subalgebra.

The canonical class is ample. A sufficiently divisible relative
canonical power is relatively ample, so relative Proj reconstructs
the original infinitesimal family from this ring. A trivial graded
ring deformation yields a trivial projective deformation. Thus
$T^1_Z$ injects into $T^1(R)_0$, which has just been proved zero.
The local-to-global sequence gives $H^1(Z,T_Z)=0$ as well.
Induction across small Artin extensions proves formal rigidity and
establishes Theorem~\ref{thm:exceptional-ordinary-rigidity46}.

\paragraph{Canonical models and terminalizations.}
The rigidity theorem concerns the canonical models Vega and Altair.
Vega is nonterminal, so the theorem makes no assertion about the
ordinary deformation space of a crepant terminalization. Altair is
terminal by Subsection~\ref{sec:altair-canonical-geometry}.



\section{Pluricanonical maps and completion of the classification}
\label{sec:stability-completion}
We determine the stable pluricanonical birationality index of both
exceptional models and then assemble the proof of
Theorem~\ref{thm:classification}. The calculations use the complete
canonical rings constructed above. Recognition transfers the conclusions
to every realization of $\mathcal B_4$ or $\mathcal B_5$.

\subsection{The distinct effective
indices}\label{ag:the-distinct-effective-indices}

Write $Z_A=\mathrm{Altair}$ and $Z_V=\mathrm{Vega}$. We use the notation \[
\begin{aligned}
e(X)&=\min\{m\geq1:P_m(X)>0\},\\
\delta(X)&=\min\{m\geq1:P_m(X)\geq2\},\\
\nu(X)&=\min\{a:P_m(X)>0\text{ for all }m\geq a\},\\
r_s(X)&=\min\{r:\varphi_m\text{ is birational for all }m\geq r\}.
\end{aligned}
\] The last is the stable pluricanonical birationality index; it is
unrelated to stable birational equivalence.

\begin{proposition}\label{ag:prop:indices}
The indices are \begin{equation}
\begin{array}{c|ccccc}
&e&\delta&\nu&\min\{a:P_m\geq2\ \forall m\geq a\}&r_s\\ \hline
\text{Altair}&5&14&12&18&23\\
\text{Vega}&4&12&14&18&23.
\end{array}
\label{ag:eq:3.1}
\end{equation}
\end{proposition}
\begin{proof}
For the numerical columns, apply Reid--Riemann--Roch \cite{Reid} to
the baskets in Theorem~\ref{ms:thm:three-baskets} with \(\chi=2\). The relevant values are \begin{equation}
\begin{array}{c|rrrrrrrrrrrr}
m&11&12&13&14&15&16&17&18&19&20&21&22\\ \hline
P_m(Z_A)&0&1&1&2&2&2&1&3&3&4&4&4\\
P_m(Z_V)&0&2&0&2&2&2&1&3&3&5&3&4 .
\end{array}
\label{ag:eq:3.2}
\end{equation} Multiplication by \(a_e\neq0\) propagates a complete block of \(e\)
consecutive degrees. The zeros in degrees \(11\) and \(13\), and
\(P_{17}=1\), prove the stated minimality. The birationality column follows from Theorems~\ref{thm:altair-stability}
and~\ref{ag:thm:vega-stability}, proved below.
\end{proof}

\subsection{The sharp birational stability of Altair}
\begin{theorem}\label{thm:altair-stability}
For Altair, the pluricanonical map $\varphi_m$ is birational if and only if
$m\ge23$. Thus its first birational degree and its birational stability
index are both $23$.
\end{theorem}

We use a simple linear-algebra observation to recover the function field.
\begin{lemma}[Recovery from linear relations]\label{lem:field-recovery}
Let $X$ be an integral projective variety. Suppose that nonzero sections
$h_0,\ldots,h_n$ of the same degree have ratios generating $\mathbb C(X)$.
Let $s_0,\ldots,s_r$ be sections of another common degree, with $s_0\ne0$. Assume that
\[
 \sum_{i=0}^n A_{ji}(s_0,\ldots,s_r)h_i=0\qquad(1\le j\le n),
\]
where the $A_{ji}$ are homogeneous of one common degree, and that the
coefficient matrix has rank $n$ over $\mathbb C(X)$. Then the map defined
by the $s_i$ is birational.
\end{lemma}
\begin{proof}
Let $F$ be the subfield generated by the ratios $s_i/s_0$. Dividing the
coefficient matrix by a nonzero monomial of its common degree puts its
entries in $F$. Its kernel is one-dimensional. Cofactors, or Gaussian
elimination over $F$, therefore recover all ratios $h_i/h_0$. These
ratios generate $\mathbb C(X)$, so $F=\mathbb C(X)$.
\end{proof}

\begin{proof}[Proof of Theorem~\ref{thm:altair-stability}]
We use the generators of $R_*$ introduced in its construction, with
$a=a_5$. On the chart $Y=(a=1)$ of the affine cone the ring is generated
by
\[
 b,u,z,e,w,p,v,y,
 \qquad\text{of residual }\mu_5\text{-weights }(3,2,4,0,4,1,3,0).
\]
The exact eliminations proving this assertion are recorded in
Appendix~\ref{app:stability-data}. This chart is smooth and finite free
of rank sixteen over $\mathbb C[b,u,z]$, by the construction and local
analysis of Altair. More explicitly, at $a=1$ the four $a$-monic module
relations eliminate $r_{21},K_{23},r_{27},r_{34}$; the remaining relation
is dependent by the second module syzygy. It is irreducible: the residual group permutes its
components transitively, since its quotient is the integral open set
$D_+(a)\subset Z_5$. Their number divides five, and their equal generic
ranks over the parameter ring show that it also divides sixteen.
Hence there is only one component.

On the dense set $b\ne0$, put
\[
 B=b^5,\quad T=bu,\quad Z=b^2z,\quad E=e,\quad
 W=b^2w,\quad P=p/b^2,\quad V=v/b,\quad Y_0=y.
\]
Writing $F=\mathbb C(B,T,Z,E,W,P,V,Y_0)$, we have
$\mathbb C(Y)=F(b)$ and $b^5=B$. The five residual scalings act
faithfully on $b$. Consequently
\[
 \mathbb C(Z_5)=F.
\]
The ordered degree-forty sections
\begin{equation}\label{eq:stability-target}
 \mathcal H=(a^8,b^5,a^4bu,a^2b^2z,a^5e,a^2b^2w,b^3p,bzv,a^4y)
\end{equation}
have ratios $(1,B,T,Z,E,W,BP,ZV,Y_0)$ on $a=1$. Since $B$ and $Z$ are
nonzero rational functions, these ratios generate the function field.

For $m=23,\ldots,27$, the complete module presentation gives the
following ordered bases:
\begin{center}
\begin{tabular}{@{}cl@{}}
\toprule
$m$ & A basis $\mathcal S_m$ of $R_m$\\
\midrule
$23$ & $a^3b,aq,be,av,K_{23}$\\
$24$ & $u^2,b^3,a^2z,af,a^2w,bp,r_{24}$\\
$25$ & $abu,a^5,a^2e,l,ay,r_{25}$\\
$26$ & $uz,a^2b^2,bq,aj,uw,a^2p,bv$\\
$27$ & $abz,a^3u,ak,ue,bf,abw,r_{27},M_{27}$\\
\bottomrule
\end{tabular}
\end{center}
The degree-$26$ relation eliminates $ar_{21}$. For each row we obtain
eight homogeneous linear relations among the nine members of
$\mathcal H$, with coefficients polynomial in $\mathcal S_m$.
Appendix~\ref{app:stability-data} specifies these rational identities
and a direct check using the defining multiplication table. The
coefficient degrees and the corresponding nonzero minors are
\begin{center}
\begin{tabular}{@{}crrr@{}}
\toprule
$m$ & $P_m$ & Coefficient degree & Minor modulo $1009$\\
\midrule
$23$&5&4&368\\
$24$&7&2&307\\
$25$&6&3&854\\
$26$&7&2&526\\
$27$&8&2&717\\
\bottomrule
\end{tabular}
\end{center}
The identities hold over $\mathbb Q$, and the indicated evaluation is
an algebra homomorphism to $\mathbb F_{1009}$ at which all denominators
are invertible. Thus each minor is nonzero in characteristic zero.
Lemma~\ref{lem:field-recovery} proves birationality in these five degrees.
Multiplication by $a\in R_5$ preserves the ratios of a subsystem, so
birationality in degree $m$ implies birationality in degree $m+5$.
It follows for every $m\ge23$.

Conversely, the Hilbert function gives
\[
 (P_{14},\ldots,P_{22})=(2,2,2,1,3,3,4,4,4),
 \qquad P_m\le1\quad(1\le m<14).
\]
For $P_m\le3$ the image has dimension less than three. If $P_m=4$ and
$\varphi_m$ were birational, its image would be all of $\mathbb P^3$,
contrary to general type. This proves the sharp lower bound.
\end{proof}

\subsection{Vega's stable index}\label{ag:sec:vega-stability}

\begin{theorem}\label{ag:thm:vega-stability}
The maps \(\varphi_{m,V}\) are birational for
every \(m\geq23\), and \(r_s(Z_V)=23\).
\end{theorem}
\begin{proof}
We use the complete canonical algebra of Subsection~\ref{sec:quartic-module46}. Write
its parameters as \[
a_4,\ u_{10},\ b_{12},\ c_{14}.
\] The invariant basis has degrees \(0,15,18,20,21,23,26,41\), with
elements \(1,h,q,r,s,k,\ell,t\). On \(a=1\), the full chart algebra is
free over \(\mathbf C[u,b,c]\) with ordered basis \begin{equation}
(1,h,q,r,s,k,\ell,t,v,w,z,d,p,g,i,j).
\label{ag:eq:3.3}
\end{equation}
Indeed, on $a=1$ the module relations of degrees $27,29,31,40$
eliminate $e_{23},f_{25},n_{27},k_{36}$, since their coefficients
are units. The relation of degree $37$ follows from the second
syzygy $(-c,b,-u,a,0)$ in (QM.4). This leaves precisely the eight
odd generators displayed above, together with the invariant basis.
Its function field is denoted by \(L\). This chart is integral: the
reduced prime divisor \((a)\) makes \(T^4-a\) Eisenstein at its
height-one valuation, and homogeneous rescaling identifies the
corresponding localized extension with a Laurent extension of the chart.

The multiplication table in Appendix~\ref{app:quartic-exact-coefficients} gives \[
v^2=u^3-\frac{108}{25}ub+\frac{432}{125}q
-\frac{972}{625}c+\frac{9983169}{100000}u.
\] Thus \(u,b,c,v\) recover \(q\). The equations for \(h^2,q^2\) recover
\(r,\ell\), with coefficient determinant \[
-12u^2+\frac{972}{25}b-\frac{29741013}{100000}.
\] Those for \(hq,hr\) recover \(s,k\), with determinant \[
-\left(b-\frac{122391}{16000}\right)^2+4uc-\frac{66}{5}u^2.
\] Both are nonzero since the parameters are algebraically independent.
The coefficient of \(t\) in \(h\ell\) is one, and the products with
\(v\) recover the remaining odd generators. Hence \[
L=\mathbf C(u,b,c,h,v).
\]

Consider the six degree-\(40\) sections \begin{equation}
\mathcal H=(a^{10},a^7b,a^5u^2,a^4uc,uh^2,uhv).
\label{ag:eq:3.4}
\end{equation} Their ratios generate \[
K_0=\mathbf C(b,u^2,uc,uh^2,uhv).
\] Adjoining \(h\) recovers \(L\), and \(h^4=(uh^2)^2/u^2\in K_0\). The
residual \(\mu_4\) acts faithfully on \(h\), so \([L:K_0]=4\) and
\(K_0=\mathbf C(Z_V)\).

It remains to recover these ratios from four consecutive systems. The
exact calculation described in Appendix~\ref{app:vega-stability-data} gives five homogeneous
relations \begin{equation}
\sum_{j=0}^{5} Q_{ij}(\mathcal S_m)\mathcal H_j=0,
\qquad 1\leq i\leq5,
\label{ag:eq:3.5}
\end{equation} where \(\mathcal S_m\) is a full basis of \(R_m\) and every
\(Q_{ij}\) has the same degree \(k\) in that basis. The data are:

\begin{center}
\begin{tabular}{@{}rrrrr@{}}
\toprule
$m$ & $P_m$ & $k$ & Matrix size & Minor modulo $1009$\\
\midrule
$23$ & $5$ & $6$ & $2072\times1260$ & $559$\\
$24$ & $7$ & $3$ & $416\times504$ & $116$\\
$25$ & $6$ & $3$ & $490\times336$ & $354$\\
$26$ & $7$ & $3$ & $562\times504$ & $434$\\
\bottomrule
\end{tabular}
\end{center}

The relations are identities over \(\mathbf Q\), verified in the
defining algebra. The modular calculation establishes only that the
selected minors do not vanish identically. Lemma~\ref{lem:field-recovery} therefore puts
every \(\mathcal H_j/\mathcal H_0\) in the image field of \(\varphi_m\).
All four maps are birational. Multiplication by \(a_4\) propagates the
conclusion to every \(m\geq23\).

Finally \(P_{22}=4\). A birational map from a threefold into
\(\mathbf P^3\) would make it rational, contrary to general type. Thus
\(\varphi_{22}\) is not birational.
\end{proof}

Among degrees below \(23\), only \(m=20\) has \(P_m>4\) for Vega. Its
first individual birational degree is consequently either \(20\) or
\(23\). The theorem does not select between these two possibilities.

\subsection{Completion of Theorem~\ref{thm:classification}}
\label{sec:completion46}
\begin{proof}[Proof of Theorem~\ref{thm:classification}]
Proposition~\ref{prop:equality-census} gives the complete equality
census and the reduction $p_g=q=0$, $\chi\in\{1,2,3\}$.
The Euler-one branch is exactly the hypersurface family by
Theorem~\ref{thm:chi1-classification}, and its moduli and ordinary
deformations are described by Theorem~\ref{thm:chi1-moduli}.
Section~\ref{sec:chi3} excludes all $283$ Euler-three baskets.

For $\chi=2$, the common tests and their applications reduce the $65$
baskets to the eleven in Table~\ref{tab:residual11}. Part~III excludes
nine of them. Part~IV constructs Vega and Altair for the two remaining
baskets, verifies that their entire normal section rings are canonical
rings, and recognizes every realization of those baskets. The ordinary
Artinian rigidity proved in the preceding section completes the
assertion about the two isolated models.
\end{proof}

The difference between the family and the two isolated models is now
visible in their canonical rings. In the Euler-one case, four
independent sections form a polynomial subring, and a single quadratic
relation retains the freedom of its coefficients. For Vega and Altair,
the nonzero two-form determines a marked surface and imposes identities
on its normal extension. The marked surface, the actual ambient
symmetry, and the associative extension of its full algebra are all
needed to recover the threefold. The deformation calculation then
establishes rigidity for the resulting canonical model.


\clearpage
\appendix

\section{The complete census and its finite witnesses}\label{app:census}
To print the complete list compactly, write $\beta_i$ for the formal
basket constituent in the following dictionary. Addition means multiset
union, and $c\beta_i$ means multiplicity $c$. Thus every row below is
an explicit full basket; the notation imposes no identification of
actual singular points. The index is recovered as $I=\operatorname{lcm}(r)$.
The companion table \texttt{data/equality\_baskets\_complete.csv} also
lists every basket in pair notation together with $I$.

\begin{center}\small
\begin{tabular}{@{}llllll@{}}
$\beta_{1}=(1,2)$ & $\beta_{2}=(1,3)$ & $\beta_{3}=(1,4)$ & $\beta_{4}=(1,5)$ & $\beta_{5}=(2,5)$ & $\beta_{6}=(1,6)$\\[2pt]
$\beta_{7}=(2,7)$ & $\beta_{8}=(3,7)$ & $\beta_{9}=(3,8)$ & $\beta_{10}=(2,9)$ & $\beta_{11}=(4,9)$ & $\beta_{12}=(3,10)$\\[2pt]
$\beta_{13}=(3,11)$ & $\beta_{14}=(4,11)$ & $\beta_{15}=(5,11)$ & $\beta_{16}=(5,12)$ & $\beta_{17}=(4,13)$ & $\beta_{18}=(6,13)$\\[2pt]
$\beta_{19}=(3,14)$ & $\beta_{20}=(5,14)$ & $\beta_{21}=(4,15)$ & $\beta_{22}=(7,15)$ & $\beta_{23}=(5,16)$ & $\beta_{24}=(7,16)$\\[2pt]
$\beta_{25}=(6,17)$ & $\beta_{26}=(7,17)$ & $\beta_{27}=(5,18)$ & $\beta_{28}=(7,18)$ & $\beta_{29}=(7,20)$ & $\beta_{30}=(9,20)$\\[2pt]
$\beta_{31}=(5,21)$ & $\beta_{32}=(8,21)$ & $\beta_{33}=(10,21)$ & $\beta_{34}=(7,22)$ & $\beta_{35}=(9,22)$ & $\beta_{36}=(7,24)$\\[2pt]
$\beta_{37}=(11,24)$ & $\beta_{38}=(11,28)$ & $\beta_{39}=(13,28)$ & $\beta_{40}=(11,30)$ & $\beta_{41}=(13,30)$\\[2pt]
\end{tabular}\end{center}
For Euler characteristic three, A denotes the empty two-visible-prime
list, B the coplanar cyclic-cover obstruction, C the commuting
cubic--quintic obstruction, and D the effective-cubic obstruction.
These are the four arguments of Section~\ref{sec:chi3}.

\subsection{Euler characteristics one and two}
For $\chi=1$ the sole row is
\[
 \basket{C1-00030}:\quad
 3\beta_1+\beta_3+\beta_6+\beta_5+\beta_8,
 \qquad I=420.
\]
It gives the family $X_{46}$. The following $65$ rows have $\chi=2$;
we suppress their common prefix \texttt{C2-}. The labels V and A in the
last column denote Vega and Altair, respectively; all other rows are
excluded in the text.

\par\noindent
\begin{minipage}[t]{.49\linewidth}\vspace{0pt}\small\setlength{\tabcolsep}{2pt}
\begin{tabular}{@{}p{9mm}p{62mm}p{3mm}@{}}\toprule Row & Basket & \\\midrule
00001 & $9\beta_{1}+3\beta_{3}+\beta_{6}+\beta_{5}+\beta_{20}$ & \\
00002 & $8\beta_{1}+\beta_{3}+\beta_{4}+\beta_{6}+2\beta_{5}+\beta_{20}$ & \\
00003 & $6\beta_{1}+\beta_{6}+\beta_{5}+2\beta_{9}+\beta_{20}$ & \\
00004 & $5\beta_{1}+2\beta_{2}+5\beta_{3}+\beta_{5}+\beta_{20}$ & \\
00005 & $5\beta_{1}+\beta_{3}+\beta_{5}+\beta_{10}+\beta_{8}+\beta_{28}$ & \\
00006 & $4\beta_{1}+5\beta_{2}+5\beta_{3}+\beta_{5}+\beta_{8}$ & \\
00007 & $4\beta_{1}+5\beta_{2}+3\beta_{3}+\beta_{5}+\beta_{20}$ & \\
00008 & $4\beta_{1}+2\beta_{2}+3\beta_{3}+\beta_{4}+2\beta_{5}+\beta_{20}$ & \\
00009 & $4\beta_{1}+3\beta_{3}+\beta_{4}+\beta_{20}+\beta_{22}$ & \\
00010 & $3\beta_{1}+5\beta_{2}+3\beta_{3}+\beta_{4}+2\beta_{5}+\beta_{8}$ & \\
00011 & $3\beta_{1}+5\beta_{2}+\beta_{3}+\beta_{4}+2\beta_{5}+\beta_{20}$ & \\
00012 & $3\beta_{1}+3\beta_{2}+3\beta_{3}+\beta_{4}+\beta_{8}+\beta_{22}$ & \\
00013 & $3\beta_{1}+3\beta_{2}+\beta_{3}+\beta_{4}+\beta_{20}+\beta_{22}$ & \\
00014 & $3\beta_{1}+\beta_{3}+2\beta_{4}+\beta_{5}+\beta_{20}+\beta_{22}$ & \\
00015 & $2\beta_{1}+4\beta_{2}+\beta_{3}+\beta_{5}+\beta_{7}+\beta_{32}$ & \\
00016 & $2\beta_{1}+3\beta_{2}+\beta_{3}+2\beta_{4}+\beta_{5}+\beta_{8}+\beta_{22}$ & \\
00017 & $2\beta_{1}+2\beta_{2}+2\beta_{3}+\beta_{5}+2\beta_{9}+\beta_{20}$ & \\
00018 & $2\beta_{1}+2\beta_{2}+2\beta_{3}+\beta_{20}+\beta_{30}$ & \\
00019 & $2\beta_{1}+2\beta_{2}+\beta_{3}+3\beta_{5}+\beta_{12}+\beta_{20}$ & \\
00020 & $2\beta_{1}+\beta_{3}+\beta_{5}+\beta_{12}+\beta_{20}+\beta_{22}$ & \\
00021 & $\beta_{1}+5\beta_{2}+2\beta_{3}+\beta_{5}+\beta_{8}+2\beta_{9}$ & \\
00022 & $\beta_{1}+4\beta_{2}+2\beta_{3}+\beta_{5}+\beta_{16}+\beta_{20}$ & \\
00023 & $\beta_{1}+2\beta_{2}+2\beta_{3}+2\beta_{5}+\beta_{8}+\beta_{29}$ & \\
00024 & $\beta_{1}+2\beta_{2}+\beta_{4}+\beta_{5}+\beta_{20}+\beta_{30}$ & \\
00025 & $\beta_{1}+2\beta_{3}+\beta_{8}+\beta_{22}+\beta_{29}$ & \\
00026 & $\beta_{1}+\beta_{4}+2\beta_{9}+\beta_{20}+\beta_{22}$ & \\
00027 & $3\beta_{2}+\beta_{4}+\beta_{8}+2\beta_{9}+\beta_{22}$ & \\
00028 & $2\beta_{2}+\beta_{4}+\beta_{16}+\beta_{20}+\beta_{22}$ & \\
00029 & $\beta_{4}+\beta_{5}+\beta_{8}+\beta_{22}+\beta_{29}$ & \\
01822 & $6\beta_{1}+\beta_{2}+2\beta_{4}+\beta_{8}+\beta_{12}+\beta_{16}$ & \\
01823 & $2\beta_{1}+8\beta_{2}+\beta_{3}+\beta_{4}+2\beta_{5}+\beta_{8}$ & \\
01824 & $2\beta_{1}+6\beta_{2}+\beta_{3}+\beta_{4}+\beta_{8}+\beta_{22}$ & \\
01825 & $2\beta_{1}+2\beta_{2}+\beta_{3}+\beta_{5}+\beta_{8}+\beta_{17}+\beta_{18}$ & \\
\bottomrule\end{tabular}\end{minipage}\hfill%
\begin{minipage}[t]{.49\linewidth}\vspace{0pt}\small\setlength{\tabcolsep}{2pt}
\begin{tabular}{@{}p{9mm}p{62mm}p{3mm}@{}}\toprule Row & Basket & \\\midrule
01826 & $\beta_{1}+5\beta_{2}+2\beta_{3}+\beta_{8}+\beta_{30}$ & \\
01827 & $\beta_{1}+5\beta_{2}+\beta_{3}+3\beta_{5}+\beta_{8}+\beta_{12}$ & \\
01828 & $\beta_{1}+3\beta_{2}+\beta_{3}+\beta_{5}+\beta_{8}+\beta_{12}+\beta_{22}$ & \\
01829 & $7\beta_{2}+2\beta_{3}+\beta_{5}+\beta_{8}+\beta_{16}$ & \\
01830 & $5\beta_{2}+\beta_{4}+\beta_{5}+\beta_{8}+\beta_{30}$ & \\
02067 & $9\beta_{1}+2\beta_{2}+\beta_{3}+2\beta_{4}+\beta_{8}+\beta_{12}$ & \\
02068 & $7\beta_{1}+3\beta_{2}+\beta_{3}+\beta_{4}+\beta_{6}+2\beta_{5}+\beta_{8}$ & \\
02069 & $6\beta_{1}+2\beta_{2}+\beta_{3}+2\beta_{4}+3\beta_{5}+\beta_{8}$ & \\
02070 & $6\beta_{1}+\beta_{3}+\beta_{6}+3\beta_{5}+\beta_{8}+\beta_{12}$ & \\
02071 & $6\beta_{1}+\beta_{3}+3\beta_{5}+\beta_{8}+\beta_{21}$ & \\
02072 & $5\beta_{1}+2\beta_{2}+2\beta_{3}+\beta_{6}+\beta_{5}+\beta_{8}+\beta_{16}$ & \\
02073 & $5\beta_{1}+2\beta_{2}+\beta_{6}+\beta_{5}+\beta_{16}+\beta_{20}$ & \\
02074 & $5\beta_{1}+\beta_{2}+\beta_{8}+2\beta_{12}+\beta_{16}$ & \\
02075 & $5\beta_{1}+\beta_{6}+2\beta_{5}+\beta_{8}+\beta_{29}$ & \\
02076 & $4\beta_{1}+\beta_{2}+2\beta_{3}+\beta_{4}+2\beta_{5}+\beta_{8}+\beta_{16}$ & \\
02077 & $4\beta_{1}+\beta_{2}+\beta_{4}+2\beta_{5}+\beta_{16}+\beta_{20}$ & \\
02078 & $3\beta_{1}+\beta_{5}+\beta_{7}+\beta_{16}+\beta_{32}$ & \\
02079 & $2\beta_{1}+\beta_{2}+3\beta_{5}+\beta_{8}+\beta_{12}+\beta_{16}$ & \\
02080 & $2\beta_{1}+\beta_{3}+\beta_{5}+2\beta_{16}+\beta_{20}$ & \\
02346 & $8\beta_{1}+3\beta_{2}+3\beta_{3}+\beta_{6}+\beta_{5}+\beta_{8}$ & \\
02347 & $7\beta_{1}+2\beta_{2}+3\beta_{3}+\beta_{4}+2\beta_{5}+\beta_{8}$ & V\\
02348 & $5\beta_{1}+3\beta_{2}+\beta_{6}+\beta_{5}+\beta_{8}+2\beta_{9}$ & \\
02349 & $5\beta_{1}+\beta_{2}+4\beta_{3}+\beta_{5}+\beta_{8}+\beta_{16}$ & \\
02350 & $5\beta_{1}+\beta_{2}+2\beta_{3}+\beta_{5}+\beta_{16}+\beta_{20}$ & \\
02351 & $4\beta_{1}+5\beta_{2}+\beta_{6}+\beta_{5}+\beta_{8}+\beta_{16}$ & \\
02352 & $3\beta_{1}+8\beta_{2}+3\beta_{3}+\beta_{5}+\beta_{8}$ & \\
02353 & $3\beta_{1}+4\beta_{2}+\beta_{4}+2\beta_{5}+\beta_{8}+\beta_{16}$ & \\
02354 & $2\beta_{1}+\beta_{2}+\beta_{3}+\beta_{5}+2\beta_{8}+\beta_{32}$ & \\
02355 & $2\beta_{1}+\beta_{2}+\beta_{3}+\beta_{5}+\beta_{8}+2\beta_{9}+\beta_{16}$ & \\
02356 & $\beta_{1}+3\beta_{2}+\beta_{3}+\beta_{5}+\beta_{8}+2\beta_{16}$ & \\
02357 & $\beta_{1}+2\beta_{2}+\beta_{3}+\beta_{5}+\beta_{8}+\beta_{11}+\beta_{28}$ & \\
02397 & $4\beta_{1}+2\beta_{2}+\beta_{4}+2\beta_{5}+\beta_{8}+2\beta_{9}$ & A\\
\bottomrule\end{tabular}\end{minipage}\par

\clearpage
\subsection{Euler characteristic three}
The common prefix is \texttt{C3-}. There are $283$ rows.
\par\noindent
\begin{minipage}[t]{.49\linewidth}\vspace{0pt}\small\setlength{\tabcolsep}{2pt}
\begin{tabular}{@{}p{9mm}p{62mm}p{3mm}@{}}\toprule Row & Basket & \\\midrule
00001 & $13\beta_{1}+2\beta_{2}+\beta_{5}+2\beta_{7}+2\beta_{9}+\beta_{19}$ & A\\
00002 & $13\beta_{1}+\beta_{2}+\beta_{3}+2\beta_{4}+\beta_{7}+\beta_{12}+\beta_{32}$ & A\\
00003 & $12\beta_{1}+3\beta_{2}+2\beta_{4}+\beta_{8}+2\beta_{9}+\beta_{21}$ & A\\
00004 & $12\beta_{1}+2\beta_{2}+2\beta_{3}+2\beta_{4}+\beta_{8}+2\beta_{9}+\beta_{12}$ & A\\
00005 & $12\beta_{1}+2\beta_{2}+\beta_{3}+\beta_{5}+\beta_{7}+\beta_{8}+\beta_{19}+\beta_{20}$ & A\\
00006 & $12\beta_{1}+2\beta_{2}+2\beta_{4}+\beta_{21}+\beta_{16}+\beta_{20}$ & A\\
00007 & $12\beta_{1}+\beta_{2}+2\beta_{3}+2\beta_{4}+\beta_{12}+\beta_{16}+\beta_{20}$ & A\\
00008 & $12\beta_{1}+\beta_{2}+\beta_{3}+\beta_{6}+3\beta_{5}+\beta_{21}+\beta_{20}$ & A\\
00009 & $12\beta_{1}+3\beta_{3}+\beta_{6}+3\beta_{5}+\beta_{12}+\beta_{20}$ & A\\
00010 & $12\beta_{1}+2\beta_{4}+\beta_{5}+\beta_{8}+\beta_{21}+\beta_{29}$ & A\\
00011 & $11\beta_{1}+2\beta_{2}+4\beta_{3}+\beta_{6}+\beta_{5}+\beta_{16}+\beta_{20}$ & A\\
00012 & $11\beta_{1}+\beta_{2}+\beta_{3}+\beta_{6}+\beta_{5}+\beta_{10}+\beta_{20}+\beta_{28}$ & A\\
00013 & $11\beta_{1}+\beta_{2}+\beta_{5}+\beta_{10}+\beta_{8}+2\beta_{9}+\beta_{27}$ & A\\
00014 & $11\beta_{1}+4\beta_{3}+\beta_{6}+2\beta_{5}+\beta_{8}+\beta_{29}$ & A\\
00015 & $11\beta_{1}+\beta_{5}+\beta_{10}+\beta_{16}+\beta_{20}+\beta_{27}$ & A\\
00016 & $10\beta_{1}+4\beta_{2}+\beta_{3}+\beta_{6}+\beta_{5}+\beta_{10}+\beta_{8}+\beta_{28}$ & A\\
00017 & $10\beta_{1}+3\beta_{2}+3\beta_{3}+\beta_{5}+\beta_{10}+\beta_{8}+\beta_{28}$ & A\\
00018 & $10\beta_{1}+2\beta_{2}+\beta_{5}+\beta_{7}+\beta_{19}+\beta_{38}$ & A\\
00019 & $10\beta_{1}+\beta_{2}+4\beta_{3}+\beta_{4}+2\beta_{5}+\beta_{16}+\beta_{20}$ & A\\
00020 & $10\beta_{1}+\beta_{2}+\beta_{3}+3\beta_{4}+\beta_{8}+\beta_{16}+\beta_{29}$ & A\\
00021 & $10\beta_{1}+\beta_{2}+\beta_{3}+\beta_{4}+\beta_{8}+\beta_{12}+\beta_{40}$ & A\\
00022 & $10\beta_{1}+2\beta_{3}+\beta_{4}+\beta_{6}+3\beta_{5}+\beta_{8}+\beta_{29}$ & A\\
00023 & $10\beta_{1}+\beta_{3}+\beta_{4}+2\beta_{5}+\beta_{11}+\beta_{20}+\beta_{27}$ & A\\
00024 & $9\beta_{1}+3\beta_{2}+\beta_{3}+\beta_{4}+2\beta_{5}+\beta_{10}+\beta_{8}+\beta_{28}$ & A\\
00025 & $9\beta_{1}+2\beta_{2}+\beta_{5}+\beta_{7}+2\beta_{8}+2\beta_{9}+\beta_{19}$ & A\\
00026 & $9\beta_{1}+\beta_{2}+\beta_{3}+2\beta_{4}+2\beta_{8}+\beta_{12}+\beta_{32}$ & A\\
00027 & $9\beta_{1}+\beta_{2}+\beta_{3}+2\beta_{4}+\beta_{8}+2\beta_{9}+\beta_{12}+\beta_{16}$ & A\\
00028 & $9\beta_{1}+\beta_{3}+\beta_{6}+\beta_{5}+\beta_{13}+\beta_{20}+\beta_{35}$ & A\\
00029 & $9\beta_{1}+\beta_{6}+3\beta_{5}+\beta_{21}+\beta_{16}+\beta_{20}$ & A\\
00030 & $8\beta_{1}+3\beta_{2}+\beta_{5}+2\beta_{10}+\beta_{8}+2\beta_{9}+\beta_{11}$ & A\\
00031 & $8\beta_{1}+2\beta_{2}+3\beta_{3}+\beta_{5}+\beta_{10}+2\beta_{11}+\beta_{20}$ & A\\
00032 & $8\beta_{1}+2\beta_{2}+\beta_{3}+2\beta_{4}+\beta_{8}+\beta_{12}+\beta_{11}+\beta_{28}$ & A\\
00033 & $8\beta_{1}+2\beta_{2}+\beta_{3}+\beta_{6}+\beta_{5}+\beta_{8}+\beta_{20}+\beta_{32}$ & A\\
00034 & $8\beta_{1}+2\beta_{2}+\beta_{3}+\beta_{5}+3\beta_{8}+\beta_{19}+\beta_{20}$ & A\\
00035 & $8\beta_{1}+2\beta_{2}+\beta_{5}+2\beta_{10}+\beta_{11}+\beta_{16}+\beta_{20}$ & A\\
00036 & $8\beta_{1}+5\beta_{3}+\beta_{5}+2\beta_{16}+\beta_{20}$ & A\\
00037 & $8\beta_{1}+3\beta_{3}+\beta_{6}+\beta_{5}+\beta_{8}+\beta_{25}+\beta_{26}$ & A\\
00038 & $8\beta_{1}+3\beta_{3}+\beta_{6}+\beta_{5}+\beta_{11}+\beta_{20}+\beta_{28}$ & A\\
00039 & $8\beta_{1}+\beta_{3}+\beta_{6}+2\beta_{5}+\beta_{8}+2\beta_{9}+\beta_{29}$ & A\\
00040 & $8\beta_{1}+\beta_{3}+\beta_{6}+\beta_{5}+\beta_{20}+\beta_{25}+\beta_{26}$ & A\\
00041 & $8\beta_{1}+2\beta_{5}+2\beta_{10}+\beta_{8}+\beta_{11}+\beta_{29}$ & A\\
00042 & $8\beta_{1}+\beta_{5}+2\beta_{9}+\beta_{11}+\beta_{20}+\beta_{27}$ & A\\
00043 & $7\beta_{1}+5\beta_{2}+\beta_{3}+\beta_{5}+\beta_{10}+2\beta_{11}+\beta_{20}$ & A\\
00044 & $7\beta_{1}+3\beta_{2}+3\beta_{3}+\beta_{6}+\beta_{5}+\beta_{8}+\beta_{11}+\beta_{28}$ & A\\
00045 & $7\beta_{1}+3\beta_{2}+\beta_{3}+\beta_{6}+\beta_{5}+\beta_{11}+\beta_{20}+\beta_{28}$ & A\\
00046 & $7\beta_{1}+3\beta_{2}+\beta_{5}+\beta_{10}+\beta_{8}+2\beta_{9}+\beta_{28}$ & A\\
00047 & $7\beta_{1}+2\beta_{2}+6\beta_{3}+2\beta_{5}+\beta_{8}+\beta_{29}$ & D\\
00048 & $7\beta_{1}+2\beta_{2}+5\beta_{3}+\beta_{5}+\beta_{8}+\beta_{11}+\beta_{28}$ & A\\
\bottomrule\end{tabular}\end{minipage}\hfill%
\begin{minipage}[t]{.49\linewidth}\vspace{0pt}\small\setlength{\tabcolsep}{2pt}
\begin{tabular}{@{}p{9mm}p{62mm}p{3mm}@{}}\toprule Row & Basket & \\\midrule
00049 & $7\beta_{1}+2\beta_{2}+4\beta_{3}+\beta_{4}+\beta_{5}+\beta_{20}+\beta_{30}$ & A\\
00050 & $7\beta_{1}+2\beta_{2}+4\beta_{3}+2\beta_{5}+\beta_{20}+\beta_{29}$ & D\\
00051 & $7\beta_{1}+2\beta_{2}+3\beta_{3}+\beta_{4}+2\beta_{5}+\beta_{8}+\beta_{17}+\beta_{18}$ & A\\
00052 & $7\beta_{1}+2\beta_{2}+3\beta_{3}+\beta_{4}+\beta_{7}+\beta_{22}+\beta_{32}$ & A\\
00053 & $7\beta_{1}+2\beta_{2}+3\beta_{3}+\beta_{5}+\beta_{11}+\beta_{20}+\beta_{28}$ & A\\
00054 & $7\beta_{1}+2\beta_{2}+2\beta_{3}+\beta_{5}+\beta_{10}+\beta_{8}+\beta_{16}+\beta_{28}$ & A\\
00055 & $7\beta_{1}+2\beta_{2}+\beta_{3}+\beta_{4}+2\beta_{5}+\beta_{10}+2\beta_{11}+\beta_{20}$ & A\\
00056 & $7\beta_{1}+2\beta_{2}+\beta_{3}+\beta_{4}+2\beta_{5}+\beta_{17}+\beta_{20}+\beta_{18}$ & A\\
00057 & $7\beta_{1}+2\beta_{2}+\beta_{5}+\beta_{10}+\beta_{16}+\beta_{20}+\beta_{28}$ & A\\
00058 & $7\beta_{1}+\beta_{2}+\beta_{3}+\beta_{4}+2\beta_{5}+\beta_{8}+\beta_{20}+\beta_{32}$ & A\\
00059 & $7\beta_{1}+2\beta_{5}+\beta_{10}+\beta_{8}+\beta_{28}+\beta_{29}$ & A\\
00060 & $6\beta_{1}+5\beta_{2}+4\beta_{3}+2\beta_{5}+\beta_{8}+\beta_{29}$ & D\\
00061 & $6\beta_{1}+5\beta_{2}+2\beta_{3}+\beta_{4}+\beta_{5}+\beta_{20}+\beta_{30}$ & A\\
00062 & $6\beta_{1}+5\beta_{2}+\beta_{3}+\beta_{4}+\beta_{7}+\beta_{22}+\beta_{32}$ & A\\
00063 & $6\beta_{1}+3\beta_{2}+4\beta_{3}+\beta_{4}+\beta_{8}+2\beta_{9}+\beta_{22}$ & A\\
00064 & $6\beta_{1}+3\beta_{2}+3\beta_{3}+\beta_{4}+2\beta_{5}+\beta_{8}+\beta_{12}+\beta_{22}$ & A\\
00065 & $6\beta_{1}+3\beta_{2}+\beta_{3}+\beta_{4}+2\beta_{5}+\beta_{12}+\beta_{20}+\beta_{22}$ & A\\
00066 & $6\beta_{1}+2\beta_{2}+4\beta_{3}+\beta_{4}+3\beta_{5}+\beta_{8}+\beta_{29}$ & C\\
00067 & $6\beta_{1}+2\beta_{2}+4\beta_{3}+\beta_{4}+\beta_{16}+\beta_{20}+\beta_{22}$ & A\\
00068 & $6\beta_{1}+2\beta_{2}+3\beta_{3}+\beta_{4}+2\beta_{5}+\beta_{8}+\beta_{11}+\beta_{28}$ & A\\
00069 & $6\beta_{1}+2\beta_{2}+2\beta_{3}+2\beta_{4}+2\beta_{5}+\beta_{20}+\beta_{30}$ & A\\
00070 & $6\beta_{1}+2\beta_{2}+2\beta_{3}+\beta_{4}+3\beta_{5}+\beta_{20}+\beta_{29}$ & C\\
00071 & $6\beta_{1}+2\beta_{2}+\beta_{3}+2\beta_{4}+3\beta_{5}+\beta_{8}+\beta_{17}+\beta_{18}$ & A\\
00072 & $6\beta_{1}+2\beta_{2}+\beta_{3}+2\beta_{4}+\beta_{5}+\beta_{7}+\beta_{22}+\beta_{32}$ & A\\
00073 & $6\beta_{1}+2\beta_{2}+\beta_{3}+\beta_{4}+2\beta_{5}+\beta_{11}+\beta_{20}+\beta_{28}$ & A\\
00074 & $6\beta_{1}+2\beta_{2}+\beta_{3}+\beta_{5}+\beta_{8}+\beta_{13}+2\beta_{14}+\beta_{15}$ & A\\
00075 & $6\beta_{1}+\beta_{2}+\beta_{3}+\beta_{5}+\beta_{7}+\beta_{17}+\beta_{18}+\beta_{32}$ & A\\
00076 & $6\beta_{1}+4\beta_{3}+\beta_{4}+\beta_{5}+\beta_{8}+\beta_{22}+\beta_{29}$ & A\\
00077 & $6\beta_{1}+2\beta_{3}+\beta_{4}+\beta_{5}+\beta_{20}+\beta_{22}+\beta_{29}$ & A\\
00078 & $5\beta_{1}+6\beta_{2}+2\beta_{3}+\beta_{4}+\beta_{8}+2\beta_{9}+\beta_{22}$ & A\\
00079 & $5\beta_{1}+5\beta_{2}+2\beta_{3}+\beta_{4}+3\beta_{5}+\beta_{8}+\beta_{29}$ & C\\
00080 & $5\beta_{1}+5\beta_{2}+2\beta_{3}+\beta_{4}+\beta_{16}+\beta_{20}+\beta_{22}$ & A\\
00081 & $5\beta_{1}+5\beta_{2}+2\beta_{4}+2\beta_{5}+\beta_{20}+\beta_{30}$ & A\\
00082 & $5\beta_{1}+4\beta_{2}+2\beta_{3}+\beta_{7}+\beta_{32}+\beta_{30}$ & A\\
00083 & $5\beta_{1}+4\beta_{2}+\beta_{3}+\beta_{4}+\beta_{10}+\beta_{8}+\beta_{22}+\beta_{28}$ & A\\
00084 & $5\beta_{1}+3\beta_{2}+2\beta_{3}+2\beta_{4}+\beta_{5}+\beta_{8}+2\beta_{9}+\beta_{22}$ & A\\
00085 & $5\beta_{1}+3\beta_{2}+2\beta_{3}+\beta_{4}+\beta_{5}+\beta_{8}+\beta_{22}+\beta_{29}$ & A\\
00086 & $5\beta_{1}+3\beta_{2}+\beta_{3}+2\beta_{4}+3\beta_{5}+\beta_{8}+\beta_{12}+\beta_{22}$ & A\\
00087 & $5\beta_{1}+3\beta_{2}+2\beta_{5}+\beta_{21}+\beta_{20}+\beta_{30}$ & A\\
00088 & $5\beta_{1}+2\beta_{2}+3\beta_{3}+2\beta_{9}+\beta_{20}+\beta_{30}$ & A\\
00089 & $5\beta_{1}+2\beta_{2}+2\beta_{3}+2\beta_{4}+\beta_{5}+\beta_{16}+\beta_{20}+\beta_{22}$ & A\\
00090 & $5\beta_{1}+2\beta_{2}+2\beta_{3}+2\beta_{5}+\beta_{12}+\beta_{20}+\beta_{30}$ & A\\
00091 & $5\beta_{1}+2\beta_{2}+2\beta_{3}+\beta_{5}+\beta_{8}+2\beta_{9}+\beta_{17}+\beta_{18}$ & A\\
00092 & $5\beta_{1}+2\beta_{2}+\beta_{3}+\beta_{5}+\beta_{7}+\beta_{12}+\beta_{22}+\beta_{32}$ & A\\
00093 & $5\beta_{1}+2\beta_{2}+\beta_{5}+\beta_{10}+2\beta_{9}+2\beta_{11}+\beta_{20}$ & A\\
00094 & $5\beta_{1}+\beta_{2}+3\beta_{3}+\beta_{6}+4\beta_{5}+\beta_{20}+\beta_{22}$ & A\\
00095 & $5\beta_{1}+\beta_{2}+2\beta_{3}+2\beta_{5}+\beta_{7}+\beta_{29}+\beta_{32}$ & C\\
00096 & $5\beta_{1}+\beta_{2}+2\beta_{3}+\beta_{5}+\beta_{17}+\beta_{16}+\beta_{20}+\beta_{18}$ & A\\
\bottomrule\end{tabular}\end{minipage}\par

\clearpage
\noindent\textbf{Euler characteristic three (continued).}
\smallskip
\par\noindent
\begin{minipage}[t]{.49\linewidth}\vspace{0pt}\small\setlength{\tabcolsep}{2pt}
\begin{tabular}{@{}p{9mm}p{62mm}p{3mm}@{}}\toprule Row & Basket & \\\midrule
00097 & $5\beta_{1}+\beta_{2}+\beta_{3}+\beta_{5}+\beta_{7}+\beta_{11}+\beta_{28}+\beta_{32}$ & A\\
00098 & $5\beta_{1}+2\beta_{3}+2\beta_{4}+2\beta_{5}+\beta_{8}+\beta_{22}+\beta_{29}$ & A\\
00099 & $5\beta_{1}+2\beta_{3}+\beta_{5}+\beta_{8}+\beta_{16}+\beta_{20}+\beta_{32}$ & A\\
00100 & $5\beta_{1}+2\beta_{4}+2\beta_{5}+\beta_{20}+\beta_{22}+\beta_{29}$ & A\\
00101 & $5\beta_{1}+4\beta_{5}+\beta_{21}+\beta_{20}+\beta_{29}$ & C\\
00102 & $4\beta_{1}+6\beta_{2}+2\beta_{3}+\beta_{5}+\beta_{7}+\beta_{16}+\beta_{32}$ & C\\
00103 & $4\beta_{1}+6\beta_{2}+2\beta_{4}+\beta_{5}+\beta_{8}+2\beta_{9}+\beta_{22}$ & A\\
00104 & $4\beta_{1}+5\beta_{2}+3\beta_{3}+\beta_{8}+2\beta_{9}+\beta_{30}$ & A\\
00105 & $4\beta_{1}+5\beta_{2}+2\beta_{4}+\beta_{5}+\beta_{16}+\beta_{20}+\beta_{22}$ & A\\
00106 & $4\beta_{1}+4\beta_{2}+3\beta_{3}+\beta_{5}+\beta_{8}+\beta_{20}+\beta_{32}$ & C\\
00107 & $4\beta_{1}+4\beta_{2}+3\beta_{3}+\beta_{5}+3\beta_{11}+\beta_{20}$ & A\\
00108 & $4\beta_{1}+4\beta_{2}+3\beta_{3}+\beta_{16}+\beta_{20}+\beta_{30}$ & A\\
00109 & $4\beta_{1}+4\beta_{2}+\beta_{3}+\beta_{6}+4\beta_{5}+\beta_{20}+\beta_{22}$ & A\\
00110 & $4\beta_{1}+4\beta_{2}+\beta_{4}+\beta_{5}+\beta_{7}+\beta_{32}+\beta_{30}$ & A\\
00111 & $4\beta_{1}+4\beta_{2}+\beta_{5}+\beta_{10}+2\beta_{11}+\beta_{16}+\beta_{20}$ & A\\
00112 & $4\beta_{1}+4\beta_{2}+\beta_{5}+\beta_{8}+2\beta_{9}+\beta_{21}+\beta_{22}$ & A\\
00113 & $4\beta_{1}+3\beta_{2}+2\beta_{3}+\beta_{5}+\beta_{8}+2\beta_{9}+\beta_{12}+\beta_{22}$ & A\\
00114 & $4\beta_{1}+3\beta_{2}+2\beta_{4}+2\beta_{5}+\beta_{8}+\beta_{22}+\beta_{29}$ & A\\
00115 & $4\beta_{1}+3\beta_{2}+4\beta_{5}+\beta_{8}+\beta_{21}+\beta_{29}$ & C\\
00116 & $4\beta_{1}+3\beta_{2}+\beta_{5}+2\beta_{9}+\beta_{31}+\beta_{33}$ & A\\
00117 & $4\beta_{1}+3\beta_{2}+\beta_{5}+\beta_{21}+\beta_{16}+\beta_{20}+\beta_{22}$ & A\\
00118 & $4\beta_{1}+3\beta_{2}+\beta_{10}+\beta_{20}+\beta_{28}+\beta_{30}$ & A\\
00119 & $4\beta_{1}+2\beta_{2}+3\beta_{3}+2\beta_{5}+\beta_{8}+2\beta_{9}+\beta_{29}$ & C\\
00120 & $4\beta_{1}+2\beta_{2}+3\beta_{3}+\beta_{5}+\beta_{8}+\beta_{29}+\beta_{30}$ & A\\
00121 & $4\beta_{1}+2\beta_{2}+3\beta_{3}+\beta_{5}+\beta_{20}+\beta_{25}+\beta_{26}$ & C\\
00122 & $4\beta_{1}+2\beta_{2}+3\beta_{3}+\beta_{8}+\beta_{22}+\beta_{40}$ & A\\
00123 & $4\beta_{1}+2\beta_{2}+2\beta_{3}+4\beta_{5}+\beta_{8}+\beta_{12}+\beta_{29}$ & C\\
00124 & $4\beta_{1}+2\beta_{2}+2\beta_{3}+\beta_{5}+\beta_{8}+2\beta_{9}+\beta_{11}+\beta_{28}$ & A\\
00125 & $4\beta_{1}+2\beta_{2}+2\beta_{3}+\beta_{5}+\beta_{12}+\beta_{16}+\beta_{20}+\beta_{22}$ & A\\
00126 & $4\beta_{1}+2\beta_{2}+\beta_{3}+\beta_{4}+\beta_{5}+2\beta_{9}+\beta_{20}+\beta_{30}$ & A\\
00127 & $4\beta_{1}+2\beta_{2}+\beta_{4}+3\beta_{5}+\beta_{12}+\beta_{20}+\beta_{30}$ & A\\
00128 & $4\beta_{1}+2\beta_{2}+\beta_{4}+2\beta_{5}+\beta_{8}+2\beta_{9}+\beta_{17}+\beta_{18}$ & A\\
00129 & $4\beta_{1}+2\beta_{2}+\beta_{4}+\beta_{7}+2\beta_{9}+\beta_{22}+\beta_{32}$ & A\\
00130 & $4\beta_{1}+2\beta_{2}+2\beta_{5}+\beta_{10}+\beta_{8}+2\beta_{11}+\beta_{29}$ & A\\
00131 & $4\beta_{1}+\beta_{2}+3\beta_{3}+\beta_{4}+2\beta_{5}+\beta_{7}+\beta_{20}+\beta_{33}$ & A\\
00132 & $4\beta_{1}+\beta_{2}+3\beta_{3}+2\beta_{5}+\beta_{16}+\beta_{20}+\beta_{29}$ & C\\
00133 & $4\beta_{1}+\beta_{2}+2\beta_{3}+\beta_{5}+\beta_{11}+\beta_{16}+\beta_{20}+\beta_{28}$ & A\\
00134 & $4\beta_{1}+\beta_{2}+\beta_{4}+2\beta_{5}+\beta_{17}+\beta_{16}+\beta_{20}+\beta_{18}$ & A\\
00135 & $4\beta_{1}+\beta_{2}+2\beta_{5}+\beta_{8}+\beta_{21}+\beta_{22}+\beta_{29}$ & A\\
00136 & $4\beta_{1}+3\beta_{3}+\beta_{8}+2\beta_{9}+\beta_{22}+\beta_{29}$ & A\\
00137 & $4\beta_{1}+2\beta_{3}+2\beta_{5}+\beta_{8}+\beta_{12}+\beta_{22}+\beta_{29}$ & A\\
00138 & $4\beta_{1}+2\beta_{5}+\beta_{31}+\beta_{29}+\beta_{33}$ & A\\
00139 & $3\beta_{1}+6\beta_{2}+\beta_{10}+\beta_{8}+\beta_{28}+\beta_{30}$ & A\\
00140 & $3\beta_{1}+5\beta_{2}+3\beta_{3}+\beta_{5}+\beta_{8}+\beta_{14}+2\beta_{15}$ & A\\
00141 & $3\beta_{1}+5\beta_{2}+3\beta_{3}+\beta_{5}+\beta_{8}+\beta_{25}+\beta_{26}$ & C\\
00142 & $3\beta_{1}+5\beta_{2}+3\beta_{3}+\beta_{5}+\beta_{11}+\beta_{20}+\beta_{28}$ & C\\
00143 & $3\beta_{1}+5\beta_{2}+\beta_{3}+\beta_{4}+\beta_{5}+\beta_{8}+2\beta_{9}+\beta_{30}$ & A\\
00144 & $3\beta_{1}+5\beta_{2}+\beta_{5}+\beta_{16}+\beta_{31}+\beta_{33}$ & A\\
\bottomrule\end{tabular}\end{minipage}\hfill%
\begin{minipage}[t]{.49\linewidth}\vspace{0pt}\small\setlength{\tabcolsep}{2pt}
\begin{tabular}{@{}p{9mm}p{62mm}p{3mm}@{}}\toprule Row & Basket & \\\midrule
00145 & $3\beta_{1}+4\beta_{2}+3\beta_{3}+\beta_{4}+2\beta_{5}+\beta_{7}+\beta_{8}+\beta_{33}$ & A\\
00146 & $3\beta_{1}+4\beta_{2}+3\beta_{3}+2\beta_{5}+\beta_{8}+\beta_{16}+\beta_{29}$ & C\\
00147 & $3\beta_{1}+4\beta_{2}+\beta_{3}+\beta_{4}+2\beta_{5}+\beta_{7}+\beta_{20}+\beta_{33}$ & A\\
00148 & $3\beta_{1}+4\beta_{2}+\beta_{3}+\beta_{4}+2\beta_{5}+3\beta_{11}+\beta_{20}$ & A\\
00149 & $3\beta_{1}+4\beta_{2}+\beta_{3}+\beta_{4}+\beta_{5}+\beta_{16}+\beta_{20}+\beta_{30}$ & A\\
00150 & $3\beta_{1}+4\beta_{2}+\beta_{4}+\beta_{7}+\beta_{16}+\beta_{22}+\beta_{32}$ & A\\
00151 & $3\beta_{1}+3\beta_{2}+\beta_{3}+\beta_{4}+\beta_{8}+4\beta_{9}+\beta_{22}$ & A\\
00152 & $3\beta_{1}+3\beta_{2}+\beta_{3}+\beta_{5}+\beta_{19}+\beta_{32}+\beta_{33}$ & A\\
00153 & $3\beta_{1}+3\beta_{2}+\beta_{4}+2\beta_{5}+\beta_{8}+2\beta_{9}+\beta_{12}+\beta_{22}$ & A\\
00154 & $3\beta_{1}+2\beta_{2}+3\beta_{3}+\beta_{8}+\beta_{16}+\beta_{22}+\beta_{29}$ & A\\
00155 & $3\beta_{1}+2\beta_{2}+\beta_{3}+\beta_{4}+2\beta_{5}+\beta_{8}+\beta_{29}+\beta_{30}$ & A\\
00156 & $3\beta_{1}+2\beta_{2}+\beta_{3}+\beta_{4}+\beta_{5}+\beta_{8}+\beta_{22}+\beta_{40}$ & A\\
00157 & $3\beta_{1}+2\beta_{2}+\beta_{3}+\beta_{4}+\beta_{8}+\beta_{20}+\beta_{22}+\beta_{32}$ & A\\
00158 & $3\beta_{1}+2\beta_{2}+\beta_{3}+\beta_{4}+2\beta_{9}+\beta_{16}+\beta_{20}+\beta_{22}$ & A\\
00159 & $3\beta_{1}+2\beta_{2}+\beta_{4}+2\beta_{5}+\beta_{12}+\beta_{16}+\beta_{20}+\beta_{22}$ & A\\
00160 & $3\beta_{1}+\beta_{2}+\beta_{3}+2\beta_{4}+3\beta_{5}+\beta_{7}+\beta_{20}+\beta_{33}$ & A\\
00161 & $3\beta_{1}+\beta_{2}+\beta_{10}+\beta_{8}+\beta_{22}+\beta_{28}+\beta_{29}$ & A\\
00162 & $3\beta_{1}+\beta_{3}+\beta_{4}+\beta_{5}+\beta_{8}+2\beta_{9}+\beta_{22}+\beta_{29}$ & A\\
00163 & $3\beta_{1}+\beta_{4}+3\beta_{5}+\beta_{8}+\beta_{12}+\beta_{22}+\beta_{29}$ & A\\
00164 & $2\beta_{1}+8\beta_{2}+3\beta_{3}+\beta_{5}+\beta_{8}+\beta_{11}+\beta_{28}$ & C\\
00165 & $2\beta_{1}+6\beta_{2}+3\beta_{3}+\beta_{5}+\beta_{8}+\beta_{9}+\beta_{37}$ & A\\
00166 & $2\beta_{1}+5\beta_{2}+\beta_{3}+\beta_{4}+2\beta_{5}+\beta_{8}+\beta_{14}+2\beta_{15}$ & B\\
00167 & $2\beta_{1}+5\beta_{2}+\beta_{3}+\beta_{4}+2\beta_{8}+\beta_{22}+\beta_{32}$ & A\\
00168 & $2\beta_{1}+5\beta_{2}+\beta_{3}+\beta_{4}+\beta_{8}+2\beta_{9}+\beta_{16}+\beta_{22}$ & A\\
00169 & $2\beta_{1}+4\beta_{2}+\beta_{3}+2\beta_{4}+3\beta_{5}+\beta_{7}+\beta_{8}+\beta_{33}$ & A\\
00170 & $2\beta_{1}+4\beta_{2}+\beta_{3}+\beta_{4}+2\beta_{16}+\beta_{20}+\beta_{22}$ & A\\
00171 & $2\beta_{1}+3\beta_{2}+\beta_{3}+\beta_{4}+\beta_{8}+\beta_{25}+\beta_{22}+\beta_{26}$ & A\\
00172 & $2\beta_{1}+3\beta_{2}+\beta_{3}+\beta_{4}+\beta_{11}+\beta_{20}+\beta_{22}+\beta_{28}$ & A\\
00173 & $2\beta_{1}+3\beta_{2}+\beta_{3}+\beta_{5}+2\beta_{7}+\beta_{32}+\beta_{33}$ & A\\
00174 & $2\beta_{1}+2\beta_{2}+2\beta_{3}+5\beta_{5}+\beta_{20}+\beta_{30}$ & A\\
00175 & $2\beta_{1}+2\beta_{2}+\beta_{3}+\beta_{4}+\beta_{5}+\beta_{8}+\beta_{16}+\beta_{22}+\beta_{29}$ & A\\
00176 & $2\beta_{1}+2\beta_{2}+\beta_{3}+4\beta_{5}+\beta_{7}+\beta_{22}+\beta_{32}$ & A\\
00177 & $2\beta_{1}+\beta_{2}+2\beta_{3}+\beta_{5}+\beta_{7}+2\beta_{9}+\beta_{20}+\beta_{33}$ & A\\
00178 & $2\beta_{1}+\beta_{2}+\beta_{3}+3\beta_{5}+\beta_{7}+\beta_{12}+\beta_{20}+\beta_{33}$ & A\\
00179 & $\beta_{1}+6\beta_{2}+\beta_{3}+\beta_{4}+2\beta_{5}+\beta_{8}+\beta_{9}+\beta_{37}$ & A\\
00180 & $\beta_{1}+6\beta_{2}+\beta_{3}+\beta_{4}+\beta_{8}+\beta_{11}+\beta_{22}+\beta_{28}$ & A\\
00181 & $\beta_{1}+5\beta_{2}+2\beta_{3}+3\beta_{5}+\beta_{12}+\beta_{39}$ & A\\
00182 & $\beta_{1}+4\beta_{2}+2\beta_{3}+\beta_{5}+\beta_{7}+\beta_{8}+2\beta_{9}+\beta_{33}$ & A\\
00183 & $\beta_{1}+3\beta_{2}+2\beta_{3}+4\beta_{5}+\beta_{8}+2\beta_{9}+\beta_{22}$ & A\\
00184 & $\beta_{1}+3\beta_{2}+2\beta_{3}+\beta_{5}+\beta_{7}+\beta_{16}+\beta_{20}+\beta_{33}$ & A\\
00185 & $\beta_{1}+2\beta_{2}+3\beta_{3}+2\beta_{5}+\beta_{29}+\beta_{39}$ & A\\
00186 & $\beta_{1}+2\beta_{2}+2\beta_{3}+4\beta_{5}+\beta_{16}+\beta_{20}+\beta_{22}$ & A\\
00187 & $\beta_{1}+\beta_{2}+3\beta_{3}+\beta_{5}+\beta_{8}+2\beta_{20}+\beta_{33}$ & A\\
00188 & $\beta_{1}+\beta_{2}+2\beta_{3}+2\beta_{5}+\beta_{7}+\beta_{8}+\beta_{29}+\beta_{33}$ & A\\
00189 & $\beta_{1}+\beta_{2}+\beta_{4}+2\beta_{5}+\beta_{7}+2\beta_{9}+\beta_{20}+\beta_{33}$ & A\\
00190 & $\beta_{1}+2\beta_{3}+5\beta_{5}+\beta_{8}+\beta_{22}+\beta_{29}$ & A\\
00191 & $7\beta_{2}+3\beta_{3}+\beta_{5}+\beta_{16}+\beta_{39}$ & A\\
00192 & $4\beta_{2}+3\beta_{3}+\beta_{5}+2\beta_{8}+\beta_{20}+\beta_{33}$ & A\\
\bottomrule\end{tabular}\end{minipage}\par

\clearpage
\noindent\textbf{Euler characteristic three (continued).}
\smallskip
\par\noindent
\begin{minipage}[t]{.49\linewidth}\vspace{0pt}\small\setlength{\tabcolsep}{2pt}
\begin{tabular}{@{}p{9mm}p{62mm}p{3mm}@{}}\toprule Row & Basket & \\\midrule
00193 & $4\beta_{2}+\beta_{4}+2\beta_{5}+\beta_{7}+\beta_{8}+2\beta_{9}+\beta_{33}$ & A\\
00194 & $3\beta_{2}+\beta_{4}+2\beta_{5}+\beta_{7}+\beta_{16}+\beta_{20}+\beta_{33}$ & A\\
00195 & $2\beta_{2}+\beta_{3}+\beta_{4}+3\beta_{5}+\beta_{29}+\beta_{39}$ & A\\
00196 & $\beta_{2}+\beta_{3}+\beta_{4}+2\beta_{5}+\beta_{8}+2\beta_{20}+\beta_{33}$ & A\\
00197 & $\beta_{2}+\beta_{4}+3\beta_{5}+\beta_{7}+\beta_{8}+\beta_{29}+\beta_{33}$ & A\\
47272 & $7\beta_{1}+2\beta_{2}+\beta_{3}+\beta_{4}+2\beta_{5}+\beta_{8}+\beta_{15}+\beta_{34}$ & B\\
47273 & $6\beta_{1}+5\beta_{2}+\beta_{3}+\beta_{4}+2\beta_{5}+\beta_{8}+\beta_{17}+\beta_{18}$ & A\\
47274 & $6\beta_{1}+\beta_{4}+\beta_{8}+\beta_{12}+2\beta_{16}+\beta_{29}$ & A\\
47275 & $5\beta_{1}+8\beta_{2}+2\beta_{3}+\beta_{4}+\beta_{5}+\beta_{8}+\beta_{30}$ & A\\
47276 & $5\beta_{1}+6\beta_{2}+\beta_{3}+\beta_{4}+2\beta_{5}+\beta_{8}+\beta_{12}+\beta_{22}$ & A\\
47277 & $5\beta_{1}+2\beta_{2}+\beta_{3}+3\beta_{5}+\beta_{8}+\beta_{12}+\beta_{17}+\beta_{18}$ & A\\
47278 & $5\beta_{1}+\beta_{2}+2\beta_{3}+\beta_{5}+\beta_{8}+\beta_{15}+\beta_{16}+\beta_{34}$ & A\\
47279 & $4\beta_{1}+8\beta_{2}+2\beta_{3}+\beta_{4}+\beta_{8}+\beta_{16}+\beta_{22}$ & A\\
47280 & $4\beta_{1}+8\beta_{2}+2\beta_{4}+2\beta_{5}+\beta_{8}+\beta_{30}$ & A\\
47281 & $4\beta_{1}+6\beta_{2}+2\beta_{5}+\beta_{8}+\beta_{21}+\beta_{30}$ & A\\
47282 & $4\beta_{1}+5\beta_{2}+2\beta_{3}+2\beta_{5}+\beta_{8}+\beta_{12}+\beta_{30}$ & A\\
47283 & $4\beta_{1}+4\beta_{2}+2\beta_{3}+\beta_{5}+\beta_{8}+\beta_{17}+\beta_{16}+\beta_{18}$ & A\\
47284 & $4\beta_{1}+4\beta_{2}+\beta_{3}+\beta_{4}+\beta_{5}+\beta_{20}+\beta_{41}$ & A\\
47285 & $4\beta_{1}+3\beta_{2}+\beta_{3}+3\beta_{5}+\beta_{8}+2\beta_{12}+\beta_{22}$ & A\\
47286 & $4\beta_{1}+\beta_{2}+\beta_{4}+2\beta_{5}+\beta_{8}+\beta_{15}+\beta_{16}+\beta_{34}$ & B\\
47287 & $3\beta_{1}+8\beta_{2}+2\beta_{4}+\beta_{5}+\beta_{8}+\beta_{16}+\beta_{22}$ & A\\
47288 & $3\beta_{1}+7\beta_{2}+\beta_{3}+\beta_{6}+4\beta_{5}+\beta_{8}+\beta_{22}$ & A\\
47289 & $3\beta_{1}+6\beta_{2}+\beta_{5}+\beta_{8}+\beta_{21}+\beta_{16}+\beta_{22}$ & A\\
47290 & $3\beta_{1}+5\beta_{2}+2\beta_{3}+\beta_{5}+\beta_{8}+\beta_{12}+\beta_{16}+\beta_{22}$ & A\\
47291 & $3\beta_{1}+5\beta_{2}+\beta_{4}+3\beta_{5}+\beta_{8}+\beta_{12}+\beta_{30}$ & A\\
47292 & $3\beta_{1}+4\beta_{2}+\beta_{4}+2\beta_{5}+\beta_{8}+\beta_{17}+\beta_{16}+\beta_{18}$ & A\\
47293 & $2\beta_{1}+7\beta_{2}+\beta_{3}+\beta_{4}+2\beta_{5}+\beta_{7}+\beta_{8}+\beta_{33}$ & A\\
47294 & $2\beta_{1}+7\beta_{2}+\beta_{3}+\beta_{4}+\beta_{5}+\beta_{8}+\beta_{16}+\beta_{30}$ & A\\
47295 & $2\beta_{1}+6\beta_{2}+\beta_{3}+\beta_{4}+5\beta_{5}+\beta_{8}+\beta_{22}$ & A\\
47296 & $2\beta_{1}+5\beta_{2}+\beta_{4}+2\beta_{5}+\beta_{8}+\beta_{12}+\beta_{16}+\beta_{22}$ & A\\
47297 & $2\beta_{1}+2\beta_{2}+\beta_{3}+6\beta_{5}+\beta_{8}+\beta_{17}+\beta_{18}$ & A\\
47298 & $\beta_{1}+7\beta_{2}+\beta_{3}+\beta_{4}+\beta_{8}+2\beta_{16}+\beta_{22}$ & A\\
47299 & $\beta_{1}+5\beta_{2}+2\beta_{3}+5\beta_{5}+\beta_{8}+\beta_{30}$ & A\\
47300 & $\beta_{1}+4\beta_{2}+\beta_{3}+3\beta_{5}+\beta_{7}+\beta_{8}+\beta_{12}+\beta_{33}$ & A\\
47301 & $\beta_{1}+3\beta_{2}+\beta_{3}+6\beta_{5}+\beta_{8}+\beta_{12}+\beta_{22}$ & A\\
47302 & $6\beta_{2}+2\beta_{3}+\beta_{5}+\beta_{7}+\beta_{8}+\beta_{16}+\beta_{33}$ & A\\
47303 & $5\beta_{2}+2\beta_{3}+4\beta_{5}+\beta_{8}+\beta_{16}+\beta_{22}$ & A\\
49933 & $11\beta_{1}+3\beta_{2}+3\beta_{3}+\beta_{6}+3\beta_{5}+\beta_{8}+\beta_{12}$ & A\\
49934 & $10\beta_{1}+3\beta_{2}+\beta_{5}+\beta_{10}+\beta_{8}+\beta_{16}+\beta_{27}$ & A\\
49935 & $10\beta_{1}+2\beta_{2}+\beta_{4}+\beta_{5}+\beta_{8}+\beta_{12}+\beta_{21}+\beta_{16}$ & A\\
49936 & $10\beta_{1}+\beta_{2}+\beta_{4}+\beta_{5}+2\beta_{12}+\beta_{16}+\beta_{20}$ & A\\
49937 & $9\beta_{1}+3\beta_{2}+\beta_{3}+\beta_{4}+2\beta_{5}+\beta_{8}+\beta_{11}+\beta_{27}$ & A\\
49938 & $9\beta_{1}+\beta_{2}+\beta_{3}+\beta_{4}+\beta_{8}+\beta_{12}+\beta_{16}+\beta_{29}$ & A\\
49939 & $9\beta_{1}+\beta_{7}+2\beta_{12}+\beta_{16}+\beta_{32}$ & A\\
49940 & $8\beta_{1}+4\beta_{2}+\beta_{6}+\beta_{5}+\beta_{7}+\beta_{16}+\beta_{32}$ & A\\
49941 & $8\beta_{1}+3\beta_{2}+\beta_{3}+\beta_{6}+\beta_{5}+\beta_{8}+\beta_{13}+\beta_{35}$ & A\\
\bottomrule\end{tabular}\end{minipage}\hfill%
\begin{minipage}[t]{.49\linewidth}\vspace{0pt}\small\setlength{\tabcolsep}{2pt}
\begin{tabular}{@{}p{9mm}p{62mm}p{3mm}@{}}\toprule Row & Basket & \\\midrule
49942 & $8\beta_{1}+2\beta_{2}+3\beta_{3}+\beta_{5}+\beta_{8}+\beta_{13}+\beta_{35}$ & A\\
49943 & $8\beta_{1}+2\beta_{2}+2\beta_{3}+\beta_{6}+3\beta_{5}+\beta_{8}+\beta_{12}+\beta_{16}$ & A\\
49944 & $8\beta_{1}+2\beta_{2}+\beta_{3}+\beta_{5}+\beta_{13}+\beta_{20}+\beta_{35}$ & A\\
49945 & $8\beta_{1}+2\beta_{2}+3\beta_{5}+\beta_{21}+\beta_{16}+\beta_{20}$ & A\\
49946 & $8\beta_{1}+2\beta_{2}+\beta_{5}+\beta_{8}+\beta_{9}+\beta_{16}+\beta_{36}$ & A\\
49947 & $8\beta_{1}+\beta_{2}+2\beta_{3}+3\beta_{5}+\beta_{12}+\beta_{16}+\beta_{20}$ & A\\
49948 & $8\beta_{1}+\beta_{2}+\beta_{3}+2\beta_{8}+2\beta_{12}+\beta_{32}$ & A\\
49949 & $8\beta_{1}+\beta_{2}+\beta_{3}+\beta_{8}+2\beta_{9}+2\beta_{12}+\beta_{16}$ & A\\
49950 & $8\beta_{1}+\beta_{2}+\beta_{5}+3\beta_{7}+\beta_{8}+\beta_{16}+\beta_{20}$ & A\\
49951 & $8\beta_{1}+\beta_{3}+2\beta_{12}+2\beta_{16}+\beta_{20}$ & A\\
49952 & $8\beta_{1}+\beta_{6}+4\beta_{5}+\beta_{8}+\beta_{12}+\beta_{29}$ & A\\
49953 & $8\beta_{1}+4\beta_{5}+\beta_{8}+\beta_{21}+\beta_{29}$ & A\\
49954 & $7\beta_{1}+3\beta_{2}+\beta_{3}+\beta_{6}+\beta_{5}+\beta_{8}+\beta_{25}+\beta_{26}$ & A\\
49955 & $7\beta_{1}+3\beta_{2}+\beta_{4}+2\beta_{5}+\beta_{7}+\beta_{16}+\beta_{32}$ & A\\
49956 & $7\beta_{1}+2\beta_{2}+\beta_{3}+\beta_{4}+2\beta_{5}+\beta_{8}+\beta_{13}+\beta_{35}$ & B\\
49957 & $7\beta_{1}+2\beta_{2}+\beta_{3}+\beta_{6}+2\beta_{5}+\beta_{8}+\beta_{16}+\beta_{29}$ & A\\
49958 & $7\beta_{1}+2\beta_{2}+\beta_{3}+\beta_{8}+2\beta_{12}+\beta_{11}+\beta_{28}$ & A\\
49959 & $7\beta_{1}+2\beta_{2}+\beta_{5}+\beta_{11}+\beta_{16}+\beta_{20}+\beta_{27}$ & A\\
49960 & $7\beta_{1}+2\beta_{5}+\beta_{8}+\beta_{11}+\beta_{27}+\beta_{29}$ & A\\
49961 & $6\beta_{1}+2\beta_{2}+\beta_{3}+\beta_{4}+2\beta_{5}+\beta_{8}+\beta_{25}+\beta_{26}$ & A\\
49962 & $6\beta_{1}+\beta_{2}+\beta_{3}+\beta_{4}+3\beta_{5}+\beta_{8}+\beta_{16}+\beta_{29}$ & A\\
49963 & $5\beta_{1}+5\beta_{2}+\beta_{5}+3\beta_{7}+\beta_{38}$ & C\\
49964 & $5\beta_{1}+2\beta_{2}+\beta_{3}+\beta_{5}+\beta_{7}+2\beta_{16}+\beta_{32}$ & A\\
49965 & $5\beta_{1}+2\beta_{2}+\beta_{5}+2\beta_{7}+2\beta_{8}+\beta_{38}$ & A\\
49966 & $5\beta_{1}+2\beta_{2}+\beta_{5}+\beta_{8}+2\beta_{9}+\beta_{13}+\beta_{35}$ & A\\
49967 & $5\beta_{1}+\beta_{2}+\beta_{3}+3\beta_{5}+2\beta_{8}+\beta_{12}+\beta_{32}$ & A\\
49968 & $5\beta_{1}+\beta_{2}+\beta_{3}+3\beta_{5}+\beta_{8}+2\beta_{9}+\beta_{12}+\beta_{16}$ & A\\
49969 & $5\beta_{1}+\beta_{2}+\beta_{5}+\beta_{13}+\beta_{16}+\beta_{20}+\beta_{35}$ & A\\
49970 & $4\beta_{1}+4\beta_{2}+\beta_{5}+\beta_{8}+\beta_{13}+\beta_{16}+\beta_{35}$ & A\\
49971 & $4\beta_{1}+3\beta_{2}+\beta_{3}+3\beta_{5}+\beta_{8}+\beta_{12}+2\beta_{16}$ & A\\
49972 & $4\beta_{1}+2\beta_{2}+\beta_{3}+3\beta_{5}+\beta_{8}+\beta_{12}+\beta_{11}+\beta_{28}$ & A\\
49973 & $4\beta_{1}+\beta_{2}+2\beta_{3}+\beta_{5}+\beta_{8}+\beta_{16}+\beta_{25}+\beta_{26}$ & A\\
49974 & $4\beta_{1}+2\beta_{3}+2\beta_{5}+\beta_{8}+2\beta_{16}+\beta_{29}$ & A\\
53332 & $5\beta_{1}+4\beta_{2}+3\beta_{3}+\beta_{20}+\beta_{41}$ & A\\
53447 & $7\beta_{1}+2\beta_{2}+3\beta_{3}+\beta_{5}+\beta_{8}+\beta_{25}+\beta_{26}$ & A\\
53448 & $7\beta_{1}+2\beta_{2}+\beta_{3}+\beta_{5}+2\beta_{7}+3\beta_{8}+\beta_{20}$ & A\\
53449 & $7\beta_{1}+\beta_{2}+3\beta_{3}+2\beta_{5}+\beta_{8}+\beta_{16}+\beta_{29}$ & A\\
53450 & $6\beta_{1}+6\beta_{2}+\beta_{3}+\beta_{6}+\beta_{5}+\beta_{8}+\beta_{11}+\beta_{28}$ & A\\
53451 & $6\beta_{1}+5\beta_{2}+3\beta_{3}+\beta_{5}+\beta_{8}+\beta_{11}+\beta_{28}$ & A\\
53452 & $6\beta_{1}+5\beta_{2}+\beta_{5}+\beta_{10}+\beta_{8}+\beta_{16}+\beta_{28}$ & A\\
53453 & $6\beta_{1}+2\beta_{2}+2\beta_{3}+\beta_{5}+\beta_{8}+\beta_{9}+\beta_{23}+\beta_{24}$ & A\\
53454 & $5\beta_{1}+5\beta_{2}+\beta_{3}+\beta_{4}+2\beta_{5}+\beta_{8}+\beta_{11}+\beta_{28}$ & A\\
53455 & $3\beta_{1}+6\beta_{2}+3\beta_{3}+4\beta_{5}+\beta_{8}+\beta_{22}$ & A\\
53456 & $3\beta_{1}+4\beta_{2}+2\beta_{3}+\beta_{5}+\beta_{8}+\beta_{11}+\beta_{16}+\beta_{28}$ & A\\
53777 & $7\beta_{1}+3\beta_{2}+\beta_{5}+\beta_{8}+2\beta_{9}+\beta_{11}+\beta_{27}$ & A\\
\bottomrule\end{tabular}\end{minipage}\par

\clearpage
\subsection{The first thirty-seven Euler-two exclusions}
The witnesses below give the initial reduction $65\to28$ used in
Section~\ref{sec:uniform-exclusions}. Their numerical criteria are proved in Part~II.

\begin{longtable}{@{}llp{.53\textwidth}@{}}
\caption{A sufficient witness for each of the first $37$ exclusions.}\label{tab:initial-cuts46}\\
\toprule Row&Test&Exact witness\\\midrule\endfirsthead
\toprule Row&Test&Exact witness\\\midrule\endhead
C2-00001 & surface profile & $m=8$, $(U_2,U_3,U_4)=(2,6,10)$, no profile\\
C2-00002 & complete system & $a=10$, $P_a=3$, lower bound $3/1210$\\
C2-00003 & surface profile & $m=8$, $(U_2,U_3,U_4)=(2,6,10)$, no profile\\
C2-00004 & complete system & $a=8$, $P_a=3$, lower bound $1/216$\\
C2-00006 & complete system & $a=12$, $P_a=4$, lower bound $5/2028$\\
C2-00007 & surface profile & $m=8$, $(U_2,U_3,U_4)=(3,7,11)$, no profile\\
C2-00008 & complete system & $a=8$, $P_a=3$, lower bound $1/216$\\
C2-00009 & complete system & $a=8$, $P_a=3$, lower bound $1/216$\\
C2-00010 & complete system & $a=12$, $P_a=4$, lower bound $5/2028$\\
C2-00011 & surface profile & $m=8$, $(U_2,U_3,U_4)=(3,7,11)$, no profile\\
C2-00012 & complete system & $a=12$, $P_a=4$, lower bound $5/2028$\\
C2-00013 & surface profile & $m=8$, $(U_2,U_3,U_4)=(2,7,10)$, no profile\\
C2-00014 & complete system & $a=8$, $P_a=3$, lower bound $1/216$\\
C2-00015 & differential & Commuting sections; see the proof of Proposition~\ref{prop:initial-cuts46}\\
C2-00016 & complete system & $a=9$, $P_a=3$, lower bound $1/300$\\
C2-00017 & complete system & $a=8$, $P_a=3$, lower bound $1/216$\\
C2-00018 & surface profile & $m=8$, $(U_2,U_3,U_4)=(4,7,11)$, no profile\\
C2-00019 & complete system & $a=10$, $P_a=3$, lower bound $3/1210$\\
C2-00020 & complete system & $a=10$, $P_a=3$, lower bound $3/1210$\\
C2-00021 & differential & Commuting sections; see the proof of Proposition~\ref{prop:initial-cuts46}\\
C2-00022 & complete system & $a=12$, $P_a=4$, lower bound $5/2028$\\
C2-00023 & complete system & $a=12$, $P_a=4$, lower bound $5/2028$\\
C2-00024 & complete system & $a=9$, $P_a=3$, lower bound $1/300$\\
C2-00025 & complete system & $a=11$, $P_a=4$, lower bound $5/1584$\\
C2-00026 & complete system & $a=8$, $P_a=3$, lower bound $1/216$\\
C2-00027 & surface profile & $m=8$, $(U_2,U_3,U_4)=(3,7,10)$, no profile\\
C2-00028 & complete system & $a=9$, $P_a=3$, lower bound $1/300$\\
C2-00029 & complete system & $a=9$, $P_a=3$, lower bound $1/300$\\
C2-01822 & complete system & $a=10$, $P_a=3$, lower bound $3/1210$\\
C2-01824 & surface profile & $m=9$, $(U_2,U_3,U_4)=(3,9,16)$, no profile\\
C2-01828 & surface profile & $m=9$, $(U_2,U_3,U_4)=(4,10,17)$, no profile\\
C2-01829 & complete system & $a=12$, $P_a=4$, lower bound $5/2028$\\
C2-01830 & complete system & $a=9$, $P_a=3$, lower bound $1/300$\\
C2-02067 & differential & Commuting sections; see the proof of Proposition~\ref{prop:initial-cuts46}\\
C2-02070 & complete system & $a=10$, $P_a=3$, lower bound $3/1210$\\
C2-02072 & complete system & $a=12$, $P_a=4$, lower bound $5/2028$\\
C2-02356 & differential & Commuting sections; see the proof of Proposition~\ref{prop:initial-cuts46}\\\bottomrule\end{longtable}


\subsection{The twenty-eight baskets after the initial tests}
This intermediate list records the destinations of every row before
the seventeen common applications. The endpoint needed in the main
proof is Table~\ref{tab:residual11}.

\begin{longtable}{@{}lrrrrp{.42\textwidth}@{}}
\caption{The remaining $28$ rows: index, first effective degree, section index, necessary torsion-state count and final argument.}\label{tab:remaining28}\\
\toprule Row&$I$&$e$&$\delta$&States&Resolution\\\midrule\endfirsthead
\toprule Row&$I$&$e$&$\delta$&States&Resolution\\\midrule\endhead
C2-00005 & 1260 & 4 & 8 & 1 & Single-visible adjunction\\
C2-01823 & 420 & 3 & 9 & 41 & Fixed-free quotient\\
C2-01825 & 5460 & 6 & 9 & 5 & K3 lattice and mixed characters\\
C2-01826 & 420 & 6 & 9 & 3 & Canonical torsion multiplication\\
C2-01827 & 420 & 3 & 9 & 87 & Fixed-free quotient\\
C2-02068 & 420 & 4 & 10 & 474 & Moving torsion Koszul\\
C2-02069 & 420 & 4 & 10 & 28 & Two-twist Koszul\\
C2-02071 & 420 & 4 & 10 & 11 & Prime restriction\\
C2-02073 & 420 & 5 & 10 & 79 & Fixed-free quotient\\
C2-02074 & 420 & 6 & 10 & 2 & Single-visible adjunction\\
C2-02075 & 420 & 5 & 10 & 9 & Complete surface classification\\
C2-02076 & 420 & 4 & 10 & 215 & Twisted quotient\\
C2-02077 & 420 & 5 & 10 & 12 & Prime restriction and invisible square\\
C2-02078 & 420 & 5 & 10 & 1 & Prime restriction\\
C2-02079 & 420 & 5 & 10 & 11 & Twisted Koszul\\
C2-02080 & 420 & 5 & 10 & 17 & Fixed-free quotient\\
C2-02346 & 420 & 4 & 12 & 51 & Twisted prime restriction\\
C2-02347 & 420 & 4 & 12 & 32 & Unique model $Z_4$\\
C2-02348 & 840 & 5 & 12 & 78 & Visible partition\\
C2-02349 & 420 & 4 & 12 & 28 & Quartic polar and mixed jet\\
C2-02350 & 420 & 8 & 12 & 4 & Complete octic conductor exclusion\\
C2-02351 & 420 & 5 & 12 & 43 & Visible partition\\
C2-02352 & 420 & 3 & 12 & 13 & Elliptic root cover\\
C2-02353 & 420 & 5 & 12 & 14 & Quintic polar and mixed jet\\
C2-02354 & 420 & 5 & 12 & 3 & Coupled prime characters\\
C2-02355 & 840 & 5 & 12 & 62 & Small marked rational curves\\
C2-02357 & 1260 & 5 & 12 & 62 & Small marked rational curves\\
C2-02397 & 840 & 5 & 14 & 11 & Unique model $Z_5$\\\bottomrule\end{longtable}



\subsection{Exact arithmetic for the Euler-one section}
\label{subsec:chi1-exact-arithmetic}

The finite statements used above are reproduced by
\texttt{verify\_chi1\_arithmetic46.py}, using only integer
arithmetic and \texttt{Fraction}.  Its inputs are the seven
explicit pairs in \eqref{eq:chi1-basket} and the five displayed
weights.  It imports no earlier data table or proof program.
The local-character computation is the following complete
finite algorithm:
\begin{enumerate}
\item For every tuple
$i\in\{0,1\}^3\times\{0,1,2,3\}\times\{0,\ldots,5\}
\times\{0,\ldots,4\}\times\{0,\ldots,6\}$, compute
\[
 j_i=\sum_Q(420/r_Q)b_Q^2i_Q\bmod420,
 \qquad
 M_i=-\sum_Q(420/r_Q)\overline{b_Qi_Q}
                  (r_Q-\overline{b_Qi_Q})\bmod840.
\]
Collect $M_i$ separately for each value of $j_i$.
\item Check that the collected set at $j$ is
$\{j^2,j^2+210,j^2+420,j^2+630\}\bmod840$ for each
$0\leq j<420$.  The four tuples with $j_i=M_i=0$ are the
ones displayed above.  Evaluate $C$ and $\delta_m$ on these
tuples; their common period is six, with the increments
specified in \eqref{eq:chi1-period-six}.
\item Pair states only when their correction masses agree.
There are ten ordered pairs: one zero pair and nine nonzero
pairs.  Formula \eqref{eq:chi1-paired-anchor} gives their
anchor sums $0$ and $1$, respectively.
\end{enumerate}
The residue check ranges over $6720$ tuples, not a search
with an unproved cutoff.  Its output includes the plurigenera,
the section bounds and the full monomial lists used in
\eqref{eq:chi1-monomial-count}.

As an independent arithmetic consistency check, it verifies
the equality of rational generating functions
\[
 H_R(s)=\frac{1+s^{23}}{(1-s^4)(1-s^5)(1-s^6)(1-s^7)}.
\]
The Riemann--Roch side is formed as
\[
 \frac{s^2(1+s)}{840(1-s)^4}
 +\frac{1-3s}{(1-s)^2}+s
 +\frac{s}{1-s}\sum_{(b,r)\in\mathcal B}
   \frac{\sum_{j=1}^r\Phi_{b,r}(j)s^j}{1-s^r}.
\]
Cross multiplication proves the identity as polynomials with
rational coefficients.  This check is not used as a substitute
for the algebraic-independence and finite-generation arguments.
The determinant-$46$ reference matrix, its integral coordinates
for all $37$ monomials, the index-two character sublattice,
and the five determinant-one function-field frames are also
verified exactly.


\section{Defining coefficients for Vega and Altair}\label{app:coefficients}

\subsection{Vega: exact quartic construction coefficients}
\label{app:quartic-exact-coefficients}
All coefficients in this appendix are rational.  The parameters and basis elements have the weights specified in the construction theorem.  Juxtaposition of a parameter and a basis element denotes the free module action.  These formulas, together with the stated reconstruction rule for the top basis element, define every invariant product.  The cover is the entire integral closure described in the text.

Use $S=A\langle1,h,q,r,s,k,\ell,t\rangle$, where $A=\mathbb C[a,u,b,c]$.  In the formula for $h\ell$, the coefficient of $t$ is one.  Write $h\ell=t+\Psi$.  First define $th=(h^2)\ell-h\Psi$ and $t\ell=h(\ell^2)-\ell\Psi$; for $x=q,r,s,k$ define $tx=h(\ell x)-\Psi x$, in that order.  Finally put $t^2=(th)\ell-t\Psi$.  Every right side reduces using products already defined.  Equivalently these are the uniquely determined products in the attached full free-module table.
\paragraph{The invariant multiplication}
\begingroup\small
\[
\begin{aligned}
h^2&=-27u^{3}+\frac{9}{4}abc+\frac{1430721}{11200}a^{2}ub\\
&\quad+\frac{702653211}{2560000}a^{4}c-\frac{27119976024987}{10035200000}a^{5}u+bq\\
&\quad-\frac{2151927}{7000}a^{3}q-6ur+\frac{81}{5}a\ell
\end{aligned}
\]

\[
\begin{aligned}
hq&=\frac{21}{20}ach-\frac{1400757}{22400}a^{2}uh+bs\\
&\quad-\frac{122391}{16000}a^{3}s-2uk
\end{aligned}
\]

\[
\begin{aligned}
hr&=-\frac{5}{2}u^{2}h-\frac{208167}{6400}a^{2}bh+\frac{1508879902059}{2508800000}a^{5}h\\
&\quad+2cs-\frac{33}{5}aus-bk\\
&\quad+\frac{122391}{16000}a^{3}k
\end{aligned}
\]

\[
\begin{aligned}
hs&=b^{3}-7ubc+\frac{263}{8}au^{2}b\\
&\quad+\frac{9477}{200}a^{2}c^{2}+\frac{5456403}{224000}a^{3}b^{2}-\frac{25928829}{51200}a^{3}uc\\
&\quad+\frac{35037437517}{20480000}a^{4}u^{2}-\frac{147285185224293}{25088000000}a^{6}b+\frac{2090907395768804409}{25088000000000}a^{9}\\
&\quad-\frac{81}{10}acq+\frac{1592253}{11200}a^{2}uq+\frac{69}{5}abr\\
&\quad-\frac{3517911}{16000}a^{4}r-6u\ell
\end{aligned}
\]

\[
\begin{aligned}
hk&=b^{2}c-6uc^{2}-\frac{3209}{240}aub^{2}\\
&\quad+\frac{927}{20}au^{2}c+\frac{2010879}{2800}a^{2}u^{3}+\frac{511524963}{3584000}a^{3}bc\\
&\quad-\frac{8780980348881}{2007040000}a^{4}ub-\frac{5395477895777961}{401408000000}a^{6}c+\frac{22875302772090097731}{224788480000000}a^{7}u\\
&\quad+6u^{2}q-\frac{137313}{3200}a^{2}bq+\frac{21556156571937}{2508800000}a^{5}q\\
&\quad+\frac{81}{10}acr+\frac{143973}{1120}a^{2}ur-b\ell\\
&\quad-\frac{47432709}{112000}a^{3}\ell
\end{aligned}
\]

\[
\begin{aligned}
h\ell&=-\frac{2801}{240}aubh-\frac{10239363}{512000}a^{3}ch-\frac{2608579173573}{2007040000}a^{4}uh\\
&\quad-2u^{2}s+\frac{503763}{22400}a^{2}bs-\frac{2480587944363}{5017600000}a^{5}s\\
&\quad+\frac{21}{20}ack-\frac{706017}{11200}a^{2}uk+t
\end{aligned}
\]

\[
\begin{aligned}
q^2&=b^{3}-3ubc+\frac{397}{24}au^{2}b\\
&\quad-\frac{2457}{400}a^{2}c^{2}-\frac{989901}{22400}a^{3}b^{2}+\frac{63177003}{256000}a^{3}uc\\
&\quad-\frac{5748444477}{4096000}a^{4}u^{2}+\frac{16894948508253}{12544000000}a^{6}b-\frac{1638634685476436667}{200704000000000}a^{9}\\
&\quad+\frac{69}{10}acq-\frac{25653}{280}a^{2}uq-\frac{12}{5}abr\\
&\quad+\frac{367173}{20000}a^{4}r+2u\ell
\end{aligned}
\]

\[
\begin{aligned}
qr&=b^{2}c-2uc^{2}+\frac{1193}{240}aub^{2}\\
&\quad+\frac{321}{40}au^{2}c-\frac{315531}{6400}a^{2}u^{3}+\frac{21502719}{512000}a^{3}bc\\
&\quad-\frac{178070154111}{286720000}a^{4}ub-\frac{43009786188027}{57344000000}a^{6}c+\frac{384469681067412309}{32112640000000}a^{7}u\\
&\quad-\frac{5}{2}u^{2}q-\frac{182373}{6400}a^{2}bq+\frac{357883635609}{627200000}a^{5}q\\
&\quad+\frac{21}{20}acr-\frac{2727}{128}a^{2}ur+b\ell\\
&\quad-\frac{122391}{16000}a^{3}\ell
\end{aligned}
\]

\[
\begin{aligned}
qs&=b^{2}h-uch-\frac{3}{5}au^{2}h\\
&\quad+\frac{2325123}{56000}a^{3}bh-\frac{674005156353}{1792000000}a^{6}h+\frac{21}{20}acs\\
&\quad-\frac{2727}{128}a^{2}us+\frac{12}{5}abk-\frac{367173}{20000}a^{4}k
\end{aligned}
\]

\[
\begin{aligned}
qk&=bch+\frac{1}{10}aubh+\frac{8835579}{89600}a^{3}ch\\
&\quad+\frac{280309165551}{250880000}a^{4}uh+2u^{2}s-\frac{297021}{11200}a^{2}bs\\
&\quad+\frac{2635278663609}{5017600000}a^{5}s+\frac{24}{5}ack+\frac{760551}{22400}a^{2}uk\\
&\quad-t
\end{aligned}
\]

\[
\begin{aligned}
q\ell&=\frac{9}{2}u^{2}b^{2}-\frac{9}{2}u^{3}c+\frac{69}{20}abc^{2}\\
&\quad-\frac{27}{10}au^{4}+\frac{145509}{6400}a^{2}b^{3}-\frac{4596013}{44800}a^{2}ubc\\
&\quad+\frac{4410131}{14336}a^{3}u^{2}b-\frac{1240548777}{10240000}a^{4}c^{2}-\frac{650099919573}{501760000}a^{5}b^{2}\\
&\quad+\frac{406964060804979}{80281600000}a^{5}uc-\frac{161444020491838881}{6422528000000}a^{6}u^{2}+\frac{19095326532455563359}{561971200000000}a^{8}b\\
&\quad-\frac{109182742903151630098371}{561971200000000000}a^{11}-bcq-\frac{617}{240}aubq\\
&\quad+\frac{350852067}{3584000}a^{3}cq-\frac{2612469023397}{2007040000}a^{4}uq+b^{2}r\\
&\quad-ucr-\frac{3}{5}au^{2}r-\frac{3738843}{56000}a^{3}br\\
&\quad+\frac{810344569059}{1792000000}a^{6}r+\frac{21}{20}ac\ell+\frac{124113}{4480}a^{2}u\ell
\end{aligned}
\]

\[
\begin{aligned}
r^2&=bc^{2}-\frac{45}{4}u^{4}+\frac{12}{5}ab^{3}\\
&\quad-\frac{553}{120}aubc-\frac{265111}{6400}a^{2}u^{2}b+\frac{5046327}{51200}a^{3}c^{2}\\
&\quad-\frac{52898084613}{401408000}a^{4}b^{2}-\frac{22907883591}{20480000}a^{4}uc+\frac{7793301012993}{1433600000}a^{5}u^{2}\\
&\quad+\frac{404770266294517881}{56197120000000}a^{7}b-\frac{1087658263162276497610713}{6294077440000000000}a^{10}+ubq\\
&\quad-\frac{109779}{2240}a^{2}cq+\frac{8628399}{56000}a^{3}uq-7u^{2}r\\
&\quad-\frac{356229}{22400}a^{2}br+\frac{1037271885669}{1254400000}a^{5}r+2c\ell\\
&\quad-\frac{33}{5}au\ell
\end{aligned}
\]

\[
\begin{aligned}
rs&=bch-\frac{23}{2}aubh-\frac{17433}{64000}a^{3}ch\\
&\quad-\frac{1220932089}{7168000}a^{4}uh-\frac{13}{2}u^{2}s-\frac{53817}{8960}a^{2}bs\\
&\quad+\frac{382481140509}{5017600000}a^{5}s+\frac{21}{20}ack-\frac{11223}{800}a^{2}uk\\
&\quad+t
\end{aligned}
\]

\[
\begin{aligned}
rk&=c^{2}h-\frac{12}{5}ab^{2}h-\frac{27}{10}auch\\
&\quad-\frac{99}{50}a^{2}u^{2}h+\frac{82050853833}{125440000}a^{4}bh-\frac{9760430709201807}{2007040000000}a^{7}h\\
&\quad-ubs-\frac{230883}{5600}a^{2}cs+\frac{3219003}{22400}a^{3}us\\
&\quad-\frac{5}{2}u^{2}k+\frac{744711}{44800}a^{2}bk+\frac{565663869279}{2508800000}a^{5}k
\end{aligned}
\]

\[
\begin{aligned}
r\ell&=ub^{3}+\frac{3}{2}u^{2}bc+\frac{117}{20}ac^{3}\\
&\quad-\frac{17053}{480}au^{3}b+\frac{41829}{6400}a^{2}b^{2}c-\frac{1928823}{22400}a^{2}uc^{2}\\
&\quad+\frac{291046939}{3584000}a^{3}ub^{2}+\frac{968173839}{7168000}a^{3}u^{2}c+\frac{77185095129}{573440000}a^{4}u^{3}\\
&\quad+\frac{87994060939269}{160563200000}a^{5}bc-\frac{72499046377747797}{12845056000000}a^{6}ub-\frac{64547919193140459}{642252800000000}a^{8}c\\
&\quad+\frac{5535845749474369075539}{71932313600000000}a^{9}u-c^{2}q+\frac{12}{5}ab^{2}q\\
&\quad+\frac{27}{10}aucq+\frac{99}{50}a^{2}u^{2}q-\frac{236961950913}{501760000}a^{4}bq\\
&\quad+\frac{27874588765101399}{8028160000000}a^{7}q+bcr-\frac{3353}{240}aubr\\
&\quad-\frac{10412667}{512000}a^{3}cr-\frac{1056232251}{286720000}a^{4}ur-\frac{5}{2}u^{2}\ell\\
&\quad+\frac{564153}{44800}a^{2}b\ell+\frac{321504614451}{1254400000}a^{5}\ell
\end{aligned}
\]

\[
\begin{aligned}
s^2&=-20u^{3}b+\frac{93}{20}ab^{2}c-\frac{189}{20}auc^{2}\\
&\quad+\frac{125213}{22400}a^{2}ub^{2}+\frac{357453}{3200}a^{2}u^{2}c+\frac{2434131}{16000}a^{3}u^{3}\\
&\quad+\frac{13039385913}{17920000}a^{4}bc-\frac{33114002760819}{10035200000}a^{5}ub-\frac{37659239734228641}{2007040000000}a^{7}c\\
&\quad+\frac{12326003449611845481}{160563200000000}a^{8}u+b^{2}q+3ucq\\
&\quad+\frac{9}{5}au^{2}q-\frac{40458141}{112000}a^{3}bq+\frac{2466566870121}{448000000}a^{6}q\\
&\quad-8ubr+\frac{486}{25}a^{2}cr+\frac{17739}{400}a^{3}ur\\
&\quad+\frac{69}{5}ab\ell-\frac{3517911}{16000}a^{4}\ell
\end{aligned}
\]

\[
\begin{aligned}
sk&=\frac{7}{2}u^{2}b^{2}-\frac{27}{2}u^{3}c+\frac{18}{5}abc^{2}\\
&\quad-\frac{81}{10}au^{4}-\frac{84027}{1280}a^{2}b^{3}+\frac{9366801}{44800}a^{2}ubc\\
&\quad-\frac{353070663}{448000}a^{3}u^{2}b-\frac{30769214283}{35840000}a^{4}c^{2}+\frac{4147666443441}{2508800000}a^{5}b^{2}\\
&\quad+\frac{18402544511307}{2007040000}a^{5}uc-\frac{11286156138822477}{321126400000}a^{6}u^{2}+\frac{50508440768985692229}{561971200000000}a^{8}b\\
&\quad-\frac{7129376146050937920674109}{4495769600000000000}a^{11}+2bcq-6aubq\\
&\quad+\frac{57429}{32000}a^{3}cq-\frac{373147772679}{125440000}a^{4}uq-b^{2}r\\
&\quad-3ucr-\frac{9}{5}au^{2}r-\frac{7319367}{28000}a^{3}br\\
&\quad+\frac{7910809765371}{1792000000}a^{6}r+\frac{81}{10}ac\ell+\frac{143973}{1120}a^{2}u\ell
\end{aligned}
\]

\[
\begin{aligned}
s\ell&=\frac{12}{5}ac^{2}h-\frac{31329}{3200}a^{2}b^{2}h-\frac{7839}{3200}a^{2}uch\\
&\quad-\frac{37341}{16000}a^{3}u^{2}h+\frac{70372428507}{50176000}a^{5}bh-\frac{407652162531294249}{40140800000000}a^{8}h\\
&\quad+bcs-\frac{3353}{240}aubs-\frac{10412667}{512000}a^{3}cs\\
&\quad-\frac{1056232251}{286720000}a^{4}us-b^{2}k+uck\\
&\quad+\frac{3}{5}au^{2}k+\frac{3738843}{56000}a^{3}bk-\frac{810344569059}{1792000000}a^{6}k
\end{aligned}
\]

\[
\begin{aligned}
k^2&=ub^{3}+\frac{99}{20}ac^{3}-\frac{2201}{40}au^{3}b\\
&\quad-\frac{428199}{3200}a^{2}b^{2}c+\frac{1462311}{4480}a^{2}uc^{2}+\frac{115164197}{89600}a^{3}ub^{2}\\
&\quad-\frac{5429609631}{1792000}a^{3}u^{2}c-\frac{15671616642027}{1003520000}a^{4}u^{3}-\frac{83585535187257}{20070400000}a^{5}bc\\
&\quad+\frac{1230839253429525963}{11239424000000}a^{6}ub+\frac{507516793995702226419}{1123942400000000}a^{8}c-\frac{1942367993972916760111419}{629407744000000000}a^{9}u\\
&\quad+3c^{2}q+\frac{12}{5}ab^{2}q-\frac{153}{5}aucq\\
&\quad-\frac{1244871}{11200}a^{2}u^{2}q+\frac{23063845923}{31360000}a^{4}bq-\frac{3187220665114915731}{14049280000000}a^{7}q\\
&\quad-2bcr+6aubr-\frac{10283193}{28000}a^{3}cr\\
&\quad-\frac{177195854043}{62720000}a^{4}ur-6u^{2}\ell+\frac{492963}{5600}a^{2}b\ell\\
&\quad+\frac{3297779067339}{313600000}a^{5}\ell
\end{aligned}
\]

\[
\begin{aligned}
k\ell&=-ub^{2}h+u^{2}ch+\frac{3}{5}au^{3}h\\
&\quad+\frac{12447}{3200}a^{2}bch+\frac{44351493}{224000}a^{3}ubh+\frac{24305375442771}{10035200000}a^{5}ch\\
&\quad+\frac{175706452040760243}{5619712000000}a^{6}uh+c^{2}s-\frac{12}{5}ab^{2}s\\
&\quad-\frac{27}{10}aucs+\frac{197697}{2240}a^{2}u^{2}s-\frac{27217166643}{50176000}a^{4}bs\\
&\quad+\frac{2116789249726434357}{112394240000000}a^{7}s-bck-\frac{617}{240}aubk\\
&\quad+\frac{181061019}{3584000}a^{3}ck+\frac{3095881754499}{2007040000}a^{4}uk-\frac{1010661}{22400}a^{2}t
\end{aligned}
\]

\[
\begin{aligned}
\ell^2&=uc^{3}+\frac{12}{5}ab^{4}-\frac{291}{20}aub^{2}c\\
&\quad+\frac{39}{4}au^{2}c^{2}+\frac{14389093}{80640}a^{2}u^{2}b^{2}-\frac{166833}{700}a^{2}u^{3}c\\
&\quad+\frac{41664663}{224000}a^{3}bc^{2}-\frac{128007}{875}a^{3}u^{4}+\frac{46941533271}{125440000}a^{4}b^{3}\\
&\quad-\frac{148336307673}{40140800}a^{4}ubc+\frac{940197688278189}{80281600000}a^{5}u^{2}b-\frac{67627404470648457}{12845056000000}a^{6}c^{2}\\
&\quad-\frac{240399359030952747}{8028160000000}a^{7}b^{2}+\frac{508794084347722069419}{3596615680000000}a^{7}uc-\frac{346592242316967258289527}{575458508800000000}a^{8}u^{2}\\
&\quad+\frac{1204084342980121941202041}{1573519360000000000}a^{10}b-\frac{430719084951939040101962880213}{100705239040000000000000}a^{13}+ub^{2}q\\
&\quad-u^{2}cq-\frac{3}{5}au^{3}q-\frac{109779}{2240}a^{2}bcq\\
&\quad+\frac{23809329}{112000}a^{3}ubq+\frac{726841160271}{250880000}a^{5}cq-\frac{176197866979881951}{5619712000000}a^{6}uq\\
&\quad+\frac{12}{5}ac^{2}r+\frac{395679}{11200}a^{2}b^{2}r-\frac{532767}{11200}a^{2}ucr\\
&\quad-\frac{329337}{11200}a^{3}u^{2}r-\frac{504848023137}{313600000}a^{5}br+\frac{65813038396551}{6422528000}a^{8}r\\
&\quad-\frac{1697}{120}aub\ell-\frac{10256283}{256000}a^{3}c\ell+\frac{43277246847}{40140800}a^{4}u\ell
\end{aligned}
\]

\endgroup
\paragraph{The quadratic field and the full anti module}
Let $v_0$ denote the first anti generator (called $v$ in the quartic construction and $w$ in the quintic construction).  The following polynomials $X_d$ are its products with the twelve anti generators.  Thus $X_{\deg v_0}=F=v_0^2$, and the complete module is generated inside $\operatorname{Frac}(S)(v_0)$ by $n_d=X_d/v_0$.  Multiplication is unambiguously $(X_d/v_0)(X_e/v_0)=X_dX_e/F$.  The polynomial product and syzygy identities are verified by the exact verifier; the proof shows that this is the full colon module, not a chosen submodule.
\begingroup\small
\[
\begin{aligned}
v_0v&=u^{3}-\frac{108}{25}a^{2}ub-\frac{972}{625}a^{4}c\\
&\quad+\frac{9983169}{100000}a^{5}u+\frac{432}{125}a^{3}q
\end{aligned}
\]

\[
\begin{aligned}
v_0w&=u^{2}b-\frac{72}{25}a^{2}b^{2}-\frac{207}{50}a^{2}uc\\
&\quad+\frac{355779}{16000}a^{3}u^{2}+\frac{55266219}{700000}a^{5}b-\frac{2559658564863}{2450000000}a^{8}\\
&\quad+\frac{6}{5}auq+\frac{216}{125}a^{3}r
\end{aligned}
\]

\[
\begin{aligned}
v_0\zeta&=u^{2}c+3au^{3}-\frac{72}{25}a^{2}bc\\
&\quad-\frac{603921}{112000}a^{3}ub-\frac{137633499}{3200000}a^{5}c+\frac{4320919034073}{12544000000}a^{6}u\\
&\quad+\frac{5738391}{140000}a^{4}q+\frac{6}{5}aur-\frac{36}{25}a^{2}\ell
\end{aligned}
\]

\[
\begin{aligned}
v_0d&=u^{2}h-\frac{72}{25}a^{2}bh-\frac{32958333}{350000}a^{5}h\\
&\quad+\frac{6}{5}aus-\frac{648}{125}a^{3}k
\end{aligned}
\]

\[
\begin{aligned}
v_0p&=\frac{977}{240}au^{2}b-\frac{27}{10}a^{2}c^{2}-\frac{136323}{22400}a^{3}b^{2}\\
&\quad+\frac{34749891}{512000}a^{3}uc-\frac{18013827951}{40960000}a^{4}u^{2}+\frac{1295195440041}{2508800000}a^{6}b\\
&\quad-\frac{14422585817323809}{4014080000000}a^{9}+\frac{6}{5}acq-\frac{109779}{4480}a^{2}uq\\
&\quad-\frac{6}{5}abr+\frac{367173}{40000}a^{4}r+u\ell
\end{aligned}
\]

\[
\begin{aligned}
v_0e&=\frac{6}{5}aub^{2}-\frac{69}{20}au^{2}c+\frac{15939}{640}a^{2}u^{3}\\
&\quad+\frac{162}{125}a^{3}bc-\frac{12582621}{280000}a^{4}ub-\frac{9913671}{1000000}a^{6}c\\
&\quad-\frac{1469228150583}{6272000000}a^{7}u+u^{2}q-\frac{72}{25}a^{2}bq\\
&\quad+\frac{1101519}{50000}a^{5}q+\frac{36}{25}a^{2}ur
\end{aligned}
\]

\[
\begin{aligned}
v_0f&=\frac{6}{5}aubc+\frac{282533}{44800}a^{2}u^{2}b+\frac{162}{125}a^{3}c^{2}\\
&\quad-\frac{147388329}{1280000}a^{4}uc+\frac{188017625331}{5017600000}a^{5}u^{2}-\frac{72}{25}a^{2}cq\\
&\quad+\frac{1428957}{56000}a^{3}uq+u^{2}r-\frac{6}{5}au\ell
\end{aligned}
\]

\[
\begin{aligned}
v_0g&=\frac{6}{5}aubh-\frac{216}{125}a^{3}ch+\frac{5160969}{280000}a^{4}uh\\
&\quad+u^{2}s-\frac{72}{25}a^{2}bs+\frac{1101519}{50000}a^{5}s\\
&\quad+\frac{36}{25}a^{2}uk
\end{aligned}
\]

\[
\begin{aligned}
v_0\eta&=\frac{69}{20}auc^{2}+\frac{524453}{44800}a^{2}ub^{2}-\frac{5175}{128}a^{2}u^{2}c\\
&\quad+\frac{37179}{400}a^{3}u^{3}-\frac{659066463}{8960000}a^{4}bc-\frac{13083555537}{100352000}a^{5}ub\\
&\quad+\frac{1104064897167711}{1003520000000}a^{7}c-\frac{37272582186665283}{11468800000000}a^{8}u-ucq\\
&\quad+\frac{9}{2}au^{2}q+\frac{703197}{56000}a^{3}bq-\frac{145282645467}{896000000}a^{6}q\\
&\quad+ubr-\frac{36}{25}a^{2}cr-\frac{18711}{16000}a^{3}ur\\
&\quad-\frac{6}{5}ab\ell+\frac{367173}{40000}a^{4}\ell
\end{aligned}
\]

\[
\begin{aligned}
v_0i&=\frac{6}{5}auch+\frac{230883}{11200}a^{2}u^{2}h-\frac{2767419}{70000}a^{4}bh\\
&\quad+\frac{338707178829}{1120000000}a^{7}h-\frac{36}{25}a^{2}cs+\frac{594}{125}a^{3}us\\
&\quad+u^{2}k-\frac{36}{25}a^{2}bk+\frac{1101519}{100000}a^{5}k
\end{aligned}
\]

\[
\begin{aligned}
v_0j&=-\frac{197283}{11200}a^{2}ubh+\frac{106083729}{1120000}a^{4}ch+\frac{70804609299}{179200000}a^{5}uh\\
&\quad+ucs+\frac{3}{10}au^{2}s-\frac{81063}{14000}a^{3}bs\\
&\quad+\frac{5407806796551}{12544000000}a^{6}s-ubk+\frac{153}{50}a^{2}ck\\
&\quad+\frac{343359}{16000}a^{3}uk-\frac{6}{5}at
\end{aligned}
\]

\[
\begin{aligned}
v_0\kappa&=-\frac{11}{5}au^{2}bh-\frac{36}{25}a^{2}c^{2}h+\frac{922473}{28000}a^{3}b^{2}h\\
&\quad-\frac{144509751}{2240000}a^{3}uch-\frac{59747975001}{179200000}a^{4}u^{2}h-\frac{37634130981}{160000000}a^{6}bh\\
&\quad-\frac{4606078924895571}{35840000000000}a^{9}h-\frac{6}{5}abcs+\frac{813}{5600}a^{2}ubs\\
&\quad+\frac{3794121}{200000}a^{4}cs-\frac{8925500783457}{25088000000}a^{5}us+\frac{6}{5}ab^{2}k\\
&\quad-\frac{27}{20}auck-\frac{71793}{4000}a^{2}u^{2}k-\frac{856737}{100000}a^{4}bk\\
&\quad-\frac{14979556881}{3200000000}a^{7}k+ut
\end{aligned}
\]

\endgroup
\paragraph{The five minimal module relations}
The relations below are relations over the polynomial parameter ring.  They must be retained together with the product equations.  Their unique second syzygy and the Hilbert-series argument in the text prove that they present the entire anti module.
\begingroup\small
\[
\begin{aligned}
0&=\frac{5}{6}bv-\frac{40797}{6400}a^{3}v-\frac{5}{6}uw\\
&\quad+ae
\end{aligned}
\]

\[
\begin{aligned}
0&=\frac{5}{6}cv+\frac{5}{2}auv-\frac{5}{6}u\zeta\\
&\quad+af
\end{aligned}
\]

\[
\begin{aligned}
0&=\frac{25}{4}abv-\frac{367173}{12800}a^{4}v+\frac{5}{6}cw\\
&\quad-\frac{15}{4}auw-\frac{5}{6}b\zeta+\frac{40797}{6400}a^{3}\zeta\\
&\quad+a\eta
\end{aligned}
\]

\[
\begin{aligned}
0&=\frac{15}{4}ubv-\frac{40797}{6400}a^{2}cv-\frac{367173}{12800}a^{3}uv\\
&\quad-\frac{15}{4}u^{2}w+\frac{40797}{6400}a^{2}u\zeta+ce\\
&\quad-bf+u\eta
\end{aligned}
\]

\[
\begin{aligned}
0&=-\frac{5}{6}cg-\frac{1}{4}aug+\frac{5}{6}bi\\
&\quad+\frac{13599}{32000}a^{3}i+\frac{5}{6}uj+a\kappa
\end{aligned}
\]

\endgroup


\subsection{Altair: exact quintic construction coefficients}
\label{app:quintic-exact-coefficients}
All coefficients in this appendix are rational.  The parameters and basis elements have the weights specified in the construction theorem.  Juxtaposition of a parameter and a basis element denotes the free module action.  These formulas, together with the stated reconstruction rule for the top basis element, define every invariant product.  The cover is the entire integral closure described in the text.

Use $S=A\langle1,q,k,t_\flat,e,f,j,\ell\rangle$, where $A=\mathbb C[a,b,u,z]$.  We retain the auxiliary homogeneous basis element $t_\flat$ so that the same exact data also define the module and deformation calculation.  It is not a further algebra generator.  The displayed formula $qk=t_\flat+\Psi$ defines it.  First define $t_\flat q=(q^2)k-q\Psi$ and $t_\flat k=q(k^2)-k\Psi$; for $x=e,f,j,\ell$ put $t_\flat x=q(kx)-\Psi x$.  Finally put $t_\flat^2=(t_\flat q)k-t_\flat\Psi$.  Thus this basis convention differs from $t=qk$ by the explicitly displayed triangular change of basis.
\paragraph{The invariant multiplication}
\begingroup\small
\[
\begin{aligned}
q^2&=-\frac{4}{121}u^{3}+bz^{2}-\frac{7}{11}b^{3}u\\
&\quad+\frac{16}{121}a^{2}uz-\frac{91258289}{19360000}a^{4}b^{2}-\frac{9583}{2200}a^{2}bq\\
&\quad-\frac{8}{11}zk-\frac{7}{11}ab^{2}e-\frac{4}{55}auf\\
&\quad-\frac{9}{605}a^{3}j
\end{aligned}
\]

\[
\begin{aligned}
qk&=-\frac{3669}{2200}a^{2}b^{2}z+\frac{61747}{1936000}a^{4}bu+\frac{9}{440}a^{2}uq\\
&\quad-\frac{9883}{4400}a^{2}bk+t_\flat+\frac{13}{44}abue\\
&\quad-\frac{1}{2}ab^{2}f-\frac{8}{11}azj
\end{aligned}
\]

\[
\begin{aligned}
qe&=\frac{8}{11}az^{2}-\frac{183}{220}ab^{2}u+\frac{73321}{968000}a^{5}b\\
&\quad+\frac{7}{220}a^{3}q-\frac{10283}{4400}a^{2}be-zf\\
&\quad+\frac{4}{11}uj+2b\ell
\end{aligned}
\]

\[
\begin{aligned}
qf&=\frac{3}{1210}abu^{2}-\frac{7}{11}ab^{4}-\frac{57}{2420}a^{3}bz\\
&\quad+\frac{553}{266200}a^{5}u-\frac{9}{605}a^{3}k-\frac{3}{11}bze\\
&\quad-\frac{16}{605}a^{2}ue-\frac{9283}{4400}a^{2}bf+\frac{14}{11}b^{2}j\\
&\quad-\frac{4}{11}u\ell
\end{aligned}
\]

\[
\begin{aligned}
qj&=-\frac{1}{40}abuz+\frac{49}{6050}a^{3}u^{2}+\frac{21}{880}a^{3}b^{3}\\
&\quad-\frac{1}{2420}a^{5}z-\frac{4}{55}auk-\frac{1}{11}u^{2}e\\
&\quad+\frac{7}{220}a^{2}ze-\frac{1}{2}buf-\frac{9}{4400}a^{4}f\\
&\quad-\frac{9283}{4400}a^{2}bj+z\ell
\end{aligned}
\]

\[
\begin{aligned}
q\ell&=-\frac{64}{605}au^{2}z-\frac{5}{11}ab^{3}z-\frac{9}{605}a^{3}z^{2}\\
&\quad+\frac{356201}{1936000}a^{3}b^{2}u-\frac{2039}{2129600}a^{7}b+\frac{27}{440}abuq\\
&\quad-\frac{1}{2420}a^{5}q+\frac{7}{11}ab^{2}k+\frac{8}{11}z^{2}e\\
&\quad-\frac{7}{22}b^{2}ue+\frac{107}{48400}a^{4}be+\frac{1}{11}u^{2}f\\
&\quad-\frac{7}{220}a^{2}zf+\frac{3}{11}bzj+\frac{16}{605}a^{2}uj\\
&\quad-\frac{10283}{4400}a^{2}b\ell
\end{aligned}
\]

\[
\begin{aligned}
k^2&=-\frac{3}{44}bu^{3}+\frac{1}{4}b^{4}u+\frac{4199}{2200}a^{2}buz\\
&\quad+\frac{323}{38720}a^{4}u^{2}-\frac{8617}{44000}a^{4}b^{3}-\frac{1}{2420}a^{6}z\\
&\quad+uzq-\frac{7}{80}a^{2}b^{2}q+\frac{17}{220}a^{2}uk\\
&\quad-\frac{2}{11}au^{2}e+\frac{1}{4}ab^{3}e+\frac{7}{220}a^{3}ze\\
&\quad-\frac{3}{20}abuf-\frac{9}{4400}a^{5}f+\frac{21}{880}a^{3}bj\\
&\quad+az\ell
\end{aligned}
\]

\[
\begin{aligned}
ke&=\frac{163}{440}abu^{2}-\frac{1}{2}ab^{4}-\frac{1141}{550}a^{3}bz\\
&\quad-\frac{511}{96800}a^{5}u-azq+\frac{7}{220}a^{3}k\\
&\quad-bze+\frac{57}{440}a^{2}ue-\frac{3}{80}a^{2}bf\\
&\quad+b^{2}j-u\ell
\end{aligned}
\]

\[
\begin{aligned}
kf&=-\frac{8}{121}au^{3}+\frac{3}{11}abz^{2}-\frac{17}{440}ab^{3}u\\
&\quad+\frac{7}{605}a^{3}uz+\frac{211963}{3872000}a^{5}b^{2}+\frac{21}{880}a^{3}bq\\
&\quad+\frac{4}{11}uze-\frac{79}{880}a^{2}b^{2}e+bzf\\
&\quad+\frac{5}{88}a^{2}uf+\frac{3}{11}buj-\frac{9}{605}a^{4}j\\
&\quad-b^{2}\ell
\end{aligned}
\]

\[
\begin{aligned}
kj&=\frac{1}{5}au^{2}z-\frac{1}{4}ab^{3}z-\frac{7}{220}a^{3}z^{2}\\
&\quad-\frac{62083}{176000}a^{3}b^{2}u-\frac{92847}{19360000}a^{7}b-\frac{3}{20}abuq\\
&\quad-\frac{9}{4400}a^{5}q-z^{2}e+\frac{1}{4}b^{2}ue\\
&\quad+\frac{93}{17600}a^{4}be-\frac{1}{4}u^{2}f+\frac{7}{80}a^{2}zf\\
&\quad+bzj+\frac{5}{88}a^{2}uj+\frac{3}{80}a^{2}b\ell
\end{aligned}
\]

\[
\begin{aligned}
k\ell&=-\frac{8}{11}az^{3}+\frac{39}{40}ab^{2}uz+\frac{68757}{193600}a^{3}bu^{2}\\
&\quad-\frac{9863}{17600}a^{3}b^{4}-\frac{75061}{968000}a^{5}bz+\frac{263}{5324000}a^{7}u\\
&\quad+\frac{2}{11}au^{2}q-\frac{1}{4}ab^{3}q-\frac{7}{220}a^{3}zq\\
&\quad+\frac{27}{440}abuk-\frac{1}{2420}a^{5}k+\frac{3}{44}bu^{2}e\\
&\quad+\frac{259}{880}a^{2}bze-\frac{19}{12100}a^{4}ue+z^{2}f\\
&\quad-\frac{1}{4}b^{2}uf-\frac{93}{17600}a^{4}bf-\frac{4}{11}uzj\\
&\quad+\frac{79}{880}a^{2}b^{2}j-bz\ell+\frac{57}{440}a^{2}u\ell
\end{aligned}
\]

\[
\begin{aligned}
e^2&=b^{2}z+\frac{8093}{4400}a^{2}bu+\frac{1}{2420}a^{6}\\
&\quad+uq-2bk+\frac{7}{220}a^{3}e\\
&\quad+a\ell
\end{aligned}
\]

\[
\begin{aligned}
ef&=-\frac{7}{11}buz-\frac{9}{1210}a^{2}u^{2}+\frac{9673}{4400}a^{2}b^{3}\\
&\quad+\frac{9}{605}a^{4}z+b^{2}q+\frac{4}{11}uk\\
&\quad+\frac{7}{220}a^{3}f-\frac{3}{11}abj
\end{aligned}
\]

\[
\begin{aligned}
ej&=-\frac{1}{11}u^{3}+bz^{2}-\frac{1}{2}b^{3}u\\
&\quad+\frac{23}{440}a^{2}uz-\frac{30369}{352000}a^{4}b^{2}-\frac{3}{80}a^{2}bq\\
&\quad-zk+\frac{1}{20}auf+\frac{7}{220}a^{3}j
\end{aligned}
\]

\[
\begin{aligned}
e\ell&=-\frac{5}{22}b^{2}u^{2}+\frac{4079}{1100}a^{2}b^{2}z-\frac{31897}{242000}a^{4}bu\\
&\quad+\frac{1}{242000}a^{8}+bzq-\frac{4}{55}a^{2}uq\\
&\quad+\frac{5}{22}a^{2}bk-t_\flat-\frac{157}{440}abue\\
&\quad+\frac{1}{2420}a^{5}e+\frac{1}{2}ab^{2}f+\frac{8}{11}azj
\end{aligned}
\]

\[
\begin{aligned}
f^2&=\frac{16}{121}u^{2}z+b^{3}z+\frac{24719}{48400}a^{2}b^{2}u\\
&\quad-\frac{92847}{2662000}a^{6}b+\frac{3}{11}buq-\frac{9}{605}a^{4}q\\
&\quad+\frac{14}{11}b^{2}k+\frac{93}{2420}a^{3}be+\frac{8}{11}azf\\
&\quad+\frac{32}{121}auj+\frac{3}{11}ab\ell
\end{aligned}
\]

\[
\begin{aligned}
fj&=\frac{4}{11}uz^{2}-\frac{5}{22}b^{2}u^{2}+bzq\\
&\quad+t_\flat
\end{aligned}
\]

\[
\begin{aligned}
f\ell&=\frac{3}{121}bu^{3}+b^{2}z^{2}-\frac{7}{22}b^{4}u\\
&\quad-\frac{17471}{24200}a^{2}buz-\frac{133}{33275}a^{4}u^{2}+\frac{808317}{3872000}a^{4}b^{3}\\
&\quad-\frac{4}{11}uzq+\frac{79}{880}a^{2}b^{2}q+\frac{3}{11}bzk\\
&\quad+\frac{16}{605}a^{2}uk+\frac{8}{121}au^{2}e-\frac{7}{22}ab^{3}e\\
&\quad+\frac{9}{605}a^{3}ze+\frac{157}{440}abuf-\frac{1}{2420}a^{5}f\\
&\quad-\frac{93}{2420}a^{3}bj
\end{aligned}
\]

\[
\begin{aligned}
j^2&=z^{3}-\frac{3}{4}b^{2}uz-\frac{9397}{17600}a^{2}bu^{2}\\
&\quad+\frac{6271}{32000}a^{4}bz+\frac{553}{1936000}a^{6}u-\frac{1}{4}u^{2}q\\
&\quad+\frac{7}{80}a^{2}zq-\frac{1}{2}buk-\frac{9}{4400}a^{4}k\\
&\quad-\frac{1}{4}abze-\frac{1}{275}a^{3}ue+\frac{3}{160}a^{3}bf\\
&\quad+\frac{1}{2}ab^{2}j-\frac{1}{20}au\ell
\end{aligned}
\]

\[
\begin{aligned}
j\ell&=\frac{1}{44}bu^{2}z-\frac{41}{4840}a^{2}u^{3}+\frac{1091}{550}a^{2}bz^{2}\\
&\quad-\frac{10353}{17600}a^{2}b^{3}u+\frac{623}{96800}a^{4}uz-\frac{962739}{77440000}a^{6}b^{2}\\
&\quad+z^{2}q-\frac{1}{4}b^{2}uq-\frac{93}{17600}a^{4}bq\\
&\quad+\frac{1}{11}u^{2}k-\frac{7}{220}a^{2}zk-\frac{6}{55}auze\\
&\quad+\frac{21}{1760}a^{3}b^{2}e+\frac{1}{4}abzf+\frac{1}{275}a^{3}uf\\
&\quad+\frac{157}{440}abuj-\frac{1}{2420}a^{5}j
\end{aligned}
\]

\[
\begin{aligned}
\ell^2&=-\frac{4}{121}u^{3}z+bz^{3}-\frac{17}{44}b^{3}uz\\
&\quad-\frac{7}{44}a^{2}b^{5}-\frac{2300377}{3872000}a^{4}b^{2}z+\frac{103929}{106480000}a^{6}bu\\
&\quad+\frac{23}{585640000}a^{10}-\frac{3}{44}bu^{2}q-\frac{259}{880}a^{2}bzq\\
&\quad+\frac{19}{12100}a^{4}uq-\frac{8}{11}z^{2}k+\frac{7}{22}b^{2}uk\\
&\quad-\frac{107}{48400}a^{4}bk-\frac{31}{44}ab^{2}ze+\frac{359}{6050}a^{3}bue\\
&\quad+\frac{1}{242000}a^{7}e+\frac{6}{55}auzf-\frac{21}{1760}a^{3}b^{2}f\\
&\quad-\frac{8}{121}au^{2}j+\frac{7}{22}ab^{3}j-\frac{9}{605}a^{3}zj
\end{aligned}
\]

\endgroup
\paragraph{The quadratic field and the full anti module}
Let $v_0$ denote the first anti generator (called $v$ in the quartic construction and $w$ in the quintic construction).  The following polynomials $X_d$ are its products with the twelve anti generators.  Thus $X_{\deg v_0}=F=v_0^2$, and the complete module is generated inside $\operatorname{Frac}(S)(v_0)$ by $n_d=X_d/v_0$.  Multiplication is unambiguously $(X_d/v_0)(X_e/v_0)=X_dX_e/F$.  The polynomial product and syzygy identities are verified by the exact verifier; the proof shows that this is the full colon module, not a chosen submodule.
\begingroup\small
\[
\begin{aligned}
v_0w&=b^{2}u+\frac{2207}{4400}a^{4}b+\frac{1}{4}a^{2}q\\
&\quad+abe
\end{aligned}
\]

\[
\begin{aligned}
v_0p&=-\frac{71}{440}a^{2}bu-\frac{1}{48400}a^{6}-bk\\
&\quad-\frac{1}{220}a^{3}e+\frac{1}{2}a\ell
\end{aligned}
\]

\[
\begin{aligned}
v_0v&=-bu^{2}+\frac{1}{2}a^{2}bz+\frac{61}{1760}a^{4}u\\
&\quad-\frac{1}{4}a^{2}k-\frac{1}{2}aue
\end{aligned}
\]

\[
\begin{aligned}
v_0y&=-buz+\frac{1}{44}a^{2}u^{2}+\frac{1}{4}a^{2}b^{3}\\
&\quad+\frac{1}{440}a^{4}z-\frac{1}{2}aze+\frac{3}{160}a^{3}f\\
&\quad-\frac{1}{2}abj
\end{aligned}
\]

\[
\begin{aligned}
v_0r_{21}&=-ab^{2}z-\frac{1091}{1100}a^{3}bu+\frac{1}{96800}a^{7}\\
&\quad-\frac{1}{2}auq+abk-bue\\
&\quad+\frac{1}{440}a^{4}e-\frac{1}{4}a^{2}\ell
\end{aligned}
\]

\[
\begin{aligned}
v_0K_{23}&=-\frac{1}{110}abu^{2}+\frac{1}{2}ab^{4}+\frac{207}{200}a^{3}bz\\
&\quad+\frac{183}{96800}a^{5}u+\frac{1}{2}azq-\frac{3}{220}a^{3}k\\
&\quad+bze-\frac{3}{110}a^{2}ue+\frac{3}{80}a^{2}bf\\
&\quad-b^{2}j
\end{aligned}
\]

\[
\begin{aligned}
v_0r_{24}&=-\frac{4}{11}a^{2}z^{2}-\frac{81}{50}a^{2}b^{2}u-\frac{8427}{176000}a^{6}b\\
&\quad-buq-\frac{9}{440}a^{4}q+\frac{123}{1760}a^{3}be\\
&\quad+\frac{1}{2}azf-\frac{2}{11}auj-\frac{3}{2}ab\ell
\end{aligned}
\]

\[
\begin{aligned}
v_0r_{25}&=\frac{7}{22}abuz+\frac{69}{2420}a^{3}u^{2}-\frac{8623}{4400}a^{3}b^{3}\\
&\quad-\frac{1}{1210}a^{5}z-ab^{2}q-\frac{2}{11}auk\\
&\quad-\frac{1}{11}a^{2}ze-buf-\frac{9}{44}a^{2}bj
\end{aligned}
\]

\[
\begin{aligned}
v_0r_{27}&=\frac{1}{22}au^{3}-abz^{2}+\frac{1}{2}ab^{3}u\\
&\quad-\frac{57}{880}a^{3}uz+\frac{1}{2}azk+\frac{3}{80}a^{2}uf\\
&\quad-buj
\end{aligned}
\]

\[
\begin{aligned}
v_0M_{27}&=-\frac{19}{22}abz^{2}+\frac{29}{88}ab^{3}u-\frac{103}{2420}a^{3}uz\\
&\quad-\frac{95363}{1548800}a^{5}b^{2}-\frac{9}{352}a^{3}bq+\frac{2}{11}azk\\
&\quad+\frac{25}{176}a^{2}b^{2}e+bzf+\frac{1}{44}a^{2}uf\\
&\quad+\frac{9}{1210}a^{4}j-b^{2}\ell
\end{aligned}
\]

\[
\begin{aligned}
v_0r_{31}&=\frac{19}{40}ab^{2}u^{2}-\frac{1589}{550}a^{3}b^{2}z+\frac{567}{8000}a^{5}bu\\
&\quad-\frac{1}{5324000}a^{9}-abzq+\frac{9}{220}a^{3}uq\\
&\quad-\frac{23}{220}a^{3}bk+\frac{1}{2}at_\flat+\frac{389}{880}a^{2}bue\\
&\quad-\frac{1}{24200}a^{6}e-\frac{1}{4}a^{2}b^{2}f-\frac{4}{11}a^{2}zj\\
&\quad-bu\ell+\frac{1}{220}a^{4}\ell
\end{aligned}
\]

\[
\begin{aligned}
v_0r_{34}&=-\frac{32}{605}a^{2}u^{2}z+\frac{3169}{2200}a^{2}b^{3}z-\frac{13}{1210}a^{4}z^{2}\\
&\quad+\frac{247347}{3872000}a^{4}b^{2}u-\frac{175271}{212960000}a^{8}b-\frac{133}{880}a^{2}buq\\
&\quad-\frac{9}{24200}a^{6}q+\frac{11283}{4400}a^{2}b^{2}k-bt_\flat\\
&\quad+\frac{4}{11}az^{2}e-\frac{5}{11}ab^{2}ue+\frac{11353}{968000}a^{5}be\\
&\quad+\frac{1}{22}au^{2}f+\frac{1}{2}ab^{3}f-\frac{1}{88}a^{3}zf\\
&\quad+\frac{19}{22}abzj+\frac{7}{605}a^{3}uj-\frac{10363}{8800}a^{3}b\ell
\end{aligned}
\]

\endgroup
\paragraph{The five minimal module relations}
The relations below are relations over the polynomial parameter ring.  They must be retained together with the product equations.  Their unique second syzygy and the Hilbert-series argument in the text prove that they present the entire anti module.
\begingroup\small
\[
\begin{aligned}
0&=2uw+\frac{1}{2}a^{2}p+2bv\\
&\quad+ar_{21}
\end{aligned}
\]

\[
\begin{aligned}
0&=-2zw-\frac{3}{55}a^{2}v-2by\\
&\quad+aK_{23}
\end{aligned}
\]

\[
\begin{aligned}
0&=2zv-2uy+ar_{27}
\end{aligned}
\]

\[
\begin{aligned}
0&=-\frac{1}{2}azp+\frac{3}{55}auv-zr_{21}\\
&\quad-uK_{23}+br_{27}
\end{aligned}
\]

\[
\begin{aligned}
0&=\frac{913}{400}a^{3}bp+\frac{16}{55}ab^{2}v+\frac{12}{11}azy\\
&\quad-\frac{1}{55}a^{3}r_{24}-2zr_{25}-2uM_{27}\\
&\quad+2br_{31}+ar_{34}
\end{aligned}
\]

\endgroup


\subsection{Altair: the global quintic derivation}
\label{app:quintic-derivation46}
The following formulas use the same generators as
Appendix~\ref{app:quintic-exact-coefficients}. They specify a derivation
of degree six on the entire normalized algebra. The value on the
auxiliary basis element $t_\flat$ follows by differentiating its
explicit expression in $qk$; thus there is no unlisted algebra generator.
\begingroup\small
\begin{align*}
\mathcal Da&=0,\qquad\mathcal Db=8w,\qquad\mathcal Du=-12v,\qquad\mathcal Dz=-14y,\\
\mathcal Dq&=-10423a^2w/550-18r_{24},\\
\mathcal Dk&=-33aK_{23}/2+9a^2v/44-22pu,\\
\mathcal De&=6ap-15r_{21},\\
\mathcal Df&=26ay/11-19r_{25},\\
\mathcal Dj&=-abw-21r_{27},\\
\mathcal D\ell&=-197a^3p/440-309abv/110+481auw/110-25r_{31},\\
\mathcal Dw&=-283a^4/880+17ae/2+14bu,\\
\mathcal Dp&=59a^2u/220-16k,\\
\mathcal Dv&=39a^2z/4-18u^2,\\
\mathcal Dy&=49a^2b^2/8-29aj/2-20uz,\\
\mathcal Dr_{21}&=173a^3u/220-31abz/2+6ak-21eu,\\
\mathcal DK_{23}&=-61a^3z/440-69au^2/55-23bj+23ez\\ &\quad -9a^2f/40+37ab^3/4,\\
\mathcal Dr_{24}&=-27a^6/24200-17a^3e/440-37a\ell/2-24qu\\ &\quad -207197a^2bu/4400,\\
\mathcal Dr_{25}&=-335397a^3b^2/8800-63a^2j/11-39abq/2-25fu\\ &\quad +26auz/11,\\
\mathcal Dr_{27}&=-31069a^5b/176000-3a^3q/40-43az^2/2-27ju\\ &\quad +a^2be/8+47ab^2u/4,\\
\mathcal DM_{27}&=-1073657a^5b/1936000-130az^2/11-27b\ell+27fz
                  \\ &\quad -9a^3q/40+589a^2be/440+521ab^2u/88,\\
\mathcal Dr_{31}&=-8807a^5u/48400-465573a^3bz/8800+17a^3k/110-31\ell u
                  \\ &\quad -9a^2bf/8+919a^2eu/220+4949abu^2/440-51aqz/2,\\
\mathcal Dr_{34}&=10423a^8/26620000+10423a^5e/121000-10423a^3\ell/1100
                 -34kq\\
                &\quad -171a^2b^2z/440-8587a^2bk/440-46a^2qu/55\\
                &\quad -9ab^2f/4+623723a^4bu/968000-63abeu/44.
\end{align*}
\endgroup
To verify the formulas, differentiate the twenty-one invariant products
and use the product table to reduce. For each of the twelve numerators
$x_g=wg$, check $w\mathcal D x_g-x_g\mathcal Dw=F\mathcal Dg$.
These tests check regularity at the branch locus as well as the Leibniz
rule; they suffice because the numerators generate the complete colon
module in all degrees. The exact verifier performs them over $\mathbb Q$.


\section{Finite verification specifications}\label{app:verification}

\subsection{Exact computations and reproducibility}
\label{app:reproducibility46}
All finite sets used in the proof have explicit bounds in the text.
The source package supplies the following independent layers.
Program and data paths below refer to the tree restored by the
computation companion; typesetting the article does not execute it.
The basket generator starts with the binary plurigenus seeds and
performs all prime packings; it takes no excluded-row list as input.
A second traversal verifies the complete equality census and its
plurigenera. The Weil-character and torsion routines enumerate
necessary local states, with multinomial weights when repeated
constituents are grouped. Actual multiplication, complete restriction
and mixed intersections are applied only in the geometric situations
proved in the text. The marked-curve routines enumerate all primitive
fractions in the stated index and degree bounds. The finite-algebra
routines check complete multiplication tables, associators, defining
ideals and homogeneous deformation matrices over exact fields.

An arithmetic output is used only for the bounded algebraic statement
it checks. In particular passing a basket or Hilbert-function test does
not imply existence, and a virtual character assignment is not treated
as an actual torsion class. The constructions have separate normality,
singularity, canonical-module and cohomology proofs. The executable
README specifies the input and output of every principal check and the
corresponding manuscript assertion. All displayed rational comparisons
can also be checked directly from the formulas in the paper.

The computation companion has five source files: the entry point
\path{anc/verification.py} and the four modules \texttt{census.py},
\texttt{exclusions.py}, \texttt{models.py} and \texttt{stability.py}.
Their data are collected in \path{anc/verification_data.json}.
The accompanying README lists the checks and their dependencies.
The command \texttt{list} shows the
available programs, \texttt{extract} restores their source and data,
and \texttt{quick} runs a specified selection in a temporary directory.
The degree-seven conductor check enumerates line-bundle degrees within
the proved bounds and verifies the degree-three and degree-six
cohomology, together with the saturated multiplication-degree
contradiction. All these conductor checks use integer monomial counts
in $\mathbb C[X_5,Y_{12}]$.


\subsection{Exact data for Altair's pluricanonical maps}
\label{app:stability-data}
The data in this subsection give a finite verification of the
function-field recovery in Theorem~\ref{thm:altair-stability}. They use the
same full canonical algebra as the construction of Altair.

\paragraph{Generators on the chart $a=1$.}
Besides $b,u,z,e,w,p,v,y$, every generator is eliminated by the following
identities. Here $t=qk$ is the degree-forty product, rather than the
adjusted degree-forty basis vector used in the homogeneous invariant
table.
\begin{align*}
q={}&-4b^2u-4be-\frac{2207}{1100}b+4w^2,\\
k={}&-4bu^2+2bz-2eu+\frac{61}{440}u-4vw,\\
j={}&-\frac1{16}b^2-4u^3+3uz+4v^2,\\
f={}&-15b^3-\frac{320}{3}bu^3+\frac{400}{3}buz
       +\frac{320}{3}bv^2+\frac{80}{3}ez
       -\frac{40}{33}u^2+\frac{160}{3}wy-\frac4{33}z,\\
l={}&-8b^2u^2+4b^2z-4beu+\frac35bu-8bvw
       +\frac1{110}e+2pw+\frac1{24200},\\
r_{21}={}&-2bv-\frac12p-2uw,\qquad
K_{23}=2by+\frac3{55}v+2wz,\\
r_{24}={}&4b^2v+3bp+4buw-2ew-\frac1{55}w,\\
r_{25}={}&-10b^2w+80bvz+40ey+40pz
       +\frac{20}{11}uv+80uwz-\frac{14}{11}y,\\
r_{27}={}&2uy-2vz,\\
M_{27}={}&-8b^3v-6b^2p+12bew-\frac{291}{440}bw
       -\frac8{11}uy+\frac8{11}vz-8w^3.
\end{align*}
The two remaining module generators satisfy $r_{31}=bK_{23}-pe$ and
$r_{34}=pq$. Each identity is checked in the full multiplication table,
so these eliminations concern the complete chart ring.

\paragraph{A finite construction of the coefficient matrices.}
Order the nine sections $\mathcal H$ as in
\eqref{eq:stability-target}, and order each basis $\mathcal S_m$ as in
the proof of Theorem~\ref{thm:altair-stability}. For
$(m,k)=(23,4),(24,2),(25,3),(26,2),(27,2)$, list all degree-$k$
monomials in $\mathcal S_m$ in the order of nondecreasing tuples of
indices, starting at zero. Multiply each by every $h_i\in\mathcal H$,
with the $h_i$ index first. On $a=1$, reduce the products in the free
$\mathbb Q[b,u,z]$ basis
\[
 1,q,k,t,e,f,j,l,w,p,v,y,r_{24},r_{25},M_{27},r_{31}.
\]
Order the polynomial coordinates lexicographically by the basis index
and then the exponents of $(b,u,z)$. The resulting columns form a
rational matrix $C_m$. A vector in its kernel is precisely a homogeneous
linear relation among the $h_i$, with coefficients of degree $k$ in
$\mathcal S_m$.

For a deterministic construction, use the reduced row echelon form of
$C_m$ and its kernel basis indexed by ascending free columns, each with
free coefficient one. The eight kernel positions used here are
\begin{center}\small
\begin{tabular}{@{}ccl@{}}
\toprule
$m$ & $k$ & Kernel positions, starting at zero\\
\midrule
23&4&$0,35,36,37,38,39,55,90$\\
24&2&$0,7,8,9,10,27,45,49$\\
25&3&$0,21,57,58,59,60,69,80$\\
26&2&$0,1,2,9,13,14,25,32$\\
27&2&$0,8,9,20,35,36,48,59$\\
\bottomrule
\end{tabular}
\end{center}
The companion files contain the complete rational coefficients of these
forty relations. Verification requires only multiplying each relation
out in the defining algebra and obtaining zero, then checking one
minor for each degree. It does not require accepting the output of a
kernel solver as a geometric premise.

\paragraph{Nonvanishing of the five minors.}
Reduce modulo $1009$, set $(a,b,u,z)=(1,843,963,118)$, and evaluate the
ordered basis above at
\[
 (1,470,846,74,34,445,389,610,
     309,164,555,981,733,785,304,629).
\]
All $256$ ordered basis products respect this evaluation. The
multiplication coefficients and the forty rational identities have
denominators prime to $1009$. Delete the first column from each of the
five evaluated $8\times9$ coefficient matrices. In increasing order of
$m$, their determinants are
\[
 368,\quad307,\quad854,\quad526,\quad717\pmod{1009}.
\]
The chart algebra is free over its polynomial parameter ring. Hence a
zero characteristic-zero normal form would reduce to zero under this
evaluation after inverting the denominators. The displayed values prove
that all five minors are nonzero in characteristic zero.

The verifier \texttt{stability/verify\_stability.py} first compares its
presentation, over $\mathbb Q$, with all invariant products, quadratic
products, mixed products and module relations in the defining data for
Altair. It then identifies the entire chart multiplication table, checks
the forty identities and the five minors, and checks the dimensions of
the canonical pieces. Only integer and rational arithmetic is used.


\subsection{Exact data for Vega's pluricanonical maps}
\label{app:vega-stability-data}
\label{ag:a.1.-the-four-recovery-certificates-for-vega}
The computation archive contains the complete rational coefficients
of the four recovery certificates, their construction inputs and
portable verifiers. They use the full canonical algebra of Vega
constructed in Part~IV and recorded in Appendix~\ref{app:coefficients}.

For clarity we restore degree subscripts on the generators. Bases of the
four section spaces used in Theorem~\ref{ag:thm:vega-stability} are \begin{equation}
\begin{aligned}
R_{23}:&\quad a^2h_{15},\ k_{23},\ a^2v_{15},\ az_{19},\ e_{23};\\
R_{24}:&\quad b^2,\ uc,\ au^2,\ a^3b,\ a^6,\ ar_{20},\ ad_{20};\\
R_{25}:&\quad uh_{15},\ as_{21},\ uv_{15},\ a^2w_{17},\ ap_{21},\ f_{25};\\
R_{26}:&\quad bc,\ aub,\ a^3c,\ a^4u,\ a^2q_{18},\ \ell_{26},\ g_{26}.
\end{aligned}
\label{ag:eq:A.1}
\end{equation} Here generator names such as \(e_{23}\) and \(g_{26}\) refer to
the full algebra of Subsection~\ref{sec:quartic-module46}; the generator
$e_{23}$ is distinct from the first nonzero canonical degree $e$.

There is a deterministic way to reconstruct each certificate. Take all
monomials of degree \(k\) in the chosen basis of \(R_m\), multiply each
by the six target functions of \eqref{ag:eq:3.4}, and reduce using the
complete multiplication table on the free chart \eqref{ag:eq:3.3}. Collect
coefficients of the parameter monomials and chart basis elements.
Compute the kernel of this rational matrix. Five suitable kernel vectors
give \begin{equation}
\sum_{j=0}^{5}Q_{ij}(S_m)H_j=0,
\qquad 1\leq i\leq5,
\label{ag:eq:A.2}
\end{equation} with \((m,k)=(23,6),(24,3),(25,3),(26,3)\). The accompanying verifier
substitutes every stored vector into the original table and checks that
the resulting rational coefficient vector is zero. It does not infer a
characteristic-zero identity from a modular nullspace.

To prove that the recovery minor is nonzero, specialize modulo \(1009\)
at \[
(u,b,c)=(3,5,7),
\] with the chart basis \eqref{ag:eq:3.3} evaluated at \[
(1,431,175,908,760,883,442,215,
89,304,437,107,14,926,292,404).
\] The six targets evaluate to \((1,5,9,21,315,51)\). The resulting
five-by-six matrices, with columns ordered as \(H_0,\ldots,H_5\), are \begin{equation}
Q_{23}=
\begin{pmatrix}
257&667&722&0&0&0\\
639&578&580&208&0&0\\
954&683&890&900&0&0\\
394&350&436&706&65&0\\
0&0&0&0&760&707
\end{pmatrix},
\label{ag:eq:A.3}
\end{equation} \begin{equation}
Q_{24}=
\begin{pmatrix}
490&911&0&0&0&0\\
0&580&911&0&0&0\\
0&0&8&429&0&0\\
692&336&768&0&246&0\\
868&416&876&84&94&221
\end{pmatrix},
\label{ag:eq:A.4}
\end{equation} \begin{equation}
Q_{25}=
\begin{pmatrix}
808&412&50&851&0&0\\
577&509&312&102&0&0\\
511&143&454&708&0&0\\
0&0&0&0&65&14\\
304&442&111&131&0&352
\end{pmatrix},
\label{ag:eq:A.5}
\end{equation} \begin{equation}
Q_{26}=
\begin{pmatrix}
497&506&0&0&0&0\\
167&96&825&0&0&0\\
0&0&497&796&0&0\\
80&317&588&337&311&0\\
467&676&654&326&131&132
\end{pmatrix}.
\label{ag:eq:A.6}
\end{equation} Deleting the first column gives determinants \(559,116,354,434\),
respectively. All denominators are invertible at the specialization.
These matrices witness nonvanishing; the complete polynomials
\(Q_{ij}\), rather than only their evaluations, are included with the
verifier.

The verification package checks the exact recovery identities in
degrees \(23,24,25,26\). These four certificates suffice, since
multiplication by \(a_4\) propagates birationality in steps of four.

\section*{AI disclosure}
OpenAI's GPT models were used throughout the mathematical development
and preparation of this paper, under the authors' direction. Their work
included developing and revising the equality classification and moduli
arguments; the explicit constructions, recognition and rigidity
calculations for Vega and Altair; the simplification of basket exclusions;
and the proofs of canonical stability. GPT also produced and revised the
computational programs, organized the accompanying data, and drafted and
edited the manuscript.

The supplied code and data allow the finite computations to be repeated
without further use of GPT. The checks use exact arithmetic to verify
bounded enumerations, algebraic identities, and rank or nonvanishing
calculations. Their geometric applicability and the hypotheses under
which they are used are justified in the mathematical proofs.

The authors are responsible for the final paper in its entirety,
including its arguments, the interpretation of the computations,
and the correctness and attribution of its references.

\section*{Acknowledgements}
We thank our supervisor, Meng Chen, for his guidance in explicit
birational geometry and for introducing us to the canonical volume
problem. His teaching on pluricanonical maps and canonical volumes
provided the foundation for this research.

\clearpage
\phantomsection

\end{document}